\documentclass[12pt, twoside, a4paper, leqno]{article}

\usepackage[all]{xy}
\usepackage[english]{babel}
\usepackage[T1]{fontenc}
\usepackage{amsthm, amsmath, amssymb, amsfonts, mathtools, graphicx}
\usepackage[hidelinks]{hyperref}
\usepackage{bookmark}
\usepackage{faktor}
\usepackage{needspace} 
\usepackage{mathrsfs} 
\usepackage[toc,page]{appendix}
\usepackage{accents} 
\usepackage{changepage} 
\usepackage{emptypage} 

\usepackage[split, makeindex, useindex]{splitidx}

\usepackage[columns=2]{idxlayout}

\usepackage{tikz}
\usetikzlibrary{tikzmark}

\usepackage{tocloft}

\usepackage[top=1in, bottom=1.5in, inner=1.5in, outer=1in]{geometry}

\newindex[Index of notation]{n}
\newindex[Index of terminology]{idx}

\let\standardsection\section%
\def\section{\cleardoublepage\standardsection}

\makeatletter
\newtheorem*{rep@thm}{\rep@title}
\newcommand{\newrepthm}[2]{%
\newenvironment{rep#1}[1]{%
 \def\rep@title{#2~\ref{##1}}%
 \begin{rep@thm}}%
 {\end{rep@thm}}}
\makeatother

\newtheorem{thm}{Theorem}[subsection]
\newrepthm{thm}{Theorem}
\newtheorem{lem}[thm]{Lemma}
\newtheorem{prop}[thm]{Proposition}
\newtheorem{cor}[thm]{Corollary}
\newtheorem{clm}[thm]{Claim}

\newtheorem{thmSec}{Theorem}[section]
\newtheorem{lemSec}[thmSec]{Lemma}
\newtheorem{propSec}[thmSec]{Proposition}
\newtheorem{corSec}[thmSec]{Corollary}
\newtheorem{clmSec}[thmSec]{Claim}

\theoremstyle{definition}
\newtheorem{dfn}[thm]{Definition}
\newtheorem{rem}[thm]{Remark}

\newtheorem{remSec}[thmSec]{Remark}

\renewcommand{\leq}{\leqslant}
\renewcommand{\geq}{\geqslant}
\renewcommand{\emptyset}{\varnothing}
\newcommand{\eps}{\varepsilon}
\newcommand{\del}{\partial}
\newcommand{\bs}{\backslash}

\expandafter\mathchardef\expandafter\varphi\number\expandafter\phi\expandafter\relax%
\expandafter\mathchardef\expandafter\phi\number\varphi%

\DeclareMathOperator{\GL}{GL}

\DeclareMathOperator{\SL}{SL}

\DeclareMathOperator{\SO}{SO}
\DeclareMathOperator{\so}{\mathfrak{so}}
\DeclareMathOperator{\SU}{SU}
\DeclareMathOperator{\su}{\mathfrak{su}}
\DeclareMathOperator{\Sp}{Sp}
\DeclareMathOperator{\smallSp}{\mathfrak{sp}}

\DeclareMathOperator{\Spin}{Spin}

\DeclareMathOperator{\Ann}{Ann}

\DeclareMathOperator{\End}{End}
\DeclareMathOperator{\Id}{Id}
\DeclareMathOperator{\R}{Re}
\DeclareMathOperator{\I}{Im}
\DeclareMathOperator{\tr}{tr}

\DeclareMathOperator{\Stab}{Stab}
\DeclareMathOperator{\Int}{Int}

\DeclareMathOperator{\Ad}{Ad}
\DeclareMathOperator{\ad}{ad}
\DeclareMathOperator{\im}{im}

\DeclareMathOperator{\supp}{supp}
\DeclareMathOperator{\Lie}{Lie}
\DeclareMathOperator{\vol}{vol}
\DeclareMathOperator{\spn}{span}
\DeclareMathOperator{\clo}{clo}
\DeclareMathOperator{\HS}{HS}
\DeclareMathOperator{\ext}{ext}
\DeclareMathOperator{\Hom}{Hom}
\DeclareMathOperator{\interior}{int}

\DeclareMathOperator{\diag}{diag}
\DeclareMathOperator{\arccosh}{arccosh} 

\DeclareMathOperator{\Residue}{\text{Res}}
\DeclareMathOperator{\rest}{res}

\newcommand{\C}{\mathit{C}}

\newcommand{\Rcal}{\cal{R}}

\newcommand{\F}{\mathcal{F}}

\newcommand{\Hcal}{\mathcal{H}}

\renewcommand{\O}{\mathcal{O}}
\renewcommand{\P}{\mathcal{P}}
\newcommand{\U}{\mathcal{U}}
\newcommand{\V}{\mathcal{V}}

\renewcommand{\S}{\mathcal{S}}

\newcommand{\CC}{\mathbb{C}}
\newcommand{\FF}{\mathbb{F}}
\newcommand{\HH}{\mathbb{H}}

\newcommand{\NN}{\mathbb{N}}
\newcommand{\OO}{\mathbb{O}}

\newcommand{\RR}{\mathbb{R}}

\newcommand{\ZZ}{\mathbb{Z}}

\renewcommand{\a}{\mathfrak{a}}
\newcommand{\e}{\mathfrak{e}}

\newcommand{\g}{\mathfrak{g}}
\newcommand{\h}{\mathfrak{h}}
\renewcommand{\k}{\mathfrak{k}}
\renewcommand{\l}{\mathfrak{l}}
\newcommand{\m}{\mathfrak{m}}
\newcommand{\n}{\mathfrak{n}}

\newcommand{\s}{\mathfrak{s}}
\newcommand{\tfrak}{\mathfrak{t}}

\newcommand{\Sfrak}{\mathfrak{S}}
\newcommand{\CS}{\mathscr{C}}

\newcommand{\dX}{\del X}
\newcommand{\Cinf}{\C^{\infty}}
\newcommand{\Ccinf}{\C_c^{\infty}}

\newcommand{\LO}[2]{\prescript{}{#2}{#1}}
\newcommand{\LOI}[2]{\prescript{#2}{}{#1}}
\newcommand{\restricted}[2]{\left.#1 \right|_{#2}}

\newcommand{\ubar}[1]{\underaccent{\bar}{#1}}

\newcommand{\Quotient}[2]{
    \begin{array}{|c|c|}
        \hline
        #1 \\ \hline
        #2 \\ \hline
    \end{array}
}

\newcommand{\DoubleQuotient}[4]{
    \begin{array}{|c|c|}
        \hline
        #1 \\ \hline
        #2 \\ \hline
        #3 \\ \hline
        #4 \\ \hline
    \end{array}
}

\newcommand{\ThreeQuotient}[3]{
    \begin{array}{|c|c|}
        \hline
        #1 \\ \hline
        #2 \\ \hline
        #3 \\ \hline
    \end{array}
}

\renewcommand{\t}[1]{ \LOI{ #1 }{t} }

\renewcommand{\c}[1]{\LOI{#1}{\circ}}

\renewcommand{\d}[1]{\ensuremath{\operatorname{d}\!{#1}}}

\newcommand*{\entry}[1]{\multicolumn{1}{|c}{#1}}

\newcommand{\blockmat}[4]{\left(\begin{array}{cc}
    #1  & \entry{#2} \\ \hline
    #3 & \entry{#4}
    \end{array}\right)
}

\def\nl{\needspace{1\baselineskip} \hfill \\}
\def\nlenum{\needspace{1\baselineskip} \hfill}

\numberwithin{equation}{section}

\DeclareFontFamily{OT1}{cmrx}{}
\DeclareFontShape{OT1}{cmrx}{m}{n}{<->cmr10}{}
\let\saveLongrightarrow\Longrightarrow
\makeatletter
\renewcommand*{\Longrightarrow}{%
    \mathrel{\rlap{\fontfamily{cmrx}\fontencoding{OT1}\selectfont=}%
    \hphantom{\saveLongrightarrow}%
    \llap{$\m@th\Rightarrow$}}}
\makeatother

\newcounter{nber}
\newcommand{\SecParNoSkip}{\thesection.\arabic{nber} \stepcounter{nber} }

\newcommand{\SecPar}{\medskip

\needspace{2\baselineskip}
\thesection.\arabic{nber} \stepcounter{nber} }

\title{Cusp forms for locally symmetric spaces of infinite volume}

\author{Gilles BECKER}

\begin{document}

\pagenumbering{Roman}

\begin{adjustwidth*}{-0.5in}{} 
    \let\newpage\relax

    \begin{center}
        \includegraphics[width=2.5cm]{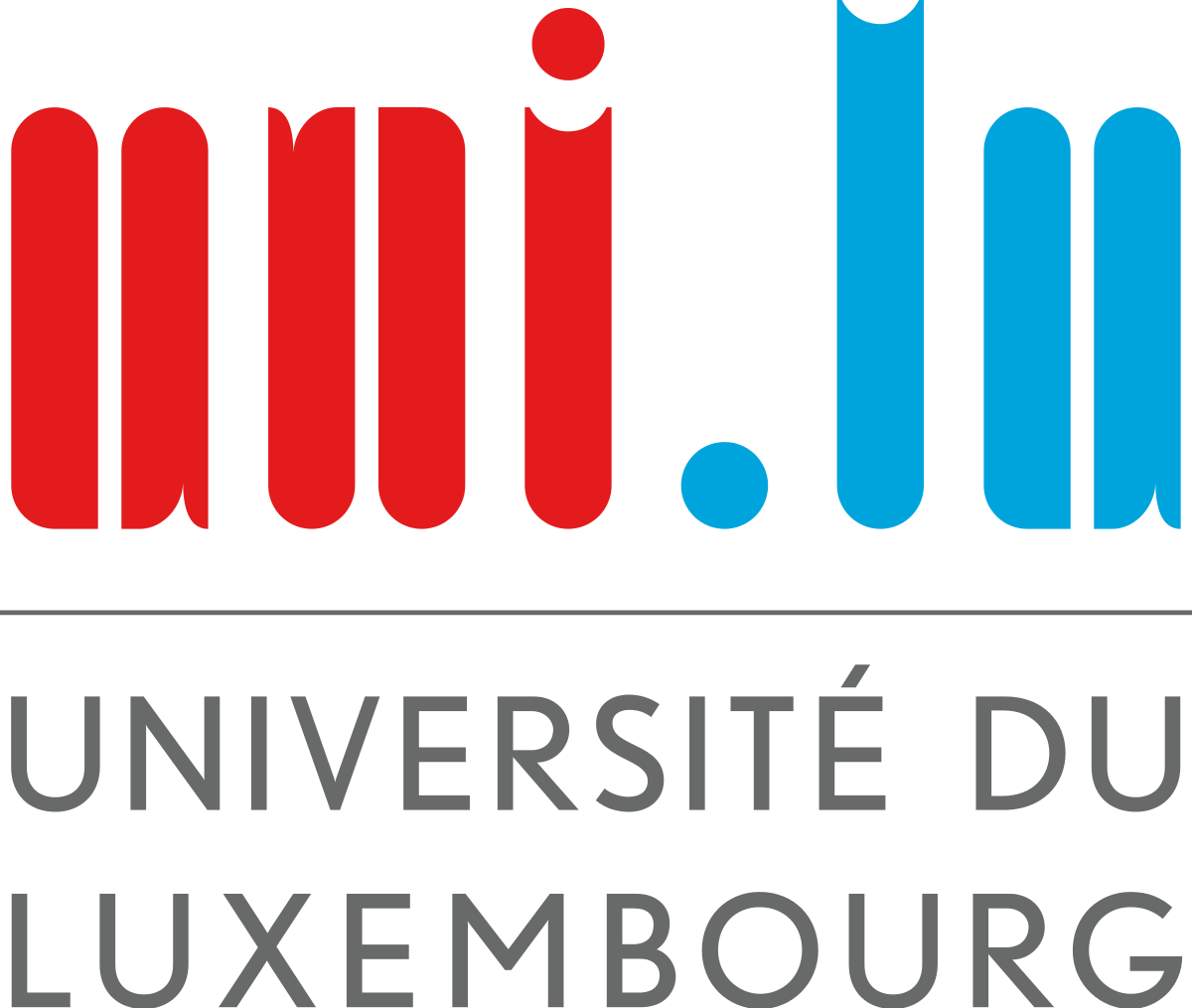}

        \bigskip

        {\small PhD-FSTC-2017-63}

        {\small The Faculty of Science, Technology and Communication}

        \vskip 2cm

        {\bf \Large DISSERTATION}

        \bigskip

        Presented on 24/10/2017 in Luxembourg

        \medskip

        to obtain the degree of

        \vskip 1cm

        \let\origbaselineskip\baselineskip
        \setlength{\baselineskip}{1.25\origbaselineskip}

        {\bf \large DOCTEUR DE L'UNIVERSITÉ DU LUXEMBOURG \\
        EN MATHÉMATIQUES}

        \medskip

        {\large by}

        \medskip

        \makeatletter

        {\large \@author}

        {\small Born on 24 March 1989 in Luxembourg}

        \vskip 1cm

        {\bf \Large \MakeUppercase{\@title}}
    \end{center}

    \vfill

    {\bf \large Dissertation defence committee}

    \smallskip

    Dr Martin OLBRICH, dissertation supervisor \\
    {\it Professor, Université du Luxembourg}

    \medskip

    Dr Ulrich BUNKE \\
    {\it Professor, Universität Regensburg}

    \medskip

    Dr Colin GUILLARMOU \\
    {\it Research Director for CNRS, \\
    Université Paris-Sud (Orsay)}

    \medskip

    Dr Gabor WIESE, chairman \\
    {\it Professor, Université du Luxembourg}

    \medskip

    Dr Salah MEHDI, vice chairman \\
    {\it Professor, Université de Lorraine}

    \thispagestyle{empty}
\end{adjustwidth*}

\cleardoublepage

\pagenumbering{roman}

\begin{abstract}
    \noindent Let $ G $ be a real simple linear connected Lie group of real rank one. Then, $ X := G/K $ is a Riemannian symmetric space with strictly negative sectional curvature. By the classification of these spaces, $X$ is a real/\-complex/\-quaternionic hyperbolic space or the Cayley hyperbolic plane.
    We define the Schwartz space $ \CS(\Gamma \bs G) $ on $ \Gamma \bs G $ for torsion-free geometrically finite subgroups $ \Gamma $ of $G$.
    We show that it has a Fréchet space structure, that the space of compactly supported smooth functions is dense in this space, that it is contained in $ L^2(\Gamma \bs G) $ and that the right translation by elements of $G$ defines a representation on $ \CS(\Gamma \bs G) $.
    Moreover, we define the space of cusp forms $ \c{\CS}(\Gamma \bs G) $ on $ \Gamma \bs G $, which is a geometrically defined subspace of $ \CS(\Gamma \bs G) $.
    It consists of the Schwartz functions which have vanishing ``constant term'' along the ordinary set $ \Omega_\Gamma \subset \dX $ and along every cusp. We show that these two constant terms are in fact related by a limit formula if the cusp is of smaller rank (not of full rank).
    The main result of this thesis consists in proving a direct sum decomposition of the closure of the space of cusp forms in $ L^2(\Gamma \bs G) $ which respects the Plancherel decomposition in the case where $ \Gamma $ is convex-cococompact and noncocompact. For technical reasons, we exclude here that $X$ is the Cayley hyperbolic plane.
\end{abstract}

\cleardoublepage

\tableofcontents

\cleardoublepage

\pagenumbering{arabic}
\setcounter{page}{1}

\section{Introduction}


\SecParNoSkip Let $ G $ be a real simple linear connected Lie group of real rank one. For these Lie groups, the abelian group $A$ in an Iwasawa decomposition $ G = NAK $ is one-dimensional. Moreover, $ X := G/K $ is then a Riemannian symmetric space with strictly negative sectional curvature. By the classification of these spaces, $X$ is a real/complex/quaternionic hyperbolic space or the Cayley hyperbolic plane. For technical reasons, we exclude $ X = \OO H^2 $ (Cayley hyperbolic plane) at several places.

\SecPar Let $ \Gamma $ be a geometrically finite subgroup of $G$.
Roughly speaking, $ \Gamma $ is a torsion-free discrete subgroup of $G$ such that $ \Gamma \bs X $ is a manifold with possibly infinite volume having finitely many cusps (of smaller rank or of full rank, rational or irrational) and funnels as ends. If the ends consist only of funnels, then $ \Gamma $ is said to be convex-cocompact.

Here is an example of the shape of a manifold $ \Gamma \bs X $ ($X$: upper-half plane) of genus 3 having one cusp (of full rank) and one funnel as ends:
\begin{center}
    \begin{tikzpicture}[scale=0.3]
        \draw plot[smooth, tension=0.55] coordinates { (-8,10) (-6, 5) (-4, 2) (-2, 1) (1, 1) (5, 3.5) };

        \draw plot[smooth, tension=0.55] coordinates { (-8.5,9.85) (-8,7) (-8,2) (-10,0) (-10,-3) (-6.5,-5.5) (-1.5,-2) (2,0) (6,0.5) };

        \draw plot[smooth, tension=1] coordinates { (-7,0)  (-6.5,-0.25) (-6,0) };

        \draw plot[smooth, tension=1] coordinates { (-6.9,-0.1)  (-6.5, 0.05) (-6.1,-0.1) };

        \draw plot[smooth, tension=1] coordinates { (-5,-3) (-4.41,-2.93) (-4,-2.5) };

        \draw plot[smooth, tension=1] coordinates { (-4.9,-3.02) (-4.59,-2.57) (-4.05,-2.6) };

        \draw plot[smooth, tension=1] coordinates { (-9,-3) (-8.45,-3.07) (-8,-2.75) };

        \draw plot[smooth, tension=1] coordinates { (-8.9,-3.05) (-8.55,-2.68) (-8.1,-2.87) };

        \draw plot[smooth, tension=1] coordinates { (-7.85,4.5) (-6.87, 4.56) (-6, 5)  };
        \draw[dashed] plot[smooth, tension=1] coordinates { (-7.85,4.5) (-6.98, 4.94) (-6, 5)  };

        \draw plot[smooth, tension=1] coordinates { (1.5,-0.2) (1.07, 0.32) (1, 1) };
        \draw[dashed] plot[smooth, tension=1] coordinates { (1.5,-0.2) (1.43, 0.48) (1, 1) };

        \draw plot[smooth, tension=1] coordinates { (5,3.5) (5.69, 2.06) (6,0.5) };
        \draw plot[smooth, tension=1] coordinates{ (5,3.5) (5.31, 1.94) (6,0.5) };

        \draw plot[smooth, tension=1] coordinates { (-8.5,9.85)  (-8.26, 9.97) (-8,10) };
        \draw plot[smooth, tension=1] coordinates{ (-8.5,9.85) (-8.24, 9.88) (-8,10) };
    \end{tikzpicture}
\end{center}

\SecPar Since the classical Fourier analysis was so successful, one generalised it to groups (not necessarily abelian or compact) and even to homogeneous spaces.
An important goal in harmonic analysis is to establish the Plancherel theorem for $ L^2(\Gamma \bs G) $, providing a ``Fourier transform''. The abstract Plancherel theorem says that we can decompose this space as a direct integral over the unitary dual $ \hat{G} $ of $G$.
Consider now the right regular representation of $G$ on $ L^2(\Gamma \bs G) $. Then, this unitary representation has a direct integral decomposition into irreducibles:
\[
    L^2(\Gamma \bs G) \simeq \int_{\hat{G}}^{\oplus} N_\pi \hat{\otimes} V_\pi \, d\kappa(\pi) \ ,
\]
where $ V_\pi $ is an irreducible unitary representation space of a representation belonging to the class of $ \pi \in \hat{G} $, $ N_\pi $ is a Hilbert space which is called the \textit{multiplicity space} of $ \pi $ in $ L^2(\Gamma \bs G) $ and $ \kappa $ is a Borel measure on $ \hat{G} $ which is called the \textit{Plancherel measure}.

M.~Olbrich and U.~Bunke determined the Plancherel decomposition quite explicitly for $ L^2(\Gamma \bs G) $ in the convex-cocompact case (cf. Theorem 11.1 of \cite[p.155]{BO00}). In the presence of cusps, it is in general not known.
If $ \pi $ is however a discrete series representation, then they showed only that $ N_\pi $, which is equal to the space $ (V_{\pi', -\infty})^\Gamma_d $ of $ \Gamma $-invariant distribution vectors having square-integrable matrix coefficients, is infinite-dimensional.
We will provide a more precise description of $ N_\pi $ (cf. \eqref{eq Intro (V_(pi, -infty))_d^Gamma = overline(D_1) oplus D_2} on p.\pageref{eq Intro (V_(pi, -infty))_d^Gamma = overline(D_1) oplus D_2}).

By Theorem 11.1 of \cite[p.155]{BO00}, $ L^2(\Gamma \bs G) $ decomposes into a continuous and a discrete part:
\[
    L^2(\Gamma \bs G) = L^2(\Gamma \bs G)_{ac} \oplus L^2(\Gamma \bs G)_d \ .
\]
Here, $ L^2(\Gamma \bs G)_{ac} $ is generated by wave packets of matrix coefficients and $ L^2(\Gamma \bs G)_d $ can be further decomposed with the help of representation theory:
\[
    L^2(\Gamma \bs G)_d = L^2(\Gamma \bs G)_{ds} \oplus L^2(\Gamma \bs G)_{res} \oplus L^2(\Gamma \bs G)_U \ ,
\]
where $ L^2(\Gamma \bs G)_{res} $ is generated by residues of Eisenstein series, $ L^2(\Gamma \bs G)_U $ is generated by matrix coefficients of ``stable'' invariant distributions supported on the limit set and $ L^2(\Gamma \bs G)_{ds} $ is generated by square-integrable matrix coefficients of discrete series representations (it is denoted by $ L^2(\Gamma \bs G)_{cusp} $ in \cite{BO00}).

Matrix coefficients give in fact the fundamental relation between representation theory and harmonic analysis.

In the group case ($ \Gamma = \{e\} $), the Plancherel theorem was established by Harish-Chandra in 1976. He finished his work in \cite{HC76}. It was Harish-Chandra who found out that the discrete series representations are the key to representation theory and harmonic analysis on connected semisimple Lie groups with finite center. Moreover, he classified the discrete series representations for these groups.
Selberg, Gel'fand and Langlands proved the Plancherel theorem in the finite volume case ($ \vol(\Gamma \bs G) < \infty $).

\SecPar Let us continue by introducing some notations. Let $ (\gamma, V_\gamma) $ be a finite-dimensional representation of $K$.
We form the homogeneous vector bundle $ V(\gamma) := G \times_K V_\gamma $ over $X$ and the locally homogeneous vector bundle $ V_Y(\gamma) := \Gamma \bs V(\gamma) $ over $ Y := \Gamma \bs X $. M.~Olbrich and U.~Bunke determined in \cite{BO00} the spectral decomposition of the space of sections $ L^2(Y, V_Y(\gamma)) $ over $Y$ with respect to the Casimir operator and other locally invariant differential operators when $ \Gamma $ is convex-cocompact as a consequence of the Plancherel decomposition of $ L^2(\Gamma \bs G) $. Spectral decompositions of $ L^2(Y, V_Y(\gamma)) $, respectively partial results, were already obtained before (sometimes cusps are allowed)
\begin{itemize}
\item by Patterson \cite{Patterson75}, for trivial $ \gamma $ and surfaces,
\item by Lax-Phillips \cite{LaxPhillips82, LaxPhillips84-1, LaxPhillips84-2, LaxPhillips85} and Perry \cite{Perry87} for trivial $ \gamma $ and higher dimensional real hyperbolic manifolds,
\item by Mazzeo-Phillips \cite{MazzeoPhillips90} and Epstein-Melrose-Mendoza \cite{EpsteinMelroseMendoza91} for differential forms on real, respectively complex, hyperbolic manifolds.
\end{itemize}
More recent results in the presence of cusps were obtained by Guillarmou-Mazzeo \cite{Guillarmou} for geometrically finite real hyperbolic manifolds and by Bunke-Olbrich \cite{BO07} for rank-one symmetric spaces.

\SecPar In the group case ($ \Gamma = \{e\} $) and the finite volume case, which are extreme cases of geometrically finite groups, the so-called cusp forms provide us a geometric characterisation of (most of) the discrete part of the direct integral decomposition.
We define the Schwartz space $ \CS(\Gamma \bs G) $ as well as the space of cusp forms $ \c{\CS}(\Gamma \bs G) \subset \CS(\Gamma \bs G) $ for geometrically finite groups (in fact we allow also a twist by a finite-dimensional unitary representation). For the analysis of the space of cusp forms, we restrict us however to the convex-cocompact case.

\SecPar The Schwartz space $ \CS(\Gamma \bs G) $ consists of $ \Gamma $-invariant functions satisfying a certain decay behaviour (involving also derivatives) at infinity, i.e. at the cusps and at the funnels.

In order to be able to define this space, we generalise first the notion of a Siegel set (cf. Section \ref{ssec: Generalised Siegel sets}) so that the following holds:
\\ Let $ \Sfrak' $ be a generalised Siegel set in $X$ at a cusp and let
\[
    \pi \colon X \cup \Omega \to \Gamma \bs (X \cup \Omega), \quad x \mapsto \Gamma x \ .
\]
Then, $ \pi $ restricted to $ \Sfrak' $ is a finite-to-one map and the closure of the image of $ \Sfrak' $ under $ \pi $ represents one end of $ \Gamma \bs (X \cup \Omega) $.

We use then these sets and sets which are ``away from the limit set'' to define some seminorms controlling the behaviour of the functions at the cusps and the funnels. We check that the Schwartz space is indeed a Fréchet space which is contained in $ L^2(\Gamma \bs G) $ and that the compactly supported smooth functions on $ \Gamma \bs G $ form a dense subspace of it (cf. Proposition \ref{Prop Ccinf(Gamma|G) is dense in the geometrically finite case}).

Let us give now the precise definition of the seminorms controlling the behaviour at the funnels: \\
Let $ \dX $ denote the boundary of $X$ and let $ \Lambda := \Lambda_\Gamma \subset \dX $ be the limit set of $ \Gamma $.
Let $ \Omega := \Omega_\Gamma = \dX \smallsetminus \Lambda_\Gamma $ be the set of ordinary points of $ \Gamma $.

Let $ U $ be an open subset of $X$ which is relatively compact in $ X \cup \Omega_\Gamma $, i.e. the closure of $U$ in $ X \cup \Omega_\Gamma $ is compact.

\medskip
\needspace{5\baselineskip}
\underline{Picture for $ G = \SL(2, \RR) $:}

\begin{center}
    \begin{tikzpicture}[scale=0.25]
        \draw[thick]  (0,0) circle (10);

        \foreach \a in {20,..., 40, 70, 71,..., 85, 105, 106,..., 130, 170, 171, ..., 190, 230, 231, ..., 255, 275, 276, ..., 290, 320, 321, ..., 340}
                      \draw[orange] (\a:9.825) -- (\a:10.175);

        \draw[thick, dashed] plot[smooth, tension=.7] coordinates {(15:10) (4,1) (4,-1) (-15:10)};

        \node at (7,0) {$U$};

        \node at (10.5,5) {\textcolor{orange}{$ \Lambda $}};
    \end{tikzpicture}
\end{center}
Each element of $ g \in G $ has a Cartan decomposition $ g = k_g a_g h_g $ with $ k_g, h_g \in K, a_g \in \bar{A}_+ $.

Let $ \U(\g) $ be the universal enveloping algebra of $ \g_\CC $ and let $ \U(\g)_l $ be the vector subspace of $ \U(\g) $ consisting of elements of order at most $l$.

Let $ l \in \NN_0 $. For $ r \geq 0 $, $ X \in \U(\g)_{l_1} $, $ Y \in \U(\g)_{l_2} $ such that $ l_1 + l_2 \leq l $ and a function $ f \in \C^l(G) $, set
\[
    \LOI{p}{U}_{r, X, Y}(f) = \sup_{gK \in U}{ (1 + \log a_g)^r a_g^\rho |L_X R_Y f(g)| } \ .
\]
For the cusps, we define similarly a seminorm $ \LO{p}{\Sfrak'}_{r, X, Y} $ on a generalised Siegel set, controlling the behaviour of the functions at a cusp, with resembling weight functions, which depend on $ \rho $ (the half-sum of positive restricted roots) and the rank of the cusp.

The Schwartz space $ \CS(\Gamma \bs G) $ on $ \Gamma \bs G $ consists then of those functions in $ \Cinf(\Gamma \bs G) $ having finite seminorms.
One can show that this definition is equivalent to Definition 8.1 in \cite[p.129]{BO00} in the convex-cocompact case.

\SecPar Let us in the following define the space of cusp forms on $ \Gamma \bs G $. It is a geometrically defined subspace of the Schwartz space. In the convex-cocompact case, it consists of the Schwartz functions which have vanishing integrals along all the horospheres in $G$ over which a nonzero Schwartz function is integrable.

Fix an Iwasawa decomposition $ G = N A K $.
Let $ M = Z_K(A) $ and $ P = M A N $ (parabolic subgroup of $G$). Then, the boundary $ \dX $ of $X$ is given by $ G/P = K/M $.
All the horospheres in $ G/K $ are of the form $ \HS_{gK, kM} := k N k^{-1} g K $ (horosphere passing through $ gK \in X $ and $ kM \in \dX $).

Let us define now the constant term along $ \Omega $ and the constant term along a cusp.

For $ g \in G(\Omega) := \{g \in G \mid gP \in \Omega\} $, $ h \in G $ and for a measurable function $ f \colon G \to \CC $, we define
\[
    f^\Omega(g, h) = \int_N{ f(gnh) \, dn }
\]
whenever the integral exists. This function is called the \textit{constant term} of $f$ along $ \Omega $.
It is well-defined if $ f $ is a Schwartz function by Proposition \ref{Prop f^Omega well-defined}.

When the geometrically finite group $ \Gamma $ is not convex-cocompact, we define in addition the constant term at a cusp.

Let $ Q $ be a parabolic subgroup of $G$ representing a cusp $ kM \in \dX $ ($ Q = k P k^{-1} $). We call $Q$ $ \Gamma $-cuspidal.

We choose a Langlands decomposition $ Q = M_Q A_Q N_Q $ so that we can construct $ N_{\Gamma, Q} \subset N_Q $ and $ M_{\Gamma, Q} \subset M_Q $ such that
\begin{enumerate}
\item $ \Gamma_Q := \Gamma \cap Q $ is contained in $ N_{\Gamma, Q} M_{\Gamma, Q} $, and
\item $ \Gamma_Q $ is a cocompact lattice in $ N_{\Gamma, Q} M_{\Gamma, Q} $.
\end{enumerate}
For $ M_{\Gamma, Q} $, we can take the smallest possible subset of $ M_Q $.

For $ g \in G $ and for a $ \Gamma $-invariant measurable function $ f \colon \Gamma \bs G \to \CC $, we define
\[
    f^Q(g) = \frac1{\vol(\Gamma_Q \bs N_{\Gamma, Q} M_{\Gamma, Q})} \int_{\Gamma_Q \bs N_Q M_{\Gamma, Q}}{ f(nmg) \, dn \, dm}
\]
whenever the integral exists. We call this function the \textit{constant term} of $f$ along $Q$.
Again, it is well-defined if $ f $ is a Schwartz function (cf. Proposition \ref{Prop f^P well-defined}).

The \textit{space of cusp forms} $ \c{\CS}(\Gamma \bs G) $ on $ \Gamma \bs G $ consists then of the functions $ f \in \CS(\Gamma \bs G) $ for which $ f^\Omega = 0 $ and $ f^Q $ vanishes  for every $\Gamma$-cuspidal parabolic subgroup $ Q $.
In Section \ref{ssec:Examples of cusp forms}, we determine the cusp forms on $ \Gamma \bs G $ which are induced from cusp forms on $G$.

\SecPar We define the little constant term of a Schwartz function along a cusp: \\
For $ g \in G $ and for a $ \Gamma $-invariant measurable function $ f \colon \Gamma \bs G \to \CC $, we define
\[
    f^{Q, lc}(g) = \frac1{ \vol(\Gamma_Q \bs N_{\Gamma, Q} M_{\Gamma, Q}) } \int_{\Gamma_Q \bs N_{\Gamma, Q} M_{\Gamma, Q}}{ f(nmg) \, dn \, dm}
\]
whenever the integral exists. This function is called the \textit{little constant term} of $f$ along $Q$.

The little constant term is equal to the constant term when we consider it along a cusp of full rank. Moreover, the constant term at a cusp of a Schwartz function is zero when the little constant term vanishes at this cusp.

We show in Theorem \ref{Thm f - f^(Q, lc) is rapidly decreasing} that a Schwartz function, with vanishing little constant term, decays rapidly on a (nongeneralised) Siegel set under some assumptions, which are always fulfilled in the real hyperbolic case.

\SecPar We prove that the right translation by elements of $G$ on the Schwartz space, respectively the space of cusp forms, is a representation of $G$ (see Theorem \ref{Thm representations geometrically finite case}). This result also implies that the Schwartz space does not depend on the choices made in its construction (cf. Proposition \ref{Prop Schwartz space is independent of choices}).

\SecPar We show that, in the real hyperbolic case, the constant term of a Schwartz function $f$ along a cusp of smaller rank (having not full rank as the cusps in the finite volume case) is zero if $ f^\Omega = 0 $ (cf. Theorem \ref{Thm f^Omega = 0 implies f^P = 0}). We conjecture that this result is also true beyond the real hyperbolic case. The main result we need in the proof is a limit formula:

Let $ G = \SO(1, n)^0 = N A K $, $ n \geq 2 $, be the identity component of $ \SO(1, n) $. Let $ Q $ be a $ \Gamma $-cuspidal parabolic subgroup with associated cusp of smaller rank. Without loss of generality, $ Q = P $. Let $ w \in K $ be a representative of the Weyl element satisfying $ w^2 = \Id $.

The group $A$ can be identified with $ (0, \infty) $ and the unipotent radical $N$ of $P$ can be identified with $ \RR^{n-1} = \RR^{n-k-1} \times \RR^k $.
The group $ N_{\Gamma, P} $ is identified with $ \{0\} \times \RR^k $. Let $ n_t = (t, 0_{1, n-2}) \in N^{\Gamma, P} \equiv \RR^{n-k-1} \times \{0\} $.

For simplicity, we state the result here only for right $K$-invariant Schwartz functions:
\begin{repthm}{Thm Conv of the constant terms}
    Let $ k = \dim(N_{\Gamma, P}) \geq 1 $ and let $ f \in \CS(\Gamma \bs G)_K $.
    Then,
    \[
        \lim_{t \to \infty} \frac1{t^k} f^\Omega(\kappa(n_t w), a_s)   \qquad (s \in (0, +\infty))
    \]
    is up to a constant equal to
    \[
        \int_0^1 u^{n - \frac{k}2 - 2}(1 - u)^{\frac{k}2 - 1} f^P(a_{\tfrac{s}{u}}) \, du \ .
    \]
\end{repthm}
\SecPar In the finite volume case, the closure of the space of cusp forms is contained in the discrete part of $ L^2(\Gamma \bs G) $ and it has a ``small'' complement (Gel'fand, Selberg, Langlands).
In the group case, Harish-Chandra proved that the closure of the space of cusp forms is equal to the discrete part of $ L^2(G) $.

Let us look now at the case, where $ \Gamma $ is convex-cocompact, noncocompact (then $ \Gamma \bs G $ has infinite volume). Famous examples, for $X$ being the upper-half plane, are:

\bigskip

\begin{minipage}{.5\textwidth}
    \begin{center}
        \begin{tikzpicture}[scale=0.3]
            \draw plot[smooth, tension=.7] coordinates {(-4.5,3.5) (-2, 3.5) (-0.25, 5) (-1,8.3)};

            \draw plot[smooth, tension=.7] coordinates {(3,8) (2,5) (3.5,3) (5.75,2.7)};

            \draw plot[smooth, tension=.7] coordinates {(4,-0.5) (2.9,1.25) (1, 2.5) (-1, 2) (-3,0)};

            \draw plot[smooth cycle, tension=.7] coordinates {(3,8) (1.03,8.55) (-1, 8.3) (0.97,7.75)};

            \draw plot[smooth cycle, tension=.7] coordinates {(5.75,2.7) (4.52,1.29) (4, -0.5) (5.23, 0.91)};

            \draw plot[smooth cycle, tension=.7] coordinates {(-4.5,3.5) (-4.09,1.54) (-3, 0) (-3.41,1.96)};
        \end{tikzpicture}
    \end{center}
\end{minipage}
\begin{minipage}{.4\textwidth}
    \begin{center}
        \begin{tikzpicture}[scale=0.4]
            \draw plot[smooth, tension=0.55] coordinates {(4,1.75) (2,1) (0,0.75) (-2,2) (-3.5,2.35) (-5,2) (-6,1.25)
            (-6.5,0)
            (-6,-1.25) (-5,-2) (-3.5,-2.35)  (-2,-2) (0,-0.75) (2,-1) (4, -1.75)};

            \draw plot[smooth, tension=1] coordinates {(-5,0) (-3.5,-0.9) (-2,0) };

            \draw plot[smooth, tension=1] coordinates {(-4.95,-0.07)(-3.5,0.6) (-2.05,-0.07)};

            \draw plot[smooth cycle, tension=1] coordinates {(4,1.75) (3.75,0) (4, -1.75) (4.25, 0)};
        \end{tikzpicture}
    \end{center}
\end{minipage}

\bigskip

For such $ \Gamma $'s, we show that the closure of the space of cusp forms admits the following decomposition:
\begin{repthm}{Thm decomposition of closure of cusp forms}
    $
        \overline{\c{\CS}(\Gamma \bs G)} = L^2(\Gamma \bs G)_{ds} \oplus L^2(\Gamma \bs G)_U \ .
    $
\end{repthm}
This result was already conjectured by M.~Olbrich in 2002 (see \cite[p.116]{Olb02}).

In order to obtain this result, it was necessary to investigate the different parts appearing in the Plancherel decomposition.

In this decomposition, matrix coefficients of $ \Gamma $-invariant distributions supported on the limit set appear and we determine which matrix coefficients belong to the space of cusp forms (see Theorem \ref{Thm c_(f, v) in 0CS(Gamma|G, phi) iff f in U_Lambda(sigma_lambda, phi)}). From this, we can conclude that $ L^2(\Gamma \bs G)_U $ is contained in $ \overline{\c{\CS}(\Gamma \bs G)} $.

We prove the following statement using the \textit{dual transform} (defined in Section \ref{sssec:The dual transform}):
\begin{repthm}{Thm a Schwartz fct in the closure of the space of cusp forms is a cusp form}
    Let $\gamma \in \hat{K}$. Then,
    \[
        \overline{\c{\CS}(\Gamma \bs G, \phi)}(\gamma) \cap \CS(\Gamma \bs G, \phi) = \c{\CS}(\Gamma \bs G, \phi)(\gamma) \ .
    \]
\end{repthm}
It follows then from this theorem and Theorem \ref{Thm c_(f, v) in 0CS(Gamma|G, phi) iff f in U_Lambda(sigma_lambda, phi)} (cf. Proposition \ref{Prop L^2(Gamma|G, phi)_res is orthogonal to the space of cusp forms}) that the space $ L^2(\Gamma \bs G)_{res} $ is orthogonal to the space of cusp forms.

For the continuous part of the decomposition, we show first that no nonzero compactly supported wave packet is a cusp form (cf. Proposition \ref{Prop no non nonzero compactly supported wave packet is cuspidal}).

Second, we show with the help of this result and Theorem \ref{Thm a Schwartz fct in the closure of the space of cusp forms is a cusp form} that $ L^2(\Gamma \bs G)_{ac} $ is actually orthogonal to the space of cusp forms (see Theorem \ref{Thm L^2(Gamma|G, phi)_ac is orthogonal to the space of cusp forms}).

A very important role in the proofs play the asymptotic expansions in the sense of N.~Wallach (see Section \ref{ssec:SomedefinitionsAndKnownResults}) and the Bunke-Olbrich extension map (cf. Definition \ref{Dfn Bunke-Olbrich extension map}).

\SecPar It remains to investigate $ L^2(\Gamma \bs G)_{ds} $:

Let $ \pi $ denote a discrete series representation. Recall that $ (V_{\pi, -\infty})_d^\Gamma $ consists of the $ \Gamma $-invariant distribution vectors $f$ for which the matrix coefficients $ c_{f, v} $ ($ v \in V_{\pi', K} $, where $ \pi' $ is the dual representation of $ \pi $) are square-integrable. We denote the space of $ \Gamma $-invariant distribution vectors $f$ for which the matrix coefficients $ c_{f, v} $ ($ v \in V_{\pi', K} $) belong to $ \CS(\Gamma \bs G) $ (resp. to $ \c{\CS}(\Gamma \bs G) $) by $ (V_{\pi, -\infty})_\CS^\Gamma $ (resp. $ (V_{\pi, -\infty})_{\c{\CS}}^\Gamma $). By Proposition \ref{Prop Every Schwartz vector is cuspidal}, $ (V_{\pi, -\infty})_\CS^\Gamma = (V_{\pi, -\infty})_{\c{\CS}}^\Gamma $.

We determine in the following a more precise description of the spaces $ (V_{\pi, -\infty})_\CS^\Gamma = (V_{\pi, -\infty})_{\c{\CS}}^\Gamma $, $ (V_{\pi, -\infty})_d^\Gamma $ and $ (V_{\pi, -\infty})^\Gamma $ from which follows that $ (V_{\pi, -\infty})_\CS^\Gamma $ is dense in $ (V_{\pi, -\infty})_d^\Gamma $.

Let us continue with introducing some more notations.
For $ \lambda \in \a_\CC^* $, set $ \sigma_\lambda(man) = \sigma(m) a^{\rho-\lambda} $. This defines a representation of $P$ on $ V_{\sigma_\lambda} := V_\sigma $. We denote by $ V(\sigma_\lambda) := G \times_P V_{\sigma_\lambda} $ the associated homogeneous bundle over $ \dX = G/P $.

We denote by $ \pi^{\sigma, \lambda} $ the $G$-representation on the space of sections of $ V(\sigma_\lambda) $ given by the left regular representation of $G$ and call it a \textit{principal series representation} of $G$.

Let us denote the dual representation of $ \sigma $ by $ \tilde{\sigma} $.
The space of distribution sections $ \C^{-\infty}(\dX, V(\sigma_\lambda)) $ is defined by $ \Cinf(\dX, V(\tilde{\sigma}_{-\lambda}))' $. Here, we take the strong dual.

We denote the space of $ \Gamma $-invariant distributions in $ \C^{-\infty}(\dX, V(\sigma_\lambda)) $ which are smooth on $ \Omega $ by $ \C^{-\infty}_\Omega(\dX, V(\sigma_\lambda))^\Gamma $ and we denote the space of $ \Gamma $-invariant distributions in $ \C^{-\infty}(\dX, V(\sigma_{-\lambda})) $ having support on the limit set $ \Lambda $ by $ \C^{-\infty}(\Lambda, V(\sigma_{-\lambda}))^\Gamma $.

Let $ (\pi, V_\pi) $ be a discrete series representation. Then, $ V_{\pi, -\infty} $ can be realised in $ \C^{-\infty}(\dX, V(\sigma_{-\lambda})) $. Moreover, $ V_{\pi, -\infty} $ (realised) is the image of a $G$-intertwining operator $ A \colon \C^{-\infty}(\dX, V(\sigma_\lambda)) \to \C^{-\infty}(\dX, V(\sigma_{-\lambda})) $.
Let
\[
    D_1 = A(\C^{-\infty}_\Omega(\dX, V(\sigma_\lambda))^\Gamma)
    \quad \text{and} \quad
    D_2 = \C^{-\infty}(\Lambda, V(\sigma_{-\lambda}))^\Gamma \cap (V_{\pi, -\infty})^\Gamma \ .
\]
By Bunke-Olbrich \cite{BO00}, it is known that $ D_1 $ is an infinite-dimensional subspace of the multiplicity space $ (V_{\pi, -\infty})_d^\Gamma $.

We show that the $K$-finite matrix coefficients $ c_{f, v} $ with $ f \in D_1 \oplus D_2 $ and $ v \in \C^K(\dX, V(\tilde{\sigma}_{\lambda})) $ ($K$-finite) are cusp forms (cf. Proposition \ref{Prop Af, smooth on Omega, invariant distributions are Schwartz}).
Furthermore, we prove that we have the following topological direct sum:
\begin{equation}\label{eq Intro (V_(pi, -infty))_d^Gamma = overline(D_1) oplus D_2}
    (V_{\pi, -\infty})_d^\Gamma = \overline{D_1} \oplus D_2
\end{equation}
(Corollary \ref{Cor decomp. of space of cuspidal vectors if pi is a DS rep}). It follows from \eqref{eq Intro (V_(pi, -infty))_d^Gamma = overline(D_1) oplus D_2} that $ (V_{\pi, -\infty})_\CS^\Gamma = (V_{\pi, -\infty})_{\c{\CS}}^\Gamma  $ is dense in $ (V_{\pi, -\infty})_d^\Gamma $. So, $ L^2(\Gamma \bs G)_{ds} $ is contained in $ \overline{\c{\CS}(\Gamma \bs G)} $.

Finally, we determine a more explicit description of $ (V_{\pi, -\infty})^\Gamma $ and of $ (V_{\pi, -\infty})_\CS^\Gamma $:
\begin{repthm}{Thm Decompositions of (V_(pi, -infty) otimes V_phi)_CS^Gamma and of (V_(pi, -infty) otimes V_phi)^Gamma}
    We have the following topological direct sum decompositions:
    \[
        (V_{\pi, -\infty})^\Gamma
        = A(\C^{-\infty}(\dX, V(\sigma_\lambda))^\Gamma) \oplus D_2
    \]
    and
    \[
        (V_{\pi, -\infty})_\CS^\Gamma = \C^{-\infty}_\Omega(\dX, V(\sigma_{-\lambda}))^\Gamma \cap \im(A)
        = D_1 \oplus D_2 \ .
    \]
\end{repthm}
To obtain this result, it was necessary to provide an explicit formula for some scalar products (cf. Theorem \ref{Thm (Af, g) = C (res(f), res(g))}).
Our investigations provide also a direct sum decomposition of the Schwartz space (cf. Proposition \ref{Prop CS(Gamma|G, phi)(gamma)}).

\needspace{3\baselineskip}
\SecPar \textit{Finally, I would like to thank Prof.~M.~Olbrich for having given me the opportunity to work in this beautiful field of research, for the many interesting discussions we had and for helping me to overcome the occurring technical difficulties. A good example is here Appendix \ref{sec: Results of independent interest about discrete series representations}, which was necessary in order to complete the work about the discrete series representations. Moreover, I would not have discovered the relation between the constant term of a Schwartz function along the ordinary set and the constant term of a Schwartz function along a cusp if my supervisor had not been convinced that there is one. Last but not least, my thanks go to my family and friends who have accompanied me along the way.}

\newpage

\section{The convex-cocompact case}

\subsection{Geometric preparations}\label{ssec:GeometricPreparations}

Let $ G $ be a real simple linear connected Lie group of real rank one. We denote its neutral element by $e$. Let $ \g $ be the Lie algebra of $ G $ and let $ \g_\CC $ be its complexification.

Let $ \l $ be a Lie algebra over $ \RR $. Let $ \U(\l) $\index[n]{UL@$ \U(\l) $} be the universal enveloping algebra of $ \l_\CC $ and let $ \U(\l)_l $ be the vector subspace of $ \U(\l) $ (filtered algebra) consisting of elements of order at most $l$.

Fix a Cartan involution $ \theta $ of $\g$, i.e. a nontrivial involutive automorphism on $ \g $. We denote the corresponding Cartan involution on $G$ also by $ \theta $.
Let $ \k $ be the $1$-eigenspace of $ \theta $ in $ \g $ and let $ \s $ be the $(-1)$-eigenspace of $ \theta $ in $ \g $.

We have: $ \g = \k \oplus \s $ (Cartan decomposition). Let $K$ be the analytic subgroup of $ G $ with Lie algebra $ \k $ (it is a maximal compact subgroup of $G$).

Fix a maximal (i.e. one-dimensional) abelian Lie subalgebra $ \a $ of $ \s $ (it exists as $ \s $ is finite-dimensional and it is unique up to the adjoint action of $K$ by Lemma 2.1.9 of \cite[p.47]{Wallach}). Let $ A $ be the analytic subgroup of $ G $ with Lie algebra $ \a $. Let $ M = Z_K(A) \subset K $ be the centraliser of $A$ in $K$. We denote the Lie algebra of $ M $ by $ \m $.

Let $ \a^{*}_{\CC} $ denote the complexified dual of $ \a $. For $ \mu \in \a^{*}_{\CC} $, we set $ a^{\mu} = e^{\mu(\log(a))} \in \CC $.

For $ \mu \in \a^* $, set $ \g_\mu = \{ X \in \g \mid [H, X] = \mu(H)X \quad \forall H \in \a \} $.

We denote by $ \Phi = \Phi(\g, \a) := \{ \mu \in \a^* \mid \mu \neq 0 \quad \text{and} \quad \g_\mu \neq \{0\} \} $ the system of restricted roots and by $ \Phi^+ $ a system of positive restricted roots. The \textit{multiplicity} of $ \mu \in \a^* $ is by definition $ m_\mu := \dim_\RR(\g_\mu) $.
Let $ \alpha \in \Phi^+ $ be the shortest root of $ \Phi $ which is in $ \Phi^+ $.

It follows from Corollary 2.17 of Chapter VII in \cite[p.291]{Helgason} that $ \Phi $ is contained in $ \{\pm \alpha, \pm 2\alpha\} $.

Set $ \n = \g_{\alpha} \oplus \g_{2\alpha} $ and $ \bar{\n} = \g_{-\alpha} \oplus \g_{-2\alpha} $.
Put $ A_{+} = \{ a \in A \mid a^{\alpha} > 1 \} $ and $ \bar{A}_{+} = \{ a \in A \mid a^{\alpha} \geq 1 \} $.
Let $ \a_{+} = \Lie(A_{+}) $ and $ \bar{\a}_+ = \Lie(\bar{A}_{+}) $.
\\ For $ s > 0 $, put $ A_{\leq s} = \{ a \in A \mid a^{\alpha} \leq s\} $ and $ A_{> s} = \{ a \in A \mid a^{\alpha} > s\} $.

We identify $ \a $ with $ \RR $ as an abelian Lie algebra via the algebra homomorphism $ H \in \a \mapsto \alpha(H) \in \RR $ and we identify $ A $ with $ (0, \infty) $ as a group via the group homomorphism $ a \in A \mapsto a^{\alpha} = e^{\alpha(\log(a))} \in (0, \infty) $. Thus, $ \bar{A}_{+} \equiv [1, \infty) $ as a set.

Hence, $ \exp \colon \a \to A $ can be identified with the exponential map on $ \RR $.

Define $ \rho \in \a^{*} $ by
$
    \rho(H) = \tfrac12 \tr(\restricted{\ad(H)}{\n}) \qquad (H \in \a) \ .
$

By construction, $ \rho $ is also equal to
\[
    \tfrac12 \sum_{\mu \in \Phi^+} m_\mu \mu = \tfrac12 (m_{\alpha} + 2m_{2\alpha}) \alpha \ .
\]
We have
\begin{center}
    \begin{tabular}{c|c|c|c|c}
                                $X$         & $\RR H^n$             & $\CC H^n$     & $\HH H^n$            & $\OO H^2$ \\
        \hline \rule{0mm}{1em}  $\rho$      & $\frac{n-1}2\alpha$ & $n \alpha$  & $(2n+1)\alpha$     & $11\alpha$
    \end{tabular} \ .
\end{center}

Let $N$ be the analytic subgroup of $ G $ with Lie algebra $ \n $. Then, $ P := MAN $ is called a \textit{parabolic subgroup} of $G$.\index{parabolic subgroup}
Let $ \bar{N} = \theta(N) $.

Each element of $ g \in G $ has a Cartan decomposition $ g = k_g a_g h_g $ with $ k_g, h_g \in K, a_g \in \bar{A}_+ $. Moreover, if $ g \neq e $, then $ a_g $ and $ k_g M $ are uniquely determined.

The quotient space $ X := G/K $ is a Riemannian symmetric space of strictly negative sectional curvature. Its geodesic boundary is $ \dX := G/P = K/M $. The Oshima compactification (cf. \cite{Oshima}) endows the set $ \bar{X} := X \cup \dX $ with a structure of a compact smooth $G$-manifold with boundary.

We know from the classification of symmetric spaces with strictly negative sectional curvature that $G$ is a linear group finitely covering the orientation-preserving isometry group of $X$ and that $ X $ is a real hyperbolic space $ \RR H^n $ ($ n \geq 1 $), a complex hyperbolic space $ \CC H^n $ ($ n \geq 2 $), a quaternionic hyperbolic space $ \HH H^n $ ($ n \geq 2 $) or the Cayley hyperbolic plane $ \OO H^2 $. 

Since the restricted root space decomposition is an orthogonal decomposition with respect to the Killing form ($\g$ is the orthogonal direct sum of simultaneous eigenspaces), $ \m \oplus \a = \g_0 $ is orthogonal to $ \n $ with respect to the Killing form. Moreover, $ \m = \k \cap \g_0 $ is orthogonal to $ \a = \s \cap \g_0 $ relative to $B$. See, e.g., Proposition 6.40 of Chapter VI in \cite[p.370]{Knapp2}. Since $ \a = \s \cap \g_0 $, $ \a $ is also orthogonal to $ \k $.

Note that $M$ is equal to $ K \cap P $ and that $ \a $ is the unique maximal abelian Lie subalgebra of $ \m \oplus \a \oplus \n $ which is orthogonal to the Lie algebra of $K$.

With respect to the Iwasawa decomposition $ G = KAN $ (resp. $ G = NAK $), we write
\[
    g = \kappa(g) a(g) n(g) \qquad \big(\text{resp. } g = \nu(g) h(g) k(g) \big)
\]
with $ \kappa(g) \in K $, $ a(g) \in A $ and $ n(g) \in N $ (resp. $ \nu(g) \in N, h(g) \in A $ and $ k(g) \in K $).

Consider the function $ G \times K \ni (g, k) \mapsto a(g^{-1}k) $. Note that this function descends to $ X \times \dX $.

\begin{lem}\label{Lem a = e^(d(eK, aK))}
    We can choose the normalisation of the distance function $ d_X $ (often simply denoted by $d$) on $ X = G/K $ such that $ a = a^{\alpha} = e^{d(eK, aK)} $ for all $ a \in \bar{A}_{+} \equiv [1, \infty) $.
    If $ a \in A_{<1} $, then $ a = e^{-d(eK, aK)} $.
\end{lem}

\begin{proof}
    Let $d$ be a distance function on $X$ (with arbitrary normalisation).

    Since $ X $ is a Riemannian symmetric with strictly negative sectional curvature (in particular it is a Hadamard space), the length of a geodesic segment connecting to points $ x, y \in X $ is equal to $ d(x, y) $.

    Let $ F(H) $ be equal to $ \pm d(eK, \exp(H)K) $, where the sign $ \pm $ is positive if $ H \in \bar{\a}_{+} $ and negative otherwise. Using that $eK$, $a_1K$ and $a_2K$ are on a geodesic for all $ a_1, a_2 \in A $ and the above fact, one easily sees that $F \colon \a \equiv \RR \to \RR $ is a continuous homomorphism.
    Since all nonzero continuous homomorphisms from $ \RR $ to $ \RR $ are of the form $ x \in \RR \mapsto \lambda x \in \RR $ for some $ \lambda \neq 0 $, $ F(H) = \lambda \alpha(H) $ for some $ \lambda \neq 0 $.
    We normalise $d$ so that $ \lambda = 1 $. Hence, $ F(H) = \alpha(H) $. So, $ e^{\pm d(eK, aK)} = e^{\alpha(\log a)} = a^\alpha = a $ for all $ a \in A $.
    The lemma follows.
\end{proof}

Consequently,
\begin{equation}\label{eq a_g = e^(d(eK, gK))}
    a_g = e^{d(eK, gK)} \qquad (g \in G)
\end{equation}
due to the $G$-invariance of the metric.

Let $B$ be the Killing form of $ \g $. Then, the Riemannian metric induces a $ \Ad(K) $-invariant inner product $ \tilde{B} $ on $ \s $.
Since $ \g $ is simple and of real rank one, the adjoint action of $K$ acts irreducibly on $ \s $ ($K$ acts transitively on $ \{ X \in \s \mid B(X, X) = 1 \} $).
It follows now from Schur's lemma that there is a nonzero $ c_{d_X} \in \RR $\index[n]{cdX@$ c_{d_X} $} such that $ \tilde{B}(X, Y) = c_{d_X} B(X, Y) $ for all $ X, Y \in \s $. Thus, $ c_{d_X} B $ is an $ \Ad(G) $-invariant nondegenerate symmetric bilinear form extending $ \tilde{B} $. We denote this form also by $ \tilde{B} $.
As both bilinear forms are positive definite on $ \s $, $ c_{d_X} > 0 $.
For $ X, Y \in \g $, set $ \langle X, Y \rangle = -\tilde{B}(X, \theta Y) $. Then, $ \langle \cdot, \cdot \rangle $ is an $ \Ad(K) $-invariant inner product on $ \g $ and $ |Z| := \sqrt{\langle Z, Z \rangle} $, $ Z \in \g $, defines a norm on $ \g $.

The map $ H \in \a \mapsto \alpha(H) \in \RR $ identifies $ \a $ with $ \RR $ isometrically. Hence, $ \a_\CC $ can be isometrically identified with $ \CC $ and $ \mu \in \a^* \mapsto \langle \alpha, \mu \rangle \in \RR $ identifies $ \a^* $ with $ \RR $ isometrically.
Moreover, $ \alpha $ is identified with 1. So, we can identify isometrically $ \a_\CC^* $ with $ \CC $.

Let us determine now the constant $ c_{d_X} $ explicitly. Let $ H \in \a $. Then,
\[
    |H|^2 = \tilde{B}(H, H) = c_{d_X} B(H, H) = c_{d_X} \tr\bigl(\ad(H) \ad(H)\bigr) \ .
\]
As $ [H, Y] = 0 $ for all $ Y \in \m \oplus \a $ and $ [H, X] = \mu(H) X $ for all $ X \in \g_\mu $,
\[
    \tr\bigl(\ad(H) \ad(H)\bigr) = 2 m_\alpha \alpha(H)^2 + 2 m_{2\alpha} \big(2 \alpha(H)\big)^2 = 2 \big( m_\alpha + 4 m_{2\alpha} \big) \alpha(H)^2 \ .
\]
Let $ H_1 \in \a $ be such that $ \alpha(H_1) = 1 $. Then, $ |H_1| = 1 $. Thus,
\[
    1 = |H_1|^2 = c_{d_X} \tr\bigl(\ad(H_1) \ad(H_1)\bigr) = 2 c_{d_X} \big( m_\alpha + 4 m_{2\alpha} \big) \ .
\]
Hence,
\begin{equation}\label{eq c_(d_X)}
    c_{d_X} = \frac1{2 (m_\alpha + 4 m_{2\alpha})} \ .
\end{equation}
If $ \mu \in \a^* $, then we define $ H_\mu \in \a $ by $ \tilde{B}(H, H_\mu) = \mu(H) $ ($ H \in \a $). This is well-defined as $\tilde{B}$ is nondegenerate.
We define a symmetric bilinear form $ \langle \cdot, \cdot \rangle $ on $ \a^* $ by $ \langle \mu, \tau \rangle = \tilde{B}(H_\mu, H_\tau) $ for all $ \mu, \tau \in \a^* $.
For $ \mu \in \a^* $, set
\[
    |\mu| := \sqrt{\langle \mu, \mu \rangle} = \sqrt{\tilde{B}(H_\mu, H_\mu)}
\]
(well-defined since the Killing form is positive definite on $ \a \subset \s $). 

The two isometries $ H \in \a \mapsto \alpha(H) \in \RR $ and $ \mu \in \a^* \mapsto \langle \alpha, \mu \rangle \in \RR $ fix the Haar measures on $A$ and $ \a^* $. We normalise the Haar measures on compact groups such that the groups have total mass 1.

The push-forward of the Lebesgue measure on $\n$ (resp. $ \bar{\n} $), normalised by the condition
\[
    \int_N a(\theta(n))^{-2\rho} \, dn = 1
    \qquad \big(\text{resp. } \int_{\bar{N}} a(\bar{n})^{-2\rho} \, d\bar{n} = 1 \big) \ ,
\]
defines a Haar measure on $N$ (resp. $ \bar{N} $).
On $G$ we consider the Haar measure $ dg = dk \, dx $, where $ dk $ is the Haar measure on $K$ and $ dx $ is the Riemannian measure on $ X = G/K $. On discrete groups, we take the counting measure.

Let $ \HS_{gK, kM} := k N k^{-1}gK $ be the horosphere passing through $ gK \in X $ and $ kM \in \dX $.\index[n]{HgKkM@$ \HS_{gK, kM} $}

\begin{lem}\label{Lem d(eK, HS_(aK, eM)) = d(eK, aK)}
    Let $ a \in A $. Then, the distance from $ eK $ to the horosphere
    \[
        \HS_{aK, eM} = \{ n a K \mid n \in N \}
    \]
    passing through $ aK \in X $ and $ eM \in \dX $ is equal to $ d(eK, aK) = |\log a| $.
\end{lem}

\begin{proof}
    Let $ a \in A $. As $ xK \mapsto a^{-1}xK $ is an isometry of $X$, it suffices to prove that the distance from $ a^{-1}K $ to the horosphere $ \HS_{eK, eM} = \{ n' K \mid n' \in N \} $ passing through $ eK \in X $ and $ eM \in \dX $ is equal to $ d(eK, aK) = |\log a| $.

    Let $ p = \frac12(\Id - \theta) $. Then, $ p \colon \g \to \g $ is the projection on $ \s \equiv T_{eK} X $.
    Since $ B(H, X) = 0 $ for all $ H \in \a $ and $ X \in \n \oplus \bar{\n} $, $ \a \subset \s $ is orthogonal to $ p(\n) \subset \s $ with respect to $ B $. Thus, $ \{a'K \mid a' \in A\} $ is orthogonal to $ \{n'K \mid n' \in N\} $ in $X$.
    As $ xK \mapsto n xK $ is an isometry of $X$ for every $n \in N$, $ \{ n a' K \mid a' \in A \} $ is a geodesic which is orthogonal to $ \{n'K \mid n' \in N\} $ for every $ n \in N $.

    Assume now that a geodesic is orthogonal to $ \{n'K \mid n' \in N\} $. Since the latter space has codimension 1 in $X$, it follows from the uniqueness property of a geodesic that it is equal to the geodesic $ \{ n a' K \mid a' \in A \} $ for some $ n \in N $.

    Since $ a^{-1}K $ must belong to the geodesic, it follows that $ \{ a'K \mid a' \in A\} $ is the unique geodesic which passes through $ a^{-1}K $ and which is orthogonal to $ \{n'K \mid n' \in N\} $.
    Hence,
    \[
        d(a^{-1}K, \HS_{eK, eM}) = d(a^{-1}K, eK) = d(eK, aK)= |\log a|
    \]
    as the geodesic $ \{ a'K \mid a' \in A\} $ intersects $ \HS_{eK, eM} = \{ n' K \mid n' \in N \} $ at $ eK $.
\end{proof}

\begin{lem}
    We have:
    \begin{equation}\label{a(g^(-1)k) <= a_g}
        a(g^{-1}k) \leq a_g \quad \text{and} \quad a(g^{-1}k)^{-1} \leq a_g \qquad ( g \in G, k \in K ) \ .
    \end{equation}
\end{lem}

\begin{proof}
    By Lemma \ref{Lem d(eK, HS_(aK, eM)) = d(eK, aK)},
    \[
        d(eK, \HS_{gK, eM}) = d(eK, \HS_{h(g)K, eM}) = d(eK, h(g)K) = d(eK, a(g^{-1})K) \ .
    \]
    So,
    \begin{multline*}
        d(eK, \HS_{gK, kM}) = \inf_{n \in N} d(eK, k n k^{-1}g K) = \inf_{n \in N} d(eK, n k^{-1}g K) \\
        = d(eK, \HS_{k^{-1}gK, eM}) = d(eK, a(g^{-1}k)K) = |\log a(g^{-1}k)| \ .
    \end{multline*}
    It follows that
    \[
        a(g^{-1}k) = e^{\pm d(eK, \HS_{gK, kM})} \qquad (g \in G, k \in K) \ ,
    \]
    where the sign $ \pm $ is positive (resp. negative) if $ eK $ lies inside (resp. outside) the corresponding horoball.
    Thus,
    \[
        |\log a(g^{-1}k)| \leq d(eK, k e k^{-1}gK) = d(eK, gK) = \log a_g \ .
    \]
    Hence,
    \[
        a(g^{-1}k) \leq a_g \quad \text{and} \quad a(g^{-1}k)^{-1} \leq a_g
    \]
    for all $ g \in G $ and $ k \in K $. The lemma follows.
\end{proof}

\begin{lem}\label{Lem a_(na) estimates}
    Let $ n \in N $ and $ a \in A $. Then,
    \[
        a_{na} \geq a_a \geq a \, , \ a_{na}^2 \geq a_n \quad \text{and} \quad \log a_{na} \geq |\log a| \ .
    \]
\end{lem}

\begin{proof}
    Let $ n \in N $ and $ a \in A $. Then, $ a_n \leq a_{na} a_{a^{-1}} = a_{na} a_a $ by the triangle inequality for $d_X$. It follows from \eqref{a(g^(-1)k) <= a_g} that $ a_{na} \geq a_a \geq a $. Thus, $ a_n \leq a_{na}^2 $ and $ \log a_{na} \geq \log a_{a} = |\log a| $.
\end{proof}

\begin{lem}\label{Lem 1 + log a_h <= (1 + log a_gh)(1 + log a_g)}
    For all $ g, h \in G $,
    \[
        1 + \log a_{gh} \leq (1 + \log a_g)(1 + \log a_h) \, , \quad 1 + \log a_g \leq (1 + \log a_{gh})(1 + \log a_h)
    \]
    and
    \[
        1 + \log a_h \leq (1 + \log a_g)(1 + \log a_{gh}) \ .
    \]
\end{lem}

\begin{proof}
    Let $ g, h \in G $. Since it follows from the triangle inequality for the Riemannian metric $d$ on $X$ that $ a_{gh} \leq a_g a_h $,
    \[
        \log a_{gh} \leq \log a_g + \log a_h \ .
    \]
    Thus,
    \begin{multline*}
        (1 + \log a_g)(1 + \log a_h) = 1 + \log a_g + \log a_h + \log a_g \log a_h \\
        \geq 1 + \log a_g + \log a_h \geq 1 + \log a_{gh} \ .
    \end{multline*}
    We get the second formula by replacing $ g $ by $ gh $ and $ h $ by $ h^{-1} $ and the third formula by replacing $ g $ by $ g^{-1} $ and $ h $ by $ gh $.
\end{proof}

\begin{lem}\label{Lem about a(g)}
    If $ S_1 $ is a relatively compact subset of $ G $ and if $ S_2 $ is a subset of $P$ such that $ \{ a(h) \mid h \in S_2 \} $ is relatively compact in $ A $ (e.g. $ S_2 $ is a relatively compact subset of $ P $), then there exist constants $ c_1, c_2 > 0 $ such that
    \[
        c_1 a(x) \leq a(g x h) \leq c_2 a(x)
    \]
    for all $ g \in S_1 $, $ h \in S_2 $ and $ x \in G $.
\end{lem}

\begin{proof}
    Let $ g \in S_1 $, $ h \in S_2 $ and $ x \in G $. Then,
    \[
        a(gxh) = a(g \kappa(x)) a(x) a(h) \ .
    \]
    The lemma follows as $ \{ g \kappa(x) \mid g \in S_1, x \in G \} $ is contained in the compact set $ \bar{S}_1K $ and as $ \overline{\{ a(h) \mid h \in S_2 \}} $ is compact.
\end{proof}

\needspace{4\baselineskip}
We write an element $ g \in \bar{N} M A N $ as $ \bar{n}_B(g) m_B(g) a_B(g) n_B(g) $ with $ \bar{n}_B(g) \in \bar{N} $, $ m_B(g) \in M $, $ a_B(g) \in A $ and $ n_B(g) \in N $ (Bruhat decomposition).

\begin{prop}\label{Theorem 3.8 of Helgason, p.414}
    Write $ \bar{n} \in \bar{N} $ as $ \bar{n} = \exp(X + Y) \ (X \in \g_{-\alpha}, Y \in \g_{-2\alpha}) $. Then,
    \begin{enumerate}
    \item $ 2 \cosh\bigl(2\log a_{\bar{n}}\bigr) = 1 + |X|^2 + (1 + \tfrac12 |X|^2)^2 + 2 |Y|^2 $ ,
    \item $ a(\bar{n}) = \sqrt{(1 + \tfrac12|X|^2)^2 + 2 |Y|^2} $ ,
    \item $ a_B(w\bar{n}) = \sqrt{\tfrac14|X|^4 + 2 |Y|^2} \qquad (\bar{n} \neq e) $ .
    \end{enumerate}
\end{prop}

\begin{proof}
    See Theorem 3.8 in Chapter IX of \cite[p.414]{Helgason} for a proof. In order to determine the appearing constant, we use that $ c_{d_X} = \frac1{2 (m_\alpha + 4 m_{2\alpha})} $ (see \eqref{eq c_(d_X)}).
\end{proof}

\begin{dfn}\label{Dfn asymp}
    Let $X$ be a set and let $f$, $g$ be two nonnegative real-valued functions on $X$. Then, we write $ f \prec g $ if there exists a positive constant $ c > 0 $ such that $ f(x) \leq c g(x) $ for all $ x \in X $ and we say that $f$ is \textit{essentially bounded} by $g$.
    We use the notation $ f \succ g $ if $ g \prec f $ and $ f \asymp g $ if $ f \prec g $ and $ g \prec f $.\index[n]{<>@$ \prec $, $ \succ $, $ \asymp $}
\end{dfn}

\begin{cor}\label{Cor a_n asymp a(theta n)}
    We have
    \[
        a_n \asymp a(\theta(n)) = a(n w)  \qquad (n \in N) \ .
    \]
    Let $ \eps > 0 $. Then,
    \[
        a_n \asymp a(\theta(n)) \asymp a_B(n w)
    \]
    for all $ n \in N $ such that $ |\log n| > \eps $.
\end{cor}

\begin{proof}
    Since
    \[
        a_B(n w) = a_B\big(w \underbrace{(w^{-1} n w)}_{ \in \bar{N} }\big) \qquad (n \in N) \ ,
    \]
    \[
        |\theta(X)| = |X| = |\Ad(w^{-1}) X| \qquad (X \in \n) \ ,
    \]
    \[
        \Ad(w^{-1})\g_\alpha \subset \g_{-\alpha} , \ \Ad(w^{-1})\g_{2\alpha} \subset \g_{-2\alpha}
    \]
    and since
    \[
        e^{\arccosh(x)} = x + \sqrt{x^2 - 1} = x \Big(1 + \sqrt{1 - \frac1{x^2}}\Big) \asymp x \qquad (x \in [1, +\infty)) \ ,
    \]
    the corollary follows immediately from Proposition \ref{Theorem 3.8 of Helgason, p.414}.
\end{proof}

\newpage

\subsection{Convex-cocompact groups}\label{ssec:Dfn convex-cocompact}

In this section, we say when a torsion-free discrete subgroup $ \Gamma $ of $G$ is convex-cocompact.

\medskip

Let us recall first of all the definition of ``torsion-free'':

\begin{dfn}
    Let $ H $ be a group. Then, we call $ x \in H $ a \textit{torsion element} if $ x $ has finite order.
    We call $ H $ \textit{torsion-free} if it has no nontrivial torsion elements.\index{torsion-free}
\end{dfn}

\begin{dfn}[{\cite[I.10, III.4]{Bourbaki71}}] \nl
    Let $ H $ be a locally compact group and let $ V $ be a locally compact topological space.
    If $ H $ acts continuously on $ V $, then we say that $ H $ acts \textit{properly} if
    \begin{equation}\label{Dfn Proper Action eq1}
        \{ h \in H \mid h.C \cap C \neq \emptyset \}
    \end{equation}
    is compact for any compact set $ C $ of $ V $.
    \\ If $ H $ is discrete and acts properly, then \eqref{Dfn Proper Action eq1} is finite and the action is called \textit{properly discontinuous}.
\end{dfn}

\begin{rem}
    If $H$ acts properly on $V$, then the space of orbits $ H \bs V $, endowed with the quotient topology, is Hausdorff.
\end{rem}

Let $ \Gamma $ be a torsion-free discrete subgroup of $G$. Then, the \textit{limit set} $ \Lambda_\Gamma $\index[n]{OmegaGamma LambdaGamma@$ \Omega_\Gamma $, $ \Lambda_\Gamma $} is by definition the set of accumulation points of the $ \Gamma $-orbit $ \bigcup_{\gamma \in \Gamma} \gamma K $ in $ \bar{X} := X \cup \dX $ (here $ X := G/K $ and $ \dX := G/P = K/M $). By Proposition 1.4 of \cite[p.88]{EO73}, this set is a closed and $ \Gamma $-invariant subset of $ \dX $. Let $ \Omega_\Gamma = \dX \smallsetminus \Lambda_\Gamma $ be the \textit{set of ordinary points} of $ \Gamma $. This set is open and clearly also a $ \Gamma $-invariant subset of $ \dX $. If no confusion is possible, then we may simply write $ \Omega $ (resp. $ \Lambda $) for $ \Omega_\Gamma $ (resp. $ \Lambda_\Gamma $).
By Proposition 8.5 of \cite[p.88]{EO73}, $ \Gamma $ acts properly discontinuously on $ X \cup \Omega_\Gamma $. As moreover $ \Gamma $ is torsion-free, $ \Gamma $ acts freely on $ X \cup \Omega_\Gamma $.

\begin{dfn}\label{Dfn convex-cocompact} \nlenum
    \begin{enumerate}
    \item We say that a left (resp. right) proper action of a group $ H $ on a topological space $ T $ is \textit{cocompact} if the quotient space $ H \bs T $ is a compact space.
    \item We say that a torsion-free, discrete subgroup $ \Gamma $ of $ G $ is \textit{convex-cocompact} if $ \Gamma $ acts cocompactly on $ X \cup \Omega $.\index{convex-cocompact}
    \end{enumerate}
\end{dfn}

\newpage

\subsection{The Schwartz space and the space of cusp forms on \texorpdfstring{$ \Gamma \bs G $}{Gamma\textbackslash G}}

\subsubsection{The Schwartz space on \texorpdfstring{$ \Gamma \bs G $}{Gamma\textbackslash G}}

First, we define the Schwartz space on $ \Gamma \bs G $. Second, we determine some basic properties of this space.

\medskip

Let $ \Gamma \subset G $ be a convex-cocompact, non-cocompact ($ \Omega \neq \emptyset $), discrete subgroup.

\begin{dfn}
    Let $ (\phi, V_\phi) $ be a finite-dimensional unitary representation of $ \Gamma $. 
    Define
    \[
        \Cinf(\Gamma \bs G, \phi)
        = \{ f \in \Cinf(G, V_\phi) \mid f(\gamma x) = \phi(\gamma)f(x) \quad \forall \gamma \in \Gamma, x \in G \} \ ,
    \]
    \[
        \Ccinf(\Gamma \bs G, \phi)
        = \{ f \in \Cinf(\Gamma \bs G, \phi) \mid f \text{ has compact support modulo $ \Gamma $} \}
    \]
    and\index[n]{LGammaGp@$ L^p(\Gamma \bs G, \phi) $}
    \begin{multline*}
        L^p(\Gamma \bs G, \phi)
        = \{ f \colon G \to V_\phi \text{ measurable} \mid \\
        \faktor{f(\gamma x) = \phi(\gamma)f(x) \quad \forall \gamma \in \Gamma, x \in G, \int_{\Gamma \bs G}{ |f(x)|^p dx < \infty } \}}{\sim} \ ,
    \end{multline*}
    with the usual identifications, for $ p \geq 1 $. 
\end{dfn}

Let $ \U_\Gamma $ be the family of open subsets of $X$ which are relatively compact in $ X \cup \Omega $. Let $ U \in \U_\Gamma $.\index[n]{UGamma@$ \U_\Gamma $}

\medskip
\needspace{5\baselineskip}
\underline{Picture for $ G = \SL(2, \RR) $:}

\begin{center}
    \begin{tikzpicture}[scale=0.25]
        \draw[thick]  (0,0) circle (10);

        \foreach \a in {20,..., 40, 70, 71,..., 85, 105, 106,..., 130, 170, 171, ..., 190, 230, 231, ..., 255, 275, 276, ..., 290, 320, 321, ..., 340}
                      \draw[orange] (\a:9.825) -- (\a:10.175);

        \draw[thick, dashed] plot[smooth, tension=.7] coordinates {(15:10) (4,1) (4,-1) (-15:10)};

        \node at (7,0) {$U$};

        \node at (10.5,5) {\textcolor{orange}{$ \Lambda $}};
    \end{tikzpicture}
\end{center}
For $Z \in \g$, set $ L_Z f(g) = \restricted{ \frac{d}{dt} }{ t = 0 } f(e^{-tZ}g) $ and $ R_Z f(g) = \restricted{ \frac{d}{dt} }{ t = 0 } f(ge^{tZ}) $. These definitions are extended to $ \U(\g) $ in the well-known way.

\index[n]{prXYf1@$ \LOI{p}{U}_{r, X, Y}(f) $}
Let $ l \in \NN_0 $. For $ r \geq 0 $, $ X \in \U(\g)_{l_1} $, $ Y \in \U(\g)_{l_2} $ such that $ l_1 + l_2 \leq l $ and a function $ f \in \C^l(G, V_\phi) $, set
\[
    \LOI{p}{U}_{r, X, Y}(f) = \sup_{gK \in U}{ (1 + \log a_g)^r a_g^\rho |L_X R_Y f(g)| } \ .
\]

\begin{dfn}\label{Dfn Schwartz space in the covex-cocompact case}
    Let $ (\phi, V_\phi) $ be a finite-dimensional unitary representation of $ \Gamma $. 
    Then, we define the Schwartz space $ \CS(\Gamma \bs G, \phi) $ on $ \Gamma \bs G $ by\index{Schwartz space}\index[n]{CGammaG@$ \CS(\Gamma \bs G, \phi) $}
    \[
        \{ f \in \Cinf(\Gamma \bs G, \phi) \mid \LOI{p}{U}_{r, X, Y}(f) < \infty \quad \forall r \geq 0, X, Y \in \U(\g), U \in \U_\Gamma \} \ .
    \]
    We equip $ \CS(\Gamma \bs G, \phi) $ with the topology induced by the seminorms.

    By definition, $ \CS(\Gamma \bs G) = \CS(\Gamma \bs G, 1) $, where $1$ is the one-dimensional trivial representation.
\end{dfn}

\begin{rem} \nlenum
    \begin{enumerate}
    \item This definition is a generalisation of Harish-Chandra's Schwartz space (see 7.1.2 of \cite[p.230]{Wallach}) as we get this space when we take $ \phi = 1 $ and $ \Gamma = \{e\} $ by Theorem 4.5.3 of \cite[p.126]{Wallach}.
    \item One can show that this definition is equivalent to Definition 8.1 in \cite[p.129]{BO00}.
    \end{enumerate}
\end{rem}

\begin{lem}\label{Lem p_(gamma U)_(r, X, Y)(f) < infty}
    Let $ U \in \U_\Gamma $, let $ r \geq 0 $ and let $ f \in \Cinf(\Gamma \bs G, \phi) $ be such that $ \LOI{p}{U}_{r, X, Y}(f) $ is finite for all $ X, Y \in \U(\g) $. For all $ \gamma \in \Gamma $, there is a constant $ C_{r, \gamma} > 0 $ such that
    \[
        \LOI{p}{\gamma U}_{r, X, Y}(f) \leq C_{r, \gamma} \LOI{p}{U}_{r, \Ad(\gamma^{-1})X, Y}(f) < \infty \qquad (X, Y \in \U(\g))  \ .
    \]
\end{lem}

\begin{proof}
    Let $ U \in \U_\Gamma $, let $ r \geq 0 $, let $ f $ be as above and let $ X, Y \in \U(\g) $. Then,
    \begin{align*}
        \LOI{p}{\gamma U}_{r, X, Y}(f) &= \sup_{gK \in U}{ (1 + \log a_{\gamma g})^r a_{\gamma g}^{\rho} |L_X R_Y f(\gamma g)| } \\
        &\leq (1 + \log a_\gamma)^r a_\gamma^{\rho} \sup_{gK \in U}{ (1 + \log a_g)^r a_g^{\rho} |L_{\Ad(\gamma^{-1})X} R_Y f(g)| } \\
        &= (1 + \log a_\gamma)^r a_\gamma^{\rho} \LOI{p}{U}_{r, \Ad(\gamma^{-1})X, Y}(f) < \infty \ ,
    \end{align*}
    where $ \Ad \colon G \to \GL(\g) $ is extended in the well-known way to $ \Ad \colon G \to \GL(\U(\g)) $. This completes the proof of the lemma.
\end{proof}

We denote by $ \clo(S) $ (resp. $ \interior(S) $) the closure (resp. interior) of a set $S$ in $X \cup \Omega$.\index[n]{clo, int@$ \clo $, $ \interior $}

\begin{lem}
    There is a set $ U_\Gamma \in \U_\Gamma $ such that
    \[
        \Gamma \bs X = \Gamma U_\Gamma \ .
    \]
\end{lem}

\begin{proof}
    As the compact set $ \Gamma \bs (X \cup \Omega) $ is covered by $ \bigcup_{U \in \U_\Gamma} { (U \cup \interior(\clo(U) \cap \dX)) } $ modulo $ \Gamma $, we can choose finitely many $ U_i $'s in $ \U_\Gamma $, say $ U_1, \dotsc, U_l $, such that the union of them covers $ X $ modulo $ \Gamma $.
    Take $ U_\Gamma := \bigcup_{i=1}^l U_i $.
\end{proof}

Fix $ U_\Gamma \in \U_\Gamma $ such that the previous lemma holds. 

\begin{lem}\label{Lem simpler version of the Schwartz space}
    $ \CS(\Gamma \bs G, \phi) $ is equal to
    \begin{equation}\label{eq LOI p U_F_(r, X, Y)(f) < infty}
        \{ f \in \Cinf(\Gamma \bs G, \phi) \mid \LOI{p}{U_\Gamma}_{r, X, Y}(f) < \infty \quad \forall r \in \NN_0, X, Y \in \U(\g) \}
    \end{equation}
    as topological spaces.
\end{lem}

\begin{proof}
    Let $ U \in \U_\Gamma $. Then, this set is covered by finitely many $ \Gamma $-translates of $U_\Gamma$: There exist $ \gamma_1, \dotsc, \gamma_s \in \Gamma $ such that
    $ U \subset \bigcup_{i=1}^s{ \gamma_i U_\Gamma } $. Let now $ \gamma \in \Gamma $, $ X, Y \in \U(\g) $, $ r \geq 0 $ and let $ f $ be in \eqref{eq LOI p U_F_(r, X, Y)(f) < infty}.

    Let $ r' $ be an integer which is greater or equal than $r$. Then, $ \LOI{p}{U_\Gamma}_{r, X, Y}(f) $ is less or equal than $ \LOI{p}{U_\Gamma}_{r', X, Y}(f) < \infty $.
    Thus, there is a constant $ C_{r', \gamma} > 0 $ such that
    \[
        \LOI{p}{\gamma U_\Gamma}_{r', X, Y}(f) \leq C_{r', \gamma} \LOI{p}{U_\Gamma}_{r', \Ad(\gamma^{-1})X, Y}(f) < \infty
    \]
    by Lemma \ref{Lem p_(gamma U)_(r, X, Y)(f) < infty}. So,
    \[
        \LOI{p}{U}_{r, X, Y}(f) \leq \sum_{i=1}^s C_{r', \gamma} \LOI{p}{U_\Gamma}_{r', \Ad(\gamma^{-1})X, Y}(f) < \infty \ .
    \]
    The lemma follows.
\end{proof}

A net $ f_\alpha $ in $ \Cinf(\Gamma \bs G, \phi) $ converges to $ f \in \Cinf(\Gamma \bs G, \phi) $ if for every compact subset $C$ of $G$, $ f_\alpha $ converges uniformly with all derivatives to $f$ on $ C $. This induces a topology on $ \Cinf(\Gamma \bs G, \phi) $.

\begin{lem}
    The space $ \Cinf(\Gamma \bs G, \phi) $ is complete.
\end{lem}

\begin{proof}
    This well-known result can be proven similarly as in 1.46 of \cite[pp.34-35]{Rudin91}.



\end{proof}

\begin{prop}\label{Prop Schwartz space is Fréchet is the convex-cocompact case}
    The Schwartz space $ \CS(\Gamma \bs G, \phi) $ is a Fréchet space.
\end{prop}

\begin{rem}
    The proof is a straightforward generalisation of that for the case $ \Gamma = \{e\} $ (see proof of Theorem 7.1.1 in \cite[pp.227-228]{Wallach}).
\end{rem}

\begin{proof}
    We can write $ U \in \U_\Gamma $ as a (possibly infinite) union of open sets in $X$ which are relatively compact in $X$.
    Indeed, $U$ is equal to $ \bigcup_{n \in \NN} (U \cap B(eK, n)) $.

    To show that $ \CS(\Gamma \bs G, \phi) $ is a Fréchet space, it suffices to prove that $ \CS(\Gamma \bs G, \phi) $ is sequentially complete as $ \U(\g) $ is countable dimensional and as it is enough to consider $ r \in \NN_0 $ and $ U = U_\Gamma $ by Lemma \ref{Lem simpler version of the Schwartz space}.

    Let $ (f_j)_j $ be a Cauchy sequence in $ \CS(\Gamma \bs G, \phi) $. Then, $ (f_j)_j $ is a Cauchy sequence in $ \Cinf(\Gamma \bs G, \phi) $.
    As $ \Cinf(\Gamma \bs G, \phi) $ is complete by the previous lemma, there exists $ f \in \Cinf(\Gamma \bs G, \phi) $ such that $ f_j $ converges uniformly with all derivatives to $f$ on compact subsets of $G$.

    Fix $ r \geq 0 $, $ X, Y \in \U(\g) $ and $ U \in \U_\Gamma $. We have:
    \[
        \LOI{p}{U}_{r, X, Y}(f) = \sup_{n \in \NN} \LOI{p}{U \cap B(eK, n)}_{r, X, Y}(f) \ .
    \]
    By convention, $ \LOI{p}{\emptyset}_{r, X, Y}(f) = 0 $.
    Let $ \eps > 0 $ be given. Let $ N \in \NN $ be sufficiently large so that $ \LOI{p}{U}_{r, X, Y}(f_i - f_j) < \frac{\eps}2 $ for all $ i, j \geq N $. Then,
    \[
        \LOI{p}{U}_{r, X, Y}(f_i) < \frac{\eps}2 + \LOI{p}{U}_{r, X, Y}(f_N)
    \]
    for all $ i \geq N $. Let $ N_n \in \NN $ be sufficiently large (maybe larger than $N$) so that
    \[
        \LOI{p}{U \cap B(eK, n)}_{r, X, Y}(f - f_i) < \frac{\eps}2
    \]
    for all $ i \geq N_n $. This is possible because for all $ n \in \NN $, there is $ c_n > 0 $ such that
    \[
        \LOI{p}{U \cap B(eK, n)}_{r, X, Y}(h) \leq c_n \sup_{gK \in U \cap B(eK, n)} |L_X R_Y h(g)| < \infty
    \]
    for all $ h \in \Cinf(\Gamma \bs G, \phi) $. Thus, $ \LOI{p}{U \cap B(eK, n)}_{r, X, Y}(f) < \eps + \LOI{p}{U \cap B(eK, n)}_{r, X, Y}(f_N) $. Hence,
    \[
        \LOI{p}{U}_{r, X, Y}(f) \leq \eps + \LOI{p}{U}_{r, X, Y}(f_N) < \infty \ .
    \]
    So, $ f \in \CS(\Gamma \bs G, \phi) $. Moreover,
    \[
        \LOI{p}{U \cap B(eK, n)}_{r, X, Y}(f - f_j) \leq \LOI{p}{U \cap B(eK, n)}_{r, X, Y}(f - f_i) + \LOI{p}{U \cap B(eK, n)}_{r, X, Y}(f_i - f_j) < \eps
    \]
    for all $ j \geq N $ and $ i \geq N_n $. It follows that $ \LOI{p}{U}_{r, X, Y}(f - f_j) \leq \eps $ for all $ j \geq N $.
    Consequently, $ \CS(\Gamma \bs G, \phi) $ is sequentially complete. The proposition follows.
\end{proof}

\begin{lem}\label{Lem int_G a_g^(-2rho) (1 + |log a_g|)^(-r) dg < oo}
    The integral $ \int_G a_g^{-2\rho} (1 + \log a_g)^{-r} dg $ is finite if and only if $ r > 1 $.
\end{lem}

\begin{proof}
    Compare with the lemma in 5.1.3 of \cite[p.139]{Wallach}.

    Since $ g \mapsto a_g $ is left and right $K$-invariant, there exists a constant $ c > 0 $ by Lemma 2.4.2 of \cite[p.61]{Wallach} such that
    \[
        \int_G{ a_g^{-2\rho} (1 + \log a_g)^{-r} \, dg } = c \int_{A_+}{ \gamma(a) a^{-2\rho} (1 + \log a)^{-r} \, da } \ ,
    \]
    where $ \gamma(a) := \prod_{\mu \in \Phi^{+}}{ \sinh(\mu(H)) } $.
    Since
    \[
        \gamma(a) \asymp \prod_{\mu \in \Phi^{+}}{ e^{\mu(\log a)} } = e^{ \sum_{\mu \in \Phi^{+}}{ \mu(\log a) } } = e^{2\rho(\log a)} = a^{2 \rho}
    \]
    on $ A $, this integral is finite if and only if
    \[
        \int_{A_+}{ (1 + \log a)^{-r} \, da } = \int_0^\infty{ (1 + t)^{-r} \, dt }
    \]
    is finite. The lemma follows as the last integral converges if and only if $ r > 1 $.
\end{proof}

\begin{prop}\label{Prop CS(Gamma|G) subset L^2(Gamma|G)}
    The Schwartz space $ \CS(\Gamma \bs G, \phi) $ is contained in $ L^2(\Gamma \bs G, \phi) $.
    \\ Moreover, the injection of $ \CS(\Gamma \bs G, \phi) $ into $ L^2(\Gamma \bs G, \phi) $ is continuous.
\end{prop}

\begin{proof}
    Let $ f \in \CS(\Gamma \bs G, \phi) $ and let $ r \geq 0 $.
    Then,
    \[
        |f(g)| \leq \LOI{p}{U_\Gamma}_{r, 1, 1}(f) (1 + \log a_g)^{-r} a_g^{-\rho}
    \]
    for all $ g \in S := U_\Gamma K \subset G $.
    We have:
    \[
        \|f\|^2_2 = \int_{\Gamma \bs G}{ |f(g)|^2 \, dg }
        \leq \int_S{ |f(g)|^2 \, dg }
        \leq \LOI{p}{U_\Gamma}_{r, 1, 1}(f)^2 \int_G{ (1 + \log a_g)^{-2r} a_g^{-2\rho} \, dg } \ .
    \]
    By Lemma \ref{Lem int_G a_g^(-2rho) (1 + |log a_g|)^(-r) dg < oo}, this is finite for $ r > \frac12 $. The proposition follows.
\end{proof}

\newpage

\subsubsection{\texorpdfstring{Density of $ \Ccinf(\Gamma \bs G, \phi) $ in $ \CS(\Gamma \bs G, \phi) $}{Density of the compactly supported smooth functions}}

\begin{sloppypar}
In the following, we do some preparation needed to prove that the space $ \Ccinf(\Gamma \bs G, \phi) $ of compactly supported smooth functions on $ \Gamma \bs G $ (viewed as a subspace of $ \CS(\Gamma \bs G, \phi) $) is dense (cf. Proposition \ref{Prop Ccinf(Gamma|G) is dense}).
\end{sloppypar}

\begin{lem}\label{Lemma 8.24 of Knapp}
    Let $ \mu $ be a restricted root $(\mu \in \{\pm \alpha, \pm 2 \alpha \})$ and let $ X \in \g_\mu $. Let $ Y = \tfrac12 (X + \theta X) \in \k $. If $ a \in A \smallsetminus \{e\} $, then
    \[
        X = \frac{ 2 a^\mu }{ a^{2\mu} - 1} (a^\mu Y - \Ad(a)^{-1} Y) \ .
    \]
\end{lem}

\begin{proof}
    This can easily be proven by direct computation. Compare with Lemma 8.24 of \cite[p.227]{Knapp}.
\end{proof}

Let $ A' = A \smallsetminus \{e\} \equiv (0, \infty) \smallsetminus \{1\} $. Let us denote by $ \Cinf_b(A_+) $ (resp. $ \Cinf_b(A') $) the smooth functions $f$ on $ A_+ $ (resp. $ A' $) having uniformly bounded derivatives on $ [1 + \eps, \infty) $ for every $ \eps > 0 $ (resp. $ (0, 1 - \eps] \cup [1 + \eps, \infty) $ for every $ \eps \in (0, 1) $).

If $ f \in \Cinf_b(A') $, then $ \restricted{f}{A_+} $ clearly belongs to $ \Cinf_b(A_+) $.

Let $ d_1, d_2 \in \NN_0 $ be such that $ d_1 \leq d_2 $. Let $ p $ (resp. $q$) be a polynomial of degree $ d_1 $ (resp. $ d_2 $). Let $ C $ be a closed subset of $ \RR $ which does not meet the zeros of $q$. Then, $ \frac{p}{q} $ has bounded derivatives of any order on $ C $.

Thus, the following functions belong to $ \Cinf_b(A') $:
\begin{equation}\label{eq functions with bounded derivatives}
    \frac1{x^2 - 1} \, , \qquad \frac{x}{x^2 - 1} \qquad \text{and} \qquad \frac{x^2}{x^2 - 1}  \ .
\end{equation}

\begin{lem}\label{Lem U(g) = U(a) oplus (Ad(a)^(-1)k U(g) + U(g)k)}
    Let $ a \in A' $. Then,
    \[
        \U(\g) = \U(\a) \oplus \bigl(\Ad(a)^{-1}(\k) \U(\g) + \U(\g) \k\bigr) \ .
    \]
\end{lem}

\begin{proof}
    We have
    \[
        \U(\k) = \CC \oplus \U(\k) \k = \CC \oplus \k \U(\k) \ .
    \]
    Let $ a \in A' $. By Lemma \ref{Lemma 8.24 of Knapp},
    \[
        \g_\CC = \n_\CC \oplus \a_\CC \oplus \k_\CC = \a_\CC \oplus (\Ad(a)^{-1}\k_\CC + \k_\CC) \ .
    \]
    Hence,
    \begin{multline*}
        \U(\g) = \bigl(\Ad(a)^{-1}\U(\k)\bigr)\U(\a)\U(\k)
        = \big(\CC \oplus \Ad(a)^{-1}(\k\U(\k))\big)\U(\a)(\CC \oplus \U(\k) \k) \\
        = \U(\a) \oplus \bigl(\Ad(a)^{-1}(\k)\U(\g) + \U(\g) \k\bigr)
    \end{multline*}
    by the Poincaré-Birkhoff-Witt theorem.
\end{proof}

Let $ \{H_i\} $ be the natural basis of $ \U(\a) $.
Let $ X \in \U(\a) $. Write
\[
    X = \sum_i c_i(X) H_i \ .
\]
Let $ \{X_{a, i}\} $ be a basis of $ \U(\g) $ induced from the above decomposition. Let $ X \in \U(\g) $. Write
\[
    X = \sum_i c_{a, i}(X) X_{a, i} \ .
\]
When we speak in the following of the coefficients of an element $ X \in \U(\g) $, then we mean with this the numbers $ c_{a, i}(X) \in \CC $.

Let $ p_a $ denote the projection of $ \U(\g) $ to $ \U(\a) $ provided by the previous lemma.

\begin{lem}\label{Lem c_i(p_a(X)) in Cinf_b(A')}
    Let $ X \in \U(\g) $. Then, $ a \in A' \mapsto c_i(p_a(X)) $ belongs to $ \Cinf_b(A') $ for all $i$.
\end{lem}

\begin{proof}
    Let $ l \in \NN_0 $. If $ X \in \CC $, there is nothing to prove. If $ X \in \g_\CC = \n_\CC \oplus \a_\CC \oplus \k_\CC $, then the lemma follows from Lemma \ref{Lemma 8.24 of Knapp}.

    Let us prove now by induction that the lemma holds also for every $ X \in \U(\g)_l $ ($ l \geq 2 $). Assume that the lemma holds for every $ X \in \U(\g)_{l-1} $. Since
    \[
        \U(\g)_l = \U(\g)_{l-1} \oplus \U(\g)_{l-1} \g \ ,
    \]
    it suffices to consider the case, where $ X = X_1 X_2 $ with $ X_1 \in \U(\g)_{l-1} $ and $ X_2 \in \g $.
    We have
    \begin{align*}
        \U(\g)_{l-1} \g &= \big(\U(\a)_{l-1} \oplus (\Ad(a)^{-1}(\k) \U(\g)_{l-2} + \U(\g)_{l-2} \k)\big)(\n \oplus \a \oplus \k) \\
        &= \underline{\U(\a)_{l-1} \n} + \Ad(a)^{-1}(\k) \U(\g)_{l-1} + \underline{\U(\g)_{l-2} \k \n} \\
        & \qquad + \U(\a)_l + \Ad(a)^{-1}(\k) \U(\g)_{l-1} + \underline{\U(\g)_{l-2} \k \a} \\
        & \qquad + \U(\g)_{l-1} \k \ .
    \end{align*}
    Here, the terms which are not under the desired form are underlined. In order to bring these terms in the form we want to have, we use that $ \n $ is contained in $ \Ad(a^{-1})\k + \k $ and the following relations:
    \begin{enumerate}
    \item $ \a \Ad(a^{-1})\k = \Ad(a^{-1})(\k) \a + [\a, \Ad(a^{-1})\k] $,
    \item $ \k \Ad(a^{-1})\k = \Ad(a^{-1})(\k) \k + [\k, \Ad(a^{-1})\k] $,
    \item $ \k \a = \a \k + [\k, \a] $.
    \end{enumerate}
    For the new terms appearing containing one of these new Lie brackets, we must check that we can apply the induction hypothesis.
    We do this by showing that all the coefficients of terms not belonging to $ \U(\g)\k $ are in $ \Cinf_b(A') $.

    Let $ q = \frac12(\Id + \theta) $. Then, $ q \colon \g \to \g $ is the orthogonal projection on $\k$.
    Since $ q(\n) = q(\bar{\n}) $ and $ q(\a) = \{0\} $, $ q(\n) \oplus \m = q(\n \oplus \m) = \k $.

    Let $ Y \in \k $ and $ H \in \a $. If $ Y \in \m $, then $ [H, \Ad(a^{-1})Y] = [H, Y] = 0 $, $ [\k, \Ad(a^{-1})Y] = [\k, Y] \subset \k $ and $ [Y, \a] = \{0\} $.

    For $ Z_k \in \g_{k\alpha} $ and $ H \in \a $, we have:
    \begin{enumerate}
    \item $ [H, \Ad(a^{-1})q(Z_k)] = \frac12 a^{-k} [H, Z_k] + \frac12 a^k \theta([H, Z_k]) $, where
        \[
            [H, Z_k] \in \n \subset \k + \Ad(a)^{\pm 1}\k \ .
        \]
    \item Let $ Y' \in \k $. Then,
        \[
            [Y', \Ad(a^{-1})q(Z_k)] = \frac12 a^{-k} [Y', Z_k] + \frac12 a^k \theta([Y', Z_k]) \ .
        \]
        If $ Y' \in \m $, then $ [Y', Z_k] \in \n \subset \k + \Ad(a)^{\pm1}\k $.

        Let $ Z'_k \in \g_{k\alpha} $. Then, $ [q(Z'_k), Z_k] \in \g_{-\alpha} \oplus \m \oplus \a \oplus \n $.
    \item $ [q(Z_k), H] \in \n \oplus \theta(\n) \subset \k + \Ad(a)^{-1}\k $.
    \end{enumerate}
    In the third case, all the coefficients are in $ \Cinf_b(A') $. In the first and the second case, the maximal growth behaviour is exactly the growth behaviour we are allowed to have in order to have that the coefficients of the elements of $ \a \n $, respectively $ \k \n $, are as desired. The lemma follows now by induction.
\end{proof}

We extend $ \alpha $ to an algebra homomorphism from $ \U(\a) $ to $ \CC $.

Define a homomorphism $ \tilde{\alpha} \colon \U(\a) \to \CC $ of Lie algebras by the following conditions:
\[
    \tilde{\alpha}(X) = 0 \text{ for } X \in \a_\CC \otimes \dotsm \otimes \a_\CC \supset \a_\CC \otimes \a_\CC
\]
and
\[
    \tilde{\alpha}(X) = \alpha(X) \qquad \text{ for } X \in \CC \oplus \a_\CC \ .
\]
\begin{prop}\label{Prop L_X R_Y a_g}
    Let $ X_1, X_2 \in \U(\g) $ and $ a \in A_+ $. Let $ f \colon G \smallsetminus K \to \CC $ be a smooth bi-$K$-invariant function. Then, there exist $ H_a \in \U(\a) $ (depending both on $ X_1 $ and $ X_2 $) such that $ a \in A_+ \mapsto \alpha(H_a) $ belongs to $ \Cinf_b(A_+) $ and
    \[
        L_{X_1} R_{X_2} \restricted{ f(g) }{g = a} = R_{H_a} \restricted{ f(g) }{g = a} \ .
    \]
    If $ f(g) = a_g $, then $ L_{X_1} R_{X_2} \restricted{ f(g) }{g = a} = \alpha(H_a) a $.
    \\ If $ f(g) = \log(a_g) $, then $ L_{X_1} R_{X_2} \restricted{ f(g) }{g = a} = \tilde{\alpha}(H_a) \log(a) $.
\end{prop}

\begin{rem}
    $ a_g $ is not differentiable at $ g \in K $: Let $ 0 \neq H \in \a $. Then,
    \[
        \restricted{ \frac{ d }{ dt } }{ t = 0^{\pm} } a_{\exp(t H)}
        = \restricted{ \frac{ d }{ dt } }{ t = 0^{\pm} } \max\{ t \alpha(H), -t \alpha(H) \} = \pm \alpha(H) \ .
    \]
    So, $ R_H a_g $ does not exist when $ g $ belongs to $ K $.
\end{rem}

\begin{proof}
    Let $ X_1, X_2 \in \U(\g) $.
    It follows from the Poincaré-Birkhoff-Witt theorem that
    \[
        \U(\g) = \U(\a) \U(\n) \U(\k) \ .
    \]
    By convention, this means that $ X \in \U(\g) $ can be written as $ \sum_i H_i Z_i Y_i $ with $ H_i \in \U(\a) $, $ Z_i \in \U(\n) $ and $ Y_i \in \U(\k) $.

    When we bring the left derivatives in $ \U(\n) $ to the right-hand side, then the corresponding right derivatives have coefficients in $ \Cinf_b(A_+) $.
    Thus, we may assume without loss of generality that $ X_1 = 1 $. The first part of the proposition follows now from Lemma \ref{Lem U(g) = U(a) oplus (Ad(a)^(-1)k U(g) + U(g)k)} and Lemma \ref{Lem c_i(p_a(X)) in Cinf_b(A')}.

    \needspace{2\baselineskip}
    The second part follows easily from the first part and the following formulas:
    \begin{enumerate}
    \item $ R_{H_1} \dotsm R_{H_k} \restricted{ a_g }{g = a} = R_{H_1} \dotsm R_{H_k} a = \prod_{j=1}^k \alpha(H_j) a $ for all $ H_j \in \a_\CC $ and all $ a \in A_+ $;
    \item $ R_H \restricted{ \log(a_g) }{g = a} = R_H \log(a) = \alpha(H) \log(a) $ for all $ H \in \a_\CC $ and all $ a \in A_+ $;
    \item $ R_{H_1} R_{H_2} \restricted{ \log(a_g) }{g = a} = R_{H_1} R_{H_2} \log(a) = 0 $ for all $ H_1, H_2 \in \a_\CC $ and all $ a \in A_+ $.
    \end{enumerate}
\end{proof}

\begin{cor}\label{Cor Estimate of L_X R_Y (1 + log a_g)}
    Let $ X_1, X_2 \in \U(\g) $ and let $V$ be a neighbourhood of $K$ in $G$. Then, there exists a constant $ c > 0 $ (depending on $X_1, X_2$) such that
    \[
        |L_{X_1} R_{X_2} \big(1 + \log(a_g)\big)| \leq c \big(1 + \log(a_g)\big)
    \]
    for all $ g \in G \smallsetminus V $.
\end{cor}

\begin{proof}
    We have
    \[
        L_{X_1} R_{X_2} a_g = L_{\Ad(k_g^{-1})X_1} R_{\Ad(h_g)X_2} \restricted{ a_x }{ x = a_g } \ .
    \]
    Let $ l \in \NN $ be sufficiently large so that $ X_1, X_2 \in \U(\g)_{l} $. Let $ \{Y_i\} $ be a basis of $ \U(\g)_{l} $. For $ k \in K $, write $ \Ad(k)X_j = \sum_i c_{i, j}(k) Y_i $. Then, $ \sup_{k \in K} |c_{i, j}(k)| $ is finite for all $i$, $j$. The corollary follows now from Proposition \ref{Prop L_X R_Y a_g}.
\end{proof}

\begin{lem}\label{Lem Right derivatives in k-direction suffice when we have all the left derivatives}
    Let $ (V, |\cdot|) $ be a finite-dimensional normed complex vector space and let $ f \in \Cinf(G, V) $ be such that $ \sup_{g \in G} |L_X R_Y f(g)| < \infty $ for all $ X \in \U(\g) $, $ Y \in \U(\k) \text{ (resp. for all } X \in \U(\k), Y \in \U(\g)) $. Then,
    \[
        \sup_{g \in G} |L_X R_Y f(g)| < \infty
    \]
    for all $ X, Y \in \U(\g) $.
\end{lem}

\begin{proof}
    Let $ f \in \Cinf(G, V) $ be such that $ \sup_{g \in G} |L_X R_Y f(g)| < \infty $ for all $ X \in \U(\g) $, $ Y \in \U(\k) $.
    Let $ X, Y \in \U(\g) $. As $ L_X f \in \Cinf(G) $, we may assume without loss of generality that $ X = 1 $.
    For $ g \in G $, we have
    \[
        R_Y f(g) = R_{\Ad(h_g)Y} (R_{h_g} f)(k_g a_g) \ .
    \]
    Thus, it remains to prove that $ \sup_{g \in G} |R_Y (R_{h_g} f)(k_g a_g)| < \infty $ for any $ Y \in \U(\g) $. Let $ Y \in \U(g) $.
    By the Poincaré-Birkhoff-Witt theorem, we have
    \[
        \U(\g) = \U(\a) \U(\bar{\n}) \U(\k) \ .
    \]
    By convention, this means that $ Y \in \U(\g) $ can be written as $ \sum_i H_i Z_i Y_i $ with $ H_i \in \U(\a) $, $ Z_i \in \U(\bar{\n}) $ and $ Y_i \in \U(\k) $.

    So, without loss of generality, we may assume that $ Y = H Z U $ with $ H \in \U(\a) $, $ Z \in \U(\bar{\n}) $ and $ U \in \U(\k) $. Then, $ Z \in \U(\bar{\n})_l $ for some $ l \geq 0 $ and $ |R_{Y} (R_{h_g} f)(k_g a_g)| $ ($ g \in G $) is less or equal than
    \begin{multline*}
        |L_{\Ad(k_g)H} R_Z R_U (R_{h_g} f)(k_g a_g)| \leq |L_{\Ad(k_g)H} L_{\Ad(a_g)Z} R_U (R_{h_g} f)(k_g a_g)| \\
        = |L_{\Ad(k_g)H} L_{\Ad(a_g)Z} R_{\Ad(h_g)^{-1}U} f(g)| \ .
    \end{multline*}
    Let $ \{ Z_j \} $ be a basis of $ \U(\bar{\n})_l $. As moreover $ \sup_{g \in G} |L_X f(g)| < \infty $ for all $ X \in \U(\g) $ and as $ \Ad(a)Z = \sum_j c_j(a) Z_j $ for some coefficients $ c_j(a) $ such that $ \sup_{a \in \bar{A}_+} |c_j(a)| $ is finite,
    $
        \sup_{g \in G} |L_X R_Y f(g)|
    $
    is finite.
    The second case of the lemma follows from the isomorphism $ \Cinf(G, V) \to \Cinf(G, V) , \ f \mapsto \check{f} $, where $ \check{f}(g) := f(g^{-1}) $.
\end{proof}

Recall that we endowed $ \bar{X} $ with a smooth manifold structure by using the Oshima compactification.
Let $ \Delta_X $ be the Laplace-Beltrami operator of $X$.

\begin{lem}\label{Lem Cut-off function in the convex-cocompact case}
    There exists a cut-off function $ \chi \in \Ccinf(X \cup \Omega, [0, 1]) $ such that
    \begin{enumerate}
    \item $ \sum_{\gamma \in \Gamma} L_\gamma \chi = 1 $,\label{Lem Cut-off function in the convex-cocompact case eq1}
    \item $ \sup_{gK \in X} a_g |d\chi(gK)| < \infty $,
        \label{Lem Cut-off function in the convex-cocompact case eq2}
    \item $ \sup_{gK \in X} a_g |\Delta_X \chi(gK)| < \infty $,\label{Lem Cut-off function in the convex-cocompact case eq3}
    \item $ \sup_{gK \in X} |L_X \chi(gK)| < \infty $ for all $ X \in \U(\g) $,\label{Lem Cut-off function in the convex-cocompact case eq4}
    \end{enumerate}
    We denote by $ \chi_\infty $ the restriction of $ \chi $ to $ \Omega $.
\end{lem}

\begin{rem}\label{Rem of Lem Cut-off function in the convex-cocompact case} \nlenum
    \begin{enumerate}
    \item It follows from Lemma \ref{Lem Right derivatives in k-direction suffice when we have all the left derivatives} that $ \sup_{gK \in X} |L_X R_Y \chi(g)| < \infty $ for all $ X, Y \in \U(\g) $.
    \item The set $ \{\gamma \in \Gamma \mid \gamma^{-1} U \cap \supp \chi \neq \emptyset \} $ is finite for every $ U \in \U_\Gamma $ since $ \chi \in \Ccinf(X \cup \Omega, [0, 1]) $ and since $ \Gamma $ acts properly discontinuous on $ X \cup \Omega $.
    \end{enumerate}
\end{rem}

\begin{proof}
    See Lemma 6.4 of \cite[p.113]{BO00}.
\end{proof}

Let us recall now how the usual topology on $ \Ccinf(\Gamma \bs G, \phi) $ is given.

A net $ f_\alpha $ in $ \Ccinf(\Gamma \bs G, \phi) $ converges to $ f \in \Ccinf(\Gamma \bs G, \phi) $ if there is a compact subset $ C $ of $G$ and $ \beta $ such that $ \supp f_\alpha \cup \supp f \subset \Gamma C $ for all $ \alpha > \beta $ and $ f_\alpha $ converges uniformly with all derivatives on $ C $.

\begin{prop}\label{Prop Ccinf(Gamma|G) is dense}
    The inclusion of $ \Ccinf(\Gamma \bs G, \phi) $ into $ \CS(\Gamma \bs G, \phi) $ is continuous with dense image.
\end{prop}

\begin{proof}
    Note first of all that $ \Ccinf(\Gamma \bs G, \phi) $ is contained in $ \CS(\Gamma \bs G, \phi) $. Let $ U \in \U_\Gamma $, $ d \geq 0 $ and $ X, Y \in \U(\g) $.
    It follows from the definition of the topology on $ \Ccinf(\Gamma \bs G, \phi) $ that the maps $ \LOI{p}{U}_{d, X, Y} $ are continuous seminorms on $ \Ccinf(\Gamma \bs G, \phi) $.
    Thus, $ \Ccinf(\Gamma \bs G, \phi) $ is continuously injected into $ \CS(\Gamma \bs G, \phi) $.

    Let us prove now that $ \Ccinf(\Gamma \bs G, \phi) $, viewed as a subspace of $ \CS(\Gamma \bs G, \phi) $, is dense in $ \CS(\Gamma \bs G, \phi) $.
    To show this we do a similar argument as in the proof of Theorem 7.1.1 in \cite[p.227]{Wallach}.

    Let $ h \in \Cinf(\RR) $ be such that $ h(x) \in [0, 1] $ for all $ x \in \RR $, $ h(x) = 1 $ if $ x \in [-1, 1] $ and $ h(x) = 0 $ if $ x \in \RR : |x| \geq 2 $. For $ r > 0 $, set $ u_r(g) = h(\tfrac{\log(a_g)}{r}) $. Then, there is a neighbourhood $V_r$ ($ r > 0 $) of $K$ such that $ u_r(g) = 1 $ for all $ g \in V_r $. Thus, $ u_r $ is smooth.

    As every smooth function is locally Lipschitz, as $ h^{(n)}(0) = 0 $ for all $ n \in \NN $ and as $ h^{(n)} $ ($ n \in \NN $) has compact support, there is a constant $ c_n > 0 $ such that
    \begin{equation}\label{eq |h^(n)(x)| <= c_n |x|}
        |h^{(n)}(x)| \leq c_n |x|
    \end{equation}
    for all $ x \in \RR $ and all $ n \in \NN $. As $ h(0) = 1 $, there is a constant $ c_0 > 0 $ such that $ |h(x) - 1| \leq c_0 |x| $ for all $ x \in [-2, 2] $. As $ |h(x)-1| = 1 $ for all $ x \in \RR : |x| \geq 2 $, the previous estimate holds also for $ x \in \RR $ such that $ |x| \geq 2 $ if we choose $ c_0 \geq \frac12 $. So,
    \begin{equation}\label{eq |h(x) - 1| <= c_0 |x|}
        |h(x) - 1| \leq c_0 |x|
    \end{equation}
    for all $ x \in \RR $. Let $ X_1, X_2 \in \U(\g) $. It follows from \eqref{eq |h^(n)(x)| <= c_n |x|}, \eqref{eq |h(x) - 1| <= c_0 |x|}, Corollary \ref{Cor Estimate of L_X R_Y (1 + log a_g)} and the chain rule that there is a constant $ c'_{X_1, X_2} > 0 $ such that
    \begin{equation}\label{eq Estimate for |L_(X_1) R_(X_2) u_r(g)|}
        |L_{X_1} R_{X_2} (u_r(g) - 1)| \leq \frac{c'_{X_1, X_2}}{r} \big( 1 + \log(a_g) \big) \qquad (g \in G) \ .
    \end{equation}
    By Remark \ref{Rem of Lem Cut-off function in the convex-cocompact case} of Lemma \ref{Lem Cut-off function in the convex-cocompact case},
    \[
        S := \{\gamma \in \Gamma \mid \gamma^{-1} U \cap \supp \chi \neq \emptyset \}
    \]
    is finite. Let $ g \in G $ such that $ gK \in U $.
    Then, $ |\sum_{\gamma \in \Gamma} (u_r \chi)(\gamma^{-1} g) - 1| $ is less or equal than
    \[
         \sum_{\gamma \in S} |u_r(\gamma^{-1} g) - 1| \chi(\gamma^{-1} g) \ .
    \]
    By \eqref{eq Estimate for |L_(X_1) R_(X_2) u_r(g)|} and as $ \chi $ takes values in $ [0, 1] $, this is again less or equal than
    \[
        c'_{1, 1} \sum_{\gamma \in S} \tfrac{1 + \log(a_{\gamma^{-1} g})}{r}
    \]
    ($ c'_{1, 1} = c_0 $). By Lemma \ref{Lem 1 + log a_h <= (1 + log a_gh)(1 + log a_g)}, we can estimate this further by
    \[
        c_{1, 1} \cdot \tfrac{1 + \log(a_g)}{r} \ ,
    \]
    where $ c_{1, 1} := c'_{1, 1} \sum_{\gamma \in S} (1 + \log(a_\gamma)) $ (this constant depends on $U$).

    It follows from \eqref{eq Estimate for |L_(X_1) R_(X_2) u_r(g)|}, Corollary \ref{Cor Estimate of L_X R_Y (1 + log a_g)}, Lemma \ref{Lem Cut-off function in the convex-cocompact case} and the Leibniz rule that for all $ X, Y \in \U(\g) $ such that $ X \in \g \U(\g) $ or $ Y \in \g \U(\g) $, there exists a positive constant $ c_{X, Y} $ (depending also on $U$) such that
    \[
        |L_X R_Y \big( \sum_{\gamma \in \Gamma} L_\gamma (u_r \chi) \big)(g)|
        = |\sum_{\gamma \in S} L_{\Ad(\gamma)^{-1}X} R_Y (u_r \chi)(\gamma^{-1} g)| \leq \frac{c_{X, Y}}{r} \big( 1 + \log(a_g) \big)
    \]
    for all $ r > 0 $ and all $ g \in G $ such that $ gK \in U $.

    Since $ u_r \chi \in \Ccinf(G) $, $ \sum_{\gamma \in \Gamma} L_\gamma (u_r \chi) \in \Ccinf(\Gamma \bs G) $.

    Let $ X, Y \in \U(\g)_{l} $ and let $ \{X_i\} $ be a basis of $ \U(\g)_{l} $.
    Then, there exists a constant $ c > 0 $ (depending on $ X $, $Y$ and $U$) such that
    \[
        \LOI{p}{U}_{d, X, Y}\bigl( \sum_{\gamma \in \Gamma} L_\gamma (u_r \chi) f - f\bigr) \leq \frac{c}{r} \sum_{i, j} \LOI{p}{U}_{d + 1, X_i, X_j}(f)
    \]
    for all $ r > 0 $. It follows that $ \sum_{\gamma \in \Gamma} L_\gamma (u_r \chi) f $ converges to $f$ in $ \CS(\Gamma \bs G, \phi) $.
\end{proof}

\newpage

\subsubsection{The space of cusp forms on \texorpdfstring{$ \Gamma \bs G $}{Gamma\textbackslash G}}

In this section, we define the constant term of a Schwartz function, we show that it is well-defined and we determine some basic properties of it.
The space of cusp forms is then the space of Schwartz functions having vanishing constant term.

For a set $U$ in $ \dX $, put\index[n]{GU@$ G(U) $, $ K(U) $}
\[
    G(U) = \{ g \in G \mid gP \in U \} \quad \text{and} \quad K(U) = \{ k \in K \mid kM \in U \} \ .
\]

\begin{dfn}
    For $ g \in G(\Omega), h \in G $ and for a measurable function $ f \colon G \to V_\phi $, we define
    \[
        f^\Omega(g, h) = \int_N{ f(gnh) \, dn }
    \]
    whenever the integral exists.
    We call this function the \textit{constant term} of $f$.\index{constant term}\index[n]{fOmega@$ f^\Omega $}
\end{dfn}

\begin{rem}
    We will show in Proposition \ref{Prop f^Omega well-defined} that the constant term is well-defined if $ f \in \CS(\Gamma \bs G, \phi) $.
\end{rem}

Let $ E = \exp(\e_1 \oplus \e_2) $ for some vector subspaces $ \e_1 \subset \g_{\alpha} $ and $ \e_2 \subset \g_{2\alpha} $.
It follows that $E$ is invariant under conjugation in $A$.
Let $ n_1 = \dim_\RR(\e_1) $, let $ n_2 = \dim_\RR(\e_2) $ and let $ \xi_E = \frac{n_1 + 2 n_2}{m_{\alpha} + 2 m_{2\alpha}} $.

Recall that $ m_{\alpha} $ (resp. $ m_{2\alpha} $) denotes the multiplicity of the root $ \alpha $ (resp. $ 2\alpha $).

For the moment we are only interested in the case where $ E = N $ (then $ \xi_E = 1 $). We state the following two results however in greater generality as we need it under this form in the geometrically finite case.

\begin{lem}\label{Lem int_E f(n) dn = int_A int_(varphi(x) = 1) f(x^a) dx a^(2 xi_E rho) da}
    For $ x \in N $, set $ \varphi(x) = a_B(x w) $. Then, there is a volume form (a nonvanishing differential form of highest degree) $dx$ on $ \{x \in E \mid \varphi(x) = 1 \} $ such that
    \[
        \int_E f(n) \, dn = \int_A \int_{\{x \in E \mid \varphi(x) = 1 \}} f(x^a) \, dx \ a^{2 \xi_E \rho} \, da
    \]
    for all $ f \in \C_c(E) $.
\end{lem}

\begin{rem}
    It follows from standard integration theory that the above formula holds also if $ f \in L^1(E) $ or if $ f $ is a nonnegative measurable function on $E$.
\end{rem}

\begin{proof}
    The lemma is standard. We provide a proof for the convenience of the reader.
    As hypersurfaces are orientable, $ \{x \in E \mid \varphi(x) = 1 \} $ has a volume form $dy$.
    Let $ n \in N $ and $ a \in A $. Then, $ n = \bar{n}_B(n w) m_B(n w) a_B(n w) n_B(n w) w $.
    Thus,
    \[
        n^a = (a \bar{n}_B(n w) a^{-1}) m_B(n w) a a_B(n w) a (a^{-1} n_B(n w) a) w \ .
    \]
    Hence,
    \[
        \varphi(n^a) = a_B(n^a w) = a^2 a_B(n w) = a^2 \varphi(n) \ .
    \]
    The map
    \[
        \Phi \colon \{x \in E \mid \varphi(x) = 1\} \times A \to E \smallsetminus \{e\} , \ (y, a) \mapsto y^{a^{-1}}
    \]
    is clearly a diffeomorphism. Its inverse sends $ n \in E \smallsetminus \{e\} $ to $ (n^a, a) $, where $ a = a_B(n w)^{-\frac12} $.

    So, $ dn = g(y, a) \, dy \, da $ for some positive measurable function $g$ on $ E \times A $ ($g$ is smooth on $ (E \smallsetminus \{e\}) \times A $).
    We have
    \begin{multline*}
        \int_E f(n^{(a')^{-1}}) \, dn = \int_A \int_{\{x \in E \mid \varphi(x) = 1 \}} f((y^{(a')^{-1}})^a) \, g(y, a) \, dy \, da \\
        = \int_A \int_{\{x \in E \mid \varphi(x) = 1 \}} f(y^a) \, g(y, a a') \, dy \, da \ .
    \end{multline*}
    Let $ y \in E $. Since $ (a')^{2 \xi_E \rho} \int_E f(n) \, dn = \int_E f(n^{(a')^{-1}}) \, dn $, $ (a')^{2 \xi_E \rho} g(y, a) = g(y, a a') $ for all $ a, a' \in A $.
    Taking $ a = 1 $ shows that $ g(y, a') = (a')^{2 \xi_E \rho} g(y, 1) $.
    The lemma follows.
\end{proof}


\begin{prop}\label{Prop N bar integral is finite}
    The integral
    \begin{equation}\label{eq int_(theta(E)) a(bar n)^(-xi_E rho) (1 + log(a(bar n)))^(-d) dbar n}
        \int_{\theta(E)}{ a(\bar{n})^{-\xi_E \rho} (1 + \log a(\bar{n}))^{-d} \, d\bar{n} }
    \end{equation}
    is finite if and only if $ d > 1 $.
    Here, $ d\bar{n} $ denotes the push-forward of the Lebesgue measure on $ \theta(\e_1) \oplus \theta(\e_2) $.
\end{prop}

\begin{rem} \nlenum
    \begin{enumerate}
    \item It follows from this proposition and Corollary \ref{Cor a_n asymp a(theta n)} that
        \[
            \int_{\theta(E)}{ a(\bar{n})^{-\xi_E \rho} (1 + \log a_{\bar{n}})^{-d} \, d\bar{n} }
        \]
        is also finite if and only if $ d > 1 $.
    \item Compare with Theorem 4.5.4 of \cite[p.126]{Wallach}.
    \end{enumerate}
\end{rem}

\begin{proof}
    Let $ S = \{ \exp(X + Y) \mid X \in \g_\alpha, Y \in \g_{2\alpha} : \tfrac14|X|^4 + 2 |Y|^2 \leq 1 \} $ (compact).

    Since $ a(\theta(n)) \asymp a_B(n w) $ on $ N \smallsetminus S $ by Corollary \ref{Cor a_n asymp a(theta n)}, \eqref{eq int_(theta(E)) a(bar n)^(-xi_E rho) (1 + log(a(bar n)))^(-d) dbar n} is finite if and only if
    \[
        \int_{E \smallsetminus S} a_B(n w)^{-\xi_E \rho} (1 + \log a_B(n w))^{-d} \, dn \ .
    \]
    By Lemma \ref{Lem int_E f(n) dn = int_A int_(varphi(x) = 1) f(x^a) dx a^(2 xi_E rho) da}, this is again finite if and only if
    \[
        \int_{\bar{A}_+} \int_{\{x \in E \mid a_B(x w) = 1 \}} a^{-2\xi_E \rho} (1 + 2 \log a)^{-d} \, dx \ a^{2 \xi_E \rho} \, da
    \]
    is finite. The above integral is again equal to
    \[
        \vol(\{x \in E \mid a_B(x w) = 1 \}) \int_{\bar{A}_+} (1 + 2 \log a)^{-d} \, da \ .
    \]
    The proposition follows as $ \int_{\bar{A}_+} (1 + 2 \log a)^{-d} \, da $ is finite if and only if $ d > 1 $.
\end{proof}

\needspace{2\baselineskip}
In order to show that the constant term is well-defined, we first need the following lemma.

\begin{lem}\label{Lem Estimate for int_N |f(g n h)| dn} \nlenum
    \begin{enumerate}
    \item Let $d > 1$, $ t \geq 1 $ and $ g, h \in G $. Then, $ n \mapsto a_{g n h}^{-t \rho} (1 + \log a_{g n h})^{-d} $ is integrable over $ N $. Moreover, if $ \eps \in [0, d - 1) $, then there exists a positive constant $ C $, depending only on $ d $ and $ \eps $, such that
        \[
            \int_{N}{ a_{g n h}^{-t \rho} (1 + \log a_{g n h})^{-d} \, dn }
            \leq C a(g)^{-t \rho} a(h^{-1})^{-t \rho} (1 + |\log a(g)a(h^{-1})^{-1}|)^{-\eps} \ .
        \]
    \item Let $d > 1$, $ t \geq \xi_E $ and $ g, h \in G $ be such that $ n(g) = n(h^{-1}) = e $. Then, $ n \mapsto a_{g n h}^{-t \rho} (1 + \log a_{g n h})^{-d} $ is integrable over $ E $. Moreover, if $ \eps \in [0, d - 1) $, then there exists a positive constant $ C $, depending only on $ d $ and $ \eps $, such that
        \[
            \int_{E}{ a_{g n h}^{-t \rho} (1 + \log a_{g n h})^{-d} \, dn }
            \leq C a(g)^{-t \rho} a(h^{-1})^{-t \rho} (1 + |\log a(g)a(h^{-1})^{-1}|)^{-\eps} \ .
        \]
    \end{enumerate}
\end{lem}

\begin{proof}
    This proof is inspired by the one of Theorem 7.2.1 in \cite[p.231]{Wallach}.

    Fix $ g, h \in G $.
    Let $ \eps \in [0, d - 1) $, $ n \in N $ and $ a \in A $. Put $ v = \theta(n) $. By formula \eqref{a(g^(-1)k) <= a_g},
    \[
        a(v) = a((va^{-1})a) = a(va^{-1}) a \leq a_{va^{-1}} a = a_{na} a \ .
    \]
    Thus, $ a_{na} \geq a^{-1} a(v) $. As in addition $ \eps \in [0, d - 1) $ and the following hold by Lemma \ref{Lem a_(na) estimates}
    \begin{enumerate}
    \item $ \log a_{na} \geq |\log a| \ $,
    \item $ 1 + \log a_{na} = 1 + \frac12 \log a_{na}^2 \geq 1 + \frac12 \log(a_n) \geq \frac12(1 + \log a_v) \ $,
    \end{enumerate}
    we have
    \begin{align}\label{Proof Lem Estimate of int_N |f(g n h)| dn}
        a_{na}^{-t \rho} (1 + \log a_{na})^{-d}
        &\leq a^{t \rho} a(v)^{-t \rho} (1 + \log a_{na})^{-\eps} (1 + \log a_{na})^{-d+\eps} \\
        &\leq 2^d a^{t \rho} a(v)^{-t \rho} (1 + |\log a|)^{-\eps} (1 + \log a_v)^{-d+\eps} \ . \nonumber
    \end{align}
    As $ g n h = \kappa(g) \Big(a(g) \big(n(g) n n(h^{-1})^{-1}\big) a(g)^{-1}\Big) a(g) a(h^{-1})^{-1} \kappa(h^{-1})^{-1} $, as $E$ is invariant under conjugation in $ A $ and as $ n(g) = n(h^{-1}) = e $ or $ E = N $,
    \[
        \int_{E}{ a_{g n h}^{-t \rho} (1 + \log a_{g n h})^{-d+\eps} \, dn }
    \]
    is equal to
    \[
        a(g)^{-2t\rho} \int_E{ a_{n a(g) a(h^{-1})^{-1}}^{-t \rho} (1 + \log a_{n a(g) a(h^{-1})^{-1}})^{-d+\eps} \, dn } \ .
    \]
    By \eqref{Proof Lem Estimate of int_N |f(g n h)| dn} and as $ \int_{\theta(E)}{ a(\bar{n})^{-\rho} (1 + \log a_{\bar{n}} )^{-d+\eps} \, d\bar{n} } $ is finite for all $ d > 1 + \eps $ by Proposition \ref{Prop N bar integral is finite}, this is less or equal than
    \[
        C a(g)^{t \rho} a(h^{-1})^{-t \rho} (1 + |\log a(g)a(h^{-1})^{-1}|)^{-\eps} a(g)^{-2t\rho}
    \]
    for some constant $C > 0$ depending only on $d$ and $ \eps $.
    So,
    \[
        \int_{E}{ a_{g n h}^{-t \rho} (1 + \log a_{g n h})^{-d} \, dn } \leq C a(g)^{-t \rho} a(h^{-1})^{-t \rho} (1 + |\log a(g)a(h^{-1})^{-1}|)^{-\eps} \ .
    \]
\end{proof}

\begin{dfn}\label{Def of U_(g, h)}\index[n]{Ugh@$ U_{g, h} $}
    Let $ g \in G(\Omega) $ and $ h \in G $. Then, we denote by $ U_{g, h} $ an open set in $X$ that is a relatively compact neighbourhood of $ gP $ in $ X \cup \Omega $ and contains the horosphere $ gNhK = gNg^{-1}(gh)K $ passing through $ ghK \in X $ and $ gP \in \dX $.
\end{dfn}

\begin{prop}\label{Prop f^Omega well-defined}
    Let $ f \in \CS(\Gamma \bs G, \phi) $, $ g \in G(\Omega) $, $ h \in G $ and $ X, Y \in \U(\g) $. Then, $ n \mapsto L_X R_Y f(gnh) $ is integrable over $N$. Thus, the constant term $ f^\Omega $ is well-defined (take $ X = Y = 1 $). Moreover, if $ U \in \U_\Gamma $ and if $ r > 1 $, then there exists $ c > 0 $ (depending only on $r$) such that
    \[
        \int_N{ |L_X R_Y f(g n h)| \, dn } \leq c \LOI{p}{U}_{2r, X, Y}(f) a(g)^{-\rho} a(h^{-1})^{-\rho} (1 + |\log a(g)a(h^{-1})^{-1}|)^{-r}
    \]
    for all $ g \in G(\Omega) $ and $ h \in G $ such that $ g n h K \in U $ for all $ n \in N $.
\end{prop}

\begin{proof}
    Let $ f \in \CS(\Gamma \bs G, \phi) $ and $ X, Y \in \U(\g) $.
    Let $ g \in G(\Omega) $ and $ h \in G $. Then, $ g n h \in U_{g, h} $ for all $ n \in N $.
    Let now $ U \in \U_\Gamma $, $ g \in G(\Omega) $ and $ h \in G $ be as above.
    For all $ d \geq 0 $, we have
    \[
        |L_X R_Y f(x)| \leq \LOI{p}{U}_{d, X, Y}(f) a_x^{-\rho} (1 + \log a_x)^{-d} \qquad (xK \in U_{g, h})
    \]
    by definition of $ \CS(\Gamma \bs G, \phi) $. Thus, for $ r > 1 $,
    \[
        \int_N{ |L_X R_Y f(g n h)| \, dn } \leq \LOI{p}{U}_{2r, X, Y}(f) \int_{N}{ a_{g n h}^{-\rho} (1 + \log a_{g n h})^{-2r} \, dn } \ .
    \]
    By Lemma \ref{Lem Estimate for int_N |f(g n h)| dn} applied to $ E = N $, there exists a constant $ c > 0 $ (depending only on $r$) such that this is again less or equal than
    \[
        c \LOI{p}{U}_{2r, X, Y}(f) a(g)^{-\rho} a(h^{-1})^{-\rho} (1 + |\log a(g)a(h^{-1})^{-1}|)^{-r} \ .
    \]
\end{proof}

\begin{cor}\label{Cor f mapsto f^Omega continuous}
    Fix $ g \in G(\Omega) $ and $ h \in G $. Then,
    \[
        f \in \CS(\Gamma \bs G, \phi) \mapsto f^\Omega(g, h) \in V_\phi
    \]
    is linear and continuous.
    Thus, $ \c{\CS}(\Gamma \bs G, \phi) $ is a closed subspace of $ \CS(\Gamma \bs G, \phi) $.
\end{cor}

\begin{proof}
    Observe that the map is linear.
    Fix $ g \in G(\Omega) $ and $ h \in G $. Let $ f \in \CS(\Gamma \bs G, \phi) $ and let $ r > 1 $. Then, by Proposition \ref{Prop f^Omega well-defined}, there exists a constant $ c > 0 $, depending on $ g $ and $ h $, such that
    \[
        |f^\Omega(g,h)| \leq c \LOI{p}{U_{g, h}}_{2r, 1, 1}(f) \ .
    \]
    Hence, $ f \mapsto f^\Omega(g, h) $ is continuous when $ g $ and $ h $ are fixed.
\end{proof}

\begin{cor}\label{Cor f^Omega(g, h) depends smoothly on g and h}
    Fix $ f \in \CS(\Gamma \bs G, \phi) $. For $ g \in G(\Omega) $, $ h \in G $ and $ X, Y \in \U(\g) $, we have
    \[
        L_X R_Y \int_N{ f(gnh) \, dn } = \int_N{ L_X R_Y f(gnh) \, dn } \ .
    \]
    In particular, $ f^\Omega(g, h) $ depends smoothly on $ g $ and $ h $.
\end{cor}

\begin{proof}
    Let us check the assumptions of the theorem about differentiation of parameter dependent integrals.
    Fix $ X, Y \in \U(\g) $ and $ r > 1 $.
    By Proposition \ref{Prop f^Omega well-defined}, we know that $ n \in N \mapsto f(g n h) $ is integrable.
    The function $ g \in G \mapsto f(g n h) $ is smooth for every $ n \in N $, $ h \in G $ as $f$ is smooth.

    Let $ V $ be an open subset of $G$ which is relatively compact in $ G(\Omega) $ and let $W$ be an open, relatively compact subset of $G$. Choose $U_{V, W} \in \U_\Gamma$ such that $ U_{g, h} \subset U_{V, W} $ for all $ g \in V $ and all $ h \in W $. By the proof of Proposition \ref{Prop f^Omega well-defined} and of Lemma \ref{Lem Estimate for int_N |f(g n h)| dn} applied to $ E = N $, there is a constant $ C > 0 $ such that
    \begin{multline}\label{Proof Cor f^Omega(g, h) depends smoothly on g and h eq1}
        |L_X R_Y f(g nh)|
        \leq C \LOI{p}{U_{V, W}}_{2r, X, Y}(f) a(\theta n(g) \theta n \theta n(h^{-1})^{-1})^{-\rho} \\
        (1 + \log a_{\theta n(g) \theta n \theta n(h^{-1})^{-1}})^{-r}
        a(g)^{-\rho} a(h^{-1})^{-\rho} (1 + |\log a(g)a(h^{-1})^{-1}|)^{-r}
    \end{multline}
    for all $ g \in V $, $ h \in W $ and $ n \in N $.
    By Lemma \ref{Lem 1 + log a_h <= (1 + log a_gh)(1 + log a_g)} and by Corollary \ref{Cor a_n asymp a(theta n)}, there exists a constant $ C' > 0 $ such that
    \[
        a(\theta n(g) \theta n \theta n(h^{-1})^{-1})^{-\rho} (1 + \log a_{\theta n(g) \theta n \theta n(h^{-1})^{-1}})^{-r}
        \leq C' a(\theta n)^{-\rho} (1 + \log a_{\theta n})^{-r}
    \]
    for all $ g \in V $, $ n \in N $ and $ h \in W $.
    It follows that there is a constant $C'' > 0$, depending only on $ V $ and $W$, such that \eqref{Proof Cor f^Omega(g, h) depends smoothly on g and h eq1} is bounded by
    \[
        C'' a(\theta n)^{-\rho} (1 + \log a_{\theta n})^{-r} \qquad (n \in N) \ .
    \]
    This is integrable over $ N $ by Proposition \ref{Prop N bar integral is finite}.
    Thus, $ |L_X R_Y f(g nh)| $ is uniformly bounded by an integrable function for all $ g \in V $ and all $ h \in W $.
    Hence, $ L_X R_Y \int_N{ f(gnh) \, dn} $ exists and is equal to $ \int_N{ L_X R_Y f(gnh) \, dn } $.
    The corollary follows.
\end{proof}

\needspace{2\baselineskip}
\begin{prop}\label{Prop equivariances in the convex-cocompact case}
    Let $ f \in \CS(\Gamma \bs G, \phi) $.
    \begin{enumerate}
    \item For $ g \in G(\Omega) $, $ h \in G $, $ \gamma \in \Gamma $, we have\label{Prop equivariances in the convex-cocompact case eq1}
        \[
            f^\Omega(\gamma g, h) = \phi(\gamma) f^\Omega(g, h) \ .
        \]
    \item For $ g \in G(\Omega) $, $ h, x \in G $, we have\label{Prop equivariances in the convex-cocompact case eq2}
        \[
            (R_x f)^\Omega(g, h) = f^\Omega(g, hx) \ .
        \]
    \item For $ g \in G(\Omega) $, $ h \in G $, $ a \in A $, $ m \in M $, we have
        \[
            f^\Omega(g am, h) = a^{-2\rho} f^\Omega(g, am h) \ .
         \]
    \item $ f^\Omega(gn_1, n_2h) = f^\Omega(g, h) $ for all $ g \in G(\Omega), h \in G, n_1, n_2 \in N $.
    \end{enumerate}
\end{prop}

\needspace{2\baselineskip}
\begin{proof}
    Let $ f \in \CS(\Gamma \bs G, \phi) $.
    \begin{enumerate}
    \item For $ g \in G(\Omega) $, $ h \in G $, $ \gamma \in \Gamma $, we have
        \[
            f^\Omega(\gamma g, h) = \int_{N}{ f(\gamma gnh) \, dn } = \phi(\gamma) \int_{N}{ f(gnh) \, dn } = \phi(\gamma) f^\Omega(g, h) \ .
        \]
    \item For $ g \in G(\Omega) $, $ h, x \in G $, we have
        \[
            (R_x f)^\Omega(g, h) = \int_N{ f(g n h x) \, dn } = f^\Omega(g, hx) \ .
        \]
    \item For $ g \in G(\Omega) $, $ h \in G $, $ a \in A $, $ m \in M $, we have
        \begin{align*}
            f^\Omega(g am, h) &= \int_{N}{ f(g am n (am)^{-1} am h) \, dn } \\
            &= a^{-2\rho} \int_{N}{ f(gn am h) \, dn } = a^{-2\rho} f^\Omega(g, am h)
         \end{align*}
         as $ n \mapsto m^{-1} n m $ is a measure-preserving diffeomorphism and as $ d(a^{-1}na) = a^{-2\rho} dn $.
    \item For $ g \in G(\Omega) $, $ h \in G $, $ n_1, n_2 \in N $,
        \[
            f^\Omega(gn_1, n_2h) = \int_{N}{ f(gn_1nn_2h) \, dn } = \int_{N}{ f(gnh) \, dn } = f^\Omega(g, h)
        \]
        as $N$ is unimodular.
    \end{enumerate}
\end{proof}

\begin{dfn}\label{Dfn Harish-Chandra space of cusp forms in the convex-cocompact finite case}
    Define the following subspace of $ \CS(\Gamma \bs G, \phi) $:
    \[
        \c{\CS}(\Gamma \bs G, \phi) = \{ f \in \CS(\Gamma \bs G, \phi) \mid f^\Omega(g, h) = 0 \quad \forall g \in G(\Omega), h \in G \} \ .
    \]
    We call it the \textit{space of cusp forms} on $ \Gamma \bs G $.\index{space of cusp forms}\index[n]{CGammaG0@$ \c{\CS}(\Gamma \bs G, \phi) $}
\end{dfn}

\begin{rem} \nlenum
    \begin{enumerate}
    \item This subspace may be considered as the analog of the Harish-Chandra space of cusp forms on $ G $.
    \item Since $ f^\Omega(\gamma g, h) = \phi(\gamma) f^\Omega(g, h) $ for all $ \gamma \in \Gamma, g \in G(\Omega), h \in G $ by Proposition \ref{Prop equivariances in the convex-cocompact case} \eqref{Prop equivariances in the convex-cocompact case eq1}, it suffices to check the condition $ f^\Omega(g, h) = 0 $ for $ g $ belonging to a fundamental set of $ G(\Omega) $ for $ \Gamma $.
    \end{enumerate}
\end{rem}

\subsubsection{The right regular representation of \texorpdfstring{$G$}{G} on the Schwartz space}

We show that the right translation by elements of $G$ on the Schwartz space, respectively the space of cusp forms, is a representation of $G$ (see Theorem \ref{Thm representations convex-cocompact case}).

\medskip

\begin{lem}\label{Lem U in U_Gamma iff clo(U) cap dX subset Omega compact}
    Let $ U $ be an open subset of $X$. Then, $ U $ is relatively compact in $ X \cup \Omega $ if and only if $ \clo(U) \cap \dX $ is a compact subset of $ \Omega $.
\end{lem}

\begin{proof}
    Let $ U $ be an open subset of $X$. If $ U $ is relatively compact in $ X \cup \Omega $, then $ \clo(U) $ is compact in $ X \cup \Omega $. Hence, $ \clo(U) \cap \dX $ is a compact subset of $ \Omega $.

    Assume now that $ \clo(U) \cap \dX $ is a compact subset of $ \Omega $. Let $ V $ be an open, relatively compact neighbourhood of $ \clo(U) \cap \dX \subset K/M $ which is also contained in $ \Omega $.
    Let $ \tilde{V} = \{ k a K \mid kM \in V , \, a \in A_+ \} \cup V $. Then, $ \tilde{V} $ is relatively compact in $ X \cup \Omega $ and a neighbourhood of $ \clo(U) \cap \dX $ in $ X \cup \Omega $. Thus, $ U \smallsetminus \tilde{V} $ is relatively compact in $X$. Hence, $ U $ is relatively compact in $ X \cup \Omega $.
\end{proof}

Let $f$ be a function on $ G $. For $ g, x \in G $, set\index[n]{Lgf Rgf@$ L_g f $, $ R_g f $}
\[
    L_g f(x) = f(g^{-1} x) \qquad \text{and} \qquad R_g f(x) = f(x g) \ .
\]
For $ f \in \CS(\Gamma \bs G, \phi) $, define
\[
    \pi(g) f(x) := R_g f(x) = f(xg)  \qquad (g, x \in G) \ .
\]
Clearly, $ \pi(g)(\CS(\Gamma \bs G, \phi)) \subset \Cinf(\Gamma \bs G, \phi) $ for all $ g \in G $.

\needspace{2\baselineskip}
For $ h \in G $ and for a subset $ U $ of $X$, put $ U^h = \{ ghK \mid gK \in U \} $.

\begin{lem}\label{Lem U^h open rel. compact in X u Omega}
    For every $ h \in G $ and $ U \in \U_\Gamma $, $ \clo(U^h) \cap \dX $ is equal to $ \clo(U) \cap \dX $ and $ U^h $ belongs also to $ \U_\Gamma $.
\end{lem}

\begin{rem}\label{Rem of Lem U^h open rel. compact in X u Omega}
    The set $ \{ hgK \mid gK \in U \} $ ($ h \in G $) belongs also to $ \U_\Gamma $ if $h$ is sufficiently close to $e$.
\end{rem}

\begin{proof}
    Fix $ h \in G $ and $ U \in \U_\Gamma $.
    Consider a sequence $ g_i K $ in $ U $ converging to some element $ kM $ in $ \dX $. Without loss of generality, we may assume that $ h_{g_i} = h_h $. Then, $ a_{g_i} $ tends to $ \infty $ and $ k_{g_i}M $ converges to $ kM \in \dX $. So, $ k_{g_i} a_{g_i} a_h^{-1} k_h^{-1} h K = g_i K $ converges also to $ kM $. Note that the limit does not change if we consider $ k_{g_i} a_{g_i} k_h^{-1} h K $.
    But $ (k_{g_i} a_{g_i} k_h^{-1})K = k_{g_i} a_{g_i} K = g_i K \in U $ for all $i$. Thus, $ kM $ belongs also to $ \clo(U^h) \cap \dX $.
    Hence, $ \clo(U) \cap \dX $ is contained in $ \clo(U^h) \cap \dX $.

    As the above is also true for $ U^h $ instead of $U$ and $ h^{-1} $ instead of $h$ and as $ U = (U^h)^{h^{-1}} $,
    \[
        \clo(U^h) \cap \dX = \clo(U) \cap \dX \subset \Omega \ .
    \]
    It follows now from Lemma \ref{Lem U in U_Gamma iff clo(U) cap dX subset Omega compact} that $ U^h $ is an open subset of $X$ which is relatively compact in $ X \cup \Omega $, i.e. $ U^h \in \U_\Gamma $.
\end{proof}

\begin{cor}\label{Cor union of U^h in U_Gamma}
    Let $ U \in \U_\Gamma $ and let $ S $ be a relatively compact subset of $ G $. Then, $ \bigcup_{h \in S} U^h $ belongs to $ \U_\Gamma $.
\end{cor}

\begin{rem}\label{Rem of Cor union of U^h in U_Gamma}
    The set $ \bigcup_{h \in S} \{ hgK \mid gK \in U \} $ only belongs to $ \U_\Gamma $ if $S$ is contained in a sufficiently small neighbourhood of $e$.
\end{rem}

\begin{proof}
    Let $ U $ and $ S $ be as above. Put $ U^S = \bigcup_{h \in S} U^h $. As $ U^h \in \U_\Gamma $ and as $ \clo(U^h) \cap \dX = \clo(U) \cap \dX $ for all $ h \in G $ by Lemma \ref{Lem U^h open rel. compact in X u Omega}, $ \clo(U) \cap \dX $ is contained in $ \clo(U^S) \cap \dX $ and $ U^S $ is open in $X$ as it is the union of open sets in $X$. Since $ \{a_h \mid h \in S \} $ is bounded, the second part of the proof of Lemma \ref{Lem U^h open rel. compact in X u Omega} can be generalised to $ U^S $. Thus, $ \clo(U^S) \cap \dX $ is equal to $ \clo(U) \cap \dX $. Hence, $ U^S \in \U_\Gamma $ by Lemma \ref{Lem U in U_Gamma iff clo(U) cap dX subset Omega compact}.
\end{proof}

\begin{prop}\label{Prop convergence of Up_(r, X, Y)(R_h f - f)}
    Let $ U \in \U_\Gamma $, let $ X, Y \in \U(\g) $, let $ f \in \CS(\Gamma \bs G, \phi) $ and let $ r \geq 0 $. Then,
    \begin{enumerate}
    \item $ \LOI{p}{U}_{r, X, Y}(R_h f) $ is finite for every $ h \in G $ and $ \LOI{p}{U}_{r, X, Y}(L_h f) $ is finite when $ h \in G $ is sufficiently close to $e$;
    \item $ \LOI{p}{U}_{r, X, Y}(R_h f - f) $, $ \LOI{p}{U}_{r, X, Y}(L_h f - f) $ converge to zero when $ h $ tends to $e$.
    \end{enumerate}
\end{prop}

\begin{proof}
    To prove this proposition, we do a similar argument as in Lemma 14 of \cite[p.21]{HC66}.

    Let $ h \in G $, $ x \in G $, $ X, Y \in \U(\g) $, $ U \in \U_\Gamma $, and let $ f \in \CS(\Gamma \bs G, \phi) $.
    We have
    \begin{align*}
        \LOI{p}{U}_{r, X, Y}(\pi(h)f) &= \sup_{gK \in U}{ (1 + \log a_g)^r a_g^{\rho} |L_X R_Y (R_h f)(g)| } \\
        &\leq (1 + \log a_h)^r a_h^{\rho} \sup_{gK \in U}{ (1 + \log a_{gh})^r a_{gh}^{\rho} |L_X R_{\Ad(h^{-1})Y} f(gh)| } \\
        &= (1 + \log a_h)^r a_h^{\rho} \LOI{p}{U^h}_{r, X, \Ad(h^{-1})Y}(f) < \infty \ .
    \end{align*}
    Thus, $ \pi(h)f $ belongs to $ \CS(\Gamma \bs G, \phi) $.
    Let $ \{X_1, \dotsc, X_n\} $ be a basis of $ \g $. Let $ H \in \g $. Write $ H = \sum_{i=1}^n c_i X_i $.
    We clearly have
    \begin{align*}
        R_{\exp H} f(x) - f(x) &= \int_0^1 \restricted{ \frac{d}{ds} }{s=0} f(x \exp(tH) \exp(sH)) \, dt \\
        &= \int_0^1 \restricted{ \frac{d}{ds} }{s=0} f(x \exp(tH) \exp\bigl(s \sum_{i=1}^n c_i X_i\bigr)) \, dt \\
        &= \sum_{i=1}^n c_i \int_0^1 \restricted{ \frac{d}{ds} }{s=0} f(x \exp(tH) \exp(s X_i)) \, dt \ .
    \end{align*}
    Hence,
    \begin{multline*}
        |L_X R_Y (R_{\exp H} f)(x) - L_X R_Y f(x)| \\
        \leq \max_{1 \leq i \leq n} |c_i| \cdot \sum_{i=1}^n \int_0^1 |\restricted{ L_X R_Y \frac{d}{ds} }{s=0} f(x \exp(tH) \exp(s X_i))| \, dt \ .
    \end{multline*}
    Since all norms are equivalent on $ \RR^n $, there exists $ c_0 > 0 $ such that
    \[
        \max_{1 \leq i \leq n} |a_i| \leq c_0 \cdot |\sum_{i = 1}^n a_i X_i|
    \]
    for any $ a_i \in \RR $. So, $ |L_X R_Y (R_{\exp H} f)(x) - L_X R_Y f(x)| $ is less or equal than
    \[
         c_0 \cdot |H| \sum_{i=1}^n \int_0^1 |L_X R_Y R_{X_i} f(x \exp(tH))| \, dt \ .
    \]
    Consequently,
    \[
         \LOI{p}{U}_{r, X, Y}(R_h f - f)
         = \sup_{xK \in U}{ (1 + \log a_x)^r a_x^\rho |L_X R_Y (R_{\exp H} f)(x) - L_X R_Y f(x)| }
    \]
    is less or equal than
    \begin{multline*}
         c_0 \cdot |H| \sum_{i=1}^n \sup_{xK \in U}{ (1 + \log a_x)^r a_x^\rho \int_0^1 |L_X R_Y R_{X_i} f(x \exp(tH))| \, dt } \\
         \leq c_0 \cdot |H| \sum_{i=1}^n \sup_{t \in [0, 1]} \sup_{xK \in U} (1 + \log a_x)^r a_x^\rho |L_X R_Y R_{X_i} f(x \exp(tH))| \ .
    \end{multline*}
    Put $ V = \bigcup_{Z \in \g \, : \, |Z| \leq 1} U^{\exp(Z)} = \{ g \exp(Z)K \mid gK \in U, Z \in \g \, : \, |Z| \leq 1 \} $. Then, $ V $ belongs to $ \U_\Gamma $ by Corollary \ref{Cor union of U^h in U_Gamma}.
    If $ |H| \leq 1 $, then the above expression is less or equal than
    \[
        c \cdot |H| \sum_{i=1}^n \LOI{p}{V}_{r, X, Y + X_i}(f) \ ,
    \]
    where $ c := c_0 \max_{Z \in \g \, : \, |Z| \leq 1} (1 + \log a_{\exp(Z)})^r a_{\exp(Z)}^\rho $.
    The proof of the other assertion works analogously. The remarks \ref{Rem of Lem U^h open rel. compact in X u Omega} and \ref{Rem of Cor union of U^h in U_Gamma} tell us to what one has to pay attention.
\end{proof}

\begin{thm}\label{Thm representations convex-cocompact case}
    $ \big(\pi, \CS(\Gamma \bs G, \phi)\big) $ and $ \big(\pi, \c{\CS}(\Gamma \bs G, \phi)\big) $ are representations of $G$.
\end{thm}

\begin{proof}
    Since the constant term is $G$-equivariant with respect to the right translation, it is enough to prove that $ \big(\pi, \CS(\Gamma \bs G, \phi)\big) $ is a representation of $G$.
    Fix $ X, Y \in \U(\g) $, $ r \geq 0 $ and $ U \in \U_\Gamma $. Let $ f \in \CS(\Gamma \bs G, \phi) $ and $ h \in G $. Then, $ U^h \in \U_\Gamma $ by Lemma \ref{Lem U^h open rel. compact in X u Omega}.

    By Proposition \ref{Prop convergence of Up_(r, X, Y)(R_h f - f)}, $ \pi(h)f $ belongs to $ \CS(\Gamma \bs G, \phi) $.

    Assume that $ h_j \to h $ in $G$ and that $ f_j \to f $ in $ \CS(\Gamma \bs G, \phi) $. 
    Let $ N^h $ be a relatively compact neighbourhood of $h$ in $G$. Put
    \[
        \tilde{U}^h = \bigcup_{\tilde{h} \in N^h} U^{\tilde{h}} = \{ g \tilde{h}K \mid gK \in U, \tilde{h} \in N^h \} \ .
    \]
    This set belongs also to $ \U_\Gamma $ by Corollary \ref{Cor union of U^h in U_Gamma}.

    Let $ \eps > 0 $ be given.
    Let $ \{X_1, \dotsc, X_n\} $ be a basis of $ \g $. Define $ X_i^a $ to be
    \[
        \underbrace{ X_i \dotsm X_i }_{a \text{ factors} }
    \]
    if $ a \in \NN := \{1, 2, 3, \dotsc \} $ and to be $1$ if $ a = 0 $. For $ m \in \NN_0 := \{0, 1, 2, \dotsc \} $, define
    \[
        I_m = \{ \alpha := (\alpha_1, \dotsc, \alpha_n) \mid \alpha_i \in \NN_0 : \sum_{i=1}^n \alpha_i \leq m \} \ .
    \]
    Let $ Y \in \U(\g)_l $. Then, $ \Ad(x^{-1})Y \in \U(\g)_l $ as, by definition, $ \Ad(x^{-1})(Y_1 \otimes Y_2) = \Ad(x^{-1})(Y_1) \otimes \Ad(x^{-1})(Y_2) $.
    As $ \Ad(x^{-1})Y \in \U(\g) $ and as $ \Ad(x^{-1})Y $ depends smoothly on $x$, there are smooth functions $ c_{\alpha_1, \dotsc, \alpha_n} $ on $G$ such that
    \[
        \Ad(x^{-1})Y = \sum_{\alpha \in I_l}{ c_{\alpha_1, \dotsc, \alpha_n}(x) X_1^{\alpha_1} \cdots X_n^{\alpha_n} }
    \]
    by the Poincaré-Birkhoff-Witt theorem.
    As $ c_{\alpha_1, \dotsc, \alpha_n}(h_j) \to c_{\alpha_1, \dotsc, \alpha_n}(h) $ when $ j \to \infty $, $ \sup_{\alpha \in I_l} |c_{\alpha_1, \dotsc, \alpha_n}(h) - c_{\alpha_1, \dotsc, \alpha_n}(h_j)| < \frac{\eps}2 $ for $ j $ sufficiently large.
    We have
    \begin{multline}\label{Proof Prop representations convex-cocompact case eq1}
        \LOI{p}{U}_{r, X, Y}(\pi(h)f - \pi(h_j)f_j) \\
        \leq \LOI{p}{U}_{r, X, Y}(\pi(h)f - \pi(h_j)f) + \LOI{p}{U}_{r, X, Y}(\pi(h_j)(f - f_j)) \ ,
    \end{multline}
    where $ \LOI{p}{U}_{r, X, Y}(\pi(h_j)(f - f_j)) $ is less or equal than
    \begin{multline*}
         \max_{\tilde{h} \in N^h}\big( (1 + \log a_{\tilde{h} })^r a_{\tilde{h} }^{\rho} \big) \LOI{p}{\tilde{U}^h}_{r, X, \Ad(h_j^{-1})Y}(f - f_j) \\
         \leq \sum_{\alpha \in I_l} \max_{\tilde{h} \in N^h}\big( (1 + \log a_{\tilde{h}})^r a_{\tilde{h}}^{\rho} |c_{\alpha_1, \dotsc, \alpha_n}(\tilde{h})| \big) \LOI{p}{\tilde{U}^h}_{r, X, X_1^{\alpha_1} \cdots X_n^{\alpha_n}}(f - f_j) < \frac{ \eps }2
    \end{multline*}
    for $ j $ sufficiently large. Here, we have used that $ h_j \in N^h $ for $ j $ sufficiently large.
    By Proposition \ref{Prop convergence of Up_(r, X, Y)(R_h f - f)} applied to $ \pi(h)f \in \CS(\Gamma \bs G, \phi) $,
    \[
        \LOI{p}{U}_{r, X, Y}(\pi(h)f - \pi(h_j)f) < \frac{ \eps }2
    \]
    for $ j $ sufficiently large.
    It follows that \eqref{Proof Prop representations convex-cocompact case eq1} tends to zero when $ j $ goes to $ \infty $. This finishes the proof of the theorem.
\end{proof}

\newpage

\subsection{Tempered distribution vectors, square-integrable distribution vectors and Schwartz vectors}

First, we recall the notion of matrix coefficients of a representation $ \pi $ of $G$ (on a reflexive Banach space) are. Second, we use the matrix coefficient map in order to define tempered distribution vectors, square-integrable distribution vectors and Schwartz vectors.

\medskip

\index[n]{Vpi-infty, Vpiinfty, VpiK@$ V_{\pi, -\infty} $, $ V_{\pi, \infty} $, $ V_{\pi, K} $}
Let $ (\pi, V_\pi) $ be a representation of $G$ on a reflexive Banach space with dual representation $ (\pi', V_{\pi'}) $. Then, we define the \textit{space of distribution vectors} to be $ V_{\pi, -\infty} := (V_{\pi', \infty})' $ (strong dual). Here the subscript $ \infty $ means that we pass to the subspace of smooth vectors and we endow $ V_{\pi', \infty} $ with the canonical Fréchet topology. We have then the following inclusions of $ G $-representations
\[
    V_{\pi, \infty} \subset V_\pi \subset V_{\pi, -\infty} \ .
\]
We denote the subspace of $K$-finite and smooth vectors of $ V_\pi $ by $ V_{\pi, K} $.
Furthermore, we denote the \textit{space of tempered distributions} by $ \CS'(\Gamma \bs G, \phi) $. It is by definition the conjugate-linear strong dual of $ \CS(\Gamma \bs G, \phi) $.\index[n]{CGammaG'@$ \CS'(\Gamma \bs G, \phi) $}

\begin{dfn}
    Let $ (\pi, V_\pi) $ be an admissible $G$-representation of finite length on a reflexive Banach space.

    We say that a $ \Gamma $-invariant distribution vector $ f \in (V_{\pi, -\infty} \otimes V_\phi)^\Gamma $ is \textit{tempered} (resp. \textit{square-integrable}) if the function
    \[
        G \ni g \mapsto c_{f, v}(g) := \langle f, \pi'(g)v \rangle \in V_\phi \ ,
    \]
    called \textit{matrix coefficient} of $ \pi $, belongs to $ \CS'(\Gamma \bs G, \phi) $ (resp. $ L^2(\Gamma \bs G, \phi) $) for all $ v \in V_{\pi', K} $. We denote the linear subspace of square-integrable (resp. tempered) $ \Gamma $-invariant distribution vectors by $ (V_{\pi, -\infty} \otimes V_\phi)_d^\Gamma $ (resp. $ (V_{\pi, -\infty} \otimes V_\phi)^\Gamma_{temp} $).
    \index{matrix coefficient}\index{tempered distribution vector}\index{square-integrable distribution vector}\index[n]{cfv@$ c_{f, v} $}
    \index[n]{Vpi-inftyGammatemp@$ (V_{\pi, -\infty} \otimes V_\phi)^\Gamma_{temp} $}\index[n]{Vpi-inftyGammad@$ (V_{\pi, -\infty} \otimes V_\phi)^\Gamma_d $}
\end{dfn}

\begin{rem} \nlenum
    \begin{enumerate}
    \item This definition is equivalent to Definition 8.3 in \cite[p.130]{BO00} by Lemma 8.4 of \cite[p.130]{BO00}.
    \item If $ f \in V_{\pi, -\infty} \otimes V_\phi $ and $ v \in V_{\pi', \infty} $, then $ c_{f, v} $ is smooth as $ G \to V_{\pi', \infty}, \; g \mapsto \pi'(g)v $ is smooth and as $ f \colon V_{\pi', \infty} \to V_\phi $ is continuous and linear.
    \item As $ L^2(\Gamma \bs G, \phi) \subset \CS'(\Gamma \bs G, \phi) $, we have
        \[
            (V_{\pi, -\infty} \otimes V_\phi)_d^\Gamma \subset (V_{\pi, -\infty} \otimes V_\phi)^\Gamma_{temp} \subset (V_{\pi, -\infty} \otimes V_\phi)^\Gamma \ .
        \]
    \end{enumerate}
\end{rem}

We equip $ (V_{\pi, -\infty} \otimes V_\phi)^\Gamma $ with the subspace topology.

\begin{lem}\label{Lem matrix coeff in L2}
    Let $ v \in (V_{\pi, -\infty} \otimes V_\phi)^\Gamma $. If $ \pi $ is irreducible and if $ c_{v, v'} \in L^2(\Gamma \bs G, \phi) $ for some $ 0 \neq v' \in V_{\pi', \infty} $, then $ c_{v, v'} \in L^2(\Gamma \bs G, \phi) $ for all $ v' \in V_{\pi', \infty} $.
\end{lem}

\begin{proof}
    The proof works completely analogously to the proof of Lemma 1.3.2 in \cite[p.23]{Wallach}.
\end{proof}

\begin{dfn}\label{Dfn Schwartz vector}
    Let $ (\pi, V_\pi) $ be an admissible $G$-representation of finite length on a reflexive Banach space.

    We say that $ f \in (V_{\pi, -\infty} \otimes V_\phi)^\Gamma $ is a \textit{Schwartz vector} (resp. \textit{cuspidal Schwartz vector}) if $ c_{f, v} \in \CS(\Gamma \bs G, \phi) $ (resp. $ c_{f, v} \in \c{\CS}(\Gamma \bs G, \phi) $) for all $ v \in V_{\pi', K} $.
    We denote the \textit{space of Schwartz vectors} (resp. \textit{space of cuspidal Schwartz vectors}) by $ (V_{\pi, -\infty} \otimes V_\phi)_\CS^\Gamma $ (resp. $ (V_{\pi, -\infty} \otimes V_\phi)_{\c{\CS}}^\Gamma $).
    \index{Schwartz vector}\index{cuspidal Schwartz vector}\index[n]{Vpi-inftyGammaC@$ (V_{\pi, -\infty} \otimes V_\phi)_\CS^\Gamma $}\index[n]{Vpi-inftyGammaC0@$ (V_{\pi, -\infty} \otimes V_\phi)_{\c{\CS}}^\Gamma $}
\end{dfn}

\begin{rem}\label{Rem Schwartz vector}\nlenum
    \begin{enumerate}
    \item Assume further that $ \pi $ is irreducible. Since $ K $ is connected, $ V_{\pi', K} = \U(\g) v $ for any nonzero $ v \in V_{\pi', K} $. It follows that $ f \in (V_{\pi, -\infty} \otimes V_\phi)^\Gamma $ is a Schwartz vector (resp. cuspidal Schwartz vector) if and only if there is a nonzero $ v \in V_{\pi', K} $ such that $ c_{f, v} \in \CS(\Gamma \bs G, \phi) $ (resp. $ c_{f, v} \in \c{\CS}(\Gamma \bs G, \phi) $).
    \item The space $ (V_{\pi, -\infty} \otimes V_\phi)_\CS^\Gamma $ is clearly contained in $ (V_{\pi, -\infty} \otimes V_\phi)_d^\Gamma $.
    \end{enumerate}
\end{rem}

Assume that $ \pi $ is irreducible.
Fix $ 0 \neq v \in V_{\pi', K} $. We equip $ (V_{\pi, -\infty} \otimes V_\phi)_\CS^\Gamma $ with the coarsest topology such that the map
$
    f \in (V_{\pi, -\infty} \otimes V_\phi)_\CS^\Gamma \mapsto c_{f, v} \in \CS(\Gamma \bs G, \phi)
$
is continuous. Analogously, we equip $ (V_{\pi, -\infty} \otimes V_\phi)_d^\Gamma $
with the coarsest topology such that the map
$
    f \in (V_{\pi, -\infty} \otimes V_\phi)_d^\Gamma \mapsto c_{f, v} \in L^2(\Gamma \bs G, \phi)
$ is continuous.

As $ \U(\g) v $ is equal to $ V_{\pi', K} $ and as $ \CS(\Gamma \bs G, \phi) \to \CS(\Gamma \bs G, \phi) , \, f \mapsto R_Y f $ and $ L^2(\Gamma \bs G, \phi) \to L^2(\Gamma \bs G, \phi) , \, f \mapsto R_Y f $ are continuous for every $ Y \in \U(\g) $, it follows that the above two maps are continuous for every $ v \in V_{\pi', K} $. This shows that the topologies on $ (V_{\pi, -\infty} \otimes V_\phi)_\CS^\Gamma $ and on $ (V_{\pi, -\infty} \otimes V_\phi)_d^\Gamma  $ do not depend on the choice of $v$.

\begin{lem}\label{Lem space of Schwartz vectors can be continuously injected into space of square-integrable vectors}
    The inclusions
    \[
        (V_{\pi, -\infty} \otimes V_\phi)^\Gamma_\CS \subset (V_{\pi, -\infty} \otimes V_\phi)^\Gamma_d , \quad
        (V_{\pi, -\infty} \otimes V_\phi)^\Gamma \subset V_{\pi, -\infty} \otimes V_\phi
    \]
    are continuous injections.
\end{lem}

\begin{proof} \nlenum
    \begin{enumerate}
    \item Consider the following commutative diagram:
        \[
            \xymatrix@R4pc@C4pc@M2.5mm{
                (V_{\pi, -\infty} \otimes V_\phi)^\Gamma_\CS \ar@{^{(}->}[d]_i \ar@{^{(}->}[r]
                & \CS(\Gamma \bs G, \phi) \ar@{^{(}->}[d]_{j} \\
                (V_{\pi, -\infty} \otimes V_\phi)^\Gamma_d \ar@{^{(}->}[r] & L^2(\Gamma \bs G, \phi)
            }
        \]
        Let us view $ (V_{\pi, -\infty} \otimes V_\phi)^\Gamma_\CS $ (resp. $ (V_{\pi, -\infty} \otimes V_\phi)^\Gamma_d $) as a subspace of $ \CS(\Gamma \bs G, \phi) $ (resp. $ L^2(\Gamma \bs G, \phi) $). The first assertion follows as $ \CS(\Gamma \bs G, \phi) \hookrightarrow L^2(\Gamma \bs G, \phi) $ is continuous.
    \item Since we consider on $ (V_{\pi, -\infty} \otimes V_\phi)^\Gamma $ the subspace topology,
        $
            (V_{\pi, -\infty} \otimes V_\phi)^\Gamma \hookrightarrow V_{\pi, -\infty} \otimes V_\phi
        $
        is also continuous.
    \end{enumerate}
    The lemma follows.
\end{proof}

\newpage

\subsection{Some definitions and known results}\label{ssec:SomedefinitionsAndKnownResults}

We consider the right regular representation of $G$ on $ L^2(\Gamma \bs G, \phi) $.
By the abstract Plancherel theorem (see for example Theorem 14.10.5 of \cite[p.336]{Wallach2}), this unitary representation has a direct integral decomposition into irreducibles:\index{abstract Plancherel theorem}\index[n]{F@$ \F $, $ \F_\pi $}
\begin{equation}\label{eq abstract Plancherel decomposition}
    L^2(\Gamma \bs G, \phi) \overset{\F}{\simeq} \int_{\hat{G}}^{\oplus} N_\pi \hat{\otimes} V_\pi \, d\kappa(\pi) \ ,
\end{equation}
where $ \hat{G} $ denotes the unitary dual of $G$, $ V_\pi $ is an irreducible unitary representation space of a representation belonging to the class of $ \pi \in \hat{G} $, $ N_\pi $ is a Hilbert space which is called the \textit{multiplicity space} of $ \pi $ in $ L^2(\Gamma \bs G, \phi) $ and $ \kappa $ is a Borel measure on $ \hat{G} $ which is called the \textit{Plancherel measure}. The isomorphism $ \F $ is a unitary equivalence of representations, called ``Fourier transform''. The action of $G$ on the right-hand side is given by $ \Id_{N_\pi} \otimes \, \pi $.

\medskip

\begin{dfn}
    We say that a bounded linear operator $ \Phi \colon L^2(\Gamma \bs G, \phi) \to L^2(\Gamma \bs G, \phi) $ is \textit{decomposable} if there is a measurable family of bounded linear maps $ L_\pi \colon N_\pi \to N_\pi $ ($ \pi \in \hat{G} $) such that
    \[
        \Phi(f) = \F^{-1}\big( (L_\pi \otimes \Id)(\F_\pi(f))_{\pi \in \hat{G}} \big)
    \]
    for all $ f \in L^2(\Gamma \bs G, \phi) $.
\end{dfn}

\begin{prop}\label{Prop G-intertwining operators are decomposable}
    Any $G$-intertwining operator on $ L^2(\Gamma \bs G, \phi) $ is decomposable.
\end{prop}

\begin{proof}
    If $G$ is a so-called type I group (e.g. an abelian (resp. compact, semisimple) connected Lie group), then a unitary representation of $G$ admits by the abstract Plancherel theorem a direct integral decomposition over $ \hat{G} $, which is similar to the one of $ (\pi, L^2(\Gamma \bs G, \phi)) $.
    Theorem 3.24 of \cite[p.77]{Führ} shows moreover that the $G$-intertwining operators are decomposable.
\end{proof}

\medskip

Bunke-Olbrich \cite{BO00} give a much more precise description of the direct integral decomposition \eqref{eq abstract Plancherel decomposition}:

By Corollary 8.7 of \cite[p.132]{BO00}, $ N_\pi $ can be realised as a subspace of $ (V_{\pi', -\infty} \otimes V_\phi)_{temp}^\Gamma $, $ \pi \in \hat{G} $, such that we have the following:
\begin{enumerate}
\item The matrix coefficient map $ c \colon (V_{\pi', -\infty} \otimes V_\phi)^\Gamma \otimes V_{\pi, \infty} \to \CS'(\Gamma \bs G, \phi) $ induces a map\index[n]{cpi@$ c_\pi $}
    \[
        c_\pi \colon N_\pi \hat{\otimes} V_\pi \to \CS'(\Gamma \bs G, \phi) \ .
    \]
    We denote the adjoint of $ c_\pi $ by $ \F_\pi \colon \CS(\Gamma \bs G, \phi) \to N_\pi \hat{\otimes} V_\pi $.
\item The maps $ \F_\pi $ satisfy $ \F(f)(\pi) = \F_\pi(f) $ for all $ f \in \CS(\Gamma \bs G, \phi) $.
\item The space $ N_\pi $ is equal to $ (V_{\pi', -\infty} \otimes V_\phi)_d^\Gamma $ if $ (V_{\pi', -\infty} \otimes V_\phi)_d^\Gamma \neq \{0\} $.
\end{enumerate}
In particular, the Plancherel measure $ \kappa $ has support on
\[
    \{ \pi \in \hat{G} \mid (V_{\pi', -\infty} \otimes V_\phi)_{temp}^\Gamma \neq \{0\} \}
\]
and $ \kappa(\{\pi\}) \neq 0 $ if and only if $ (V_{\pi', -\infty} \otimes V_\phi)_d^\Gamma \neq \{0\} $.
If $ (V_{\pi', -\infty} \otimes V_\phi)_d^\Gamma \neq \{0\} $, then we choose the Plancherel measure and the scalar product on $ N_\pi = (V_{\pi', -\infty} \otimes V_\phi)_d^\Gamma $ such that
\begin{enumerate}
\item $ \kappa(\{\pi\}) = 1 $,
\item the matrix coefficient map $ c_\pi $ induces an isometric embedding of $ N_\pi \hat{\otimes} V_\pi $ into $ L^2(\Gamma \bs G, \phi) $.
\end{enumerate}
In particular, $ \F_\pi $ extends then to a map from $ L^2(\Gamma \bs G, \phi) $ to $ (V_{\pi', -\infty} \otimes V_\phi)_d^\Gamma \hat{\otimes} V_\pi $.

The unitary dual $ \hat{G} $ is the disjoint union of the discrete series representations, the unitary principal series representations and the complementary series representations:
\[
    \hat{G} = \hat{G}_d \cup \hat{G}_u \cup \hat{G}_c \ ,
\]
\needspace{4\baselineskip}
where\index[n]{Gd@$ \hat{G}_d $, $ \hat{G}_u $, $ \hat{G}_c $}
\begin{align*}
    & \hat{G}_d := \{ \pi \in \hat{G} \mid V_\pi \text{ is square-integrable} \} \ ; \\
    & \hat{G}_u := \{ \pi \in \hat{G} \mid V_\pi \text{ is tempered} \} \smallsetminus \hat{G}_d \ ; \\
    & \hat{G}_c := \{ \pi \in \hat{G} \mid V_\pi \text{ is not tempered} \} \ .
\end{align*}
By Theorem 11.1 of \cite[p.155]{BO00}, $ L^2(\Gamma \bs G, \phi) $ decomposes into a continuous and a discrete part:\index[n]{LGammaG2ac@$ L^2(\Gamma \bs G, \phi)_{ac/ds/res/U} $}
\begin{equation}\label{eq Plancherel decomposition}
    L^2(\Gamma \bs G, \phi) = L^2(\Gamma \bs G, \phi)_{ac} \oplus L^2(\Gamma \bs G, \phi)_d \ .
\end{equation}
Here, the space $ L^2(\Gamma \bs G, \phi)_{ac} $ is generated by wave packets of matrix coefficients and $ L^2(\Gamma \bs G)_d $ can be further decomposed with the help of representation theory:
\[
    L^2(\Gamma \bs G, \phi)_d = L^2(\Gamma \bs G, \phi)_{ds} \oplus L^2(\Gamma \bs G, \phi)_{res} \oplus L^2(\Gamma \bs G, \phi)_U \ ,
\]
where $ L^2(\Gamma \bs G, \phi)_{res} $ is generated by residues of Eisenstein series, $ L^2(\Gamma \bs G, \phi)_U $ is generated by matrix coefficients of ``stable'' invariant distributions supported on the limit set and $ L^2(\Gamma \bs G, \phi)_{ds} $ is generated by square-integrable matrix coefficients of discrete series representations (it is denoted by $ L^2(\Gamma \bs G, \phi)_{cusp} $ in \cite{BO00}).

The space $ L^2(\Gamma \bs G, \phi)_{ac} $ corresponds essentially to $ \hat{G}_u $ and $ L^2(\Gamma \bs G, \phi)_d $ corresponds essentially to $ \hat{G}_d \cup \hat{G}_c $: $ L^2(\Gamma \bs G, \phi)_{ds} $ corresponds to $ \hat{G}_d $ and $ L^2(\Gamma \bs G, \phi)_{res} \oplus L^2(\Gamma \bs G, \phi)_U $ corresponds essentially to $ \hat{G}_c $.

In the following sections, we determine the contribution of the space of cusp forms to each of the appearing spaces.

Let us denote the orthogonal projection on $ L^2(\Gamma \bs G, \phi)_{ac} $ (resp. $ L^2(\Gamma \bs G, \phi)_{res} $, $ L^2(\Gamma \bs G, \phi)_U $, $ L^2(\Gamma \bs G, \phi)_{ds} $ by $ p_{ac} $ (resp. $ p_{res} $, $ p_U $, $ p_{ds} $).\index[n]{pac, pres, pU, pds@$ p_{ac} $ $ p_{res} $, $ p_U $, $ p_{ds} $}

\begin{prop}\label{Prop a G-int op commutes with p_ds, p_ac and p_res + p_U}
    Let $ I $ be a $G$-intertwining operator on $ L^2(\Gamma \bs G, \phi) $. Then, $ I $ commutes with $ p_{ac} $, $ p_{ds} $ and $ p_{res} + p_U $.
\end{prop}

\begin{proof}
    As $I$ is decomposable by Proposition \ref{Prop G-intertwining operators are decomposable},
    \[
        p(f) = \F^{-1}\big( (L_\pi \otimes \Id)(\F_\pi(f))_{\pi \in \hat{G}} \big) \ .
    \]
    Let $ p = p_{ac} $ (resp. $ p = p_{ds} $, $ p = p_{res} + p_U $). By the proof of Theorem 11.1 of \cite[p.155]{BO00},
    \[
        p(f) = \F^{-1}\big( (L'_\pi \otimes \Id)(\F_\pi(f))_{\pi \in \hat{G}} \big)
    \]
    with $ L'_\pi = 0 $ or $ L'_\pi = \Id $ for almost all $ \pi \in \hat{G} $. As $ L_\pi \circ L'_\pi = L'_\pi \circ L_\pi $ for almost all $ \pi \in \hat{G} $, the proposition follows.
\end{proof}

Let us continue with some elementary consequences of the Plancherel decomposition.

\begin{lem}\label{Lem (c_pi o F_pi)^2 = c_pi o F_pi, F_pi o c_pi = Id on V_(pi', -infty)^Gamma otimes V_pi}
    Let $ \pi \in \hat{G} $ be such that $ (V_{\pi', -\infty} \otimes V_\phi)_d^\Gamma \neq \{0\} $. Then,
    $ \F_{\pi} \circ c_{\pi} $ is the identity on $ (V_{\pi', -\infty} \otimes V_\phi)^\Gamma_d \hat{\otimes} V_{\pi} $ and
    \[
        (c_{\pi} \circ \F_{\pi})^2 = c_{\pi} \circ \F_{\pi} \ .
    \]
\end{lem}

\begin{proof} Let $ \pi $ be as above.
    \begin{enumerate}
    \item For all $ f, g \in \F_\pi(L^2(\Gamma \bs G, \phi)) = (V_{\pi', -\infty} \otimes V_\phi)^\Gamma_d \hat{\otimes} V_{\pi} $, we have
        \[
            (\F_{\pi}(c_{\pi}(f)), g) = (c_{\pi}(f), c_{\pi}(g)) = (f, g) \ .
        \]
        Since $ f $ and $ \F_{\pi}(c_{\pi}(f)) $ belong to the Hilbert space $ (V_{\pi', -\infty} \otimes V_\phi)^\Gamma_d \hat{\otimes} V_{\pi} $, it follows that $ f = \F_{\pi}(c_{\pi}(f)) $.
    \item Let $ f, g \in L^2(\Gamma \bs G, \phi) $. As $ (V_{\pi', -\infty} \otimes V_\phi)_d^\Gamma \neq \{0\} $, $ (c_{\pi} \circ \F_{\pi})(f) \in L^2(\Gamma \bs G, \phi) $. Thus, $ (c_{\pi} \circ \F_{\pi})^2(f) $ is well-defined. Moreover, we have
        \[
            ((c_{\pi} \circ \F_{\pi})^2(f), g) = (\F_{\pi}((c_{\pi} \circ \F_{\pi})(f)), \F_{\pi}(g))_{N_{\pi} \otimes V_{\pi}} = ((c_{\pi} \circ \F_{\pi})(f), g) \ .
        \]
        As $ (\cdot, \cdot) $ is nondegenerate, $ (c_{\pi} \circ \F_{\pi})^2 = c_{\pi} \circ \F_{\pi} $.
    \end{enumerate}
\end{proof}

\begin{cor}\label{Cor c_pi o F_pi is the orth. proj. on V_(pi', -infty)^Gamma otimes V_pi}
    Let $ \pi \in \hat{G} $ be such that $ (V_{\pi', -\infty} \otimes V_\phi)_d^\Gamma \neq \{0\} $. Then, $ c_{\pi} \circ \F_{\pi} $ is the orthogonal projection on $ c_{\pi}\big((V_{\pi', -\infty} \otimes V_\phi)_d^\Gamma \hat{\otimes} V_{\pi}\big) $. Thus,
    \[
        p_d \colon f \in L^2(\Gamma \bs G, \phi) \mapsto \sum_{\pi \in \hat{G} \, : \, (V_{\pi', -\infty} \otimes V_\phi)_d^\Gamma \neq \{0\}} c_{\pi}(\F_{\pi}(f)) \in L^2(\Gamma \bs G, \phi)_d
    \]
    is the orthogonal projection on $ L^2(\Gamma \bs G, \phi)_d $.
\end{cor}

\begin{proof}
    The corollary follows from the previous lemma and the equality $ (c_{\pi} \circ \F_{\pi})^* = c_{\pi} \circ \F_{\pi} $.
\end{proof}

\begin{cor}\label{Cor c_pi(N_pi hat(otimes) V_pi), c_pi(N_pi otimes F) closed}
    Let $ \pi \in \hat{G} $ be such that $ (V_{\pi', -\infty} \otimes V_\phi)_d^\Gamma \neq \{0\} $, let $ \gamma \in \hat{K} $ and let $ F $ be a subspace of $ V_{\pi}(\gamma) $. Then,
    \begin{equation}\label{Cor c_pi(N_pi hat(otimes) V_pi), c_pi(N_pi otimes Cv) closed eq1}
        c_\pi\big((V_{\pi', -\infty} \otimes V_\phi)_d^\Gamma \hat{\otimes} V_{\pi}\big)
    \end{equation}
    and
    \begin{equation}\label{Cor c_pi(N_pi hat(otimes) V_pi), c_pi(N_pi otimes Cv) closed eq2}
        S_{\pi, F} := c_\pi((V_{\pi', -\infty} \otimes V_\phi)^\Gamma_d \otimes F)
    \end{equation}
    are closed in $ L^2(\Gamma \bs G, \phi) $.
\end{cor}

\begin{proof}
    Let $ p = \Id - c_\pi \circ \F_\pi $. Then, $p$ is a continuous orthogonal projection on the orthogonal complement of \eqref{Cor c_pi(N_pi hat(otimes) V_pi), c_pi(N_pi otimes Cv) closed eq1}. Thus, \eqref{Cor c_pi(N_pi hat(otimes) V_pi), c_pi(N_pi otimes Cv) closed eq1} is equal to $ \ker p $ and hence closed.
    So, $ S_{\pi, V_{\pi}(\gamma)} = c_\pi\big((V_{\pi', -\infty} \otimes V_\phi)_d^\Gamma \otimes V_{\pi}(\gamma)\big) = c_\pi\big((V_{\pi', -\infty} \otimes V_\phi)_d^\Gamma \hat{\otimes} V_{\pi}\big)(\gamma) $ is closed, too.

    Let $ \gamma \in \hat{K} $ and let $ F $ be a subspace of $ V_{\pi}(\gamma) $ (finite-dimensional). Let $ F^\perp $ denote the orthogonal complement of $ F $ in $V_{\pi}(\gamma)$. Then, $ S_{\pi, F} $ is the orthogonal complement of $ S_{\pi, F^\perp} $ in $ S_{\pi, V_{\pi}(\gamma)} $ and is hence closed.
\end{proof}

Fix $ 0 \neq v \in V_{\pi, K} $. It follows from Corollary \ref{Cor c_pi(N_pi hat(otimes) V_pi), c_pi(N_pi otimes F) closed} and the open mapping theorem that the topology on $ N_\pi $ is the coarsest such that
\[
    f \in N_\pi \mapsto c_{\pi}(f \otimes v)
\]
is continuous. Recall that this topology does not depend on the choice of $v$.

\medskip

Let us define now some representations and some bundles playing an important role in the following.

Let $ (\sigma, V_\sigma) $ be a finite-dimensional unitary representation of $M$. Then, we define its Weyl-conjugate $ \sigma^w $\index[n]{sigmaw@$ \sigma^w $} by $ \sigma^w(m) := \sigma(w^{-1} m w) $, where $ w \in N_K(\a) $ is a representative of the nontrivial Weyl group element.

We say that $ \sigma $ is \textit{Weyl-invariant}\index{Weyl-invariant} if $ \sigma $ is equivalent to $ \sigma^w $. We assume that $ \sigma $ denotes a Weyl-invariant representation of $M$ from now on.
It is either irreducible or of the form $ \sigma' \oplus \sigma'^w $ with $ \sigma' $ irreducible and not Weyl-invariant and it can be extended in both cases to a representation of $ N_K(\a) $. We denote this extension also by $ \sigma $. As $ \sigma(w) $ is only unique up to a sign, we fix such an extension.
We denote the dual representation of $ \sigma $ by $ \tilde{\sigma} $.

Let $ (\gamma, V_\gamma) $ be a finite-dimensional representation of $K$.
We form the homogeneous vector bundle $ V(\gamma) := G \times_K V_\gamma $ over $X$ and the locally homogeneous vector bundle $ V_Y(\gamma) := \Gamma \bs G \times_K V_\gamma $ over $ Y := \Gamma \bs X $.
Set $ V(\gamma, \phi) = V(\gamma) \otimes V_\phi $ and $ V_Y(\gamma, \phi) = V_Y(\gamma) \otimes V_\phi $.\index[n]{Vgamma VYgamma VBsigmalambda@$ V(\gamma) $, $ V_Y(\gamma) $, $ V_B(\sigma_\lambda) $}

Let $ W $ be a vector space which is equipped with a left $K$-action. Let\index[n]{WK@$ W^K $}
\[
    W^K = \{ v \in W \mid k.v = v \quad \forall k \in K \} \ .
\]
Consider the tensor product action of $K$ on $ L^2(\Gamma \bs G, \phi) \otimes V_\gamma $. Using the isomorphism $ L^2(Y, V_Y(\gamma, \phi)) \simeq [L^2(\Gamma \bs G, \phi) \otimes V_\gamma]^K $, we define\index[n]{CY@$ \CS(Y, V_Y(\gamma, \phi)) $}
\[
    \CS(Y, V_Y(\gamma, \phi)) := [\CS(\Gamma \bs G, \phi) \otimes V_\gamma]^K \subset L^2(Y, V_Y(\gamma, \phi)) \ .
\]
For $ \lambda \in \a_\CC^* $, set $ \sigma_\lambda(man) = \sigma(m) a^{\rho-\lambda} $. This defines a representation of $P$ on $ V_{\sigma_\lambda} := V_\sigma $. We denote by $ V(\sigma_\lambda) := G \times_P V_{\sigma_\lambda} $\index[n]{Vsigmalambda@$ V(\sigma_\lambda) $} the associated homogeneous bundle over $ \dX = G/P $. Let $ B = \Gamma \bs \Omega $. Then, $ V_B(\sigma_\lambda) := \restricted{\Gamma \bs V(\sigma_\lambda)}{\Omega} $ defines a bundle over $B$.

Recall that $ (\phi, V_\phi) $ is a finite-dimensional unitary representation of $ \Gamma $. We denote its dual by $ (\tilde{\phi}, V_{\tilde{\phi}}) $.

On the space of sections of $ V(\sigma_\lambda) $, we consider the left regular representation of $G$. We denote it by $ \pi^{\sigma, \lambda} $ and call it a \textit{principal series representation} of $G$. We define the bundle $ V(\sigma_\lambda, \phi) := V(\sigma_\lambda) \otimes V_\phi $ on $ \dX $ carrying the tensor product action of $ \Gamma $ and we define $ V_B(\sigma_\lambda, \phi) := \restricted{\Gamma \bs (V(\sigma_\lambda) \otimes V_\phi)}{\Omega} $.

The space of distribution sections $ \C^{-\infty}(\dX, V(\sigma_\lambda)) $\index[n]{C-inftydXVsigmalambda@$ \C^{-\infty}(\dX, V(\sigma_\lambda)) $} (resp. $ \C^{-\infty}(\dX, V(\sigma_\lambda, \phi)) $) is defined by $ \Cinf(\dX, V(\tilde{\sigma}_{-\lambda}))' $ (resp. $ \Cinf(\dX, V(\tilde{\sigma}_{-\lambda}, \phi))' $). Here, we take the strong dual.
Analogously, we define $ \C^{-\infty}(B, V_B(\sigma_\lambda)) $\index[n]{C-inftyBVBsigmalambda@$ \C^{-\infty}(B, V_B(\sigma_\lambda)) $} (resp. $ \C^{-\infty}(B, V_B(\sigma_\lambda, \phi)) $) by $ \Cinf(B, V_B(\tilde{\sigma}_{-\lambda}))' $ (resp. $ \C^{-\infty}(B, V_B(\tilde{\sigma}_{-\lambda}, \phi))' $). Here, we take also the strong dual.
We denote the space of $ \Gamma $-invariant distribution sections by $ \C^{-\infty}(\dX, V(\sigma_\lambda))^\Gamma $ (resp. $ \C^{-\infty}(\dX, V(\sigma_\lambda, \phi))^\Gamma $).\index[n]{C-inftydXVsigmalambdaGamma@$ \C^{-\infty}(\dX, V(\sigma_\lambda))^\Gamma $}

We denote the subspace of $K$-finite functions of $ \Cinf(\dX, V(\sigma_\lambda)) $ by $ \C^K(\dX, V(\sigma_\lambda)) $.\index[n]{CKdXVsigmalambda@$ \C^K(\dX, V(\sigma_\lambda)) $}

Note that we can identify distribution sections of $ V_B(\sigma_\lambda, \phi) $ with $ \Gamma $-invariant sections of $ \restricted{ V(\sigma_\lambda, \phi) }{\Omega} $. By a slight abuse of notation, we denote the subspace of $ \C^{-\infty}(\dX, V(\sigma_\lambda, \phi)) $ consisting of $\Gamma$-invariant distributions supported on the limit set by $ \C^{-\infty}(\Lambda, V(\sigma_\lambda, \phi))^\Gamma $.
\index[n]{C-inftyLambdaVsigmalambdaGamma@$ \C^{-\infty}(\Lambda, V(\sigma_\lambda, \phi))^\Gamma $}

Let $ W $ be a vector space and let $ A \in \End(W) $. If we write $ A(f \otimes v) $ ($ f \in W $, $ v \in V_\phi $), then this means that we consider simply the trivial action on $ V_\phi $:
\[
    A(f \otimes v) := A(f) \otimes v   \qquad (f \in W, v \in V_\phi) \ .
\]
In contrast to bilinear parings $ \langle \cdot, \cdot \rangle $, we will always denote sesquilinear pairings by $ (\cdot, \cdot) $. By convention, these pairings are complex linear with respect to the first variable.\index[n]{( . , . ), langle . , . rangle@$ (\cdot, \cdot) $, $ \langle \cdot, \cdot \rangle $}

\begin{dfn}[{10.1.2 of \cite[p.5]{Wallach2}, \cite[p.96]{BO00}}] \nlenum \\
    Let us consider $ f \in \Cinf(\dX, V(\sigma_\lambda)) $ as a function on $G$ taking values in $ V_{\sigma_\lambda} $ and satisfying\index{Knapp-Stein intertwining operator}
    \[
        f(gp) = \sigma_\lambda(p)^{-1} f(g)
    \]
    for all $ g \in G $ and $ p \in P $. If $ \R(\lambda) < 0 $, then the \textit{Knapp-Stein intertwining operator}
    \[
        \hat{J}^w_{\sigma, \lambda} \colon \Cinf(\dX, V(\sigma_\lambda)) \to \Cinf(\dX, V(\sigma^w_{-\lambda}))
    \]
    is defined by the convergent integral
    \[
        (\hat{J}^w_{\sigma, \lambda} f)(g) := \int_N{ f(g n w) \, dn } \ .
    \]
    Its continuous extension to distributions is obtained by duality.

    Set $ \hat{J}_{\sigma, \lambda} = \sigma(w) \hat{J}^w_{\sigma, \lambda} $. Then, this operator from $ \Cinf(\dX, V(\sigma_\lambda)) $ to $ \Cinf(\dX, V(\sigma_{-\lambda})) $ does not depend on the choice of $w$.\index[n]{Jf@$ \hat{J}^w_{\sigma, \lambda} $, $ \hat{J}_{\sigma, \lambda} $}
\end{dfn}

\begin{rem} \nlenum
    \begin{enumerate}
    \item The convergence of the integral follows easily from Proposition \ref{Prop N bar integral is finite} applied to $ E = N $.
    \item We denote $ \hat{J}^w_{\sigma, \lambda} \otimes \Id \colon \C^*(\dX, V(\sigma_\lambda, \phi)) \to \C^*(\dX, V(\sigma^w_{-\lambda}, \phi)) $ also by $ \hat{J}^w_{\sigma, \lambda} $.
    \item If $ \R(\lambda) \geq 0 $, then it is defined by meromorphic continuation (see \cite{KS71} or Lemma 5.2 of \cite[p.97]{BO00}).
    \end{enumerate}
\end{rem}

\needspace{2\baselineskip}
The following is a well-known result, which we use in several proofs:
\begin{lem}\label{Lem int_K f(k) dk = int_(bar N) f(kappa(bar n) a(bar n)^(-2rho) dbar(n)}
    Let $ f $ be a continuous right $M$-invariant function on $K$. Then,
    \[
        \int_K f(k) \, dk = \int_{\bar{N}} f(\kappa(\bar{n})) a(\bar{n})^{-2\rho} \, d\bar{n}
        = \int_N f(\kappa(n w)) a(n w)^{-2\rho} \, dn \ .
    \]
\end{lem}

\begin{rem}
    The map $ \bar{n} \mapsto \kappa(\bar{n})M $ is a diffeomorphism of $ \bar{N} $ onto $ K/M \smallsetminus \{wM\} $.
\end{rem}

\begin{proof}
    Let $ f $ be a continuous right $M$-invariant function on $K$. Then,
    \[
        \int_K f(k) \, dk = \int_{\bar{N}} f(\kappa(\bar{n})) a(\bar{n})^{-2\rho} \, d\bar{n}
    \]
    by Proposition 8.46 of \cite[p.542]{Knapp2}. The normalisation is correct as
    \[
        \int_K 1 \, dk = 1 = \int_{\bar{N}} a(\bar{n})^{-2\rho} \, d\bar{n} \ .
    \]
    Thus,
    \[
        \int_K f(k) \, dk = \int_K f(w k) \, dk = \int_{\bar{N}} f(w \kappa(\bar{n})) a(\bar{n})^{-2\rho} \, d\bar{n} \ .
    \]
    Since $ n \in N \mapsto w^{-1} n w \in \bar{N} $ is a diffeomorphism transforming $ dn $ into $ d\bar{n} $, this is again equal to
    \[
        \int_N f(\kappa(n w)) a(n w)^{-2\rho} \, dn \ .
    \]
\end{proof}

Let $ \dX = U \cup Q $, where $U$ is open and $Q$ is the complement of $U$. Let $ \lambda \in \a^*_\CC $ and $ f \in \C^{-\infty}(\dX, V(\sigma_\lambda, \phi)) $ with support in $Q$. Fix $ gP \in U $.
Let $ 1_N \colon N \to \RR, n \mapsto 1 $. As $ f $ has support in $Q$ and as the open set $ G \smallsetminus \{gP\} = g N w P $ is diffeomorphic to $N$, we can identify $ f $ with a distribution $ f_N $ in $ \C^{-\infty}(N, V_\sigma \otimes V_\phi) $ such that $ \int_N f(g n w) \, dn := \langle f_N, 1_N \rangle_N $ is equal to $ \hat{J}^w_{\sigma, \lambda} f(g) $.

We show next that matrix coefficients have under certain assumptions an asymptotic expansion in the sense of N.~Wallach.
The following two lemmas are proven in \cite[p.110]{BO00} in the setting of the Poisson transform. For the convenience of the reader, we redo here the proof for matrix coefficients.

\begin{lem}\label{Lemma about Asymptotic expansions}
    Let $ \dX = U \cup Q $, where $U$ is open and $Q$ is the complement of $U$. Let $ \lambda \in \a^*_\CC $.
    Let $ f \in \C^{-\infty}(\dX, V(\sigma_\lambda, \phi)) $ with support in $Q$, let $ v \in \C^K(\dX, V(\tilde{\sigma}_{-\lambda})) $ be a $K$-finite function, let $ a \in A_+ $ and let $ k \in K(U), h \in K $. 
    Then, there exist smooth functions $ \psi_{h, j} $ ($ j \in \NN $) on $ K(U) $ 
    such that
    \begin{equation}\label{eqAsymptotic expansion}
        c_{f, v}(kah) = a^{-(\lambda + \rho)} \langle (\hat{J}^w_{\sigma, \lambda} f)(k), v(h^{-1}w) \rangle
            + \sum_{j=1}^\infty{ a^{-(\lambda + \rho + j\alpha)} \psi_{h, j}(k) } \ .
    \end{equation}
    For $ a \gg 0 $ and $ k $ in a compact subset of $K(U)$, the series converges uniformly. In particular, for $ a \to \infty $, we have
    \begin{equation}\label{eqAsymptotic expansion2}
        c_{f, v}(kah) = a^{-(\lambda + \rho)} \langle (\hat{J}^w_{\sigma, \lambda} f)(k), v(h^{-1}w) \rangle
            + O(a^{-\R(\lambda)-\rho-\alpha})
    \end{equation}
    uniformly for $ k $ varying in a compact subset of $K(U)$.
\end{lem}

\begin{proof}
    Let $ U $, $ Q $, $ v $ and $ f $ as above.
    For $ k \in K(U) $, $ h \in K $ and $ a \in A_+ $, we have 
    \[
        c_{f, v}(kah) = \int_K{ \langle f(l), \pi^{\tilde{\sigma}, -\lambda}(kah) v(l) \rangle \, dl }
        = \int_K{ \langle f(l), v(h^{-1} a^{-1} k^{-1} l) \rangle \, dl } \ ,
    \]
    where the integral is a formal notation meaning that the distribution $f$ has to be applied to the smooth integral kernel. Since $K$ is unimodular, this yields
    \[
        \int_K{ \langle f(k w l), v(h^{-1} a^{-1} w l) \rangle \, dl } \ .
    \]
    By Lemma \ref{Lem int_K f(k) dk = int_(bar N) f(kappa(bar n) a(bar n)^(-2rho) dbar(n)}, this is equal to
    \[
        \int_{\bar{N}}{ \langle f(k w \kappa(\bar{n})), v(h^{-1} a^{-1} w \kappa(\bar{n})) \rangle a(\bar{n})^{-2\rho} \, d\bar{n} } \ .
    \]
    As $ \bar{n} = \kappa(\bar{n}) a(\bar{n}) n(\bar{n}) $, $ \kappa(\bar{n}) = \bar{n} n(\bar{n})^{-1} a(\bar{n})^{-1} $.
    So,
    \begin{multline*}
        \kappa(h^{-1} a^{-1} w \kappa(\bar{n})) = h^{-1} \kappa(a^{-1} w \bar{n} n(\bar{n})^{-1} a(\bar{n})^{-1}) \\
        = h^{-1} w \kappa(\underbrace{(a \bar{n} a^{-1})}_{\in \bar{N}} \underbrace{(a a(\bar{n})^{-1})}_{\in A} \underbrace{(a(\bar{n}) n(\bar{n})^{-1} a(\bar{n})^{-1})}_{\in N})
        = h^{-1} w \kappa(a \bar{n} a^{-1})
    \end{multline*}
    and
    \begin{multline*}
        a(h^{-1} a^{-1} w \kappa(\bar{n})) = a(a \kappa(\bar{n})) = a(a \bar{n} n(\bar{n})^{-1} a(\bar{n})^{-1}) \\
        = a((a \bar{n} a^{-1}) a n(\bar{n})^{-1}) a(\bar{n})^{-1}
        = a(a \bar{n} a^{-1}) a a(\bar{n})^{-1} \ .
    \end{multline*}
    Thus, $ c_{f, v}(kah) $ is equal to
    \begin{equation}\label{Proof U_Lambda eq1}
        a^{-(\lambda + \rho)} \int_{\bar{N}}{ a(\bar{n})^{\lambda -\rho} a(a \bar{n} a^{-1})^{-(\lambda + \rho)} \langle f(k w \kappa(\bar{n})), v(h^{-1} w \kappa(a \bar{n} a^{-1})) \rangle \, d\bar{n} } \ .
    \end{equation}
    For $ t > 0 $, define $ a_t \in A $ by $ t = a_t^{-\alpha} $. Define
    \[
        \Phi \colon (0, \infty) \times \bar{N} \to \bar{N}, \quad (t, \bar{n}) \mapsto a_t \bar{n} a_t^{-1} \ .
    \]
    As $ \bar{\n} = \g_{-\alpha} \oplus \g_{-2\alpha} $, this can be written as
    \[
        \Phi(t, \exp(X+Y)) = \exp(tX + t^2Y) \qquad (X \in \g_{-\alpha}, \; Y \in \g_{-2\alpha}) \ .
    \]
    Thus, we can extend $ \Phi $ and hence
    \[
        (t, \bar{n}) \mapsto a(a_t \bar{n} a_t^{-1})^{-(\lambda + \rho)} v(h^{-1}w \kappa(a_t \bar{n} a_t^{-1})) = v(h^{-1}w a_t \bar{n} a_t^{-1})
    \]
    to a function on $ \RR \times \bar{N} $ that has a Taylor expansion relative to $t$ at $ t=0 $ converging in the space of smooth functions on $ \bar{N} $ with values in $ V(\tilde{\sigma}_{-\lambda}) $.
    As $ \lim_{t \to 0}{a_t \bar{n} a_t^{-1}} = e \in G $, it is given by
    \[
        a(a_t \bar{n} a_t^{-1})^{-(\lambda + \rho)} v(h^{-1} w \kappa(a_t \bar{n} a_t^{-1})) = v(h^{-1}w) + \sum_{j=1}^\infty{ A_{h, j}(\bar{n}) t^j } \ ,
    \]
    where $ A_{h, j} \colon \bar{N} \to V(\tilde{\sigma}_{-\lambda}) $ are analytic. Inserting this expansion into \eqref{Proof U_Lambda eq1} yields
    \begin{multline*}
        c_{f, v}(kah)
        = a^{-(\lambda + \rho)} \int_{\bar{N}} \langle { \overbrace{a(\bar{n})^{\lambda -\rho} f(k w \kappa(\bar{n}))}^{ = f(k w \bar{n}) }, v(h^{-1}w) \rangle \, d\bar{n} } \\
        + \sum_{j=1}^\infty{ a^{-(\lambda + \rho + j\alpha)} \underbrace{ \int_{\bar{N}}{ a(\bar{n})^{\lambda -\rho} \langle f(k w \kappa(\bar{n})), A_{h, j}(\bar{n}) \rangle \, d\bar{n} } } }_{ =: \psi_{h, j}(k) } \ .
    \end{multline*}
    Since $K(U)$ is open and since $ k \in K(U) $, there exists a neighbourhood $ N_e $ of $e$ in $K$ such that $ k N_e $ is still contained in $ K(U) $.
    As
    $
        k w \kappa(\bar{n}) \in K(Q)
    $
    if $ f(k w \kappa(\bar{n})) \neq 0 $, $ w \kappa(\bar{n}) \in k^{-1} K(Q) $ for all $ \bar{n} \in \bar{N} $ such that $ f(k w \kappa(\bar{n})) \neq 0 $.
    Note that the intersection of $ N_e $ with $ k^{-1} K(Q) $ is empty since $ U \cap Q = \emptyset $ and that the set $ \{ \bar{n} \in \bar{N} \mid w \kappa(\bar{n}) \not \in N_e \} $ is relatively compact.
    Thus, $ k \mapsto f(k w \kappa(\cdot)) $ is a smooth family of distributions with compact support in $ \bar{N} $. Hence, the $ \psi_{h, j} $'s are smooth.
    Since $ \int_{\bar{N}}{ \langle f(k w \bar{n}), v(h^{-1}w) \rangle \, d\bar{n} } $ is equal to
    \[
        \int_N{ \langle f(knw), v(h^{-1}w) \rangle \, dn }
        = \langle (\hat{J}^w_{\sigma, \lambda} f)(k), v(h^{-1}w) \rangle \ ,
    \]
    we have finally
    \[
        c_{f, v}(kah) = a^{-(\lambda + \rho)} \langle (\hat{J}^w_{\sigma, \lambda} f)(k), v(h^{-1}w) \rangle
        + \sum_{j=1}^\infty{ a^{-(\lambda + \rho + j\alpha)} \psi_{h, j}(k) } \ .
    \]
\end{proof}

Let $ U $ be an open subset of $ \dX $. Let $ V $ be a relatively compact subset of $U$.
For any $ l \in \NN_0 $, any bounded subset $ A $ of $ \U(\g)_l $ and any $ f \in \Cinf(U, V(\sigma_\lambda, \phi)) $ (viewed as a function on $G$ with values in $ V_\sigma \otimes V_\phi $), set
\[
    \rho_{l, A}(f) = \sup_{X \in A, kM \in V} |R_X f(\kappa(k \, .))| \ .
\]
Then, $ \rho_{l, A} $ is a seminorm on $ \Cinf(U, V(\sigma_\lambda, \phi)) $. As it suffices to take countably many seminorms, the seminorms induce a Fréchet space structure on $ \Cinf(U, V(\sigma_\lambda, \phi)) $.

Define\index[n]{C-inftydXVsigmalambdaU@$ \C^{-\infty}_U(\dX, V(\sigma_\lambda, \phi)) $}
\[
    \C^{-\infty}_U(\dX, V(\sigma_\lambda, \phi)) = \{ f \in \C^{-\infty}(\dX, V(\sigma_\lambda, \phi)) \mid \restricted{f}{U} \in \C^\infty(U, V(\sigma_\lambda, \phi)) \} \ .
\]
We put on $ \C^{-\infty}_U(\dX, V(\sigma_\lambda, \phi)) $ the coarsest topology such that the inclusion map $ \C^{-\infty}_U(\dX, V(\sigma_\lambda, \phi)) \to \C^{-\infty}(\dX, V(\sigma_\lambda, \phi)) $ and the restriction map
\[
    \rest_U \colon \C^{-\infty}_U(\dX, V(\sigma_\lambda, \phi)) \to \Cinf(U, V(\sigma_\lambda, \phi))
\]
are continuous.

\begin{lem}\label{Lemma2 about Asymptotic expansions} \nlenum
    \begin{enumerate}
    \item Let $ f \in \Cinf(\dX, V(\sigma_\lambda, \phi)) $, $ v \in \C^K(\dX, V(\tilde{\sigma}_{-\lambda})) $ and $ k, h \in K $.

        If $ \R(\lambda) > 0 $, then there is $ \eps > 0 $ (depending on $ \lambda $ but not on $f$) such that, for $ a \to \infty $, $ c_{f, v}(kah) $ is equal to
        \[
            a^{\lambda - \rho} \langle f(k), \hat{J}^w_{\tilde{\sigma}, -\lambda} v(h^{-1}w^{-1}) \rangle + O(a^{\R(\lambda) - \rho -\eps})
        \]
        uniformly in $ k \in K $.

        If $ |\R(\lambda)| < \tfrac12 $ and $ \lambda \neq 0 $, then there is $ \eps > 0 $ (depending on $ \lambda $ but not on $f$) such that, for $ a \to \infty $, $ c_{f, v}(kah) $ is equal to
        \[
            a^{\lambda - \rho} \langle f(k), \hat{J}^w_{\tilde{\sigma}, -\lambda} v(h^{-1}w^{-1}) \rangle
            + a^{-\lambda - \rho} \langle (\hat{J}^w_{\sigma, \lambda} f)(k), v(h^{-1}w) \rangle
            + O(a^{-\tfrac{\alpha}2 - \rho -\eps})
        \]
        uniformly in $ k \in K $. The remainder depends jointly on $ \lambda $ and $f$.\label{Lemma2 about Asymptotic expansions eq1}
    \item Let $ f \in \Cinf(\dX, V(\sigma_{-\lambda}, \phi)) $ with $ \R(\lambda) > 0 $, $ v \in \C^K(\dX, V(\tilde{\sigma}_{\lambda})) $ and $ k, h \in K $. Then, there is $ \eps > 0 $ (depending on $ \lambda $ but not on $f$) such that, for $ a \to \infty $, $ c_{f, v}(kah) $ is equal to
        \[
            a^{\lambda - \rho} \langle (\hat{J}^w_{-\lambda} f)(k), v(h^{-1}w) \rangle + O(a^{\R(\lambda) - \rho -\eps})
        \]
        uniformly in $ k \in K $.
    \item Let $ \dX = U \cup Q $, where $U$ is open and $Q$ is the complement of $U$, and let $ f \in \C^{-\infty}_U(\dX, V(\sigma_\lambda, \phi)) $ and let $ v \in \C^K(\dX, V(\tilde{\sigma}_{-\lambda})) $.

        If $ \R(\lambda) > 0 $, then there exists $ \eps > 0 $ (depending on $ \lambda $ but not on $f$) such that, for $ a \to \infty $, $ c_{f, v}(kah) $ is equal to
        \[
            a^{\lambda - \rho} \langle f(k), \hat{J}^w_{\tilde{\sigma}, -\lambda} v(h^{-1}w^{-1}) \rangle + O(a^{\R(\lambda) - \rho -\eps})
        \]
        uniformly in $ kM \in \dX $ varying in a compact subset of $U$.\label{Lemma2 about Asymptotic expansions eq2}

        If $ |\R(\lambda)| < \tfrac12 $ and $ \lambda \neq 0 $, then there is $ \eps > 0 $ (depending on $ \lambda $ but not on $f$) such that, for $ a \to \infty $, $ c_{f, v}(kah) $ is equal to
        \[
            a^{\lambda - \rho} \langle f(k), \hat{J}^w_{\tilde{\sigma}, -\lambda} v(h^{-1}w^{-1}) \rangle
            + a^{-\lambda - \rho} \langle (\hat{J}^w_{\sigma, \lambda} f)(k), v(h^{-1}w) \rangle
            + O(a^{-\tfrac{\alpha}2 - \rho -\eps})
        \]
        uniformly in $ kM \in \dX $ varying in a compact subset of $U$. The remainder depends jointly on $ \lambda $ and $f$.
    \end{enumerate}
\end{lem}

\begin{proof} \nlenum
    \begin{enumerate}
    \item Let $ f \in \Cinf(\dX, V(\sigma_\lambda, \phi)) $ with $ \R(\lambda) > 0 $, $ v \in \C^K(\dX, V(\tilde{\sigma}_{-\lambda})) $ and $ k, h \in K $.
        By Theorem 5.3.4 of \cite[p.146]{Wallach},
        \[
            \lim_{a \to \infty} a^{\rho - \lambda} c_{f, v}(a)
            = \int_{\bar{N}}{ \langle f(e), v(\bar{n}) \rangle \, d\bar{n} }
            = \int_N{ \langle f(e), v(w^{-1} n w) \rangle \, dn } \ .
        \]
        Thus,
        \[
            \lim_{a \to \infty} a^{\rho - \lambda} c_{f, v}(kah)
            = \langle f(k), \hat{J}^w_{\tilde{\sigma}, -\lambda} v(h^{-1}w^{-1}) \rangle
        \]
        as $ c_{f, v}(kah) = c_{\pi^{\sigma, \lambda}(k^{-1})f, \pi^{\tilde{\sigma}, -\lambda}(h)v}(a) $.
        The above integral exists by Lemma 5.3.1 of \cite[p.144]{Wallach}.
        The first part of the assertion now follows by Theorem 4.4.3 of \cite[p.119]{Wallach}.
    \item Let $ f \in \Cinf(\dX, V(\sigma_{-\lambda}, \phi)) $ with $ \R(\lambda) > 0 $, $ v \in \C^K(\dX, V(\tilde{\sigma}_{\lambda})) $ and $ k, h \in K $.
        By Theorem 5.3.4 of \cite[p.146]{Wallach},
        \begin{multline*}
            \lim_{a \to \infty} a^{\rho - \lambda} c_{f, v}(kah)
            = \lim_{a \to \infty} a^{\rho - \lambda} c_{v, f}(h^{-1} w a w^{-1} k^{-1}) \\
            = \langle v(h^{-1}w), (\hat{J}^w_{-\lambda} f)(k) \rangle
            = \langle (\hat{J}^w_{-\lambda} f)(k), v(h^{-1}w) \rangle \ .
        \end{multline*}
        The above integral exists by Lemma 5.3.1 of \cite[p.144]{Wallach}.
        The second assertion follows now from Theorem 4.4.3 of \cite[p.119]{Wallach}.

        The second part of the first assertion follows from the two assertions we have just shown and the fact that the leading term in an asymptotic expansion depends holomorphically on $ \lambda \in \{ \mu \in \a^*_\CC \mid |\R(\mu)| < \tfrac12 , \, \mu \neq 0 \} $ (see sections 12.4, 12.5 and 12.6 of \cite{Wallach2}).
    \item Let $ W, W' $ be compact subsets of $U$ such that $ W \subset \interior(W') $. Let $ \chi \in \Ccinf(U) $ be such that $ \restricted{\chi}{W_1} \equiv 1 $. Then, $ f = \chi f + (1 - \chi)f $ with $ \chi f $ smooth and $ \supp\big((1-\chi)f\big) \subset \dX \smallsetminus \interior(W') $.
        The third assertion follows by applying \eqref{Lemma2 about Asymptotic expansions eq1} to $ \chi f $ and Lemma \ref{Lemma about Asymptotic expansions} to $ (1-\chi)f $ for $ kM \in W $.
    \end{enumerate}
\end{proof}

\begin{dfn}
    We call the smallest element $ \delta_{\Gamma} \in \a^{*} $ such that $ \sum_{\gamma \in \Gamma}{ a_{\gamma}^{-(\lambda + \rho) } } $ converges for all $ \lambda \in \a^{*} $ with $ \lambda > \delta_{\Gamma} $ the \textit{critical exponent} of $ \Gamma $. We set $ \delta_{\Gamma} = -\infty $ if $ \Gamma $ is the trivial group.\index{critical exponent}\index[n]{deltaGamma@$ \delta_{\Gamma} $}
\end{dfn}

\begin{rem} \nlenum
    \begin{enumerate}
    \item If $\Gamma$ is nontrivial, then $ \delta_\Gamma\in [-\rho,\rho) $. 
    \item It follows from \eqref{eq a_g = e^(d(eK, gK))} that this definition differs only by a $ \rho $-shift from the one defined in \cite{Corlette}.
    \end{enumerate}
\end{rem}

Let us now introduce the push-down\index[n]{pi@$ \pi_* $}
\[
    \pi_* \colon \C^{\infty}(\dX, V(\tilde{\sigma}_{-\lambda}, \tilde{\phi})) \to \C^{\infty}(B, V_B(\tilde{\sigma}_{-\lambda}, \tilde{\phi}))
\]
needed to define the extension map that we will define in a moment. Let $ \pi(\gamma) = \pi^{\tilde{\sigma}, -\lambda}(\gamma) \otimes \tilde{\phi}(\gamma) $ ($ \gamma \in \Gamma $). By using the identification
\[ \C^\infty(B, V_B(\tilde{\sigma}_{-\lambda}, \tilde{\phi})) = \Cinf(\Omega, V(\tilde{\sigma}_{-\lambda}, \tilde{\phi}))^\Gamma \ , \]
we define $ \pi_* $ by
\[
    \pi_*(f)(kM) = \sum_{\gamma \in \Gamma}{ (\pi(\gamma)f)(kM) } \qquad (f \in \C^{\infty}(\dX, V(\tilde{\sigma}_{-\lambda}, \tilde{\phi})), kM \in \Omega)
\]
if the sum converges \cite[p.88]{BO00}. By Lemma 4.2 of \cite[p.89]{BO00}, this is the case when $ \R(\lambda) > \delta_\Gamma $.

\begin{dfn}[{\cite[p.91]{BO00}}]\label{Dfn Bunke-Olbrich extension map}
    For $ \R(\lambda) > \delta_\Gamma $, the \textit{Bunke-Olbrich extension map}\index{Bunke-Olbrich extension map}\index[n]{ext@$ \ext $}
    \[
        \ext \colon \C^{-\infty}(B, V_B(\sigma_\lambda, \phi)) \to \C^{-\infty}(\dX, V(\sigma_\lambda, \phi))^\Gamma
    \]
    is defined to be the adjoint of $ \pi_* $.
\end{dfn}

\begin{rem}
    When we want to make our notation more precise, we denote $ \ext $ by $ \ext_\lambda $.
\end{rem}

If $ X \neq \OO H^2 $, then $ \ext $ has a meromorphic continuation to all of $ \a_\CC^* $ (see Theorem 5.10 of \mbox{\cite[p.103]{BO00}}).

\begin{dfn}[{\cite[p.92]{BO00}}]
    We define the restriction map\index{restriction map}\index[n]{rest@$ \rest $, $ \rest_\Omega $}
    \[
        \rest \colon \C^{-\infty}(\dX, V(\sigma_\lambda, \phi))^\Gamma \to \C^{-\infty}(B, V_B(\sigma_\lambda, \phi))
    \]
    as the composition of the identification map
    \[
        \C^{-\infty}(\Omega, V(\sigma_\lambda, \phi))^\Gamma \to \C^{-\infty}(B, V_B(\sigma_\lambda, \phi))
    \]
    with the restriction
    \[
        \rest_\Omega \colon \C^{-\infty}(\dX, V(\sigma_\lambda, \phi)) \to \C^{-\infty}(\Omega, V(\sigma_\lambda, \phi)) \ .
    \]
\end{dfn}

By Lemma 4.5 of \cite[p.92]{BO00}, it is the left-inverse of $ \ext $:
\begin{equation}\label{eq res o ext = Id}
    \rest \circ \ext = \Id \ .
\end{equation}

\newpage

\subsection{Invariant distributions supported on the limit set}

Recall that we denote the space of $\Gamma$-invariant distributions supported on the limit set by $ \C^{-\infty}(\Lambda, V(\sigma_\lambda, \phi))^\Gamma $. In the following, we provide some information about these spaces. For further details, see \cite{BO00}. They are important because $ L^2(\Gamma \bs G, \phi)_U \oplus L^2(\Gamma \bs G, \phi)_{res} $ is generated by matrix coefficients coming from $\Gamma$-invariant distributions supported on the limit set.

\medskip

By Theorem 4.7 of \cite[p.93]{BO00}, $ \C^{-\infty}(\Lambda, V(\sigma_\lambda, \phi))^\Gamma = \{0\} $ for all $ \lambda \in \a^*_\CC $ such that $ \R(\lambda) > \delta_\Gamma $.

Assume that $ X \neq \OO H^2 $. Then, $ \C^{-\infty}(\Lambda, V(\sigma_\lambda, \phi))^\Gamma $ is finite-dimensional for all $ \lambda \in \a^*_\CC $ by Theorem 6.1 of \cite[p.109]{BO00}.

\begin{dfn}[{\cite[p.126]{BO00}}]\label{Dfn U_Lambda(sigma_lambda, phi)}
    For $ \lambda \in \a_\CC^* $, define\index[n]{ULambda@$ U_\Lambda(\sigma_\lambda, \phi) $}
    \[
        U_\Lambda(\sigma_\lambda, \phi) := \{ \varphi \in \C^{-\infty}(\Lambda, V(\sigma_\lambda, \phi))^\Gamma \mid \rest \circ \hat{J}^w_{\sigma, \lambda}(\varphi) = 0 \} \ .
    \]
    We call it the \textit{space of ``stable'' $ \Gamma $-invariant distributions supported on the limit set}.
\end{dfn}

\begin{rem}
    The operator $ \rest \circ \hat{J}^w_{\sigma, \lambda} $ is well-defined on $ \C^{-\infty}(\Lambda, V(\sigma_\lambda, \phi))^\Gamma $ by Lemma 5.3 of \cite[p.98]{BO00}.
\end{rem}

We denote the space of germs at $ \lambda $ ($ \lambda \in \a_\CC^* $) of holomorphic families
\[
    \mu \mapsto f_\mu \in \C^{-\infty}(\dX, V(\sigma_\mu, \phi))  \qquad (\text{resp. } \mu \mapsto f_\mu \in \C^{-\infty}(B, V_B(\sigma_\mu, \phi)))
\]
by
\[
    \O_\lambda \C^{-\infty}(\dX, V(\sigma_{.}, \phi)) \qquad (\text{resp. } \O_\lambda \C^{-\infty}(B, V_B(\sigma_{.}, \phi))) \ .
\]
Let $ \O^0_\lambda \C^{-\infty}(B, V_B(\sigma_{.}, \phi)) $ be defined by
\[
     \{ f_\mu \in \O_\lambda \C^{-\infty}(B, V_B(\sigma_{.}, \phi)) \mid (\mu - \lambda) \ext f_\mu \in \O_\lambda \C^{-\infty}(\dX, V(\sigma_{.}, \phi)) \} \ .
\]

\begin{dfn} [{\cite[p.119]{BO00}}]\label{Dfn E_Lambda(sigma_lambda, phi)}
    Let\index[n]{ELambda@$ E_\Lambda(\sigma_\lambda, \phi) $}
    \[
        E_\Lambda(\sigma_\lambda, \phi) = \{ \Residue_{\mu = \lambda} \ext_\mu(f_\mu) \mid f_\mu \in \O^0_\lambda \C^{-\infty}(B, V_B(\sigma_{.}, \phi)) \}
    \]
    ($ \lambda \in \a^*_\CC $) be the space of $ \Gamma $-invariant distributions on the limit set generated by the singular parts of $ \ext $.
\end{dfn}

\begin{prop}\label{Prop res ext(f) supported on Lambda}
    For any $ \lambda \in \a^*_\CC $, $ E_\Lambda(\sigma_\lambda, \phi) $ is contained in $ \C^{-\infty}(\Lambda, V(\sigma_\lambda, \phi))^\Gamma $.
\end{prop}

\begin{proof}
    This follows from the identity $ \rest \circ \ext = \Id $. See Proposition 6.11 of \cite[p.119]{BO00} for a detailed proof.
\end{proof}

By Proposition 7.6 of \cite[p.126]{BO00}, we have
\[
    \C^{-\infty}(\Lambda, V(\sigma_\lambda, \phi))^\Gamma
    = E_\Lambda(\sigma_\lambda, \phi) \oplus U_\Lambda(\sigma_\lambda, \phi)
\]
for all $ \R(\lambda) \geq 0 $.

\begin{prop}\label{Prop C^(-infty)(Lambda, V(sigma_lambda))^Gamma neq 0 implies lambda in a^*}
    Let $ \lambda \in \a^*_\CC $ be such that $ \R(\lambda) \geq 0 $. Then,
    \[
        \C^{-\infty}(\Lambda, V(\sigma_\lambda, \phi))^\Gamma \neq \{0\}
    \]
    implies that $ \lambda \in \a^* $.
\end{prop}

\begin{rem}
    Proposition 7.8 of \cite[p.127]{BO00} provides even more information. It shows for example that $ \{ \lambda \in \a^*_\CC \mid \R(\lambda) \geq 0, \C^{-\infty}(\Lambda, V(\sigma_\lambda, \phi))^\Gamma \neq \{0\} \} $ is a finite subset of $ [0, \delta_\Gamma] $.
\end{rem}

\begin{proof}
    The proposition follows from Lemma 7.2 and Proposition 7.3 of \cite[p.123]{BO00}.
\end{proof}

For $ \sigma \in \hat{M} $ and $ \lambda \in \a_\CC^* $ with $ \R(\lambda) > 0 $, we denote by $ (\bar{\pi}^{\sigma, \lambda}, I^{\sigma,\lambda}_{\pm \infty}) $ the unique irreducible subrepresentation of the principal series representation $ \pi^{\sigma, \lambda} $ acting on $ C^{\pm\infty}(\dX, V(\sigma_\lambda)) $. We denote the underlying $(\g, K)$-module of $ I^{\sigma, \lambda}_{\infty} $ by $ I^{\sigma, \lambda} $.

The Langlands classification (cf. Chapter 5 of \cite{Wallach}), a unique Langlands parameter $ (\sigma, \lambda) $, $ \sigma \in \hat{M} $, $ \lambda \in \a^*_\CC : \R(\lambda) > 0 $, is associated to any irreducible nontempered representation $ (\pi, V_\pi) $ such that $ V_{\pi, \pm \infty} $ is equivalent to $ (\bar{\pi}^{\sigma, \lambda}, I^{\sigma,\lambda}_{\pm \infty}) $.

Let $ \sigma \in \hat{M} $, $ \lambda \in \a^*_\CC $ with $ \R(\lambda) > 0 $. By Proposition 9.2 of \cite[p.135]{BO00}, we have
\begin{multline*}
    (V_{\pi, -\infty} \otimes V_\phi)^\Gamma_d = (I^{\sigma, \lambda}_{-\infty} \otimes V_\phi)_d^\Gamma
    = (I^{\sigma, \lambda}_{-\infty} \otimes V_\phi)_{temp}^\Gamma \\
    = \C^{-\infty}(\Lambda, V(\sigma_\lambda, \phi))^\Gamma = E_\Lambda(\sigma_\lambda, \phi) \oplus U_\Lambda(\sigma_\lambda, \phi) \ ,
\end{multline*}
where $ \pi := \bar{\pi}^{\sigma, \lambda} $. If one of these spaces is nontrivial, then $ \pi \in \hat{G}_c $.

For $ \lambda \in \a^*_\CC $, set
\[
    \C^{-\infty}(\Lambda, V(\sigma_\lambda, \phi))_d^\Gamma
    = \C^{-\infty}(\dX, V(\sigma_\lambda, \phi))_d^\Gamma
    \cap \C^{-\infty}(\Lambda, V(\sigma_\lambda, \phi))^\Gamma \ .
\]

\begin{lem}\label{Lem C^(-infty)(Lambda, V(sigma_0, phi))_d^Gamma = U_Lambda(sigma_0, phi)}
    We have
    \[
        \C^{-\infty}(\dX, V(\sigma_0, \phi))_d^\Gamma
        = \C^{-\infty}(\Lambda, V(\sigma_0, \phi))_d^\Gamma = U_\Lambda(\sigma_0, \phi) \ .
    \]
\end{lem}

\begin{proof}
    It follows from Lemma \ref{Lemma about Asymptotic expansions} that $ f \in \C^{-\infty}(\Lambda, V(\sigma_0, \phi))^\Gamma $ is square-integrable if and only if $ \rest \circ \hat{J}_0 f = 0 $. Thus, $ \C^{-\infty}(\Lambda, V(\sigma_0, \phi))_d^\Gamma = U_\Lambda(\sigma_0, \phi) \, $.

    See Proposition 9.5 of \cite[p.138]{BO00} for a complete proof.
\end{proof}

If $ \pi = \bar{\pi}^{\sigma, \lambda} \in \hat{G}_c $, then $ \sigma = \sigma^w $ and $ \lambda \in \a^* $ (see e.g. Theorem 16.6 in \cite[p.656]{Knapp}). 

Let $ (\cdot, \cdot) $ be the scalar product on $ (I^{\sigma, \lambda}_{-\infty} \otimes V_\phi)_d^\Gamma $ induced by the matrix coefficient map to $ L^2(\Gamma \bs G, \phi) $.

By Corollary 10.5 of \cite[p.145]{BO00}, $ E_\Lambda(\sigma_\lambda, \phi) \oplus U_\Lambda(\sigma_\lambda, \phi) $ is orthogonal with respect to $ (\cdot, \cdot) $ if $ \R(\lambda) > 0 $. By Proposition 4.23 of \cite[p.45]{Olb02}, $ \C^{-\infty}(\Lambda, V(\sigma_0, \phi))^\Gamma $ is either $ E_\Lambda(\sigma_0, \phi) $ or $ U_\Lambda(\sigma_0, \phi) $.

\newpage

\subsection{Determination of the cuspidal invariant distributions supported on the limit set}

In this section, we determine the cuspidal invariant distribution vectors among the invariant distributions supported on the limit set.

\medskip

Recall that we assume that $ X \neq \OO H^2 $. Because of Proposition \ref{Prop C^(-infty)(Lambda, V(sigma_lambda))^Gamma neq 0 implies lambda in a^*}, we consider in the following only $ \lambda \in \bar{\a}^*_+ $.

\begin{lem}\label{Lem c_(f, v) in the Schwartz space}
    If $ \lambda \in \a^*_+ $, $ f \in \C^{-\infty}(\Lambda, V(\sigma_\lambda, \phi))^\Gamma $ or if $ f \in \C^{-\infty}(\Lambda, V(\sigma_0, \phi))_d^\Gamma $, then $f$ is a Schwartz vector.
\end{lem}

\begin{rem}
    If $ \lambda \in \a^*_+ $, then the lemma shows in particular that $ E_\Lambda(\sigma_\lambda, \phi) $ is contained in the space of Schwartz vectors by Proposition \ref{Prop res ext(f) supported on Lambda}.
\end{rem}

\begin{proof}
    By assumption, $ \lambda \geq 0 $.
    Let $ f $ be as above and let $ v \in \C^K(\dX, V(\tilde{\sigma}_{-\lambda})) $.
    \\ Let $ U \in \U_\Gamma $, $ X, Y \in \U(\g) $, $ k \in K(\Omega) $, $ h \in K $ and $ a \in \bar{A}_+ $.
    There exists an open relatively compact subset $ W $ of $ \Omega $ such that $ U $ is contained in $ \{ kaK \mid k \in K(W), a \in \bar{A}_+ \} $ modulo a relatively compact subset of $X$. Indeed, we can choose any $ W $ strictly containing $U$. Thus, we may assume without loss of generality that $U$ is of the form $ \{ kaK \mid k \in K(W), a \in \bar{A}_+ \} $.
    Hence, there exist smooth functions $ \psi_{h, X, Y, j} $ ($ j \in \NN $) on $ K(\Omega) $ such that
    \begin{multline*}
        L_X R_Y c_{f, v}(kah) = c_{\pi^{\sigma, \lambda}(X) f, \pi^{\tilde{\sigma}, -\lambda}(Y) v}(kah) \\
        = a^{-(\lambda + \rho)} \langle (\hat{J}^w_{\sigma, \lambda} \pi^{\sigma, \lambda}(X) f)(k), \pi^{\tilde{\sigma}, -\lambda}(Y) v(h^{-1}w) \rangle
            + \sum_{j=1}^\infty{ a^{-(\lambda + \rho + j\alpha)} \psi_{h, X, Y, j}(k) }
    \end{multline*}
    by Lemma \ref{Lemma about Asymptotic expansions}. If $ \lambda = 0 $, then $ f \in U_\Lambda(\sigma_0, \phi) $ by Lemma \ref{Lem C^(-infty)(Lambda, V(sigma_0, phi))_d^Gamma = U_Lambda(sigma_0, phi)}. So, if $ \lambda = 0 $, then the real part of the leading exponent of the asymptotic expansion for $ a \to \infty $ of $ L_X R_Y c_{f, v}(kah) $ is less or equal than $ -\rho - \alpha < -\rho $ as $ f \in U_\Lambda(\sigma_0, \phi) $ implies that $ \rest \circ \hat{J}^w_{\sigma, \lambda} \pi^{\sigma, \lambda}(X) f = \rest \circ \, \pi^{\sigma^w, -\lambda}(X) \hat{J}^w_{\sigma, \lambda} f = 0 $. If $ \lambda > 0 $, then it is less or equal than $ -\lambda - \rho < -\rho $.
    As in any case the real part of the leading exponent is less than $ \rho $, $ c_{f, v} \in \CS(\Gamma \bs G, \phi) $.
\end{proof}

\begin{prop}\label{Prop c_(f, v)^Omega formula} \nl
    If $ f \in \C^{-\infty}(\Lambda, V(\sigma_\lambda, \phi))^\Gamma $ with $ \lambda \in \a^*_+ $, then
    \begin{equation}\label{Prop c_(f, v)^Omega formula eq}
        c_{f, v}^\Omega(g, h) = \langle \hat{J}^w_{\sigma, \lambda} f(g), \hat{J}^w_{\tilde{\sigma}, -\lambda} v(h^{-1}) \rangle
    \end{equation}
    for all $ g \in G(\Omega) $, $ h \in G $ and $ v \in \C^K(\dX, V(\tilde{\sigma}_{-\lambda})) $.

    If $ f \in \C^{-\infty}(\Lambda, V(\sigma_0, \phi))_d^\Gamma = U_\Lambda(\sigma_0, \phi) $, then $ c_{f, v}^\Omega = 0 $ for all $ v \in \C^K(\dX, V(\tilde{\sigma}_0)) $.
\end{prop}

\begin{rem}
    The expression $ \hat{J}^w_{\sigma, \lambda} f(g) $ is well-defined for all $ g \in G(\Omega) $ and all $ f \in \C^{-\infty}(\Lambda, V(\sigma_\lambda, \phi))^\Gamma $ by Lemma 5.3 of \cite[p.98]{BO00}.
\end{rem}

\begin{proof}
    Fix $ g \in G(\Omega) $ and $ h \in G $. Let $ f, v$ be as above.
    By Lemma \ref{Lem c_(f, v) in the Schwartz space}, $ c_{f, v} $ belongs to $ \CS(\Gamma \bs G, \phi) $. So, $ (c_{f, v})^\Omega(g, h) $ is well-defined.

    Without loss of generality, we may assume that $ h = e $ (replace $v$ by $ \pi^{\tilde{\sigma}, -\lambda}(h)v $ in order to get the general formula).
    Then, $ c_{f, v}(gn) $ ($ n \in N $) is equal to
    \begin{equation}\label{Proof Lem c_(f, v)^Omega formula eq1}
        \langle f, \pi^{\tilde{\sigma}, -\lambda}(gn) v \rangle
        = \int_{K}{ \langle f(g n k), v(k) \rangle \, dk } \ .
    \end{equation}
    By Lemma \ref{Lem int_K f(k) dk = int_(bar N) f(kappa(bar n) a(bar n)^(-2rho) dbar(n)}, this yields
    \begin{equation}\label{Proof Prop c_(f, v)^Omega formula eq2}
        \int_N \langle f(g n n' w), v(n' w) \rangle \, dn'
        = \int_N \langle f(g n' w), v(n^{-1} n' w) \rangle \, dn' \ .
    \end{equation}
    Let $ 1_N \colon N \to \RR, \, n \mapsto 1 $. Since $ f $ has support in $Q$ and since the open set $ G \smallsetminus \{gP\} = g N w P $ is diffeomorphic to $N$, we can identify $ f $ with a distribution $ f_N $ in $ \C_c^{-\infty}(N, V_\sigma \otimes V_\phi) $ such that $ \int_N f(g n w) \, dn := f_N(1_N) $ is equal to $ \hat{J}^w_{\sigma, \lambda} f(g) $.

    Let $ u(n) = v(n^{-1} w) $ ($ n \in N $). Then, \eqref{Proof Prop c_(f, v)^Omega formula eq2} is again equal to
    \[
        (f_N * \check{u})(n) \ ,
    \]
    where $ \check{u}(n) := u(n^{-1}) $ and $ * $ denotes the convolution on $N$.
    We have
    \[
        |u(n)| = a(nw)^{-\lambda-\rho} |v(\kappa(n w))| \leq C a(nw)^{-\lambda-\rho}  \ ,
    \]
    where $ C := \|v\|_{\infty, K} = \sup_{k \in K} |v(k)| $.
    Let $ \{X_j\} $ be a basis of $ \U(\n)_l $, where $l$ denotes the order of $f_N$.
    As $f_N$ is a distribution on $N$ having compact support, there exists a compact subset $ S $ of $N$ such that
    \[
        |(f_N * \check{u})(n)| \prec \sum_j \|L_{X_j} u\|_{\infty, S} \prec \sup_{n' \in S} a(n^{-1} n' w)^{-\lambda-\rho}
        \prec a(n w)^{-\lambda-\rho}
    \]
    by Proposition \ref{Theorem 3.8 of Helgason, p.414}. Thus, $ f_N * \check{u} $ is integrable over $N$ if $ \lambda \in \a^*_+ $.

    For $ \lambda \in \a^*_+ $, it remains to prove the following:

    \begin{clm}
        If $ \lambda \in \a^*_+ $, then
        \[
            \int_N (f_N * \check{u})(n) \, dn = \langle f_N(1_N) , \int_N \check{u}(n) \, dn \rangle \ .
        \]
        In other words, $ c_{f, v}^\Omega(g, e) $ is equal to $ \langle \hat{J}^w_{\sigma, \lambda} f(g), \hat{J}^w_{\tilde{\sigma}, -\lambda} v(e) \rangle $.
    \end{clm}

    \begin{proof}
        If $f_N$ were a compactly supported regular distribution of $N$ with values in $ V(\sigma_\lambda, \phi) $, then the claim would directly follow by Fubini's theorem. Since this is not the case, we take the convolution with a test function in order to be able to apply Fubini.

        Let $ \phi \in \Ccinf(N) $ be such that $ \int_N \phi(n) \, dn = 1 $. Since $ f_N $ and $ \phi $ have compact support and since $ \phi * f_N $ is a smooth compactly supported function on $N$ with values in $ V_\sigma \otimes V_\phi $, 
        \begin{multline*}
            \int_N (f_N * \check{u})(n) \, dn = \int_N \phi(n) \, dn \cdot \int_N (f_N * \check{u})(n) \, dn = \int_N \big(\phi * (f_N * \check{u})\big)(n) \, dn \\
            = \int_N \big( (\phi * f_N) * u\big)(n) \, dn
            = \langle \int_N (\phi * f_N)(n) \, dn,  \int_N \check{u}(n) \, dn \rangle = \langle (\phi * f_N)(1_N), \int_N \check{u}(n) \, dn \rangle
        \end{multline*}
        and
        \begin{multline*}
            (\phi * f_N)(1_N) = \big((\phi * f_N) * 1_N\big)(e) = \big(\phi * (f_N * 1_N)\big)(e)
            = \big(\phi * f_N(1_N) 1_N\big)(e) \\
            = f_N(1_N) \big(\phi * 1_N\big)(e)
            = f_N(1_N) \int_N \phi(n) \, dn = f_N(1_N) \ .
        \end{multline*}
        Thus,
        \[
            \int_N (f_N * \check{u})(n) \, dn = \langle f_N(1_N) , \int_N \check{u}(n) \, dn \rangle \ .
        \]
        Here, $ f_N(1_N) = \hat{J}^w_{\sigma, \lambda} f(g) $ and
        \[
            \int_N \check{u}(n) \, dn = \int_N u(n^{-1}) \, dn
            = \int_N v(n w) \, dn
            = \hat{J}^w_{\tilde{\sigma}, -\lambda} v(e) \ .
        \]
    \end{proof}

    \needspace{2\baselineskip}
    Let us now look at the case $ \lambda = 0 $.

    De Rham theory asserts that we have the following assertions:

    \begin{lem}
        Let $ M $ be a connected, oriented smooth manifold of dimension $n$.
        Let $ \Omega^{*}_{c, \infty}(M) $ denote the de Rham complex of compactly supported smooth complex-valued functions on $M$.
        Then,
        \[
            H^n(\Omega^{*}_{c, \infty}(M)) \simeq \CC \ .
        \]
        The isomorphism is given by $ [\omega] \mapsto \int \omega \in \CC $.
    \end{lem}

    \begin{proof}
        For a proof, see for example Corollary 5.8 of \cite[p.47]{BottTu}.
    \end{proof}

    \begin{lem}
        Let $ M $ be a connected, oriented smooth manifold of dimension $n$ and let $ \Omega^{*}_{c, -\infty}(M) $ be the de Rham complex of distributional forms (of compact support) on $M$. Then,
        \[
            \Omega^{*}_{c, \infty}(M) \hookrightarrow \Omega^{*}_{c, -\infty}(M)
        \]
        induces an isomorphism in cohomology.
    \end{lem}

    \begin{proof}
        See Theorem 14 of \cite[p.94]{DeRham}.
    \end{proof}

    It follows from the two lemmas that $ H^n(\Omega^{*}_{c, -\infty}(M)) \simeq \CC $. The isomorphism is given by
    \[
        [\omega] \mapsto \omega(1_M) \qquad (\omega \in \Omega^{*}_{c, -\infty}(M)) \ .
    \]
    In particular, if $ \omega(1_M) = 0 $ for some $ \omega \in \Omega^n_{c, -\infty}(M) $, then
    \begin{equation}\label{eq omega = d eta}
        \omega = d\eta
    \end{equation}
    for some $ \eta \in \Omega^{n-1}_{c, -\infty}(M) $.

    \begin{clm}\label{Clm Estimate of |(L_X v)(n w)|}
        Let $ X \in \n $. Then,
        \[
            |(L_X v)(n w)| \prec a(nw)^{-(\frac{\alpha}2 + \rho)}
        \]
        for all $ n \in N $.
    \end{clm}

    \begin{proof}
        Let $ n \in N $. Write $ n $ as $ a \tilde{n} a^{-1} $ with $ \tilde{n} \in N $ and $ a \in A $ such that $ a_B(\tilde{n} w) = 1 $.
        \\ Since
        \[
            a \tilde{n} a^{-1} w
            = a \tilde{n} w a = (a \bar{n}_B(\tilde{n} w) a^{-1}) m_B(\tilde{n} w) (a_B(\tilde{n} w) a^2) (a^{-1} n_B(\tilde{n} w) a) \ ,
        \]
        $ \bar{n}_B(n w) = \bar{n}_B(\tilde{n}^a w) = a \bar{n}_B(\tilde{n} w) a^{-1} $ and $ a_B(n w) = a_B(\tilde{n}^a w) = a_B(\tilde{n} w) a^2 = a^2 $.
        So,
        \[
            a = a_B(n w)^{\frac12} \ .
        \]
        Let $ S = \{ n \in N \mid |\log(n)| \geq 1 \} $. Let $ X \in \n $.
        Since $ \bar{n} \in \bar{N} \mapsto v(\bar{n}) $ belongs to $ \C^1(\bar{N}, V_{\tilde{\sigma}}) $,
        \[
            v(\bar{n}) = v(e) + O(|\log \bar{n}|)
        \]
        and $ L_X v(\bar{n}) = O(|\log \bar{n}|) $. 
        Let
        \[
            c = \sup \{|\log \bar{n}_B(n w)| \mid n \in N : a_B(nw) = 1 \} < \infty \ .
        \]
        Then, there is $ c' > 0 $ such that
        \begin{multline*}
            |L_X v(n w)| = |L_X v(\bar{n}_B(n w))| a_B(n w)^{-\rho} = |L_X v(a \bar{n}_B(\tilde{n} w) a^{-1} )| a^{-2\rho} \\
            \leq c' a^{-2\rho} |\log a \bar{n}_B(\tilde{n} w) a^{-1}|
            \leq c c' a^{-\alpha - 2\rho} = c c' a_B(n w)^{-(\frac{\alpha}2 + \rho)}
        \end{multline*}
        for every $ n \in S $.
        The claim follows as $ \{ n \in N \mid |\log(n)| \leq 1 \} $ is compact and as $ a_B(n w) \asymp a(n w) $ for every $ n \in S $, by Proposition \ref{Theorem 3.8 of Helgason, p.414}.
    \end{proof}

    By Lemma \ref{Lem C^(-infty)(Lambda, V(sigma_0, phi))_d^Gamma = U_Lambda(sigma_0, phi)},$ \C^{-\infty}(\Lambda, V(\sigma_0, \phi))_d^\Gamma = U_\Lambda(\sigma_0, \phi) $.
    Then, $ f_N(1_N) = 0 $.

    Let $ \{Z_j\} $ be a basis of $ \n $ respecting the decomposition $ \n = \g_\alpha \oplus \g_{2\alpha} $. By \eqref{eq omega = d eta}, there exist finitely many $ S_j \in \C_c^{-\infty}(N, V_\sigma \otimes V_\phi) $ such that
    \[
        f_N = \sum_j R_{Z_j} S_j \ .
    \]
    Thus,
    \[
        f_N * \check{u} = \sum_j (R_{Z_j} S_j) * \check{u} = -\sum_j S_j * (L_{Z_j} \check{u}) \ .
    \]
    Hence,
    \[
        \int_N (f_N * \check{u})(n) \, dn = -\sum_j \langle S_j(1_N) , \int_N (L_{Z_j} v)(n w) \, dn \rangle \ .
    \]
    It remains to show that $ \int_N (L_{Z_j} v)(n w) \, dn = 0 $ for all $ j $.

    Consider first the case $ \dim \n = 1 $.
    Since $ |(L_{Z_j} v)(n w)| \prec a(nw)^{-(\frac{\alpha}2 + \rho)} $ and since $ v(nw) $ converges to zero when $ |\log n| $ tends to $ \infty $, it follows that $ \int_N (L_{Z_j} v)(n w) \, dn = 0 $ for any $j$, by partial integration.

    Assume in the following that $ \dim \n \geq 2 $.
    Let $ \n_j = \RR Z_j $ and $ \tilde{\n}_j = \bigoplus_{i \neq j} \RR Z_i $. Then, $ \n = \n_j \oplus \tilde{\n}_j $ is an $A$-invariant decomposition with $ \dim \tilde{\n}_j \geq 1 $.

    Let $ N_j = \exp(\n_j) $ and $ \tilde{N}_j = \exp(\tilde{n}_j) $. We normalise the Haar measure $ dn_j $ on $ N_j \equiv \RR $ such that
    \[
        \int_{N_j} F(n_j) \, dn_j = \int_{\RR} F(\exp(x Z_j)) \, dx
    \]
    for all $ F \in L^1(N_j) $. By 4.A.2.1 of \cite[p.133]{Wallach}, we can normalise the measure $ d\tilde{n}_j $ on $ \tilde{N}_j $, given by the push-forward of the Lebesgue measure on $ \tilde{\n}_j $, such that
    \[
        \int_N F(n) \, dn = \int_{N_j} \int_{\tilde{N}_j} F(n_j \tilde{n}_j) \, d\tilde{n}_j \, dn_j
    \]
    for all $ F \in L^1(N) $. Thus, if $ L_{Z_j} F $ is integrable over $N$, then $ \int_N L_{Z_j} F(n) \, dn $ is equal to
    \[
        \int_\RR \int_{\tilde{N}_j} L_{Z_j} F(\exp(x Z_j) \tilde{n}_j) \, d\tilde{n}_j \, dx
        = \lim_{r \to \infty} \int_{-r}^r \int_{\tilde{N}_j} L_{Z_j} F(\exp(x Z_j) \tilde{n}_j) \, d\tilde{n}_j \, dx \ .
    \]
    Since $ |(L_{Z_j} v)(n w)| \prec a(nw)^{-(\frac{\alpha}2 + \rho)} $, $ (L_{Z_j} v)(\cdot \, w) $ is integrable over $N$ by Claim \ref{Clm Estimate of |(L_X v)(n w)|}.
    As moreover $ F(n_j \tilde{n}_j) $ is integrable over $ \tilde{N}_j $ when $ n_j $ varies in a compact subset of $ N_j $, the above is again equal to
    \[
        -\lim_{r \to \infty} \int_{-r}^r \restricted{ \frac{d}{dt} }{t=0} \int_{\tilde{N}_j} F(\exp((x + t) Z_j) \tilde{n}_j) \, d\tilde{n}_j \, dx
    \]
    by the theorem about differentiation of parameter dependent integrals. This yields
    \[
        -\lim_{r \to \infty} \int_{\tilde{N}_j} F(\exp(r Z_j) \tilde{n}_j) \, d\tilde{n}_j \, dx
        + \lim_{r \to \infty} \int_{\tilde{N}_j} F(\exp(-r Z_j) \tilde{n}_j) \, d\tilde{n}_j \, dx
    \]
    by partial integration. This limit is equal to zero as $ \int_{\tilde{N}_j} v(\exp(x Z_j) \tilde{n}_j w) \, d\tilde{n}_j $ converges to zero when $ |x| $ tends to $ \infty $ by the following claim.

    \begin{clm}
        $ \int_{\tilde{N}_j} a(n_j \tilde{n}_j w)^{-\rho} \, d\tilde{n}_j $ converges to zero when $ |\log \tilde{n}_j| $ converges to $ \infty $.
    \end{clm}

    \begin{proof}
        Let $ E_j = \{ \exp(X + Y) \mid X \in \g_\alpha, Y \in \g_{2\alpha} : \tfrac14|X|^4 + 2 |Y|^2 \leq 1 \} \cap \tilde{N}_j $ (compact).
        Since $ \int_{E_j} a(n_j \tilde{n}_j w)^{-\rho} \, d\tilde{n}_j $ converges to zero when $ |\log \tilde{n}_j| $ converges to $ \infty $ and since
        $ a(n w) \asymp a_B(n w) $ on $ N \smallsetminus E_j $ by Corollary \ref{Cor a_n asymp a(theta n)}, it suffices to show that
        \begin{equation}\label{eq int_(tilde(N)_j - E_j) a_B(n_j tilde n_j w)^(-rho) dn}
            \lim_{t \to \pm \infty} \int_{\tilde{N}_j \smallsetminus E_j} a_B(\exp(t Z_j) \tilde{n}_j w)^{-\rho}  \, d\tilde{n}_j = 0 \ .
        \end{equation}
        As $ a_B(n_j \tilde{n}_j w) \asymp \sqrt{a_B(n_j w)^2 + a_B(\tilde{n}_j w)^2} $ for all $ n_j \in N_j $ with $ |\log n_j| \geq 1 $ and all $ \tilde{n}_j \in E_j $.
        Let $ m_j = 1 $ if $ Z_j \in \g_\alpha $ and $ m_j = 2 $ if $ Z_j \in \g_{2\alpha} $. By Lemma \ref{Lem int_E f(n) dn = int_A int_(varphi(x) = 1) f(x^a) dx a^(2 xi_E rho) da}, we have \eqref{eq int_(tilde(N)_j - E_j) a_B(n_j tilde n_j w)^(-rho) dn} if and only if
        \begin{multline*}
            \lim_{t \to \pm \infty} \int_{\bar{A}_+} \int_{\{x \in \tilde{N}_j \mid a_B(x w) = 1\}} \big(a_B(\exp(t Z_j) w)^2 + a^4 a_B(x w)^2\big)^{-\rho} \, dx
            \ a^{2 \rho - m_j} \, da = 0 \\
            \iff \lim_{t \to \infty} \int_{\bar{A}_+} \big(t^{\frac2{m_j}} + a^4\big)^{-\frac{\rho}2} a^{2\rho - m_j} \, da = 0
            \iff \lim_{t \to \infty} \int_{\bar{A}_+} \frac{a^{2\rho - m_j}}{t^{\frac{\rho}{m_j}} + a^{2\rho}} \, da = 0 \ .
        \end{multline*}
        But
        \begin{multline}\label{eq int_(bar(A)_+) a^(2rho - m_j)/(t^(rho/m_j) + a^(2rho))}
            \int_{\bar{A}_+} \frac{a^{2\rho - m_j}}{t^{\frac{\rho}{m_j}} + a^{2\rho}} \, da
            = \int_1^\infty \frac{x^{2\rho - m_j - 1}}{t^{\frac{\rho}{m_j}} + x^{2\rho}} \, dx
            = \int_1^\infty x^{- m_j - 1} \cdot \frac{x^{2\rho}}{t^{\frac{\rho}{m_j}} + x^{2\rho}} \, dx \\
            = \int_1^\infty x^{- m_j - 1} \, dx
            - \int_1^\infty x^{- m_j - 1} \cdot \frac{t^{\frac{\rho}{m_j}}}{t^{\frac{\rho}{m_j}} + x^{2\rho}} \, dx \ .
        \end{multline}
        Since
        \[
            x^{- m_j - 1} \cdot \frac{t^{\frac{\rho}{m_j}}}{t^{\frac{\rho}{m_j}} + x^{2\rho}} \leq x^{- m_j - 1} \leq \frac1{x^2}
        \]
        for all $ x \geq 1 $ and all $ t \geq 0 $, \eqref{eq int_(bar(A)_+) a^(2rho - m_j)/(t^(rho/m_j) + a^(2rho))} converges to
        \[
            \int_1^\infty x^{- m_j - 1} \, dx - \int_1^\infty x^{- m_j - 1} \, dx =  0
        \]
        by Lebesgue's theorem of dominated convergence.
    \end{proof}
    So, $ c_{f, v}^\Omega(g, e) = \int_N (T * \check{u})(n) \, dn $ vanishes.
    This finishes the proof of the proposition.
\end{proof}

\begin{prop}\label{Prop U_Lambda(sigma_lambda, phi) consists of cusp forms}
    If $ f \in U_\Lambda(\sigma_\lambda, \phi) $ with $ \lambda \in \a^*_+ $ or $ \lambda = 0 $, then $f$ is a cuspidal Schwartz vector.
\end{prop}

\begin{proof}
    Let $ v \in \C^K(\dX, V(\tilde{\sigma}_{-\lambda})) $ ($ \lambda \in \bar{\a}^*_+ $). By Lemma \ref{Lem c_(f, v) in the Schwartz space}, $ c_{f, v} \in \CS(\Gamma \bs G, \phi) $ and, by Proposition \ref{Prop c_(f, v)^Omega formula}, $ c_{f, v}^\Omega(g, h) $ is equal to zero as $ f \in U_\Lambda(\sigma_\lambda, \phi) $ implies that $ \hat{J}^w_{\sigma, \lambda} f(g) = 0 $ for all $ g \in G(\Omega) $. Hence, $ c_{f, v} \in \c{\CS}(\Gamma \bs G, \phi) $.
\end{proof}

\begin{thm}\label{Thm c_(f, v) in 0CS(Gamma|G, phi) iff f in U_Lambda(sigma_lambda, phi)} \nl
    If $ f \in \C^{-\infty}(\Lambda, V(\sigma_\lambda, \phi))^\Gamma $ with $ \lambda \in \a^*_+ $ or if $ f \in \C^{-\infty}(\Lambda, V(\sigma_0, \phi))_d^\Gamma $, then $f$ is a cuspidal Schwartz vector if and only if $ f $ belongs to $ U_\Lambda(\sigma_\lambda, \phi) $.
\end{thm}

\begin{proof}
    The implication from the right to the left follows by Proposition \ref{Prop U_Lambda(sigma_lambda, phi) consists of cusp forms}.
    As the implication in the other direction is trivial for $ \lambda = 0 $ by Lemma \ref{Lem C^(-infty)(Lambda, V(sigma_0, phi))_d^Gamma = U_Lambda(sigma_0, phi)}. Let now $ \lambda \in \a^*_+ $.
    Without loss of generality, we may assume that $ \C^{-\infty}(\Lambda, V(\sigma_\lambda, \phi))^\Gamma \neq \{0\} $. Let $ f $ be as above such that $ c_{f, v} \in \c{\CS}(\Gamma \bs G, \phi) $.
    Let us assume now that there is $ x \in G(\Omega) $ such that $ \hat{J}^w_{\sigma, \lambda} f(x) \neq 0 $.
    Then, there is $ v' \in V(\tilde{\sigma}^w_{\lambda}) $ such that
    \begin{equation}\label{Proof Prop c_(f, v) in °CS(Gamma|G, phi) iff f in U_Lambda(sigma_lambda, phi) eq1}
        \langle \hat{J}^w_{\sigma, \lambda} f(x), v' \rangle \neq 0 \ .
    \end{equation}
    Let $ C = \int_{\bar{N}} a(\bar{n})^{-\lambda-\rho} \, d\bar{n} > 0 $ (finite by Proposition \ref{Prop N bar integral is finite} applied to $ E = N $) and let $ g \in \C^K(\dX, V(\tilde{\sigma}^w_{-\lambda})) $ be such that $ g(k) = \frac1C v' $ for all $ k \in K $. Then, $ \hat{J}^w_{\tilde{\sigma}, -\lambda} g(k) = v' $ for all $ k \in K $ by construction.
    By Proposition \ref{Prop c_(f, v)^Omega formula}, we have
    \[
        \langle \hat{J}^w_{\sigma, \lambda} f(x), v' \rangle
        = \langle \hat{J}^w_{\sigma, \lambda} f(x), \hat{J}^w_{\tilde{\sigma}, -\lambda} g(k) \rangle = 0
    \]
    for all $ k \in K $. This contradicts \eqref{Proof Prop c_(f, v) in °CS(Gamma|G, phi) iff f in U_Lambda(sigma_lambda, phi) eq1}. So, $ \hat{J}^w_{\sigma, \lambda} f(x) = 0 $ for all $ x \in G(\Omega) $. In other words, $ f \in U_\Lambda(\sigma_\lambda, \phi) $. This finishes the proof of the theorem.
\end{proof}

\newpage

\subsection{Discrete series representations}\label{ssec: Discrete series representations}

We recall how a discrete series representation can be embedded into a principal series representation and we show how the embedding map is related to the Poisson transform, which can be viewed as sum of matrix coefficients.

\medskip

\begin{dfn} \nlenum
    \begin{enumerate}
    \item We say that a representation $ (\pi, V_\pi) $ is \textit{square-integrable} if all the matrix coefficients $ c_{v, \tilde{v}} $ ($ v \in V_{\pi, \infty} , \tilde{v} \in V_{\pi', K} $) belong to $ L^2(G) $.
    \item A representation $ (\pi, V_\pi) $ is called a \textit{discrete series representation} if it is an irreducible square-integrable representation of $G$.
    \item A representation $ (\pi, V_\pi) $ is said to be \textit{tempered} if all the matrix coefficients $ c_{v, \tilde{v}} $ ($ v \in V_{\pi, \infty}, \tilde{v} \in V_{\pi', K} $) belong to $ \CS'(G) $.
    \end{enumerate}
\end{dfn}

\begin{rem}\label{Rem discrete series rep and tempered rep} \nlenum
    \begin{enumerate}
    \item It follows from the previous definition that a discrete series representation is tempered.
    \item By Lemma 9.4 of \cite[p.137]{BO00}, we know the following:\label{Rem discrete series rep and tempered rep, Lemma 9.4}
        \\ If $ (\pi, V_\pi) \in \hat{G} $ is tempered, then $ (V_{\pi, -\infty} \otimes V_\phi)^\Gamma_{temp} = (V_{\pi, -\infty} \otimes V_\phi)^\Gamma $. Moreover,
        \[
            (V_{\pi, -\infty} \otimes V_\phi)^\Gamma \ni f \mapsto c_{f, v} \in \CS'(\Gamma \bs G, \phi)
        \]
        is continuous.
    \item $ \pi $ is a discrete series representation (resp. tempered representation) if and only if $ \pi' $ is a discrete series representation (resp. tempered representation).
    \end{enumerate}
\end{rem}

Let $ (\gamma, V_\gamma) $ be a finite-dimensional representation of $K$. Recall that $ V(\gamma) := G \times_K V_\gamma $ (homogeneous vector bundle over $X$) and that $ V(\gamma, \phi) := V(\gamma) \otimes V_\phi $.
Sections of $V(\gamma)$ will be viewed as functions from $G$ to $V_\gamma$ satisfying
\[
    f(gk) = \gamma(k)^{-1} f(g)
\]
for all $ g \in G $ and $ k \in K $.

\begin{dfn}[{\cite[p.93]{BO00}}]
    Let $ T \in \Hom_M(V_\sigma, V_\gamma) $ and $ \lambda \in \a^*_\CC $.
    Then, the \textit{Poisson transform}\index{Poisson transform}\index[n]{PTlambdaf@$ P^T_\lambda f $}
    \[
        P^T_\lambda \colon \C^{-\infty}(\dX, V(\sigma_\lambda)) \to \Cinf(X, V(\gamma))
    \]
    is defined by
    \[
        (P^T_\lambda f)(g) :=
        \int_K a(g^{-1}k)^{-(\lambda+\rho)} \gamma(\kappa(g^{-1}k)) T f(k) \, dk \ ,
    \]
    where the integral is a formal notation meaning that the distribution $f$ has to be applied to the smooth integral kernel.
\end{dfn}

\begin{rem}
    We denote $ P^{T \otimes \Id}_\lambda \otimes \Id \colon \C^{-\infty}(\dX, V(\sigma_\lambda, \phi)) \to \Cinf(X, V(\gamma, \phi)) $ also by $ P^T_\lambda $.
\end{rem}

\begin{dfn}
    Let $ (\pi, V_\pi) $ be an admissible representation of $G$ of finite length on a reflexive Banach space. Let $ \Lambda_{V_{\pi, K}} \in \a^* $ be the real part of the highest leading exponent appearing in an asymptotic expansion for $ a \to \infty $ of $ c_{v, \tilde{v}}(ka) $ ($ v \in V_{\pi, \infty} $, $ \tilde{v} \in V_{\pi', K} $).
\end{dfn}

\begin{rem}
    Alternatively, one could define $ \Lambda_{V_{\pi, K}} $ as in 4.3.5 of \cite[p.116]{Wallach}.
\end{rem}

For any $ v \in V_{\tilde{\gamma}} $, $ \mu \in \a^*_\CC $ and $ T \in \Hom_M(V_\sigma, V_\gamma) $, define $ v_{T, \mu} \in \Cinf(\dX, V_{\tilde{\sigma}_\mu}) $ by
\begin{equation}\label{eq Def of v_T}
    v_{T, \mu}(k) = \LOI{T}{t} \tilde{\gamma}(k^{-1}) v \qquad (k \in K) \ .
\end{equation}
Here, $ \tilde{\gamma} $ denotes the dual representation of $ \gamma $.\index[n]{vT@$ v_{T, \mu} $}

For the rest of this section, $ (\pi, V_\pi) $ will denote a discrete series representation unless stated otherwise. Let us recall in the following the construction of an embedding map $ \beta \colon V_{\pi, -\infty} \to \C^{-\infty}(\dX, V(\sigma_{-\lambda})) $.

Casselman's subrepresentation theorem provides an injective $G$-intertwining operator from $ V_{\pi, \infty} $ to $ \Cinf(\dX, V(\sigma_{-\lambda})) $. Let us give the explicit construction.

Let $ \lambda \in \a^* $ be such that $ -(\lambda + \rho) = \Lambda_{V_{\pi, K}} $. As $ (\pi, V_\pi) $ is a discrete series representation, $ 0 < \lambda \in \a^* $.

Let $ \gamma \in \hat{K} $ and $ V_\gamma \subset V_\pi(\gamma) $ be such that $ V_\gamma $ is a nonzero irreducible submodule of $ V_\pi(\gamma) $ (e.g. take the minimal $K$-type $(\gamma, V_\gamma) $ of $ V_\pi $) and choose an embedding $ 0 \neq \tilde{t} \in \Hom_K(V_{\tilde{\gamma}}, V_{\pi'}) $. 
For $ v \in V_{\pi, -\infty} $, define $ F_v \in \Cinf(G, V_\gamma) $ by
\[
    \langle F_v(g), \tilde{v} \rangle = \langle \pi(g^{-1}) v, \tilde{t} \tilde{v} \rangle \qquad (\tilde{v} \in V_{\tilde{\gamma}}) \ .
\]
Then, $ F_v(ka) $ has an asymptotic expansion for $ v \in V_{\pi, \infty} $ by Theorem 4.4.3 of \cite[p.119]{Wallach}.
So, there exists a nonzero map $ \tilde{\beta} \colon V_{\pi, \infty} \to \Cinf(\dX, V((\restricted{\gamma}{M})_{-\lambda}) $ (which is a $G$-intertwining operator) such that
\[
    F_v(ka) = a^{-\lambda - \rho} (\tilde{\beta} v)(k) + O(a^{-\lambda - \rho - \eps}) \ .
\]
Choose $ \sigma \in \hat{M} $ such that $ \Hom_M(V_{\gamma}, V_{\sigma}) \neq \{0\} $.
Let $ t \in \Hom_M(V_{\gamma}, V_{\sigma}) $ be such that $ t \tilde{\beta} v $ is nonzero for some $ v \in V_{\pi, \infty} $.
\\ For $ v \in V_{\pi, \infty} $, set $ \beta v = t \tilde{\beta} v \in \Cinf(\dX, V(\sigma_{-\lambda})) $.

It follows from \cite{Collingwood} that $ V_\pi $ appears with multiplicity one in $ L^2(\dX, V(\sigma_{-\lambda})) $.
Hence, $ \beta $ is unique up to a constant when $ \sigma $ is fixed.

\begin{dfn}\label{Dfn topological embedding}
    Let $ X_1 $ and $ X_2 $ be two topological spaces.
    We say that $ f \colon X_1 \to X_2 $ is a \textit{topological embedding} if it is a homeomorphism onto its image (equipped with the subspace topology).\index{topological embedding}
\end{dfn}

It follows from Casselman-Wallach globalisation theory (cf. Chapter 11 in \cite{Wallach2}) that the image of $ \beta $ is closed in $ \Cinf(\dX, V(\sigma_{-\lambda})) $. Thus, $ \beta $ is a topological embedding.
We get a projection
\[
    q \colon \C^{-\infty}(\dX, V(\sigma_\lambda)) \to V_{\pi, -\infty}
\]
by forming the adjoint with respect to Hermitian scalar products.
We can extend $ \beta $ to a map between the corresponding spaces of distribution vectors by functoriality and get a $G$-intertwining operator
\[
    A := \beta \circ q \colon \C^{-\infty}(\dX, V(\sigma_\lambda)) \to \C^{-\infty}(\dX, V(\sigma_{-\lambda})) \ .
\]
Thus, we have the following commutative diagram:
\begin{eqnarray*}
    \xymatrix@R-1pc@C-2pc{
        \C^{-\infty}(\dX, V(\sigma_\lambda)) \ar@{->}[rr]^{\LOI{ A := \beta \circ q }{ }} \ar@{->>}[rd]_{q}
        & & \C^{-\infty}(\dX, V(\sigma_{-\lambda})) \\
        & V_{\pi, -\infty} \ar@{^{(}->}[ur]_{\beta}
    }
\end{eqnarray*}
Let us denote the operator $ \beta \otimes \Id $ (resp. $ q \otimes \Id $, $ A \otimes \Id $) also by $ \beta $ (resp. $ q $, $ A $).

\begin{lem}\label{Lem beta is a top embedding}
    The maps $ \beta \colon V_{\pi, -\infty} \otimes V_\phi \to \C^{-\infty}(\dX, V(\sigma_{-\lambda}, \phi)) $, $ \beta \colon (V_{\pi, -\infty} \otimes V_\phi)^\Gamma \to \C^{-\infty}(\dX, V(\sigma_{-\lambda}, \phi))^\Gamma $ are topological embeddings.
\end{lem}

\begin{proof}
    Since $ \beta \colon V_{\pi, -\infty} \otimes V_\phi \to \C^{-\infty}(\dX, V(\sigma_{-\lambda}, \phi)) $ is a continuous injection and since it follows from the globalisation theory of Casselman and Wallach that $ \im(\beta) $ is closed in $ \C^{-\infty}(\dX, V(\sigma_{-\lambda}, \phi)) $, $ \beta \colon V_{\pi, -\infty} \otimes V_\phi \to \im(\beta) $ is a homeomorphism when we endow $ \im(\beta) $ with the subspace topology.
    The first assertion follows.

    The second assertion follows as $ \beta(v) $ ($ v \in V_{\pi, -\infty} \otimes V_\phi $) is $ \Gamma $-invariant if and only if $ v $ is $ \Gamma $-invariant.
\end{proof}

Let us recall now the relation between the Poisson transform $ P^T_\lambda $ and the matrix coefficients of principal series representations. It is proven by a standard computation using integral formulas.

\begin{lem}\label{Lem BO00, 9.1}
    Let $ \sigma $ (resp. $ \gamma $) be a finite-dimensional representation of $M$ (resp. $K$) and let $ T \in \Hom_M(V_\sigma, V_\gamma) $.
    Consider the Poisson transform
    \[
        P^T_\lambda \colon \C^{-\infty}(\dX, V(\sigma_\lambda, \phi)) \to \Cinf(X, V(\gamma, \phi)) \simeq [\Cinf(G, V_\phi) \otimes V_\gamma]^K \ .
    \]
    Then, for any $ v \in V_{\tilde{\gamma}} $, $ f \in \C^{-\infty}(\dX, V(\sigma_\lambda, \phi)) $ and $ g \in G $, we have
    \[
        \langle P^T_\lambda f(g), v \rangle = c_{f, v_{T, -\lambda}}(g) \in V_\phi \ .
    \]
\end{lem}

Let us show in the following how the Poisson transform and $ \tilde{\beta} \circ q $ are related.

Let $ v \in V_{\pi, -\infty} \otimes V_\phi $ and $ f \in \C^{-\infty}(\dX, V(\sigma_\lambda, \phi)) $ be such that $ qf = v $. Then,
\begin{equation}\label{eq <F_v(g), tilde v> = c_(v, tilde(t) tilde(v))(g)}
    \langle F_v(g), \tilde{v} \rangle := \langle \pi(g^{-1}) v, \tilde{t} \tilde{v} \rangle = c_{v, \tilde{t} \tilde{v}}(g)
    = c_{f, \t{q}(\tilde{t} \tilde{v})}(g)
\end{equation}
for all $ g \in G $ and $ \tilde{v} \in V_{\tilde{\gamma}} $.
Define $ T \in \Hom_M(V_\sigma, V_\gamma) $ by
\begin{equation}\label{eq Def of T}
    \langle T(v), \tilde{v} \rangle = \langle v, [\t{q} \circ \tilde{t}(\tilde{v})](e) \rangle
    \qquad (v \in V_\sigma, \tilde{v} \in V_{\tilde{\gamma}}) \ .
\end{equation}
Observe that $ \tilde{v}_{T, -\lambda} = \t{q}(\tilde{t}(\tilde{v})) $ for all $ \tilde{v} \in V_{\tilde{\gamma}} $.
Thus,
\begin{equation}\label{eq F_v(g) = P^T f(g)}
    F_v(g) = P^T_\lambda f(g)   \qquad (g \in G)
\end{equation}
by Lemma \ref{Lem BO00, 9.1}. Hence, we have the following:

Let $ f \in \C^{-\infty}(\dX, V(\sigma_\lambda, \phi)) $ be such that $ q f \in V_{\pi, \infty} \otimes V_\phi $ and let $ k \in K $. Then,
\begin{equation}\label{eq DS leading term of the Poisson transform}
    \lim_{a \to \infty} a^{\lambda + \rho} P^T_\lambda f(ka) = (\tilde{\beta} (q f))(k) \ .
\end{equation}

\begin{lem}\label{Lem Poisson transform injective ds rep}
    The map
    \[
        \C^{-\infty}(\dX, V(\sigma_\lambda, \phi))/\ker(q) \to \Cinf(X, V(\gamma, \phi)) , \quad f \mapsto P^T f
    \]
    is well-defined and injective.
\end{lem}

\begin{proof}
    By \eqref{eq F_v(g) = P^T f(g)}, $ F_{q f} = P^T_\lambda f $. Thus, the above map is well-defined.
    As $ \langle F_{q f}(g), \tilde{v} \rangle = c_{q f, \tilde{t} \tilde{v}}(g) $ for all $ g \in G $ and $ \tilde{v} \in V_{\tilde{\gamma}} $, by \eqref{eq <F_v(g), tilde v> = c_(v, tilde(t) tilde(v))(g)} and as
    $ v \in V_{\pi, -\infty} \otimes V_\phi \mapsto c_{v, \tilde{t} \tilde{v}} $ is injective for every nonzero $ \tilde{v} \in V_{\tilde{\gamma}} $, the lemma follows.
\end{proof}

\newpage

\subsubsection{Schwartz vectors are cuspidal}

The main result, which we prove in this section, is that a Schwartz vector of a discrete series representation is cuspidal (see Proposition \ref{Prop Every Schwartz vector is cuspidal}).

\medskip

It follows from Definition \ref{Dfn Schwartz vector} that $ qf \in (V_{\pi, -\infty} \otimes V_\phi)^\Gamma $ is a Schwartz vector (resp. cuspidal Schwartz vector) if $ c_{f, v} \in \CS(\Gamma \bs G, \phi) $ (resp. $ c_{f, v} \in \c{\CS}(\Gamma \bs G, \phi) $) for all $ v \in \C^K(\dX, V(\tilde{\sigma}_{-\lambda})) \cap \im(\t{q}) $.

\begin{dfn}[{\cite[p.14]{Borel}}]\label{DfnDiracSequence}
    We say that a sequence $ (\phi_n)_{n \in \NN} $ in $ \Ccinf(G, \RR) $ is a \textit{Dirac sequence} if\index{Dirac sequence}
    \begin{enumerate}
    \item $ \phi_n \geq 0 \quad \forall n \in \NN $;\label{DiracProp1}
    \item $ \supp(\phi_n) \underset{n \to \infty}{\longrightarrow} \{ e \} $: For any neighbourhood $U$ of $e$, there is $ N \in \NN $ such that $ \supp(\phi_n) \subset U $ for all $ n \geq N $;\label{DiracProp2}
    \item $ \int_G{ \phi_n(g) \, dg } = 1 \quad \forall n \in \NN $.\label{DiracProp3}
    \end{enumerate}
\end{dfn}

\begin{rem}
    Let $ (\phi_n)_{n \in \NN} $ be a Dirac sequence in $ \Ccinf(G, \RR) $ and let $ f \in \Cinf(G, V_\phi) $. Then, $ \phi_n * f $ converges to $f$.
\end{rem}

\begin{lem}\label{LemDirac}
    Let $ f \in \CS(\Gamma \bs G, \phi) $. Let $ U, U' \in \U_\Gamma $ such that $ U' \subset U $. Let $ V $ be a symmetric, open, relatively compact neighbourhood of the neutral element $ e \in G $. We assume that $ V U' \subset U $. Let $ (\phi_n)_{n \in \NN} $ be a Dirac sequence, let $ X, Y \in \U(\g) $ and let $ r \geq 0 $. Then, $ \LOI{p}{U'}_{r, X, Y} (\phi_n * f - f) $ converges to zero when $n$ goes to $ \infty $.
\end{lem}

\begin{proof} Let $ f \in \CS(\Gamma \bs G, \phi) $ and let $ U, U'$ and $V$ be as above.

    Let $ \eps > 0 $ be given. Let $ X \in \U(\g)_l $ and $ Y \in \U(\g) $. Let $ \{X_i\} $ be a basis of $ \U(\g)_l $. By the proof of Proposition \ref{Prop convergence of Up_(r, X, Y)(R_h f - f)}, there exists a symmetric, open, relatively compact neighbourhood $W \subset V$ of $e$ such that
    \begin{equation}\label{Proof of LemDirac eq1}
        \sup_{x \in W} \LOI{p}{U'}_{r, X, Y}(L_x f - f) \leq \frac{\eps}2 \ , \qquad \sup_{x \in W} \LOI{p}{U'}_{r, X_i, Y}(L_x f - f) \leq \frac{\eps}2
    \end{equation}
    for all $i$. Because of Property \eqref{DiracProp2} of Definition \ref{DfnDiracSequence}, there exists $ N \in \NN $ such that $ \supp(\phi_j) \subset W $ for all $ j \geq N $.
    \\ Consequently, for $ j \geq N $ and $ gK \in U' $, we have by Property \eqref{DiracProp3} of Definition \ref{DfnDiracSequence} that
    \begin{align*}
        & |L_X R_Y (\phi_j * f)(g) - L_X R_Y f(g)| \\
        =& |\int_G{ \phi_j(x) L_{\Ad(x^{-1}) X} R_Y (L_x f)(g) \, dx } - \int_G{ L_X R_Y f(g) \phi_j(x) \, dx }| \\
        \leq& \int_G{ |L_{\Ad(x^{-1}) X} R_Y (L_x f)(g) - L_X R_Y f(g)| \phi_j(x) \, dx } \\
        \leq& \int_G{ |L_{\Ad(x^{-1}) X - X} R_Y (L_x f)(g)| \phi_j(x) \, dx } \\
        & \qquad + \int_G{ |L_X R_Y (L_x f)(g) - L_X R_Y f(g)| \phi_j(x) \, dx } \ .
    \end{align*}
    \needspace{2\baselineskip}
    Let us look now at the two integrals separately:
    \begin{enumerate}
    \item The first integral is less or equal than
        \begin{multline*}
            \sup_{x \in \supp(\phi_j)} |L_{\Ad(x^{-1}) X - X} R_Y (L_x f)(g)| \\
            \leq \sum_i (\sup_{x \in \supp(\phi_j)} |c_i(x)|) \sup_i \sup_{x \in W} |L_{X_i} R_Y (L_x f)(g)| \ ,
        \end{multline*}
        where $ c_i $ are smooth functions on $G$ such that
        \[
            \Ad(x^{-1}) X - X = \sum_i c_i(x) X_i \ .
        \]
        As moreover $ c_i(e) = 0 $ for any $i$, $ \sup_{x \in \supp(\phi_j)} |c_i(x)| $ converges to zero when $ j $ tends to $ \infty $, for every $i$.
        So, there is $ N' \geq N $ such that
        \begin{equation}\label{Proof of LemDirac eq2}
            \sum_i (\sup_{x \in \supp(\phi_j)} |c_i(x)|) \sup_i \sup_{x \in W} \LOI{p}{U'}_{r, X_i, Y}(L_x f - f) < \frac{\eps}2
        \end{equation}
        for every $ j \geq N' $.
    \item By Property \eqref{DiracProp1} of Definition \ref{DfnDiracSequence} and as $ \supp(\phi_j) \subset W $,
        the second integral is less or equal than
        \begin{align*}
            & \int_W{|L_X R_Y (L_x f)(g) - L_X R_Y f(g)|\phi_j(x) \, dx} \\
            \leq & \sup_{x \in W} |L_X R_Y (L_x f)(g) - L_X R_Y f(g)| \ .
        \end{align*}
    \end{enumerate}
    Thus,
    \[
        \LOI{p}{U'}_{r, X, Y}(\phi_j * f - f)
        \leq \sup_{x \in W} \LOI{p}{U'}_{r, X, Y}(L_x f - f) < \eps
    \]
    for every $ j \geq N' $, by \eqref{Proof of LemDirac eq1} and \eqref{Proof of LemDirac eq2}. The lemma follows.
\end{proof}

Let $ g \in G(\Omega) $ and $ h \in G $. Recall that $ U_{g, h} $ denotes an open set in $X$ that is a relatively compact neighbourhood of $ gP $ in $ X \cup \Omega $ and contains the horosphere $ gNhK = gNg^{-1}(gh)K $ passing through $ ghK \in X $ and $ gP \in \dX $. Let $ U $ be an open subset of $X$ containing $ U_{g, h} $ that is relatively compact in $ X \cup \Omega $.
Let $ V $ be a symmetric, open, relatively compact neighbourhood of the neutral element $ e \in G $ such that $ V U_{g, h} \subset U $.

The following two results are due to Harish-Chandra.
\begin{thm}\label{Thm 7.2.2}
    Let $ f \in \CS(G) $ and assume that $ \dim Z_G(\g)f < \infty $. Then, $ f \in \c{\CS}(G) $.
\end{thm}

\begin{proof}
    For a proof, see Theorem 7.2.2 of \cite[p.233]{Wallach}.
\end{proof}

\begin{cor}\label{Cor c_(f, v) is a cusp form if pi is a discrete series rep}
    Let $ f \in V_{\pi, \infty} $ and $ v \in V_{\pi', K} $. Then, $ c_{f, v} \in \c{\CS}(G) $.
\end{cor}

\begin{proof}
    Let $ k, h \in K $, $ a \in A $ and $ X, Y \in \U(\g) $. Since $ \pi $ is a discrete series representation of $G$, $ \dim Z_G(\g)f < \infty $ and the leading exponent of the asymptotic expansion of $ L_X R_Y c_{f, v}(k a h) = c_{\pi(X) f, \pi'(Y) v}(k a h) $ for $ a \to \infty $ is less or equal than $ -\lambda - \rho $. Thus, $ c_{f, v} \in \CS(G) $. Hence, $ c_{f, v} \in \c{\CS}(G) $ by Theorem \ref{Thm 7.2.2}.
\end{proof}

\begin{dfn}\label{Dfn pi(f)v}
    Let $ (\pi, V_\pi) $ be a representation of $G$ on a reflexive Banach space. Let $ f \in \C_c(G) $, $ v \in V_\pi $. Then,
    \[
        \pi(f)v := \int_G f(g) \pi(g)v \, dg \in V_\pi \ .
    \]
    For $ f \in \C_c(G) $ and $ v \in V_\pi \otimes V_\phi $, we set $ \pi(f)v = (\pi(f) \otimes \Id)v $.
\end{dfn}

\begin{rem}
    Let $ f \in \Ccinf(G) $ and $ v \in V_{\pi, -\infty} $. Then, $ \pi(f)v \in V_{\pi, \infty} $ (see, e.g., \cite[p.136]{Dijk}).
\end{rem}

\begin{prop}\label{Prop Every Schwartz vector is cuspidal}
    If $ f \in (V_{\pi, -\infty} \otimes V_\phi)^\Gamma $ is a Schwartz vector, then $f$ is cuspidal.
\end{prop}

\begin{proof}
    Let $ (\phi_j)_j $ be a Dirac sequence in $ \Ccinf(G, \RR) $. Let $ v \in V_{\pi', K} $.
    As $ c_{f, v} \in \CS(\Gamma \bs G, \phi) \subset \Cinf(\Gamma \bs G, \phi) \subset \Cinf(G, V_\phi) $, $ c_{\pi(\phi_j)f, v} = \phi_j * c_{f, v} $ converges to $ c_{f, v} $ in $ \Cinf(G, V_\phi) $ and in particular pointwise. Moreover, as $ \phi_j \in \Ccinf(G) $, $ \pi(\phi_j)f \in V_{\pi, \infty} \otimes V_\phi $.
    Since furthermore $ \pi $ is a discrete series representation, $ \phi_j * c_{f, v} = c_{\pi(\phi_j)f, v} $ belongs to $ \c{\CS}(G, V_\phi) $ by Corollary \ref{Cor c_(f, v) is a cusp form if pi is a discrete series rep}. So,
    $
        \int_N{ \phi_j * c_{f, v}(gnh) \, dn } = 0 
    $.
    \\ Fix $ r > 1 $. Then, $ \int_N{ |\phi_j * c_{f, v}(gnh) - c_{f, v} (gnh)| \, dn } $ is less or equal than
    \[
        \LOI{p}{U}_{r, 1, 1} (\phi_j * c_{f, v} - c_{f, v}) \int_N{ a_{g n h}^{-\rho} (1 + \log a_{g n h})^{-r} \, dn } \ ,
    \]
    where $ \int_N{ a_{g n h}^{-\rho} (1 + \log a_{g n h})^{-r} \, dn } $ is finite by Lemma \ref{Lem Estimate for int_N |f(g n h)| dn} applied to $ E = N $. By Lemma \ref{LemDirac}, this converges to zero. Hence,
    \[
        (c_{f, v})^\Omega(g, h)
        = \int_N{ c_{f, v}(gnh) \, dn }
        = \lim_{j \to \infty} \int_N{ \phi_j * c_{f, v}(gnh) \, dn } = 0 \ .
    \]
\end{proof}

\begin{lem}\label{Lem space of Schwartz vector is top. isom. to its image}
    Fix $ 0 \neq v \in V_{\pi, K} $. Then, $ (V_{\pi', -\infty} \otimes V_\phi)^\Gamma_\CS $ is a Fréchet space which is topologically isomorphic to $ \{ c_{\pi}(f \otimes v) \mid f \in (V_{\pi', -\infty} \otimes V_\phi)^\Gamma_\CS \} $.
    \\ Moreover, $ (V_{\pi', -\infty} \otimes V_\phi)^\Gamma_\CS \otimes V_{\pi}(\gamma) $ ($ \gamma \in \hat{K} $) is topologically isomorphic to
    \begin{equation}\label{eq Lem space of Schwartz vector is top. isom. to its image}
        \spn\{ c_{\pi}(f \otimes v)
        \mid f \in (V_{\pi', -\infty} \otimes V_\phi)^\Gamma_\CS , \ v \in V_{\pi}(\gamma) \} \ .
    \end{equation}
\end{lem}

\begin{proof}
    It follows from Corollary \ref{Cor c_pi(N_pi hat(otimes) V_pi), c_pi(N_pi otimes F) closed} combined with the first part of Lemma \ref{Lem V = (S cap V) oplus (S^perp cap V)} that $ \{ c_{\pi}(f \otimes v) \mid f \in (V_{\pi', -\infty} \otimes V_\phi)^\Gamma_\CS \} $ and \eqref{eq Lem space of Schwartz vector is top. isom. to its image} are closed.

    Fix $ 0 \neq v \in V_{\pi, K} $. Since moreover we take the coarsest topology on $ (V_{\pi', -\infty} \otimes V_\phi)^\Gamma_\CS $ such that $ f \in (V_{\pi', -\infty} \otimes V_\phi)^\Gamma_\CS \mapsto c_{\pi}(f \otimes v) \in \CS(\Gamma \bs G, \phi) $ is continuous and since $ f \mapsto c_{\pi}(f \otimes v) $ is injective, $ (V_{\pi', -\infty} \otimes V_\phi)^\Gamma_\CS $ is topologically isomorphic to $ \{ c_{\pi}(f \otimes v) \mid f \in (V_{\pi', -\infty} \otimes V_\phi)^\Gamma_\CS \} $. Hence, $ (V_{\pi', -\infty} \otimes V_\phi)^\Gamma_\CS $ is a Fréchet space, too.
    This shows the first part. As it follows from this part that $ (V_{\pi', -\infty} \otimes V_\phi)^\Gamma_\CS \otimes V_\pi(\gamma) $ is a Fréchet space and as
    \[
        c_{\pi} \colon V_{\pi', -\infty} \otimes V_\phi \to \CS'(G, V_\phi), \ f \mapsto c_{\pi}(f \otimes v)
    \]
    is injective, one can prove the second part similarly.
\end{proof}

\newpage

\subsubsection{The \texorpdfstring{$ (Af, g) $}{(Af, g)}-formula}

The $ (Af, g) $-formula (see Theorem \ref{Thm (Af, g) = C (res(f), res(g))}) is an explicit formula for the pairing between certain invariant distributions. This formula is absolutely necessary in order to be, for example, able to describe the space of Schwartz vectors more explicitly.

\medskip

Assume that $ X \neq \OO H^2 $ or that $ \delta_\Gamma < 0 $ and $ \lambda \in (0, -\delta_\Gamma) $ for the remaining part about the discrete series representations. Then, $ \ext_\lambda $ and $ \ext_{-\lambda} $ are defined. Recall that $ \ext_\mu $ is only defined if $ \R(\mu) > \delta_\Gamma $, when $ X = \OO H^2 $.

Define $ \check{{ }} \, \colon \Cinf(K) \to \Cinf(K) $ by $ \check{\psi}(k) = \psi(k^{-1}) $.

Recall that $ T \in \Hom_M(V_\sigma, V_\gamma) $ is defined by \eqref{eq Def of T}.

\begin{lem}\label{Lem Averaged asymptotics of P^T_(-lambda) g(ka)}
    Let $ f \in \Cinf(\dX, V((\restricted{\gamma}{M})_\lambda, \phi)) $ and $ g \in \C^{-\infty}(\dX, V(\sigma_\lambda, \phi)) $.
    Then, there exists $ \eps > 0 $ such that for $ a \to \infty $ we have
    \[
        \int_K (f(k), (P^T_\lambda g)(ka)) \, dk
        = a^{-\lambda - \rho} \int_K (f(k), (\tilde{\beta}(qg))(k)) \, dk + O(a^{-\lambda -\rho -\eps}) \ .
    \]
\end{lem}

\begin{proof}
    We prove this lemma by doing a similar argument as in the proof of Corollary 6.3 of \cite[p.112]{BO00}, which itself is adapted from \cite[5.1]{Schlichtkrull}.
    \\ Let $f$ and $g$ be as above.
    Let $ \{v_i\} $ be a basis of $ V_\gamma \otimes V_\phi $. So, there are smooth functions $ f_i $ on $K$ such that $ f(k) = \sum_i f_i(k) v_i $.
    Then, $ \int_K (f(k), (P^T_\lambda g)(ka)) \, dk $ is equal to
    \begin{multline*}
        \sum_i (v_i, \int_K P^T_\lambda g(ka) \overline{f_i(k)} \, dk) \\
        = \sum_i (v_i, \int_K \Big(\int_K a(a^{-1}h)^{\lambda-\rho} \gamma(\kappa(a^{-1}h)) Tg(kh) \, dh \Big) \overline{f_i(k)} \, dk)
    \end{multline*}
    as $K$ is unimodular. By Fubini's theorem, this yields
    \[
        \sum_i (v_i, \int_K a(a^{-1}h)^{\lambda-\rho} \gamma(\kappa(a^{-1}h)) T \int_K \overline{f_i(k)} g(kh) \, dk \, dh)
        = \sum_i (v_i, P^T_\lambda(\check{\bar{f}}_i *_K g)(a)) \ ,
    \]
    where $ *_K $ denotes the convolution on $K$.
    Since $ \check{\bar{f}}_i *_K g \in \Cinf(\dX, V(\sigma_\lambda, \phi)) $, the above is equal to
    \begin{equation}\label{Proof of Lem Averaged asymptotics of P^T_(-lambda) g(ka) eq1}
        a^{-\lambda - \rho} \sum_i (v_i, \tilde{\beta}(q(\check{\bar{f}}_i *_K g))(e)) + O(a^{-\lambda -\rho -\eps})
    \end{equation}
    by \eqref{eq DS leading term of the Poisson transform}. Here, $e$ denotes the neutral element of $G$.
    The lemma follows by doing similar computations as before backwards:
    Since
    \[
        \tilde{\beta}(q(\check{\bar{f}}_i *_K g))(e) 
        = \int_K \overline{f_i(k)} \big(\pi^{\sigma, -\lambda}(k^{-1})\tilde{\beta}(qg)\big)(e) \, dk = \int_K \overline{f_i(k)} (\tilde{\beta}(qg))(k) \, dk \ ,
    \]
    \eqref{Proof of Lem Averaged asymptotics of P^T_(-lambda) g(ka) eq1} is again equal to
    \begin{align*}
        & a^{-\lambda - \rho} \sum_i \int_K f_i(k) (v_i, (\tilde{\beta}(qg))(k)) \, dk + O(a^{-\lambda -\rho -\eps}) \\
        =& a^{-\lambda - \rho} \int_K (f(k), (\tilde{\beta}(qg))(k)) \, dk + O(a^{-\lambda -\rho -\eps}) \ .
    \end{align*}
\end{proof}

\begin{lem}\label{Lem A off-diagonally smoothing}
    The map $ A $ is off-diagonally smoothing:
    \[
        A\big(\C^{-\infty}_\Omega(\dX, V(\sigma_\lambda, \phi))\big) \subset \C^{-\infty}_\Omega(\dX, V(\sigma_{-\lambda}, \phi)) \ .
    \]
\end{lem}

\begin{proof}
    This is shown in the proof of Lemma 9.6 of \cite[p.139]{BO00}.
\end{proof}

\begin{lem}\label{Lem Asymptotics of c (Af, v)(kah) for -lambda} \nl
    Let $ f \in \C^{-\infty}(\dX, V(\sigma_\lambda, \phi)) $ be such that $ Af \in \C^{-\infty}_\Omega(\dX, V(\sigma_{-\lambda}, \phi)) $, let $ v \in \C^K(\dX, V(\tilde{\sigma}_{-\lambda})) \cap \im(\t{q}) $ and let $ h \in K $.
    Then, there exist $ \eps > 0 $ and $ g_h \in \Cinf(K(\Omega), V_\phi) $ such that for $ a \to \infty $ we have
    \[
        c_{f, v}(kah) = a^{-\lambda - \rho} g_h(k) + O(a^{-\lambda -\rho -\eps})
    \]
    uniformly as $ kM $ varies in a compact subset of $ \Omega $.
\end{lem}

\begin{proof}
    Fix $ f $, $ v = \t{A}\tilde{v} $ and $ h $ as above. Without loss of generality, we may assume that there is $ k \in K(\Omega) $ such that $ c_{f, v}(k a h) $ is nonzero for $ a \in A_+ $ sufficiently large (otherwise the assertion of the lemma is trivial).

    Let $W$ be a compact subset of $ \Omega $ and let $ W_1 $ be a compact subset of $ \Omega $ such that $ W \subset \interior(W_1) $. Let $ \chi \in \Ccinf(\Omega) $ be such that $ \restricted{\chi}{W_1} \equiv 1 $.
    Then, $ Af = \chi Af + (1-\chi)Af $ with $ (1 - \chi) Af $ having support in $ \dX \smallsetminus \interior(W_1) $ and with $ \chi Af $ belonging to $ \Cinf(\dX, V(\sigma_{-\lambda}, \phi)) $ by Lemma \ref{Lem A off-diagonally smoothing}.
    By Lemma \ref{Lemma2 about Asymptotic expansions} and by Lemma \ref{Lemma about Asymptotic expansions}, $ c_{\chi Af, \tilde{v}}(kah) $ and $ c_{(1 - \chi) Af, \tilde{v}}(kah) $ have an asymptotic expansion for $ a \to \infty $ and $kM \in W $.
    Since moreover we can choose $W$ arbitrarily, $ c_{f, v}(kah) = c_{Af, \tilde{v}}(kah) $, $ kM \in \Omega $, has an asymptotic expansion for $ a \to \infty $, too. Let $ \mu - \rho $ be its leading exponent.
    As discrete series representations have regular infinitesimal characters, there are functions $ 0 \not \equiv g_{1, h} \in \Cinf(K(\Omega), V_\phi), g_2 \in \Cinf(K(\Omega)A_{> t}K, V_\phi) $ for some $ t \geq 1 $ and $ \eps > 0 $ such that 
    \[
        c_{f, v}(kah) = a^{\mu - \rho} g_{1, h}(k) + g_2(kah) \ ,
    \]
    where $ |g_2(kah)| \leq C a^{\mu -\rho -\eps} $ ($ kM \in W $, $ a \in A_{> t} $) for some constant $ C > 0 $. Choose $W$ sufficiently big so that $ g_{1, h} $ is not identically zero.

    Let us now prove that $ \mu \leq -\lambda $.
    We denote the convolution on $K$ by $*_K$.
    Let $ \psi \colon K \to \RR $ be a smooth function with compact support in $ WM \subset K $ such that $ \int_K g_{1, h}(k) \psi(k) \, dk \neq 0 $. Then, on one hand $ \check{\psi} *_K c_{f, v}(ah) $ is equal to
    \[
        \int_K \psi(k) c_{f, v}(kah) \, dk = a^{\mu - \rho} \int_K \psi(k) g_{1, h}(k) \, dk + \int_K \psi(k) g_2(kah) \, dk
    \]
    and $ \int_K |\psi(k) g_2(kah)| \, dk \leq \sup_{k \in WM} |g_2(kah)| \leq C a^{\mu -\rho -\eps} $.
    \\ On the other hand, $ \check{\psi} *_K c_{f, v}(ah) $ is equal to
    $
        c_{\pi^{\sigma, \lambda}(\check{\psi})f, v}(ah).
    $
    \\ As $ \psi \in \Cinf(K) $, $ \pi(\check{\psi})q(f) \in V_{\pi, \infty} \otimes V_\phi $. Thus, $ c_{\pi^{\sigma, \lambda}(\check{\psi})f, v} = c_{\pi(\check{\psi})q(f), \t{\beta}(\tilde{v})} $ is a matrix coefficient of the discrete series representation $ V_\pi $.
    It follows now from the definition of $ \lambda $ that $ \mu \leq -\lambda $.
    This completes the proof of the lemma.
\end{proof}

\begin{cor}\label{Cor of Lem Asymptotics of c (Af, v)(kah) for -lambda}
    We also have the following analogous statement: \\
    There exists $ \eps > 0 $ such that for $ a \to \infty $ we have
    \[
        (P^T_\lambda f)(ka) = a^{-\lambda - \rho} (\tilde{\beta} (qf))(k) + O(a^{-\lambda -\rho -\eps})
    \]
    uniformly as $ kM $ varies in a compact subset of $ \Omega $.
\end{cor}

\begin{proof}
    The corollary follows from the previous lemma.
    That the leading term of $ (P^T_\lambda f)(ka) $ is equal to $ (\tilde{\beta} (qf))(k) $ follows from Lemma \ref{Lem Averaged asymptotics of P^T_(-lambda) g(ka)}.
\end{proof}

\begin{prop}\label{Prop Af, smooth on Omega, invariant distributions are Schwartz}
    For every distribution $ f \in \C^{-\infty}(\dX, V(\sigma_\lambda, \phi)) $ such that $ Af \in \C^{-\infty}_\Omega(\dX, V(\sigma_{-\lambda}, \phi))^\Gamma $, $ qf$ is a cuspidal Schwartz vector.
\end{prop}

\begin{proof}
    Let $ f $ be as above and let $ v = \t{A}\tilde{v} \in \C^K(\dX, V(\tilde{\sigma}_{-\lambda})) \cap \im(\t{q}) $. Then, $ c_{f, v} = c_{Af, \tilde{v}} $ is left $ \Gamma $-equivariant since $ Af \in \C^{-\infty}(\dX, V(\sigma_{-\lambda}, \phi))^\Gamma $.
    For all $ X, Y \in \U(\g) $,
    \[
        L_X R_Y c_{f, v}(kah) = c_{\pi^{\sigma, \lambda}(X) f, \pi^{\tilde{\sigma}, -\lambda}(Y) v}(kah) \ .
    \]
    \begin{sloppypar}
        Since the $G$-intertwining operator $ A \colon \C^{-\infty}(\dX, V(\sigma_\lambda, \phi)) \to \C^{-\infty}(\dX, V(\sigma_{-\lambda}, \phi)) $ induces a $ \g $-intertwining operator between the same spaces and since the $G$-operator $ \t{A} \colon \Cinf(\dX, V(\tilde{\sigma}_\lambda)) \to \Cinf(\dX, V(\tilde{\sigma}_{-\lambda})) $ induces a $ (\g, K) $-intertwining operator between $ \C^K(\dX, V(\tilde{\sigma}_\lambda)) $ and $ \C^K(\dX, V(\tilde{\sigma}_{-\lambda})) $,
    \end{sloppypar}
    \[
        A(\pi^{\sigma, \lambda}(X) f) = \pi^{\sigma, -\lambda}(X) Af \in \im(A) \cap \C^{-\infty}_\Omega(\dX, V(\sigma_{-\lambda}, \phi))
    \]
    and
    \[
        \pi^{\tilde{\sigma}, -\lambda}(Y) v = \pi^{\tilde{\sigma}, -\lambda}(Y) \t{A} \tilde{v} = \t{A} \big(\pi^{\tilde{\sigma}, \lambda}(Y) \tilde{v}\big) \in \C^K(\dX, V(\tilde{\sigma}_{-\lambda})) \cap \im(\t{q}) \ .
    \]
    Thus, $ c_{f, v} \in \CS(\Gamma \bs G, \phi) $ by Lemma \ref{Lem Asymptotics of c (Af, v)(kah) for -lambda}. Hence, $ qf $ is a Schwartz vector. It follows by Proposition \ref{Prop Every Schwartz vector is cuspidal} that $qf$ is cuspidal.
\end{proof}

Let $ (\sigma, V_\sigma) $ be a finite-dimensional unitary representation of $M$ and let $ \mu \in \a^*_\CC $.

As $ \sigma $ and $ \phi $ are unitary representations, we have a natural positive definite sesquilinear pairing on $ V_\sigma \otimes V_\phi $, which we denote by $ (\cdot, \cdot) $.

For $ p \in P $, $ v_1 \in V_{\sigma_\mu} \otimes V_\phi $, $ v_2 \in V_{\sigma_{-\bar{\mu}}} \otimes V_\phi $, set $ p.(v_1, v_2) = (\sigma_\mu(p) v_1, \sigma_{-\bar{\mu}}(p) v_2) $. This defines a left $P$-action on $ (V_{\sigma_\mu} \otimes V_\phi) \times (V_{\sigma_{-\bar{\mu}}} \otimes V_\phi) $.

Then, $ (\cdot, \cdot) $ is a $P$-equivariant map from $ (V_{\sigma_\mu} \otimes V_\phi) \otimes (V_{\sigma_{-\bar{\mu}}} \otimes V_\phi) $ to $ \CC_{1_{-\rho}} $.
\\ Indeed, for $ p \in P $, $ v_1 \in V_{\sigma_\mu} \otimes V_\phi $ and $ v_2 \in V_{\sigma_{-\bar{\mu}}} \otimes V_\phi $, we have
\[
    (\sigma_\mu(p) v_1, \sigma_{-\bar{\mu}}(p) v_2) = a(p)^{\rho - \mu} \overline{a(p)^{\rho + \bar{\mu}}} (v_1, v_2) = a(p)^{2\rho} (v_1, v_2)
    = 1_{-\rho}(p) (v_1, v_2) \ .
\]
Let $ f_1 \in \Cinf(B, V_B(\sigma_\mu, \phi)) $ and $ f_2 \in \Cinf(B, V_B(\sigma_{-\bar{\mu}}, \phi)) $. Then, $ (f_1(x), f_2(x)) $ belongs to $ V_B(1_{-\rho}) \simeq \Lambda^{max} T^* B $. In other words, it is a volume form on $B = \Gamma \bs \Omega$.
So,
\begin{equation}\label{eq Def of (.,.)_B}
    (f_1, f_2)_B := \int_{\Gamma \bs \Omega}{ (f_1(x), f_2(x)) }
\end{equation}
defines a natural pairing between $ \Cinf(B, V_B(\sigma_\mu, \phi)) $ and $ \Cinf(B, V_B(\sigma_{-\bar{\mu}}, \phi)) $.

Recall that $ \chi_\infty $ denotes the restriction of $ \chi $ to $ \Omega $.

\begin{lem}\label{Lem int_B (f_1(x), f_2(x)) = int_K chi_infty(kM) (f_1(k), f_2(k)) dk}
    Let $ f_1 \in \Cinf(B, V_B(\sigma_\mu, \phi)) $ and $ f_2 \in \Cinf(B, V_B(\sigma_{-\bar{\mu}}, \phi)) $. Then,
    \[
        \int_{\Gamma \bs \Omega} (f_1(x), f_2(x)) = \int_K \chi_\infty(kM) (f_1(k), f_2(k)) \, dk \ .
    \]
\end{lem}

\begin{proof}
    We have:
    \begin{multline*}
        \int_K \chi_\infty(kM) (f_1(k), f_2(k)) \, dk = \int_\Omega \chi_\infty(kM) (f_1(k), f_2(k)) \, dk_M \\
        = \int_{\Gamma \bs \Omega} \sum_{\gamma \in \Gamma} \gamma^*(\chi_\infty(kM) (f_1(k), f_2(k)) \, dk_M) \ .
    \end{multline*}
    Since $ (f_1(k), f_2(k)) \, dk $ is a $ \Gamma $-invariant form and since $ \sum_{\gamma \in \Gamma} \chi_\infty(\kappa(\gamma k)M) = 1 $ for all $ kM \in \Omega $, this is again equal to
    \[
        \int_{\Gamma \bs \Omega} (f_1(k), f_2(k)) \, dk_M = \int_{\Gamma \bs \Omega} (f_1(x), f_2(x)) \ .
    \]
\end{proof}

For $ f_1 \in \C^{-\infty}(B, V_B(\sigma_\mu, \phi)) $, $ f_2 \in \Cinf(B, V_B(\sigma_{-\bar{\mu}}, \phi)) $, we define $ (f_1, f_2)_B $ by
\[
    \langle f_1, I(f_2) \rangle := f_1(I(f_2)) \ ,
\]
where $ I \colon \Cinf(B, V_B(\sigma_{-\bar{\mu}}, \phi)) \to \Cinf(B, V_B(\tilde{\sigma}_{-\mu}, \phi)) $ is the canonical conjugate-linear isomorphism.
Then, $ (\cdot, \cdot)_B $ is a sesquilinear pairing between
\[
    \C^{-\infty}(B, V_B(\sigma_\mu, \phi)) \quad \text{and} \quad \Cinf(B, V_B(\sigma_{-\bar{\mu}}, \phi)) \ .
\]
We define $ (f_2, f_1)_B = \overline{(f_1, f_2)_B} $.
Note that if the pairing between two distributions $ f_1 $ and $ f_2 $ is well-defined for two of the three pairings, then they give the same result.
By abuse of notation, we write always $ (f_1, f_2)_B = \int_{\Gamma \bs \Omega}{ (f_1(x), f_2(x)) } $.

By Corollary 8.7 of \cite[p.132]{BO00} (see Section \ref{ssec:SomedefinitionsAndKnownResults} for details),
\[
    (V_{\pi, -\infty} \otimes V_\phi)_d^\Gamma \ \hat{\otimes} \ V_{\pi'}
    \hookrightarrow L^2(\Gamma \bs G, \phi), \quad f \otimes v \mapsto c_{f, v}
\]
induces an isometric embedding of $ (V_{\pi, -\infty} \otimes V_\phi)_d^\Gamma \ \hat{\otimes} \ V_{\pi'} $ into $ L^2(\Gamma \bs G, \phi) $. Thus,
\[
    (f_1, f_2) (v_1, v_2) =
    (f_1 \otimes v_1, f_2 \otimes v_2)
    = (c_{f_1, v_1}, c_{f_2, v_2})_{L^2(\Gamma \bs G, \phi)}
\]
for all $ f_1, f_2 \in (V_{\pi, -\infty} \otimes V_\phi)_d^\Gamma $ and $ v_1, v_2 \in V_{\pi'} $. Let $ \zeta \in V_{\pi', K} $ be of norm 1. Then,
\[
    (f_1, f_2) = (c_{f_1, \zeta}, c_{f_2, \zeta})_{L^2(\Gamma \bs G, \phi)}
    \qquad (f_1, f_2 \in (V_{\pi, -\infty} \otimes V_\phi)_d^\Gamma)
\]
induces a sesquilinear pairing $ (\cdot, \cdot) $ between the tempered (thus all, by Remark \ref{Rem discrete series rep and tempered rep}) $ \Gamma $-invariant distribution vectors and the Schwartz vectors:
\[
    (f_1, f_2) := (c_{f_1, \zeta}, c_{f_2, \zeta})
    \qquad (f_1 \in (V_{\pi, -\infty} \otimes V_\phi)^\Gamma, f_2 \in (V_{\pi, -\infty} \otimes V_\phi)_\CS^\Gamma) \ ,
\]
where the pairing on the right-hand side is simply the pairing between $ \CS'(\Gamma \bs G, \phi) $ and $ \CS(\Gamma \bs G, \phi) $.

This new sesquilinear pairing does not depend on the choice of $ \zeta $. Indeed, if $ f_1 $ and $ f_2 $ are fixed, then $ F(v_1, v_2) := c_{f_1, v_1}(c_{f_2, v_2}) $ is a $ (\g, K) $-invariant sesquilinear form on $ V_{\pi', K} $. It follows from Schur's lemma that there exists a constant $ C(f_1, f_2) $, depending on $ f_1 $ and $ f_2 $, such that $ c_{f_1, v_1}(c_{f_2, v_2}) = C(f_1, f_2) (v_1, v_2) $ for all $ v_1, v_2 \in V_{\pi', K} $. Thus, $ c_{f_1, v}(c_{f_2, v}) = C(f_1, f_2) $ for all $ v \in V_{\pi', K} $ of norm 1.
For $ f_1 \in (V_{\pi, -\infty} \otimes V_\phi)^\Gamma $ and $ f_2 \in (V_{\pi, -\infty} \otimes V_\phi)_\CS^\Gamma $, define
\[
    (f_2, f_1) = \overline{(f_1, f_2)} \ .
\]
Note again that there is no ambiguity.

Fix $ \zeta \in V_{\pi', K} \cap \im \tilde{t} $ of norm 1. Fix $ 0 \neq u \in V_{\tilde{\gamma}} $ such that $ \tilde{t}(u) = \zeta $.
Then, $ \t{q}(\zeta) = u_{T, -\lambda} $.

It follows from the definition of $ (\cdot, \cdot) $, the topologies on $ (V_{\pi, -\infty} \otimes V_\phi)^\Gamma $ and $ (V_{\pi, -\infty} \otimes V_\phi)_\CS^\Gamma $ and Remark \ref{Rem discrete series rep and tempered rep} \eqref{Rem discrete series rep and tempered rep, Lemma 9.4} that $ (f_1, f_2) $ depends continuously on $ f_1 $ (resp. $ f_2 $) when $ f_2 $ (resp. $ f_1 $) is fixed.

\begin{sloppypar}
    If $ f \in \C^{-\infty}(\dX, V(\sigma_\lambda, \phi)) $ with $ Af \in \C^{-\infty}(\dX, V(\sigma_{-\lambda}, \phi))^\Gamma $ and $ g \in \C^{-\infty}(\dX, V(\sigma_{-\lambda}, \phi))^\Gamma \cap \im(\beta) $, then, as $ \beta $ is injective, set
    \begin{equation}
        (Af, g) = (q(f), \beta^{-1}(g))
    \end{equation}
    whenever the right-hand side is well-defined.
\end{sloppypar}

\begin{lem}\label{Lem ext regular}
    The Bunke-Olbrich extension map
    \[
        \ext \colon \C^{-\infty}(B, V_B(\sigma_\lambda, \phi)) \to \C^{-\infty}(\dX, V(\sigma_\lambda, \phi))^\Gamma
    \]
    is regular at $ \lambda $.
\end{lem}

\begin{proof}
    We must show that $ \C^{-\infty}(\dX, V(\sigma_\lambda)) $ is reducible. 
    Then, $ \ext $ is regular at $ \lambda $ by Proposition 4.21 of \cite[p.44]{Olb02}.
    \\ It follows from Casselman's subrepresentation theorem and the functorial properties of smooth globalisation that $ V_{\pi, -\infty} \otimes V_\phi $ is a quotient of $ \C^{-\infty}(\dX, V(\sigma_{-\lambda}, \phi)) $ and a nontrivial $G$-submodule of $ \C^{-\infty}(\dX, V(\sigma_{-\lambda}, \phi)) $ (By a $G$-submodule we always mean a representation on a closed subspace.). As moreover the Langlands quotient is contained in $ \ker q $, the latter space is a nontrivial $G$-invariant submodule of $ \C^{-\infty}(\dX, V(\sigma_\lambda, \phi)) $. The lemma follows.
\end{proof}

\begin{lem}\label{Lem C^(-infty)(Lambda, V(sigma_lambda, phi))^Gamma is contained in ker q}
    $ \C^{-\infty}(\Lambda, V(\sigma_\lambda, \phi))^\Gamma $ is contained in $ \ker q $.
\end{lem}

\begin{proof}
    $ \C^{-\infty}(\Lambda, V(\sigma_\lambda, \phi))^\Gamma $ is contained in the Langlands quotient by Proposition 4.21 of \cite[p.44]{Olb02}, which is again contained in $ \ker(q) $.
\end{proof}

\begin{thm}\label{Thm Collingwood result} \nlenum
    \begin{enumerate}
    \item The space $ A(\C^{-\infty}_\Omega(\dX, V(\sigma_\lambda, \phi))^\Gamma) \, (\text{resp. } A(\C^{-\infty}(\dX, V(\sigma_\lambda, \phi))^\Gamma) ) $ has finite codimension in
        \[
            \im(A) \cap \C^{-\infty}_\Omega(\dX, V(\sigma_{-\lambda}, \phi))^\Gamma \quad (\text{resp. } \im(A) \cap \C^{-\infty}(\dX, V(\sigma_{-\lambda}, \phi))^\Gamma ) \ .
        \]
    \item The space $ \im(A) \cap \C^{-\infty}_\Omega(\dX, V(\sigma_{-\lambda}, \phi))^\Gamma $ is dense in
        \[
            \im(A) \cap \C^{-\infty}(\dX, V(\sigma_{-\lambda}, \phi))^\Gamma \ .
        \]
    \end{enumerate}
\end{thm}

\begin{proof}
    This is the content of Theorem \ref{Thm V_Omega subset V_(-infty)^Gamma is dense} in the appendix.
\end{proof}

\begin{thm}\label{Thm (Af, g) = C (res(f), res(g))} 
    There exists a real constant $ C \neq 0 $ such that
    \[
        (Af, g) = C (\rest(f), \rest(g))_B \in \CC
    \]
    for all $ f \in \C^{-\infty}(\dX, V(\sigma_\lambda, \phi))^\Gamma $ and 
    $ g \in \C^{-\infty}(\dX, V(\sigma_{-\lambda}, \phi))^\Gamma \cap \im(\beta) $ such that $f$ or $g$ is smooth on $ \Omega $.
\end{thm}

\begin{rem}
    $ (Af, g) $ is well-defined by Proposition \ref{Prop Af, smooth on Omega, invariant distributions are Schwartz} and as $A$ is off-diagonally smoothing by Lemma \ref{Lem A off-diagonally smoothing}.
\end{rem}

\begin{proof}
    The argument is adapted from Proposition 10.4 of \cite[p.144]{BO00}.
    Let $ f \in \C^{-\infty}_\Omega(\dX, V(\sigma_\lambda, \phi))^\Gamma $ and $ g = A \tilde{g} \in \C^{-\infty}_\Omega(\dX, V(\sigma_{-\lambda}, \phi))^\Gamma \cap \im(\beta) $.
    Since $ \ext $ is regular at $ \lambda $ by Lemma \ref{Lem ext regular}, there is a holomorphic family $ \mu \mapsto f_\mu \in \C^{-\infty}_\Omega(\dX, V(\sigma_\mu, \phi))^\Gamma $ defined in a neighbourhood of $ \lambda $ such that $ \rest(f) = \rest(f_\lambda) $.
    Since $ \C^{-\infty}(\Lambda, V(\sigma_\lambda, \phi))^\Gamma $ is contained in $ \ker q $ by Lemma \ref{Lem C^(-infty)(Lambda, V(sigma_lambda, phi))^Gamma is contained in ker q}, we may assume without loss of generality that $ f = f_\lambda $.

    Let $ B_R $ be the ball in $X$ with center $ eK $ and radius $R$.
    Let $ \chi \in \Ccinf(X \cup \Omega) $ be as in Lemma \ref{Lem Cut-off function in the convex-cocompact case}.

    Let us denote the scalar product in $ L^2(B_R, V(\gamma, \phi)) $ by $ (\cdot, \cdot)_{B_R} $.
    Set
    \[
        S_R(\mu) = (P^T_\mu f_\mu, \chi P^T_\lambda \tilde{g})_{B_R} \ .
    \]
    Then, $ P^T_\lambda f $ and $ P^T_\lambda \tilde{g} $ belongs to $ L^2(Y, V_Y(\gamma, \phi)) $ by Lemma 6.2, 3., of \cite[p.110]{BO00} (asymptotic expansions for the Poisson transform, compare with Lemma \ref{Lemma2 about Asymptotic expansions}) and by Corollary \ref{Cor of Lem Asymptotics of c (Af, v)(kah) for -lambda}. Thus,
    \[
        (P^T_\lambda f, P^T_\lambda \tilde{g})_{L^2(Y, V_Y(\gamma, \phi))} = \lim_{R \to \infty} \lim_{\mu \to \lambda} S_R(\mu) \ .
    \]
    By \cite{BO95}, there exists $ c(\sigma) \in \RR $ such that $ (-\Omega_G + c(\sigma) + \mu^2) \circ P_\mu^S = 0 $ for all $ S \in \Hom_M(V_\sigma, V_\gamma) $ and all $ \mu \in \a^*_\CC $.
    Let $ D = -\Omega_G + c(\sigma) + \lambda^2 $ be the shifted Casimir operator.
    Then, $ D P_\mu^T f_\mu = (\lambda^2 - \mu^2) P_\mu^T f_\mu $. So,
    \[
        S_R(\mu) = \frac1{\lambda^2 - \mu^2} (D P^T_\mu f_\mu, \chi P^T_\lambda \tilde{g})_{B_R}
    \]
    for all $ \mu \in \a^*_\CC $ such that $ \R(\mu) > 0 $ and $ \mu \neq \lambda $. 
    Observe that there exists a selfadjoint endomorphism $ \Rcal $ such that $ D = \nabla^* \nabla + \Rcal $, where $ \nabla^* \nabla $ is the Bochner Laplacian associated to the invariant connection $ \nabla $ of $ V(\gamma, \phi) $. As $ \mu \mapsto P^T_\mu f_\mu $ is holomorphic, we can write
    \[
        P^T_\mu f_\mu = F_0 + F_1(\mu - \lambda) + \dotsb \,
    \]
    with $ F_i \in \Cinf(X, V(\gamma, \phi))^\Gamma $. 
    We have:
    \begin{multline*}
        D F_0 + D F_1(\mu - \lambda) + \dotsb = D P^T_\mu f_\mu  \\
        = (\lambda^2 - \mu^2) P^T_\mu f_\mu
        = -(\mu + \lambda) (F_0(\mu - \lambda) + F_1(\mu - \lambda)^2 + \dotsb) \ .
    \end{multline*}
    Thus, $ D F_0 = 0 $ and $ \lim_{\mu \to \lambda} D\big(\frac{P^T_\mu f_\mu}{\mu - \lambda}\big) = D F_1 $.
    Hence, $ \lim_{\mu \to \lambda} S_R(\mu) $ is equal to
    \begin{multline*}
        \frac1{2\lambda} (DF_1, \chi P^T_\lambda \tilde{g})_{B_R}
        = \frac1{2\lambda} \Big( (DF_1, \chi P^T_\lambda \tilde{g})_{B_R}
        - (F_1, \chi D P^T_\lambda \tilde{g})_{B_R} \Big) \\
        = \frac1{2\lambda} \Big( (DF_1, \chi P^T_\lambda \tilde{g})_{B_R}
        - (F_1, D \chi P^T_\lambda \tilde{g})_{B_R}
        + (F_1, [D, \chi] P^T_\lambda \tilde{g})_{B_R} \Big) \ .
    \end{multline*}
    By applying Green's formula, this yields
    \[
        \frac1{2\lambda} \Big( (\nabla_n F_1, \chi P^T_\lambda \tilde{g})_{\del B_R}
        - (F_1, \nabla_n \chi P^T_\lambda \tilde{g})_{\del B_R}
        + (F_1, [D, \chi] P^T_\lambda \tilde{g})_{B_R} \Big) \ ,
    \]
    where $ n $ is the exterior unit normal vector field at $ \del B_R $.
    By Lemma 6.2, 3., of \cite[p.110]{BO00} (compare with Lemma \ref{Lemma2 about Asymptotic expansions}), there exists $ \eps > 0 $ such that for $ a \to \infty $ we have
    \[
        P^T_\mu f_\mu(ka) = a^{\mu - \rho} c_\gamma(\mu) T f(k) + O(a^{\R(\mu) -\rho - \eps})
    \]
    uniformly as $kM $ varies in compact subsets of $\Omega$.
    Here, $ c_\gamma(\mu) \in \End_M(V_\gamma) $ is the Harish-Chandra $c$-function. This meromorphic function is given by
    \[
        c_\gamma(\mu) := \int_{\bar{N}}{ a(\bar{n})^{-(\mu+\rho)} \gamma(\kappa(\bar{n})) \, d\bar{n} }
    \]
    for $ \R(\mu) > 0 $.
    Since the above formula also holds for $ \mu = \lambda $ and any $ f \in \C^{-\infty}_\Omega(\dX, V(\sigma_\lambda, \phi))^\Gamma $, it follows from Corollary \ref{Cor of Lem Asymptotics of c (Af, v)(kah) for -lambda} that $ c_\gamma(\lambda) = 0 $.
    As
    \[
        a^{\mu - \rho} = a^{\mu - \lambda} a^{\lambda - \rho} = a^{\lambda - \rho} e^{(\mu - \lambda) \log(a)} = a^{\lambda - \rho} \sum_{j=0}^\infty{ \frac{\log(a)^j}{j!} (\mu - \lambda)^j }
    \]
    and as
    \[
        c_\gamma(\mu) = \sum_{j=0}^\infty (-1)^j (\mu - \lambda)^j
        \int_{\bar{N}}{ \frac{\log(a(\bar{n}))^j}{j!} a(\bar{n})^{-\lambda-\rho} \gamma(\kappa(\bar{n})) \, d\bar{n} } =: \sum_{j=0}^\infty (\mu - \lambda)^j c_\gamma^j(\lambda) \ ,
    \]
    we have for $ a \to \infty $
    \[
        P^T_\mu f_\mu(ka) = a^{\lambda - \rho} c_\gamma(\lambda) + (\mu - \lambda)(a^{\lambda - \rho} (\log(a) c_\gamma(\lambda) + c_\gamma^1(\lambda)) T f(k)) + \dotsb + O(a^{\lambda -\rho - \eps})
    \]
    uniformly as $ kM $ varies in compact subsets of $ \Omega $, for all $ \mu \in \a^*_\CC $ such that $ \R(\mu) < \lambda $. 
    Assume from now on that $ \mu \in \a^*_\CC $ satisfies $ 0 < \R(\mu) < \lambda $.

    By Cauchy's integral formula and as $ c_\gamma(\lambda) = 0 $, we have for $ a \to \infty $
    \begin{equation}\label{Proof Lem (Af, g) = C <res(f), res(g)> eq1}
        F_1(ka) = a^{\lambda - \rho} c_\gamma^1(\lambda) T f(k)
        + O(a^{\lambda -\rho - \eps})
    \end{equation}
    uniformly as $ kM $ varies in compact subsets of $ \Omega $, as $ c_\gamma(\lambda) = 0 $.
    \\ Let $ \omega_X = \frac{\omega_n}{2^r} $, where $ n = \dim X $, $ \omega_n = \vol(S^{n-1}) = \frac{2\pi^{\frac{n}2}}{\Gamma(\frac{n}2)} $ and $ r \in \NN $ is such that $ r \alpha = 2\rho $.
    It follows from \eqref{Proof Lem (Af, g) = C <res(f), res(g)> eq1} and Corollary \ref{Cor of Lem Asymptotics of c (Af, v)(kah) for -lambda} combined with properties \eqref{Lem Cut-off function in the convex-cocompact case eq2} and \eqref{Lem Cut-off function in the convex-cocompact case eq3} of Lemma \ref{Lem Cut-off function in the convex-cocompact case} that $ | (F_1, [D, \chi] P^T_\lambda \tilde{g})| $ is integrable over $X$. Thus, $ (F_1, [D, \chi] P^T_\lambda \tilde{g})_X $ is equal to
    \[
        \int_{\Gamma \bs X} \sum_{\gamma \in \Gamma} (F_1(\gamma x), D \big(\chi(\gamma x) P^T_\lambda \tilde{g}(\gamma x)\big) ) \, dx
        = \int_{\Gamma \bs X} (F_1(x), D\big(\sum_{\gamma \in \Gamma} \chi(\gamma x) P^T_\lambda \tilde{g}(x)\big) ) \, dx
    \]
    as $ \sum_{\gamma \in \Gamma} \chi(\gamma x) = 1 $, as $ D \circ P^T_\lambda = 0 $ and as $ F_1 $ and $ P^T_\lambda \tilde{g} $ belong to $ \Cinf(X, V(\gamma, \phi))^\Gamma $. Hence,
    \begin{equation}\label{Proof Lem (Af, g) = C <res(f), res(g)> eq2}
        (F_1, [D, \chi] P^T_\lambda \tilde{g})_X = 0 \ .
    \end{equation}
    Recall that $ \chi_\infty $ denotes the restriction of $ \chi $ to $ \Omega $. Then, $ \chi_\infty(k) = \lim_{a \to \infty} \chi(kaK) $ as $ \chi $ is continuous.
    As $ \chi $ has compact support in $ X \cup \Omega $, $ \chi_\infty $ has compact support in $ \Omega $.
    So, by \eqref{Proof Lem (Af, g) = C <res(f), res(g)> eq1}, \eqref{Proof Lem (Af, g) = C <res(f), res(g)> eq2} and by Lemma \ref{Lem Averaged asymptotics of P^T_(-lambda) g(ka)}, $ (P^T_\lambda f, P^T_\lambda \tilde{g})_{L^2(Y, V_Y(\gamma, \phi))} $ is equal to
    \begin{multline*}
        \frac{\omega_X}{2\lambda} \Big( (\lambda - \rho)
        \int_{\dX}{ (c_\gamma^1(\lambda) T f(k), \chi_\infty(k) (\tilde{\beta} (q\tilde{g}))(k)) \, dk} \\
        - (-\lambda - \rho) \int_{\dX}{ (c_\gamma^1(\lambda) T f(k), \chi_\infty(k) (\tilde{\beta} (q\tilde{g}))(k)) \, dk} \Big)
    \end{multline*}
    as $ \lim_{R \to \infty} e^{-2\rho R} \vol(\del B_R) = \omega_X $.
    This is again equal to
    \[
        \frac{ \omega_X }{ (t, t) } \int_{\dX}{ \chi_\infty(k) ((t \circ c_\gamma^1(\lambda) T) f(k), g(k)) \, dk} \ ,
    \]
    where $ (T_1, T_2) \Id_{V_\sigma} = T_2^* \circ T_1 $ (since $ \sigma \in \hat{M} $ this is well-defined by Schur's lemma). By Lemma \ref{Lem int_B (f_1(x), f_2(x)) = int_K chi_infty(kM) (f_1(k), f_2(k)) dk}, this yields
    \[
        \frac{ (c_\gamma^1(\lambda) T, t^*) \omega_X }{ (t, t) } (\rest(f), \rest(g))_B \ .
    \]
    We have just shown that
    \begin{equation}\label{Proof Lem (Af, g) = C <res(f), res(g)> eq4}
        (P^T_\lambda f, P^T_\lambda \tilde{g})_{L^2(Y, V_Y(\gamma, \phi))} = \frac{ (c_\gamma^1(\lambda) T, t^*) \omega_X }{ (t, t) } (\rest(f), \rest(g))_B \ .
    \end{equation}
    Thus, by Lemma \ref{Lem BO00, 9.1}, we have
    \[
        (c_{f, u_{T, -\lambda}}, c_{\tilde{g}, u_{T, -\lambda}})_{L^2(\Gamma \bs G)}
        = \int_{\Gamma \bs G} (\langle P^T_\lambda f(x), u \rangle, \langle P^T_\lambda \tilde{g}(x), u \rangle) \, dx \ .
    \]
    By Fubini's theorem, this yields
    \begin{align*}
        & \int_{\Gamma \bs G} \int_K \langle P^T_\lambda f(xk), u \rangle \overline{\langle P^T_\lambda \tilde{g}(xk), u \rangle} \, dk \, dx \\
        =& \int_{\Gamma \bs G} \int_K \langle \gamma(k^{-1}) P^T_\lambda f(x), u \rangle \overline{\langle \gamma(k^{-1}) P^T_\lambda \tilde{g}(x), u \rangle} \, dk \, dx \ .
    \end{align*}
    By Schur's lemma, this is again equal to
    \[
        \frac{\|u\|^2}{\dim(V_{\tilde{\gamma}})} \int_{\Gamma \bs G} (P^T_\lambda f(x), P^T_\lambda \tilde{g}(x)) \, dx
        = \frac{ \|u\|^2 (c_\gamma^1(\lambda) T, t^*) \omega_X }{ (t, t) \dim(V_{\tilde{\gamma}}) } \cdot (\rest(f), \rest(g))_B \ .
    \]
    Thus, by \eqref{Proof Lem (Af, g) = C <res(f), res(g)> eq4}, $ (Af, g) $ is equal to
    \[
        C (\rest(f), \rest(g))_B \ ,
    \]
    where
    \[
        C := \frac{ \|u\|^2 (c_\gamma^1(\lambda) T, t^*) \omega_X }{ (t, t) \dim(V_{\tilde{\gamma}}) }
        = \frac{ (c_\gamma^1(\lambda) T, t^*) \omega_X }{ (t, t) (\tilde{t}, \tilde{t}) \dim(V_{\tilde{\gamma}}) }
    \]
    as $ \tilde{t}(u) = \xi $ and as $ \xi $ has norm 1.
    Since $ (\cdot, \cdot) $ is positive definite,
    \[
        0 < (Af, Af) = C (\rest(f), \rest(Af))_B \quad (f \in \C^{-\infty}_\Omega(\dX, V(\sigma_\lambda, \phi))^\Gamma \smallsetminus \ker A)
    \]
    shows that $ C \in \RR \smallsetminus \{0\} $.

    This proves the formula in the case where $f$ and $g$ are smooth on $ \Omega $. Let us show now that we can extend this formula.

    Since $ \beta^{-1} \colon \im(A) \to V_{\pi, -\infty} \otimes V_\phi $ is continuous ($ \beta $ is a topological embedding),
    \[
        \{ v \in (V_{\pi, -\infty} \otimes V_\phi)^\Gamma \mid \beta v \in \C^{-\infty}_\Omega(\dX, V(\sigma_{-\lambda}, \phi))^\Gamma \cap \im(A) \}
    \]
    is dense in $ (V_{\pi, -\infty} \otimes V_\phi)^\Gamma $ by Theorem \ref{Thm Collingwood result}.

    Let $ g = \beta(\tilde{g}) \in \C^{-\infty}(\dX, V(\sigma_{-\lambda}, \phi))^\Gamma \cap \im(A) $. Then, there exist $ \tilde{g}_j \in (V_{\pi, -\infty} \otimes V_\phi)^\Gamma $ with $ \beta(\tilde{g}_j) \in \C^{-\infty}_\Omega(\dX, V(\sigma_{-\lambda}, \phi))^\Gamma \cap \im(A) $ such that $ \tilde{g}_j $ converges to $ \tilde{g} $ in $ (V_{\pi, -\infty} \otimes V_\phi)^\Gamma $.
    Thus,
    \begin{multline*}
        (Af, g) = (qf, \tilde{g}) = \lim_{j \to \infty} (qf, \tilde{g}_j) = \lim_{j \to \infty} (Af, \beta(\tilde{g}_j)) \\
        = C \lim_{j \to \infty} (\rest(f), \rest(\beta(\tilde{g}_j)))_B = C (\rest(f), \rest(g))_B \ .
    \end{multline*}
    \begin{sloppypar}
        Since $ \C^{-\infty}_\Omega(\dX, V(\sigma_\lambda, \phi))^\Gamma $ is dense in $ \C^{-\infty}(\dX, V(\sigma_\lambda, \phi))^\Gamma $ by Lemma 6.6 of \cite[p.114]{BO00}, one can show similarly that the formula holds when $ f \in \C^{-\infty}(\dX, V(\sigma_\lambda, \phi))^\Gamma $ and $ g \in \C^{-\infty}_\Omega(\dX, V(\sigma_{-\lambda}, \phi))^\Gamma \cap \im(A) $.
    \end{sloppypar}
    This completes the proof of the theorem.
\end{proof}

\begin{prop}\label{Prop Schwartz implies Af smooth on Omega}
    If $ q f $ $ (f \in \C^{-\infty}(\dX, V(\sigma_\lambda, \phi)) )$ is a Schwartz vector, then $ Af \in \C^{-\infty}_\Omega(\dX, V(\sigma_{-\lambda}, \phi))^\Gamma $.
\end{prop}

\begin{proof}
    Let $ f $ be as above. By Theorem \ref{Thm (Af, g) = C (res(f), res(g))}, there exists $ C \neq 0 $ such that
    \[
        (q(\ext \varphi), qf) = (A(\ext \varphi), Af) = C (\varphi, \rest(Af))_B
    \]
    for all $ \varphi \in \Cinf(B, V_B(\sigma_\lambda, \phi)) $.
    Recall that $ (V_{\pi, -\infty} \otimes V_\phi)^\Gamma \ni f \mapsto c_{f, v} \in \CS'(\Gamma \bs G, \phi) $ is continuous by Remark \ref{Rem discrete series rep and tempered rep} \eqref{Rem discrete series rep and tempered rep, Lemma 9.4}.
    \\ Thus, the left-hand side is also well-defined for $ \varphi \in \C^{-\infty}(B, V_B(\sigma_\lambda, \phi)) $ and depends continuously on $ \varphi $. Hence, we can define a continuous conjugate-linear map
    \[
        g \colon \C^{-\infty}(B, V_B(\sigma_\lambda, \phi)) \to \CC
    \]
    by $ g(\varphi) = \frac1C \overline{(q(\ext \varphi), qf)} $. Since $B$ is compact, we can identify $ \C^{-\infty}(B, V_B(\sigma_\lambda, \phi))' $ with $ \Cinf(B, V_B(\sigma_{-\lambda}, \phi)) $. So, $ g $ is given by a smooth function. By construction, $ g = \rest(Af) $ as distributions. Consequently, $ \rest(Af) $ is smooth.
\end{proof}

\needspace{2\baselineskip}
For the definition of the maps $ \F_\pi $ ($ \pi \in \hat{G} $), see Section \ref{ssec:SomedefinitionsAndKnownResults}.

\begin{prop}\label{Prop c_pi(F_pi(f)) for discrete series} \nlenum
    \begin{enumerate}
    \item Let $ \gamma \in \hat{K} $. If $ f \in \CS(\Gamma \bs G, \phi)(\gamma) $, then $ \F_\pi(f) \in (V_{\pi', -\infty} \otimes V_\phi)^\Gamma_{\c{\CS}} \otimes V_{\pi}(\gamma) $ and $ c_\pi(\F_\pi(f)) \in \c{\CS}(\Gamma \bs G, \phi)(\gamma) $.
    \item If $ f \in \CS(\Gamma \bs G, \phi)_K $ ($K$-finite element), then $ c_\pi(\F_\pi(f)) \in \c{\CS}(\Gamma \bs G, \phi)_K $.
    \end{enumerate}
\end{prop}

\begin{rem} \nlenum
    \begin{enumerate}
    \item $ c_{\pi} \circ \F_{\pi} $ is the orthogonal projection on $ (V_{\pi', -\infty} \otimes V_\phi)_d^\Gamma \hat{\otimes} V_{\pi} $ by Corollary \ref{Cor c_pi o F_pi is the orth. proj. on V_(pi', -infty)^Gamma otimes V_pi}.
    \item It follows from this proposition and Lemma \ref{Lem (c_pi o F_pi)^2 = c_pi o F_pi, F_pi o c_pi = Id on V_(pi', -infty)^Gamma otimes V_pi} that
        \[
            \F_\pi(\CS(\Gamma \bs G, \phi)(\gamma)) \qquad (\gamma \in \hat{K})
        \]
        is equal to $ (V_{\pi', -\infty} \otimes V_\phi)^\Gamma_\CS \otimes V_{\pi}(\gamma) $ and that $ (V_{\pi', -\infty} \otimes V_\phi)^\Gamma_\CS \otimes V_{\pi, K} $ is equal to $ \F_{\pi}(\CS(\Gamma \bs G, \phi)_K) $.
    \end{enumerate}
\end{rem}

\begin{proof}
    Let $ \pi $ be as above and $ f \in \CS(\Gamma \bs G, \phi)_K $. The construction of the embedding map $ \beta_{\pi'} $ provides $ \sigma \in \hat{M} $ and $ \lambda > 0 $.
    We may assume without loss of generality that $ f $ is nonzero and belongs to a $K$-isotypic component $ L^2(\Gamma \bs G, \phi)(\gamma) $ for some $ \gamma \in \hat{K} $.
    Let $ \{v_1, \dotsc, v_l\} $ be an orthonormal basis of $ V_{\pi}(\gamma) $.
    Write $ \F_\pi(f) = \sum_{j=1}^l \tilde{f}_j \otimes v_j $ with $ \tilde{f}_j \in N_\pi $.
    Let $ \varphi \in \Cinf(B, V_B(\sigma_{\lambda}, \phi)) $. Then, for all $ i \in \{1, \dotsc, l\} $, we have
    \begin{multline*}
        (c_\pi(q_{\pi'}(\ext(\varphi)) \otimes v_i), f)
        = (q_{\pi'}(\ext(\varphi)) \otimes v_i, \F_\pi(f))_{N_\pi \hat{\otimes} V_\pi} \\
        = \sum_{j=1}^l (q_{\pi'}(\ext(\varphi)), \tilde{f}_j) (v_i, v_j)
        = (q_{\pi'}(\ext(\varphi)), \tilde{f}_i) (v_i, v_i) \ .
    \end{multline*}
    By Theorem \ref{Thm (Af, g) = C (res(f), res(g))} and as $ (v_i, v_i) = 1 $, there exists a constant $ C \neq 0 $ such that the last part is equal to
    \[
        C (\varphi, \rest(\beta_{\pi'} \tilde{f}_i))_B \ .
    \]
    Recall that $ (V_{\pi, -\infty} \otimes V_\phi)^\Gamma \ni f \mapsto c_{f, v} \in \CS'(\Gamma \bs G, \phi) $ is continuous by Remark \ref{Rem discrete series rep and tempered rep} \eqref{Rem discrete series rep and tempered rep, Lemma 9.4}. Thus, one can show similarly as in the proof of Proposition \ref{Prop Schwartz implies Af smooth on Omega} that $ \beta_{\pi'}(\tilde{f}_i) $ is smooth on $ \Omega $ for all $i$.
    Hence, $ \F_\pi(f) \in (V_{\pi', -\infty} \otimes V_\phi)^\Gamma_{\c{\CS}} \otimes V_{\pi}(\gamma) $ by Proposition \ref{Prop Af, smooth on Omega, invariant distributions are Schwartz}. So, $ c_\pi(\F_\pi(f)) \in \c{\CS}(\Gamma \bs G, \phi) $ .
\end{proof}

\begin{prop}\label{Prop Schwartz vectors are dense in the space of Gamma-invariant vectors}
    The space $ (V_{\pi, -\infty} \otimes V_\phi)^\Gamma_{\c{\CS}} $ is dense in $ (V_{\pi, -\infty} \otimes V_\phi)^\Gamma $.
\end{prop}

\begin{proof}
    It follows from Proposition \ref{Prop Af, smooth on Omega, invariant distributions are Schwartz}, Proposition \ref{Prop Schwartz implies Af smooth on Omega} and Proposition \ref{Prop Every Schwartz vector is cuspidal} that $ \beta\big((V_{\pi, -\infty} \otimes V_\phi)^\Gamma_{\c{\CS}}\big) $ is equal to $ \im(A) \cap \C^{-\infty}_\Omega(\dX, V(\sigma_{-\lambda}, \phi))^\Gamma $.

    The proposition follows now immediately from Theorem \ref{Thm Collingwood result} and the fact that
    \[
        \beta \colon (V_{\pi, -\infty} \otimes V_\phi)^\Gamma \to \im(A) \cap \C^{-\infty}(\dX, V(\sigma_{-\lambda}, \phi))^\Gamma
    \]
    is a homeomorphism by Lemma \ref{Lem beta is a top embedding}. 
\end{proof}

\newpage

\subsubsection{The contribution of the discrete series representations to the closure of the space of cusp forms}

In this section, we show that cuspidal Schwartz vectors are dense in the multiplicity space of the discrete series representation $ \pi $. For this, we compute first the annihilator of certain subspaces of the space of Schwartz vectors. We will reuse these results in the following sections.

\medskip


From now on we identify $ V_\pi $ (resp. $ V_{\pi'} $) and $ \beta(V_\pi) $ (resp. $ \t{q}(V_{\pi'}) $), equipped with the topology induced by $ V_\pi $, as topological spaces.
So, we view $ \beta $ and $ \t{q} $ as the identity map and $ V_{\pi, -\infty} $ as a subspace of $ \C^{-\infty}(\dX, V(\sigma_{-\lambda})) $.

Let\index[n]{D1 D2@$ D_1 $, $ D_2 $}
\[
    D_1 = A(\C^{-\infty}_\Omega(\dX, V(\sigma_\lambda, \phi))^\Gamma)
\]
and
\[
    D_2 = \C^{-\infty}(\Lambda, V(\sigma_{-\lambda}, \phi))^\Gamma \cap (V_{\pi, -\infty} \otimes V_\phi)^\Gamma \ .
\]
Then, $ D_1 $ and $ D_2 $ are contained in $ (V_{\pi, -\infty} \otimes V_\phi)^\Gamma \cap \C^{-\infty}_\Omega(\dX, V(\sigma_{-\lambda}, \phi))^\Gamma $, consisting of cuspidal Schwartz vectors by Proposition \ref{Prop Af, smooth on Omega, invariant distributions are Schwartz}. We endow $ D_1 $ and $ D_2 $ with the subspace topology induced from $ \C^{-\infty}_\Omega(\dX, V(\sigma_{-\lambda}, \phi))^\Gamma $.


It follows from the proof of Proposition 9.8 of \cite[p.141]{BO00} that $ D_1 $ is infinite-dimensional, contrary to $ D_2 $ which is finite-dimensional by Theorem 6.1 of \cite[p.109]{BO00} if $ X \neq \OO H^2 $. There is also an alternative proof which does not use the meromorphy of $ \ext $ and which is hence also valid in the exceptional case (see \cite[pp.121-122]{BO00}).

Let $ D_{1, -\infty} = A(\C^{-\infty}(\dX, V(\sigma_\lambda, \phi))^\Gamma) \subset (V_{\pi, -\infty} \otimes V_\phi)^\Gamma $.\index[n]{D1-infty@$ D_{1, -\infty} $}

By Theorem \ref{Thm Collingwood result},
\[
    D_1 \subset \im(A) \cap \C^{-\infty}_\Omega(\dX, V(\sigma_{-\lambda}, \phi))^\Gamma
\]
and
\[
    D_{1, -\infty} \subset \im(A) \cap \C^{-\infty}(\dX, V(\sigma_{-\lambda}, \phi))^\Gamma
\]
have finite codimension.

The intersection of $ D_{1, -\infty} $ with $ D_2 $ is trivial by Theorem \ref{Thm (Af, g) = C (res(f), res(g))}.

\begin{dfn} \nlenum\index{annihilator}
    \begin{enumerate}
    \item For a vector subspace $ S $ of $ (V_{\pi, -\infty} \otimes V_\phi)^\Gamma $, we define the \textit{annihilator} of $S$ by
        \[
            \Ann_\CS(S) := \{ v \in (V_{\pi, -\infty} \otimes V_\phi)^\Gamma_\CS \mid (v, v') = 0 \quad \forall v' \in S \} \ .
        \]
    \item For a vector subspace $ S $ of $ (V_{\pi, -\infty} \otimes V_\phi)_\CS^\Gamma $, we define the \textit{annihilator} of $S$ by
        \[
            \Ann(S) := \{ v \in (V_{\pi, -\infty} \otimes V_\phi)^\Gamma \mid (v, v') = 0 \quad \forall v' \in S \} \ .
        \]
    \item For a vector subspace $ S $ of $ \im(A) \cap \C^{-\infty}(\dX, V(\sigma_{-\lambda}, \phi))^\Gamma $, we define the \textit{annihilator} of $S$ by
        \[
            \Ann_\Omega(S) := \{ g \in \im(A) \cap \C^{-\infty}_\Omega(\dX, V(\sigma_{-\lambda}, \phi))^\Gamma \mid (g, h) = 0 \quad \forall h \in S \} \ .
        \]
    \end{enumerate}
\end{dfn}

\begin{rem}\label{Rem Ann(S)} \nlenum
    \begin{enumerate}
    \item We endow $ \Ann_\CS(S) $, $ \Ann(S) $ and $ \Ann_\Omega(S) $ with the subspace topology.
    \item If $ \Ann(S) $ is contained in $ (V_{\pi, -\infty} \otimes V_\phi)_d^\Gamma $, then it is equal to the orthogonal complement of $S$ in $ (V_{\pi, -\infty} \otimes V_\phi)_d^\Gamma $.
    \item Let $ S $ be a vector subspace of $ (V_{\pi, -\infty} \otimes V_\phi)^\Gamma $. Then, $ \Ann_\CS(S) = \Ann_\Omega(S) $ (as vector spaces) by Proposition \ref{Prop Af, smooth on Omega, invariant distributions are Schwartz} and Proposition \ref{Prop Schwartz implies Af smooth on Omega}.
    \end{enumerate}
\end{rem}

\begin{lem}\label{Lem Ann(D_1) = D_2}
    The annihilator of $ D_1 $ is equal to $ D_2 $.
\end{lem}

\begin{rem}\label{Rem Ann(D_1 oplus D_2) = 0} \nlenum
    \begin{enumerate}
    \item It follows that the annihilator of $ D_1 \oplus D_2 $ is zero. Indeed,
        \[
            \Ann(D_1 \oplus D_2) \subset \Ann(D_1) \cap \Ann(D_2) = D_2 \cap \Ann(D_2) = D_2 \cap (D_2)^\perp = \{0\} \ .
        \]
    \item It follows that $ \Ann_\CS(D_1) = D_2 $ and $ \Ann_\Omega(D_1) = D_2 $.
        \\ The full strength of Theorem \ref{Thm (Af, g) = C (res(f), res(g))} is not needed to prove that
        $ \Ann_\Omega(D_1) $ is equal to $ D_2 $.
    \end{enumerate}
\end{rem}

\begin{proof} \nl
    \underline{$ \supset $:}
    For all $ A f \in D_1 $ and $ v \in D_2 $,
    \[
        (A f, v) = C (\rest(f), \rest(v)) = 0
    \]
    by Theorem \ref{Thm (Af, g) = C (res(f), res(g))}. Thus, $ D_2 \subset \Ann(D_1) $.

    \underline{$ \subset $:}
    The annihilator $ \Ann(D_1) $ of $ D_1 $ is equal to
    \[
        \{ v \in (V_{\pi, -\infty} \otimes V_\phi)^\Gamma \mid (A f, v) = 0 \quad \forall f \in \C^{-\infty}_\Omega(\dX, V(\sigma_\lambda, \phi))^\Gamma \} \ .
    \]
    By Theorem \ref{Thm (Af, g) = C (res(f), res(g))}, this space is equal to
    \begin{equation}\label{Proof Prop cuspidal vectors are dense if pi is a DS rep eq1}
        \{ v \in (V_{\pi, -\infty} \otimes V_\phi)^\Gamma \mid (\rest(f), \rest(v))_B = 0 \quad \forall f \in \C^{-\infty}_\Omega(\dX, V(\sigma_\lambda, \phi))^\Gamma \} \ .
    \end{equation}
    By Lemma \ref{Lem ext regular}, $ \ext $ is regular at $ \lambda $. Thus, the continuous map $ \rest $ is surjective. Hence, \eqref{Proof Prop cuspidal vectors are dense if pi is a DS rep eq1} is equal to
    \[
        \{ v \in (V_{\pi, -\infty} \otimes V_\phi)^\Gamma \mid \rest(v) = 0 \}
        = (V_{\pi, -\infty} \otimes V_\phi)^\Gamma \cap \C^{-\infty}(\Lambda, V(\sigma_{-\lambda}, \phi))^\Gamma = D_2 \ .
    \]
    So, the annihilator of $D_1$ is contained in $ D_2 $. Consequently, $ \Ann(D_1) = D_2 $.
\end{proof}

\begin{cor}\label{Cor decomp. of space of cuspidal vectors if pi is a DS rep}
    We have the following topological direct sum:
    \[
        (V_{\pi, -\infty} \otimes V_\phi)_d^\Gamma = \overline{D_1} \oplus D_2 \ .
    \]
    In other words, the cuspidal Schwartz vectors are dense in $ (V_{\pi, -\infty} \otimes V_\phi)_d^\Gamma $ (with respect to the Hilbert space topology on $ (V_{\pi, -\infty} \otimes V_\phi)_d^\Gamma $).
    In particular, the multiplicity space $ N_{\pi'} := (V_{\pi, -\infty} \otimes V_\phi)_d^\Gamma $ of $ \pi' $ is the completion of $ D_1 \oplus D_2 $ with respect to the norm on $ N_{\pi'} $. It can be realised in $ (V_{\pi, -\infty} \otimes V_\phi)^\Gamma $.
\end{cor}

\begin{rem}\nlenum
    \begin{enumerate}
    \item It follows already from the proof of Proposition 9.8 of \cite[p.141]{BO00} that $ D_1 $ is an infinite-dimensional subspace of $ (V_{\pi, -\infty} \otimes V_\phi)_d^\Gamma $.
        This corollary provides us with more precise information about the multiplicity space $ N_{\pi'} $.
    \item The above decomposition depends on the choice of the principal series in which the discrete series is embedded into (choice of $ \sigma \in \hat{M} $).
    \end{enumerate}
\end{rem}

\begin{proof}
    It follows from Lemma \ref{Lem Ann(D_1) = D_2} that $ D_2 $ is the orthogonal complement of $ D_1 $ in $ (V_{\pi, -\infty} \otimes V_\phi)_d^\Gamma $. Thus,
    \[
        (V_{\pi, -\infty} \otimes V_\phi)_d^\Gamma = \overline{D_1} \oplus D_1^\perp = \overline{D_1} \oplus D_2 \ .
    \]
    The last assertion is now obvious as $ D_1 $ and $ D_2 $ are contained in $ (V_{\pi, -\infty} \otimes V_\phi)_{\c{\CS}}^\Gamma $.
\end{proof}

\begin{cor}
    The space $ c_\pi((V_{\pi', -\infty} \otimes V_\phi)^\Gamma_\CS \otimes V_{\pi, K}) $ is dense in $ c_\pi(N_\pi \hat{\otimes} V_\pi) $.
\end{cor}

\begin{proof}
    By Corollary \ref{Cor decomp. of space of cuspidal vectors if pi is a DS rep}, $ (V_{\pi', -\infty} \otimes V_\phi)^\Gamma_\CS $ is dense in $ N_\pi $.

    Since $ c_\pi $ is continuous and since we have the following dense inclusions
    \[
        (V_{\pi', -\infty} \otimes V_\phi)^\Gamma_\CS \otimes V_{\pi, K} \subset N_\pi \otimes V_{\pi, K} \subset N_\pi \hat{\otimes} V_\pi \ ,
    \]
    $ c_\pi((V_{\pi', -\infty} \otimes V_\phi)^\Gamma_\CS \otimes V_{\pi, K}) $ is dense in $ c_\pi(N_\pi \hat{\otimes} V_\pi) $.
\end{proof}

\newpage

\subsubsection[\texorpdfstring{Topological decomposition of $ (V_{\pi, -\infty} \otimes V_\phi)_\CS^\Gamma $ and of \\ \mbox{$ (V_{\pi, -\infty} \otimes V_\phi)^\Gamma $}}{Topological decompositions}]{\texorpdfstring{Topological decomposition of $ (V_{\pi, -\infty} \otimes V_\phi)_\CS^\Gamma $ and of \mbox{$ (V_{\pi, -\infty} \otimes V_\phi)^\Gamma $}}{Topological decompositions}}

In this section, we show that we have the following topological decompositions:
\[
    (V_{\pi, -\infty} \otimes V_\phi)^\Gamma = D_{1, -\infty} \oplus D_2
\]
and
\[
    (V_{\pi, -\infty} \otimes V_\phi)_\CS^\Gamma = D_1 \oplus D_2
\]
(see Theorem \ref{Thm Decompositions of (V_(pi, -infty) otimes V_phi)_CS^Gamma and of (V_(pi, -infty) otimes V_phi)^Gamma}).




\begin{lem}\label{Lem TopologicalDirectSum D_1^(-infty) oplus D_2}
    We have the following topological decompositions:
    \[
        \C^{-\infty}(\dX, V(\sigma_{-\lambda}, \phi))^\Gamma \cap \im(A) = \Ann(D_2) \oplus D_2
    \]
    and
    \[
        \C^{-\infty}_\Omega(\dX, V(\sigma_{-\lambda}, \phi))^\Gamma \cap \im(A) = \Ann_\Omega(D_2) \oplus D_2 \ .
    \]
\end{lem}

\begin{proof}
    Since $ \Ann(D_2) \oplus D_2 \supset D_{1, -\infty} \oplus D_2 $ is closed and finite-codimensional, $ \Ann(D_2) \oplus D_2 $ is a topological direct sum by Lemma \ref{Lem X = Y oplus F topological direct sum if Y is closed and if F is finite-dimensional}.
    \\ Assume that $ \big(\C^{-\infty}(\dX, V(\sigma_{-\lambda}, \phi))^\Gamma \cap \im(A)\big)/(\Ann(D_2) \oplus D_2) $ (endowed with the quotient topology) is nontrivial. It follows from the Hahn-Banach theorem that there exists a nonzero continuous linear functional $ T $ from
    \[
        \bigl(\C^{-\infty}(\dX, V(\sigma_{-\lambda}, \phi))^\Gamma \cap \im(A)\bigr)/\bigl(\Ann(D_2) \oplus D_2\bigr)
    \]
    to $ \CC $. It induces a nonzero functional on $ \C^{-\infty}(\dX, V(\sigma_{-\lambda}, \phi))^\Gamma \cap \im(A) $. This is impossible since the annihilator of $ D_1 \oplus D_2 $ is zero by Remark \ref{Rem Ann(D_1 oplus D_2) = 0}. Thus,
    \[
        \C^{-\infty}(\dX, V(\sigma_{-\lambda}, \phi))^\Gamma \cap \im(A) = \Ann(D_2) \oplus D_2 \ .
    \]
    The second assertion is proven completely similarly.
\end{proof}

Let $ \star $ be $ \infty $ (or $ -\infty $, resp.) and let $ \sharp $ be $ \emptyset $ (or $ \Omega $, resp.).
Fix a trivialisation
\[
    \C^{-\infty}(B, V_B(\sigma_{-\lambda}, \phi)) \hookrightarrow \O_{-\lambda} \C^{-\infty}(B, V_B(\sigma_{.}, \phi)), \quad \varphi \mapsto \varphi_\mu \ .
\]
Recall that here $ \O_{-\lambda} \C^{-\infty}(B, V_B(\sigma_{.}, \phi)) $ is the space of germs at $ -\lambda $ of holomorphic families $ \mu \mapsto \varphi_\mu \in \C^{-\infty}(B, V_B(\sigma_\mu, \phi))) $.
Define $ E_\sharp(\sigma_{-\lambda}, \phi) $ by
\[
     \{ \varphi \in \rest(\C^{-\infty}_\sharp(\dX, V(\sigma_{-\lambda}, \phi))^\Gamma \cap \im(A)) \mid \ext_\mu \varphi_\mu \text{ is regular in } \mu = -\lambda \} \ .
\]
Since $ \ext $ has at most finite-dimensional singularities by Theorem 5.10 of \cite[p.103]{BO00}, $ E_\sharp(\sigma_{-\lambda}, \phi) $ is finite-codimensional in $ \rest(\C^{-\infty}_\sharp(\dX, V(\sigma_{-\lambda}, \phi))^\Gamma \cap \im(A)) $. Let $S_\sharp(\sigma_{-\lambda}, \phi)$ be a finite-dimensional complement.

Define $ \ext(\varphi) = \ext_{-\mu} \restricted{ \varphi_{-\mu} }{ \mu = \lambda } $ for $ \varphi \in E_\sharp(\sigma_{-\lambda}, \phi) $.
Then, $ (\rest \circ \ext)(\varphi) = \varphi $ for all $ \varphi \in E_\sharp(\sigma_{-\lambda}, \phi) $ as $ \rest \circ \ext_{-\mu} = \Id $ by \eqref{eq res o ext = Id}.
We have:
\begin{equation}\label{eq rest(C^(-infty)(dX, V(sigma_mu, phi))^Gamma) = E_mu oplus S_mu}
    \rest(\C^{-\infty}_\sharp(\dX, V(\sigma_{-\lambda}, \phi))^\Gamma \cap \im(A)) = E_\sharp(\sigma_{-\lambda}, \phi) \oplus S_\sharp(\sigma_{-\lambda}, \phi) \ .
\end{equation}

\needspace{2\baselineskip}
The following is the crucial result needed in the proof of Lemma \ref{Lem res induces an embedding}.

\begin{lem}\label{Lem res restricted to E_sharp is an isomorphism}
    The map
    \[
        \restricted{ \ext }{ E_\sharp(\sigma_{-\lambda}, \phi) } \colon E_\sharp(\sigma_{-\lambda}, \phi) \to \ext(E_\sharp(\sigma_{-\lambda}, \phi))
    \]
    is a topological isomorphism with continuous linear inverse $ \rest $.
\end{lem}

\begin{proof}
    Let $ f \in E_\sharp(\sigma_{-\lambda}, \phi) $. Then, $ (\rest \circ \ext)(f) = f $ by \eqref{eq res o ext = Id}.
    Let now $ g \in \ext(E_\sharp(\sigma_{-\lambda}, \phi)) $. Then, there exists $ f \in E_\sharp(\sigma_{-\lambda}, \phi) $ such that $ g = \ext(f) $. Thus, $ \rest(g) = f \in E_\sharp(\sigma_{-\lambda}, \phi) $ and $ (\ext \circ \rest)(g) = \ext(f) = g $. Hence, $ \ext \circ \rest = \Id $ on $ \ext(E_\sharp(\sigma_{-\lambda}, \phi)) $. So,
    \[
        \restricted{ \ext }{ E_\sharp(\sigma_{-\lambda}, \phi) } \colon E_\sharp(\sigma_{-\lambda}, \phi) \to \ext(E_\sharp(\sigma_{-\lambda}, \phi))
    \]
    is a linear isomorphism and a homeomorphism with continuous linear inverse $ \rest $ as $ \rest $ is continuous and linear and as $ \ext $ is continuous and linear on $E_\sharp(\sigma_{-\lambda}, \phi)$ by Theorem 5.10 of \cite[p.103]{BO00}.
\end{proof}

Let $ \{s_1, \dotsc, s_n\} $ be a basis of $S_\sharp(\sigma_{-\lambda}, \phi)$ (here $n$ may depend on $ \sharp $).
\\ Choose $ f_j \in \C^{-\infty}_\sharp(\dX, V(\sigma_{-\lambda}, \phi))^\Gamma \cap \im(A) $ such that $ s_j = \rest(f_j) $. For $ s = \sum_{j=1}^n a_j s_j \in S_\sharp(\sigma_{-\lambda}, \phi) $, set $ \imath_\sharp(s) = \sum_{j=1}^n a_j f_j $. Then, $ \imath_\sharp $ is a continuous linear bijective map from $S_\sharp(\sigma_{-\lambda}, \phi)$ to $ \imath_\sharp(S_\sharp(\sigma_{-\lambda}, \phi)) \subset \C^{-\infty}_\sharp(\dX, V(\sigma_{-\lambda}, \phi))^\Gamma \cap \im(A) $ having $ \rest $ as continuous inverse.

\begin{lem}\label{Lem res induces an embedding}
    The maps
    \[
        \restricted{ \rest }{ \Ann(D_2) }  \colon \Ann(D_2) \to \C^{-\infty}(B, V_B(\sigma_{-\lambda}, \phi))
    \]
    and
    \[
        \restricted{ \rest }{ \Ann_\Omega(D_2) }  \colon \Ann_\Omega(D_2) \to \Cinf(B, V_B(\sigma_{-\lambda}, \phi))
    \]
    are topological embeddings. 
    It follows that
    \begin{enumerate}
    \item $ \rest(\C^{-\infty}(\dX, V(\sigma_{-\lambda}, \phi))^\Gamma \cap \im(A)) $ is closed in $ \C^{-\infty}(B, V_B(\sigma_{-\lambda}, \phi)) $, and that
    \item $ \rest(\C^{-\infty}_\Omega(\dX, V(\sigma_{-\lambda}, \phi))^\Gamma \cap \im(A)) $ is closed in $ \Cinf(B, V_B(\sigma_{-\lambda}, \phi)) $.
    \end{enumerate}
\end{lem}

\begin{proof}
    As $ \bigl(\im(A) \cap \C^{-\infty}_\sharp(\dX, V(\sigma_{-\lambda}, \phi))^\Gamma\bigr)/D_2 \simeq \Ann_\sharp(D_2) $, 
    it remains to prove the following:
    \\ The map $ \rest $ induces a topological isomorphism from $ \Ann_\sharp(D_2) $ to
    \[
        \rest\bigl(\C^{-\infty}_\sharp(\dX, V(\sigma_{-\lambda}, \phi))^\Gamma \cap \im(A)\bigr) \ .
    \]
    Let us prove this now. By \eqref{eq rest(C^(-infty)(dX, V(sigma_mu, phi))^Gamma) = E_mu oplus S_mu}, we have
    \[
        \rest\bigl(\C^{-\infty}_\sharp(\dX, V(\sigma_{-\lambda}, \phi))^\Gamma \cap \im(A)\bigr) = E_\sharp(\sigma_{-\lambda}, \phi) \oplus S_\sharp(\sigma_{-\lambda}, \phi) \ .
    \]
    As the sum of two continuous linear maps between topological vector spaces is a continuous linear map,
    \[
        \restricted{\ext}{ \tilde{E}_\sharp(\sigma_{-\lambda}, \phi) } \oplus \, \imath_\sharp \colon \rest(\C^{-\infty}_\sharp(\dX, V(\sigma_{-\lambda}, \phi))^\Gamma \cap \im(A)) \to W_\sharp \ ,
    \]
    where
    \[
        W_\sharp := \C^{-\infty}(\Lambda, V(\sigma_{-\lambda}, \phi))^\Gamma + \C^{-\infty}_\sharp(\dX, V(\sigma_{-\lambda}, \phi))^\Gamma \cap \im(A) \ ,
    \]
    is a continuous linear map when we equip the space $ W_\sharp $ with the subspace topology.

    As $ \Ann_\sharp(D_2) $ is a closed finite-codimensional subspace of $ W_\sharp $ (this is a consequence of Lemma \ref{Lem TopologicalDirectSum D_1^(-infty) oplus D_2}), there is a continuous projection $ p_\sharp $ from $ W_\sharp $ to $ \Ann_\sharp(D_2) $ by Lemma \ref{Lem X = Y oplus F topological direct sum if Y is closed and if F is finite-dimensional}.
    It follows from this and Lemma \ref{Lem res restricted to E_sharp is an isomorphism} that
    \begin{equation}\label{eq explicit description of (rest|Ann(D_2))^(-1)}
        p_\sharp \circ (\restricted{\ext}{ \tilde{E}_\sharp(\sigma_{-\lambda}, \phi) } \oplus \, \imath_\sharp)
    \end{equation}
    is the continuous linear inverse map of $ \rest $.
\end{proof}

\begin{prop}\label{Prop about structure} \nlenum
    \begin{enumerate}
    \item\label{Prop about structure eq1}
        The linear map
        \[
            v + D_2 \in \big(\im(A) \cap \C^{-\infty}(\dX, V(\sigma_{-\lambda}, \phi))^\Gamma\big)/D_2 \mapsto (v' \in D_1 \mapsto (v', v)) \in D'_1
        \]
        is well-defined and injective.
    \item\label{Prop about structure eq2}
        For every linear functional $T$ on $ \Ann(D_2) \, (\text{resp. }\Ann_\Omega(D_2) )$ that is continuous with respect to the strong topology, there exists $ \psi_T \in \Cinf(B, V_B(\sigma_\lambda, \phi)) $ $ (\text{resp. } \psi_T \in \C^{-\infty}(B, V_B(\sigma_\lambda, \phi)) ) $, depending on $T$, such that
        \[
            T = (\rest(\cdot), \psi_T)_B \ .
        \]
    \item If $ v \in \Ann_\Omega(D_2) \, (\text{resp. } v \in \Ann(D_2) )$, then $ v = A(\ext(\psi_T)) $, where $ T \colon g \in \Ann(D_2) \mapsto (g, v) \, (\text{resp. } T \colon g \in \Ann_\Omega(D_2) \mapsto (g, v) )$.
    \item\label{Prop about structure eq4}
        The space $ \Ann(D_2) \, (\text{resp. } \Ann_\Omega(D_2) )$ is equal to $ D_{1, -\infty} \, (\text{resp. } D_1) $.
        \\ In particular, $ D_{1, -\infty} \, (\text{resp. } D_1) $ is closed in
        \[
            \im(A) \cap \C^{-\infty}(\dX, V(\sigma_{-\lambda}, \phi))^\Gamma \quad (\text{resp. } \im(A) \cap \C^{-\infty}_\Omega(\dX, V(\sigma_{-\lambda}, \phi))^\Gamma) \ .
        \]
        Thus, $ D_{1, -\infty} $ is the closure of $ D_1 $ in $ \C^{-\infty}(\dX, V(\sigma_{-\lambda}, \phi))^\Gamma $.
    \end{enumerate}
\end{prop}

\begin{proof}
    Recall that $ \star $ is $ \infty $ (or $ -\infty $, resp.) and that $ \sharp $ is $ \emptyset $ (or $ \Omega $, resp.).

    Let $ E $ be $ \Ann(D_2) $ (resp. $ \Ann_\Omega(D_2) $) and let $ F $ be $ \Ann_\Omega(D_2) $ (resp. $ \Ann(D_2) $).
    \begin{enumerate}
    \item It follows from Theorem \ref{Thm (Af, g) = C (res(f), res(g))} that the map considered in \eqref{Prop about structure eq1} is well-defined.
        Let now $ v_1, v_2 \in (V_{\pi, -\infty} \otimes V_\phi)^\Gamma $ be such that
        \[
            (g, v_1) = (g, v_2)
        \]
        for all $ g \in D_1 $. Thus, $ v_1 - v_2 \in \Ann(D_1) = D_2 $ by Lemma \ref{Lem Ann(D_1) = D_2}.
    \item Let $ T \colon E \to \CC $ be a continuous linear functional.

        For $ \varphi = \rest(f) \in \rest(E) $, define $ \tilde{T}(\varphi) = T(f) $.
        Since
        \[
            \rest \colon E \to \rest(E) \subset \C^{-\star}(B, V_B(\sigma_{-\lambda}, \phi))
        \]
        is a topological embedding map by Lemma \ref{Lem res induces an embedding}, this defines a continuous linear functional from $ \rest(E) $ to $ \CC $.

        By the Hahn-Banach theorem, there exists $ \tilde{\tilde{T}} $ a continuous linear functional from $ \C^{-\star}(B, V_B(\sigma_{-\lambda}, \phi)) $ to $ \CC $ such that $ \restricted{ \tilde{\tilde{T}} }{ \rest(E) } = \tilde{T} $.
        Thus, $ \tilde{\tilde{T}} \in \C^\star(B, V_B(\tilde{\sigma}_{-\lambda}, \phi)) $.
        Hence, by using the canonical conjugate-linear isomorphism $ I \colon \C^\star(B, V_B(\sigma_{\lambda}, \phi)) \to \C^\star(B, V_B(\tilde{\sigma}_{\lambda}, \phi)) $, one gets
        \[
            \psi_T \in \C^\star(B, V_B(\sigma_\lambda, \phi))
        \]
        (depending on $T$) such that
        \[
            \tilde{\tilde{T}}(\varphi) = (\varphi, \psi_T)_B \qquad (\varphi \in \rest(E)) \ .
        \]
        So, for $ g \in E $, we have
        \[
            T(g) = \tilde{T}(\rest(g)) = \tilde{\tilde{T}}(\rest(g)) = (\rest(g), \psi_T)_B \ .
        \]
    \item Let $ v $ and $ T $ be as above. Then, $ (v, v') = (\psi_T, \rest(v'))_B $ for all $ v' \in E $. Since $ \beta v \in \Ann(D_2) $, this holds even for all $ v' \in \im(A) \cap \C^{-\infty}(\dX, V(\sigma_{-\lambda}, \phi))^\Gamma $ (resp. $ v' \in \im(A) \cap \C^{-\infty}_\Omega(\dX, V(\sigma_{-\lambda}, \phi))^\Gamma $) by Lemma \ref{Lem TopologicalDirectSum D_1^(-infty) oplus D_2}.
        By Theorem \ref{Thm (Af, g) = C (res(f), res(g))}, there exists a real constant $ C \neq 0 $ such that
        \[
            (v', v) = (\rest(v'), \psi_T)_B = \frac1C (v', A(\ext(\psi_T)))
        \]
        for all $ v' \in D_1 $. Thus, $ \frac1C A(\ext(\psi_T)) - v \in \Ann(D_1) = D_2 $.
        But on the other hand, $ \frac1C A(\ext(\psi_T)) - v \in \Ann(D_2) $. Thus, $ v = \frac1C A(\ext(\psi_T)) $.
    \item By Theorem \ref{Thm (Af, g) = C (res(f), res(g))}, $ D_{1, -\infty} $ is contained in $ \Ann(D_2) $ and $ D_1 $ is contained in $ \Ann_\Omega(D_2) $.
        Let $ v \in \Ann(D_2) $ (resp. $ \Ann_\Omega(D_2) $). Then, $ (\cdot, v) \in \Ann_\Omega(D_2)' $ (resp. $ \Ann(D_2)' $).

        It follows from \eqref{Prop about structure eq1} and \eqref{Prop about structure eq2} that there exist $ \varphi \in \C^{-\star}(B, V_B(\sigma_\lambda, \phi)) $ and $ g \in D_2 $ such that $ v = A(\ext(\varphi)) + g $.
        Then,
        \[
            g = v - A(\ext(\varphi)) \in D_2 \cap \Ann(D_2) = \{0\} \ .
        \]
        So, $ \Ann(D_2) $ (resp. $ \Ann_\Omega(D_2) $) is equal to $ D_{1, -\infty} $ (resp. $ D_1 $).
        In particular, $ D_{1, -\infty} $ (resp. $ D_1 $) is closed in $ \C^{-\infty}_\sharp(\dX, V(\sigma_{-\lambda}, \phi))^\Gamma $. As moreover $ D_1 $ is dense in $ D_{1, -\infty} $, the assertion follows. 
    \end{enumerate}
\end{proof}

\begin{cor}\label{Cor D_1' = D_(1, -infty)}
    The conjugate-linear strong dual $ D_1' $ $(\text{resp. } (D_{1, -\infty})' )$ of $ D_1 $ $(\text{resp. } D_{1, -\infty} )$ is equal to $ D_{1, -\infty} $ $(\text{resp. } D_1 )$ as topological spaces.
\end{cor}

\begin{proof}
    It follows from Proposition \ref{Prop about structure} that $ D_1' $ is equal to $ D_{1, -\infty} $ as a vector space. Moreover, $ D_{1, -\infty} \to D_1' , \, g \mapsto (g, \cdot) $ is continuous.

    Let $T$ be a continuous linear functional on $ D_1 $. For $ \varphi \in \Cinf(B, V_B(\sigma_\lambda, \phi)) $, set
    \[
        \Psi_T(\varphi) = T\big(A(\ext(\varphi))\big) \ .
    \]
    Then, $\Psi_T$ is a conjugate-linear functional on $ \Cinf(B, V_B(\sigma_\lambda, \phi)) $.
    \\ Since $ \Cinf(B, V_B(\sigma_\lambda, \phi))' $ is topologically isomorphic to $ \C^{-\infty}(B, V_B(\sigma_{-\lambda}, \phi)) $, we can view $ \Psi_T $ as an element in $ \C^{-\infty}(B, V_B(\sigma_{-\lambda}, \phi)) $. But it follows from Proposition \ref{Prop about structure} that $ \Psi_T \in \rest(D_{1, -\infty}) $. Set $ f_T = \rest^{-1}(\Psi_T) \in D_{1, -\infty} $. It follows from Lemma \ref{Lem res induces an embedding} that $ T \mapsto f_T $ is continuous. This map is by construction the inverse of the above map. Thus, $ D_1' = D_{1, -\infty} $ as topological spaces. One can show similarly that $ (D_{1, -\infty})' = D_1 $ as topological spaces.
\end{proof}

\begin{cor}\label{Cor Schwartz topology = smooth on Omega topology}
    $ (V_{\pi, -\infty} \otimes V_\phi)^\Gamma_\CS $ is equal to $ \im(A) \cap \C^{-\infty}_\Omega(\dX, V(\sigma_{-\lambda}, \phi))^\Gamma $ as topological spaces.
\end{cor}

\begin{proof}
    The corollary follows by the open mapping theorem as $ (V_{\pi, -\infty} \otimes V_\phi)^\Gamma_\CS $ can be continuously injected into $ ((V_{\pi, -\infty} \otimes V_\phi)^\Gamma)' $ and as the latter space is topologically isomorphic to $ D_1 \oplus D_2 = \im(A) \cap \C^{-\infty}_\Omega(\dX, V(\sigma_{-\lambda}, \phi))^\Gamma $ as topological spaces, by Corollary \ref{Cor D_1' = D_(1, -infty)} 
    and the open mapping theorem.

\end{proof}

\begin{prop}\label{Prop space of Schwartz vectors top. isom to space of Gamma-inv dist vectors}
    \begin{sloppypar}
        The conjugate-linear strong dual $ \big( (V_{\pi, -\infty} \otimes V_\phi)^\Gamma_\CS \big)' $ of $ (V_{\pi, -\infty} \otimes V_\phi)^\Gamma_\CS $ is topologically isomorphic to $ (V_{\pi, -\infty} \otimes V_\phi)^\Gamma $.
    \end{sloppypar}
\end{prop}

\begin{rem}
    It follows that $ (V_{\pi, -\infty} \otimes V_\phi)^\Gamma $ is a DF-space.
\end{rem}

\begin{proof}
    By Lemma \ref{Lem TopologicalDirectSum D_1^(-infty) oplus D_2} and Proposition \ref{Prop about structure}, $ (V_{\pi, -\infty} \otimes V_\phi)^\Gamma_\CS = D_1 \oplus D_2 $ (topological direct sum).
    As $ D_2 $ is finite-dimensional, the conjugate-linear strong dual $ D_2' $ of $ D_2 $ is topologically isomorphic to $ D_2 $.
    By Corollary \ref{Cor D_1' = D_(1, -infty)} combined with Corollary \ref{Cor Schwartz topology = smooth on Omega topology}, the conjugate-linear dual $ D_1' $ of $ D_1 $ is topologically isomorphic to
    \[
        D_{1, -\infty} = A(\C^{-\infty}(\dX, V(\sigma_\lambda, \phi))^\Gamma) \ .
    \]
    The proposition follows now from Lemma \ref{Lem X = Y oplus F topological direct sum if Y is closed and if F is finite-dimensional}.
\end{proof}

\begin{prop}
    The inclusions
    \[
        (V_{\pi, -\infty} \otimes V_\phi)^\Gamma_\CS \subset (V_{\pi, -\infty} \otimes V_\phi)^\Gamma_d \subset (V_{\pi, -\infty} \otimes V_\phi)^\Gamma \subset V_{\pi, -\infty} \otimes V_\phi
    \]
    are continuous injections.
\end{prop}

\begin{proof}
    By Lemma \ref{Lem space of Schwartz vectors can be continuously injected into space of square-integrable vectors}, it remains to show that
    \[
        (V_{\pi, -\infty} \otimes V_\phi)^\Gamma_d \hookrightarrow (V_{\pi, -\infty} \otimes V_\phi)^\Gamma
    \]
    is continuous. This follows from Proposition \ref{Prop space of Schwartz vectors top. isom to space of Gamma-inv dist vectors} and the fact that dual map is continuous (see Lemma \ref{Lem dual map is continuous}).
\end{proof}

\needspace{2\baselineskip}

\begin{thm}\label{Thm Decompositions of (V_(pi, -infty) otimes V_phi)_CS^Gamma and of (V_(pi, -infty) otimes V_phi)^Gamma}
    We have the following topological direct sums
    \[
        (V_{\pi, -\infty} \otimes V_\phi)^\Gamma 
        = A(\C^{-\infty}(\dX, V(\sigma_\lambda, \phi))^\Gamma) \oplus (\C^{-\infty}(\Lambda, V(\sigma_{-\lambda}, \phi))^\Gamma \cap \im(A))
    \]
    and
    \begin{align*}
        (V_{\pi, -\infty} \otimes V_\phi)_\CS^\Gamma &= \C^{-\infty}_\Omega(\dX, V(\sigma_{-\lambda}, \phi))^\Gamma \cap \im(A) \\
        &= A(\C^{-\infty}_\Omega(\dX, V(\sigma_\lambda, \phi))^\Gamma) \oplus (\C^{-\infty}(\Lambda, V(\sigma_{-\lambda}, \phi))^\Gamma \cap \im(A)) \ .
    \end{align*}
\end{thm}

\begin{rem} \nlenum
    \begin{enumerate}
    \item The space $ A(\C^{-\infty}_\Omega(\dX, V(\sigma_\lambda, \phi))^\Gamma) $ is infinite-dimensional by Lemma 9.6 of \cite[p.139]{BO00} combined with Lemma 6.6 of \cite[p.114]{BO00} and it is contained in $ \C^{-\infty}_\Omega(\dX, V(\sigma_{-\lambda}, \phi))^\Gamma $. Thus, $ (V_{\pi, -\infty} \otimes V_\phi)_\CS^\Gamma $ is infinite-dimensional. Hence, $ N_{\pi'} $ is also infinite-dimensional --- a result that was already proven by M.~Olbrich and U.~Bunke (see Proposition 9.8 of \cite[p.141]{BO00}).
    \item By Theorem 4.7 of \cite[p.93]{BO00}, $ \C^{-\infty}(\Lambda, V(\sigma_\mu, \phi))^\Gamma = \{0\} $ for $ \R(\mu) > \delta_\Gamma $. Hence, $ \C^{-\infty}(\Lambda, V(\sigma_{-\lambda}, \phi))^\Gamma \cap \im(A) = \{0\} $ for $ \lambda \in (0, -\delta_\Gamma) $.
    \item The above decompositions depends on the choice of the principal series in which the discrete series is embedded into (choice of $ \sigma \in \hat{M} $).
    \end{enumerate}
\end{rem}

\begin{proof}
    By Corollary \ref{Cor decomp. of space of cuspidal vectors if pi is a DS rep}, we have
    \[
        (V_{\pi, -\infty} \otimes V_\phi)_\CS^\Gamma
        \underset{\text{dense}}{\subset} N_{\pi'} = (V_{\pi, -\infty} \otimes V_\phi)_d^\Gamma
        \subset \C^{-\infty}(\dX, V(\sigma_{-\lambda}, \phi))^\Gamma \ .
    \]
    The theorem follows now directly from Lemma \ref{Lem TopologicalDirectSum D_1^(-infty) oplus D_2}, Proposition \ref{Prop about structure} and Corollary \ref{Cor Schwartz topology = smooth on Omega topology}.
\end{proof}

Let $ (\pi, V_\pi) $ be an admissible representation of $G$ of finite length on a reflexive Banach space. Recall that for $ f \in \C_c(G) $ and $ v \in V_\pi $, we defined
\begin{equation}\label{eq pi(f)v = int_G f(g) pi(g)v dg}
    \pi(f)v = \int_G f(g) \pi(g)v \, dg \in V_\pi
\end{equation}
(cf. Definition \ref{Dfn pi(f)v}).
Let $ f \in \Ccinf(G) $ and $ v \in V_{\pi, -\infty} $. Then, $ \pi(f)v \in V_{\pi, \infty} $ (see, e.g., \cite[p.136]{Dijk}).

If $ V $ is a complex vector space, then we denote by $ \bar{V} $ the vector space having the same elements and the same additive group structure as $V$ but which is equipped with the following scalar multiplication:
\[
    \lambda \cdot \ubar{v} := \bar{\lambda} \ubar{v}
\]
for all $ \lambda \in \CC $ and $ \ubar{v} \in \bar{V} $. Let $ \ubar{\cdot} \colon V \to \bar{V} , \, v \mapsto \ubar{v} $ be the natural isomorphism between these two vector spaces.

For $ \varphi \in V_{\bar{\pi}', -\infty} $, $ f \in \CS(G) $ and $ v \in V_{\pi, \infty} $, set
\begin{equation}\label{eq (varphi, pi(f)v) = (c_(varphi, v), f)_G}
    (\varphi, \pi(f)v) = (c_{\varphi, \ubar{v}}, f)_G \ .
\end{equation}
Here, $ (\cdot, \cdot)_G $ denotes the pairing between $ \CS'(G) $ and $ \CS(G) $. The right-hand side is well-defined as $ c_{\varphi, \ubar{v}} \in \CS'(G) $ by Lemma 8.4 of \cite[p.130]{BO00}. By Lemma 8.5 of \cite[p.131]{BO00}, the matrix coefficient map
\[
    \varphi \in V_{\bar{\pi}', -\infty} \mapsto c_{\varphi, \ubar{v}} \in \CS'(G)
\]
is continuous. Thus,
\[
    (\cdot, \pi(f)v) := (c_{\cdot, \ubar{v}}, f)_G
\]
is a continuous functional on $ V_{\bar{\pi}', -\infty} $. If $ f \in \Ccinf(G) $ and $ \varphi \in V_{\bar{\pi}', \infty} $, then
\[
    (c_{\varphi, \ubar{v}}, f)_G = \int_G c_{\varphi, \ubar{v}}(g) \overline{f(g)} \, dg
    = \int_G \int_K \langle \varphi(k), \bar{\pi}(g)\ubar{v} \rangle \, dk \, \overline{f(g)} \, dg \ .
\]
By Fubini's theorem, this is again equal to
\[
    \int_K \langle \varphi(k), \int_G \overline{f(g) \pi(g)v} \, dg \rangle \, dk = (\varphi, \pi(f)v) \ .
\]
Since $ (c_{\cdot, \ubar{v}}, f)_G $ and $ (\cdot, \pi(f)v) $ are continuous functionals on $ V_{\bar{\pi}', -\infty} $, this holds also for every $ \varphi \in V_{\bar{\pi}', -\infty} $.
So, \eqref{eq pi(f)v = int_G f(g) pi(g)v dg} is compatible with \eqref{eq (varphi, pi(f)v) = (c_(varphi, v), f)_G}.

As $ (F, \cdot)_G \, $ ($ F \in \CS'(G) $ fixed) is a continuous functional on $ \CS(G, V_\phi) $, $ (\varphi, \pi(f)v) $ depends continuously on $f$ when $ \varphi $ and $v$ are fixed.
So, if $ (f_j)_j $ is a sequence in $ \Ccinf(G) $ converging to $f$, then $ \pi(f_j)v $ converges to $ \pi(f)v $ in $ V_{\pi, -\infty} $.

For $ f \in \CS(G) $ and $ v \in V_{\pi, \infty} \otimes V_\phi $, we set $ \pi(f)v = (\pi(f) \otimes \Id)v $.

\medskip

Let us give now a formula for the projection of a Schwartz function onto $ D_1 \otimes V_{\pi'}(\gamma) $.

\begin{lem}
    Let $f$ be a Schwartz vector and let $ v \in V_{\pi', K} $ be of norm 1. Then, $ \rest(f) $ is up to a constant equal to
    \[
        \overline{\pi_*(\pi'(\chi \overline{c_{f, v}})v)} \in \Cinf(B, V_B(\sigma_{-\lambda}, \phi)) \ .
    \]
\end{lem}

\begin{proof}
    Let $ v \in V_{\pi', K} $ be of norm 1.
    \\ For $ \varphi \in \Cinf(B, V_B(\sigma_\lambda, \phi)) $ and a Schwartz vector $f$, we have
    \begin{multline}\label{eq Proof of lemma, C(varphi, res(Af))_B = ...}
        C (\varphi, \rest(f))_B = (A(\ext(\varphi)), f) = \int_{\Gamma \bs G} c_{q(\ext(\varphi)), v}(g) \overline{c_{f, v}(g)} \, dg \\
        = \int_G c_{\ext(\varphi), v}(g) \chi(gK) \overline{c_{f, v}(g)} \, dg \ .
    \end{multline}
    Since $ \chi\overline{c_{f, v}} \in \CS(G, V_\phi) $ by Lemma \ref{Lem chi f in CS(G, V_phi)}, $ \Ccinf(G, V_\phi) $ is dense in $ \CS(G, V_\phi) $ by Proposition \ref{Prop Ccinf(Gamma|G) is dense}, there is a sequence $ (F_j)_j $ in $ \Ccinf(G, V_\phi) $ converging to $ \chi\overline{c_{f, v}} $ in $ \CS(G, V_\phi) $. Thus, \eqref{eq Proof of lemma, C(varphi, res(Af))_B = ...} is again equal to
    \[
        \lim_{j \to \infty} \int_G c_{\ext(\varphi), v}(g) F_j(g) \, dg
        = \lim_{j \to \infty} \int_G \int_B \langle \varphi(x), \pi_*(\pi'(g) v)(x) \rangle \, F_j(g) \, dg
    \]
    by Theorem \ref{Thm (Af, g) = C (res(f), res(g))}. 
    Let $ F \in \Ccinf(G, V_\phi) $. Then, $ \pi'(F)v $ is well-defined, belongs to $ V_{\pi', \infty} \otimes V_\phi $ and $ g \in G \mapsto F(g) \pi'(g) v $ is integrable. Thus, $ \pi_*(\pi'(F)v) $ is well-defined and belongs to $ \Cinf(B, V_B(\sigma_{-\lambda}, \phi)) $.
    Moreover, the above integral is again equal to
    \[
        \lim_{j \to \infty} \int_B \int_G \langle \varphi(x), \pi_*(F_j(g) \pi'(g) v)(x) \rangle \, dg
        = \lim_{j \to \infty} (\varphi, \overline{\pi_*(\pi'(F_j) v)})_B
    \]
    by Fubini's theorem.
    Since $(\varphi, \cdot)_B$ is continuous on $ \C^{-\infty}(B, V_B(\sigma_{-\lambda}, \phi)) $, since $ \pi'(F_j) v $ converges to $ \pi'(\chi\overline{c_{f, v}})v $ in $ V_{\pi, -\infty} \otimes V_\phi $ and since $ \pi_* $ is continuous on $ \C^{-\infty}(\dX, V(\sigma_{-\lambda}, \phi)) $ ($ \ext $ is continuous by Theorem 5.10 of \cite[p.103]{BO00}), this is again equal to
    \[
        (\varphi, \overline{\pi_*(\pi'(\chi\overline{c_{f, v}}) v)(x)})_B \ .
    \]
    So, $ C \rest(f) = \overline{\pi_*(\pi'(\chi \overline{c_{f, v}})v)} \in \Cinf(B, V_B(\sigma_{-\lambda}, \phi)) $.
\end{proof}

\begin{cor}
    Let $ \{ v_i \} $ be an orthonormal basis of $ V_{\pi'}(\gamma) $, let $ f \in \CS(\Gamma \bs G, \phi) $ and let $ \varphi_i = \frac1C \overline{\pi_*(\pi^{\sigma, -\lambda}(\chi \bar{f})v_i)} $.
    Let $ p_{D_1 \otimes V_{\pi'}(\gamma)} $ be the projection onto $ D_1 \otimes V_{\pi'}(\gamma) $ provided by Theorem \ref{Thm Decompositions of (V_(pi, -infty) otimes V_phi)_CS^Gamma and of (V_(pi, -infty) otimes V_phi)^Gamma}.
    Then, $ p_{D_1 \otimes V_{\pi'}(\gamma)}\big(\F_{\pi'}(f)\big) $ is equal to
    \[
        \sum_i (\restricted{\rest}{D_1})^{-1}(\varphi_i) \otimes v_i \ .
    \]
\end{cor}

\begin{rem} \nlenum
    \begin{enumerate}
    \item When $ \lambda > \delta_\Gamma $, then the above formula is quite explicit as $ \pi_* $ is then also defined without meromorphic extension.
    \item Since $ \ext $ is regular at $ \mu = \lambda $ by Lemma \ref{Lem ext regular}, $ \pi_* $ is regular at $ \mu = -\lambda $.
    \item As we have seen, one can construct $ (\restricted{\rest}{D_1})^{-1} $ with the help of $ \ext $, see \eqref{eq explicit description of (rest|Ann(D_2))^(-1)} in the proof of Lemma \ref{Lem res induces an embedding} (p.\pageref{Lem res induces an embedding}).
    \end{enumerate}
\end{rem}

\begin{proof}
    Let $ \{ v_i \} $ be an orthonormal basis of $ V_{\pi'}(\gamma) $, let $ f \in \CS(\Gamma \bs G, \phi) $ and let $ \varphi_i = \frac1C \overline{\pi_*(\pi'(\chi \bar{f})v_i)} $.
    Then, it follows from the proof of the previous lemma that $ p_{D_1 \otimes V_{\pi'}(\gamma)}\big(\F_{\pi'}(f)\big) $ is equal to $ \sum_i (\restricted{\rest}{D_1})^{-1}(\varphi_i) \otimes v_i $.
\end{proof}

\newpage

\subsubsection[\texorpdfstring{Topological decomposition of $ \CS(\Gamma \bs G, \phi)_{ds}(\gamma) $ and of \\ \mbox{$ \CS'(\Gamma \bs G, \phi)_{ds}(\gamma) $}}{Topological decompositions II}]{\texorpdfstring{Topological decomposition of $ \CS(\Gamma \bs G, \phi)_{ds}(\gamma) $ and of \mbox{$ \CS'(\Gamma \bs G, \phi)_{ds}(\gamma) $}}{Topological decompositions II}}

For $ \gamma \in \hat{K} $, set\index[n]{CGammaG_ds@$ \CS(\Gamma \bs G, \phi)_{ds}(\gamma) $} 
\[
    \CS(\Gamma \bs G, \phi)_{ds}(\gamma) = \CS(\Gamma \bs G, \phi)(\gamma) \cap L^2(\Gamma \bs G, \phi)_{ds}
\]
and
\[
    \CS'(\Gamma \bs G, \phi)_{ds}(\gamma) = (\CS(\Gamma \bs G, \phi) \cap L^2(\Gamma \bs G, \phi)_{ds})'(\gamma) \ .
\]

\begin{prop}\label{Prop Decomposition of C(Gamma|G) cap L^2(Gamma|G)_ds(gamma)}
    Let $ \gamma \in \hat{K} $. Then, $ c_\pi\big(\F_\pi( \CS(\Gamma \bs G, \phi)(\gamma) )\big) $ is equal to
    \begin{multline*}
        \spn\{ c_{\pi}(f \otimes v) \mid f \in (V_{\pi', -\infty} \otimes V_\phi)^\Gamma , \, v \in V_{\pi}(\gamma) \} \cap \CS(\Gamma \bs G, \phi) \\
        = \spn\{ c_{\pi}(f \otimes v) \mid f \in (V_{\pi', -\infty} \otimes V_\phi)^\Gamma_\CS , \, v \in V_{\pi}(\gamma) \} \ .
    \end{multline*}
    Moreover, we have
    \begin{multline}\label{eq Decomposition of CS(Gamma|G, phi)(gamma)}
        \CS(\Gamma \bs G, \phi)_{ds}(\gamma)  \\
        = \bigoplus_{ \pi \in \hat{G} \, : \, (V_{\pi', -\infty} \otimes V_\phi)^\Gamma_d \neq \{0\} } \spn\{ c_{\pi}(f \otimes v) \mid f \in (V_{\pi', -\infty} \otimes V_\phi)^\Gamma_\CS , \, v \in V_{\pi}(\gamma) \}
    \end{multline}
    (finite topological direct sum).
\end{prop}

\begin{proof}
    Let
    \[
        S_\pi = c_\pi((V_{\pi, -\infty} \otimes V_\phi)^\Gamma_d \otimes V_{\pi}(\gamma)) \ .
    \]
    This space is closed in $ L^2(\Gamma \bs G, \phi) $ by Lemma \ref{Cor c_pi(N_pi hat(otimes) V_pi), c_pi(N_pi otimes F) closed}.

    Set $ \tilde{S}_\pi = S_{\pi} \cap \CS(\Gamma \bs G, \phi) $ (closed in $ \CS(\Gamma \bs G, \phi) $ by Lemma \ref{Lem V = (S cap V) oplus (S^perp cap V)}).
    
    Let $ p_{S_{\pi}} $ denote the orthogonal projection of $ L^2(\Gamma \bs G, \phi) $ onto $ S_\pi $.
    Then, $ p_{S_{\pi}}(g) \in \CS(\Gamma \bs G, \phi) $ for all $ g \in \CS(\Gamma \bs G, \phi) $, $ p_{S_{\pi}} = c_\pi \circ \F_\pi $ on $ L^2(\Gamma \bs G, \phi)(\gamma) $ and
    \[
        p_{S_{\pi}}(\CS(\Gamma \bs G, \phi)(\gamma))
        = \tilde{S}_\pi = \spn\{ c_{\pi}(f \otimes v)
        \mid f \in (V_{\pi', -\infty} \otimes V_\phi)^\Gamma_\CS , \, v \in V_{\pi}(\gamma) \}
    \]
    by Proposition \ref{Prop c_pi(F_pi(f)) for discrete series}. The first part follows now from the above combined with the following fact:
    \begin{multline*}
        \spn\{ c_{\pi}(f \otimes v) \mid f \in (V_{\pi', -\infty} \otimes V_\phi)^\Gamma , \, v \in V_{\pi}(\gamma) \} \cap L^2(\Gamma \bs G, \phi) \\
        = \spn\{ c_{\pi}(f \otimes v) \mid f \in (V_{\pi', -\infty} \otimes V_\phi)^\Gamma_d , \, v \in V_{\pi}(\gamma) \} \ .
    \end{multline*}
    Indeed, one can easily show this by using an orthonormal basis of $ V_{\pi}(\gamma) $ and by using the natural topological isomorphism between $ (V_{\pi', -\infty} \otimes V_\phi)^\Gamma_d $ and $ ((V_{\pi', -\infty} \otimes V_\phi)^\Gamma_d)' $ (compare with the proof of Proposition \ref{Prop c_pi(F_pi(f)) for discrete series}).

    As for a fixed $K$-type only finitely many terms are nonzero in the above direct sum, \eqref{eq Decomposition of CS(Gamma|G, phi)(gamma)} follows by Proposition \ref{Prop CS(Gamma|G) subset L^2(Gamma|G)} and Lemma \ref{Lem V = (S cap V) oplus (S^perp cap V) generalisation} applied to $ H = L^2(\Gamma \bs G, \phi)_{ds}(\gamma) $.
\end{proof}

\begin{lem}\label{Lem V' = Vbar_1 oplus ... oplus Vbar_d}
    Let $V$ be a closed subspace of $ \CS(\Gamma \bs G, \phi) $ and let $ V_1, \dotsc, V_d $ ($ d \in \NN $) be closed vector subspaces of $V$ such that $ V = \bigoplus_{j=1}^d V_j $ is a (topological) orthogonal direct sum decomposition. Let $ \bar{V}_j $ be the closure of $ V_j $ in $ \CS'(\Gamma \bs G, \phi) $. Then,
    \[
        V' = \bigoplus_{j=1}^d \bar{V}_j
    \]
    (topological direct sum). Moreover, the space $ V_i $, viewed as a subspace of $ V_i' $, is dense in $ V_i' $ and $ \bar{V}_i $ $ (i \in \{1, \dotsc, d\}) $ is topologically isomorphic to $ V_i' $. Hence, $ \bar{V}_i $ is a DF-space if $V$ is a Fréchet space.
\end{lem}

\begin{proof}
    Since $ \CS(\Gamma \bs G, \phi) $ is dense in $ \CS'(\Gamma \bs G, \phi) $, $V$ is also dense in $V'$.
    Indeed, this follows by the arguments below applied to
    \[
        \CS(\Gamma \bs G, \phi) = \underbrace{V}_{=: V_1} \oplus \underbrace{(V^\perp \cap \CS(\Gamma \bs G, \phi))}_{=: V_2}
    \]
    (topological orthogonal direct sum by Lemma \ref{Lem V = (S cap V) oplus (S^perp cap V)}).
    
    As $ V = \bigoplus_{j=1}^d V_j $ is a topological direct sum,
    \[
        V' = \bar{V} = \bar{V}_1 + \dotsb + \bar{V}_d \ .
    \]
    Assume now that $ f \in \bar{V}_k \cap \bar{V}_l $ ($ k \neq l $). Then, there exist $ f_j^{(k)} \in V_k $ converging to $f$ and $ f_j^{(l)} \in V_l $ converging to $f$, too.
    Let $ v \in \bigoplus_{j \neq k} V_j $. Then, $ (f, v) = \lim_{j \to \infty} (f_j^{(k)}, v) = 0 $. So, $ f \in \Ann\bigl(\bigoplus_{j \neq k} V_j\bigr) $. Using the other sequence, one shows that $ f $ belongs also to the annihilator of $ V_k $. Thus, $ f \in \Ann(V) $. Hence, $ f = 0 $.
    The first part follows.

    Let $ \iota_j \colon V'_j \to V' $ be the natural topological embeddings. It follows from Lemma \ref{Lem V' simeq V'_1 oplus V'_2} that $ V' = \bigoplus_{j=1}^d \iota_j(V'_j) $ (topological direct sum). Thus, $ \bar{V}_i = \iota_i(V'_i) \simeq V'_i $ and $ V_i $, viewed as a subspace of $ V_i' $, is dense in $ V_i' $.
    Hence, $ \bar{V}_i $ is a DF-space if $V$ is a Fréchet space.
\end{proof}

\begin{prop}\label{Prop Decomposition of (C(Gamma|G) cap L^2_ds)'(gamma)}
    Let $ \gamma \in \hat{K} $. Then, the space
    \[
        \spn\{ c_\pi(f \otimes v) \mid f \in (V_{\pi', -\infty} \otimes V_\phi)^\Gamma , \, v \in V_{\pi}(\gamma) \}
    \]
    is closed in $ \CS'(\Gamma \bs G, \phi) $ and topologically isomorphic to $ (V_{\pi', -\infty} \otimes V_\phi)^\Gamma \otimes V_{\pi}(\gamma) $.
    Moreover, we have
    \begin{multline}\label{eq Decomposition of (C(Gamma|G) cap L^2_ds)'(gamma)}
        \CS'(\Gamma \bs G, \phi)_{ds}(\gamma) \\
        = \bigoplus_{ \pi \in \hat{G} \, : \, (V_{\pi', -\infty} \otimes V_\phi)^\Gamma_d \neq \{0\} } \spn\{ c_{\pi}(f \otimes v) \mid f \in (V_{\pi', -\infty} \otimes V_\phi)^\Gamma , \, v \in V_{\pi}(\gamma) \}
    \end{multline}
    (topological direct sum).
\end{prop}

\begin{proof}
    By Proposition \ref{Prop space of Schwartz vectors top. isom to space of Gamma-inv dist vectors}, Lemma \ref{Lem space of Schwartz vector is top. isom. to its image} and as $ c_\pi $ is injective, we have
    \begin{align}\label{eq Proof of Prop Decomposition of (C(Gamma|G) cap L^2_ds)'(gamma)}
        & \spn\{ c_\pi(f \otimes v) \mid f \in (V_{\pi', -\infty} \otimes V_\phi)^\Gamma , \, v \in V_{\pi}(\gamma) \} \\
        \simeq & (V_{\pi', -\infty} \otimes V_\phi)^\Gamma \otimes V_{\pi}(\gamma) 
        \simeq \big( (V_{\pi', -\infty} \otimes V_\phi)^\Gamma_\CS \otimes V_{\pi}(\gamma) \big)' \nonumber \\
        \simeq & ( \{ c_\pi(f \otimes v) \mid f \in (V_{\pi', -\infty} \otimes V_\phi)^\Gamma_\CS , \, v \in V_{\pi}(\gamma) \} )' \nonumber
    \end{align}
    as vector spaces.
    It follows from Proposition \ref{Prop Schwartz vectors are dense in the space of Gamma-invariant vectors} that $ \{ c_\pi(f \otimes v) \mid f \in (V_{\pi', -\infty} \otimes V_\phi)^\Gamma_\CS , \, v \in V_{\pi}(\gamma) \} $ is dense in $ \{ c_\pi(f \otimes v) \mid f \in (V_{\pi', -\infty} \otimes V_\phi)^\Gamma , \, v \in V_{\pi}(\gamma) \} $.
    Thus, $ \spn\{ c_\pi(f \otimes v) \mid f \in (V_{\pi', -\infty} \otimes V_\phi)^\Gamma , \, v \in V_{\pi}(\gamma) \} $ is closed and \eqref{eq Decomposition of (C(Gamma|G) cap L^2_ds)'(gamma)} holds by Proposition \ref{Prop Decomposition of C(Gamma|G) cap L^2(Gamma|G)_ds(gamma)} and Lemma \ref{Lem V' = Vbar_1 oplus ... oplus Vbar_d}.
    Moreover,
    \[
        \spn\{ c_\pi(f \otimes v) \mid f \in (V_{\pi', -\infty} \otimes V_\phi)^\Gamma , \,
        v \in V_{\pi}(\gamma) \}
    \]
    is a DF-space when we equip the space with the subspace topology.
    It follows now from the open mapping theorem that the isomorphisms in \eqref{eq Proof of Prop Decomposition of (C(Gamma|G) cap L^2_ds)'(gamma)} are topological.
\end{proof}

\newpage

\subsubsection{Equivalent topologies}

We finish our investigations on the discrete series representations by determining an equivalent topology on the space of invariant distribution vectors, respectively on the space of Schwartz vectors.

\begin{lem}\label{Lem space of tempered Gamma-invariant vectors is equipped with the coarsest top. s.t. the matrix coeff. map is continuous}
     Let $ 0 \neq v \in V_{\pi', K} $. Then, the subspace topology on $ (V_{\pi, -\infty} \otimes V_\phi)^\Gamma $ is the coarsest topology such that the map
     \[
        f \in (V_{\pi, -\infty} \otimes V_\phi)^\Gamma \mapsto c_{f, v} \in \CS'(\Gamma \bs G, \phi)
     \]
     is continuous.
\end{lem}

\begin{proof}
    Denote the subspace topology on $ (V_{\pi, -\infty} \otimes V_\phi)^\Gamma $ by $ \tau_1 $ and the coarsest topology such that the above map is continuous by $ \tau_2 $.

    Similarly as for $ V_{\pi}(\gamma) $ in the proof of Proposition \ref{Prop Decomposition of (C(Gamma|G) cap L^2_ds)'(gamma)}, one can show that $ \{ c_{f, v} \mid f \in (V_{\pi, -\infty} \otimes V_\phi)^\Gamma \} $ is closed and a DF-space when we equip the space with the subspace topology. So, $ (V_{\pi, -\infty} \otimes V_\phi)^\Gamma $ is topologically isomorphic to $ \{ c_{f, v} \mid f \in (V_{\pi, -\infty} \otimes V_\phi)^\Gamma \} $ for any $ 0 \neq v \in V_{\pi', K} $ by the open mapping theorem.

    It follows that $ (V_{\pi, -\infty} \otimes V_\phi)^\Gamma $, equipped with the topology $ \tau_2 $, is topologically isomorphic to $ \{ c_{f, v} \mid f \in (V_{\pi, -\infty} \otimes V_\phi)^\Gamma \} $.
    But $ (V_{\pi, -\infty} \otimes V_\phi)^\Gamma $, equipped with the topology $ \tau_1 $, is also topologically isomorphic to $ \{ c_{f, v} \mid f \in (V_{\pi, -\infty} \otimes V_\phi)^\Gamma \} $. The lemma follows.
\end{proof}

\begin{cor}\label{Cor Equivalent topology on (V_(pi, -infty) otimes V_phi)^Gamma}
     The topology on $ (V_{\pi, -\infty} \otimes V_\phi)^\Gamma $ is the coarsest topology such that the map
     \[
        f \in (V_{\pi, -\infty} \otimes V_\phi)^\Gamma \mapsto (f, \cdot) \in \big( (V_{\pi, -\infty} \otimes V_\phi)^\Gamma_\CS \big)'
     \]
     is continuous.
\end{cor}

\begin{proof}
    The corollary follows from Proposition \ref{Prop Decomposition of C(Gamma|G) cap L^2(Gamma|G)_ds(gamma)}, Lemma \ref{Lem space of tempered Gamma-invariant vectors is equipped with the coarsest top. s.t. the matrix coeff. map is continuous} and the fact that
    \[
        \CS(\Gamma \bs G, \phi)(\gamma) = \CS(\Gamma \bs G, \phi)_{ds}(\gamma)
            \oplus (\CS(\Gamma \bs G, \phi) \cap L^2(\Gamma \bs G, \phi)_{ds}(\gamma)^\perp)
    \]
    (topological direct sum), by Proposition \ref{Prop c_pi(F_pi(f)) for discrete series} and by Lemma \ref{Lem V = (S cap V) oplus (S^perp cap V)}.
\end{proof}

\begin{cor}
     The topology on $ (V_{\pi, -\infty} \otimes V_\phi)_\CS^\Gamma $ is the coarsest topology such that the map
     \[
        f \in (V_{\pi, -\infty} \otimes V_\phi)_\CS^\Gamma \mapsto (f, \cdot) \in \big( (V_{\pi, -\infty} \otimes V_\phi)^\Gamma \big)'
     \]
     is continuous.
\end{cor}

\begin{proof}
    The corollary follows from Proposition \ref{Prop Decomposition of (C(Gamma|G) cap L^2_ds)'(gamma)} and the fact that
    \[
        \CS'(\Gamma \bs G, \phi)(\gamma) = \CS'(\Gamma \bs G, \phi)_{ds}(\gamma)
            \oplus (\CS(\Gamma \bs G, \phi) \cap L^2(\Gamma \bs G, \phi)_{ds}^\perp)'(\gamma)
    \]
    (topological direct sum), by the proof of Corollary \ref{Cor Equivalent topology on (V_(pi, -infty) otimes V_phi)^Gamma} and by Lemma \ref{Lem V' simeq V'_1 oplus V'_2}.
\end{proof}

\newpage

\subsubsection{Decomposition of the Schwartz space}

\begin{prop}\label{Prop CS(Gamma|G, phi)(gamma)}
    Fix $ \gamma \in \hat{K} $. Then,
    \begin{multline*}
        \CS(\Gamma \bs G, \phi)(\gamma) = L^2(\Gamma \bs G, \phi)_{ds}(\gamma) \cap \c{\CS}(\Gamma \bs G, \phi)
        \oplus L^2(\Gamma \bs G, \phi)_U(\gamma) \\
        \oplus L^2(\Gamma \bs G, \phi)_{res}(\gamma)
        \oplus L^2(\Gamma \bs G, \phi)_{ac}(\gamma) \cap \CS(\Gamma \bs G, \phi)
    \end{multline*}
    (topological direct sum).
\end{prop}

\begin{rem}\label{Rem of Prop CS(Gamma|G, phi)(gamma)}
    We have for the appearing spaces a precise description except for $ L^2(\Gamma \bs G, \phi)_{ac}(\gamma) \cap \CS(\Gamma \bs G, \phi) $.
    The space $ L^2(\Gamma \bs G, \phi)_{ds}(\gamma) \cap \c{\CS}(\Gamma \bs G, \phi) $ is described in Proposition \ref{Prop Decomposition of C(Gamma|G) cap L^2(Gamma|G)_ds(gamma)}. The spaces $ L^2(\Gamma \bs G, \phi)_U $ and $ L^2(\Gamma \bs G, \phi)_{res} $ are described in \cite{BO00}.
\end{rem}

\begin{proof}
    Fix $ \gamma \in \hat{K} $. Note first that the right-hand side is contained in the left-hand side. Let now $ f \in \CS(\Gamma \bs G, \phi)(\gamma) $. By the Plancherel decomposition (see \eqref{eq Plancherel decomposition} on p.\pageref{eq Plancherel decomposition}), we can write
    \[
        f = f_1 + f_2 + f_3 + f_4
    \]
    with $ f_1 \in L^2(\Gamma \bs G, \phi)_{ds}(\gamma) $, $ f_2 \in L^2(\Gamma \bs G, \phi)_U(\gamma) $, $ f_3 \in L^2(\Gamma \bs G, \phi)_{res}(\gamma) $ and $ f_4 \in L^2(\Gamma \bs G, \phi)_{ac}(\gamma) $. By Proposition \ref{Prop c_pi(F_pi(f)) for discrete series} and by Proposition \ref{Prop U_Lambda(sigma_lambda, phi) consists of cusp forms}, $ f_1 $ and $ f_2 $ belong to $ \c{\CS}(\Gamma \bs G, \phi) $. As $ f_3 \in \CS(\Gamma \bs G, \phi) $ by Lemma \ref{Lem c_(f, v) in the Schwartz space},
    \[
        f_4 = f - f_1 - f_2 - f_3 \in \CS(\Gamma \bs G, \phi) \ .
    \]
    The proposition follows now from Lemma \ref{Lem V = (S cap V) oplus (S^perp cap V) generalisation}.
\end{proof}

\newpage

\subsection{The dual transform, Eisenstein series, Eisenstein integrals and wave packets}

\subsubsection{The dual transform}\label{sssec:The dual transform}

In this section we define the dual transform and we show how it is related to the constant term.
This allows us then to show that
\[
    \overline{\c{\CS}(\Gamma \bs G, \phi)}(\tau) \cap \CS(\Gamma \bs G, \phi) = \c{\CS}(\Gamma \bs G, \phi)(\tau)
\]
for every $ \tau \in \hat{K} $ (cf. Theorem \ref{Thm a Schwartz fct in the closure of the space of cusp forms is a cusp form}).

\begin{lem}\label{Lem (k in K | gk in G(Omega)) is open and dense}
    The set $ \{ k \in K \mid gk \in G(\Omega) \} $ is open and dense in $K$.
\end{lem}

\begin{proof}
    Let $ g \in G $. Then, $ S := \{ k \in K \mid gk \in G(\Omega) \} = K \cap g^{-1}G(\Omega) $ is open in $K$. Let $ k \in K $. Since $ G(\Omega) $ is dense in $G$, there is a sequence $ (h_j)_j $ in $G(\Omega)$ such that $ h_j $ converges to $ gk $.

    Thus, $ g^{-1} h_j $ converges to $ k $. Hence, $ \kappa(g^{-1} h_j) $ converges to $ k $, too.
    Since
    \[
        g \kappa(g^{-1} h_j) P = h_j P \in \Omega \ ,
    \]
    $ g \kappa(g^{-1} h_j) \in G(\Omega) $. So, $ \kappa(g^{-1} h_j) \in S $. The lemma follows.
\end{proof}

Let $ \tau \in \hat{K} $.
Let us denote the space of functions $ \varphi \in \Cinf(\Gamma \bs G(\Omega)/N \times N \bs G, V_\tau \otimes V_\phi) $ which satisfy
\[
    \varphi(\gamma g m a, a^{-1} m^{-1} h k) = a^{-2\rho} \tau(k)^{-1} \otimes \phi(\gamma) \varphi(g, h)
\]
($\gamma \in \Gamma $, $ g \in G(\Omega) $, $ m \in M $, $ a \in A $, $h \in G $, $k \in K $) and which are compactly supported modulo the $A$-equivariance by $ \Ccinf(\Gamma \bs G(\Omega)/N \times_{M A} N \bs G/K, V_{\Hcal}(\tau, \phi)) $.

Set $ \Hcal = \Gamma \bs G(\Omega)/N \times_{M A} N \bs G/K $.
For $ \varphi \in \Ccinf(\Hcal, V_{\Hcal}(\tau, \phi)) $ and $ g \in G $, set
\[
    R^* \varphi(g) = \int_{\{k \in K \, \mid \, gk \in G(\Omega)\}} \varphi(gk, k^{-1}) \, dk \ .
\]
Then, $ R^* \varphi(gh) = \tau(h)^{-1} R^* \varphi(g) $ for all $ g \in G $ and $ h \in K $.

For measurable functions $ \varphi, \psi $ from $ \Hcal $ to $ V_{\Hcal}(\tau, \phi) $ (such functions are defined analogously as above), define
\[
    (\varphi, \psi) = \int_{\Gamma \bs G(\Omega)/MN} (\varphi(g, e), \psi(g, e)) \, dg \ ,
\]
where $ \Gamma \bs G(\Omega)/MN $ is equipped with its invariant volume form, whenever the integral exists.

\begin{lem}\label{Lem (varphi, f^Omega) = (R^* varphi, f)} \nlenum
    \begin{enumerate}
    \item $ R^* \varphi $ is well-defined and belongs to $ \Cinf(Y, V_Y(\tau, \phi)) $.
    \item Let $ f \in \Ccinf(Y, V_Y(\tau, \phi)) $. Then,
        \[
            (\varphi, f^\Omega) = (R^* \varphi, f) \ .
        \]
    \end{enumerate}
\end{lem}

\begin{proof}
    Let $ f $ be as above. Then,
    \[
        (\varphi, f^\Omega)
        = \int_{\Gamma \bs G(\Omega)/MN} (\varphi(g, e), \int_N f(gn) \, dn) \, dg \ .
    \]
    Since $ \Gamma \bs G(\Omega)/M $ fibers over $ \Gamma \bs G(\Omega)/MN $, this is again equal to
    \[
        \int_{\Gamma \bs G(\Omega)/M} (\varphi(g, e), f(g)) \, dg \ .
    \]
    Since $ \Gamma \bs G(\Omega)/M $ fibers over $ \Gamma \bs G/K $ with fiber $ { \{ k \in K \, \mid \, gk \in G(\Omega) \} } $, this yields
    \begin{align*}
        & \int_{\Gamma \bs G/K} \int_{ \{ k \in K \, \mid \, gk \in G(\Omega) \} } (\varphi(g k, e), \tau(k)^{-1} f(g)) \, dk \, dg \\
        =& \int_{\Gamma \bs G/K} (\int_{ \{ k \in K \, \mid \, gk \in G(\Omega) \} } \varphi(gk, k^{-1}) \, dk, f(g)) \, dg
        = (R^* \varphi, f) \ .
    \end{align*}
    As $ \int_{\Gamma \bs G(\Omega)/M} |(\varphi(g, e), f(g))| \, dg $ is clearly finite for all $ f \in \Ccinf(Y, V_Y(\tau, \phi))$,
    \[
        k \in \{ h \in K \mid gh \in G(\Omega) \} \mapsto \varphi(gk, k^{-1})
    \]
    is integrable for almost all $ g \in G $, by Fubini's theorem and the above computations.
    It follows by the theorem about differentiation of parameter dependent integrals that $ R^* \varphi(g) $ is well-defined for all $ g \in G $ and that $ R^* \varphi \in \Cinf(Y, V_Y(\tau, \phi)) $.
\end{proof}

Let $ \sigma $ be a finite-dimensional representation of $M$. Define $ \sigma_{MN}(m n) = \sigma(m) $.\index[n]{sigma_MN@$\sigma_{MN}$}
We denote $ V_\sigma $ endowed with this action by $ V_{\sigma_{MN}} $. Set $ V_{MN}(\sigma) = G \times_{M N} V_{\sigma_{MN}} $, $ V_{MN}(\sigma, \phi) = V_{MN}(\sigma) \otimes V_\phi $ (carrying the tensor product action of $ \Gamma $) and
\[
    V_{B'}(\sigma, \phi) = \Gamma \bs \restricted{V_{MN}(\sigma, \phi)}{G(\Omega)/MN}
\]
(bundle on $ B' := \Gamma \bs G(\Omega)/MN $). Then, $ \Ccinf(\Gamma \bs G(\Omega)/MN, V_{B'}(\sigma, \phi)) $ is equal to
\begin{multline*}
    \{ \varphi \in \Ccinf(\Gamma \bs G(\Omega)/N, V_\sigma \otimes V_\phi) \mid \\
    \varphi(\gamma g m) = (\sigma(m)^{-1} \otimes \phi(\gamma)) \varphi(g)
    \quad \forall \gamma \in \Gamma , \, g \in G(\Omega) , \, m \in M \} \ .
\end{multline*}
Let $ T \in \Hom_M(V_\sigma, V_\tau) $ ($ \sigma \in \hat{M} $, $ \tau \in \hat{K} $) and let $ \varphi \in \Ccinf(\Gamma \bs G(\Omega)/MN, V_{B'}(\sigma, \phi)) $.
\\ For $ g \in G(\Omega) $, $ n \in N $, $ a \in A $ and $ k \in K $, define
\[
    F(g, n a k) = a^{2\rho} \tau(k)^{-1} T \varphi(g a) \ .
\]
Then, $F$ is well-defined and belongs to $ F \in \Ccinf(\Hcal, V_{\Hcal}(\tau, \phi)) $.

By abuse of notation, we denote $ F $ also by $ T \varphi $.

Let $ e \neq \bar{n} \in \bar{N} $.  Write $ \bar{n} \in \bar{N} $ as $ \bar{n} = \exp(X + Y) \ (X \in \g_{-\alpha},  Y \in \g_{-2\alpha} ) $.

By Proposition \ref{Theorem 3.8 of Helgason, p.414}, $ a_B(w\bar{n}) = \sqrt{\tfrac14|X|^4 + 2 |Y|^2} = \frac{ \sqrt{|X|^4 + 8 |Y|^2} }2 $. Let $ \eps > 0 $. Then,
\[
    a_B(w\bar{n}) \asymp a(\bar{n})
\]
for all $ \bar{n} \in \bar{N} $ such that $ |\log \bar{n}| > \eps $ by Corollary \ref{Cor a_n asymp a(theta n)}.

For $ t > 0 $, set $ a_t = \exp(\log(t) H_\alpha) $.

\begin{lem}\label{Lem a a((bar n)^(a^(-1/2)))^(-1) converges} \nlenum
    \begin{enumerate}
    \item Let $ e \neq \bar{n} \in \bar{N} $ and $ a \in A $. Then,
        \[
            a a(\bar{n}^{a^{-\tfrac12}})^{-1}
                \leq a_B(w\bar{n})^{-1}
            \ .
        \]
    \item Let $ r > 0 $, $ R > 0 $ and let $ \eps > 0 $. Then, there is a constant $ c > 0 $ (depending on $R$ and $ \eps $) such that
        \[
            |a^r a(\bar{n}^{a^{-\tfrac12}})^{-r} - a_B(w\bar{n})^{-r}| \leq \frac{c}a \cdot a(\bar{n})^{-r - 1}
        \]
        for all $ r \in [0, R] $, $ a \in \bar{A}_+ $ and $ \bar{n} \in \bar{N} $ such that $ |\log \bar{n}| > \eps $.
    \item Let $ r > 0 $, $ R > 0 $ and let $ \eps > 0 $. Then, there is a constant $ c > 0 $ (depending on $r$ and $ \eps $) such that
        \[
            |a(n^a w)^{-r} - a_B(n^a w)^{-r}| \leq \frac{c}{a^{2(r+1)}} \cdot a(n w)^{-r - 1}
        \]
        for all $ r \in [0, R] $, $ a \in \bar{A}_+ $ and $ n \in N $ such that $ |\log n| > \eps $.
    \end{enumerate}
\end{lem}

\begin{proof}
    Let $ e \neq \bar{n} \in \bar{N} $ and $ t > 0 $. By Proposition \ref{Theorem 3.8 of Helgason, p.414} (we use the notation preceding the lemma),
    \[
        a_t a(\bar{n}^{a_t^{-\tfrac12}})^{-1} = \frac{t}{ \sqrt{(1 + \tfrac{t}{2}|X|^2)^2 + 2 t^2 |Y|^2} }
        = \frac1{ \sqrt{(\frac1{t} + \tfrac12|X|^2)^2 + 2 |Y|^2} } \ .
    \]
    Thus, the first assertion holds.
    
    Let $ r \geq 0 $, $ t > 0 $ and let $ e \neq \bar{n} \in \bar{N} $. Since
    \begin{multline*}
        \frac{ a^r a(\bar{n}^{a^{-\tfrac12}})^{-r} }{ a_B(w\bar{n})^{-r} } - 1
        = \frac{ \big(\frac14|X|^4 + 2 |Y|^2\big)^{\frac{r}2} }{ \big((\frac1{t} + \tfrac12|X|^2)^2 + 2 |Y|^2\big)^{\frac{r}2} } - 1 \\
        = \frac{ a_B(w\bar{n})^r }{ \big(\frac1{t^2} + \frac1t |X|^2 + a_B(w\bar{n})^2 \big)^{\frac{r}2} } - 1
        = \frac{ 1 }{ \big( (\frac1{t^2} + \frac1t |X|^2) a_B(w\bar{n})^{-2} + 1 \big)^{\frac{r}2} } - 1 \ .
    \end{multline*}
    Let $ t \geq 1 $ and let $ \bar{n} \in \bar{N} $ be such that $ |\log \bar{n}| > \eps $. Then,
    \[
        (\tfrac1{t^2} + \tfrac1t |X|^2) a_B(w\bar{n})^{-2} \leq \frac1t \cdot \frac{ 1 + |X|^2 }{ a_B(w\bar{n})^{2} }
    \]
    is uniformly bounded by a constant $ C > 0 $.
    As $ s \in (-1, \infty) \mapsto \frac{ 1 }{ (s + 1)^{\frac{r}2}} - 1 $ is analytic, there are coefficients $ a_k(r) $ such that
    \[
        \frac{ 1 }{ (s + 1)^{\frac{r}2} } - 1 = \sum_{k=0}^\infty a_k(r) s^k  \ .
    \]
    As the function vanishes at $ s = 0 $, $ a_0(r) = 0 $. Hence, $ \frac1s |\frac{ 1 }{ (s + 1)^{\frac{r}2} } - 1| $ is uniformly bounded when $s$ varies in $ (0, C] $.
    Moreover,
    \[
        0 \leq \frac1s \big(1 - \frac{ 1 }{ (s + 1)^{\frac{r_1}2} }\big) \leq \frac1s \big(1 - \frac{ 1 }{ (s + 1)^{\frac{r_2}2} }\big)
    \]
    for all $ s \in (0, \infty) $ and all $ 0 \leq r_1 \leq r_2 $. Let $ R > 0 $. Then, the above estimate is also uniform when $r$ varies in $ [0, R] $.
    So,
    \[
        |\tfrac{ a^r a(\bar{n}^{a^{-\tfrac12}})^{-r} }{ a_B(w\bar{n})^{-r} } - 1| \prec (\tfrac1{t^2} + \tfrac1t |X|^2) a_B(w\bar{n})^{-2} \prec \tfrac1t a_B(w\bar{n})^{-1}
    \]
    for $ r \in [0, R] $, $ t \geq 1 $ and $ \bar{n} \in \bar{N} $ such that $ |\log \bar{n}| > \eps $.
    The second assertion follows.

    We have
    \begin{multline*}
        a(n^a w)^{-r} - a_B(n^a w)^{-r} = a(n^a w)^{-r} - a^{-r} a_B(n^{a^{\frac12}} w)^{-r} \\
        = a^{-r} \big(a^r a( (w n^{a^{\frac12}} w)^{a^{-\frac12}})^{-r} - a_B(w (w n^{a^{\frac12}} w))^{-r}\big) \ .
    \end{multline*}
    It follows from the previous assertion that there is $ c > 0 $, depending only on $R$, such that this is again less or equal than
    \[
        \frac{c \, a^{-r}}a \cdot a_B(n^{a^{\frac12}} w)^{-r - 1} \leq \frac{c}{a^{2(r+1)}} \cdot a_B(n w)^{-r - 1}
    \]
    for all $ r \in [0, R] $, $ a \in \bar{A}_+ $ and $ n \in N $ such that $ |\log n| > \eps $. The third assertion follows.
\end{proof}

Let $ N_0 A_0 K_0 $ be the standard Iwasawa decomposition of $ \SU(1, 2) $.
Write $ g = \kappa_0(g) a_0(g) n_0(g) $ ($ g \in \SU(1, 2) $) with $ \kappa_0(g) \in K_0 $, $ a_0(g) \in A_0 $ and $ n_0(g) \in N_0 $.
Set
\[
    w_0 = \begin{pmatrix} 1&0&0 \\0&-1&0 \\0&0&-1 \end{pmatrix} \in K_0 \cap \SO(1, 2)^0 \ .
\]
Let $ N_0^1 = \{ n_{v, 0} \mid v \in \RR \} $ and $ N_0^2 = \{ n_{0, r} \mid r \in \RR \} $ (we use here the notations of Appendix \ref{sec:Sp(1,n)}).

\begin{lem}\label{Lem for SU(1, 2)-reduction}  \nlenum
    \begin{enumerate}
    \item Let $ X_1 \in \g_\alpha $ be nonzero. Then, there is a Lie group homomorphism
        \[
            \Phi \colon \Spin(1, 2) \to G
        \]
        respecting the Iwasawa decompositions such that $ \exp(X_1) \in \im \Phi $.
    \item Assume that $ \g_{2\alpha} \neq \{0\} $. Let $ X_1 \in \g_\alpha $ and $ X_2 \in \g_{2\alpha} $ be nonzero. Then, there is a Lie group homomorphism $ \Phi \colon \SU(1, 2) \to G $ respecting the Iwasawa decompositions such that $ \exp(X_1) \in \Phi(N_0^1) $, $ \exp(X_2) \in \Phi(N_0^2) $.
    \end{enumerate}
\end{lem}

\begin{rem}
    Let $ \Phi $ be as above. Then, a standard representative of the Weyl element of $ \Spin(1, 2) $ (resp. $ \SU(1, 2) $) provides a representative $ w \in K $ of the Weyl element of $G$ satisfying $ w^2 \in Z(G) $ and $ w^4 = e $ when $ \Phi $ is provided by the first assertion (resp. $ w^2 = e $ when $ \Phi $ is provided by the second assertion).
\end{rem}

\begin{proof}
    Assume that $ \g_{2\alpha} \neq \{0\} $. Let $ X_1 \in \g_\alpha $ and $ X_2 \in \g_{2\alpha} $ be nonzero. Let $ \theta_0 $ be the Cartan involution on $ \su(1, 2) $ associated to $ \Lie(K_0) $.
    It follows from Theorem 3.1, Chapter IX, of \cite[p.409]{Helgason} ($ \SU(2, 1) $-reduction) and Lemma 3.7, Chapter IX, of \cite[p.413]{Helgason} that the Lie subalgebra $ \g^* \subset \g $ generated by $ X_1 $, $ X_2 $, $ \theta(X_1) $, $ \theta(X_2) $ is isomorphic to $ \su(1, 2) $ via a map $ \varphi $ preserving Iwasawa decompositions (compare with the proof of Theorem \ref{Thm explicit formulas}).
    Let $ G^* $ be the analytic subgroup of $G$ with Lie algebra $ \g^* $. Since $ G^* \subset G $ is linear, there is a Lie group homomorphism
    \[
        \Phi \colon \SU(1, 2) \to G
    \]
    respecting the Iwasawa decompositions such that $ \exp(X_1) \in \Phi(N_0^1) $, $ \exp(X_2) \in \Phi(N_0^2) $.
    The first assertion is proven similarly. One may use that $ \Spin(1, 2) $ is isomorphic to $ \SL(2, \RR) $.
\end{proof}

Let us continue with some preparations for Proposition \ref{Prop lim_(a to infty) a^rho (R^* varphi)(ka) converges}.

\begin{lem}\label{Lem kappa_0(n^( a_(1/u) ) w_0) is analytic}
    Let $ G = \SU(1, 2) $. Let $ n \in N_0 $. Then, $ \lim_{a \to \infty} \kappa_0(n^a w_0) $ exists and belongs to $M_0 := Z_{K_0}(A_0)$. Moreover,
    \[
        (0, \infty) \times (N_0 \smallsetminus \{e\}) \to K_0 , \quad (u, n) \mapsto \kappa_0(n^{ a_{\frac1u} } w_0)
    \]
    extends analytically to $ \RR \times (N_0 \smallsetminus \{e\}) $.
\end{lem}

\begin{proof}
    Let $ G = \SU(1, 2) $. Let $ e \neq n_{v, r} \in N_0 $.
    Direct computation shows that
    \[
        \kappa_0(n_{v, r} w_0) =
        \begin{pmatrix}
            \frac{ v^2 + 1 + 2 i r }{ \sqrt{ (v^2 + 1)^2 + 4 r^2 } } & 0 & 0 \\
            0 & \frac{ v^2 - 1 + 2 i r }{ \sqrt{ (v^2 + 1)^2 + 4 r^2 } } & -\frac{ 2 v }{ v^2 + 1 + 2 i r } \\
            0 & \frac{ 2 v }{ \sqrt{ (v^2 + 1)^2 + 4 r^2 } } & \frac{ v^2 - 1 - 2 i r }{ v^2 + 1 + 2 i r }
        \end{pmatrix} \in K_0 \ .
    \]
    Thus,
    \[
        \kappa_0(n_{v, r}^{a_{\frac1u}} w_0) = \begin{pmatrix}
            \frac{ v^2 + u^2 + 2 i r }{ \sqrt{ (v^2 + u^2)^2 + 4 r^2 } } & 0 & 0 \\
            0 & \frac{ v^2 - u^2 + 2 i r }{ \sqrt{ (v^2 + u^2)^2 + 4 r^2 } } & -\frac{ 2 u v }{ v^2 + u^2 + 2 i r } \\
            0 & \frac{ 2 u v }{ \sqrt{ (v^2 + u^2)^2 + 4 r^2 } } & \frac{ v^2 - u^2 - 2 i r }{ v^2 + u^2 + 2 i r }
        \end{pmatrix}
    \]
    (here $ n^a := a n a^{-1} $) and
    \[
        \lim_{a \to \infty} \kappa_0(n_{v, r}^a w_0) = \lim_{u \to 0^+} \kappa_0(n_{v, r}^{ a_{\frac1u} } w_0)
        = \begin{pmatrix}
            \frac{ v^2 + 2 i r }{ \sqrt{v^4 + 4 r^2} } & 0 & 0 \\
            0 & \frac{ v^2 + 2 i r }{ \sqrt{v^4 + 4 r^2} } & 0 \\
            0 & 0 & \frac{ v^2 - 2 i r }{ v^2 + 2 i r }
        \end{pmatrix} \in M_0 \ .
    \]
    \underline{Remark:} This matrix is equal to the identity when $ r = 0 $.
    
    The lemma follows.
\end{proof}

Set
\[
    N_1 = \{ n_x := \begin{pmatrix} 1 & x \\ 0 & 1 \end{pmatrix} \mid x \in \RR \} , \qquad
    A_1 = \{ a_t := \begin{pmatrix} \sqrt{t} & 0 \\ 0 & \frac1{\sqrt{t}} \end{pmatrix} \mid t > 0 \}
\]
and
\[
    K_1 = \big\{ \begin{pmatrix} \cos \phi & \sin \phi \\ -\sin \phi & \cos \phi \end{pmatrix} \mid \phi \in [0, 2\pi] \big\} \ .
\]
Then, $ N_1 A_1 K_1 $ is an Iwasawa decomposition of $ \SL(2, \RR) $. This induces an Iwasawa decomposition $ N_2 A_2 K_2 $ of $ \Spin(1, 2) $.

Put $ w_1 = \begin{pmatrix} 0 & 1 \\ -1 & 0 \end{pmatrix} \in K_1 $.

\begin{lem}\label{Lem kappa_0(n^( a_(1/u) ) w_0) is analytic for SL(2, R)}
    Let $ G = \SL(2, \RR) \simeq \Spin(1, 2) $. Let $ n \in N_1 $. Then, $ \lim_{a \to \infty} \kappa(n^a w_1) $ exists and belongs to $M_1 := Z_{K_1}(A_1)$. Moreover,
    \[
        (0, \infty) \times (N_1 \smallsetminus \{e\}) \to K_1 , \quad (u, n) \mapsto \kappa(n^{ a_{\frac1u} } w_1)
    \]
    extends analytically to $ \RR \times (N_1 \smallsetminus \{e\}) $.
\end{lem}

\begin{proof}
    Let $ G = \SL(2, \RR) $.
    Direct computation shows that
    \[
        \kappa(n_x w_1) =
        \begin{pmatrix}
            -\frac{ x }{ \sqrt{x^2 + 1} } & \frac{ 1 }{ \sqrt{x^2 + 1} } \\
            -\frac{ 1 }{ \sqrt{x^2 + 1} } & -\frac{ x }{ \sqrt{x^2 + 1} }
        \end{pmatrix} \in K_1 \ .
    \]
    Thus,
    \[
        \kappa(n_x^{a_{\frac1u}} w_0) = \begin{pmatrix}
            -\frac{ x }{ \sqrt{x^2 + u^2} } & \frac{ u }{ \sqrt{x^2 + u^2} } \\
            -\frac{ u }{ \sqrt{x^2 + u^2} } & -\frac{ x }{ \sqrt{x^2 + u^2} }
        \end{pmatrix}
    \]
    and
    \[
        \lim_{a \to \infty} \kappa(n_x^a w_1) = \lim_{u \to 0^+} \kappa(n_x^{ a_{\frac1u} } w_1)
        = \begin{pmatrix}
            -1 & 0\\
            0 & -1
        \end{pmatrix} \in M_1 \ .
    \]
    The lemma follows.
\end{proof}

Let $ X_1 \in \g_\alpha $ be nonzero and let $ X_2 \in \g_{2\alpha} $. If $ \g_{2\alpha} $ is nonzero, then we choose also $ X_2 \neq 0 $.
Let $ H = \Spin(1, 2) \simeq \SL(2, \RR) $ if $ X_2 = 0 $ and let $ H = \SU(1, 2) $ otherwise.
Let $ \Phi \colon H \to G $ the Lie group homomorphism provided by Lemma \ref{Lem for SU(1, 2)-reduction}.
Let $ N' = \Phi(N_2 \cap H) $ if $ H = \Spin(1, 2) $ and let $ N' = \Phi(N_0 \cap H) $ if $ H = \SU(1, 2) $.

\begin{lem}\label{Lem cup_(m in M) m N' m^(-1) = N}
    We have
    \[
        \bigcup_{m \in M} m N' m^{-1} = N \ .
    \]
\end{lem}

\begin{proof}
    If $ X = \OO H^2 $, then $ \g_\alpha = \OO $ and $ \g_{2\alpha} = \I(\OO) $ and $ M = \Spin(7) $.

    By Lemma 8.8 of \cite[p.30]{Wolf}, $ (\restricted{\Ad}{M}, \OO) $ is isomorphic to the real 8-dimensional spin representation $ \Delta $ and $ (\restricted{\Ad}{M}, \I(\OO)) $ is isomorphic to the standard representation on $ \RR^7 $.

    By Corollary 5.4 of \cite[p.32]{Adams}, $ M = \Spin(7) $ acts transitively on
    \[
        \{ X + Y \mid X \in \g_\alpha, Y \in \g_{2\alpha} : |X| = |Y| = 1 \} \ .
    \]
    Let $ r_1 > 0 $ and $ r_2 > 0 $. By the proof of the mentioned corollary, $ M $ acts also transitively on
    \[
        \{ X + Y \mid X \in \g_\alpha, Y \in \g_{2\alpha} : |X| = r_1 , \, |Y| = r_2 \} \ .
    \]
    Thus,
    \[
        \bigcup_{m \in M} \Ad(m) \n' = \n \ ,
    \]
    where $ \n' = \Lie(N') $. The lemma follows in the exceptional case.

    If $ X \neq \OO H^2 $, then we may assume that $ G $ is equal to $ \SO(1, n)^0 $ (resp. $ \SU(1, n) $, $ \Sp(1, n) $) for some $ n \geq 2 $. Indeed, let $ \Psi \colon \tilde{G} \to G $ be a covering of $G$. Then, $ \Psi_* = \Id $ and $ \tilde{G} = N A \tilde{K} $ as $ N $ and $A$ are simply connected. Moreover,
    \[
        \Ad(m)X = \Psi_*(\Ad(m)X) = \Ad(\Psi(m))\Psi_*(X) = \Ad(\Psi(m))(X)
    \]
    for all $ m \in \tilde{M} $ and $ X \in \n $.
    
    Let $ \FF = \RR $ (resp. $ \FF = \CC $, $ \FF = \HH $).
    If $ \FF = \RR $, then
    \[
        m_A n_v m_A^{-1} = n_{Av}
    \]
    for all $ A \in \SO(n-1) $ and $ v \in \RR^{n-1} $, by \eqref{eq m_A n_v m_A^(-1) = n_(Av)} of Appendix \ref{sec:Sp(1,n)}. Thus, $ M = \SO(n-1) $ acts transitively on spheres of $ \n = \RR^{n-1} $ in this case.

    Let us consider now the two other cases:
    For all $ v \in \FF^{n-1} $, $ r \in \I(\FF) $ ($ \I(\CC) = i\RR $ by convention), $ A \in \SU(n-1) $ (resp. $ \Sp(n-1) $), $ q \in \SU(1) $ (resp. $ q \in \Sp(1) $),
    \[
        m_{A, q} n_{v, r} m_{A, q}^{-1} = n_{Av \bar{q}, q r \bar{q}} = n_{Av \bar{q}, q r q^{-1}} = n_{Av \bar{q}, \Ad(q)r}
    \]
    by \eqref{eq m_(A, q) n_(v, r) m_(a, q)^(-1) = n_(Av bar(q), q r bar(q))} of Appendix \ref{sec:Sp(1,n)}.
    For $ R > 0 $, set $ S_1(0, R) = \{v \in \FF^{n-1} \mid \|v\| = R\} $ and $ S_2(0, R) = \{r \in \I(\FF) \mid |r| = R\} $.
    If $ \FF = \CC $, then $ (q, r) \in \SU(1) \times S_2(0, R) \mapsto q r \bar{q} \in S_2(0, R) $ is transitive.
    As moreover
    \[
        (A, v) \in \SU(n-1) \times S_1(0, R) \mapsto Av \in S_1(0, R)
    \]
    is also transitive, $M$ acts transitive on $ S_1(0, R) \times S_2(0, R) $.

    Similarly, one shows that $ M $ acts  transitively on $ S_1(0, R) \times S_2(0, R) $ if $ \FF = \HH $.
    Thus,
    \[
        \bigcup_{m \in M} m N' m^{-1} = N \ .
    \]
\end{proof}

For $ x \in N $, set $ m_x = \lim_{a \to \infty} \kappa(x^a w) $. It follows from the previous corollary that $ x \in N \smallsetminus \{e\} \mapsto m_x \in M $ is a well-defined function. Note that $ m_x $ has the following two properties:
\begin{enumerate}
\item $ m_{x^a} = m_x $ for all $ x \in N $ and $ a \in A $.
\item $ m_{m_1 x m_2} = m_1 m_x w^{-1} m_2 w $ for all $ x \in N $ and $ m_1, m_2 \in M $.
\end{enumerate}

\begin{cor}\label{Cor t |tau(kappa(n^(a_t) w)) - tau(m_n)| is uniformly bounded}
    Let $ \tau \in \hat{K} $ and $ C $ be a compact subset of $ N \smallsetminus \{e\} $. Then,
    \[
        t |\tau(\kappa(n^{a_t} w)) - \tau(m_n)|
    \]
    is uniformly bounded when $ t \geq 1 $ and $ n \in C $.
\end{cor}

\begin{proof}
    Let $ \tilde{w} $ be the Weyl element provided by Lemma \ref{Lem for SU(1, 2)-reduction}.
    Let $ m \in M $ be such that $ \tilde{w} = w m $. Since moreover $ \kappa(n w) = \kappa(n w m) m^{-1} = \kappa(n \tilde{w}) m^{-1} $ and $ |\tau(m^{-1})| = 1 $, we may assume without loss of generality that $ w = \tilde{w} $.

    For $ G = \SU(1, 2) $, the corollary follows now immediately from Lemma \ref{Lem kappa_0(n^( a_(1/u) ) w_0) is analytic} and the fact that a finite-dimensional representation on an arbitrary Lie group is real analytic.

    For the general case, we need also Lemma \ref{Lem cup_(m in M) m N' m^(-1) = N}, Lemma \ref{Lem kappa_0(n^( a_(1/u) ) w_0) is analytic for SL(2, R)} and the fact that
    \[
        \tau(\kappa((m n m^{-1})^a w)) - \tau(m_{m n m^{-1}}) = \tau(m) (\tau(\kappa(n^a w)) - \tau(m_n)) \tau(w^{-1} m w)
    \]
    for all $ n \in N $ and $ m \in M $.
\end{proof}

\begin{lem}\label{Lem Res_0 int_A a(x^a w)^(-lambda - rho) a^(2 rho) da = 1/2}
    $ \Residue_0 \int_A a(x^a w)^{-\lambda - \rho} a^{2\rho} \, da $ is equal to $ \Residue_0 \int_A a^{-2\lambda} \, da = \frac12 $ for every $ x \in N $ with $ a_B(x w) = 1 $.
\end{lem}

\begin{proof}
    Let $ r \in \RR $. Then, $ \Residue_0 \int_A a(x^a w)^{-\lambda - \rho} a^{2\rho} \, da = \Residue_0 \int_{A_{\geq e^r}} a(x^a w)^{-\lambda - \rho} a^{2\rho} \, da $.

    By Lemma \ref{Lem a a((bar n)^(a^(-1/2)))^(-1) converges}, there is a constant $ c > 0 $ such that
    \[
        |a(x^a w)^{-(\lambda + \rho)} - a_B(x^a w)^{-(\lambda + \rho)}| \leq \frac{c}{a^{2(\lambda + \rho + \alpha)}}
    \]
    for all $ \lambda \in [0, 2] $, $ a \in \bar{A}_+ $ and $ n \in N $ such that $ a_B(x w) = 1 $. But
    \[
        \Residue_0 \int_{A_{\geq e^r}} a^{-2(\lambda + \alpha)} \, da = \lim_{\lambda \to 0^+} \frac{ \lambda }{2(\lambda + 1)}e^{-2(\lambda + 1) r} = 0 \ .
    \]
    So,
    \[
        \Residue_0 \int_A a(x^a w)^{-\lambda - \rho} a^{2\rho} \, da = \Residue_0 \int_{A_{\geq e^r}} a_B(x^a w)^{-\lambda - \rho} a^{2\rho} \, da
        = \Residue_0 \int_A a^{-2\lambda} \, da = \frac12
    \]
    for every $ x \in N $ with $ a_B(x w) = 1 $. The lemma follows.
\end{proof}

Recall that we denote the Harish-Chandra $c$-function by $ c_\gamma(\mu) \in \End_M(V_\gamma) $ ($ \gamma \in \hat{K} $). This meromorphic function is given by
\[
    c_\gamma(\mu) := \int_{\bar{N}}{ a(\bar{n})^{-(\mu+\rho)} \gamma(\kappa(\bar{n})) \, d\bar{n} }
\]
for $ \R(\mu) > 0 $.

\begin{prop}\label{Prop Other expression for Res_0 c_tau(lambda) o T}
    Let $ T \in \Hom_M(V_\sigma, V_\tau) $ ($ \sigma \in \hat{M} $, $ \tau \in \hat{K} $). Then, $ \Residue_0 c_\tau(\lambda) \circ T $ is equal to $ \frac12 \tau(w) \int_{\{x \in N \, \mid \, a_B(x w) = 1\}} \tau(m_x) \, dx \circ T $.
\end{prop}

\begin{proof}
    Let $ \eps > 0 $. Let $ r \in \RR $ be sufficiently large so that $ |\tau(\kappa(x^a w)) - \tau(m_x)| < \eps $ for every $ a \in A_{e^r} $ and $ x \in N $ with $ a_B(x w) = 1 $.
    By Lemma \ref{Lem int_E f(n) dn = int_A int_(varphi(x) = 1) f(x^a) dx a^(2 xi_E rho) da} applied to $ E = N $,
    \begin{multline*}
        \int_{\bar{N}} a(\bar{n})^{-(\lambda + \rho)} \tau(\kappa(\bar{n})) \, d\bar{n} \circ T
        = \tau(w) \int_N a(n w)^{-(\lambda + \rho)} \tau(\kappa(n w)) \, dn \circ T \\
        = \tau(w) \int_{\{x \in N \, \mid \, a_B(x w) = 1\}} \int_A a(x^a w)^{-(\lambda + \rho)} \tau(\kappa(x^a w)) a^{2\rho} \, da \, dx \circ T
    \end{multline*}
    for every $ \lambda > 0 $. Thus, $ \Residue_0 c_\tau(\lambda) \circ T $ is equal to
    \[
        \lim_{\lambda \to 0^+} \lambda \tau(w) \int_{\{x \in N \, \mid \, a_B(x w) = 1\}}
            \int_{A_{\geq e^r}} a(x^a w)^{-(\lambda + \rho)} \tau(\kappa(x^a w)) a^{2\rho} \, da \, dx \circ T \ .
    \]
    Let $ M_r = \sup_{x \in N \, : \, a_B(x w) = 1, \, a \in A_{\geq e^r}, \, \lambda \in [0, 2]} \frac{a(x^a w)^{-\lambda-\rho}}{a_B(x^a w)^{-\lambda-\rho}} < \infty $.
    Then,
    \[
        \lambda a(x^a w)^{-(\lambda + \rho)} a^{2\rho} \leq \lambda M_r a_B(x^a w)^{-\lambda-\rho} a^{2\rho} = \lambda M_r a^{-2\lambda}
    \]
    for every $ x \in N $ with $ a_B(x w) = 1 $, $ a \in A_{\geq e^r} $ and $ \lambda \in (0, 2] $. This is integrable on $ \{x \in N \, \mid \, a_B(x w) = 1\} \times A_{\geq e^r} $.
    Since moreover $ \int_{A_{\geq e^r}} \lambda a^{-2\lambda} \, da = \frac12 e^{-2\lambda r} \leq \frac12 $ for every $ \lambda \in (0, 2] $ and since we can choose $ \eps > 0 $ arbitrarily small, $ \Residue_0 c_\tau(\lambda) \circ T $ is equal to
    \[
        \lim_{\lambda \to 0^+} \lambda \tau(w) \int_{\{x \in N \, \mid \, a_B(x w) = 1\}} \tau(m_x) \int_A a(x^a w)^{-(\lambda + \rho)} a^{2\rho} \, da \, dx \circ T \ .
    \]
    It follows from the proof of Lemma \ref{Lem Res_0 int_A a(x^a w)^(-lambda - rho) a^(2 rho) da = 1/2} that this is again equal to
    \[
        \frac12 \tau(w) \int_{\{x \in N \, \mid \, a_B(x w) = 1\}} \tau(m_x) \, dx \circ T \ .
    \]
    The proposition follows.
\end{proof}

\begin{prop}\label{Prop Res_0 c_tau(lambda) T injective if T is nonzero and Res_0 c_tau(lambda) is nonzero on V_tau(sigma)}
    Let $ T \in \Hom_M(V_\sigma, V_\tau) $ ($ \sigma \in \hat{M} $, $ \tau \in \hat{K} $). Then:
    \begin{enumerate}
    \item $ \Residue_0 c_\tau(\lambda) T $ is up to a constant equal to $ T^w $. Hence, $ \Residue_0 c_\tau(\lambda) T $ is injective if and only if $ T $ is nonzero and $ \Residue_0 c_\tau(\lambda) $ is nonzero on $ V_\tau(\sigma) $.
    \item $ \Residue_0 c_\tau(\lambda) T = 0 \iff T = 0 \text{ or } p_\sigma(0) \neq 0 $.
    \end{enumerate}
\end{prop}

\begin{proof} \nlenum
    \begin{enumerate}
    \item By Theorem 14.16 of \cite[p.541]{Knapp}, $ \C^{-\infty}(\dX, V(\sigma_0)) $ is irreducible if and only if
        \begin{enumerate}
        \item $ \sigma $ is not Weyl-invariant, or
        \item $ \sigma $ is Weyl-invariant and $ p_\sigma(0) = 0 $.
        \end{enumerate}
        If $ p_\sigma(0) \neq 0 $, then $ \hat{J}_{\sigma, 0} $ is regular by Proposition 7.3, 3., of \cite[p.123]{BO00}.

        Since moreover $ \Residue_0 \hat{J}_{\sigma, \lambda} $ is an intertwining operator, it follows from Schur's lemma and the above that
        \[
            \Residue_0 \hat{J}_{\sigma, \lambda} = c \Id
        \]
        for some constant $ c \in \CC $. 
        Let $ T \in \Hom_M(V_\sigma, V_\tau) $ and $ f \in \C^{-\infty}(\dX, V(\sigma_\lambda)) $. Then,
        \[
            \int_K \tau(k) T \hat{J}_{\sigma, \lambda} f(k) \, dk = \int_K \tau(k) \tau(w) c_\tau(\lambda) T \sigma(w)^{-1} f(k) \, dk
        \]
        by Lemma 5.5, 1., of \cite[p.100]{BO00}.
        Thus,
        \[
            c \int_K \tau(k) T f(k) \, dk = \int_K \tau(k) \tau(w) \Residue_0 c_\tau(\lambda) T \sigma(w)^{-1} f(k) \, dk \ .
        \]
        Taking $ f = \delta_e v $ ($ v \in V_\sigma $) on $K$ yields
        \[
            c \, T^w = c \, \tau(w) T \sigma(w)^{-1} = \Residue_0 c_\tau(\lambda) T \ .
        \]
        Thus, $ c $ is nonzero if and only if $ \Residue_0 c_\tau(\lambda) $ is nonzero on $ V_\tau(\sigma) $.

        Assume that $ T \neq 0 $. Then, it follows from Schur's lemma that $ T $ is injective. Hence, $ T^w $ is also injective.
        The first assertion follows.
    \item Let $ \mu \in \a^*_\CC $. By $(20)$ of \cite[p.100]{BO00},
        \[
            c_\tau(-\mu)^w c_\tau(\mu) T = \frac1{p_\sigma(\mu)} T \ .
        \]
        By $ (21) $ of  \cite[p.100]{BO00}, $ c_\tau(\mu)^* = c_\tau(\bar{\mu})^w $.

        Let $ \lambda \in \a^*_+ $. Then, it follows from the above that
        \[
            c_\tau(i\lambda)^* c_\tau(i\lambda) T = \frac1{p_\sigma(i\lambda)} T \ .
        \]
        Since $ c_\tau(\mu) $ is meromorphic and has poles of order at most 1, there exist $ F_i \in \End_M(V_\tau) $ such that
        \[
            c_\tau(i\lambda) = \frac1{i\lambda} F_{-1} + F_0 + i\lambda F_1 + \dotsb
        \]
        on $ V_\tau(\sigma) $, for all $ \lambda $ in a punctured neighbourhood of $0$ in $ \a^* $.
        Assume that $ T \neq 0 $. Let $ v \in V_\sigma $ be such that $ T v \neq 0 $. Then, $ T^w v $ is also nonzero as $ T^w $ is injective.
        We have
        \[
            \|c_\tau(i\lambda) T v\|^2 = \frac1{p_\sigma(i\lambda)} \|T v\|^2 \ .
        \]
        The left-hand side is again equal to
        \begin{multline*}
            \frac1{\lambda^2} \|F_{-1}(T v)\|^2 + \frac1{i\lambda} \big((F_{-1}(T v), F_0(T v)) - (F_0(T v), F_{-1}(T v))\big) \\
            + \big(\|F_0(T v)\|^2 - (F_{-1}(T v), F_1(T v)) - (F_1(T v), F_{-1}(T v)\big) + O(\lambda) \ .
        \end{multline*}
        Hence, if $ \Residue_{\mu = 0} c_\tau = F_{-1} = 0 $, then
        \[
            \frac1{p_\sigma(i\lambda)} \|T v\|^2 = \|F_0(T v)\|^2 + O(\lambda) \ .
        \]
        So, $ p_\sigma(0) \neq 0 $. If $ \Residue_{\mu = 0} c_\tau = F_{-1} \neq 0 $, then $ p_\sigma(0) = 0 $.
        Consequently, $ \Residue_{\mu = 0} c_\tau T = 0 $ if and only if $ p_\sigma(0) \neq 0 $.
    \end{enumerate}
    The proposition follows.
\end{proof}

For $ \sigma \in \hat{M} $, set
\begin{multline*}
    W_{\sigma} = \{ \varphi \in \Ccinf(\Gamma \bs G(\Omega)/MN, V_{B'}(\sigma, \phi)) \mid \\
         p_\sigma(0) \neq 0 \quad \text{or} \quad \int_A \varphi(k a) a^\rho \, da = 0 \quad \forall k \in K(\Omega) \} \ .
\end{multline*}
For a finite-dimensional representation $ \sigma $ of $M$, set
\[
    \tilde{W}_\sigma = \{ \varphi \in \Ccinf(\Gamma \bs G(\Omega)/MN, V_{B'}(\sigma, \phi)) \mid \int_A \varphi(k a) a^\rho \, da = 0 \quad \forall k \in K(\Omega) \} \ .
\]
Let $ 0 \neq T \in \Hom_M(V_\sigma, V_\tau) $ ($ \sigma \in \hat{M} $, $ \tau \in \hat{K} $). Then, $ \varphi \in W_{\sigma} $ if and only if
\[
    \Residue_0 c_\tau(\lambda) \int_A T \varphi(k a, e) a^\rho \, da = 0
\]
for all $ k \in K(\Omega) $, by Proposition \ref{Prop Res_0 c_tau(lambda) T injective if T is nonzero and Res_0 c_tau(lambda) is nonzero on V_tau(sigma)}. By the Peter-Weyl theorem,
\begin{equation} \label{eq Decomposition of Ccinf(Gamma|G(Omega)/N x (M A) N|G/K, V_B'(tau, phi))}
    \Ccinf(\Hcal, V_{\Hcal}(\tau, \phi))
    = \bigoplus_{\sigma \in \hat{M}} \Hom_M(V_\sigma, V_\tau) \otimes \Ccinf(\Gamma \bs G(\Omega)/MN, V_{B'}(\sigma, \phi)) \ .
\end{equation}
Let
\[
    W_\tau = \bigoplus_{\sigma \in \hat{M}} \Hom_M(V_\sigma, V_\tau) \otimes W_{\sigma}
\]
and
\[
    \tilde{W}_\tau = \{ \varphi \in \Ccinf(\Hcal, V_{\Hcal}(\tau, \phi)) \mid \int_A \varphi(k a, e) a^\rho \, da = 0 \quad \forall k \in K(\Omega) \} \ .
\]
Then, $ \tilde{W}_\tau $ is contained in $ W_\tau $. Let $ \sigma \in \hat{M} $. Let $ T_i $ be a basis of $ \Hom_M(V_\sigma, V_\tau) $. Since the sum of $ \I(T_i) $ over $i$ is direct, $ \tilde{W}_\tau $ is equal to $ \bigoplus_{\sigma \in \hat{M}} \Hom_M(V_\sigma, V_\tau) \otimes \tilde{W}_\sigma $.

\begin{prop}\label{Prop lim_(a to infty) a^rho (R^* varphi)(ka) converges}
    Let $ T \in \Hom_M(V_\sigma, V_\tau) $ $ (\sigma \in \hat{M}, \tau \in \hat{K}) $.
    \begin{enumerate}
    \item Let $ k \in K(\Omega) $ and $ \varphi \in \Ccinf(\Gamma \bs G(\Omega)/MN, V_{B'}(\sigma, \phi)) $. Then, \\
        $ \lim_{t \to \infty} a_t^\rho (R^* (T \varphi))(ka_t) $ converges to
        \[
            \int_N T \varphi\big(k a_B(y w)^{-1}, m_y^{-1} \big) a_B(y w)^{-2\rho} \, dy = \Residue_0 c_\tau(\lambda) \int_A T \varphi(k a, e) a^\rho \, da \ .
        \]
    \item Let $ V $ be a relatively compact subset of $ K(\Omega) $ and let $ \varphi \in W_{\sigma} $. Then,
        \[
            a^{\rho + \frac{\alpha}2} (R^* (T \varphi))(k a)
        \]
        is uniformly bounded when $ a \in \bar{A}_+ $ and $ k \in V $.
    \end{enumerate}
\end{prop}

\begin{proof}
    Let $ k \in K(\Omega) $, $ t \geq 1 $, let $ T $ and $ \varphi $ be as above. We have
    \begin{multline*}
        a_t^\rho (R^* (T \varphi))(ka_t)
        = a_t^\rho \int_{ \{ h \in K \, \mid \, k a_t h \in G(\Omega) \} } T \varphi(k a_t h, h^{-1}) \, dh \\
        = a_t^\rho \int_{ \{ \bar{n} \in \bar{N} \, \mid \, k a_t \kappa(\bar{n}) \in G(\Omega) \} } T \varphi\big(ka_t \kappa(\bar{n}), \kappa(\bar{n})^{-1}\big) a(\bar{n})^{-2\rho} \, d\bar{n} \ .
    \end{multline*}
    Since $ \kappa(\bar{n}) = \bar{n} n(\bar{n})^{-1} a(\bar{n})^{-1} = \bar{n} a(\bar{n})^{-1} (a(\bar{n}) n(\bar{n})^{-1} a(\bar{n})^{-1}) $, this is again equal to
    \[
        a_t^\rho \int_{ \{ \bar{n} \in \bar{N} \, \mid \, k a_t \kappa(\bar{n}) \in G(\Omega) \} }
        T \varphi\big(ka_t \bar{n} a(\bar{n})^{-1}, \kappa(\bar{n})^{-1} \big) a(\bar{n})^{-2\rho} \, d\bar{n} \ .
    \]
    By the change of variables $ \bar{n} = y^{a_t^{-\tfrac12}} $ ($ d\bar{n} = a_t^\rho dy $), this yields
    \begin{equation} \label{Proof of Lem lim_(a to infty) a^rho (R^* varphi)(ka) converges eq1}
        \int_{ \{ y \in \bar{N} \, \mid \, k a_t \kappa(y^{a_t^{-\tfrac12}}) \in G(\Omega) \} } T \varphi\big(k y^{a_t^{\tfrac12}} a_t a(y^{a_t^{-\tfrac12}})^{-1}, \kappa(y^{a_t^{-\tfrac12}})^{-1} \big) a_t^{2\rho} a(y^{a_t^{-\tfrac12}})^{-2\rho} \, dy \ .
    \end{equation}
    Let $ V_k $ be an open relatively compact neighbourhood of $ kM $ in $ \Omega $. Let
    \[
        U_k = \{ h a K \mid hM \in V_k, a \in A_+ \} \in \U_\Gamma \ .
    \]
    As $ \varphi $ has compact support in $ \Gamma \bs G(\Omega)/N $, as $ \{\gamma \in \Gamma \mid \gamma U_k \cap U_k \neq \emptyset \} $ is finite and as there is $ s > 0 $ such that $ h a n K \in U_k $ for all $ hM \in V_k $, $ a \in A_{\geq s} $ and $ n \in N $, there is $ s' > 0 $ (depending only on $ V_k $) such that $ \varphi(h a n) = 0 $ for all $ h \in V_k $, $ a \in A_{\geq s'} $ and $ n \in N $.

    Thus, $ \varphi $ is bounded on $ \{ k_1 a k_2 \mid k_1M \in V_k, a \in \bar{A}_+, k_2 \in K \} $. When $ a $ tends to $ \infty $, then $ y^{a^{\tfrac12}} $ converges pointwise to $ e \in \bar{N} $.

    Let $ \eps > 0 $. Choose $ s'' \geq s' > 0 $ sufficiently large so that we also have $ k \kappa(y^{a^{\tfrac12}})M \in V_k $ for all $ y \in \bar{N} $ such that $ |\log y| \leq \eps $  and all $ a \in A_{\geq s''} $.

    Let $ y \in \bar{N} $ be such that $ |\log y| \leq \eps $.
    As $ a(\bar{n}) \geq 1 $ for all $ \bar{n} \in \bar{N} $,
    \[
        k \bar{n} a = k \kappa(\bar{n}) \underbrace{\big(a(\bar{n}) a\big)}_{ \in A_{\geq s''} } \big(a^{-1} n(\bar{n}) a\big)
    \]
    for all $ \bar{n} \in \bar{N} $ and $ a \in A_{\geq s''} $.
    Let $ \eps > 0 $ be sufficiently small so that
    \[
        \inf_{a \in \bar{A}_+ , \, y \in \bar{N} \, : \, |\log y| \leq \eps} a a(y^{a^{-\tfrac12}})^{-1} \geq \inf_{y \in \bar{N} \, : \, |\log y| \leq \eps} a(y)^{-1}
    \]
    is greater or equal than $ s'' $. Thus,
    \[
        \varphi\big(k y^{a^{\tfrac12}} a a(y^{a^{-\tfrac12}})^{-1}\big) = 0
    \]
    for all $ y \in \bar{N} $ such that $ |\log y| \leq \eps $ and $ a \in \bar{A}_+ $.
    Hence, there is a constant $ c > 0 $ such that
    \[
        |T \varphi\big(k y^{a^{\tfrac12}} a a(y^{a^{-\tfrac12}})^{-1}, \kappa(y^{a^{-\tfrac12}})^{-1}\big)| a^{2\rho} a(y^{a^{-\tfrac12}})^{-2\rho}
    \]
    is less or equal than $ c \, a(y)^{-2\rho} $ for all $ a \in \bar{A}_+ $, by Lemma \ref{Lem a a((bar n)^(a^(-1/2)))^(-1) converges}.
    
    By Lemma \ref{Lem for SU(1, 2)-reduction}, there is a representative of the Weyl element $ w \in K $ such that $ w^2 = e $.

    It follows from Lebesgue's theorem of dominated convergence that
    \[
        \lim_{t \to \infty} a_t^\rho (R^* (T \varphi))(ka_t)
    \]
    is equal to
    \[
        \int_N T \varphi\big(k a_B(y w)^{-1}, m_y^{-1} w\big) a_B(y w)^{-2\rho} \, dy \ .
    \]
    By Lemma \ref{Lem int_E f(n) dn = int_A int_(varphi(x) = 1) f(x^a) dx a^(2 xi_E rho) da} applied to $ E = N $, we obtain
    \[
        \int_A \int_{\{x \in N \, \mid \, a_B(x w) = 1 \}} T \varphi\big(k a_B(x^a w)^{-1}, m_y^{-1} w\big) a_B(x^a w)^{-2\rho} a^{2\rho} \, dx \, da \ .
    \]
    Since $ a_B(y^a w) = a^2 a_B(y w) = a^2 $, this is again equal to
    \[
        \int_A \int_{\{x \in N \, \mid \, a_B(x w) = 1 \}} T \varphi\big(k a^{-2}, m_x^{-1} w\big) a^{-2\rho} \, dx \, da \ .
    \]
    By the change of variables $ a \mapsto a^{-\frac12} $, this yields
    \[
        \frac12 \int_A \int_{\{x \in N \, \mid \, a_B(x w) = 1 \}} T \varphi\big(k a, m_x^{-1} w\big) a^{\rho} \, dx \, da \ .
    \]
    The first assertion follows now by Proposition \ref{Prop Other expression for Res_0 c_tau(lambda) o T}.

    Let $ V $ be a relatively compact subset of $ K(\Omega) $ and let $ 0 < p < q $. For $ y \in \bar{N} \smallsetminus \{e\} $, set $ \bar{m}_y = \lim_{a \to \infty} \tau\big(\kappa(y^{a^{-1}})\big) $.
    It remains to show that
    \begin{multline*}
        t^{\frac12} |\int_{ \{ y \in \bar{N} \, \mid \, |\log y| \in [p, q], \, k a_t \kappa(y^{a_t^{-\tfrac12}}) \in G(\Omega) \} } \tau\big(\kappa(y^{a_t^{-\tfrac12}})\big) T \varphi\big(k y^{a_t^{\tfrac12}} a_t a(y^{a_t^{-\tfrac12}})^{-1}, e \big) a_t^{2\rho} a(y^{a_t^{-\tfrac12}})^{-2\rho} \, dy \\
        - \int_{ \{ y \in \bar{N} \, \mid \, |\log y| \in [p, q] \} } \tau(\bar{m}_y) T \varphi\big(k a_B(w y)^{-1}, e\big) a_B(w y)^{-2\rho} \, dy|
    \end{multline*}
    is uniformly bounded when $ t \geq 1 $ and $ k \in V $. Indeed, by the above and as $ a_t a(y^{a_t^{-\tfrac12}})^{-1} \leq a_B(w y)^{-1} $ for all $ t > 0 $ and $ y \in \bar{N} \smallsetminus \{e\} $, we can choose $p$ and $q$ so that the integrand is zero on $ \{ y \in \bar{N} \mid |\log y| \not \in [p, q] \} $. Note that we could have used this argument already above in order to show that we can apply Lebesgue's theorem of dominated convergence.


    For $ Y \in \bar{n} $, $ H \in \a $, define
    \[
        F_k(Y, H) = T \varphi\big(k \exp(Y) \exp(H), e \big) \ .
    \]
    Then, $ \lim_{t \to \infty} F_k(\log(y^{a_t^{\tfrac12}}), \log(a_t a(y^{a_t^{-\tfrac12}})^{-1})) = F_k(0, \log(a_B(w y)^{-1})) $.

    Since a $ \C^1 $-function is locally Lipschitz, there is a constant $ L > 0 $ such that
    \begin{multline*}
        |F_k(\log(y^{a_t^{\tfrac12}}), \log(a_t a(y^{a_t^{-\tfrac12}})^{-1})) - F_k(0, \log(a_B(w y)^{-1}))| \\
        \leq L \big(|\log(y^{a_t^{\tfrac12}})| + |\log(a_t a(y^{a_t^{-\tfrac12}})^{-1}) - \log(a_B(w y)^{-1})|\big)
    \end{multline*}
    for all $ y \in \bar{N} $ with $ |\log y| \in [p, q] $ and all $ t $ sufficiently large (then $ y^{a_t^{\tfrac12}} a_t a(y^{a_t^{-\tfrac12}})^{-1} \in G(\Omega) $).

    For $ t \geq 1 $, $ |\log y^{a_t^{\tfrac12}}| 
    \leq \frac1{\sqrt{t}} |\log y| $.

    Using the Taylor expansion of $ \log(1 + x) $ at $ x = 0 $, one can show that
    \[
        |\log\bigl(a_t a(y^{a_t^{-\tfrac12}})^{-1}\bigr) - \log\bigl(a_B(w y)^{-1}\bigr)| \prec \frac1t \qquad (t \geq 1, y \in \bar{N} : |\log y| > \eps) \ .
    \]
    By Corollary \ref{Cor t |tau(kappa(n^(a_t) w)) - tau(m_n)| is uniformly bounded},
    \[
        \sqrt{t} |\tau\big(\kappa(y^{a_t^{-\tfrac12}})\big) - \tau(\bar{m}_y)| = \sqrt{t} |\tau\big(\kappa((w y w^{-1})^{a_t^{\tfrac12}})\big) - \tau(w \bar{m}_y)|
    \]
    is uniformly bounded when $ t \geq 1 $ and $ y \in \bar{N} $ such that $ |\log y| \in [p, q] $.

    The proposition follows now easily from these estimates and Lemma \ref{Lem a a((bar n)^(a^(-1/2)))^(-1) converges}.
\end{proof}

\begin{lem}\label{Lem R^* varphi is square-integrable}
    If $ \varphi \in W_\tau $, then $ R^* \varphi \in L^2(Y, V_Y(\tau, \phi)) $.
\end{lem}

\begin{rem}
    It follows that $ (R^* \varphi, f) $ is well-defined for all $ \varphi \in W_\tau $ and $ f \in L^2(Y, V_Y(\tau, \phi)) $.

    If, in addition, $ f \in \CS(Y, V_Y(\tau, \phi)) $, then $ (\varphi, f^\Omega) $ is also well-defined (recall that $ \varphi $ has compact support) and
    \[
        (R^* \varphi, f) = (\varphi, f^\Omega)
    \]
    by Lemma \ref{Lem (varphi, f^Omega) = (R^* varphi, f)}, as $ \Ccinf(Y, V_Y(\tau, \phi)) $ is dense in $ \CS(Y, V_Y(\tau, \phi)) $ and as
    \[
        f \in \CS(Y, V_Y(\tau, \phi)) \mapsto (\varphi, f^\Omega)
    \]
    is continuous. 
\end{rem}

\begin{proof}
    Let $ \varphi \in \Ccinf(\Hcal, V_{\Hcal}(\tau, \phi)) $. Since $ (R^* \varphi(g), R^* \varphi(g)) \geq 0 $ for all $ g \in G $,
    \[
        \int_{\Gamma \bs G} (R^* \varphi(g), R^* \varphi(g)) \, dg
    \]
    is equal to
    \begin{align*}
        & \int_G \chi(gK) (\int_{\{k \in K \, \mid \, gk \in G(\Omega)\}} \varphi(gk, k^{-1}) \, dk, \int_{\{k \in K \, \mid \, gk \in G(\Omega)\}} \varphi(gk, k^{-1}) \, dk) \, dg \\
        =& \int_K \int_{A_+} \chi(haK) (a^\rho \int_K \varphi(h a k, k^{-1}) \, dk, a^\rho \int_K \varphi(h a k, k^{-1}) \, dk) \underbrace{\gamma(a) a^{-2\rho}}_{\text{bounded}} \, da \, dh \ ,
    \end{align*}
    where $ \gamma(a) := \prod_{\mu \in \Phi^{+}}{ \sinh(\mu(H)) } $, by Fubini's theorem.
    If $ \varphi \in W_\tau $, then this is finite by Proposition \ref{Prop lim_(a to infty) a^rho (R^* varphi)(ka) converges} combined with \eqref{eq Decomposition of Ccinf(Gamma|G(Omega)/N x (M A) N|G/K, V_B'(tau, phi))}.
\end{proof}

\begin{lem}\label{Lem < f, R^* varphi > = 0 if f in clo(0CS)}
    Let $ f \in \overline{\c{\CS}(Y, V_Y(\tau, \phi))} $. Then,
    \[
        (f, R^* \varphi) = 0
    \]
    for all $ \varphi \in W_\tau $.
\end{lem}

\begin{rem}
    The pairing is well-defined as $ R^* \varphi $ is square-integrable by Lemma \ref{Lem R^* varphi is square-integrable} combined with \eqref{eq Decomposition of Ccinf(Gamma|G(Omega)/N x (M A) N|G/K, V_B'(tau, phi))}.
\end{rem}

\begin{proof}
    Let $ f \in \c{\CS}(Y, V_Y(\tau, \phi)) $ and let $ \varphi $ be as above. Then,
    \[
        (f, R^* \varphi) = (f^\Omega, \varphi) = 0 \ .
    \]
    Let now $ f \in \overline{\c{\CS}(Y, V_Y(\tau, \phi))} $. Then, there is a sequence $ (f_j)_j $ in $ \c{\CS}(Y, V_Y(\tau, \phi)) $ such that $ f_j $ converges to $f$ in $ L^2(Y, V_Y(\tau, \phi)) $. Since $ R^* \varphi $ is square-integrable,
    \[
        (f, R^* \varphi) = \lim_{j \to \infty} (f_j, R^* \varphi) = 0 \ .
    \]
\end{proof}

\begin{prop}\label{Prop < f, R^* varphi > = 0 => f^Omega = 0}
    Let $ f \in \CS(Y, V_Y(\tau, \phi)) $ be such that
    \[
        (f, R^* \varphi) = 0
    \]
    for all $ \varphi \in \tilde{W}_\tau $. Then, $ f^\Omega = 0 $.
\end{prop}

\begin{rem}
    Again the pairing is well-defined as $ R^* \varphi $ is square-integrable by Lemma \ref{Lem R^* varphi is square-integrable} combined with \eqref{eq Decomposition of Ccinf(Gamma|G(Omega)/N x (M A) N|G/K, V_B'(tau, phi))}.
\end{rem}

\begin{proof}
    Note that it suffices to prove the following: \\
    \textit{Let $ \sigma = \restricted{\tau}{M} $ (finite-dimensional representation of $M$). Let $ f \in \CS(Y, V_Y(\tau, \phi)) $ be such that
    \[
        \int_{\Gamma \bs G(\Omega)/MN} (f^\Omega(g, e), \varphi(g)) \, dg = 0
    \]
    for all $ \varphi \in \tilde{W}_\sigma $. Then, $ f^\Omega = 0 $.}

    Let $ \varphi \in \Ccinf(G(\Omega)/MN, V(\sigma_{MN}) \otimes V_\phi) $ (smooth section of $ V(\sigma_{MN}) \otimes V_\phi $ with compact support in $ G(\Omega)/MN $) be such that $ \int_A \varphi(k a) a^\rho \, da = 0 $ for all $ k \in K(\Omega) $. Then, $ \int_A \varphi(g a) a^\rho \, da $ vanishes also for all $ g \in G(\Omega) $.
    Let $ \pi^{\sigma_{MN}} $ denote the $G$-representation on the space of sections $ V(\sigma_{MN}) $. Let $ \pi(\gamma) $ be the induced action of $ \pi^{\sigma_{MN}}(\gamma) \otimes \phi(\gamma) $ ($ \gamma \in \Gamma $).
    Then,
    \[
        \Phi(g) := \sum_{\gamma \in \Gamma} (\pi(\gamma) \varphi)(gMN) = \sum_{\gamma \in \Gamma} \varphi(\gamma^{-1} g) \qquad (g \in G(\Omega))
    \]
    is well-defined and $ \Phi $ belongs to $ \tilde{W}_\sigma $. Indeed, since $ \Gamma $ acts properly discontinuously on $ \Omega $, it acts also properly continuously on $ G(\Omega)/MN $. So, the above sum is locally finite as $ \varphi $ has compact support in $ G(\Omega)/MN $.
    Thus,
    \begin{align}\label{eq int_(K(Omega)) int_A (f^Omega(k, a), varphi(k, a)) da dk = 0}
        0 &= \int_{\Gamma \bs G(\Omega)/MN} (f^\Omega(g, e), \Phi(g)) \, dg \\ \nonumber
        &= \int_{G(\Omega)/MN} (f^\Omega(g, e), \varphi(g)) \, dg \\ \nonumber
        &= \int_{K(\Omega)/M} \int_A (f^\Omega(k a, e), \varphi(k a)) a^{2\rho} \, da \, dk \\ \nonumber
        &= \int_{K(\Omega)} \int_A (f^\Omega(k, a), \varphi(k a)) \, da \, dk \ . \nonumber
    \end{align}
    Here, we view $ \varphi $ also as a function on $ K(\Omega)/M \times A $.

    Let $ V(\sigma) = K \times_M V_{\sigma} $.\index[n]{Vsigma@$ V(\sigma) $}
    Let us denote the space of conjugate-linear continuous functionals on $ \Cinf(K(\Omega)/M \times A, V(\sigma) \otimes V_\phi) $ with compact support by $ \C_c^{-\infty}(K(\Omega)/M \times A, V(\sigma) \otimes V_\phi) $.

    For $ h \in K(\Omega) $, $ a \in A $ and $ v \in V_\sigma \otimes V_\phi $, define
    \[
        \delta_{(h, a)} v \in \C_c^{-\infty}(K(\Omega)/M \times A, V(\sigma) \otimes V_\phi)
    \]
    via
    \[
        (\delta_{(h, a)} v, \psi) = (v, \psi(h, a))
    \]
    ($ \psi \in \Cinf(K(\Omega)/M \times A, V(\sigma) \otimes V_\phi) $).

    Define the partial Fourier transform $ \hat{T}(0) $ of $ T \in \C_c^{-\infty}(K(\Omega)/M \times A, V(\sigma) \otimes V_\phi) $ at $ \lambda = 0 $ by
    \[
        (\hat{T}(0), \psi) = (T, a^{\rho} \psi)
    \]
    ($ \psi \in \Cinf(K(\Omega)/M, V(\sigma) \otimes V_\phi) $). Let $ f \in \Ccinf(K(\Omega)/M \times A, V(\sigma) \otimes V_\phi) $. Then,
    \begin{multline*}
        (\hat{T}_f(0), \psi) = (T_f, a^{\rho} \psi) = \int_{K(\Omega)} \int_A (f(k, a), \psi(k)) a^{\rho} \, da \, dk \\
        = \int_{K(\Omega)} (\int_A f(k, a) a^{\rho} \, da, \psi(k)) \, dk =: (\hat{f}(0), \psi)
    \end{multline*}
    for all $ \psi \in \Cinf(K(\Omega)/M, V(\sigma) \otimes V_\phi) $.
    Let $ h \in K(\Omega) $, $ a_1, a_2 \in A $ and $ v \in V_\sigma \otimes V_\phi $. Let
    \[
        S = a_1^{-\rho} \delta_{(h, a_1)} v - a_2^{-\rho} \delta_{(h, a_2)} v \in \C_c^{-\infty}(K(\Omega)/M \times A, V(\sigma) \otimes V_\phi) \ .
    \]
    Then, $ \hat{S}(0) $ vanishes. Indeed, $ \hat{S}(0)(\psi) $ ($ \psi \in \Cinf(K(\Omega)/M, V(\sigma) \otimes V_\phi) $) is equal to
    \[
        a_1^{-\rho} (\delta_{(h, a_1)} v, a^{\rho} \psi)
        - a_2^{-\rho} (\delta_{(h, a_2)} v, a^{\rho} \psi)
        = a_1^{-\rho} (v, a_1^{\rho} \psi(h))
        - a_2^{-\rho} (v, a_2^{\rho} \psi(h))
        = 0 \ .
    \]
    \begin{clm}
        There is a sequence $ (f_j)_j $ in $ \Ccinf(K(\Omega)/M \times A, V(\sigma) \otimes V_\phi) $ such that
        \begin{enumerate}
        \item $ f_j $ converges to $ S = a_1^{-\rho} \delta_{(h, a_1)} v - a_2^{-\rho} \delta_{(h, a_2)} v $ in $ \C_c^{-\infty}(K(\Omega)/M \times A, V(\sigma) \otimes V_\phi) $, and
        \item $ \int_A f_j(k, a) \, a^{\rho} \, da = 0 $ for all $ k \in K(\Omega) $.
        \end{enumerate}
    \end{clm}

    \begin{proof}
        Choose $ (\chi_n)_{n \in \NN} $ in $ \Ccinf(K(\Omega)/M, [0, 1]) $ such that
        \[
            \supp(\chi_n) \underset{n \to \infty}{\longrightarrow} \{ hM \} \qquad \text{and} \qquad \int_K{ \chi_n(kM) \, dk } = 1  \quad \text{for all } n \in \NN \ .
        \]
        Let $ \{v_i\} $ be an orthonormal basis of $ V_\sigma \otimes V_\phi $.
        Let $ g^{(i)}_{0, 1}(k) = \tau(k^{-1}h) v_i $. Then, $ g^{(i)}_{0, 1} \in \Cinf(K(\Omega)/M, V(\sigma) \otimes V_\phi) $.
        Choose $ g_{0, 2} \in \Ccinf(A, \RR) $ such that $ \int_A g_{0, 2}(a) a^{-\rho} \, da = 1 $.

        For $ k \in K(\Omega) $ and $ a \in A $, set $ g^{(i)}_0(k, a) = g^{(i)}_{0, 1}(k) g_{0, 2}(a) $. Then,
        \[
            g^{(i)}_0 \in \Cinf(K(\Omega)/M \times A, V(\sigma) \otimes V_\phi)
        \]
        is such that
        \[
            \int_A (g^{(i)}_0(k, a), \tau(k^{-1}h) v_i) \, a^{\rho} \, da = 1
        \]
        for all $ k \in K(\Omega) $. Let $ (g_{j, 1})_j $ (resp. $ (g_{j, 2})_j $) be a Dirac sequence in $ \Ccinf(A, \RR) $ such that $ \supp(g_{j, 1}) $ (resp. $ \supp(g_{j, 2}) $) converges to $ \{a_1\} $ (resp. $ \{a_2\} $).
        Set
        \[
            g_j(k, a) = \chi_j(kM) \big(a_1^{-\rho} g_{j, 1}(a) \tau(k^{-1}h) v - a_2^{-\rho} g_{j, 2}(a) \tau(k^{-1}h) v\big) \ .
        \]
        Then, $ g_j \in \Ccinf(K(\Omega)/M \times A, V(\sigma) \otimes V_\phi) $ converges to $ S $ in $ \C_c^{-\infty}(K(\Omega)/M \times A, V(\sigma) \otimes V_\phi) $.
        For $ v' \in V_\tau \otimes V_\phi $, set
        \[
            \lambda_{i, j}(k) = \int_A (g_j(k, a), \tau(k^{-1}h) v_i) a^{\rho} \, da \ .
        \]
        Then, $ k \in K(\Omega) \mapsto \lambda_{i, j}(k) $ belongs to $ \Ccinf(K(\Omega)/M, \CC) $.
        For $ k \in K(\Omega) $ and $ a \in A $, set
        \[
            f_j(k, a) = g_j(k, a) - \sum_l \lambda_{l, j}(k) g^{(l)}_0(k, a) \ .
        \]
        Then,
        \[
            f_j \in \Ccinf(K(\Omega)/M \times A, V(\sigma) \otimes V_\phi) \qquad (\supp f_j \subset \supp \chi_j) \ .
        \]
        For all $ k \in K(\Omega) $ and all $i$, we have
        \[
            \int_A (f_j(k, a), \tau(k^{-1}h) v_i) a^{\rho} \, da = \lambda_{i, j}(k) - \sum_l \lambda_{l, j}(k) (v_l, v_i) = 0 \ .
        \]
        Thus, $ \tau(h^{-1}k) \int_A f_j(k, a) a^{\rho} \, da = 0 $ for all $ k \in K(\Omega) $. Hence, $ \int_A f_j(k, a) a^{\rho} \, da $ vanishes for all $ k \in K(\Omega) $, too.
        Let $ \psi \in \Cinf(K(\Omega)/M \times A, V(\sigma) \otimes V_\phi) $. Then,
        \begin{align*}
            & \int_{K(\Omega)} \int_A (f_j(k, a), \psi(k, a)) \, da \, dk \\
            =& \underbrace{ \int_{K(\Omega)} \int_A (g_j(k, a), \psi(k, a)) \, da \, dk }_{ \to S(\psi) }
            - \sum_l \underbrace{\int_{K(\Omega)} \lambda_j(k, v_l) \int_A (g^{(l)}_0(k, a), \psi(k, a)) \, da \, dk}_{ =: I_{j, l}(\psi) }
        \end{align*}
        and
        \begin{align*}
            I_{j, l}(\psi) &= \int_{K(\Omega)} \int_A (g_j(k, a), a^{\rho} \int_A (\psi(k, a'), g^{(l)}_0(k, a')) \, da' \tau(k^{-1}h) v_l) \, da \, dk \\
            &= (\hat{g_j}(0), \int_A (\psi(k, a'), g^{(l)}_0(k, a')) \, da' \tau(k^{-1}h) v_l)
        \end{align*}
        converges to 0 as $ \hat{g_j}(0) $ converges to $ \hat{S}(0) $ in $ \C_c^{-\infty}(K(\Omega)/M \times A, V(\sigma) \otimes V_\phi) $. The claim follows.
    \end{proof}

    It follows from the claim, applied to $ \{v_i\} $, and \eqref{eq int_(K(Omega)) int_A (f^Omega(k, a), varphi(k, a)) da dk = 0} that
    \begin{equation} \label{eq a_1^rho f^Omega(k, a_1) = a_2^rho f^Omega(k, a_2)}
        a_1^{-\rho} f^\Omega(k, a_1) = a_2^{-\rho} f^\Omega(k, a_2) \ .
    \end{equation}
    Consequently, the right-hand side of \eqref{eq a_1^rho f^Omega(k, a_1) = a_2^rho f^Omega(k, a_2)} converges to zero when $ a_2 $ tends to $ \infty $. So, $ f^\Omega(k, a_1) = 0 $. The proposition follows.
\end{proof}

\begin{thm}\label{Thm a Schwartz fct in the closure of the space of cusp forms is a cusp form}
    Let $ \tau \in \hat{K} $. Then,
    \[
        \overline{\c{\CS}(Y, V_Y(\tau, \phi))} \cap \CS(Y, V_Y(\tau, \phi)) = \c{\CS}(Y, V_Y(\tau, \phi)) \ .
    \]
    Thus, we also have
    \[
        \overline{\c{\CS}(\Gamma \bs G, \phi)}(\tau) \cap \CS(\Gamma \bs G, \phi) = \c{\CS}(\Gamma \bs G, \phi)(\tau) \ .
    \]
\end{thm}

\begin{proof}
    As $ \c{\CS}(Y, V_Y(\tau, \phi)) $ is clearly contained in the left-hand side, it suffices to prove the other inclusion.

    Let $ f \in \overline{\c{\CS}(Y, V_Y(\tau, \phi))} \cap \CS(Y, V_Y(\tau, \phi)) $.
    Then, $ (f, R^* \varphi) = 0 $ for all $ \varphi \in W_\tau $, by Lemma \ref{Lem < f, R^* varphi > = 0 if f in clo(0CS)}.
    It follows by Proposition \ref{Prop < f, R^* varphi > = 0 => f^Omega = 0} that $ f^\Omega = 0 $.
\end{proof}

\newpage

\subsubsection{About the orthogonality of Eisenstein series to the space of cusp forms}

We show in this section that the Eisenstein series are orthogonal to the space of cusp forms when the critical exponent $ \delta_\Gamma $ is negative. The proof without making this assumption is much more complicated (see Theorem \ref{Thm L^2(Gamma|G, phi)_ac is orthogonal to the space of cusp forms} in Section \ref{ssec:The contrib of the WPs to the closure of the space of cusp forms}).

\begin{dfn}[{\cite[p.143]{BO00}}] \nl
    For $ \varphi \in \C^{-\infty}(B, V_B(\sigma_\lambda, \phi)) $ and $ v \in \C^K(\dX, V(\tilde{\sigma}_{-\lambda})) $, we define the \textit{Eisenstein series} $ E(\lambda, \varphi, v) \in \Cinf(\Gamma \bs G, \phi) $ by\index{Eisenstein series}\index[n]{ElambdaphivQ@$ E(\lambda, \varphi, v) $}
    \[
        E(\lambda, \varphi, v) = c_{\ext(\varphi), v}
    \]
    (well-defined as a meromorphic function).
\end{dfn}

\begin{lem}\label{Lem sum_(gamma in Gamma) a(gamma x)^{-Re(lambda) -rho} < oo iff Re(lambda) > delta_Gamma}
    Let $ \lambda \in \a^* $ and $ x \in G(\Omega) $. Then,
    \[
        \sum_{\gamma \in \Gamma} a(\gamma^{-1} x)^{-\lambda -\rho} < \infty \iff \lambda > \delta_\Gamma \ .
    \]
\end{lem}

\begin{proof}
    Let us assume that $ \lambda \geq -\rho $. The proof for $ \lambda < -\rho $ works completely similarly and is left to the reader.
    Let $ U \subset K $ be an open neighbourhood of $ \kappa(x) $ such that $ \bar{U} $ is a compact subset of $ K(\Omega) $.
    Let $W$ be a compact subset of $ (\dX \smallsetminus \bar{U}M)M $ containing $ \Lambda M $ in its interior. Then, by Lemma 2.3 of \cite[p.84]{BO00}, there is a constant $ C > 0 $ such that
    $
        C a_\gamma \leq a(\gamma^{-1} \kappa(x))
    $
    for all $ \gamma \in \Gamma $ with $ k_\gamma \in W $.
    Since the closure of $ \dX \smallsetminus WM $ is a compact subset of $ \Omega $ and since $ \Gamma $ acts properly discontinuously on $ X \cup \Omega $, $ S := \{ \gamma K \mid \gamma \in \Gamma : k_\gamma M \in \dX \smallsetminus WM \} $ is finite.
    Since moreover $ \sum_{\gamma \in \Gamma} a(\gamma^{-1} x)^{-\lambda -\rho}  = \sum_{gK \in \{ \gamma K \mid \gamma \in \Gamma \}} a(g^{-1} x)^{-\lambda -\rho} $, we may assume without loss of generality that $ k_\gamma \in W $. Since on one hand
    \[
        \sum_{\gamma \in \Gamma} a(\gamma^{-1} x)^{-\lambda -\rho} = a(x)^{-\lambda -\rho} \sum_{\gamma \in \Gamma} a(\gamma^{-1} \kappa(x))^{-\lambda -\rho}
        \leq C a(x)^{-\lambda -\rho} \sum_{\gamma \in \Gamma} a_\gamma^{-\lambda -\rho}
    \]
    and on the other hand
    \[
        \sum_{\gamma \in \Gamma} a(\gamma^{-1} x)^{-\lambda -\rho}
         = a(x)^{-\lambda -\rho} \sum_{\gamma \in \Gamma} a(\gamma^{-1} \kappa(x))^{-\lambda -\rho}
         \geq a(x)^{-\lambda -\rho} \sum_{\gamma \in \Gamma} a_\gamma^{-\lambda -\rho} \ ,
    \]
    the lemma follows by definition of $ \delta_\Gamma $ and by the fact that $ \sum_{\gamma \in \Gamma} a_\gamma^{-\lambda -\rho} $ diverges for $ \lambda = \delta_\Gamma $.
\end{proof}

\begin{dfn}
    Define the Schwartz space on $ G $ by
    \[
        \CS(G, V_\phi) = \{ f \in \Cinf(G, V_\phi) \mid p_{r, X, Y}(f) < \infty \quad \forall r \geq 0, X, Y \in \U(\g) \} \ ,
    \]
    where $ p_{r, X, Y}(f) = \LOI{p}{X}_{r, X, Y}(f) $.
    We equip $ \CS(G, V_\phi) $ with the topology induced by the seminorms.
\end{dfn}

\begin{rem}
    This Schwartz space is just the vector-valued version of $ \CS(G, 1) $. So, instead of the absolute value, we take here the norm on $ V_\phi $.
\end{rem}

\needspace{3\baselineskip}
The following results are of independent interest.

\begin{lem}\label{Lem chi f in CS(G, V_phi)}
    Let $ \chi \in \Ccinf(X \cup \Omega) $ be as in Lemma \ref{Lem Cut-off function in the convex-cocompact case}.
    \\ Let $ f \in \Cinf(\Gamma \bs G, \phi) $. Then,
    \[
        f \in \CS(\Gamma \bs G, \phi) \iff \chi f \in \CS(G, V_\phi) \ .
    \]
\end{lem}

\begin{rem}
    Here $ \chi f(g) $ ($g \in G$) is by definition equal to $ \chi(gK)f(g) $.
\end{rem}

\begin{proof}
    Let $ X, Y \in \g, r \geq 0 $ and let $ f \in \Cinf(\Gamma \bs G, \phi) $.

    $\underline{\Rightarrow:}$
    Assume that $ f \in \CS(\Gamma \bs G, \phi) $. Since $ \Lambda_{\{e\}} = \emptyset $, $ \Omega_{\{e\}} = \dX $. So, $ U_{\F_{\{e\}}} $ is equal to $ \bar{X} $.
    Let $ U $ be an open subset of $X$, containing $ \supp(\chi) \cap X $, that is relatively compact in $ X \cup \Omega $. Then, $ U \in \U_\Gamma $. By the Leibniz rule,
    \[
        \sup_{gK \in X} (1 + \log a_g)^r a_g^\rho |L_X \chi f(g)|
        = \sup_{gK \in U} (1 + \log a_g)^r a_g^\rho |L_X \chi f(g)|
    \]
    is less or equal than
    \begin{multline*}
        \sup_{gK \in X} |L_X \chi(gK)| \sup_{gK \in U} (1 + \log a_g)^r a_g^\rho |f(g)| \\
        \qquad + \sup_{gK \in X} |\chi(gK)| \sup_{gK \in U} (1 + \log a_g)^r a_g^\rho |L_X f(g)| \ .
    \end{multline*}
    This is finite by Lemma \ref{Lem Cut-off function in the convex-cocompact case} and as $ f \in \CS(\Gamma \bs G, \phi) $. The assertion for $ Y \in \g $, respectively $ X, Y \in \g $, can be shown similarly.

    $\underline{\Leftarrow:}$ Assume that $ \chi f \in \CS(G, V_\phi) $. As $ \clo(U) $ is compact in $ X \cup \Omega $ and since $ \supp \chi $ contains by construction an open set $ V \subset X \cup \Omega $ such that $ \bigcup_{\gamma \in \Gamma} \gamma V = X \cup \Omega $, there exist $ \gamma_1 := e, \dotsc, \gamma_n \in \Gamma $ such that $ U \subset \bigcup_{i=1}^n \gamma_i \supp \chi =: C $. Since $ \Gamma $ acts properly discontinuously on $ X \cup \Omega $,
    \[
        \{\gamma \in \Gamma \mid \gamma U \cap \supp \chi \neq \emptyset \} \subset \{\gamma \in \Gamma \mid \gamma C \cap C \neq \emptyset \}
    \]
    is finite. Thus, we have
    \[
        |f(g)| = \sum_{\gamma \in \Gamma} |\chi f(\gamma g)| \leq p_{r, 1, 1}(\chi f) \sum_{i=1}^n a_{\gamma_i g}^{-\rho} (1 + \log a_{\gamma_i g})^{-r}
    \]
    for all $ gK \in U $. By Lemma \ref{Lem 1 + log a_h <= (1 + log a_gh)(1 + log a_g)}, this is less or equal than
    \[
        p_{r, 1, 1}(\chi f) \sum_{i=1}^n a_{\gamma_i}^{\rho} (1 + \log a_{\gamma_i})^r a_g^{-\rho} (1 + \log a_g)^{-r} \ .
    \]
    Hence, $ \LOI{p}{U}_{r, 1, 1}(f) $ is finite.
    Since $ \{\gamma \in \Gamma \mid \gamma U \cap \supp \chi \neq \emptyset \} $ is finite, we have
    \begin{multline*}
        |L_X f(g)| = \sum_{i=1}^n |\chi L_X f(\gamma_i g)| = \sum_{i=1}^n |L_X (\chi f)(\gamma_i g) - (L_X \chi(\gamma_i gK)) f(\gamma_i g)| \\
        \leq \sum_{i=1}^n |L_X (\chi f)(\gamma_i g)| + \big(\sum_{i=1}^n |L_X \chi(\gamma_i gK)|\big) |f(g)|
    \end{multline*}
    for all $ gK \in U $. As $ \chi f \in \CS(G, V_\phi) $,
    \[
        \sum_{i=1}^n |L_X (\chi f)(\gamma_i g)| \leq p_{r, X, 1}(\chi f) \sum_{i=1}^n a_{\gamma_i}^{\rho} (1 + \log a_{\gamma_i})^r a_g^{-\rho} (1 + \log a_g)^{-r}
    \]
    with $ p_{r, X, 1}(\chi f) < \infty $.
    As in addition $ \sup_{xK \in X} |L_X \chi(xK)| < \infty $ by \eqref{Lem Cut-off function in the convex-cocompact case eq4} of Lemma \ref{Lem Cut-off function in the convex-cocompact case},
    \[
        \LOI{p}{U}_{r, X, 1}(f) \leq p_{r, X, 1}(\chi f) \sum_{i=1}^n a_{\gamma_i}^{\rho} (1 + \log a_{\gamma_i})^r + \big(\sum_{i=1}^n \sup_{gK \in X} |L_X \chi(\gamma_i gK)|\big) \LOI{p}{U}_{r, 1, 1}(f)
    \]
    is finite. The assertion for $ Y \in \g $, respectively $ X, Y \in \g $, is proven similarly. In that case, one must use Remark \ref{Rem of Lem Cut-off function in the convex-cocompact case} of Lemma \ref{Lem Cut-off function in the convex-cocompact case} instead of \eqref{Lem Cut-off function in the convex-cocompact case eq4}. The assertion for $ X, Y \in \U(\g) $ follows by induction.
    This completes the proof of the lemma.
\end{proof}

\begin{lem}\label{Lem chi f in Ccinf(G, V_phi)}
    Let $ \chi \in \Ccinf(X \cup \Omega) $ be as in Lemma \ref{Lem Cut-off function in the convex-cocompact case}.
    \\ Let $ f \in \Cinf(\Gamma \bs G, \phi) $. Then,
    \[
        f \in \Ccinf(\Gamma \bs G, \phi) \iff \chi f \in \Ccinf(G, V_\phi) \ .
    \]
\end{lem}

\begin{proof}
    Let $ f \in \Cinf(\Gamma \bs G, \phi) $.

    $\underline{\Rightarrow:}$
    Assume that $ f \in \Ccinf(\Gamma \bs G, \phi) $. Then, there is a compact subset $ C $ of $X$ such that $ \supp(f)K \subset \Gamma C $.
    As in the proof of the previous lemma, we can show that there exist $ \gamma_1 := e, \dotsc, \gamma_n \in \Gamma $ such that $ C \subset \bigcup_{i=1}^n \gamma_i \supp \chi =: C' $. Since $ \Gamma $ acts properly discontinuously on $ X \cup \Omega $,
    \[
        \Gamma_{C'} := \{\gamma \in \Gamma \mid \gamma \supp(\chi) \cap \gamma C \neq \emptyset \} \subset \{\gamma \in \Gamma \mid \gamma C' \cap C' \neq \emptyset \}
    \]
    is finite. So, $ \supp(\chi f)K \subset \supp(\chi) \cap \supp(f)K \subset \supp(\chi) \cap \Gamma C \subset \bigcup_{\gamma \in \Gamma_{C'}} \gamma C $ and this set is compact.

    $\underline{\Leftarrow:}$ Assume that $ \chi f \in \Ccinf(G, V_\phi) $.
    For $ g \in G $, we have
    \[
        f(g) = \sum_{\gamma \in \Gamma} \chi f(\gamma g)
    \]
    Thus, $ \supp f \subset \Gamma \supp(\chi f) $.
    This completes the proof of the lemma.
\end{proof}

\begin{lem}\label{Lem about converges in Schwartz spaces}
    Let $ \chi \in \Ccinf(X \cup \Omega) $ be as in Lemma \ref{Lem Cut-off function in the convex-cocompact case}, let $ f \in \CS(\Gamma \bs G, \phi) $ and let $ (f_j)_j $ be a sequence in $ \CS(\Gamma \bs G, \phi) $. Then,
    \[
        f_j \to f \text{ in } \CS(\Gamma \bs G, \phi) \quad \iff \quad \chi f_j \to \chi f \text{ in } \CS(G, V_\phi) \ .
    \]
\end{lem}

\begin{proof}
    The lemma follows from the proof of Lemma \ref{Lem chi f in CS(G, V_phi)}.
\end{proof}

\begin{dfn}
    Let $ f \in \CS(\Gamma \bs G, \phi) $ and $ \lambda \in \a_\CC^* $. Define
    \[
        \widehat{f^\Omega}(\lambda, x, k) := \int_A a^{\lambda - \rho} f^\Omega(x, a k) \, da \qquad (x \in G(\Omega), k \in K)
    \]
    whenever the integral converges.
\end{dfn}

Let us show now that $ \widehat{f^\Omega}(\lambda, x, k) $ ($ x \in G(\Omega), k \in K $) is well-defined for $ \lambda \in i\a^* $ and $ \delta_\Gamma < 0 $.

\begin{lem}\label{Lem widehat{f^Omega}(mu, . , k) well-defined}
    Let $ f \in \CS(\Gamma \bs G, \phi) $. Then,
    \[
        \int_A \int_N |a^{\lambda - \rho} f(x n a k)| \, dn \, da < \infty \qquad (x \in G(\Omega), k \in K)
    \]
    for $ \lambda \in i\a^* $ and $ \delta_\Gamma < 0 $.
\end{lem}

\begin{proof}
    Let $ \chi \in \Ccinf(X \cup \Omega) $ be as in Lemma \ref{Lem Cut-off function in the convex-cocompact case} and let $ \lambda \in \a_\CC^* $. Since
    \[
        |f(g)| = \sum_{\gamma \in \Gamma} \chi(\gamma g K) |f(g)| = \sum_{\gamma \in \Gamma} |\chi f(\gamma g)|
    \]
    and since $ \chi \geq 0 $, $ \int_N |a^{\lambda - \rho} f(xnak)| \, dn $ ($ x \in G(\Omega) $) is equal to
    \begin{equation}\label{Proof Lem widehat{f^Omega}(mu, . , k) well-defined eq1}
        \sum_{\gamma \in \Gamma} a^{\R(\lambda) - \rho} \int_N |\chi f(\gamma x n a k)| \, dn
    \end{equation}
    by Fubini's theorem. 
    By Lemma \ref{Lem chi f in CS(G, V_phi)}, $ \chi f \in \CS(G, V_\phi) $.
    Fix $ d > 1 $.
    By Proposition \ref{Prop f^Omega well-defined} with $ \Gamma = \{e\} $, \eqref{Proof Lem widehat{f^Omega}(mu, . , k) well-defined eq1} is less or equal than
    \[
        c \, p_{2d, 1, 1}(f) \sum_{\gamma \in \Gamma} a^{\R(\lambda) - \rho} a(\gamma x)^{-\rho}
            a^\rho (1 + |\log a(\gamma x) a|)^{-d}
    \]
    for some constant $ c > 0 $. Thus,
    \begin{align*}
        & \sum_{\gamma \in \Gamma} \int_A a^{\R(\lambda)} a(\gamma x)^{-\rho} (1 + |\log a(\gamma x) a|)^{-d} \, da \\
        = & \sum_{\gamma \in \Gamma} a(\gamma x)^{-\rho} a(\gamma x)^{-\R(\lambda)} \int_A a^{\R(\lambda)} (1 + |\log a|)^{-d} \, da \\
        = & \sum_{\gamma \in \Gamma} a(\gamma x)^{-\R(\lambda) -\rho} \Big(\int_{-\infty}^0 e^{\R(\lambda)t} (1 - t)^{-d} \, dt
            + \int_0^\infty e^{\R(\lambda)t} (1 + t)^{-d} \, dt \Big) \ .
    \end{align*}
    The integrals converge if and only if $ \R(\lambda) = 0 $ and, 
    by Lemma \ref{Lem sum_(gamma in Gamma) a(gamma x)^{-Re(lambda) -rho} < oo iff Re(lambda) > delta_Gamma}, the series converges if and only if $ \R(\lambda) > \delta_\Gamma $.
    Thus, $ \int_A \int_N |a^{\lambda - \rho} f(xnak)| \, dn \, da $ is finite by Fubini's theorem if $ \R(\lambda) = 0 $ and $ \delta_\Gamma < 0 $.
\end{proof}

Recall that the pairing between $ \Cinf(B, V_B(\sigma_\mu, \phi)) $ ($ \mu \in \a^*_\CC $) and $ \Cinf(B, V_B(\sigma_{-\bar{\mu}}, \phi)) $ is defined by
\[
    \int_B (f_1(x), f_2(x)) \ ,
\]
where $ (\cdot, \cdot) $ denotes the scalar product on $ V_\sigma \otimes V_\phi $ (see \eqref{eq Def of (.,.)_B}). This pairing induces in a natural way a pairing between $ \C^{-\infty}(B, V_B(\sigma_\mu, \phi)) $ and $ \Cinf(B, V_B(\sigma_{-\bar{\mu}}, \phi)) $.

\begin{prop}\label{Prop with widehat f^Omega}
    Let $ \lambda \in i \a^* $. Let $ f \in \CS(\Gamma \bs G, \phi) $, $ \varphi \in \C^{-\infty}(B, V_B(\sigma_\lambda, \phi)) $ and $ v \in \C^K(\dX, V(\tilde{\sigma}_{-\lambda})) $. Assume $ \delta_\Gamma < 0 $. Then, $ \int_K{ \overline{v(k^{-1})} \otimes \widehat{f^\Omega}(-\lambda, \cdot, k) \, dk } \in \Cinf(B, V_B(\sigma_\lambda, \phi)) $ and
    \begin{equation*}
        (E(\lambda, \varphi, v), f)_{L^2(\Gamma \bs G, \phi)} := \int_{\Gamma \bs G}{ (E(\lambda, \varphi, v)(g), f(g)) \, dg }
    \end{equation*}
    is equal to
    \begin{equation*}
        (\varphi, \int_K{ \overline{v(k^{-1})} \otimes \widehat{f^\Omega}(-\lambda, \cdot, k) \, dk })_B \ .
    \end{equation*}
    In particular,
    $
        \int_{\Gamma \bs G}{ (E(\lambda, \varphi, v)(g), f(g)) \, dg } = 0
    $
    if $ f \in \c{\CS}(\Gamma \bs G, \phi) $. In other words, the Eisenstein series are orthogonal to the space of cusp forms.
\end{prop}

\begin{rem}
    If $ \varphi \not \in \Cinf(B, V_B(\sigma_\lambda, \phi)) $, then the integral
    \[
        \int_{\Gamma \bs G}{ (E(\lambda, \varphi, v)(g), f(g)) \, dg }
    \]
    is in general only a formal notation meaning that the tempered distribution $ E(\lambda, \varphi, v) \in \CS'(\Gamma \bs G, \phi) $ (see Remark \ref{Rem discrete series rep and tempered rep} \eqref{Rem discrete series rep and tempered rep, Lemma 9.4}) has to be applied to the Schwartz function $f$.
\end{rem}

\begin{proof}
    Let $ \lambda \in i \a^* $ and $\varphi \in \Cinf(B, V_B(\sigma_\lambda, \phi)) $. Let $f$ and $v$ be as above. As $ \Cinf(B, V_B(\sigma_\lambda, \phi)) $ is dense in $ \C^{-\infty}(B, V_B(\sigma_\lambda, \phi)) $ and as the pairings depend continuously on $\varphi$, we may assume without loss of generality that $ \varphi \in \Cinf(B, V_B(\sigma_\lambda, \phi)) $.

    Let $ \{v_i\} $ be an orthonormal basis of $ V_\phi $ and let $ \tilde{v}_j \in V_{\tilde{\phi}} $ be such that $ \langle v_i, \tilde{v}_j \rangle = \delta_{i, j} $.
    We identify $ V_\phi $ with $ \CC^{\dim(V_\phi)} $ via the basis $ \{v_i\} $. Then,
    \[
        (v, v') = \sum_j \overline{\langle v, \tilde{v}_j \rangle} v'_j
    \]
    for all $ v, v' \in V_\phi $.
    Note that $ \int_K{ \overline{v(k^{-1})} \otimes \widehat{f^\Omega}(-\lambda, \cdot, k) \, dk } \in \Cinf(B, V_B(\sigma_\lambda, \phi)) $.

    Let $ \chi_\infty \in \Ccinf(\Omega) $ be as in Lemma \ref{Lem Cut-off function in the convex-cocompact case}. Assume $ \delta_\Gamma < 0 $.
    By Lemma \ref{Lem widehat{f^Omega}(mu, . , k) well-defined},
    \begin{multline*}
        \int_{\Gamma \bs \Omega}{ (\varphi(x), \int_K{ \overline{v(k^{-1})} \otimes \widehat{f^\Omega}(-\lambda, x, k) \, dk }) } \\
        = \int_{\Gamma \bs \Omega}{ \int_K \int_A \int_N (\varphi(x), \overline{v(k^{-1})} \otimes f(x n a k) a^{-2\rho} a^{-\lambda+\rho}) \, dn \, da \, dk }
    \end{multline*}
    is well-defined and it is equal to
    \begin{equation}\label{Proof Prop with widehat f^Omega eq1}
        \int_{\Gamma \bs \Omega} \int_G (\varphi(x), \overline{v(\kappa(g^{-1}))} \otimes f(x g) a(g^{-1})^{\lambda-\rho}) \, dg \ .
    \end{equation}
    Moreover, we can apply Fubini's theorem. Thus, \eqref{Proof Prop with widehat f^Omega eq1} is equal to
    \begin{align*}
        & \int_G{ \int_K{ \chi_\infty(x M) (a(g^{-1}x)^{-\lambda-\rho} \varphi(x), \overline{v(\kappa(g^{-1}x))} \otimes f(g)) \, dx } \, dg } \\
        = & \int_G{ \int_K{ \chi_\infty(x M) (\langle \varphi(x), \pi^{\tilde{\sigma}, -\lambda}(g)v(x) \rangle, f(g)) \, dx } \, dg } \\
        = & \sum_j \int_{\Gamma \bs G}{ \sum_{\gamma \in \Gamma} \int_K \chi_\infty(x M) \overline{\langle \varphi(x), \pi^{\tilde{\sigma}, -\lambda}(\gamma g)v(x) \otimes \tilde{v}_j \rangle} f_j(\gamma g) \, dx \, dg } \ .
    \end{align*}
    Since $ (\cdot, \cdot) $ is unitary, this is equal to
    \[
        \sum_j \int_{\Gamma \bs G}{ \sum_{\gamma \in \Gamma} \int_K \chi_\infty(x M) \overline{\langle \varphi(x), \pi^{\tilde{\sigma}, -\lambda}(\gamma g)v(x) \otimes \tilde{\phi}(\gamma) \tilde{v}_j \rangle} f_j(g) \, dx \, dg } \ .
    \]
    By definition of $ \pi(\gamma) $, this yields
    \[
        \sum_j \int_{\Gamma \bs G}{ \sum_{\gamma \in \Gamma} \int_K \chi_\infty(x M) \overline{\langle \varphi(x), \pi(\gamma)(\pi^{\tilde{\sigma}, -\lambda}(g)v(x) \otimes \tilde{v}_j) \rangle} f_j(g) \, dx \, dg } \ .
    \]
    By Lemma \ref{Lem int_B (f_1(x), f_2(x)) = int_K chi_infty(kM) (f_1(k), f_2(k)) dk}, this is again equal to
    \begin{align*}
        & \sum_j \int_{\Gamma \bs G}{ \overline{\langle \varphi, \pi_*(\pi^{\tilde{\sigma}, -\lambda}(g)v \otimes \tilde{v}_j) \rangle_B} f_j(g) \, dg } \\
        = & \sum_j \int_{\Gamma \bs G}{ \overline{\langle \ext(\varphi), \pi^{\tilde{\sigma}, -\lambda}(g)v \otimes \tilde{v}_j \rangle_B} f_j(g) \, dg } \\
        = & \int_{\Gamma \bs G}{ (\langle \ext(\varphi), \pi^{\tilde{\sigma}, -\lambda}(g)v \rangle_{\dX}, f(g)) \, dg } \ .
    \end{align*}
    The proposition follows.
\end{proof}

\newpage

\subsubsection{Wave packets of Eisenstein series}

We show that a compactly supported wave packet belongs to the Schwartz space.

\medskip

Until the end of this section, we assume again that $ X \neq \OO H^2 $.

By \cite[p.86]{BO00}, we can consider holomorphic, smooth or continuous families of sections $ \a^*_\CC \ni \mu \mapsto f_\mu \in \C^{\pm \infty}(B, V_B(\sigma_\mu, \phi)) $.
We denote the linear space of smooth families $ \a^*_+ \ni \mu \mapsto \varphi_{i\mu} \in \Cinf(B, V_B(\sigma_{i\mu}, \phi)) $ with compact support in $ \a^*_+ $ by $ \Hcal^\sigma_0(\phi) $\index[n]{Hsigma@$ \Hcal^\sigma_0(\phi) $, $ \Hcal^{\tilde{\sigma}, K} $} and the linear space of smooth families $ \a^*_+ \ni \mu \mapsto v_{-i\mu} \in \C^K(\dX, V(\tilde{\sigma}_{-i\mu})) $ by $ \Hcal^{\tilde{\sigma}, K} $.

\begin{dfn}\label{Dfn wave packet}
    Let $ \varphi \in \Hcal^\sigma_0(\phi) $ and $ v \in \Hcal^{\tilde{\sigma}, K} $. Then, the \textit{wave packet transform}
    \[
        E(\varphi, v) \colon \Hcal^\sigma_0(\phi) \times \Hcal^{\tilde{\sigma}, K} \to \Cinf(\Gamma \bs G, \phi)
    \]
    is defined by
    \[
        E(\varphi, v) := \int_{\a^*_+} c_{\ext(\varphi_{i\mu}), v_{-i\mu}} \ p_\sigma(i\mu) \, d\mu \ .
    \]
    The section $ E(\varphi, v) $ is called a \textit{compactly supported wave packet (of Eisenstein series)}.\index{compactly supported wave packet}\index[n]{Ephiv@$ E(\varphi, v) $}
\end{dfn}

\begin{prop}\label{Prop Image of wave packet transform in Schwartz space}
    If $ \varphi \in \Hcal^\sigma_0(\phi) $ and $ v \in \Hcal^{\tilde{\sigma}, K} $, then $ E(\varphi, v) $ belongs to $ \CS(\Gamma \bs G, \phi) $.
\end{prop}

\begin{proof}
    Compare with Lemma 10.7 of \cite[p.147]{BO00}. Let $ U \in \U_\Gamma $ and let $ gK \in U $. Fix $ X, Y \in \U(\g) $. Then, $ L_X R_Y E(\varphi, v)(g) $ is equal to
    \[
         \int_{\a^*_+} c_{\pi^{\sigma, i\mu}(X) \ext(\varphi_{i\mu}), \pi^{\tilde{\sigma}, -i\mu}(Y) v_{-i\mu}}(g) p_\sigma(i\mu) \, d\mu \ .
    \]
    By Lemma \ref{Lemma2 about Asymptotic expansions}, there exists $ \eps > 0 $ such that we have
    \begin{align}\label{Proof of Prop Image of wave packet transform in Schwartz space eq1}
        & c_{\pi^{\sigma, i\mu}(X) \ext(\varphi_{i\mu}), \pi^{\tilde{\sigma}, -i\mu}(Y) v_{-i\mu}}(kah) \\
        = & a^{i\mu - \rho} \langle \ext(\varphi_{i\mu})(k), (\hat{J}^w_{\tilde{\sigma}, -i\mu} \pi^{\tilde{\sigma}, -i\mu}(Y) v_{-i\mu})(h^{-1}w^{-1}) \rangle
        \nonumber \\
        & \qquad + a^{-(i\mu + \rho)} \langle (\hat{J}^w_{\sigma, i\mu} \pi^{\sigma, i\mu}(X) \ext(\varphi_{i\mu}))(k), \pi^{\tilde{\sigma}, -i\mu}(Y) v_{-i\mu}(h^{-1}w) \rangle \nonumber \\
        & \qquad + a^{-(\rho + \eps)} R(i\mu, \pi^{\sigma, i\mu}(X) \ext(\varphi_{i\mu}), kah) \nonumber
    \end{align}
    for $ kM $ varying in a compact subset of $ \Omega $ and for $a$ sufficiently large and $ h \in K $, where the remainder function $ (\mu, kah) \mapsto R(i\mu, \pi^{\tilde{\sigma}, -i\mu}(Y) v_{-i\mu}, kah) $ is uniformly bounded. Let us write $ U $ as the union of a relatively compact set $ U_1 $ in $X$ and a set $ U_2 \subset X $ that is relatively compact in $ X \cup \Omega $ such that \eqref{Proof of Prop Image of wave packet transform in Schwartz space eq1} holds for any $ kaK \in U_2 $ and $ h \in K $.
    As $ \mu \mapsto \langle \ext(\varphi_{i\mu})(k), (\hat{J}^w_{\tilde{\sigma}, -i\mu} \pi^{\tilde{\sigma}, -i\mu}(Y) v_{-i\mu})(h^{-1}w^{-1}) \rangle $
    has compact support with respect to $ \mu $ and is smooth with respect to $ (\mu, k, h) $ (as long as $ kM $ varies in a compact subset of $ \Omega $), for all $ N \in \NN_0 $ there exists a constant $ C_N $ such that
    \begin{multline*}
        |\int_{\a^*_+} e^{i\mu \log a} \langle \ext(\varphi_{i\mu})(k), (\hat{J}^w_{\tilde{\sigma}, -i\mu} \pi^{\tilde{\sigma}, -i\mu}(Y) v_{-i\mu})(h^{-1}w^{-1}) \rangle \ p_\sigma(i\mu) \, d\mu| \\
        \leq C_N (1 + |\log a|)^{-N}
    \end{multline*}
    as the Fourier transform of a Schwartz function on $ \RR^n $ (here $ n = 1 $) is again a Schwartz function on $ \RR^n $.
    Let $W$ be a compact subset of $ \Omega $. Let $ \chi \in \Ccinf(\Omega) $ be a cut-off function which is equal to one on $W$. Then,
    \[
        \chi \hat{J}^w_{\sigma, i\mu} \pi^{\sigma, i\mu}(X) (1 - \chi) \ext(\varphi_{i\mu})
    \]
    is a smoothing operator by Lemma 5.3 of \cite[p.98]{BO00}. Since $ \ext(\varphi_{i\mu})(x) = \varphi_{i\mu}(x) $ for all $ x \in K(\Omega) $,
    \[
        \chi \hat{J}^w_{\sigma, i\mu} \pi^{\sigma, i\mu}(X) \ext(\varphi_{i\mu})
    \]
    is also a smoothing operator. So,
    \[
        \mu \mapsto \langle (\hat{J}^w_{\sigma, i\mu} \pi^{\sigma, i\mu}(X) \ext(\varphi_{i\mu}))(k), \pi^{\tilde{\sigma}, -i\mu}(Y) v_{-i\mu}(h^{-1}w) \rangle
    \]
    is smooth with respect to $k$ as long as $ kM $ varies in a compact subset of $ \Omega $.
    Since in addition it has compact support with respect to $ \mu $ and is also smooth with respect to $ (\mu, k, h) $ (as long as $ kM $ varies in a compact subset of $ \Omega $), for all $ N \in \NN_0 $ there exists a constant $ C'_N $ such that
    \begin{multline*}
        |\int_{\a^*_+} e^{-i\mu \log a} \langle (\hat{J}^w_{\sigma, i\mu} \pi^{\sigma, i\mu}(X) \ext(\varphi_{i\mu}))(k), \pi^{\tilde{\sigma}, -i\mu}(Y) v_{-i\mu}(h^{-1}w) \rangle p_\sigma(i\mu) \, d\mu| \\
        \leq C'_N (1 + |\log a|)^{-N} \ .
    \end{multline*}
    Thus,
    \[
        \LOI{p}{U}_{N, X, Y}(E(\varphi, v))
        \leq \LOI{p}{U_1}_{N, X, Y}(E(\varphi, v)) + \LOI{p}{U_2}_{N, X, Y}(E(\varphi, v)) < \infty
    \]
    for all $ N \in \NN_0 $.
    Hence, $ E(\varphi, v) $ belongs to $ \CS(\Gamma \bs G, \phi) $ as the seminorms $ \LOI{p}{V}_{N, X, Y} $ ($ V \in \U_\Gamma $, $ X, Y \in \U(\g) $, $ N \in \NN_0 $) are sufficient to define the topology on $ \CS(\Gamma \bs G, \phi) $.
\end{proof}

Recall the definition of $ v_{T, \lambda} $: For $ v \in V_{\tilde{\gamma}} $ ($ \gamma \in \hat{K} $), $ \lambda \in \a^*_\CC $ and $ T \in \Hom_M(V_\sigma, V_\gamma) $, we defined $ v_{T, \lambda} \in \Cinf(\dX, V(\tilde{\sigma}_{\lambda})) $ ($ \lambda \in \a^*_\CC $) by
\[
    v_{T, \lambda}(k) = \t{T} \tilde{\gamma}(k^{-1})v  \qquad (k \in K) \ .
\]
Let $ v_T \colon \a^*_+ \ni \mu \mapsto v_{T, -i\mu} $.
Set
\[
    \|v_T\| = \big(\int_K (\t{T} \tilde{\gamma}(k^{-1})v(k), \t{T} \tilde{\gamma}(k^{-1})v(k)) \, dk\big)^{\tfrac12} \ .
\]
Let $ \sigma \in \hat{M} $ and let $ \lambda \in i\a^* \smallsetminus \{0\} $.
Recall that we defined $ \pi^{\sigma, \lambda}(f) $ ($ f \in \CS(G, V_\phi) $) such that
\[
    \pi^{\sigma, \lambda}(F)v = \int_G F(g) \pi^{\sigma, \lambda}(g)v \, dg
\]
for all $ F \in \C_c(G, V_\phi) $ and $ v \in \Cinf(\dX, V(\sigma_\lambda)) $ (cf. \eqref{eq pi(f)v = int_G f(g) pi(g)v dg} and \eqref{eq (varphi, pi(f)v) = (c_(varphi, v), f)_G}).

\begin{lem}\label{Lem wave packet varphi_imu formula}
    Let $ \sigma \in \hat{M} $, $ \gamma \in \hat{K} $, $ \a^*_+ \ni \mu \mapsto \varphi_{i\mu} \in L^2(B, V_B(\sigma_{i\mu}, \phi)) $ be a measurable family, $ v \in V_{\tilde{\gamma}} $ and $ T \in \Hom_M(V_\sigma, V_\gamma) $. Let $ E(\varphi, v_T) \in \CS(\Gamma \bs G, \phi) $. Then,
    \[
        \varphi_{i\mu} = \frac1{\|v_T\|^2} \overline{\pi_*\big(\pi^{\tilde{\sigma}, -i\mu}(\chi \overline{E(\varphi, v_T)})(v_{T, -i\mu})\big)} \in \Cinf(B, V_B(\sigma_{i\mu}, \phi)) \ .
    \]
    Thus, $ \a^*_+ \ni \mu \mapsto \varphi_{i\mu} $ is a smooth family.
\end{lem}

\begin{proof}
    Let $ \sigma \in \hat{M} $, $ \psi \in \Hcal^\sigma_0(\phi) $. Let $ \varphi $ and $ v_T $ be as above. By \cite[p.156]{BO00}, $ (E(\psi, v_T), E(\varphi, v_T))_{L^2(\Gamma \bs G, \phi)} $ is equal to
    \[
        \|v_T\|^2 \cdot \int_{\a^*_+} (\psi_{i\mu}, \varphi_{i\mu})_B \ p_\sigma(i\mu) \, d\mu \ .
    \]
    But it is also equal to
    \begin{align*}
        & \int_{\Gamma \bs G} E(\psi, v_T)(g) \overline{E(\varphi, v_T)(g)} \, dg \\
        =& \int_G E(\psi, v_T)(g) \overline{(\chi E(\varphi, v_T))(g)} \, dg \\
        =& \int_G \int_{\a^*_+} \int_K \langle \ext(\psi_{i\mu})(k), \pi^{\tilde{\sigma}, -i\mu}(g) v_{T, -i\mu}(k) \rangle \, dk \ p_\sigma(i\mu) \, d\mu \ \overline{(\chi E(\varphi, v_T))(g)} \, dg \\
        =& \int_{\a^*_+} \int_K \langle \ext(\psi_{i\mu})(k), \int_G \overline{(\chi E(\varphi, v_T))(g)} \pi^{\tilde{\sigma}, -i\mu}(g) v_{T, -i\mu}(k) \, dg \rangle \, dk \ p_\sigma(i\mu) \, d\mu \\
        =& \int_{\a^*_+} \int_{\Gamma \bs \Omega} \langle \psi_{i\mu}(x), \pi_*\big(\pi^{\tilde{\sigma}, -i\mu}(\chi \overline{E(\varphi, v_T)})(v_{T, -i\mu})\big)(x) \rangle
        \ p_\sigma(i\mu) \, d\mu \ .
    \end{align*}
    Since $ \Hcal^\sigma_0(\phi) $ is dense in the Hilbert space $ \int^\oplus_{\a^*_+} L^2(B, V_B(\sigma_{i\lambda}, \phi)) p_\sigma(i\lambda) \, d\lambda $,
    \[
        \varphi_{i\mu} = \frac1{\|v_T\|^2} \overline{\pi_*\big(\pi^{\tilde{\sigma}, -i\mu}(\chi \overline{E(\varphi, v_T)})(v_{T, -i\mu})\big)} \ .
    \]
    The lemma follows.
\end{proof}

\begin{cor}\label{Cor wave packet varphi_imu formula for p_ac(f)}
    Let $ \gamma \in \hat{K} $ and $ f \in \CS(\Gamma \bs G, \phi)(\gamma) $. Then, $ p_{ac}(f) $, where $ p_{ac} $ denotes the orthogonal projection on $ L^2(\Gamma \bs G, \phi)_{ac} $, is a finite sum of wave packets of the form $ E(\varphi, v_T) \in L^2(\Gamma \bs G, \phi) $ with $ v \in V_{\tilde{\gamma}} $, $ \sigma \in \hat{M} $, $ T \in \Hom_M(V_{\tilde{\sigma}}, V_{\tilde{\gamma}}) $ and
    \[
        \varphi_{i\mu} := \frac1{\|v_T\|^2} \overline{\pi_*\big(\pi^{\tilde{\sigma}, -i\mu}(\chi \bar{f})(v_{T, -i\mu})\big)} \in \Cinf(B, V_B(\sigma_{i\mu}, \phi)) \ .
    \]
\end{cor}

\begin{proof}
    Let $ f \in \CS(\Gamma \bs G, \phi)(\gamma) $. Then, $ p_{ac}(f) $ is a finite sum of wave packets of the form
    \[
        E(\varphi, v_T) := \int_{\a^*_+} c_{\ext(\varphi_{i\mu}), v_{T, -i\mu}} \ p_\sigma(i\mu) \, d\mu \in L^2(\Gamma \bs G, \phi) \ ,
    \]
    where $ v \in V_{\tilde{\gamma}} $, $ \sigma \in \hat{M} $, $ T \in \Hom_M(V_{\sigma}, V_{\gamma}) $ and $ \a^*_+ \ni \mu \mapsto \varphi_{i\mu} \in L^2(B, V_B(\sigma_{i\mu}, \phi)) $ is a measurable family. Write
    \[
        p_{ac}(f) = \sum_i E(\varphi^{(i)}, v^{(i)}_{T_i})
    \]
    with $ T_i \in \Hom_M(V_{\sigma^{(i)}}, V_{\gamma}) $ ($ \sigma^{(i)} \in \hat{M} $), $ (v^{(k)}_{T_k}, v^{(l)}_{T_l}) = 0 $ if $ k \neq l $.
    Let $ \psi^{(j)} \in \Hcal^{\sigma^{(j)}}_0(\phi) $.

    By \cite[p.156]{BO00}, $ (E(\psi^{(j)}, v^{(j)}_{T_j}), p_{ac}(f))_{L^2(\Gamma \bs G, \phi)} $ is equal to
    \[
        \|v^{(j)}_{T_j}\|^2 \cdot \int_{\a^*_+} (\psi^{(j)}_{i\mu}, \varphi^{(j)}_{i\mu})_B \ p_\sigma(i\mu) \, d\mu \ .
    \]
    But it is also equal to
    \[
        \int_{\Gamma \bs G} E(\psi^{(j)}, v^{(j)}_{T_j})(g) \overline{f(g)} \, dg \ .
    \]
    As in the proof of the previous lemma, one can show that this is again equal to
    \[
        \int_{\a^*_+} \int_{\Gamma \bs \Omega} \langle \psi^{(j)}_{i\mu}(x), \pi_*\big(\pi^{\tilde{\sigma}^{(j)}, -i\mu}(\chi \bar{f})(v^{(j)}_{T_j, -i\mu})\big)(x) \rangle
        \ p_{\sigma^{(j)}}(i\mu) \, d\mu \ .
    \]
    The corollary follows as $ \Hcal^\sigma_0(\phi) $ is dense in the Hilbert space
    \[
        \int^\oplus_{\a^*_+} L^2(B, V_B(\sigma_{i\lambda}, \phi)) p_\sigma(i\lambda) \, d\lambda \ .
    \]
\end{proof}

\newpage

\subsubsection{Eisenstein integrals and its relation to the Poisson transform}

In this section, we recall the definition of \textit{Eisenstein integrals} and we show how the Poisson transform and the Eisenstein integrals are related. Moreover, we state several known results about Eisenstein integrals.

\begin{dfn}
    For $ A \in \Hom_M(V_\gamma, V_\tau(\sigma)) $ ($ \sigma \in \hat{M} $, $ \gamma, \tau \in \hat{K} $) and $ \lambda \in \a^*_\CC $, set
    \[
        E_{\gamma, \tau, \sigma, \lambda, A}(g) = \int_K a(g^{-1}k)^{-(\lambda + \rho)} \tau(\kappa(g^{-1}k)) A \gamma(k^{-1}) \, dk \ .
    \]
    These are the so-called \textit{Eisenstein integrals}.\index{Eisenstein integrals}\index[n]{EgammatausigmalambdaA@$ E_{\gamma, \tau, \sigma, \lambda, A} $}
\end{dfn}


Let $ \gamma \in \hat{K} $, $ d_\gamma = \dim(V_\gamma) $ and $ \{v^{(\gamma)}_1, \dotsc, v^{(\gamma)}_{d_\gamma}\} $ be an orthonormal basis of $ V_\gamma $.
For $ k \in K $, set $ \chi_\gamma(k) = \tr \gamma(k) $ ($ \chi_\gamma $ is called the \textit{character} of $ \gamma $). Let\index[n]{Egamma@$E_\gamma$}
\[
    E_\gamma = d_\gamma \int_K \overline{\chi_\gamma(k)} \pi(k) \, dk \ .
\]

\begin{lem}\label{Lem E_gamma w = sum_j F_j v_j}
    Let $ (\pi, W) $ be a Fréchet representation of $K$ and let $ \gamma \in \hat{K} $. Then, there are $ F_{\gamma, j, w} \in \Hom_K(V_\gamma, W) $ such that
    \[
        E_\gamma w = \sum_{j=1}^{d_\gamma} F_{\gamma, j, w}(v^{(\gamma)}_j) \ .
    \]
\end{lem}

\begin{proof}
    We identify $ V_\gamma $ with $ \CC^{d_\gamma} $ via this basis. Then, $ \overline{v^{(\gamma)}_j} = v^{(\gamma)}_j $ for all $j$.

    Let $ \{(v^{(\gamma)}_1)^*, \dotsc, (v^{(\gamma)}_{d_\gamma})^*\} $ be the dual basis of $ \{v^{(\gamma)}_1, \dotsc, v^{(\gamma)}_{d_\gamma}\} $.

    Then, $ (v^{(\gamma)}_i)^*( v^{(\gamma)}_j ) = (v^{(\gamma)}_i, v^{(\gamma)}_j) = \delta_{i, j} $.

    Let $ (\pi, W) $ be a Fréchet representation of $K$.
    For $ \gamma \in \hat{K} $, $ j = 1, \dotsc, d_\gamma $, $ w \in W $ and $ v \in V_{\tilde{\gamma}} $, set
    \[
        F_{\gamma, j, w}(v) = d_\gamma \int_K \pi(k)w \, (v^{(\gamma)}_j)^*\big(\gamma(k)^{-1} v\big) \, dk \in W
    \]
    ($W$ is a complete, locally convex, Hausdorff topological vector space).
    Since
    \[
        (v^{(\gamma)}_j)^*\big(\gamma(k)^{-1} v^{(\gamma)}_j\big) = (v^{(\gamma)}_j, \gamma(k)^{-1} v^{(\gamma)}_j)
        = \overline{(v^{(\gamma)}_j, \gamma(k) v^{(\gamma)}_j)} = \overline{\tr \gamma(k)} \ ,
    \]
    \[
        \sum_{j=1}^{d_\gamma} F_{\gamma, j, w}(v^{(\gamma)}_j) = d_\gamma \int_K \overline{\tr \gamma(k)} \pi(k)w \, dk = E_\gamma w \ .
    \]
    The lemma follows.
\end{proof}

For $ \lambda \in \a^*_\CC $, $ \sigma \in \hat{M} $, let $ H_\infty^{\sigma, \lambda} = \Cinf(\dX, V(\sigma_\lambda)) $.

\begin{prop}\label{Prop E_gamma P^T f in terms of Eisenstein integrals}
    Let $ \lambda \in \a^*_\CC $, $ \sigma \in \hat{M} $, $ \gamma, \tau \in \hat{K} $, $ f \in H_\infty^{\sigma, \lambda} $ and $ T \in \Hom_M(V_\sigma, V_\tau) $. Then,
    there are $ A_j \in \Hom_M(V_\gamma, V_\tau(\sigma)) $ such that
    \[
        E_\gamma (P^T_\lambda f)(g) = \sum_{j=1}^{d_\gamma} E_{\gamma, \tau, \sigma, \lambda, A_j}(g) v^{(\gamma)}_j \ .
    \]
\end{prop}

\begin{proof}
    Let $ g \in G $. Define $ F_{\gamma, j, P^T_\lambda f}(g) $ by
    \begin{align*}
        & d_\gamma \int_K (P^T_\lambda f)(k^{-1} g) (v^{(\gamma)}_j)^*\big( \gamma(k)^{-1} \cdot\big) \, dk \\
        = & d_\gamma \int_K a(g^{-1} h)^{-(\lambda+\rho)} \tau(\kappa(g^{-1} h)) \int_K T f(k^{-1} h) (v^{(\gamma)}_j)^*\big( \gamma(k)^{-1} \cdot\big) \, dh \, dk \\
        = & d_\gamma \int_K a(g^{-1} h)^{-(\lambda+\rho)} \tau(\kappa(g^{-1} h)) \int_K T f(k^{-1}) (v^{(\gamma)}_j)^*\big( \gamma(h k)^{-1} \cdot\big) \, dk \, dh \ .
    \end{align*}
    Set $ A_j := d_\gamma \int_K T f(k^{-1}) (v^{(\gamma)}_j)^*\big( \gamma(k)^{-1} \cdot\big) \, dk $.
    Let $ m \in M $ and $ v \in V_\gamma $. Then,
    \begin{align*}
        A_j(\gamma(m)v) &= d_\gamma \int_K T f(k^{-1}) (v^{(\gamma)}_j)^*\big( \gamma(k^{-1} m) v\big) \, dk \\
        &= d_\gamma \int_K T \sigma(m) f(k^{-1}) (v^{(\gamma)}_j)^*\big( \gamma(k)^{-1} v\big) \, dk
        = \tau(m) A_j(v) \ .
    \end{align*}
    Thus, $ A_j \in \Hom_M(V_\gamma, V_\tau(\sigma)) $. Hence, $ E_{\gamma, \tau, \sigma, \lambda, A_j} = F_{\gamma, j, P^T_\lambda f} $ is an Eisenstein integral and
    \[
        E_\gamma (P^T_\lambda f)(g) = \sum_{j=1}^{d_\gamma} E_{\gamma, \tau, \sigma, \lambda, A_j}(g) v^{(\gamma)}_j
    \]
    by the proof of Lemma \ref{Lem E_gamma w = sum_j F_j v_j}.
\end{proof}

Let $ \Phi $ be the canonical isomorphism between the principal series representation $ L^2(K \times_M V_\sigma) $ and $ H^{\sigma, \lambda} $. Then, $ \Phi(f)(k a n) = a^{\lambda - \rho} f(k) $ ($ f \in L^2(K \times_M V_\sigma) $) and $ \Phi^{-1} $ sends $ f \in H^{\sigma, \lambda} $ to $ \restricted{f}{K} $.

Let us give in the following an explicit description of the isomorphism
\[
    \Cinf(K \times_M V_\tau(\sigma)) \simeq \Hom_M(V_\sigma, V_\tau) \otimes H_\infty^{\sigma, \lambda} \ .
\]
Let $ T \in \Hom_M(V_\sigma, V_\tau) $ and let $ f \in H_\infty^{\sigma, \lambda} $. Then,
\[
    \restricted{T f}{K} \in \Cinf(K \times_M V_\tau(\sigma)) \ .
\]
Let now $ F \in \Cinf(K \times_M V_\tau(\sigma)) $.
Let $ V_\tau(\sigma) = V_{\sigma_1} \oplus \dotsm \oplus V_{\sigma_k} $ be an orthogonal direct sum decomposition of irreducible $M$-representations.
Let $ p_i $ be the canonical projection from $ V_\tau $ to $ V_{\sigma_i} $ and let $ \iota_i $ be the canonical inclusion from $ V_{\sigma_i} $ to $ V_\tau $.
Since $ V_\sigma \simeq V_{\sigma_i} $, there is a $M$-intertwining operator $ S_i $ from $ V_\sigma $ to $ V_{\sigma_i} $.
Let $ T_i = \iota_i \circ S_i \in \Hom_M(V_\sigma, V_\tau) $.
Define $ f_i \in H_\infty^{\sigma, \lambda} $ by $ f_i(k) = S_i^{-1}\big(p_i\big(F(k)\big)\big) $ ($ k \in K $). Then,
\[
    F = \sum_{i=1}^k \restricted{T_i f_i}{K}
\]
as $ \Id = \sum_{i=1}^k \iota_i \circ p_i $ on $ V_\tau(\sigma) $. So, the isomorphism maps $F$ to $ \sum_{i=1}^k T_i \otimes f_i $.

\begin{lem}
    Let $ E_{\gamma, \tau, \sigma, \lambda, A} $ be an Eisenstein integral and $ v \in V_\gamma $. Then, there are $ T_i \in \Hom_M(V_\sigma, V_\tau) $ and $ f_{i, v} \in H_\infty^{\sigma, \lambda} $ such that
    \[
        E_{\gamma, \tau, \sigma, \lambda, A}(m a) v = \sum_i P_\lambda^{T_i} f_{i, v}(ma)  \qquad (m \in M, a \in A) \ .
    \]
\end{lem}

\begin{proof}
    For $ v \in V_\gamma $ and $ k \in K $, set $ F(k) = A \gamma(k^{-1}) v $. Then,
    \[
        F \in \Cinf(K \times_M V_\tau(\sigma)) \simeq \Hom_M(V_\sigma, V_\tau) \otimes H_\infty^{\sigma, \lambda} \ .
    \]
    Let $ f_{i, v} = S_i^{-1} \circ p_i \circ F $ and let $ T_i = \iota_i \circ S_i \in \Hom_M(V_\sigma, V_\tau) $. It follows from the above that
    \[
        \sum_i T_i f_{i, v}(k) = A(\gamma(k)^{-1} v) \ .
    \]
    The lemma follows now easily from this.
\end{proof}

Let $ (\gamma, V_\gamma) $ be a finite-dimensional representation of $K$.
Let $ c_\gamma \colon \a^*_\CC \to \End_M(V_\gamma) $ be the meromorphic function which is given by\index[n]{cgamma@$ c_\gamma(\lambda) $}
\[
    c_\gamma(\lambda) := \int_{\bar{N}}{ a(\bar{n})^{-(\lambda+\rho)} \gamma(\kappa(\bar{n})) \, d\bar{n} }
\]
for $ \R(\lambda) > 0 $. It is the Harish-Chandra $c$-function.

\begin{prop}\label{Prop Asymptotic expansion of Eisenstein integrals}
    Let $ E_{\gamma, \tau, \sigma, i\lambda, A} $ be an Eisenstein integral with $ \lambda \in \a^* \smallsetminus \{0\} $ and let $ v \in V_\gamma $. Then, there is $ \eps > 0 $ such that, for $ a \to \infty $,
    \[
        E_{\gamma, \tau, \sigma, i\lambda, A}(m a) v = a^{i\lambda - \rho} c_\tau(i\lambda) A(\gamma(m)^{-1} v)
        + a^{-i\lambda - \rho} A^w\big(c_\gamma(-i\lambda) \gamma(m)^{-1} v\big)
        + O(a^{-\alpha - \rho -\eps})
    \]
    uniformly in $ m \in M $.
\end{prop}

\begin{proof}
    For $ \gamma \in \hat{K} $ and $ \R(\mu) < 0 $, define
    \[
        j_{\gamma, \mu} = \int_{\bar{N}} a(\bar{n})^{\mu - \rho} \gamma(\kappa(\bar{n}))^{-1} \, d\bar{n}
        \in \End_M V_\gamma
    \]
    (this corresponds to Definition 2.17 in \cite[p.31]{Olb94}). By analytic continuation, we get a meromorphic family on $ \a^*_\CC $.
    By the previous lemma,
    \[
        E_{\gamma, \tau, \sigma, \lambda, A}(m a) v = \sum_i P_\lambda^{T_i} f_{i, v}(ma)  \qquad (m \in M, a \in A) \ .
    \]
    Let $ \lambda \in \a^* \smallsetminus \{0\} $. Then, there is $ \eps > 0 $ such that, for $ a \to \infty $,
    \[
        \sum_j P^{T_j} f_{j, v}(ka) = a^{i\lambda - \rho} c_\tau(i\lambda) \sum_j T_j f_{j, v}(k)
        + a^{-i\lambda - \rho} \sum_j T_j^w (\hat{J}_{\sigma, i\lambda} f_{j, v})(k)
        + O(a^{-\alpha - \rho -\eps})
    \]
    uniformly in $ kM \in \dX $, by Lemma 6.2, 1., of \cite[p.110]{BO00} (compare with Lemma \ref{Lemma2 about Asymptotic expansions}). If $ k = m \in M $, then
    \[
        \sum_j T_j^w (\hat{J}_{\sigma, i\lambda} f_{j, v})(m)
    \]
    is equal to
    \[
        \tau(m)^{-1} \tau(w) \sum_j T_j (\hat{J}^w_{\sigma, i\lambda} f_{j, v})(e)
        = \tau(m)^{-1} \tau(w) \sum_j \iota_j (\hat{J}^w_{\sigma_j, i\lambda} \big(p_j(F)\big))(e) \ .
    \]
    This is again equal to
    \[
        \tau(m)^{-1} \tau(w) A\big(j_{\gamma, i\lambda} \gamma(w)^{-1} v\big)
        = A^w\big(j_{\gamma, i\lambda}^w \gamma(m)^{-1} v\big) \ .
    \]
    By Proposition 3.13 of \cite[p.31]{Olb94},
    \begin{equation}\label{eq int_(bar N) a(bar n)^(lambda - rho) gamma(kappa(bar n))^(-1) d bar n = gamma(w)^(-1) c_gamma(-lambda) gamma(w)}
        j_{\gamma, -\mu}^w = c_\gamma(\mu) \qquad (\mu \in \a^*_\CC) \ .
    \end{equation}
    The proposition follows.
\end{proof}

\begin{cor}\label{Cor Comparison of terms in the asymptotic expansions}
    Let $ \gamma, \tau \in \hat{K} $, $ \sigma \in \hat{M} $. Let $ \{v^{(\gamma)}_j\} $ is an orthonormal basis of $ V_{\gamma} $. Let $ \lambda \in \a^* \smallsetminus \{0\} $, $ \psi_{i\lambda} \in H_\infty^{\sigma, i\lambda} $ and $ T \in \Hom_M(V_\sigma, V_\tau) $. Let $ k \in K $. Then, there are $ A_{k, \gamma, j} \in \Hom_M(V_{\gamma}, V_\tau(\sigma)) $ such that
    \[
        c_\tau(i\lambda) T \big(E_{\gamma} \pi^{\sigma, i\lambda}(k^{-1}) \psi_{i\lambda} \big)(e) = c_\tau(i\lambda) \sum_{j=1}^{d_{\gamma}} A_{k, \gamma, j} v^{(\gamma)}_{j}
    \]
    and
    \[
        T^w \hat{J}_{\sigma, -i\lambda} \big(E_{\gamma} \pi^{\sigma, i\lambda}(k^{-1}) \psi_{i\lambda} \big)(e)
        = \sum_{j=1}^{d_{\gamma}} A_{k, \gamma, j}^w c_{\gamma}(i\lambda) v^{(\gamma)}_{j} \ .
    \]
\end{cor}

\begin{rem}
    Let $ E_{\gamma, \tau, \sigma, \lambda, A} $ be an Eisenstein integral. Then, one can show that
    \[
        E_{\gamma, \tau, w \sigma, -\lambda, c_\tau(\lambda)^w A^w} = E_{\gamma, \tau, \sigma, \lambda, A c_\gamma(\lambda)}
    \]
    whenever $ c_\tau(\mu) $ and $ c_\gamma(\mu) $ are regular in $ \mu = \lambda $. Compare with Theorem 13.2.9 of \cite[p.232]{Wallach2}.
\end{rem}

\begin{proof}
    Let $ \lambda \in \a^* \smallsetminus \{0\} $ and let $ \psi_{i\lambda} \in H_\infty^{\sigma, i\lambda} $. Let $ k \in K $. Then, there are $ A_{k, \gamma, j} \in \Hom_M(V_{\gamma}, V_\tau(\sigma)) $ such that
    \[
        E_{\gamma} (P^T_{i\lambda} \pi^{\sigma, i\lambda}(k^{-1}) \psi_{i\lambda})(g) = \sum_{j=1}^{d_{\gamma}} E_{\gamma, \tau, \sigma, i\lambda, A_{k, \gamma, j}}(g) v^{(\gamma)}_{j}  \qquad (g \in G)
    \]
    by Proposition \ref{Prop E_gamma P^T f in terms of Eisenstein integrals}.
    By Lemma 6.1, 1., of \cite[p.110]{BO00}, Proposition \ref{Prop Asymptotic expansion of Eisenstein integrals} (applied twice) and by uniqueness of asymptotic expansions,
    \[
        c_\tau(i\lambda) T \big(E_{\gamma} \pi^{\sigma, i\lambda}(k^{-1}) \psi_{i\lambda} \big)(e) = c_\tau(i\lambda) \sum_{j=1}^{d_{\gamma}} A_{k, \gamma, j} v^{(\gamma)}_{j}
    \]
    and
    \[
        T^w \hat{J}_{\sigma, -i\lambda} \big(E_{\gamma} \pi^{\sigma, i\lambda}(k^{-1}) \psi_{i\lambda} \big)(e)
        = \sum_{j=1}^{d_{\gamma}} A_{k, \gamma, j}^w c_{\gamma}(i\lambda) v^{(\gamma)}_{j} \ .
    \]
\end{proof}

\begin{thm}\label{Thm 13.3.2 for Eisenstein integrals}
    Let $ E_{\gamma, \tau, \sigma, i\nu, A} $ be an Eisenstein integral and let $ m \in M $. Then,
    \[
        \int_A \int_N \int_{\a^*_+} E_{\gamma, \tau, \sigma, i\nu, A}(nma) p_\sigma(i\nu) \, d\nu \, a^{-i\mu-\rho} \, dn \, da \qquad (\mu \in \a^*)
    \]
    is equal to
    \[
        \frac1{2\pi} c_\tau(i\mu)^{-1} \big( a^{i\mu - \rho} c_\tau(i\mu) A \gamma(m)^{-1} + a^{-i\mu - \rho} A^w c_\gamma(i\mu) \gamma(m)^{-1} \big) \ .
    \]
\end{thm}

\begin{rem}
    It follows from the proof (look in particular at the formula for $ \Phi(P, \alpha)^Q(ma) $ on p.239 of \cite{Wallach2}) that the above formula is deduced from the following equality:
    \begin{multline*}
        a^{-\rho} \int_N \int_{\a^*_+} E_{\gamma, \tau, \sigma, i\mu, A}(nma) p_\sigma(i\mu) \, d\mu \, dn \\
        = \int_{\a^*_+} \Big( c_\tau(i\mu)^{-1} \big( c_\tau(i\mu) A \gamma(m)^{-1} + A^w c_\gamma(i\mu) \gamma(m)^{-1} \big) a^{i\mu} \Big) d\mu  \qquad (a \in A) \ .
    \end{multline*}
\end{rem}

\begin{proof}
    For $ m \in M $, set $ \phi(m) = A(\gamma(m)^{-1} \, \cdot) $.
    For $ \varphi \in \Hom(V_\gamma, V_\tau) $, set
    \begin{enumerate}
    \item $ C_{P|P}(e, \lambda)\varphi = c_\tau(\lambda) \varphi $,
    \item $ C_{P|P}(w, \lambda)\varphi = \varphi^w c_\gamma(-\lambda) $,
    \item $ C_{P|P}(e, w^{-1}\lambda)\varphi = c_\tau(-\lambda) \varphi $, and
    \item $ C_{P|P}(w, w^{-1}\lambda)\varphi = \varphi^w c_\gamma(\lambda) $
    \end{enumerate}
    whenever this is well-defined. This notation makes it easier to compare N.~Wallach's results with ours.

    Assume now that $ \lambda \in i \a^* \smallsetminus \{0\} $. By Proposition \ref{Prop Asymptotic expansion of Eisenstein integrals}, there is $ \eps > 0 $ such that, for $ a \to \infty $,
    \begin{multline}\label{Proof of Prop Functional equation of Eisenstein integrals eq1}
        E_{\gamma, \tau, \sigma, \lambda, A}(m a) v
        = a^{\lambda - \rho} \big(C_{P|P}(e, \lambda) \phi(m)\big) v \\
        + a^{-\lambda - \rho} \big(C_{P|P}(w, \lambda) \phi(m)\big) v + O(a^{-\alpha - \rho -\eps})
    \end{multline}
    and
    \begin{multline}\label{Proof of Prop Functional equation of Eisenstein integrals eq2}
        E_{\gamma, \tau, \sigma, -\lambda, A}(m a) v
        = a^{\lambda - \rho} \big(C_{P|P}(w, w^{-1}\lambda) \phi(m)\big) v \\
        + a^{-\lambda - \rho} \big(C_{P|P}(e, w^{-1}\lambda) \phi(m)\big) v + O(a^{-\alpha - \rho -\eps})
    \end{multline}
    uniformly in $ m \in M $. Compare this with Theorem 13.2.6 of \cite[p.230]{Wallach2}.

    Set $ \LOI{C}{0}_{P|P}(e, \lambda) = C_{P|P}(e, \lambda)^{-1} C_{P|P}(e, \lambda) = \Id $ and, for $ \varphi \in \Hom(V_\gamma, V_\tau) $, set
    \[
        \LOI{C}{0}_{P|P}(w, w\lambda)\varphi = C_{P|P}(e, \lambda)^{-1} C_{P|P}(w, w^{-1}\lambda) \varphi = c_\tau(\lambda)^{-1} \varphi^w c_\gamma(\lambda) \ .
    \]
    The theorem follows now from Theorem 13.3.2 of \cite[p.236]{Wallach2}.

    Since $ \int_{\bar{N}} a(\bar{n})^{-2\rho} \, d\bar{n} = 1 $ and
    \[
        2\pi \int_A \int_{\a^*} f(\nu) a^{i\nu} \, d\nu \, a^{-i\mu} \, da = f(\mu)
    \]
    for all $ \mu \in \a^* $ and $ f \in \Ccinf(\a^*) $, the appearing multiplicative constant is $ \frac1{2\pi} $.
\end{proof}

\newpage

\subsubsection{The contribution of the wave packets of Eisenstein series to the closure of the space of cusp forms}\label{ssec:The contrib of the WPs to the closure of the space of cusp forms}

We show in this section that the wave packets of Eisenstein series are orthogonal to the space of cusp forms (cf. Theorem \ref{Thm L^2(Gamma|G, phi)_ac is orthogonal to the space of cusp forms}). We conclude from this and our previous results that
\[
    \overline{\c{\CS}(\Gamma \bs G, \phi)} = L^2(\Gamma \bs G, \phi)_{ds} \oplus L^2(\Gamma \bs G, \phi)_U
\]
(Theorem \ref{Thm decomposition of closure of cusp forms}).

\medskip

Assume that $ X \neq \OO H^2 $. 

\begin{lem}\label{Lem about Knapp-Stein intertwining operator}
    Fix $ g \in G(\Omega) $, $ \sigma \in \hat{M} $ and $ \lambda \in \a^*_\CC $. Then,
    \[
        f \in \C^{-\infty}_\Omega(\dX, V(\sigma_\lambda, \phi))^\Gamma \mapsto \hat{J}_{\sigma, \lambda}f(g)
    \]
    is continuous.
\end{lem}

\begin{proof}
    Fix $ g \in G(\Omega) $ and $ \lambda \in \a^*_\CC $. Let $Q$ be an open neighbourhood of $ \Lambda $ that does not contain $ gP $ and let $ U = \dX \smallsetminus Q $. Let $ U' $ be a compact subset of $U$ containing $ gP $ and let $ \chi \in \Ccinf(U) $ be such that $ \restricted{ \chi }{ U' } \equiv 1 $.
    Then,
    \[
        f \in \C^{-\infty}_\Omega(\dX, V(\sigma_\lambda, \phi))^\Gamma \mapsto \hat{J}_{\sigma, \lambda}(\chi f)(g) 
    \]
    is continuous as multiplication with $ \chi $ is a smooth operation and as
    \[
        \hat{J}_{\sigma, \lambda} \colon \Cinf(\dX, V(\sigma_\lambda, \phi)) \to \Cinf(\dX, V(\sigma^w_{-\lambda}, \phi))
    \]
    is continuous by Lemma 5.1 and Lemma 5.2 of \cite[pp.96-97]{BO00}. The map
    \[
        f \in \C^{-\infty}_\Omega(\dX, V(\sigma_\lambda, \phi))^\Gamma \mapsto \hat{J}_{\sigma, \lambda}((1 - \chi)f)(g) 
    \]
    is also continuous by Lemma 5.3 of \cite[p.98]{BO00} and as multiplication with $ 1 - \chi $ is a continuous operation.
    Thus,
    \[
        f \in \C^{-\infty}_\Omega(\dX, V(\sigma_\lambda, \phi))^\Gamma \mapsto \hat{J}_{\sigma, \lambda}f(g) 
    \]
    is continuous. This completes the proof of the lemma.
\end{proof}

For $ T \in \Hom_M(V_\sigma, V_\gamma) $, define $ T^w = \gamma(w) T \sigma(w)^{-1} $. Then, $ T^w $ belongs also to $ \Hom_M(V_\sigma, V_\gamma) $.
Note that $ T^w $ does not depend on the choice of $w$.\index[n]{Tw@$ T^w $}

\begin{thm}\label{Theorem 13.3.2}
    Let $ \varphi \in \Hcal^{\sigma}_0(\phi) $ and let $ k \in K(\Omega) $. 
    Let $ T \in \Hom_M(V_\sigma, V_\gamma) $. Then,
    \begin{equation}\label{Theorem 13.3.2 eq1}
        \int_A \int_N \int_{\a^*_+} P^T_{i\nu} \ext(\varphi_{i\nu})(kna) p_\sigma(i\nu) \, d\nu \, a^{-i\mu-\rho} \, dn \, da \qquad (\mu \in \a^*)
    \end{equation}
    is equal to
    \begin{equation}\label{Theorem 13.3.2 eq2}
        \frac1{2\pi} c_\gamma(i\mu)^{-1} \big(c_\gamma(i\mu) T\varphi_{i\mu}(k) + T^w \hat{J}_{\sigma, -i\mu} \ext(\varphi_{-i\mu})(k) \big) \ .
    \end{equation}
\end{thm}

\begin{proof}
    Fix $ \mu \in \a^* $ and $ k \in K(\Omega) $. Put $ \psi_{i\nu} = \ext(\varphi_{i\nu}) $. It suffices to show that
    \begin{multline}\label{Proof Theorem 13.3.2 eq1}
        a^{-\rho} \int_N \int_{\a^*_+} P^T_{i\nu} \psi_{i\nu}(kna) p_\sigma(i\nu) \, d\nu \, dn \\
        = \int_{\a^*} c_\gamma(i\nu)^{-1} \big(c_\gamma(i\nu) T \varphi_{i\nu}(k) + T^w \hat{J}_{\sigma, -i\nu} \psi_{-i\nu}(k) \big) a^{i\nu} \, d\nu \ .
    \end{multline}
    Indeed,
    \begin{equation}\label{Proof Theorem 13.3.2 eq2}
        2\pi \int_A \int_{\a^*} f(\nu) a^{i\nu} \, d\nu \, a^{-i\mu} \, da = f(\mu)
    \end{equation}
    for all $ \mu \in \a^* $ and $ f \in \Ccinf(\a^*) $. It follows from \eqref{Proof Theorem 13.3.2 eq1} that \eqref{Theorem 13.3.2 eq1} is equal to
    \[
        \int_A \int_{\a^*} c_\gamma(i\nu)^{-1} \big(c_\gamma(i\nu) T \varphi_{i\nu}(k) + T^w \hat{J}_{\sigma, -i\nu} \psi_{-i\nu}(k) \big) a^{i\nu} \, d\nu \, a^{-i\mu} \, da \ .
    \]
    By \eqref{Proof Theorem 13.3.2 eq2} and as $ \varphi \in \Hcal^{\sigma}_0(\phi) $, this is equal to \eqref{Theorem 13.3.2 eq2}.

    Let us prove first \eqref{Proof Theorem 13.3.2 eq1} for $ \Gamma = \{e\} $. Then, $ \psi_{i\nu} \in \Cinf(\dX, V(\sigma_{i\nu}, \phi)) $ for all $ \nu \in \a^*_+ $.
    Let $ \gamma_1, \gamma_2, \dotsc $ be the elements of $ \hat{K} $. Then,
    \begin{multline*}
        \int_N |\int_{\a^*_+} P^T_{i\nu} \big(\sum_{j=1}^l E_{\gamma_j} \pi^{\sigma, i\nu}(k^{-1}) \psi_{i\nu} \big)(na) p_\sigma(i\nu) \, d\nu \\
        - \int_{\a^*_+} P^T_{i\nu} \pi^{\sigma, i\nu}(k^{-1}) \psi_{i\nu}(na) p_\sigma(i\nu) \, d\nu| \, dn
    \end{multline*}
    is equal to
    \[
        \int_N |\sum_{j=1}^l E_{\gamma_j} \int_{\a^*_+} P^T_{i\nu} \psi_{i\nu}(kna) p_\sigma(i\nu) \, d\nu
        - \int_{\a^*_+} P^T_{i\nu} \psi_{i\nu}(kna) p_\sigma(i\nu) \, d\nu| \, dn \ .
    \]
    By Theorem 7.1.1 of \cite[p.227]{Wallach} and by the proof of Proposition \ref{Prop f^Omega well-defined}, there exist constants $ C_{r, l} $ converging to zero when $ r > 1 $ is fixed and $l$ goes to $ \infty $ such that this is less or equal than
    \[
        C_{r, l} a^\rho (1 + |\log a|)^{-r} \ .
    \]
    Thus, the left-hand side of \eqref{Proof Theorem 13.3.2 eq1} is equal to
    \begin{equation}\label{Proof Theorem 13.3.2 eq3}
        \lim_{l \to \infty} a^{-\rho} \int_N \int_{\a^*_+} P^T_{i\nu} \big(\sum_{j=1}^l E_{\gamma_j} \pi^{\sigma, i\nu}(k^{-1}) \psi_{i\nu} \big)(na) p_\sigma(i\nu) \, d\nu \ .
    \end{equation}
    Let $ \{v^{(\gamma_j)}_i\} $ be an orthonormal basis of $ V_{\gamma_j} $.
    By Corollary \ref{Cor Comparison of terms in the asymptotic expansions}, there are $ A_{k, j, j'} \in \Hom_M(V_{\gamma_j}, V_\gamma(\sigma)) $ such that
    \[
        c_\gamma(i\nu) T \big(E_{\gamma_j} \pi^{\sigma, i\nu}(k^{-1}) \psi_{i\nu} \big)(e) = c_\gamma(i\nu) \sum_{j'=1}^{d_{\gamma_j}} A_{k, j, j'} v^{(\gamma_j)}_{j'}
    \]
    and
    \[
        T^w \hat{J}_{\sigma, -i\nu} \big(E_{\gamma_j} \pi^{\sigma, i\nu}(k^{-1}) \psi_{i\nu} \big)(e) = \sum_{j'=1}^{d_{\gamma_j}} A_{k, j, j'}^w c_{\gamma_j}(i\nu) v^{(\gamma_j)}_{j'}
    \]
    for every $ j \in \NN $, $ \nu \in \a^* \smallsetminus \{0\} $ and $ k \in K $.
    By Theorem \ref{Thm 13.3.2 for Eisenstein integrals}, \eqref{Proof Theorem 13.3.2 eq3} yields
    \begin{multline*}
        \lim_{l \to \infty} \int_{\a^*} c_\gamma(i\nu)^{-1} \big(c_\gamma(i\nu) T \big(\sum_{j=1}^l E_{\gamma_j} \pi^{\sigma, i\nu}(k^{-1}) \psi_{i\nu} \big)(e) \\
        + T^w \hat{J}_{\sigma, -i\nu} \big(\sum_{j=1}^l E_{\gamma_j} \pi^{\sigma, i\nu}(k^{-1}) \psi_{i\nu} \big)(e) \big) a^{i\nu} \, d\nu \ .
    \end{multline*}
    \eqref{Proof Theorem 13.3.2 eq1} now follows.
    \\ Let us now show \eqref{Proof Theorem 13.3.2 eq1} for $ \Gamma \neq \{e\} $.
    Let $ (\phi_j)_j $ be a Dirac sequence in $ \Ccinf(G, \RR) $. As $ P^T_{i\mu} \pi^{\sigma, i\mu}(\phi_j) \psi_{i\mu} $ converges to $ P^T_{i\mu} \psi_{i\mu} $ and as $ \pi^{\sigma, i\mu}(\phi_j)\psi_{i\mu} \in \Cinf(\dX, V(\sigma_{i\mu}, \phi)) $ (since $ \phi_j \in \Ccinf(G) $) and as it depends smoothly on $ \mu \in \a^*_+ $, the function
    \[
        \int_{\a^*_+} P^T_{i\mu} \pi^{\sigma, i\mu}(\phi_j) \psi_{i\mu} p_\sigma(i\mu) \, d\mu
    \]
    is a wave packet in $X$ and belongs therefore to $ \CS(X, V_X(\gamma, \phi)) $. By the proven case $ \Gamma = \{e\} $,
    \begin{equation}\label{Proof Theorem 13.3.2 eq4}
        a^{-\rho} \int_N \int_{\a^*_+} P^T_{i\nu} \pi^{\sigma, i\nu}(\phi_j) \psi_{i\nu}(kna) p_\sigma(i\nu) \, d\nu \, dn
    \end{equation}
    is equal to
    \begin{equation}\label{Proof Theorem 13.3.2 eq5}
        \int_{\a^*} c_\gamma(i\nu)^{-1} \big(c_\gamma(i\nu) T \pi^{\sigma, i\nu}(\phi_j) \varphi_{i\nu}(k) + T^w \hat{J}_{\sigma, -i\nu} \pi^{\sigma, -i\nu}(\phi_j) \psi_{-i\nu}(k) \big) a^{i\nu} \, d\nu \ .
    \end{equation}
    One can show that \eqref{Proof Theorem 13.3.2 eq4} converges to
    \[
        a^{-\rho} \int_N \int_{\a^*_+} P^T_{i\nu} \psi_{i\nu}(kna) p_\sigma(i\nu) \, d\nu \, dn
    \]
    by using Lemma \ref{LemDirac} and the proof of Proposition \ref{Prop f^Omega well-defined}.
    As $ \varphi \in \Hcal^{\sigma}_0(\phi) $, \eqref{Proof Theorem 13.3.2 eq5} converges to
    \[
        \int_{\a^*} c_\gamma(i\nu)^{-1} \big(c_\gamma(i\nu) T \varphi_{i\nu}(k) +  T^w \hat{J}_{\sigma, -i\nu} \psi_{-i\nu}(k) \big) a^{i\nu} \, d\nu
    \]
    by Lemma \ref{Lem about Knapp-Stein intertwining operator}. The theorem follows.
\end{proof}

The following is the crucial result needed to prove that $ L^2(\Gamma \bs G, \phi)_{ac} $ is orthogonal to the space of cusp forms (Theorem \ref {Thm L^2(Gamma|G, phi)_ac is orthogonal to the space of cusp forms}).

\begin{prop}\label{Prop no non nonzero compactly supported wave packet is cuspidal}
    No non nonzero compactly supported wave packet is cuspidal.
\end{prop}

\begin{proof}
    The main result we need in order to prove this proposition is Theorem \ref{Theorem 13.3.2}, providing an explicit formula for the Fourier transform of the constant term of a compactly supported wave packet.

    Let $ f \in L^2(\Gamma \bs G, \phi)_{ac} \cap \c{\CS}(\Gamma \bs G, \phi) $ be a cuspidal, compactly supported wave packet.
    Since the projection of a function $ g \in \c{\CS}(\Gamma \bs G, \phi) $ to $ L^2(\Gamma \bs G, \phi)(\gamma) $ belongs still to $ \c{\CS}(\Gamma \bs G, \phi) $ for any $ \gamma \in \hat{K} $, we may assume without loss of generality that $ f \in (L^2(\Gamma \bs G, \phi)_{ac})(\gamma) \cap \c{\CS}(\Gamma \bs G, \phi) $ for some $ \gamma \in \hat{K} $.
    We have
    \[
        \c{\CS}(\Gamma \bs G, \phi)(\gamma) \simeq [\c{\CS}(\Gamma \bs G, \phi) \otimes V_{\tilde{\gamma}}]^K \otimes V_\gamma
        = \c{\CS}(Y, V_Y(\tilde{\gamma}, \phi)) \otimes V_\gamma \ ,
    \]
    where $ V_Y(\tilde{\gamma}, \phi) := V_Y(\tilde{\gamma}) \otimes V_\phi $.
    The isomorphism
    \[
        \c{\CS}(Y, V_Y(\tilde{\gamma}, \phi)) \otimes V_\gamma \to \c{\CS}(\Gamma \bs G, \phi)(\gamma)
    \]
    maps $ F \otimes v $ to $ \langle F, v \rangle $.
    Thus, $f$ is a finite sum 
    of wave packets of the form $ \langle E(\varphi, T), v \rangle $ with $ T \in \Hom_M(V_{\tilde{\sigma}}, V_{\tilde{\gamma}}) $ for some $ \sigma \in \hat{M} $, $ v \in V_\gamma $, $ \varphi \colon \lambda \mapsto \varphi_{i\lambda} \in L^2(B, V_B(\tilde{\sigma}_{i\lambda}, \phi)) $ is a measurable family with compact support in $ \lambda $ and $  E(\varphi, T) \in \c{\CS}(Y, V_Y(\tilde{\gamma}, \phi)) $.

    Let $ \{T_1^{(r)}, \dotsc, T_{n_r}^{(r)}\} $ be an orthonormal basis of $ \Hom_M(V_{\tilde{\sigma}^{(r)}}, V_{\tilde{\gamma}}) $ and let $ v_1, \dotsc, v_s $ be a basis of $ V_\gamma $. Write
    \[
        f = \sum_{r=1}^s \sum_{q=1}^{n_r} \langle \int_{\a^*_+} P^{T_q^{(r)}}_{i\nu} \ext\bigl(\varphi^{(q, r)}_{i\nu}\bigr) p_{\sigma^{(r)}}(i\nu) \, d\nu, v_r \rangle \ .
    \]
    Since the sum of $ \I(T_i^{(r)}) $ over $i$ is direct for $r$ fixed,
    \begin{equation}\label{Proof of Prop L^2(Gamma|G, phi)_ac cap c(CS)(Gamma|G, phi) = 0 eq1}
        \sum_{q=1}^{n_r} \int_{\a^*_+} P^{T_q^{(r)}}_{i\nu} \ext\bigl(\varphi^{(q, r)}_{i\nu}\bigr) p_{\sigma^{(r)}}(i\nu) \, d\nu
    \end{equation}
    belongs also to $ \c{\CS}(Y, V_Y(\tilde{\gamma}, \phi)) $ for $ r = 1, \dotsc, s $.
    It follows from Lemma \ref{Lem wave packet varphi_imu formula} that $ \varphi^{(q, r)}_{i\nu} \in \Cinf(B, V_B(\tilde{\sigma}^{(r)}_{i\nu}, \phi)) $ for all $ \nu \in \a^*_+ $ and that $ \nu \in \a^*_+ \mapsto \varphi^{(q, r)}_{i\nu} $ is smooth.
    \\ By Theorem \ref{Theorem 13.3.2}, the Fourier transform of the constant term of \eqref{Proof of Prop L^2(Gamma|G, phi)_ac cap c(CS)(Gamma|G, phi) = 0 eq1} at $ i\mu $ ($ \mu \in \a^*_+ $) is equal to
    \[
        \frac1{2\pi} \sum_{q=1}^{n_r} T_q^{(r)} \varphi^{(q, r)}_{i\mu}(k)
    \]
    for all $ k \in K(\Omega) $, since $ \varphi^{(q, r)}_{-i\mu} = 0 $ for all $ \mu \in \a^*_+ $.
    By \eqref{Proof of Prop L^2(Gamma|G, phi)_ac cap c(CS)(Gamma|G, phi) = 0 eq1},
    \[
        \sum_{q=1}^{n_r} T_q^{(r)} \varphi^{(q, r)}_{i\mu}(k) = 0
    \]
    for all $ k \in K(\Omega) $.
    Applying $ (T_q^{(r)})^* $ yields $ \varphi^{(q, r)}_{i\mu} = 0 $ for $ q = 1, \dotsc, n_r $.
    Hence,
    \[
        f = \sum_{r=1}^s \sum_{q=1}^{n_r} \langle \int_{\a^*_+} P^{T_q^{(r)}}_{i\nu} \ext\bigl(\varphi^{(q, r)}_{i\nu}\bigr) p_{\sigma^{(r)}}(i\nu) \, d\nu, v_r \rangle = 0 \ .
    \]
    The proposition follows.
\end{proof}

Let $ \sigma \in \hat{M} $. Define $ N^{\sigma}_{\pi} $ ($ \pi \in \hat{G} $) by $ \begin{cases}
    N_\pi   & \text{if } \pi = \pi^{\sigma, i \lambda} \text{ for some } \lambda \in \a^*_+ \\
    \{0\}   & \text{otherwise}
\end{cases} \ . $

Let $ \F_\sigma $ be the restriction of $ \F $ to the closed $G$-invariant space $ \F^{-1}\big(\int^\oplus_{\hat{G}} N^{\sigma}_\pi \hat{\otimes} V_\pi \, d\kappa(\pi)\big) $.

Let $ \sigma \in \hat{M} $ and $ \chi_\sigma \in \Ccinf(\a^*_+) $ be a cut-off function.
For $ f \in L^2(\Gamma \bs G, \phi) $, define
\[
    \Psi_\sigma(f) = \F_\sigma^{-1}\big( (\chi_\sigma(\lambda) \F_{\pi^{\sigma, i\lambda}}(f))_{\lambda \in \a^*_+} \big) \ .
\]
Note that $ \Psi_\sigma $ is a $G$-intertwining operator. 
For $ f \in L^2(\Gamma \bs G, \phi)(\gamma) $ ($ \gamma \in \hat{K} $), define
\[
    \Psi(f) = \sum_{\sigma \in \hat{M}} \Psi_\sigma(f) \ .
\]
Then, the map $ \Psi $ is well-defined as only finitely many maps $ \Psi_\sigma $ are nonzero.
We call such a function a \textit{spectral cut-off function}.

Let $ \sigma \in \hat{M} $ and let $W$ be a closed $G$-invariant subspace of $ L^2(\Gamma \bs G, \phi)_{ac} $.

Let us denote the orthogonal projection of $ L^2(\Gamma \bs G, \phi) $ to $ W $ by $ p_W $. This is clearly a $G$-intertwining operator.

\begin{lem}\label{Lem spectral cut-off maps W(gamma) to W(gamma)}
    If $ f \in W(\gamma) \ (\gamma \in \hat{K}) $, then $ \Psi(f) \in W(\gamma) $.
\end{lem}

\begin{proof}
    Since $ p_W $ is decomposable by Proposition \ref{Prop G-intertwining operators are decomposable}, there is a measurable family of maps $ L_\pi \colon N_\pi \to N_\pi $ ($ \pi \in \hat{G} $) such that
    \[
        p_W(f) = \F^{-1}\big( (L_\pi \otimes \Id)(\F_\pi(f))_{\pi \in \hat{G}} \big)
    \]
    for every $ f \in L^2(\Gamma \bs G, \phi) $. For $ f \in W $, we have
    \begin{align*}
         \Psi_\sigma(f)
         &= \Psi_\sigma(p_W(f))
         = \F_\sigma^{-1}\big( (\chi_\sigma(\lambda) \F_{\pi^{\sigma, i\lambda}}(p_W(f)))_{\lambda \in \a^*_+} \big) \\
         &= \F_\sigma^{-1}\big( (\chi_\sigma(\lambda) \F_{\pi^{\sigma, i\lambda}}(\F_\sigma^{-1}\big( (L_{\pi^{\sigma, i\mu}} \otimes \Id)(\F_{\pi^{\sigma, i\mu}}(f))_{\mu \in \a^*_+} \big)))_{\lambda \in \a^*_+} \big) \\
         &= \F_\sigma^{-1}\big( (L_{\pi^{\sigma, i\lambda}} \otimes \Id)(\chi_\sigma(\lambda) \F_{\pi^{\sigma, i\lambda}}(f))_{\lambda \in \a^*_+} \big) \\
         &= \F_\sigma^{-1}\big( (L_{\pi^{\sigma, i\lambda}} \otimes \Id)(\F_{\pi^{\sigma, i\lambda}}( \F_\sigma^{-1}\big( (\chi_\sigma(\mu) \F_{\pi^{\sigma, i\mu}}(f))_{\mu \in \a^*_+} \big) ))_{\lambda \in \a^*_+} \big) \\
         &= p_W(\Psi_\sigma(f)) \in W \ .
    \end{align*}
    The lemma follows.
\end{proof}

\begin{thm}\label{Thm L^2(Gamma|G, phi)_ac is orthogonal to the space of cusp forms}
    $ L^2(\Gamma \bs G, \phi)_{ac} $ is orthogonal to the space of cusp forms.
\end{thm}

\begin{rem}
    If $ \delta_\Gamma < 0 $, then the orthogonality of $ L^2(\Gamma \bs G, \phi)_{ac} $ to $ \c{\CS}(\Gamma \bs G, \phi) $ follows directly from Proposition \ref{Prop with widehat f^Omega}.
\end{rem}

\begin{proof}
    Assume that $ L^2(\Gamma \bs G, \phi)_{ac} $ is not orthogonal to the space of cusp forms. Then, there is $ f \in \c{\CS}(\Gamma \bs G, \phi)(\tau) $ ($ \tau \in \hat{K} $) and $ g \in L^2(\Gamma \bs G, \phi)_{ac} $ such that
    \[
        (p_{ac}(f), g) = (f, g) \neq 0 \ .
    \]
    In particular, $ p_{ac}(f) $ is nonzero. Let $ q $ be the orthogonal projection on $ \overline{\c{\CS}(\Gamma \bs G, \phi)} $.
    Since $ p_{ac} \circ q = q \circ p_{ac} $ by Proposition \ref{Prop a G-int op commutes with p_ds, p_ac and p_res + p_U}, $ p_{ac}(f) \in \overline{\c{\CS}(\Gamma \bs G, \phi)}(\tau) $.

    Let $ \Psi $ be a spectral cut-off function. Then, $ \Psi(p_{ac}(f)) $ belongs to $ \overline{\c{\CS}(\Gamma \bs G, \phi)}(\tau) $ by Lemma \ref{Lem spectral cut-off maps W(gamma) to W(gamma)} and it is a finite sum of compactly supported wave packets by construction and by Corollary \ref{Cor wave packet varphi_imu formula for p_ac(f)} and hence a Schwartz function by Proposition \ref{Prop Image of wave packet transform in Schwartz space}.
    It follows now from Theorem \ref{Thm a Schwartz fct in the closure of the space of cusp forms is a cusp form} that $ \Psi(p_{ac}(f)) $ is a cusp form.
    Thus, $ \Psi(p_{ac}(f)) = 0 $ by Proposition \ref{Prop no non nonzero compactly supported wave packet is cuspidal}. As we can choose $ \Psi $ arbitrarily, it follows that $ p_{ac}(f) = 0 $. This is a contradiction! The theorem follows.
\end{proof}

\begin{prop}\label{Prop L^2(Gamma|G, phi)_res is orthogonal to the space of cusp forms}
    $ L^2(\Gamma \bs G, \phi)_{res} $ is orthogonal to the space of cusp forms.
\end{prop}

\begin{proof}
    Let $p$ be the orthogonal projection on $ L^2(\Gamma \bs G, \phi)_{res} \oplus L^2(\Gamma \bs G, \phi)_U $ and let $q$ be the orthogonal projection on $ \overline{\c{\CS}(\Gamma \bs G, \phi)} $. Since $ p \circ q = q \circ p $ by Proposition \ref{Prop a G-int op commutes with p_ds, p_ac and p_res + p_U},
    \[
        p(\overline{\c{\CS}(\Gamma \bs G, \phi)}) = (L^2(\Gamma \bs G, \phi)_{res} \oplus L^2(\Gamma \bs G, \phi)_U) \cap \overline{\c{\CS}(\Gamma \bs G, \phi)} \ .
    \]
    As $ L^2(\Gamma \bs G, \phi)_U $ is contained in $ \overline{\c{\CS}(\Gamma \bs G, \phi)} $ by Theorem \ref{Thm c_(f, v) in 0CS(Gamma|G, phi) iff f in U_Lambda(sigma_lambda, phi)}, this is again equal to
    \[
        (L^2(\Gamma \bs G, \phi)_{res} \cap \overline{\c{\CS}(\Gamma \bs G, \phi)}) \oplus L^2(\Gamma \bs G, \phi)_U \ .
    \]
    Let $ f \in L^2(\Gamma \bs G, \phi)_{res} $ and $ g \in \overline{\c{\CS}(\Gamma \bs G, \phi)} $. Then,
    \[
        (f, g)_{L^2(\Gamma \bs G, \phi)} = (f, p(g))_{L^2(\Gamma \bs G, \phi)}
    \]
    As moreover $ L^2(\Gamma \bs G, \phi)_U $ is orthogonal to $ L^2(\Gamma \bs G, \phi)_{res} $, it suffices to show that
    \[
        L^2(\Gamma \bs G, \phi)_{res} \cap \overline{\c{\CS}(\Gamma \bs G, \phi)} = \{0\} \ .
    \]
    Since the orthogonal projection on $ (L^2(\Gamma \bs G, \phi)_{res} \oplus L^2(\Gamma \bs G, \phi)_U) \cap \overline{\c{\CS}(\Gamma \bs G, \phi)} $ is decomposable by Proposition \ref{Prop G-intertwining operators are decomposable}, there are subspaces $ S_\pi $ ($ \pi \in \hat{G}_c $) of $ N_\pi $ such that
    \[
        L^2(\Gamma \bs G, \phi)_{res} \cap \overline{\c{\CS}(\Gamma \bs G, \phi)} = \bigoplus_{\pi \in \hat{G}_c} c_\pi(S_\pi \hat{\otimes} V_\pi)
    \]
    Since in particular $ c_\pi(S_\pi \hat{\otimes} V_\pi) \subset L^2(\Gamma \bs G, \phi)_{res} $, every $ S_\pi $ must be contained in some $ E_\Lambda(\sigma_\lambda, \phi) $.
    Assume now that some space $ S_\pi $ is nontrivial. Then, it contained also a nonzero $K$-finite matrix coefficient $ c_\pi(f \otimes v) $.
    Thus, $ c_\pi(f \otimes v) \in \overline{\c{\CS}(\Gamma \bs G, \phi)} \cap \CS(\Gamma \bs G, \phi) $ by Lemma \ref{Lem c_(f, v) in the Schwartz space}.
    Hence, $ c_\pi(f \otimes v) \in \c{\CS}(\Gamma \bs G, \phi) $ by Theorem \ref{Thm a Schwartz fct in the closure of the space of cusp forms is a cusp form}.
    This is a contradiction by Theorem \ref{Thm c_(f, v) in 0CS(Gamma|G, phi) iff f in U_Lambda(sigma_lambda, phi)}.
\end{proof}

\begin{thm}\label{Thm decomposition of closure of cusp forms}
    We have:
    \begin{equation}\label{Thm decomposition of closure of cusp forms eq}
        \overline{\c{\CS}(\Gamma \bs G, \phi)} = L^2(\Gamma \bs G, \phi)_{ds} \oplus L^2(\Gamma \bs G, \phi)_U \ .
    \end{equation}
\end{thm}

\begin{rem}
    The decomposition \eqref{Thm decomposition of closure of cusp forms eq} was already conjectured by M.~Olbrich in 2002 (see \cite[p.116]{Olb02}).
\end{rem}

\begin{proof}
    As the spaces $ L^2(\Gamma \bs G, \phi)_U $ and $ L^2(\Gamma \bs G, \phi)_{ds} $ are contained in $ \overline{\c{\CS}(\Gamma \bs G, \phi)} $ by Proposition \ref{Prop U_Lambda(sigma_lambda, phi) consists of cusp forms} and Corollary \ref{Cor decomp. of space of cuspidal vectors if pi is a DS rep}, we have
    \begin{multline*}
        \overline{\c{\CS}(\Gamma \bs G, \phi)} = (L^2(\Gamma \bs G, \phi)_{ac} \oplus L^2(\Gamma \bs G, \phi)_{res}) \cap \overline{\c{\CS}(\Gamma \bs G, \phi)} \\
        \oplus L^2(\Gamma \bs G, \phi)_{ds} \oplus L^2(\Gamma \bs G, \phi)_U \ .
    \end{multline*}
    The theorem follows now from Theorem \ref{Thm L^2(Gamma|G, phi)_ac is orthogonal to the space of cusp forms} and from Proposition \ref{Prop L^2(Gamma|G, phi)_res is orthogonal to the space of cusp forms}.
\end{proof}

\newpage

\addtocontents{toc}{\protect\newpage}

\section{The geometrically finite case}

\subsection{Decompositions for varying parabolic subgroups}\label{ssec:DecompositionsForVaryingParabolics}

For notations, see Section \ref{ssec:GeometricPreparations}.

Since we may have several cusps, it is in general no longer sufficient to work only with one fixed parabolic subgroup of $G$. We investigate in the following what happens when we pass from one parabolic subgroup of $G$ to another.

For a Cartan involution $ \theta' $ of $\g$, $ h \in G $ and $ X \in \g $, set $ \LOI{\theta'}{h}(X) = \Ad(h) \theta'(\Ad(h)^{-1} X) $.
\\ For a Cartan involution $ \theta' $ of $G$ and $ g, h \in G $, set $ \LOI{\theta'}{h}(g) = h \theta'(h^{-1} g h) h^{-1} $.

Note that $ \LOI{\theta'}{h} $ is again a Cartan involution if $ \theta' $ is a Cartan involution.

We denote the maximal compact subgroup of $G$ associated to a Cartan involution $ \theta' $ by $ K_{\theta'} $.
Let $ \theta' $ be a Cartan involution of $G$. Then,
\[
    \{ g \in G \mid \LOI{\theta'}{h}(g) = g \} = \{ g \in G \mid h\theta'(h^{-1}gh)h^{-1} = g \} = h K_{\theta'} h^{-1} \qquad (h \in G) \ .
\]
Moreover, there exists $ h \in G $ such that $ \theta' = \LOI{\theta}{h} $ by, e.g., Lemma 2.3.2 of \cite[p.57]{Wallach}.
Thus, the maximal compact subgroups of $G$ are conjugate.

For $ X, Y \in \g $, set $ \langle X, Y \rangle_{\theta'} = \langle \Ad(h)^{-1} X, \Ad(h)^{-1} Y \rangle $. Note that this definition is independent of the choice of $h$ and that $ \langle \cdot, \cdot \rangle_{\theta'} $ is an $ \Ad(K_{\theta'}) $-invariant inner product on $ \g $.

By abuse of notation, we denote the norm $ \g \to [0, \infty) , \, X \mapsto \sqrt{\langle X, X \rangle_{\theta'}} $ also by $ |\cdot| $.
It will be clear from the context which norm is meant.

\begin{dfn} \nlenum
    \begin{enumerate}
    \item A subgroup $ Q $ of $G$ is called a \textit{parabolic subgroup} of $G$ if there exists $ g \in G $ such that $ Q = g P g^{-1} $.\index{parabolic subgroup}
    \item Let $Q$ be a parabolic subgroup of $G$. We say that $ Q = M_Q A_Q N_Q $ is a \textit{Langlands decomposition} if there exists $ g \in G $ such that $ M_Q = g M g^{-1} $, $ A_Q = g A g^{-1} $ and $ N_Q = g N g^{-1} $.\index{Langlands decomposition}
    \end{enumerate}
\end{dfn}

\begin{rem}
    This definition does not depend on the choice of $P$. It follows from the uniqueness property of the Cartan involutions, we have seen above, the fact that all the maximal abelian Lie subalgebras in $ \s $ are conjugate under the adjoint action of $K$ by Lemma 2.1.9 of \cite[p.47]{Wallach} and the construction of $P$ that the above definition does not depend on the choice of $P$.
\end{rem}

Let $Q$ be a parabolic subgroup of $G$ with Langlands decomposition $ M_Q A_Q N_Q $. Let $ \m_Q $ (resp. $ \a_Q $, $ \n_Q $) be the Lie algebra of $ M_Q $ (resp. $A_Q$, $ N_Q $).
Since the Killing form $B$ is invariant under the adjoint action of $G$, it follows from the previous definition and Section \ref{ssec:GeometricPreparations} that $ \m_Q \oplus \a_Q \oplus \n_Q $ is an orthogonal decomposition relative to $B$.

\begin{lem}\label{Lem g K g^(-1) cap K}
    Let $ g \in G $. Then, $ g K g^{-1} \cap K $ is equal to $K$ if $ g \in K $ and equal to a conjugate of $M$ in $K$ (which is strictly contained in $K$) otherwise.

    Moreover, if $ K' $ is a maximal compact subgroup of $G$, then there exists a unique element $ hK \in X $ such that $ K' = h K h^{-1} $.
\end{lem}

\begin{proof}
    If $ g \in K $, then $ g K g^{-1} \cap K $ is clearly equal to $K$.
    Assume that there exists $ g \in G \smallsetminus K $ such that $ g K g^{-1} \cap K \neq \{e\} $. Let $ h $ be a nontrivial element of $ g K g^{-1} \cap K $.
    Then, $ a_g \neq e $ and $ h g K = g K $. Thus, $ h k_g a_g K = k_g a_g K $. Hence, $ h k_g M = k_g M $, by the uniqueness property of the Cartan decomposition. So, $ k_g^{-1} h k_g \in M \iff h \in k_g M k_g^{-1} $.
    The first assertion follows.

    Let $ h_1, h_2 \in G $ be such that $ K' = h_1 K h_1^{-1} = h_2 K h_2^{-1} $. Then, $ h_2^{-1} h_1 K h_1^{-1} h_2 = K $.
    It follows from the first assertion that $ h_2^{-1} h_1 $ belongs to $K$. Thus, $ h_1 K = h_2 K $.
    The lemma follows.
\end{proof}

\begin{lem}\label{Lem A_Q = k A k^{-1}}
    Let $ \theta' $ be a Cartan involution of $\g$ and let $ Q = M_Q A_Q N_Q $ be a parabolic subgroup of $G$.
    The following conditions are equivalent:
    \begin{enumerate}
    \item $ M_Q $ is contained in $ K_{\theta'} $;
    \item $ \a_Q $ is orthogonal to the Lie algebra of $ K_{\theta'} $.
    \end{enumerate}
    Assume now that one of two conditions is satisfied.
    Let $ k \in K $ be such that $ Q = kPk^{-1} $ ($ kM $ is unique). Then, there exists a unique $ h \in N_Q A_Q $ such that $ K_{\theta'} = h K h^{-1} $. Moreover, $ N_Q = k N k^{-1} $, $ A_Q = h k A k^{-1} h^{-1} $ and $ M_Q = h k M k^{-1} h^{-1} = Z_{K_{\theta'}}(A_Q) \subset K_{\theta'} $.
\end{lem}

\begin{rem} \nlenum
    \begin{enumerate}
    \item When we take the Cartan involution $ \LOI{\theta}{h} $ and $ \Lie(h k A k^{-1} h^{-1}) $ for $ \a $ in Section \ref{ssec:GeometricPreparations} (note that this is indeed possible), then the constructed parabolic subgroup of $G$ in Section \ref{ssec:GeometricPreparations} is equal to $Q$.
    \item Let $ g \in G $ and let $ M_Q = g M g^{-1} $, $ A_Q = g A g^{-1} $, $ N_Q = g N g^{-1} $ and $ Q = M_Q A_Q N_Q $. Let $ k, h $ be as in the statement of the lemma. Then, it follows from the lemma that $ gM = h k M $. In particular, $ g K g^{-1} = h K h^{-1} $.
    \end{enumerate}
\end{rem}

\begin{proof}
    Let $ \theta' $, $ k $ and $ Q = M_Q A_Q N_Q $ be as above. Then, $ N_Q = k N k^{-1} $.

    Let $ g \in G $ be such that $ A_Q = g A g^{-1} $ and $ M_Q = g M g^{-1} $.

    \begin{enumerate}
    \item Assume that $ M_Q $ is contained in $ K_{\theta'} $.

        As $ m a = a m $ for all $ a \in A $ and $ m \in M $, $ m a = a m $ for all $ a \in A_Q $ and $ m \in M_Q $. So, $ M_Q = Z_{K_{\theta'}}(A_Q) $.

        Let $ h \in G $ be such that $ K_{\theta'} = h K h^{-1} $.

        Since $ K_{\theta'} A_Q N_Q = (h k K k^{-1} h^{-1}) (h k A k^{-1} h^{-1}) (h k N k^{-1} h^{-1}) $ is an Iwasawa decomposition of $G$, we can choose $ h \in N_Q A_Q = A_Q N_Q $ such that $ K_{\theta'} = h K h^{-1} $. $h$ is unique by Lemma \ref{Lem g K g^(-1) cap K}.
        
        As the Killing form is invariant under the adjoint action of $G$, as $ \a $ is orthogonal to $ \k $ and as $M$ is equal to $ K \cap P $, $ \Lie(h k A k^{-1} h^{-1}) $ is orthogonal to the Lie algebra of $ K_{\theta'} = h k K k^{-1} h^{-1} $ and $ M_Q = h k M k^{-1} h^{-1} \subset K_{\theta'} \cap Q $.

        Thus, $ A_Q = h k A k^{-1} h^{-1} $. Hence, $ \Lie(A_Q) = \Lie(h k A k^{-1} h^{-1}) $ is also orthogonal to the Lie algebra of $ K_{\theta'} $.
    \item Assume that $ \a_Q $ is orthogonal to the Lie algebra of $ K_{\theta'} $. As the Killing form is invariant under the adjoint action of $G$, as $ \Lie(h k A k^{-1} h^{-1}) $ is also orthogonal to the Lie algebra of $ K_{\theta'} $ and as $ \a $ is the unique maximal abelian Lie subalgebra of $ \m \oplus \a \oplus \n $ whose Lie algebra is orthogonal to that of $K$, $ A_Q = h k A k^{-1} h^{-1} $. Hence, $ M_Q = h k M k^{-1} h^{-1} $.
    \end{enumerate}

    The lemma follows.
\end{proof}

Let us also introduce some more notation.

Let $ \theta' $ be a Cartan involution of $ \g $.
Let $ Q = \tilde{M}_Q \tilde{A}_Q N_Q $ be a parabolic subgroup of $G$. Let $ \tilde{K}_Q $ be a maximal compact subgroup of $G$ containing the compact group $ \tilde{M}_Q $ (other choices are $ a \tilde{K}_Q a^{-1} $ ($ e \neq a \in \tilde{A}_Q $)).

Since $ G = NAK $ is an Iwasawa decomposition, it follows from Lemma \ref{Lem A_Q = k A k^{-1}} that $ G = N_Q \tilde{A}_Q \tilde{K}_Q $ is an other Iwasawa decomposition.

Let $ h $ be the unique element in $ N_Q \tilde{A}_Q $ such that $ K_{\theta'} = h \tilde{K}_Q h^{-1} $.
Then, $ M_{\theta', Q} := h \tilde{M}_Q h^{-1} $ is contained in $ K_{\theta'} $. Let $ A_{\theta', Q} = h \tilde{A}_Q h^{-1} $.

When it is clear which Cartan involution is meant, then we often use the simplified notations $ M_Q $ and $ A_Q $. We do the same for the notations which we introduce below.

Let $ k \in K $ be such that $ Q = kPk^{-1} $ and let $ h' \in N_Q A_{\theta', Q} $ be the element provided by Lemma \ref{Lem A_Q = k A k^{-1}}.
Then, $ K_{\theta'} = h' K {h'}^{-1} $.

Let $ \bar{N}_{\theta', Q} = \theta'(N_Q) $. We denote the Lie algebra of $ M_{\theta', Q} $ (resp. $ A_{\theta', Q} $, $ N_Q $, $ \bar{N}_{\theta', Q} $) by $ \m_{\theta', Q} $ (resp. $ \a_{\theta', Q} $, $ \n_Q $, $ \bar{\n}_{\theta', Q} $).

Define $ \alpha_{\theta', Q} \in \a_{\theta', Q}^{*} $ (resp. $ \rho_{\theta', Q} \in \a_{\theta', Q}^{*} $) by
\[
    \alpha_{\theta', Q}(H) = \alpha(\Ad(h'k)^{-1}(H)) \quad (\text{resp. } \rho_{\theta', Q}(H) = \rho(\Ad(h'k)^{-1}(H))) \; (H \in \a_{\theta', Q}) \ .
\]
Put $ A_{\theta', Q, +} = \{ a \in A_{\theta', Q} \mid a^{\alpha_{\theta', Q}} > 1 \} $ and $ \bar{A}_{\theta', Q, +} = \{ a \in A_{\theta', Q} \mid a^{\alpha_{\theta', Q}} \geq 1 \} $.
\\ For $ s > 0 $, put $ A_{\theta', Q, \leq s} = \{ a \in A_{\theta', Q} \mid a^{\alpha_{\theta', Q}} \leq s\} $ and $ A_{\theta', Q, > s} = \{ a \in A_{\theta', Q} \mid a^{\alpha_{\theta', Q}} > s\} $.

If $ \mu \in \a_{\theta', Q}^* $, then we define $ H_{\theta', Q, \mu} \in \a_{\theta', Q} $ by
\[
    H_{\theta', Q, \mu} = (\Ad(h'k)^{-1}(H))_\mu    \qquad ( H \in \a_{\theta', Q} ) \ .
\]
For $ H \in \a $ and $ \mu \in \a_{\theta', Q}^* $, set $ (h'k)^{-1}.\mu(H) = \mu(\Ad(h'k)H) $. Then, $ (h'k)^{-1}.\mu \in \a^* $ for all $ \mu \in \a_{\theta', Q}^* $.

We define a symmetric bilinear form $ \langle \cdot, \cdot \rangle_{\theta', Q} $ on $ \a_{\theta', Q}^* $ by
\[
    \langle \mu, \tau \rangle_{\theta', Q} = \langle (h'k)^{-1}.\mu, (h'k)^{-1}.\tau \rangle
\]
for all $ \mu, \tau \in \a_{\theta', Q}^* $.
For $ \mu \in \a_{\theta', Q}^* $, set $ |\mu| := \sqrt{\langle \mu, \mu \rangle_{\theta', Q}} $.

Thus, $ H \in \a_{\theta', Q} \mapsto \alpha_{\theta', Q}(H) \in \RR $ identifies $ \a_{\theta', Q} $ with $ \RR $ isometrically and that $ \mu \in \a_{\theta', Q}^* \mapsto \langle \alpha, \mu \rangle_{\theta', Q} \in \RR $ identifies $ \a_{\theta', Q}^* $ with $ \RR $ isometrically. Moreover, $ \alpha_{\theta', Q} $ is identified with 1.

For $ g \in G $, we set
\[
    \kappa_{\theta', Q}(g) = h'k\kappa(k^{-1}{h'}^{-1}gh'k)k^{-1}{h'}^{-1} \ .
\]
Define $ a_{\theta', Q} $, $ a_{\theta', Q} $, $ n_{\theta', Q} $, $ \nu_{\theta', Q} $, $ h_{\theta', Q} $ and $ k_{\theta', Q} $, $ k_{\theta', Q, \cdot} $, $ a_{\theta', Q, \cdot} $, $ h_{\theta', Q, \cdot} $ analogously.

By abuse of notation, we will often write for $ a_{\theta', Q, g}^{\alpha_{\theta', Q}} $ simply $ a_g $.

Then,
\[
    g = \kappa_{\theta', Q}(g) a_{\theta', Q}(g) n_{\theta', Q}(g) = \nu_{\theta', Q}(g) h_{\theta', Q}(g) k_{\theta', Q}(g) \qquad (g \in G)
\]
(Iwasawa decomposition of $G$ with respect to $ N_Q A_{\theta', Q} K_{\theta'} $) and
\[
    g = k_{\theta', Q, g} a_{\theta', Q, g} h_{\theta', Q, g} \qquad (g \in G)
\]
(Cartan decomposition of $G$). These decompositions are independent of the choice of $k$.
Since
\[
    a_{k^{-1}{h'}^{-1} g h' k} = a_{{h'}^{-1} g h'}   \qquad (g \in G, k \in K) \ ,
\]
\begin{equation}\label{eq Relation to standard a_g}
    a_{\theta', Q, g}^{\alpha_{\theta', Q}} = a_{{h'}^{-1} g h'}^\alpha = a_{{h'}^{-1} g h'} \ .
\end{equation}
Note that
\begin{equation}\label{eq a_Q(g)^(alpha_Q) = a(g)^alpha}
    a_{\theta', Q}(g)^{\alpha_{\theta', Q}} = a(k^{-1}{h'}^{-1} g h'k)^{\alpha} \qquad (g \in G) \ .
\end{equation}
The Haar measure on $N$ (resp. $ \bar{N} $, $A$) induces a Haar measure on $ N_Q $ (resp. $ \bar{N}_{\theta', Q} $, $A_{\theta', Q}$) so that
\[
    \int_{N_Q} a_{\theta', Q}(\theta'(n))^{-2\rho_{\theta', Q}} \, dn = 1
    \qquad \big(\text{resp. } \int_{\bar{N}_{\theta', Q}} a_{\theta', Q}(\bar{n})^{-2\rho_{\theta', Q}} \, d\bar{n} = 1 \big) \ .
\]
Let $ \theta'' $ be an other Cartan involution.
Let $ h'' \in N_Q A_{\theta', Q} $ be the element provided by Lemma \ref{Lem A_Q = k A k^{-1}} for $ \theta' = \theta'' $.
Let $ \tilde{h} = h'' {h'}^{-1} \in N_Q A_{\theta', Q} $. Then, $ K_{\theta''} = h'' K {h''}^{-1} = \tilde{h} K_{\theta'} \tilde{h}^{-1} $.
Since
\begin{multline*}
    a_{\theta'', Q}(g)^{\alpha_{\theta'', Q}} = a(k^{-1}{h''}^{-1} g h''k)^{\alpha} \\
    = a(k^{-1} {h'}^{-1} (h' {h''}^{-1} g h'' {h'}^{-1}) h'k)^{\alpha} = a_{\theta', Q}(\tilde{h}^{-1} g \tilde{h})^{\alpha_{\theta', Q}}
\end{multline*}
by \eqref{eq a_Q(g)^(alpha_Q) = a(g)^alpha}, it follows that
\begin{equation}\label{eq Compare a(g)}
    a_{\theta'', Q}(g)^{\alpha_{\theta'', Q}} = a_{\theta', Q}(\tilde{h}^{-1} g \tilde{h})^{\alpha_{\theta', Q}} \ .
\end{equation}
Similarly, one shows, using \eqref{eq Relation to standard a_g}, that
\begin{equation}\label{eq Compare a_g}
    a_{\theta', Q, g}^{\alpha_{\theta'', Q}} = a_{\tilde{h}^{-1} g \tilde{h}}^{\alpha_{\theta', Q}} \ .
\end{equation}
Let us show now how we can compare spaces constructed from different Cartan involutions.

Let $ K' $ be a maximal compact subgroup of $G$.
Let $ h \in G $ be such that $ K' = h K h^{-1} $. Then, $ K' $ is the stabiliser of $ hK $. So,
\[
    X' := G/K' \to X = G/K , \quad g K' \mapsto g h K
\]
is a natural isomorphism between $ X' $ and $X$.
Let $ P' $ be a parabolic subgroup of $G$. 
Let $ h \in G $ be such that $ P' = h P h^{-1} $. Then, $ P' $ is the stabiliser of $ hP $. So,
\[
    \dX' := G/P' \to \dX = G/P , \quad g P' \mapsto g h P
\]
is a natural isomorphism between $ \dX' $ and $ \dX $.

If $ h \in G $ is such that $ K' = h K h^{-1} $ and $ P' = h P h^{-1} $, then
\[
    \bar{X} = X \cup \dX = G/K \cup G/P
\]
is canonically isomorphic to
\[
    \bar{X}' := X' \cup \dX' = G/K' \cup G/P' \ .
\]
To shorten notations, we simply write $ g K' = g h K $ and $ g P' = g h P $.

Let us come now to the impact on functions of the change of the Cartan involution.

If $ f \colon G \to \CC $ is a right $K$-invariant function, then $ R_h f $ is a right $K'$-invariant function on $ G $ as
\[
    R_h f(g k) = f((g h) (h^{-1} k h)) = f(g h) = R_h f(g)
\]
for all $ g \in G $ and all $ h \in K' $.

\newpage

\subsection{Geometrically finite groups and \texorpdfstring{$\Gamma$}{Gamma}-cuspidal parabolic subgroups}\label{ssec:Geometrically finite groups}

We define geometrically finite groups and $\Gamma$-cuspidal parabolic subgroups of $G$, which are the representatives of the cusps of the space $ \Gamma \bs X = \Gamma \bs G / K $.

\medskip


Let $ \Gamma $ be a torsion-free discrete subgroup of $G$. Recall that $ \Lambda_\Gamma $ denotes the limit set (closed and $ \Gamma $-invariant set) and $ \Omega_\Gamma $ the set of ordinary points (open and $ \Gamma $-invariant set). Moreover, $ \dX = \Omega_\Gamma \cup \Lambda_\Gamma $ and $ \Gamma $ acts properly discontinuously on $ X \cup \Omega_\Gamma $.
As moreover $ \Gamma $ is torsion-free, $ \Gamma $ acts freely on $ X \cup \Omega_\Gamma $. See Section \ref{ssec:Dfn convex-cocompact} for further details.

\begin{dfn}\label{Dfn Gamma-cuspidal}\index{Gamma-cuspidal@$\Gamma$-cuspidal parabolic subgroup}
    We say that a parabolic subgroup $ Q = M_Q A_Q N_Q $ of $ G $ is $ \mathit{\Gamma} $\textit{-cuspidal} if $ \Gamma_Q := \Gamma \cap Q $ is an infinite subgroup of $G$ which is contained in $ N_Q M_Q $.\index[n]{GammaQ@$ \Gamma_Q $}
\end{dfn}

\begin{rem} \nlenum
    \begin{enumerate}
    \item This definition does not depend on the chosen Langlands decomposition as $ N_Q M_Q = N_Q n M_Q n^{-1} $ for all $ n \in N_Q $.
    \item If $ Q $ is $ \Gamma $-cuspidal, then $ \gamma Q \gamma^{-1} $ is $ \Gamma $-cuspidal for all $ \gamma \in \Gamma $. 
    \end{enumerate}
\end{rem}

Let $ Q_1, Q_2 $ be two $\Gamma$-cuspidal parabolic subgroups of $ G $. Then, we say that $ Q_1 $ is $ \Gamma $-equivalent to $ Q_2 $ if there exists $ \gamma \in \Gamma $ such that $ Q_2 = \gamma Q_1 \gamma^{-1} $. We write $ Q_1 \sim_\Gamma Q_2 $ in this case.

Let $ \P_\Gamma $ be the set of $\Gamma$-conjugacy classes of $ \Gamma $-cuspidal parabolic subgroups of $G$.\index[n]{PGamma@$ \P_\Gamma $}

\begin{dfn}
    We say that $ u \in \dX $ is a \textit{cusp} for $ \Gamma $ if its stabiliser is a $ \Gamma $-cuspidal parabolic subgroup of $G$.
\end{dfn}

Let $ [Q]_\Gamma \in \P_\Gamma $. Then, there exists $ k \in K $ such that $ Q = kPk^{-1} $.
Let $ \tilde{M}_Q = k M k^{-1} $ and let $ \tilde{A}_Q = k A k^{-1} $.

For a $ \Gamma $-cuspidal parabolic subgroup $ Q = \tilde{M}_Q \tilde{A}_Q N_Q $ of $ G $, we consider the exact sequence
\begin{equation}
    0   \rightarrow     N_Q   \rightarrow     N_Q \tilde{M}_Q   \stackrel{l_{\tilde{M}_Q}}{\rightarrow}   \tilde{M}_Q    \rightarrow    0 \ ,
\end{equation}
where $ l_{\tilde{M}_Q} \colon N_Q \tilde{M}_Q \to \tilde{M}_Q $ is the projection onto $\tilde{M}_Q$. Put $ \tilde{M}_{\Gamma, Q} = \overline{l_{\tilde{M}_Q}(\Gamma_Q)} \subset \tilde{M}_Q $.

\needspace{2\baselineskip}\index[n]{N_GammaQ M_Gamma@$ N_{\Gamma, Q} $, $ M_{\Gamma, Q} $}
By Auslander's theorem (Theorem 3.1 of \cite[p.25]{BO07}), we have the following:
\begin{enumerate}
\item There exist a connected subgroup $ \tilde{N}_{\Gamma, Q} \subset N_Q $ and $ n \in N_Q $ such that $ n \Gamma_Q n^{-1} $ normalises $ \tilde{N}_{\Gamma, Q} $ and such that $ n \Gamma_Q n^{-1} $ is a cocompact lattice in the unimodular group $ \tilde{Q}_\Gamma := \tilde{N}_{\Gamma, Q} \tilde{M}_{\Gamma, Q} $.
\item Fix $ n \in N_Q $ as above and $ a \in \tilde{A}_Q $. 
    Let $ N_{\Gamma, Q} = n^{-1} \tilde{N}_{\Gamma, Q} n \subset N_Q $, $ M_{\Gamma, Q} = n^{-1} \tilde{M}_{\Gamma, Q} n $ and $ Q_\Gamma = N_{\Gamma, Q} M_{\Gamma, Q} $.
    \\ Let us introduce the following notations according to the appearing case: 
    \begin{enumerate}
    \item $ \tilde{M}_{\Gamma, Q} \neq \{e\} $: Let $ \theta_Q = n^{-1} a^{-1} \theta a n $, $ K_Q = n^{-1} a^{-1} K a n $, $ A_Q = n^{-1} \tilde{A}_Q n $ and $ M_Q = n^{-1} \tilde{M}_Q n $. Then, $ M_{\Gamma, Q} \subset M_Q \subset K_Q $.
    \item $ \tilde{M}_{\Gamma, Q} = \{e\} $ ($ \Gamma_Q \subset N_Q $): We can define the groups as previously (we want to allow a maximum of freedom) or we use the following definitions:
        Let $ \theta_Q = \theta $, $ K_Q = K $, $ A_Q = \tilde{A}_Q $, $ M_Q = \tilde{M}_Q $. Then, $ M_{\Gamma, Q} = \{e\} \subset M_Q \subset K_Q $ and $ Q_\Gamma = N_{\Gamma, Q} $.
    \end{enumerate}
    For $ e \neq \gamma \in \Gamma $, we set $ A_{\gamma Q \gamma^{-1}} = \gamma A_Q \gamma^{-1} $. We define $ N_{\gamma Q \gamma^{-1}} $, $ N_{\Gamma, \gamma Q \gamma^{-1}} $, $ K_{\gamma Q \gamma^{-1}} $, $ A_{\gamma Q \gamma^{-1}} $, $ M_{\gamma Q \gamma^{-1}} $, $ M_{\Gamma, \gamma Q \gamma^{-1}} $ and $ (\gamma Q \gamma^{-1})_\Gamma $ analogously.
    Then,
    \[
        \theta_{\gamma Q \gamma^{-1}} = \LOI{\theta}{\gamma} \ .
    \]
    These definitions fix uniquely the above $ n $ and $a$ for $ \gamma Q \gamma^{-1} $.

    Moreover, there exists a normal subgroup of finite index $ \Gamma' $ of $ \Gamma_Q $ such that $ M_{\Gamma', Q} := \overline{l_{M_Q}(\Gamma')} \subset M_Q $ is a torus (commutative, compact, connected Lie group).
\item We can write any element $ \gamma \in \Gamma' $ (resp. $ \gamma \in \Gamma_Q $) as $ n_\gamma m_\gamma $ with $ n_\gamma \in N_{\Gamma, Q} $ and $ m_\gamma \in M_{\Gamma, Q} $ such that $ m_\gamma $ centralises (resp. normalises) $ N_{\Gamma, Q} $.
\item The group $ \{ n_\gamma \mid \gamma \in \Gamma' \} $ is a lattice in $ N_{\Gamma, Q} $.
\end{enumerate}
As Definition \ref{Dfn Gamma-cuspidal} does not depend on the chosen Langlands decomposition, $ \Gamma_Q $ is also an infinite subgroup of $G$ which is contained in $ N_Q M_Q $.

When we write in the following $ Q = M_Q A_Q N_Q $, then we mean with this that we consider the Langlands decomposition for which $ M_Q $ and $ A_Q $ are given as above.

Let $ w_Q \in N_{K_Q}(A_Q) $ denote a representative of the nontrivial element of the Weyl group of $ (\g, \a_Q) $. 

By abuse of notation, we will write for $ a_{Q, g}^{\alpha_Q} $ simply $ a_g $.

Let $ \varphi \colon \{ n_\gamma \mid \gamma \in \Gamma' \} \to M_{\Gamma', Q} , \, n_\gamma \mapsto m_\gamma $. This defines a group homomorphism from $ \{ n_\gamma \mid \gamma \in \Gamma' \} $ to $ M_{\Gamma', Q} $. Indeed:
\begin{enumerate}
\item This map is well-defined: Let $ \gamma_1, \gamma_2 \in \Gamma' $ be such that $ n_{\gamma_1} = n_{\gamma_2} $. Then,
    \[
        \gamma_1 \gamma_2^{-1} = n_{\gamma_1} m_{\gamma_1} m_{\gamma_2}^{-1} n_{\gamma_1}^{-1}
        = m_{\gamma_1} m_{\gamma_2}^{-1} \in \Gamma \cap K_Q
    \]
    since $ m_{\gamma_1} $ and $ m_{\gamma_2} $ centralise $ N_{\Gamma, Q} $. As $ \Gamma $ is torsion-free, it follows that $ \gamma_1 \gamma_2^{-1} = e $. Thus, $ \gamma_1 = \gamma_2 $. Hence, $ m_{\gamma_1} = m_{\gamma_2} $.
\item It is a group homomorphism: The neutral element is clearly mapped to the neutral element.
    Let $ \gamma_1, \gamma_2 \in \Gamma' $. Then,
    \[
        \gamma_1 \gamma_2 = n_{\gamma_1} m_{\gamma_1} n_{\gamma_2} m_{\gamma_2}
        = (n_{\gamma_1} n_{\gamma_2}) (m_{\gamma_1} m_{\gamma_2})
    \]
    as $ m_{\gamma_1} $ centralises $ N_{\Gamma, Q} $. Thus, $ n_{\gamma_1} n_{\gamma_2} = n_{\gamma_1 \gamma_2} $ is mapped to $ m_{\gamma_1 \gamma_2} = m_{\gamma_1} m_{\gamma_2} $.
\end{enumerate}

By Lemma 3.2 of \cite[p.26]{BO07}, there exists a subgroup $V$ of finite index in $ \{ n_\gamma \mid \gamma \in \Gamma' \} $ such that the restriction of $ \varphi $ to $V$ extends to a group homomorphism from $ N_{\Gamma, Q} $ to $ M_{\Gamma, Q} $, which we denote by abuse of notation also by $ \varphi $.

By 3.1.5 of \cite[p.26]{BO07}, we can choose $ \Gamma' $ with the following properties:
\begin{enumerate}
\item $ \{ n_\gamma \mid \gamma \in \Gamma' \} $ is contained in $V$;
\item $ M_{\Gamma', Q} $ coincides with $ \overline{\varphi(N_{\Gamma, Q})} $;
\item $ M_{\Gamma', Q} $ has finite index in $ M_{\Gamma, Q} $.
\end{enumerate}

Let $ Y_{\Gamma} $ (resp. $ Y_{\Gamma_Q} $) denote the manifold without boundary $ \Gamma \bs X $ (resp. $ \Gamma_Q \bs X $). 

Let $ g \in G $ be such that $ Q = g P g^{-1} $.
Let $ \Lambda_Q = \{eQ\} = \{gP\} \subset \dX $ (limit set of $ \Gamma_Q $) and let $ \Omega_Q = \dX \smallsetminus \Lambda_Q \subset \dX $ (ordinary set of $ \Gamma_Q $).

Let $ \bar{Y}_{\Gamma} $ (resp. $ \bar{Y}_{\Gamma_Q} $, $ \bar{Y}_{\Omega, \Gamma_Q} $) denote the manifold with boundary $ \Gamma \bs (X \cup \Omega_\Gamma) $ (resp. $ \Gamma_Q \bs (X \cup \Omega_Q) $, $ \Gamma_Q \bs (X \cup \Omega_\Gamma) $). 

\begin{dfn}\label{Dfn geometrically finite}
    We say that a torsion-free discrete subgroup $ \Gamma $ of $G$ is \textit{geometrically finite} if the following conditions hold:\index{geometrically finite}
    \begin{enumerate}
    \item The set $ \P_\Gamma $ of $\Gamma$-conjugacy classes of $ \Gamma $-cuspidal parabolic subgroups of $G$ is finite. We denote the elements of $ \P_{\Gamma} $ by $ [P_1]_\Gamma, \dotsc, [P_m]_\Gamma $.
    \item There exists a compact subset $ C $ of $ \bar{Y}_\Gamma $ such that the set of connected components of $ \bar{Y}_\Gamma \smallsetminus C $ is in bijection with $ \P_{\Gamma} $.
    \item For each connected component $ Y_i \subset \bar{Y}_\Gamma \smallsetminus C $ corresponding to $ [P_i]_\Gamma $ ($ i = 1, \dotsc, m $), there exists a subset $ \tilde{Y}_i $ of $ \bar{Y}_{\Omega, \Gamma_{P_i}} $ such that
        \begin{enumerate}
        \item $ \bar{Y}_{\Gamma_{P_i}} \smallsetminus \tilde{Y}_i $ is compact in $ \bar{Y}_{\Gamma_{P_i}} $,
        \item the map $ \pi_i \colon \bar{Y}_{\Omega, \Gamma_{P_i}} \to \bar{Y}_\Gamma , \, \Gamma_{P_i} x \mapsto \Gamma x $ is injective on $ \tilde{Y}_i $ and $ Y_i = \pi_i(\tilde{Y}_i) $.
        \end{enumerate}
    \end{enumerate}
\end{dfn}

\begin{rem} \nlenum
    \begin{enumerate}
    \item If $ \Gamma \subset G $ is geometrically finite, then the spaces $ Y_\Gamma $, $ Y_{\Gamma_{P_i}} $ are Riemannian manifolds in a natural way and the maps $ \pi_i $ restricted to $ \tilde{Y}_i \cap Y_{\Gamma_{P_i}} $ are isometric embeddings.
    \item If $ \Gamma \subset G $ is geometrically finite and if $ \P_{\Gamma} = \emptyset $, then $ \Gamma $ is convex-cocompact.
    \item Our definition is equivalent to Bowditch's definition \cite{Bowditch} (he gives several definitions and proves that they are all equivalent).
    \end{enumerate}
\end{rem}

From now on we will assume $ \Gamma \subset G $ is geometrically finite unless stated otherwise.

Let $ e_i \colon Y_i \to \tilde{Y}_i \subset \bar{Y}_{\Gamma_{P_i}} $ be the inverse map of $ \restricted{\pi_i}{\tilde{Y}_i} $.
As $ \pi_i $ is an open map, $ e_i $ is continuous.

Let $Q$ be a $ \Gamma $-cuspidal parabolic subgroup of $G$.

Let $ w_Q \in N_{K_Q}(\a_Q) $ be a representative of the nontrivial Weyl group element.
Let $ w_i = w_{P_i} $.

\begin{lem}\label{Lem Gamma_(P_i)|Lambda_i compact}
    $ \Gamma_Q \bs (\Lambda_\Gamma \smallsetminus \{ e Q \}) $ is a compact subset of $ \Gamma_Q \bs \Omega_Q $.
\end{lem}

\begin{proof}
    Since $ \Lambda_\Gamma $ is closed in $ \dX $, $ \Lambda_i = \Lambda_\Gamma \cap \Omega_{P_i} $ is closed in $ \Omega_{P_i} $. Thus, $ \Gamma_{P_i} \bs \Lambda_i $ is closed in $ \Gamma_{P_i} \bs \Omega_{P_i} $ by definition of the quotient topology.

    As $ e_i(Y_i) $ is contained in $ \bar{Y}_{\Omega, \Gamma_{P_i}} $, $ e_i(Y_i) \cap \Gamma_{P_i} \bs \Lambda_i = \emptyset $. Moreover, it follows from Definition \ref{Dfn geometrically finite} that $ C_i := \bar{Y}_{\Gamma_{P_i}} \smallsetminus e_i(Y_i) $ is a compact subset of $ \bar{Y}_{\Gamma_{P_i}} $.
    Thus, $ \Gamma_{P_i} \bs \Lambda_i $ is contained in $ C_i \cap (\Gamma_{P_i} \bs \Omega_{P_i}) $, which is a compact subset of $ \Gamma_{P_i} \bs \Omega_{P_i} $. 
    The lemma follows for $ Q = P_i $. Hence, it holds also for $ Q = \gamma P_i \gamma^{-1} $ ($ \gamma \in \Gamma $).
\end{proof}

\begin{lem}\label{Lem Gamma_Q is geometrically finite}
    $ \Gamma_Q $ is a geometrically finite subgroup of $G$. Moreover, $ \P_{\Gamma_Q} = \{ Q \} $.
\end{lem}

\begin{proof}
    Let $ Q' $ be another parabolic subgroup of $G$ (different from $Q$). As all parabolic subgroups of $G$ are conjugate, there exists $ g \in G $ such that $ Q' = g \bar{Q} g^{-1} $. Since $ \bar{Q} $ is the stabiliser of $ w_Q Q $, $ Q' $ is the stabiliser of $ g w_Q Q $. 

    It follows from the Bruhat decomposition ($ G = N_Q w_Q Q \cup Q $), that there exists also $ n \in N_Q $ such that $ Q' $ is the stabiliser of $ n w_Q Q $.
    Thus, $ Q' = n \bar{Q} n^{-1} $.
    Hence,
    \begin{multline*}
        Q \cap Q' = M_Q A_Q N_Q \cap (n M_Q n^{-1}) (n A_Q n^{-1}) (n \bar{N}_Q n^{-1}) \\
        = (n M_Q n^{-1}) (n A_Q n^{-1}) =: \tilde{M} \tilde{A} \ .
    \end{multline*}
    So, $ \Gamma_Q \cap Q' \subset \tilde{M} \tilde{A} \cap M_Q N_Q = \tilde{M} \tilde{A} \cap \tilde{M} N_Q = \tilde{M} $ as the Langlands decomposition of $ \tilde{M} \tilde{A} N_Q $ is unique. Since $ \tilde{M} $ is compact and $ \Gamma_Q $ is discrete, it follows that $ \Gamma_Q \cap Q' $ is finite (in fact it is trivial as $ \Gamma $ is torsion-free).
    Consequently, $ \P_{\Gamma_Q} = \{ [Q]_{\Gamma_Q} \} = \{ Q \} $.
    The lemma follows now easily from this.
\end{proof}



\newpage

\subsection{\texorpdfstring{$ N_\Gamma N^\Gamma $}{N\_Gamma N\string^Gamma} decomposition of the unipotent radical \texorpdfstring{$N$}{N} of a cuspidal parabolic subgroup of \texorpdfstring{$G$}{G}}

We show in this section how the unipotent radical $ N_Q $ of a $ \Gamma $-cuspidal parabolic subgroup $Q$ of $G$ can be decomposed.

\medskip

From now on, we assume that $ Q $ is a $ \Gamma $-cuspidal parabolic subgroup of $G$ unless stated otherwise.

Let $ \n_{\Gamma, Q} $ be the Lie algebra of $ N_{\Gamma, Q} $.
Since
\[
    0 \to \g_{2\alpha_Q} \cap \n_{\Gamma, Q} \mathrel{\mathop{\longrightarrow}_{\iota}} \n_{\Gamma, Q}
    \mathrel{\mathop{\longrightarrow}_{p}} p_{\g_{\alpha_Q}}(\n_{\Gamma, Q}) \to 0 \ ,
\]
where $ \iota $ is the canonical injection and $p$ is the canonical projection, is a short exact sequence,
\[
    \dim(\n_{\Gamma, Q}) = \dim(p_{\g_{\alpha_Q}}(\n_{\Gamma, Q})) + \dim(\g_{2\alpha_Q} \cap \n_{\Gamma, Q}) \ .
\]
Let $ \n^{\Gamma, Q}_1 $ (resp. $ \n^{\Gamma, Q}_2 $) be the orthogonal complement of $ p_{\g_{\alpha_Q}}(\n_{\Gamma, Q}) $ in $ \g_{\alpha_Q} $ (resp. $ \g_{2\alpha_Q} \cap \n_{\Gamma, Q} $ in $ \g_{2\alpha_Q} $) with respect to $ \langle \cdot, \cdot \rangle_{\theta_Q} $.
Set $ \n^{\Gamma, Q} = \n^{\Gamma, Q}_1 \oplus \n^{\Gamma, Q}_2 $. This space is invariant under the adjoint action of $A_Q$ as so are $ \n^{\Gamma, Q}_1 $ and $ \n^{\Gamma, Q}_2 $.
Assume that $ X \in \n_{\Gamma, Q} \cap \n^{\Gamma, Q} $. Write $ X = X_1 + X_2 $ with $ X_j \in \g_{j\alpha_Q} $. Then, $ X_1 \in p_{\g_{\alpha_Q}}(\n_{\Gamma, Q}) \cap \n^{\Gamma, Q}_1 = \{0\} $. Thus, $ X = X_2 \in \g_{2\alpha_Q} \cap \n_{\Gamma, Q} \cap \n^{\Gamma, Q}_2 = \{0\} $. Hence, $ X = 0 $. Since in addition $ \dim(\n_{\Gamma, Q}) + \dim(\n^{\Gamma, Q}) = \dim(\n_Q) $ by construction,
\begin{equation}
    \n_Q = \n_{\Gamma, Q} \oplus \n^{\Gamma, Q} \ .
\end{equation}
Set $ N^{\Gamma, Q} = \exp(\n^{\Gamma, Q}) $ (invariant under conjugation in $ A_Q $).\index[n]{N^GammaQ@$ N^{\Gamma, Q} $}

Let $ N^{\Gamma, Q}_1 = \exp(\n^{\Gamma, Q}_1) $ and $ N^{\Gamma, Q}_2 = \exp(\n^{\Gamma, Q}_2) $.

We put as measure on all measurable subsets $S$ of $ N_Q $ the push-forward of the Lebesgue measure on $ \log(S) $ (with an appropriate normalisation).
Here, $ \log $ is the inverse map of the exponential map $ \exp \colon \n_Q \to N_Q $.

Since $ \n_{\Gamma, Q} $ and $ \g_{2\alpha_Q} $ are invariant under the adjoint action of $ M_{\Gamma, Q} $ and since $ p_{\g_{\alpha_Q}} $ is clearly $ M_Q $-equivariant, $ \g_{2\alpha_Q} \cap \n_{\Gamma, Q} $ and $ p_{\g_{\alpha_Q}}(\n_{\Gamma, Q}) $ are also invariant under the adjoint action of $ M_{\Gamma, Q} $.

As $ \langle \cdot, \cdot \rangle_{\theta_Q} $ is invariant under the adjoint action of $K_Q$, it follows that $ \n^{\Gamma, Q}_1 $, $ \n^{\Gamma, Q}_2 $ and $ \n^{\Gamma, Q} $ are invariant under the adjoint action of $ M_{\Gamma, Q} $, too. Thus, $ N^{\Gamma, Q} $ is invariant under conjugation in $ M_{\Gamma, Q} $.

\begin{lem}\label{Lem using the Baker-Campbell-Hausdorff formula}
    For all $ X, Y \in \n_Q $, $ \exp(X) \exp(Y) = \exp(X+Y)\exp(\frac12 [X, Y]) $.
\end{lem}

\begin{proof}
    By the Baker-Campbell-Hausdorff formula, we have
    \[
        \log(\exp(X) \exp(Y)) = X + Y + \frac12 [X, Y] \qquad (X, Y \in \n_Q)
    \]
    and
    \[
        \log(\exp(X) \exp(Y)) = X + Y \qquad (X \in \n_Q, Y \in \g_{2\alpha_Q}) \ .
    \]
    Thus, $ \exp(X) \exp(Y) = \exp(X+Y) \exp(\frac12 [X, Y]) $ for all $ X, Y \in \n_Q $. The lemma follows.
\end{proof}

\begin{prop}\label{Prop Decomposition of N}
    We have: $ N_Q = N_{\Gamma, Q} N^{\Gamma, Q} $. Moreover, every $ n \in N_Q $ writes uniquely as $ n_1 n_2 $ with $ n_1 \in N_{\Gamma, Q} $ and $ n_2 \in N^{\Gamma, Q} $.
\end{prop}

\begin{rem}
    As $ N^{\Gamma, Q} $ may be no Lie group, the $ N_{\Gamma, Q} N^{\Gamma, Q} $ decomposition of $ n n' $ ($ n, n' \in N^{\Gamma, Q} $) can have a nontrivial $ N_{\Gamma, Q} $-part!
\end{rem}

\begin{proof}
    As $ N_{\Gamma, Q} $ (resp. $ N^{\Gamma, Q} $) is isomorphic to $ \n_{\Gamma, Q} $ (resp. $ \n^{\Gamma, Q} $), as $ \n_Q = \n_{\Gamma, Q} \oplus \n^{\Gamma, Q} $ and as $ \g_{2\alpha_Q} = (\n_{\Gamma, Q} \cap \g_{2\alpha_Q}) \oplus (\n^{\Gamma, Q} \cap \g_{2\alpha_Q}) $,
    \[
        N_{\Gamma, Q} \times N^{\Gamma, Q} \to N_Q , \; (n_1, n_2) \mapsto n_1 n_2
    \]
    is a diffeomorphism by Proposition 17 of Chapter III of \cite[p.237]{Bourbaki}. The proposition follows.
\end{proof}

Let $ n \in N_Q $. Write $ n = n_1 n_2 $ with $ n_1 \in N_{\Gamma, Q} $ and $ n_2 \in N^{\Gamma, Q} $.
\\ Set
\[
    p_{N_{\Gamma, Q}}(n) = n_1 \in N_{\Gamma, Q}
\]
and
\[
    p_{N^{\Gamma, Q}}(n) = n_2 \in N^{\Gamma, Q} \ .
\]
Then, $ p_{N_{\Gamma, Q}} \colon N_Q \to N_{\Gamma, Q} $ and $ p_{N^{\Gamma, Q}} \colon N_Q \to N^{\Gamma, Q} $ are analytic maps by the proof of the previous proposition.

\newpage

\subsection{Generalised Siegel sets}\label{ssec: Generalised Siegel sets}

In the following, we generalise the notions of a Siegel set.
Moreover, we use these sets to get a cover of the quotient space $ \Gamma \bs X $.

\begin{dfn}\index{at most $k$-to-one map}
    Let $ S, T $ be sets, $ k \in \NN $ and $ f $ a map from $ S $ to $ T $. Then, $f$ is said to be \textit{at most $k$-to-one} if, for all $ y \in T $, $ f^{-1}(y) $ contains at most $k$ elements.
\end{dfn}

\index[n]{SiegelSet@$\Sfrak$}
\begin{dfn}\label{Dfn Generalised Siegel set}\index{generalised Siegel set} \nlenum
    \begin{enumerate}
    \item 
        Let $ Q = M_Q A_Q N_Q $ be a $ \Gamma $-cuspidal parabolic subgroup of $G$. Recall that $ w_Q \in N_{K_Q}(A_Q) $ denotes a representative of the nontrivial element of the Weyl group of $ (\g, \a_Q) $. 
        Let $ V $ be a connected subset of $N_Q$, which is the closure of an open subset of $N_Q$, such that 
        \[
            \Phi \colon V M_{\Gamma, Q} \to \Gamma_Q \bs N_Q M_{\Gamma, Q} \ , \quad n m \mapsto \Gamma_Q nm
        \]
        is surjective and at most $k$-to-one for some $ k \in \NN $.
        Let $ V' $ be a subset of $ V \subset N_Q $, which is relatively compact in $ N_Q $, such that $ V \smallsetminus V' $ is closed and such that
        \[
            \{ n \in N_Q \mid n w_Q Q \in N_{\Gamma, Q} (\Lambda_\Gamma \smallsetminus \{eQ\}) \} M_{\Gamma, Q}
        \]
        is contained in $ \Gamma_Q V' M_{\Gamma, Q} $.
        \\ We say that a subset $ \Sfrak_{Q, t, V, V'} $ of $G$ is a \textit{generalised Siegel set} with respect to $ Q $ if it is of the form
        \[
            V A_Q K_Q \smallsetminus V' A_{Q, <t} K_Q \ ,
        \]
        where $ A_{Q, < t} := \{ a \in A_Q \mid a^{\alpha_Q} < t \} $ (here $ t > 0 $ is fixed).
    \item We denote often $ \Sfrak_{Q, t, V, V'} $ simply by $ \Sfrak_Q $ or $ \Sfrak $.
    \item We call $ \Sfrak' := \{ g K_Q \mid g \in \Sfrak \} $ a \textit{generalised Siegel set in $\mathit{G/K_Q}$}.
    \item Let $K'$ be a maximal compact subgroup of $G$ and let $ h \in G $ be such that $ K_Q = h K' h^{-1} $. Set $ X' = G/K' $. Then, we say that
        \[
            \Sfrak'_{X'} := \{g h K' \mid g K_Q \in \Sfrak' \}
        \]
        is a \textit{generalised Siegel set in $\mathit{X'}$}.
    \item We say that $ \Sfrak' $ is a generalised Siegel set at $ eQ $.
        Moreover, if $u$ is $ \Gamma $-cuspidal, then $ \Sfrak' $ is also called $ \Gamma $-cuspidal.
    \end{enumerate}
\end{dfn}

\begin{rem}\label{Rem of Dfn Generalised Siegel set} \nlenum
    \begin{enumerate}
    \item We use here the convention that $ \emptyset A_{Q, <t} K_Q = \emptyset $.
        We can take $ V' = \emptyset $ if $ \Gamma = \Gamma_Q $.
    \item As $ \Gamma_Q $ is a cocompact lattice in $ N_{\Gamma, Q} M_{\Gamma, Q} $ and as $ N_Q = N_{\Gamma, Q} N^{\Gamma, Q} $ by Proposition \ref{Prop Decomposition of N}, it is indeed possible to choose $ V $ as above.
    \item It follows from Lemma \ref{Lem Gamma_(P_i)|Lambda_i compact} that the set $ \Gamma_Q \bs \{ n \in N_Q \mid n w_Q Q \in \Lambda_\Gamma \smallsetminus \{eQ\} \} $ is compact. Thus, $ V' $ exists since moreover $ \Phi $ is surjective and at most $k$-to-one.
    \item If $ \Sfrak $ is a generalised Siegel set at a cusp of full rank, then $ N^{\Gamma, Q} = \{e\} $. Thus, we can take $ V' = V $. So, $ \Sfrak $ is of the form
        \[
            V A_{Q, \geq t} K_Q \ ,
        \]
        where $ A_{Q, \geq t} := \{ a \in A_Q \mid a^{\alpha_Q} \geq t \} $ (here $ t > 0 $ is fixed). We call such a set a \textit{Siegel set} with respect to $ Q $.
    \item Let $ \Sfrak_1 $ and $ \Sfrak_2 $ be two generalised Siegel sets with respect to $ Q $. Then, $ \Gamma_Q \Sfrak'_1 $ is equal to $ \Gamma_Q \Sfrak'_2 $ up to a relatively compact subset of $ \Gamma_Q \bs (X \cup \Omega_\Gamma) $.\label{Rem of Dfn Generalised Siegel set, difference in Gamma_Q U}
    \end{enumerate}
\end{rem}

\index[n]{SiegelSet@$\Sfrak$}
\begin{dfn}\label{Dfn Generalised standard Siegel set}\index{generalised standard Siegel set} \nlenum
    \begin{enumerate}
    \item 
        Let $ Q = M_Q A_Q N_Q $ be a $ \Gamma $-cuspidal parabolic subgroup of $G$.
        Let $ \omega $ be a compact subset of $N_{\Gamma, Q}$ such that the map
        \[
            \omega M_{\Gamma, Q} \to \Gamma_Q \bs N_{\Gamma, Q} M_{\Gamma, Q} \ , \quad n m \mapsto \Gamma_Q nm
        \]
        is surjective and at most $k$-to-one for some $ k \in \NN $.
        Let $ \omega' $ be an open relatively compact subset of $ N^{\Gamma, Q} $ such that the compact set
        \[
            \{ n \in N^{\Gamma, Q} \mid n w_Q Q \in N_{\Gamma, Q} (\Lambda_\Gamma \smallsetminus \{eQ\}) \}
        \]
        is contained in $ \omega' $. 
        \\ We say that a subset $ \Sfrak_{std, Q, t, \omega, \omega'} $ of $G$ is a \textit{generalised standard Siegel set} with respect to $ Q $ if it is of the form
        \[
            \omega N^{\Gamma, Q} A_Q K_Q \smallsetminus \omega \omega' A_{Q, <t} K_Q \ ,
        \]
        where $ A_{Q, < t} := \{ a \in A_Q \mid a^{\alpha_Q} < t \} $ (here $ t > 0 $ is fixed).
    \item We call $ \Sfrak' := \{ g K_Q \mid g \in \Sfrak \} $ a \textit{generalised standard Siegel set in $\mathit{G/K_Q}$}.
    \item Let $K'$ be a maximal compact subgroup of $G$ and let $ h \in G $ be such that $ K_Q = h K' h^{-1} $. Set $ X' = G/K' $. Then, we say that
        \[
            \Sfrak'_{X'} := \{g h K' \mid g K_Q \in \Sfrak' \}
        \]
        is a \textit{generalised standard Siegel set in $\mathit{X'}$}.
    \item We say that $ \Sfrak' $ is a generalised standard Siegel set at $ eQ $.
        Moreover, if $u$ is $ \Gamma $-cuspidal, then $ \Sfrak' $ is also called $ \Gamma $-cuspidal.
    \end{enumerate}
\end{dfn}

\begin{rem} \nlenum
    \begin{enumerate}
    \item We use here the convention that $ \omega \emptyset A_{Q, <t} K_Q = \emptyset $.
        We can take $ \omega' = \emptyset $ if $ \Gamma = \Gamma_Q $.
    \item As $ \Gamma_Q $ is a cocompact lattice in $ N_{\Gamma, Q} M_{\Gamma, Q} $ and as $ N_Q = N_{\Gamma, Q} N^{\Gamma, Q} $ by Proposition \ref{Prop Decomposition of N}, it is indeed possible to choose $ \omega $ as above.
    \item It follows from Lemma \ref{Lem Gamma_(P_i)|Lambda_i compact} and the fact that the image of a compact set by a continuous map between metric spaces (here $ p_{N^{\Gamma, Q}} \colon \Gamma_Q \bs N_Q \to N^{\Gamma, Q} $, the canonical projection) is compact that the set $ \{ n \in N^{\Gamma, Q} \mid n w_Q Q \in N_{\Gamma, Q} (\Lambda_\Gamma \smallsetminus \{eQ\}) \} = p_{N^{\Gamma, Q}}(\Gamma_Q \bs \{ n \in N_Q \mid n w_Q Q \in \Lambda_\Gamma \smallsetminus \{eQ\} \}) $ is compact. Thus, $ \omega' $ exists.
    \item A generalised standard Siegel set is of course a generalised Siegel set and every generalised Siegel set in $X$ is clearly contained in the union of a generalised standard Siegel set in $X$ and an open subset of $X$ which is relatively compact in $ X \cup \Omega $.
    \item If $ \Sfrak $ is a generalised standard Siegel set at a cusp of full rank, then $ N^{\Gamma, Q} = \{e\} $. Thus, we can take $ \omega' = \{e\} $. So, $ \Sfrak $ is of the form
        \[
            \omega A_{Q, \geq t} K_Q \ ,
        \]
        where $ A_{Q, \geq t} := \{ a \in A_Q \mid a^{\alpha_Q} \geq t \} $ (here $ t > 0 $ is fixed). We call such a set a \textit{standard Siegel set} with respect to $ Q $.
    \item If $ \gamma \in \Gamma_Q $, then $ \gamma \omega $ is of the same form than $ \omega $.
    \end{enumerate}
\end{rem}

\medskip

\needspace{7\baselineskip}
\underline{Example of a generalised standard Siegel set when $X$ is the upper half-space:}

\begin{center}
    \includegraphics[width = 12cm]{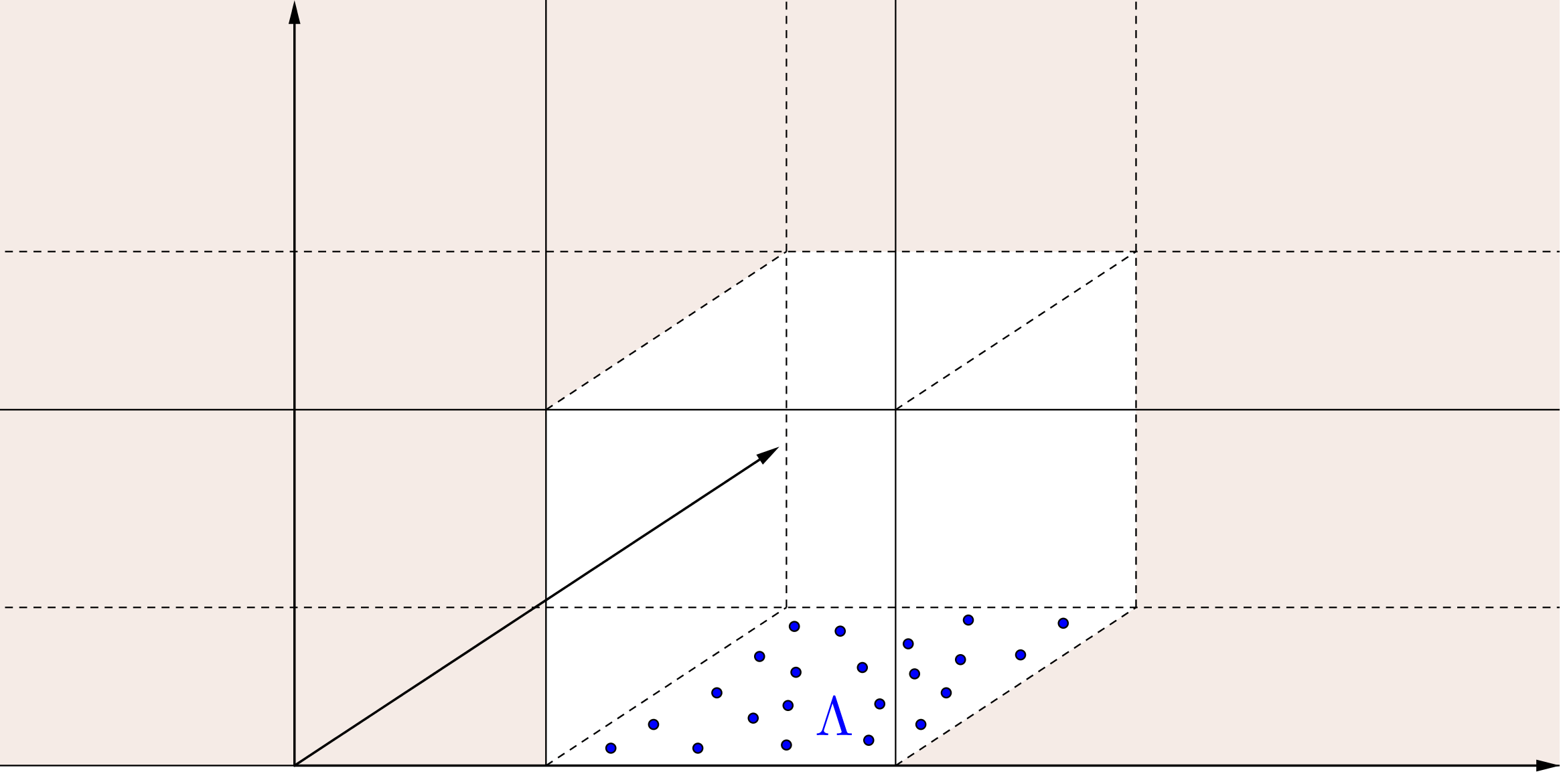}
\end{center}
N.B. Only $ \omega N^{\Gamma, Q} A_Q K_Q $ is drawn on the picture.

\medskip

Let $ \Sfrak' $ be a generalised (standard) Siegel set in $ G/K_Q $ with respect to $ Q $ and let $ \gamma \in \Gamma $. Then, $ \gamma \Sfrak' \gamma^{-1} $ is a generalised (standard) Siegel set in $ G/(\gamma K_Q \gamma^{-1}) $ with respect to $ \gamma Q \gamma^{-1} $ and $ \gamma \Sfrak' $ is a generalised Siegel set in $G/K_Q$ with respect to $ Q $.

Let $ \omega_{std, i} $ be a subset of $ N_{\Gamma, P_i} $ and let $ \omega'_{std, i} $ be a subset of $ N^{\Gamma, P_i} $ such that Definition \ref{Dfn Generalised standard Siegel set} for a generalised standard Siegel set with respect to $P_i$ holds for $ \omega = \omega_{std, i} $ and $ \omega' = \omega'_{std, i} $.

Set $ \Sfrak_{std, i} = \omega_{std, i} N^{\Gamma, P_i} A_i K_i \smallsetminus \omega_{std, i} \omega'_{std, i} (A_i)_{< 1} K_i $ (generalised standard Siegel set in $G$ with respect to $ P_i $).\index[n]{SiegelSetstdi@$ \Sfrak_{std, i} $}

Let $ h_i \in N_i A_i $ be such that $ K_i = h_i K h_i^{-1} $. Set
\[
    \Sfrak'_{X, std, i} = \{ g h_i K \in X \mid gK_i \in \Sfrak'_{std, i} \} \ .
\]
Recall that $ e_i \colon Y_i \to \bar{Y}_{\Gamma_{P_i}} $ is the inverse map of $ \restricted{\pi_i}{\tilde{Y}_i} $ (see Definition \ref{Dfn geometrically finite} for notations).

\begin{lem}\label{Lem e_i(Y_i) subset Gamma_(P_i) clo(Sfrak'_i)}
    There are generalised standard Siegel sets $ \Sfrak'_{X, i, 1}, \Sfrak'_{X, i, 2} $ in $X$ with respect to $ P_i $ such that
    \[
        \Gamma_{P_i} \clo(\Sfrak'_{X, i, 1}) \subset e_i(Y_i) \subset \Gamma_{P_i} \clo(\Sfrak'_{X, i, 2}) \ .
    \]
\end{lem}

\begin{proof}
    Let $ C_i = \bar{Y}_{\Gamma_{P_i}} \smallsetminus e_i(Y_i) $. By the proof of Lemma \ref{Lem Gamma_(P_i)|Lambda_i compact}, $ \Gamma_{P_i} \bs \Lambda_i $ is contained in $ C_i $. Thus, $ e_i(Y_i) $ is relatively compact in $ \Gamma_{P_i} \bs (X \cup \Omega_\Gamma \cup \{eP_i\}) $.

    Since $ \Gamma_{P_i} \clo(\Sfrak'_{X, std, i}) \cup \{eP_i\} $ is a compact subset of $ \Gamma_{P_i} \bs (X \cup \Omega_\Gamma \cup \{eP_i\}) $,
    \[
        C_{i, 1} := \Gamma_{P_i} \clo(\Sfrak'_{X, std, i}) \smallsetminus e_i(Y_i)
        \quad \text{and} \quad
        C_{i, 2} := e_i(Y_i) \smallsetminus \Gamma_{P_i} \clo(\Sfrak'_{X, std, i})
    \]
    are compact in $ \Gamma_{P_i} \bs (X \cup \Omega_\Gamma) $. Thus, $ C_{i, 1} \cap \Gamma_{P_i} \bs \Omega_\Gamma $ and $ C_{i, 2} \cap \Gamma_{P_i} \bs \Omega_\Gamma $ are compact in $ \Gamma_{P_i} \bs \Omega_\Gamma $.

    Let $ h_i \in G $ be such that $ K_i = h_i K h_i^{-1} $ and $ P_i = h_i P h_i^{-1} $. As moreover
    \[
        \Omega_{P_i} = \{ n w_i P_i \mid n \in N_i \} = \{ n w_i h_i P \mid n \in N_i \} \subset \dX \ ,
    \]
    $ \Gamma_{P_i} \clo(\Sfrak'_{X, std, P_i, t_1, \omega_i, \omega'_{i, 1}}) $ is contained in $ \bar{Y}_{\Gamma_{P_i}} \smallsetminus C_{i, 1} $ for some $ t_1 \geq 1 $ and some $ \omega'_{i, 1} $ containing $ \omega'_i $ and $ C_{i, 2} $ is contained in
    \[
        \Gamma_{P_i} \clo(\Sfrak'_{X, std, P_i, t_2, \omega_i, \omega'_{i, 2}})
    \]
    for some $ t_2 \leq 1 $ and $ \omega'_{i, 2} $ which is contained in $ \omega'_i $.
    Set $ \Sfrak'_{X, i, 1} = \Sfrak'_{X, std, P_i, t_1, \omega_i, \omega'_{i, 1}} $ and $ \Sfrak'_{X, i, 2} = \Sfrak'_{X, std, P_i, t_2, \omega_i, \omega'_{i, 2}} $. The lemma follows.
\end{proof}

\begin{prop}\label{Prop Decomp. of Gamma|(X cup Omega) in the geometrically finite case}\nlenum
    \begin{enumerate}
    \item Let $ \Sfrak'_{X, i} $ be generalised Siegel sets in $X$ with respect to $ P_i $ ($ i = 1, \dotsc, m $). Then, there is an open relatively compact subset $ U $ of $ X \cup \Omega_\Gamma $ such that
        \begin{equation}\label{Prop Decomp. of Gamma|(X cup Omega) in the geometrically finite case eq1}
            \Gamma \bs (X \cup \Omega_\Gamma) = \Gamma U \cup \bigcup_{i=1}^m \Gamma \clo(\Sfrak'_{X, i}) \ .
        \end{equation}
    \item There are generalised Siegel sets $ \Sfrak_i $ at $ P_i $ ($ i = 1, \dotsc, m $) and an open relatively compact subset $ U $ of $ X \cup \Omega_\Gamma $ such that
        \begin{equation}\label{Prop Decomp. of Gamma|(X cup Omega) in the geometrically finite case eq2}
            \Gamma \bs (X \cup \Omega_\Gamma) = \Gamma U \cup \Gamma \clo(\Sfrak'_{X, 1}) \cup \Gamma \clo(\Sfrak'_{X, 2}) \cup \dotsm \cup \Gamma \clo(\Sfrak'_{X, m})
        \end{equation}
        (disjoint union),
        \begin{equation}\label{Prop Decomp. of Gamma|(X cup Omega) in the geometrically finite case eq3}
            \{ \gamma \in \Gamma \mid \gamma \clo(\Sfrak'_{X, k}) \cap \clo(\Sfrak'_{X, l}) \neq \emptyset \} = \emptyset
        \end{equation}
        for all $ k \neq l $, and
        \begin{equation}\label{Prop Decomp. of Gamma|(X cup Omega) in the geometrically finite case eq4}
            \{ \gamma \in \Gamma \mid \gamma \clo(\Sfrak'_{X, i}) \cap \clo(\Sfrak'_{X, i}) \neq \emptyset \} \subset \Gamma_{P_i}
        \end{equation}
        for all $i$.
    \end{enumerate}
\end{prop}

\begin{proof} \nlenum
    \begin{enumerate}
    \item By Definition \ref{Dfn geometrically finite}, there is a compact $C$ in $ \bar{Y}_\Gamma $ such that
        \begin{equation}\label{Proof of Prop Decomp. of Gamma|(X cup Omega) in the geometrically finite case eq1}
            \bar{Y}_\Gamma = C \cup \bigcup_{i=1}^m Y_i \ .
        \end{equation}
        Let $ \Sfrak_i $ be generalised Siegel sets with respect to $ P_i $ ($ i = 1, \dotsc, m $). Without loss of generality, we may assume that $ \Sfrak_i = \Sfrak_{std, i} $.
        By Lemma \ref{Lem e_i(Y_i) subset Gamma_(P_i) clo(Sfrak'_i)}, there is an open relatively compact subset $ V $ of $ X \cup \Omega $ such that
        \[
            e_i(Y_i) \subset \Gamma_{P_i} V \cup \Gamma_{P_i} \clo(\Sfrak'_{X, std, i}) \subset \bar{Y}_{\Omega, \Gamma_{P_i}}
        \]
        for all $ i $. Applying $ \pi_i $ to this inclusion shows that $ Y_i $ is contained in $ \Gamma V \cup \Gamma \clo(\Sfrak'_{X, std, i}) $.
        Let $ U \in \U_\Gamma $ be such that $ C \cup \Gamma V $ is contained in $ \Gamma U $. The first assertion follows.
    \item By Lemma \ref{Lem e_i(Y_i) subset Gamma_(P_i) clo(Sfrak'_i)}, there are generalised Siegel sets $ \Sfrak'_{X, i} $ in $X$ with respect to $ P_i $ such that $ e_i(Y_i) $ is contained in $ \Gamma_{P_i} \clo(\Sfrak'_{X, i}) $. Thus, $ Y_i $ is contained in $ \Gamma \clo(\Sfrak'_{X, i}) $.
        Hence, $ V := \bar{Y}_\Gamma \smallsetminus \bigcup_{i=1}^m \Gamma \clo(\Sfrak'_{X, i}) $ is relatively compact in $ \bar{Y}_\Gamma $ by \eqref{Proof of Prop Decomp. of Gamma|(X cup Omega) in the geometrically finite case eq1}. Let $ U' \in \U_\Gamma $ be such that $ V $ is contained in $ \Gamma U' $.
        Set $ U = U' \smallsetminus \bigcup_{i=1}^m \clo(\Sfrak'_{X, i}) $ (open). The disjoint decomposition in \eqref{Prop Decomp. of Gamma|(X cup Omega) in the geometrically finite case eq2} follows as $ Y_i \cap Y_j = \emptyset $ if $ i \neq j $.

        As $ \Gamma \clo(\Sfrak'_{X, i}) $ is contained in $ Y_i $ for all $i$ and as $ Y_k \cap Y_l = \emptyset $ for $ k \neq l $,
        \[
            \Gamma \clo(\Sfrak'_{X, k}) \cap \Gamma \clo(\Sfrak'_{X, l}) = \emptyset
        \]
        for all $ k \neq l $. Thus, \eqref{Prop Decomp. of Gamma|(X cup Omega) in the geometrically finite case eq3} holds.

        Let $ p_i \colon X \cup \Omega_\Gamma \to \Gamma_{P_i} \bs (X \cup \Omega_\Gamma), \, x \mapsto \Gamma_{P_i} x $.

        The set $ \{ \gamma \in \Gamma \mid \gamma \clo(\Sfrak'_{X, i}) \cap \clo(\Sfrak'_{X, i}) \neq \emptyset \} $ is contained in
        \[
            \{ \gamma \in \Gamma \mid \gamma p_i^{-1}( e_i(Y_i) ) \cap p_i^{-1}( e_i(Y_i) ) \neq \emptyset \}
        \]
        and this set is contained in $ \Gamma_{P_i} $ as $ \pi_i \colon \bar{Y}_{\Omega, \Gamma_{P_i}} \to \bar{Y}_\Gamma , \, \Gamma_{P_i} x \mapsto \Gamma x $ is injective on $ e_i(Y_i) = \tilde{Y}_i $. This shows \eqref{Prop Decomp. of Gamma|(X cup Omega) in the geometrically finite case eq4}.
    \end{enumerate}
\end{proof}

\begin{cor}\label{Cor Decomp. of Gamma|X in the geometrically finite case}\nlenum
    \begin{enumerate}
    \item Let $ \Sfrak_{X, i} $ be generalised Siegel sets in $X$ with respect to $ P_i $ ($ i = 1, \dotsc, m $). Then, there is $ U \in \U_\Gamma $ such that
        \begin{equation}\label{Cor Decomp. of Gamma|X in the geometrically finite case eq1}
            \Gamma \bs X = \Gamma U \cup \bigcup_{i=1}^m \Gamma \Sfrak'_{X, i} \ .
        \end{equation}
    \item There are generalised Siegel sets $ \Sfrak_i $ at $ P_i $ ($ i = 1, \dotsc, m $) and $ U \in \U_\Gamma $ such that
        \begin{equation}
            \Gamma \bs X = \Gamma U \cup \Gamma \Sfrak'_{X, 1} \cup \dotsm \cup \Gamma \Sfrak'_{X, m}
        \end{equation}
        (disjoint union),
        \[
            \{ \gamma \in \Gamma \mid \gamma \Sfrak'_{X, k} \cap \Sfrak'_{X, l} \neq \emptyset \} = \emptyset
        \]
        for all $ k \neq l $, and
        \[
            \{ \gamma \in \Gamma \mid \gamma \Sfrak'_{X, i} \cap \Sfrak'_{X, i} \neq \emptyset \} \subset \Gamma_{P_i}
        \]
        for all $i$.
    \end{enumerate}
\end{cor}

\begin{proof}
    By the previous proposition, we have a decomposition of the form
    \[
        \Gamma \bs (X \cup \Omega_\Gamma) = \Gamma V \cup \Gamma \clo(\Sfrak'_{X, 1}) \cup \dotsm \cup \Gamma \clo(\Sfrak'_{X, m})
    \]
    with $ V $ an open relatively compact subset of $ X \cup \Omega_\Gamma $ and $ \Sfrak'_{X, i} $ generalised Siegel sets in $ X $ with respect to $P_i$.

    Set $ U = V \cap X \in \U_\Gamma $. Since $ \Sfrak'_{X, i} $ is closed in $ X $, $ \clo(\Sfrak'_{X, i}) \cap X = \Sfrak'_{X, i} $. Thus,
    \[
        \Gamma \bs X = \Gamma U \cup \bigcup_{i=1}^m \Gamma \Sfrak'_{X, i} \ .
    \]
    The corollary follows.
\end{proof}

\begin{prop}\label{Prop [gamma in Gamma | gamma (U cap cup_i Sfrak'_i) cap (U cap cup_i Sfrak'_i) neq emptyset]}
    Let $ U \in \U_\Gamma $ and let $ \Sfrak'_{X, i} $ be generalised Siegel sets in $X$ with respect to $ P_i $ ($ i = 1, \dotsc, m $). Then,
    \begin{equation}\label{eq [gamma in Gamma | gamma (U cap cup_i Sfrak'_i) cap (U cap cup_i Sfrak'_i) neq emptyset]}
        \{ \gamma \in \Gamma \mid \gamma \bigl(\clo(U) \cup \bigcup_{i=1}^m \clo(\Sfrak'_{X, i})\bigr) \cap \bigl(\clo(U) \cup \bigcup_{i=1}^m \clo(\Sfrak'_{X, i})\bigr) \neq \emptyset \}
    \end{equation}
    is finite.
\end{prop}

\begin{proof}
    By Proposition \ref{Prop Decomp. of Gamma|(X cup Omega) in the geometrically finite case}, there are generalised Siegel sets $ \Sfrak_i $ at $ P_i $ ($ i = 1, \dotsc, m $) and an open relatively compact subset $ V $ of $ X \cup \Omega_\Gamma $ such that
    \[
        \Gamma \bs (X \cup \Omega_\Gamma) = \Gamma V \cup \Gamma \clo(\Sfrak'_{X, 1}) \cup \Gamma \clo(\Sfrak'_{X, 2}) \cup \dotsm \cup \Gamma \clo(\Sfrak'_{X, m})
    \]
    (disjoint union),
    \begin{equation}\label{Proof of Prop [gamma in Gamma | gamma (U cap cup_i Sfrak'_i) cap (U cap cup_i Sfrak'_i) neq emptyset] eq1}
        \{ \gamma \in \Gamma \mid \gamma \clo(\Sfrak'_{X, i}) \cap \clo(\Sfrak'_{X, j}) \neq \emptyset \} = \emptyset
    \end{equation}
    for all $ i \neq j $, and
    \begin{equation}\label{Proof of Prop [gamma in Gamma | gamma (U cap cup_i Sfrak'_i) cap (U cap cup_i Sfrak'_i) neq emptyset] eq2}
        \{ \gamma \in \Gamma \mid \gamma \clo(\Sfrak'_{X, i}) \cap \clo(\Sfrak'_{X, i}) \neq \emptyset \} \subset \Gamma_{P_i}
    \end{equation}
    for all $i$.

    Let $ U \in \U_\Gamma $. It follows from Definition \ref{Dfn geometrically finite}, Definition \ref{Dfn Generalised Siegel set} and Remark \ref{Rem of Dfn Generalised Siegel set} \eqref{Rem of Dfn Generalised Siegel set, difference in Gamma_Q U} that we may assume without loss of generality that the $ \Sfrak'_{X, i} $'s are the given generalised Siegel sets by adding $V$ to $U$.
    Thus, \eqref{eq [gamma in Gamma | gamma (U cap cup_i Sfrak'_i) cap (U cap cup_i Sfrak'_i) neq emptyset]} is equal to
    \begin{multline*}
        \{ \gamma \in \Gamma \mid \gamma \clo(U) \cap \clo(U) \neq \emptyset \} \cup \bigcup_{i=1}^m \{ \gamma \in \Gamma \mid \gamma \clo(\Sfrak'_{X, i}) \cap \clo(\Sfrak'_{X, i}) \neq \emptyset \} \\
        \cup \bigcup_{i=1}^m \{ \gamma \in \Gamma \mid \gamma \clo(\Sfrak'_{X, i}) \cap \clo(U) \neq \emptyset \} \ .
    \end{multline*}
    By the proof of Lemma \ref{Lem e_i(Y_i) subset Gamma_(P_i) clo(Sfrak'_i)}, there is a generalised Siegel set $ \Sfrak'_{X, i, 2} $ of the form $ \Sfrak'_{X, std, P_i, t_2, \omega_{i, 2}, \omega'_{i, 2}} $ in $X$ such that
    \begin{equation}\label{eq inclusion of Siegel sets}
        \clo(\Sfrak'_{X, i}) \subset \interior\bigl(\clo(\Sfrak'_{X, i, 2})\bigr) \subset e_i(Y_i) \ .
    \end{equation}
    Here, $ \interior\bigl(\clo(\Sfrak'_{X, i, 2})\bigr) $ is equal to
    \[
        \omega_{i, 2}^\circ N^{\Gamma, P_i} A_i h_i K \smallsetminus \omega_{i, 2}^\circ \overline{\omega'_{i, 2}} A_{i, \leq t_2} h_i K
        \cup \{ n_1 n_2 w_i h_i P \mid n_1 \in \omega_{i, 2}^\circ , \, n_2 \in N^{\Gamma, P_i} \smallsetminus \overline{\omega'_{i, 2}} \} \ ,
    \]
    where $ h_i \in G $ is such that $ K_i = h_i K h_i^{-1} $ and $ P_i = h_i P h_i^{-1} $, $ \omega_{i, 2}^\circ \supset \omega_i $ is the interior of $ \omega_{i, 2} $ in $ N_{\Gamma, P_i} $ and $ \overline{\omega'_{i, 2}} $ denotes the closure of $ \omega'_{i, 2} $ in $ N^{\Gamma, P_i} $.

    Thus, $ \clo(U) $ is contained in
    $
        \bigcup_{\gamma \in \Gamma} \gamma (V \cup \bigcup_{j=1}^m \interior\bigl(\clo(\Sfrak'_{X, j, 2})\bigr))
    $.
    Hence, there are $ \gamma_1, \dotsc, \gamma_r \in \Gamma $ such that $ \clo(U) $ is contained in
    \[
        \bigcup_{k=1}^r \gamma_k (V \cup \bigcup_{j=1}^m \interior\bigl(\clo(\Sfrak'_{X, j, 2})\bigr)) \ .
    \]
    By \eqref{Proof of Prop [gamma in Gamma | gamma (U cap cup_i Sfrak'_i) cap (U cap cup_i Sfrak'_i) neq emptyset] eq1} and since moreover $ \Gamma V \cap \Gamma \clo(\Sfrak'_{X, j}) = \emptyset $ for all $j$,
    \[
        \{\gamma \in \Gamma \mid \clo(U) \cap \gamma \clo(\Sfrak'_{X, i}) \neq \emptyset \}
    \]
    is contained in
    \[
        \{\gamma \in \Gamma \mid \big(\bigcup_{k=1}^r \gamma_k \clo(\Sfrak'_{X, i, 2})\big) \cap \gamma \clo(\Sfrak'_{X, i}) \neq \emptyset \} \ .
    \]
    It follows from \eqref{eq inclusion of Siegel sets} and Definition \ref{Dfn geometrically finite} that this set is again contained in $ \Gamma_{P_i} $. Hence, the above set is finite.

    Since $ \Gamma $ acts properly continuous on $ X \cup \Omega $, $ \{ \gamma \in \Gamma \mid \gamma \clo(U) \cap \clo(U) \neq \emptyset \} $ is finite for all $i$. Thus, the proposition follows by \eqref{Proof of Prop [gamma in Gamma | gamma (U cap cup_i Sfrak'_i) cap (U cap cup_i Sfrak'_i) neq emptyset] eq2}.
\end{proof}

\begin{lem}\label{Lem gamma g h in Siegel set}
    Let $ \Sfrak_{std, Q, t, \omega, \omega'} $ be a generalised standard Siegel set with respect to $ Q $ and let $ S $ be a relatively compact subset of $G$. Then, there exists a generalised standard Siegel set $ \Sfrak_{std, Q, \tilde{t}, \omega, \tilde{\omega}'} $ such that
    \[
        \gamma_{g, h} g h \in \Sfrak_{std, Q, \tilde{t}, \omega, \tilde{\omega}'}
    \]
    for some $ \gamma_{g, h} \in \Gamma_Q $, for all $ g \in \Sfrak_{std, Q, t, \omega, \omega'} $ and all $ h \in S $.
\end{lem}

\begin{proof}
    Let $ \Sfrak := \Sfrak_{std, Q, t, \omega, \omega'} $ and $ S $ be as above and let $ h \in S $. Let $ g \in G $.
    \\ We have
    \begin{align*}
        g h &= \nu_Q(g h) h_Q(gh) k_Q(gh) \\
        &= \big( \nu_Q(g) h_Q(g) \nu_Q(k_Q(g)h) h_Q(g)^{-1} \big) h_Q(g) h_Q(k_Q(g)h) k_Q(gh) \ . \nonumber
    \end{align*}
    Let $ \tilde{\omega}' $ be such that its closure in $ N^{\Gamma, Q} $ is strictly contained in $ \omega' $ and such that $ \Sfrak_{std, Q, t, \omega, \tilde{\omega}'} $ is a generalised standard Siegel set with respect to $ Q $.

    Let $ D = \overline{\tilde{\omega}'} $ (compact) and $ U = N_{\Gamma, Q} \omega' $ (open).
    \\ Then, $ V = \{ n \in N_Q \mid D n \subset U , \, D n^{-1} \subset U \} $ is an open neighbourhood of $ e $. Moreover, $ N_{\Gamma, Q} D V \subset U $.

    By definition of $V$, $ N_{\Gamma, Q} D $ is contained in $ U n $ for all $ n \in V $. Hence, $ N_Q \smallsetminus U n $ is contained in $ N_Q \smallsetminus N_{\Gamma, Q} D $ for all $ n \in V $.

    Let $ n \in V $. Then, $ N_Q $ is the disjoint union of $ U n $ and $ (N_Q \smallsetminus U) n $.

    Thus, $ (N_Q \smallsetminus U) n = N_Q \smallsetminus U n $. Hence, $ V (N_Q \smallsetminus U) = \bigcup_{n \in V} (N_Q \smallsetminus U n) \subset N_Q \smallsetminus N_{\Gamma, Q} D $.

    Let $ \delta \in (0, t) $ be such that
    \[
        \nu_Q(g) a \nu_Q(k h) a^{-1} \in V (N_Q \smallsetminus N_{\Gamma, Q} \omega') \subset N_Q \smallsetminus N_{\Gamma, Q} \overline{\tilde{\omega}'}
    \]
    for all $ a \in A_{Q, \leq \delta} $, $ k \in K_Q $, $ h \in S $ and $ g \in G $ such that $ p_{N^{\Gamma, Q}}(\nu_Q(g)) \in N^{\Gamma, Q} \smallsetminus \omega' $.

    Let $ \tilde{t} = \delta \inf_{k \in K_Q, h \in S} h_Q(k h)^{\alpha_Q} > 0 $. Then,
    \[
        g h \in N_Q A_{Q, \geq \tilde{t}} K_Q \subset \Gamma_Q \Sfrak_{std, Q, \tilde{t}, \omega, \tilde{\omega}'}
    \]
    for all $ g \in \omega N^{\Gamma, Q} A_{Q, \geq \delta} K_Q $ and all $ h \in S $.
\end{proof}

\begin{prop}\label{Prop Gamma|G = Gamma U K cup Cup_(i=1)^m Sfrak_i}
    There is $ U_\Gamma \in \U_\Gamma $ such that
    \[
        \Gamma \bs G = \Gamma U_\Gamma K \cup \bigcup_{i=1}^m \Sfrak_{std, i} \ .
    \]
\end{prop}

\begin{proof}
    By Corollary \ref{Cor Decomp. of Gamma|X in the geometrically finite case}, there is $ U \in \U_\Gamma $ such that
    \[
        \Gamma \bs X = \Gamma U \cup \bigcup_{i=1}^m \Gamma \Sfrak'_{X, std, i} \ .
    \]
    Let $ h_i \in G $ be such that $ K_i = h_i K h_i^{-1} $ and $ P_i = h_i P h_i^{-1} $. Thus,
    \[
        \Gamma \bs G = \Gamma U K \cup \bigcup_{i=1}^m \Gamma \{ g h_i k \mid g \in \Sfrak_{std, i}, k \in K \} \ .
    \]
    By Lemma \ref{Lem gamma g h in Siegel set}, there exists a generalised standard Siegel set $ \Sfrak_{std, P_i, \tilde{t}_i, \omega_i, \tilde{\omega}'_i} $ such that
    \[
        \{ g h_i k \mid g \in \Sfrak_{std, i}, k \in K \} \subset \bigcup_{\gamma \in \Gamma_{P_i}} \gamma \Sfrak_{std, P_i, \tilde{t}_i, \omega_i, \tilde{\omega}'_i} \ .
    \]
    Hence,
    \[
        \Gamma \bs G = \Gamma U K \cup \bigcup_{i=1}^m \Gamma \Sfrak_{std, P_i, \tilde{t}_i, \omega_i, \tilde{\omega}'_i} \ .
    \]
    Let $ U_i $ be an open subset of $ G/K_i $ such that $ U_i $ is relatively compact in
    \[
        G/K_i \cup \{gP_i \mid g h_i P \in \Omega_\Gamma\}
    \]
    and such that $ \Sfrak_{std, P_i, \tilde{t}_i, \omega_i, \tilde{\omega}'_i} $ is contained in
    \[
        U_i K_i \cup \Sfrak_{std, i} \subset \{gk \in G \mid gK_i \in U_i, k \in K\} \cup \Sfrak_{std, i} \ .
    \]
    Then, $ U'_i := \{ g h_i K \mid g K_i \in U_i \} $ is an open subset of $X$ which is relatively compact in $ X \cup \Omega_\Gamma $ ($ U'_i \in \U_\Gamma $).
    By Corollary \ref{Cor union of U^h in U_Gamma},
    \[
        \{gK \mid gK_i \in U_i\} = \{ g h_i^{-1} K \mid gK \in U'_i \}
    \]
    belongs also to $ \U_\Gamma $. The proposition follows.
\end{proof}

\newpage

\subsection{Computation of the critical exponent of \texorpdfstring{$ \Gamma_Q $}{Gamma\_Q} and definition of the rank of a cusp}

We define the rank of a cusp and we determine an explicit formula for the critical exponent of $ \Gamma_Q $ (see Proposition \ref{Prop Critical exponent of delta_(Gamma_Q)}).

\medskip

\begin{dfn}\nlenum
    \begin{enumerate}
    \item The \textit{rank} of the cusp associated to $ [Q]_\Gamma \in \P_\Gamma $ is by definition $ \dim N_{\Gamma, Q} $.\index{rank of a cusp}
    \item We say that the cusp associated to $ [Q]_\Gamma \in \P_\Gamma $ is \textit{of full rank} if its rank is equal to $ \dim N_Q $, or equivalently, if $ \Gamma_Q \bs N_Q M_{\Gamma, Q} $ is compact. Otherwise we say that the cusp associated to $ [Q]_\Gamma $ is \textit{of smaller rank}.
    \end{enumerate}
\end{dfn}

\needspace{2\baselineskip}
Let $ r_1 = \dim p_{\g_{\alpha_Q}}(\n_{\Gamma, Q}) $ and $ r_2 = \dim(\g_{2\alpha_Q} \cap \n_{\Gamma, Q}) $.

\begin{prop}\label{Prop Critical exponent of delta_(Gamma_Q)}
    The critical exponent $ \delta_{\Gamma_Q} $ of $ \Gamma_Q $ is equal to
    \[
        \Big(\frac{r_1}2 + r_2\Big)\alpha_Q - \rho_Q = \Big(\frac{r_1 - m_{\alpha_Q}}2 + r_2 - m_{2\alpha_Q}\Big)\alpha_Q \ .
    \]
\end{prop}

\begin{rem}\label{Rem delta_(Gamma_Q) <= 0}
    A direct consequence of this proposition is that the critical exponent $ \delta_{\Gamma_Q} $ of $ \Gamma_Q $ belongs to $ (-\rho_Q, 0] $.
    \\ If $ eQ $ is of smaller rank, then $ \delta_{\Gamma_Q} $ belongs even to $ (-\rho_Q, 0) $.
\end{rem}

\begin{proof}
    Let $ \lambda_Q \in \a^*_Q \equiv \RR $. If $ \gamma \in \Gamma_Q := \Gamma \cap N_Q M_Q $, then $ a_\gamma = a_{\nu_Q(\gamma)} $.
    Thus, $ \sum_{\gamma \in \Gamma_Q}{ a_\gamma^{-(\lambda_Q+\rho_Q)} } $ is finite if and only if
    \begin{equation}\label{Proof of Prop Critical exponent of delta_(Gamma_Q) eq1}
        \sum_{\gamma \in \Gamma_Q}{ a(\theta_Q \nu_Q(\gamma))^{-(\lambda_Q+\rho_Q)} }
    \end{equation}
    is finite, by Corollary \ref{Cor a_n asymp a(theta n)}.
    \\ Let $ n \in N_Q $. Write $ n = \exp(X+Y) $ ($ X \in \g_{\alpha_Q} $, $ Y \in \g_{2\alpha_Q} $). By Proposition \ref{Theorem 3.8 of Helgason, p.414},
    \[
        a(\theta_Q n) = \sqrt{(1 + \tfrac12|X|^2)^2 + 2 |Y|^2} \ .
    \]
    Since $ (x, y) \mapsto \big((1 + \tfrac12x^2)^2 + 2 y^2\big)^{-\frac{\lambda_Q+\rho_Q}{2}} $ is a continuous, nonnegative function that is decreasing in both variables, \eqref{Proof of Prop Critical exponent of delta_(Gamma_Q) eq1} is finite if and only if
    \[
        \int_{N_{\Gamma, Q}} a(\theta_Q n)^{-(\lambda_Q+\rho_Q)} \, dn
    \]
    is finite.
    Let $ N'_1 = \exp\bigl(p_{\g_{\alpha_Q}}(\n_{\Gamma, Q})\bigr) $ and $ N'_2 = \exp(\g_{2\alpha_Q} \cap \n_{\Gamma, Q}) $. Then, it follows from the explicit formulas that the above integral is finite if and only if
    \[
        \int_{N'_1 N'_2} a(\theta_Q n)^{-(\lambda_Q+\rho_Q)} \, dn
    \]
    is finite (we will prove a much more general formula in Lemma \ref{Lem a_(n_(v_1, phi(v_1) + r_1) n_(v_2, r_2) a) asymp a_(n_(v_1 + v_2, r_1 + r_2) a)}). By Proposition \ref{Prop N bar integral is finite}, this is case if and only if $ \lambda_Q > \big(\frac{r_1}2 + r_2\big)\alpha_Q - \rho_Q $.
\end{proof}

\index[n]{rhoGammaQ@$ \rho^{\Gamma_Q} $, $ \rho_{\Gamma_Q} $}
For $ H \in \a_Q $, set
\[
    \rho^{\Gamma_Q}(H) = \frac12 \tr(\restricted{\ad(H)}{\n^{\Gamma, Q}})
\]
and
\[
    \rho_{\Gamma_Q}(H) = \rho_Q(H) - \rho^{\Gamma_Q}(H) \ .
\]

\begin{lem}\label{Lem rho_(Gamma_Q) = delta_(Gamma_Q)}
    We have:
    \[
        \rho_{\Gamma_Q} = \delta_{\Gamma_Q} \ .
    \]
\end{lem}

\begin{proof}
    For $ H \in \a_Q $, we have
    \begin{align*}
        \rho^{\Gamma_Q}(H) &= \tfrac12 \big( \dim(\n^{\Gamma, Q}_1) \alpha_Q(H) + 2 \dim(\n^{\Gamma, Q}_2) \alpha_Q(H) \big) \\
        &= \tfrac12 \big( (m_{\alpha_Q} - r_1) \alpha_Q(H) + 2 (m_{2\alpha_Q} - r_2) \alpha_Q(H) \big) \\
        &= \rho_Q(H)  - (\tfrac12 r_1 + r_2) \alpha_Q(H) \ .
    \end{align*}
    Thus, $ \rho_{\Gamma_Q} = \delta_{\Gamma_Q} $ by the previous proposition.
\end{proof}

\newpage

\subsection{The Schwartz space and the space of cusp forms on \texorpdfstring{$ \Gamma \bs G $}{Gamma\textbackslash G}}

\subsubsection{The Schwartz space on \texorpdfstring{$ \Gamma \bs G $}{Gamma\textbackslash G}}

In this section, we define the Schwartz space on $ \Gamma \bs G $ (in the geometrically finite case) and we show that a Schwartz function is indeed square-integrable.

\medskip

Let $ E = \exp(\e_1 \oplus \e_2) $ for some vector subspaces $ \e_1 \subset \g_{\alpha_Q} $ and $ \e_2 \subset \g_{2\alpha_Q} $. Let $ n_1 = \dim_\RR(\e_1) $, let $ n_2 = \dim_\RR(\e_2) $ and let $ \xi_{Q, E} = \frac{n_1 + 2 n_2}{m_{\alpha_Q} + 2 m_{2\alpha_Q}} $.\index[n]{xiQE@$ \xi_{Q, E} $}

Recall that we put as measure on all measurable subsets $S$ of $ N_Q $ the push-forward of the Lebesgue measure on $ \log(S) $ (with an appropriate normalisation).

The following two lemmas are needed to ensure that the Schwartz space on $ \Gamma \bs G $, which we will define in a moment, lies in $ L^2(\Gamma \bs G, \phi) $.

\begin{lem}\label{Lem Estimate of int_CE a_na^(-xi rho) (1 + log a_na)^(-2d)dn}
    Let $ t \geq \xi_{Q, E} $ and let $ C $ be a compact subset of $N_Q$.
    For $ a \in A_Q $, $d > 1$, $ n \mapsto a_{n a}^{-t \rho_Q} (1 + \log a_{n a})^{-d} $ is integrable over $ CE $ and for all $ \eps \in [0, d-1) $, there exists a positive constant $ c $ depending only on $ d $ such that
    \[
        \int_{CE}{ a_{n a}^{-t \rho_Q} (1 + \log a_{n a})^{-d} \, dn } \leq c \, a^{t \rho_Q} (1 + |\log a|)^{-\eps} \ .
    \]
\end{lem}

\begin{rem}\label{Rem of Lem Estimate of int_CE a_na^(-xi rho) (1 + log a_na)^(-2d)dn}
    $ \rho^{\Gamma_Q} $ is equal to $ \xi_{Q, N^{\Gamma, Q}} \rho_Q $. Indeed:
    \[
        \xi_{Q, N^{\Gamma, Q}} \rho_Q = \frac{\dim(\n^{\Gamma, Q}_1) + 2 \dim(\n^{\Gamma, Q}_2)}{m_{\alpha_Q} + 2 m_{2\alpha_Q}} \rho_Q
        = \frac12\bigl(\dim(\n^{\Gamma, Q}_1) + 2 \dim(\n^{\Gamma, Q}_2)\bigr) \alpha_Q = \rho^{\Gamma_Q} \ .
    \]
\end{rem}

\begin{proof}
    Let $ \eps \in [0, d-1) $ and $ a \in A_Q $. We have
    \[
        \int_{C E}{ a_{na}^{-t \rho_Q} (1 + \log a_{na})^{-d} \, dn}
        =\int_{C E}{ a_{xya}^{-t \rho_Q} (1 + \log a_{xya})^{-d} \, dy \, dx } \ .
    \]
    By the proof of Lemma \ref{Lem Estimate for int_N |f(g n h)| dn}, this is less or equal than
    \[
        c \, a^{t \rho_Q} (1 + |\log a|)^{-\eps} \int_C \int_E a_Q(\theta_Q(x y))^{-t \rho_Q} (1 + \log a_{\theta_Q(x y)} )^{-d+\eps} \, dy \, dx
    \]
    for some constant $ c > 0 $. The lemma follows now from Lemma \ref{Lem about a(g)}, Lemma \ref{Lem 1 + log a_h <= (1 + log a_gh)(1 + log a_g)} and Proposition \ref{Prop N bar integral is finite}.
\end{proof}

\begin{lem}\label{Lem int_CEAK( a(g^(-1))^(-2(1-xi) rho)(1 + log a_g|)^(-2r) a_g^(-2xi rho) dg < infty}
    Let $ r > 1 $ and $ \lambda_Q \in \a_Q^*  $. Let $ C $ be a compact subset of $N_Q$. Then,
    \[
        \int_{C E A_Q K_Q}{ a_Q(g^{-1})^{-2(1-\xi_{Q, E})\rho_Q}(1 + \log a_g)^{-2r} a_g^{-2\xi_{Q, E} \rho_Q} \, dg } < \infty \ .
    \]
\end{lem}

\begin{proof}
    It follows from our choice of $ \xi_{Q, E} $ that $ \xi_{Q, E} \in [0, 1] $.
    We have
    \begin{multline*}
        \int_{C E A_Q K_Q}{ a_Q(g^{-1})^{-2(1-\xi_{Q, E})\rho_Q}(1 + \log a_g)^{-2r} a_g^{-2\xi_{Q, E}\rho_Q} \, dg } \\
        = \int_{C E A_Q}{ a^{2(1-\xi_{Q, E})\rho_Q} (1 + \log a_{na})^{-2r} a_{na}^{-2\xi_{Q, E} \rho_Q} a^{-2\rho_Q} \, dn \, da } \ .
    \end{multline*}
    By Lemma \ref{Lem Estimate of int_CE a_na^(-xi rho) (1 + log a_na)^(-2d)dn}, this is again less or equal than
    \[
        c \int_{A_Q}{ a^{-2\xi_{Q, E}\rho_Q} a^{2\xi_{Q, E} \rho_Q} (1 + |\log a|)^{-r} \, da } = 2c \int_0^\infty{ (1 + t)^{-r} \, dt } \ .
    \]
    The lemma follows as the last integral is finite for $ r > 1 $.
\end{proof}

As in the convex cocompact case, we denote the family of open subsets of $X$ which are relatively compact in $ X \cup \Omega_\Gamma $ by $ \U_\Gamma $ and for $ f \in \Cinf(G, V_\phi), X, Y \in \U(\g), r \geq 0 $, $ U \in \U_\Gamma $, we set
\[
    \LOI{p}{U}_{r, X, Y}(f) = \sup_{gK \in U}{ (1 + \log a_g)^r a_g^\rho |L_X R_Y f(g)| } \ .
\]
Let $ P_1, \dotsc, P_m $ be representatives of the $ \Gamma $-conjugacy classes of $ \Gamma $-cuspidal parabolic subgroups of $G$.

Let $ \theta_i = \theta_{P_i} $, $ K_i = K_{P_i} $, $ N_i = N_{P_i} $, $ \bar{N}_i = \bar{N}_{P_i} $, $ M_i = M_{P_i} $ and $ A_i = A_{P_i} $.
\index[n]{K_i M_i A_i N_i@$ K_i $, $ M_i $, $ A_i $, $ N_i $}

\index[n]{VGammaVstdGamma@$ \V_\Gamma $, $ \V_{std, \Gamma} $}
Let us denote the family of generalised Siegel sets (resp. generalised standard Siegel sets) in $G/K_i$ with respect to $P_i$ ($ i = 1, \dotsc, m $) by $ \V_\Gamma $ (resp. $ \V_{std, \Gamma} $).
Let $ l \in \NN_0 $. For a generalised Siegel set $ \Sfrak' \in \V_\Gamma $ in $ G/K_Q $ with respect to $ Q $, $ r \geq 0 $, $ X \in \U(\g)_{l_1} $, $ Y \in \U(\g)_{l_2} $ such that $ l_1 + l_2 \leq l $ and a function $ f \in \C^l(G, V_\phi) $, put\index[n]{prXYf2@$ \LO{p}{\Sfrak'}_{r, X, Y}(f) $}
\[
    \LO{p}{\Sfrak'}_{r, X, Y}(f) = \sup_{gK_Q \in \Sfrak'}{ a_Q(g^{-1})^{\rho_{\Gamma_Q}} (1 + \log a_g)^r a_g^{\rho^{\Gamma_Q}} |L_X R_Y f(g)| } \ .
\]
Set $ \xi_i := \xi_{P_i, N^{\Gamma, P_i}} $.\index[n]{xii@$ \xi_i $}

\begin{rem}\label{Rem Schwartz space well-defined} \nlenum
    \begin{enumerate}
    \item We recall that we write for $ a_{Q, g}^{\alpha_Q} $ simply $ a_g $ by abuse of notation.
    \item Let $ \xi = \xi_{Q, N^{\Gamma, Q}} $. Then, $ \rho_{\Gamma_Q} = (1-\xi)\rho_Q = \delta_{\Gamma_Q} $ and $ \rho^{\Gamma_Q} = \xi \rho_Q = \rho_Q - \delta_{\Gamma_Q} $ by definition of $ \rho_{\Gamma_Q} $, by Remark \ref{Rem of Lem Estimate of int_CE a_na^(-xi rho) (1 + log a_na)^(-2d)dn} and by Lemma \ref{Lem rho_(Gamma_Q) = delta_(Gamma_Q)}.
    \end{enumerate}
\end{rem}

\begin{lem}\label{Lem c_1 a_n a^(-1) <= a_(na) <= c_2 a_n a^(-1)}
    \nlenum
    \begin{enumerate}
    \item There exists $ c \in (0, 1] $ such that
        \[
            c \, a_n a^{-1} \leq a_{na}
        \]
        for all $ n \in N_Q $ and $ a \in A_Q $.
    \item Let $ r > 0 $. Then, there exists $ c \geq 1 $ such that
        \[
            a_{na} \leq c \, a_n a^{-1}
        \]
        for all $ n \in N_Q $ and $ a \in A_{Q, \leq r} $.
        Moreover, if $ r \leq 1 $, then we can take $ c = 1 $.
    \end{enumerate}
\end{lem}

\begin{proof} \nlenum
    \begin{enumerate}
    \item Let $ n \in N_Q $ and $ a \in A_Q $. Let $W$ be a compact subset of $ \dX $ containing all $ k_n M_Q $, $ n \in N_Q $, but not $ w_Q M_Q $. This is possible as $N_Q$ is the horosphere passing through $ eK \in X $ and $ eM_Q \in \dX $. Let $U$ be a neighbourhood of $ k_0 := w_Q $ in $K$ satisfying $ \bar{U}M_Q \cap W = \emptyset $. Then, by Corollary 2.4 of \cite[p.85]{BO00}, there exists $ c > 0 $ such that
        \[
            a_{na} = a_{n w_Q a^{-1}} \geq c \, a_n a^{-1} \ .
        \]
    \item Let $ r > 0 $. For $ a \in A_{Q, \leq r} $, we have
        \[
            a_{na} \leq a_n a_a \leq \max\{1, r^2\} a_n a^{-1}
        \]
        as $ a_h h \leq h^2 \leq r^2 $ for all $ h \in A_{Q, \geq 1} \cap A_{Q, \leq r} $ and $ a_h = h^{-1} $ for all $ h \in A_{Q, \leq 1} $. In particular, $ \max\{1, r^2\} $ is equal to 1 if $ r \leq 1 $.
    \end{enumerate}
    This completes the proof of the lemma.
\end{proof}

\begin{prop}\label{Prop seminorms are equivalent}
    Let $ f \in \Cinf(\Gamma \bs G, \phi) $, $ \Sfrak' $ be a generalised Siegel set with respect to $ Q $, $ X, Y \in \U(\g) $, $ r \geq 0 $ and $ U \in \U_\Gamma $ be such that $ U $ is also contained in $ \Sfrak' $. Then, there exists a constant $ c \in (0, 1] $ such that
    \[
        c \LOI{p}{U}_{r, X, Y}(f) \leq \LO{p}{\Sfrak'}_{r, X, Y}(f) \ .
    \]
\end{prop}

\begin{proof}
    As $ eM_Q \not \in \clo(U) \cap \dX \subset G/P $, $ \sup_{gK \in U} a_{n_Q(g^{-1})} $ is bounded. Thus, there exists a constant $ c \in (0, 1] $ such that
    \[
        c^{\rho_{\Gamma_Q}} a_g^{\rho_{\Gamma_Q}} \leq a_Q(g^{-1})^{\rho_{\Gamma_Q}}
    \]
    by Lemma \ref{Lem c_1 a_n a^(-1) <= a_(na) <= c_2 a_n a^(-1)} applied to $ n = n_Q(g^{-1})^{-1} $ and $ a = a_Q(g^{-1})^{-1} $. Hence,
    \[
        c^{\rho_{\Gamma_Q}} \LOI{p}{U}_{r, X, Y}(f) \leq \LO{p}{\Sfrak'}_{r, X, Y}(f)  \ .
    \]
    This finishes the proof of the proposition.
\end{proof}

\begin{lem}\label{Lem X in U(n) enough}
    Let $ (V, |\cdot|) $ be a finite-dimensional normed complex vector space, let $ f \in \Cinf(G, V) $, $ X \in \U(\g)_{l_1} $, $ Y \in \U(\g)_{l_2} $ and $ U \in \U_\Gamma $.

    Let $ \{Z_1, \dotsc, Z_{d_{l_1}}\} $ be a basis of $ \U(\n_Q)_{l_1} $, $ \{Y_1, \dotsc, Y_{d_{l_2}}\} $ be a basis of $ \U(\g)_{l_2} $.
    \\ If $ \Gamma \bs X $ has a cusp, then we can write
    \[
        L_X R_Y f(g) = \sum_{i=1}^{d_{l_1}} \sum_{j=1}^{d_{l_2}} c_{i, j}(g) L_{Z_i} R_{Y_j} f(g)   \qquad (gK \in U)
    \]
    with $ \sup_{gK \in U} |c_{i, j}(g)| < \infty $.
\end{lem}

\begin{proof}
    Let $ f \in \Cinf(G, V) $, $ X \in \g, Y \in \U(\g) $ and $U$ as above. Without loss of generality, we may assume that $ Y = 1 $.

    Since $ U \in \U_\Gamma $ and since $ \Gamma \bs X $ has a cusp, there is a generalised standard Siegel set $ \Sfrak' $ with respect to $ Q $ such that $ U $ is contained in $ \Sfrak' $.

    It follows that there exist $ c_1 > 0 $ and $ c_2 > 0 $ such that $ |\log(\nu_Q(g))| \leq c_1 $ and $ h_Q(g)^{\alpha_Q} \leq c_2 $ for all $ gK \in U $.

    Let $ \{X_1, \dotsc, X_d\} $ be a basis of $ \g = \n_Q \oplus \m_Q \oplus \a_Q \oplus \bar{\n}_Q $ induced by the inclusions
    \[
        0 \subset \n_Q \subset \n_Q \oplus \m_Q \oplus \a_Q \subset \n_Q \oplus \m_Q \oplus \a_Q \oplus \g_{-\alpha_Q} \subset \g \ .
    \]
    Write $ \Ad(n)X = \sum_{j=1}^d a_j(n) X_j $ ($ n \in N_Q $). Since $ \Ad(n)X $ depends smoothly on $ n \in N_Q $, all the functions $ a_j $ are smooth. Thus,
    \[
        \sup_j \sup_{gK \in U} |a_j(\nu_Q(g)^{\pm1})| \leq \sup_j \sup_{n \in N_Q \, : \, |\log(n)| \leq c_1} |a_j(n)| < \infty \ .
    \]
    Since
    \begin{multline*}
        \Ad(g^{-1})X = \Ad(k_Q(g)^{-1}) \Ad(h_Q(g)^{-1}) \Ad(\nu_Q(g)^{-1})X \\
        = \sum_{j=1}^d a_j(\nu_Q(g)^{-1}) \Ad(k_Q(g)^{-1}) \Ad(h_Q(g)^{-1}) X_j \ ,
    \end{multline*}
    \[
        L_X f(g) = \restricted{ \frac{d}{dt} }{ t = 0 } f(\exp(-t X) g) = \restricted{ \frac{d}{dt} }{ t = 0 } f(g \exp(-t \Ad(g^{-1})X))
    \]
    is equal to
    \begin{align*}
        & \sum_{j=1}^{\dim(\n_Q)} a_j(\nu_Q(g)^{-1}) \restricted{ \frac{d}{dt} }{ t = 0 } f(\exp(-t \Ad(\nu_Q(g))X_j) g) \\
        & \qquad + \sum_{j=\dim(\n_Q)+1}^d a_j(\nu_Q(g)^{-1}) \restricted{ \frac{d}{dt} }{ t = 0 } f(g \exp(-t \Ad(k_Q(g)^{-1}) \Ad(h_Q(g)^{-1})  X_j)) \ .
    \end{align*}
    Since $ \restricted{\Ad(a)}{ \m_Q \oplus \a_Q } = \Id $ for all $ a \in A_Q $ and $ \Ad(a)Z = a^{r \alpha_Q} Z $ for all $ a \in A_Q $ and $ Z \in g_{r \alpha_Q} $ ($ r \in \{\pm 1, \pm 2\}) $, this is again equal to
    \begin{align*}
        & \sum_{j \, : \, X_j \in \n_Q} a_j(\nu_Q(g)^{-1}) a_j(\nu_Q(g)) L_{X_j} f(g)
        - \sum_{j \, : \, X_j \in \m_Q \oplus \a_Q} a_j(\nu_Q(g)^{-1}) R_{\Ad(k_Q(g)^{-1}) X_j } f(g) \\
        & \qquad - \sum_{j \, : \, X_j \in \g_{-\alpha_Q}} a_j(\nu_Q(g)^{-1}) h_Q(g)^{\alpha_Q} R_{\Ad(k_Q(g)^{-1}) X_j } f(g) \\
        & \qquad - \sum_{j \, : \, X_j \in \g_{-2\alpha_Q}} a_j(\nu_Q(g)^{-1}) h_Q(g)^{2\alpha_Q} R_{\Ad(k_Q(g)^{-1}) X_j } f(g) \ .
    \end{align*}
    Since $ \Ad(g)X $ ($ g \in G $, $ X \in \g $) depends smoothly on $g$, $K_Q$ is compact, \\
    $ \sup_{gK \in U} h_Q(g)^{\alpha_Q} \leq c_2 $ and $ \sup_{gK \in U} h_Q(g)^{2\alpha_Q} \leq c_2^2 $, the assertion of the lemma holds for $ X \in \g $. The lemma follows now by induction.
\end{proof}

\begin{cor}\label{Cor X in U(n) enough for seminorms}
    Let $ f \in \Cinf(\Gamma \bs G, \phi) $, $ X, Y \in \U(\g) $, $ r \geq 0 $ and $ U \in \U_\Gamma $.
    If $ \Gamma \bs X $ has a cusp, then
    $
        \LOI{p}{U}_{r, X, Y}(f)
    $
    is finite if $ \LOI{p}{U}_{r, Z_1, Z_2}(f) $ is finite for all $ Z_1 \in \U(\n_Q) $ and $ Z_2 \in \U(\g) $.
\end{cor}

\begin{proof}
    The corollary follows directly from Lemma \ref{Lem X in U(n) enough}.
\end{proof}

\begin{dfn}\label{Dfn Schwartz space in the geometrically finite case}
    Let $ (\phi, V_\phi) $ be a finite-dimensional unitary representation of $ \Gamma $. 
    Then, we define the Schwartz space on $ \Gamma \bs G $ by\index{Schwartz space}\index[n]{CGammaG@$ \CS(\Gamma \bs G, \phi) $}
    \begin{multline*}
        \CS(\Gamma \bs G, \phi) = \{ f \in \Cinf(\Gamma \bs G, \phi) \mid \LOI{p}{U}_{r, X, Y}(f) < \infty \; \forall r \geq 0, X, Y \in \U(\g), \\
        \LO{p}{\Sfrak'_Q}_{r, X, Y}(f) < \infty \; \forall \, U \in \U_\Gamma, \Sfrak'_Q \in \V_{std, \Gamma}, X \in \U(\n_Q), Y \in \U(\g), r \geq 0 \} \ .
    \end{multline*}
    We equip $ \CS(\Gamma \bs G, \phi) $ with the topology induced by the seminorms.

    By definition, $ \CS(\Gamma \bs G) = \CS(\Gamma \bs G, 1) $, where $1$ is the one-dimensional trivial representation.
\end{dfn}

\begin{rem}\label{Rem of Dfn Schwartz space in the geometrically finite case} \nlenum
    \begin{enumerate}
    \item It will follow from Proposition \ref{Prop Schwartz space is independent of choices} that the Schwartz space does not depend on made choices.
    \item By Lemma \ref{Lem p_(gamma U)_(r, X, Y)(f) < infty}, Proposition \ref{Prop Gamma|G = Gamma U K cup Cup_(i=1)^m Sfrak_i}, Proposition \ref{Prop [gamma in Gamma | gamma (U cap cup_i Sfrak'_i) cap (U cap cup_i Sfrak'_i) neq emptyset]}, Proposition \ref{Prop seminorms are equivalent} and Corollary \ref{Cor X in U(n) enough for seminorms}, the space $ \CS(\Gamma \bs G, \phi) $ is also equal to (resp. topologically isomorphic to)
        \begin{multline*}
            \{ f \in \Cinf(\Gamma \bs G, \phi) \mid \LOI{p}{U_\Gamma}_{r, X, Y}(f) < \infty \; \forall r \in \NN_0, X, Y \in \U(\g), \\
            \LO{p}{\Sfrak'_{std, i}}_{r, X, Y}(f) < \infty \; \forall r \geq 0, X \in \U(\n_{P_i}), Y \in \U(\g), i \in \{1, \dotsc, m\} \} \ ,
        \end{multline*}
        where $U_\Gamma$ is provided by Proposition \ref{Prop Gamma|G = Gamma U K cup Cup_(i=1)^m Sfrak_i}. Compare with Lemma \ref{Lem simpler version of the Schwartz space}.
    \item The $ \CS(\Gamma \bs G, \phi) $ has a Fréchet space structure. This follows by the proof of Proposition \ref{Prop Schwartz space is Fréchet is the convex-cocompact case}. Indeed, one uses in the proof the seminorms we have here and the previous remark instead of Lemma \ref{Lem simpler version of the Schwartz space}.
    \item This definition is an extension of Definition \ref{Dfn Schwartz space in the covex-cocompact case} as $ m = 0 $ gives the old definition.
    \end{enumerate}
\end{rem}

\begin{prop}
    The Schwartz space $ \CS(\Gamma \bs G, \phi) $ is contained in $ L^2(\Gamma \bs G, \phi) $.
    \\ Moreover, the injection of $ \CS(\Gamma \bs G, \phi) $ into $ L^2(\Gamma \bs G, \phi) $ is continuous.
\end{prop}

\begin{proof}
    Let $ f \in \CS(\Gamma \bs G, \phi) $.
    Fix $ r > 1 $. Set $ S_0 = U_\Gamma K $, where $ U_\Gamma $ is provided by Proposition \ref{Prop Gamma|G = Gamma U K cup Cup_(i=1)^m Sfrak_i}.
    Put $ l_0(f) = \LOI{p}{U_\Gamma}_{r, 1, 1}(f) $ and $ l_i(f) = \LO{p}{\Sfrak'_{std, i}}_{r, 1, 1}(f) $. Then,
    \[
        |f(g)| \leq l_0(f) (1 + \log a_g)^{-r} a_g^{-\rho}
    \]
    for all $ g \in S_0 \subset G $ and
    \[
        |f(g)| \leq l_i(f) a_{P_i}(g^{-1})^{-(1-\xi_i)\rho_{P_i}} (1 + \log a_g)^{-r} a_g^{-\xi_i \rho_{P_i}}
    \]
    for all $ g \in \Sfrak_{std, i} \subset G $.
    We have:
    \begin{align*}
        \int_{\Gamma \bs G}{ |f(g)|^2 \, dg }
        &\leq \int_{S_0}{ |f(g)|^2 \, dg } + \sum_{i=1}^m \int_{\Sfrak_{std, i}}{ |f(g)|^2 \, dg } \\
        &\leq l_0(f)^2 \int_G{ (1 + \log a_g)^{-2r} a_g^{-2\rho} \, dg } \\
        &\qquad + \sum_{i=1}^m l_i(f)^2 \int_{\Sfrak_{std, i}}{ a_{P_i}(g^{-1})^{-2(1-\xi_i)\rho_{P_i}} (1 + \log a_g)^{-2r} a_g^{-2 \xi_i \rho_{P_i}} \, dg } \ .
    \end{align*} 
    This is finite by Lemma \ref{Lem int_G a_g^(-2rho) (1 + |log a_g|)^(-r) dg < oo} and by Lemma \ref{Lem int_CEAK( a(g^(-1))^(-2(1-xi) rho)(1 + log a_g|)^(-2r) a_g^(-2xi rho) dg < infty} with $ E = N^{\Gamma, P_i} $ and $ A = A_i $ ($ i = 1, \dotsc, m $). The proposition follows.
\end{proof}

\newpage

\subsubsection{The space of cusp forms on \texorpdfstring{$ \Gamma \bs G $}{Gamma\textbackslash G}}

We define the constant term of a Schwartz function along a cuspidal parabolic subgroup of $G$. We show that it is well-defined and that it has some basic properties. Furthermore, we define the space of cusp forms.

\medskip

Define\index[n]{VphiGammaQ@$ V_\phi^{\Gamma_Q} $}
\[
    V_\phi^{\Gamma_Q} = \{ v \in V_\phi \mid \phi(\gamma)v = v \quad \forall \gamma \in \Gamma_Q \} \ .
\]
Let $ f \colon G \to V_\phi $ be a measurable function. Let $ p_Q $ be the orthogonal projection on $ V_\phi^{\Gamma_Q} $. Put\index[n]{f0@$ f^0 $}\index[n]{pQ@$ p_Q $}
\[
    f^0(g) = p_Q f(g) \qquad (g \in G) \ .
\]
The function $ f^0 $ is a left $ \Gamma_Q $-invariant function from $G$ to $ V_\phi^{\Gamma_Q} $.
Since $ p_Q $ has operator norm 1, $ |f^0(g)| \leq |f(g)| $ for all $ g \in G $.

As $ p_Q $ is continuous, $ L_X R_Y f^0(g) = (L_X R_Y f)^0(g) $ for all $ g \in G $ and $ X, Y \in \U(\g) $.

We will use the projection $ p_Q $ in the definition of the constant term along a cusp (cf. Definition \ref{Dfn constant term of f along Q}).

\begin{lem}\label{Lem about V_phi^(Gamma_Q)}
    Let $ Q' = \gamma Q \gamma^{-1} $ ($ \gamma \in \Gamma $). Then,
    \begin{enumerate}
    \item $ V_\phi^{\Gamma_{Q'}} = \phi(\gamma) V_\phi^{\Gamma_Q} $;
    \item the orthogonal complement of $ V_\phi^{\Gamma_{Q'}} $ is equal to $ \phi(\gamma) (V_\phi^{\Gamma_Q})^\perp $;
    \item $ p_{Q'} \circ \phi(\gamma) = \phi(\gamma) \circ p_Q $.
    \end{enumerate}
\end{lem}

\begin{proof} \nlenum
    \begin{enumerate}
    \item We have:
        \begin{align*}
            v \in V_\phi^{\Gamma_{Q'}}
            \iff & \phi(\gamma \gamma_Q \gamma^{-1})v = v \quad \forall \gamma_Q \in \Gamma_Q \\
            \iff & \phi(\gamma)^{-1} v \in V_\phi^{\Gamma_Q} \ .
        \end{align*}
        Thus, $ V_\phi^{\Gamma_{Q'}} = \phi(\gamma) V_\phi^{\Gamma_Q} $.
    \item Since $ V_\phi = V_\phi^{\Gamma_Q} \oplus (V_\phi^{\Gamma_Q})^\perp $ and since $ \phi $ is unitary, we have
        \begin{align*}
            (V_\phi^{\Gamma_{Q'}})^\perp &= \{ v \in V_\phi \mid (v, v') = 0 \quad \forall v' \in V_\phi^{\Gamma_{Q'}} \} \\
            &= \{ v \in V_\phi \mid (\phi(\gamma)^{-1}v, v') = 0 \quad \forall v' \in V_\phi^{\Gamma_Q} \} \\
            &= \phi(\gamma) \{ v \in V_\phi \mid (v, v') = 0 \quad \forall v' \in V_\phi^{\Gamma_Q} \}
            = \phi(\gamma) (V_\phi^{\Gamma_Q})^\perp
        \end{align*}
        by the previous assertion.
    \item Let $ v \in V_\phi $. Write $ v = v_1 + v_2 $ with $ v_1 \in V_\phi^{\Gamma_Q} $ and $ v_2 \in (V_\phi^{\Gamma_Q})^\perp $.
        Then, $ \phi(\gamma) v = \phi(\gamma) v_1 + \phi(\gamma) v_2 $ with $ \phi(\gamma)v_1 \in V_\phi^{\Gamma_{Q'}} $ and $ \phi(\gamma)v_2 \in (V_\phi^{\Gamma_{Q'}})^\perp $ by the two previous assertions. Thus, $ p_{Q'}(\phi(\gamma)v) = \phi(\gamma) p_Q(v) $.
    \end{enumerate}
\end{proof}

\begin{dfn}\label{Dfn constant term of f along Q}
    For $ g \in G $ and for a $ \Gamma $-equivariant measurable function $ f \colon G \to V_\phi $, we define
    \[
        f^Q(g) = \frac1{\vol(\Gamma_Q \bs N_{\Gamma, Q} M_{\Gamma, Q})} \int_{\Gamma_Q \bs N_Q M_{\Gamma, Q}}{ f^0(nmg) \, dn \, dm}
    \]
    whenever the integral exists. We call this function the \textit{constant term} of $f$ along $Q$. Recall that we defined $ f^0(x) $ ($x \in G$) by $ p_Q f(x) $, where $ p_Q $ is the orthogonal projection on $ V_\phi^{\Gamma_Q} $.\index{constant term}\index[n]{fQ@$ f^Q $}
\end{dfn}

\begin{prop}\label{Prop f^P well-defined}
    Let $ f \in \CS(\Gamma \bs G, \phi) $, $ Y \in \U(\g) $ and $ g \in G $. Then, $ nm \mapsto R_Y f^0(nmg) $ is integrable over $ \Gamma_Q \bs N_Q M_{\Gamma, Q} $. Thus, the constant term $ f^Q $ along $Q$ is well-defined (take $ Y = 1 $). Moreover, if $ \Sfrak := \Sfrak_{std, Q, t, \omega, \omega'} $ is a generalised standard Siegel set in $G$ with respect to $ Q $ and if $ r > 1 $, then there exists a constant $ C > 0 $ (depending on $ \omega $ and $r$) such that
    \begin{multline*}
        \frac1{\vol(\Gamma_Q \bs N_{\Gamma, Q} M_{\Gamma, Q})} \int_{\Gamma_Q \bs N_Q M_{\Gamma, Q}} { |R_Y f(nmg)| \, dn \, dm} \\
        \leq C \LO{p}{\Sfrak'}_{2r, 1, Y}(f) a_Q(g^{-1})^{-\rho_Q} (1 + |\log a_Q(g^{-1})|)^{-r}
    \end{multline*}
    for all $ g \in G $ with $ h_Q(g) \geq t $.
\end{prop}

\begin{proof}
    Let $ \Sfrak := \Sfrak_{std, Q, t, \omega, \omega'} $ be a generalised standard Siegel set in $G$ with respect to $ Q $ and let $ g \in G $ be as above.
    \\ Note that there is a generalised standard Siegel set containing all $ g \in G $ with $ h_Q(g) \geq t $.

    By definition, the map
    \[
        \omega M_{\Gamma, Q} \to \Gamma_Q \bs N_{\Gamma, Q} M_{\Gamma, Q} \, , \quad n m \mapsto \Gamma_Q nm
    \]
    is surjective and at most $k$-to-one for some $ k \in \NN $.
    For $ g = \nu_Q(g) h_Q(g) k_Q(g) \in G $, we have
    \[
        \int_{\Gamma_Q \bs N_Q M_{\Gamma, Q}}{ |R_Y f(nmg)| \, dn \, dm} = \int_{\Gamma_Q \bs N_Q M_{\Gamma, Q}}{ |R_Y f(n m h_Q(g) k_Q(g))| \, dn \, dm}
    \]
    as $ N_Q $ is unimodular. This divided by $ \vol(\Gamma_Q \bs N_{\Gamma, Q} M_{\Gamma, Q}) $ is less or equal than
    \[
        c \int_{\omega N^{\Gamma, Q} M_{\Gamma, Q}}{ |R_Y f(n m h_Q(g) k_Q(g))| \, dn \, dm}
    \]
    for some constant $ c > 0 $.
    \\ Let $ r > 1 $. As $ f \in \CS(\Gamma \bs G, \phi) $, $ \int_{M_{\Gamma, Q}} \int_{\omega N^{\Gamma, Q}}{ |R_Y f(n h_Q(g) m k_Q(g))| \, dn \, dm} $ is less or equal than
    \[
          \LO{p}{\Sfrak'}_{2r, 1, Y}(f)
          \int_{M_{\Gamma, Q}} \int_{\omega N^{\Gamma, Q}}{ h_Q(g)^{(1 - \xi) \rho_Q} (1 + \log a_{n h_Q(g)})^{-2r} a_{n h_Q(g)}^{-\xi \rho_Q} \, dn \, dm} \ ,
    \]
    where $ \xi := \xi_{Q, N^{\Gamma, Q}} $.
    But as $ \vol(M_{\Gamma, Q}) = 1 $, this is less or equal than
    \[
        \LO{p}{\Sfrak'}_{2r, 1, Y}(f) h_Q(g)^{(1 - \xi)\rho_Q} \int_{\omega N^{\Gamma, Q}} { (1 + \log a_{n h_Q(g)})^{-2r} a_{n h_Q(g)}^{-\xi \rho_Q} \, dn } \ .
    \]
    So, $ I := \frac1{\vol(\Gamma_Q \bs N_{\Gamma, Q} M_{\Gamma, Q})} \int_{\Gamma_Q \bs N_Q M_{\Gamma, Q}}{ |R_Y f(nmg)| \, dn \, dm} $ is finite by Fubini's theorem as the last integral is finite by Lemma \ref{Lem Estimate of int_CE a_na^(-xi rho) (1 + log a_na)^(-2d)dn}. Moreover, there exists a constant $ c' > 0 $ (depending on $ \omega $ and $ r $) such that
    \[
        \int_{\omega N^{\Gamma, Q}}{ a_{n h_Q(g)}^{-\xi \rho_Q} (1 + \log a_{n h_Q(g)})^{-2r} \, dn } \leq c' \LO{p}{\Sfrak'}_{2r, 1, Y}(f) h_Q(g)^{\xi \rho_Q} (1 + |\log h_Q(g)|)^{-r} \ .
    \]
    Hence,
    \[
        I \leq c c' \LO{p}{\Sfrak'}_{2r, 1, Y}(f) h_Q(g)^{\rho_Q} (1 + |\log h_Q(g)|)^{-r} \ .
    \]
    The proposition follows.
\end{proof}

\begin{cor}\label{Cor f^Q is continuous}
    Fix $ g \in G $. Then, $ f \in \CS(\Gamma \bs G, \phi) \mapsto f^Q(g) $ is linear and continuous.
\end{cor}

\begin{proof}
    Observe that the map is linear.
    Fix $ g \in G $ and let $ f \in \CS(\Gamma \bs G, \phi) $. Then, by Proposition \ref{Prop f^P well-defined}, we have
    \[
        |f^Q(g)| \leq C_g(f) \ ,
    \]
    where $ C_g(f) $ is a constant depending on $ g $ and linearly on $f$.
    Hence, $ f \mapsto f^Q(g) $ is continuous.
\end{proof}

\begin{cor}\label{Prop f^P depends smoothly on g}
    Let $ f \in \CS(\Gamma \bs G, \phi) $ and let $ Y \in \U(\g) $. Then, $ f^Q $ is smooth and for $ g \in G $, we have
    \[
        R_Y (f^Q)(g) = \frac1{\vol(\Gamma_Q \bs N_{\Gamma, Q} M_{\Gamma, Q})} \int_{\Gamma_Q \bs N_Q M_{\Gamma, Q}}{ R_Y f^0(nmg) \, dn \, dm} \ .
    \]
\end{cor}

\begin{proof}
    The proof is analogue to the proof of Corollary \ref{Cor f^Omega(g, h) depends smoothly on g and h}.
\end{proof}

\begin{prop}\label{Prop equivariances in the geometrically finite case}
    Let $ f \in \CS(\Gamma \bs G, \phi) $. Then, the constant of $f$ along $Q$ has the following properties:
    \begin{enumerate}
    \item $ f \mapsto f^Q $ is equivariant with respect to right translation:
        \[
            (R_x f)^Q = R_x f^Q \qquad (x \in G) \ ;
        \]
    \item $ f^Q $ is left $ N_Q M_{\Gamma, Q} $-invariant;\label{Prop equivariances in the geometrically finite case eq2}
    \item Let $ Q' = \gamma Q \gamma^{-1} $ $ (\gamma \in \Gamma) $.
        Then, $ f^{Q'} $ is well-defined and\label{Prop equivariances in the geometrically finite case eq3}
        \[
            f^{Q'}(g) = \phi(\gamma) f^Q(\gamma^{-1} g) \ .
        \]
    \end{enumerate}
\end{prop}

\needspace{2\baselineskip}
\begin{proof}
    Let $ f \in \CS(\Gamma \bs G, \phi) $.
    \begin{enumerate}
    \item For $ g \in G $, $ x \in G $, we have
        \[ (R_x f)^Q(g) = \frac1{\vol(\Gamma_Q \bs N_{\Gamma, Q} M_{\Gamma, Q})} \int_{\Gamma_Q \bs N_Q M_{\Gamma, Q}}{ f^0(nmgx) \, dn \, dm} = R_x f^Q(g) \ . \]
    \item For $ g \in G $, $ n' \in N_Q $, $ m' \in M_{\Gamma, Q} $, we have
        \begin{align*}
            f^Q(n'm'g) &= \frac1{\vol(\Gamma_Q \bs N_{\Gamma, Q} M_{\Gamma, Q})} \int_{\Gamma_Q \bs N_Q M_{\Gamma, Q}}{ f^0(nmn'm'g) \, dn \, dm} \\
            &= \frac1{\vol(\Gamma_Q \bs N_{\Gamma, Q} M_{\Gamma, Q})} \int_{\Gamma_Q \bs N_Q M_{\Gamma, Q}}{ f^0(nmg) \, dn \, dm} = f^Q(g)
        \end{align*}
        as $N_Q M_{\Gamma, Q}$ is unimodular.
    \item This part is similarly proven than the corresponding assertion in \cite[p.74]{Borel}. \\
        Let $ \gamma \in \Gamma $ and let $ Q' = \gamma Q \gamma^{-1} $.
        The unipotent radical $ N_{Q'} $ of $ Q' $ is clearly equal to $ \gamma N_Q \gamma^{-1} $, $ N_{\Gamma, Q'} = \gamma N_{\Gamma, Q} \gamma^{-1} $ and $ M_{\Gamma, Q'} = \gamma M_{\Gamma, Q} \gamma^{-1} $.

        By Remark \ref{Rem Schwartz space well-defined} and by the proof of Proposition \ref{Prop f^P well-defined}, $ f^{Q'} $ is well-defined.

        Consider now
        \[
            \Phi \colon \Gamma_Q \bs N_{\Gamma, Q} M_{\Gamma, Q} \to \Gamma_{Q'} \bs N_{\Gamma, Q'} M_{\Gamma, Q'} , \; \Gamma_Q n m \mapsto \Gamma_{Q'} \Int \gamma(n m) \ ,
        \]
        where $ \Int g(x) = gxg^{-1} $.
        This map is well-defined as for any $ \gamma_Q \in \Gamma_Q $ and $ n \in N_Q $, $ m \in M_{\Gamma, Q} $, we have
        \begin{multline*}
            \Phi(\Gamma_Q \gamma_Q n m) = \Gamma_{Q'} \Int \gamma(\gamma_Q n m) \\
            = \Gamma_{Q'} \Int \gamma(\gamma_Q) \Int \gamma(n m)
            = \Gamma_{Q'} \Int \gamma(n m) = \Phi(\Gamma_Q n m)
        \end{multline*}
        since $ \Int \gamma(\gamma_Q) \in \Gamma \cap N_{\Gamma, Q'} M_{Q'} = \Gamma_{Q'} $.
        The map $ \Phi $ is clearly an isomorphism as $ \Int \gamma $ is an automorphism transforming $ \frac1{\vol(\Gamma_Q \bs N_{\Gamma, Q} M_{\Gamma, Q})} dn dm $ into $ \frac1{\vol(\Gamma_{Q'} \bs N_{\Gamma, Q'} M_{\Gamma, Q'})} dn' dm' $. Here, $ dn dm $ (resp. $ dn' dm' $) is the Haar measure on $ N_{\Gamma, Q} M_{\Gamma, Q} $ (resp. on $ N_{\Gamma, Q'} M_{\Gamma, Q'} $).
        Since moreover $ N_Q = N_{\Gamma, Q} N^{\Gamma, Q} $,
        \begin{multline*}
            f^{Q'}(g) = \frac1{\vol(\Gamma_{Q'} \bs N_{\Gamma, Q'} M_{\Gamma, Q'})} \int_{\Gamma_{Q'} \bs N_{Q'} M_{\Gamma, Q'}}{ p_{Q'} f(n'm'g) \, dn' \, dm'} \\
            = \frac1{\vol(\Gamma_Q \bs N_{\Gamma, Q} M_{\Gamma, Q})} \cdot \int_{\Gamma_Q \bs N_Q M_{\Gamma, Q}}{ p_{Q'} f(\gamma nm \gamma^{-1} g) \, dn \, dm}
        \end{multline*}
        for all $ g \in G $. By Lemma \ref{Lem about V_phi^(Gamma_Q)}, this is equal to
        \[
            \frac1{\vol(\Gamma_Q \bs N_{\Gamma, Q} M_{\Gamma, Q})} \cdot \phi(\gamma)
            \int_{\Gamma_Q \bs N_Q M_{\Gamma, Q}}{ f^0(nm \gamma^{-1} g) \, dn \, dm}
            = \phi(\gamma) f^Q(\gamma^{-1} g) \ .
        \]
    \end{enumerate}
\end{proof}

\begin{dfn}
    For $ g \in G(\Omega), h \in G $ and for a measurable function $ f \colon G \to V_\phi $, we define
    \[
        f^\Omega(g, h) = \int_N{ f(gnh) \, dn }
    \]
    whenever the integral exists.
    We call this function the \textit{constant term} of $f$ along $ \Omega $.\index{constant term}
\end{dfn}

\begin{rem} \nlenum
    \begin{enumerate}
    \item If $ f \in \CS(\Gamma \bs G, \phi) $, then the constant term of $f$ along $ \Omega $ is well-defined by Proposition \ref{Prop f^Omega well-defined}.
    \item You can use your favorite parabolic subgroup of $G$ to define $ G(\Omega) $. If $ \Gamma \bs X $ has a cusp, then you can define $ G(\Omega) $ for example by $ \{ g \in G \mid g P_1 \in \Omega \} $.
    \end{enumerate}
\end{rem}

\begin{dfn}\label{Dfn Harish-Chandra space of cusp forms in the geometrically finite case}
    Define the following subspace of $ \CS(\Gamma \bs G, \phi) $:
    \[
        \c{\CS}(\Gamma \bs G, \phi) = \{ f \in \CS(\Gamma \bs G, \phi) \mid f^\Omega = 0 , \,
        f^Q = 0 \quad \forall [Q]_\Gamma \in \P_\Gamma \} \ . 
    \]
    We call it the \textit{space of cusp forms} on $ \Gamma \bs G $.\index{space of cusp forms}\index[n]{CGammaG0@$ \c{\CS}(\Gamma \bs G, \phi) $}
\end{dfn}

\begin{rem} \nlenum
    \begin{enumerate}
    \item This subspace may be considered as the analog of the Harish-Chandra space of cusp forms on $ G $ and is an extension of Definition \ref{Dfn Harish-Chandra space of cusp forms in the geometrically finite case}.
    \item By Proposition \ref{Prop equivariances in the geometrically finite case} \eqref{Prop equivariances in the geometrically finite case eq3}, it suffices to check the condition $ f^Q = 0 $ for $ Q = P_i $ ($ i = 1, \dotsc, m $). 
    \end{enumerate}
\end{rem}

\newpage

\subsubsection{The little constant term}

In the following, we define the little constant term, we determine its basic properties and we relate it to the constant term.

\begin{dfn}\label{Dfn little constant term}
    For $ g \in G $ and for a $ \Gamma $-equivariant measurable function $ f \colon G \to V_\phi $, we define
    \[
        f^{Q, lc}(g) = \frac1{ \vol(\Gamma_Q \bs N_{\Gamma, Q} M_{\Gamma, Q}) } \int_{\Gamma_Q \bs N_{\Gamma, Q} M_{\Gamma, Q}}{ f^0(nmg) \, dn \, dm}
    \]
    whenever the integral exists. We call this function the \textit{little constant term} of $f$ along $Q$.
    Recall that we defined $ f^0(x) $ ($x \in G$) by $ p_Q f(x) $.\index{little constant term}\index[n]{fQlc@$ f^{Q, lc} $}
\end{dfn}

\begin{rem} \nlenum
    \begin{enumerate}
    \item Recall that $ \Gamma_Q \bs N_{\Gamma, Q} M_{\Gamma, Q} $ is compact.
    \item The little constant term exists for example if $ f \colon G \to V_\phi $ is a $ \Gamma $-equivariant continuous function.
    \end{enumerate}
\end{rem}

\begin{prop}\label{Prop Little constant term properties}
    Let $ f \colon G \to V_\phi $ be a $ \Gamma $-equivariant measurable function such that $ |f|^{Q, lc} $ is finite.
    Then, we have the following:
    \begin{enumerate}
    \item $ f \mapsto f^{Q, lc} $ is equivariant with respect to right translation:
        \[
            (R_x f)^{Q, lc} = R_x f^{Q, lc} \qquad (x \in G) \ ;
        \]
    \item $ f^{Q, lc} $ is left $ N_{\Gamma, Q} M_{\Gamma, Q} $-invariant;\label{eq f^(Q, lc) N_Gamma-invariant}
    \item Let $ Q' = \gamma Q \gamma^{-1} $ $ (\gamma \in \Gamma) $.
        Then, $ |f|^{Q', lc} $ is finite and
        \[
            f^{Q', lc}(g) = \phi(\gamma) f^{Q, lc}(\gamma^{-1} g) \qquad (g \in G) \ .
        \]
    \item\label{eq (f^(Q, lc))^(Q, lc)}
        $
            (f^{Q, lc})^{Q, lc} = f^{Q, lc}
        $;
    \item $ (f * \phi)^{Q, lc} = f^{Q, lc} * \phi $ for all $ \phi \in \C_c(G) $;
    \item $ R_Y (f^{Q, lc}) = (R_Y f)^{Q, lc} $ for all $ f \in \Cinf(\Gamma \bs G, \phi) $ and $ Y \in \U(\g) $.
    \end{enumerate}
\end{prop}

\begin{proof}
    The first three assertions can be proven similarly as the corresponding assertions in Proposition \ref{Prop equivariances in the geometrically finite case}.
    The remaining assertions can be easily checked by direct computation.
\end{proof}

\begin{lem}\label{Lemma omega p_(N^Gamma)(n(g)) h(g) K in Sfrak'}
    Let $ \Sfrak := \Sfrak_{std, Q, t, \omega, \omega'} $ be a generalised standard Siegel set in $G$ with respect to $ Q $. Then,
    \[
        n m p_{N^\Gamma, Q}(\nu_Q(g)) h_Q(g) K_Q \in \Sfrak'
    \]
    for all $ n \in \omega $ and $ m \in M_{\Gamma, Q} $.
\end{lem}

\begin{proof}
    Let $ n \in \omega $ and $ m \in M_{\Gamma, Q} $. Let $ \Sfrak $ be as above. Then,
    \begin{equation}\label{eq omega p_(N^Gamma)(n(g)) h(g) K in Sfrak'}
        n m p_{N^\Gamma, Q}(\nu_Q(g)) h_Q(g) K_Q \in \Sfrak'
    \end{equation}
    for all $ g \in \Sfrak $ with $ p_{N^\Gamma, Q}(\nu_Q(g)) \in N^{\Gamma, Q} \smallsetminus \omega' $ as $ \omega' $ is invariant under conjugation in $M_{\Gamma, Q}$. For $ g \in \Sfrak $ with $ p_{N^\Gamma, Q}(\nu_Q(g)) \in \omega' $, $ h_Q(g) \geq t $. Hence, \eqref{eq omega p_(N^Gamma)(n(g)) h(g) K in Sfrak'} holds also for $ g \in \Sfrak $ with $ p_{N^\Gamma, Q}(\nu_Q(g)) \in \omega' $.
\end{proof}

\begin{dfn}\label{Dfn Schwartz function near eQ}
    We say that $ f \in \Cinf(\Gamma \bs G, \phi) $ is a \textit{Schwartz function near $eQ$} if there exists $ s' \geq 0 $ such that\index{Schwartz function near a cusp}
    \[
        \LO{p}{\Sfrak'}_{r, X, Y}(f)
    \]
    is finite for every $ r \geq 0 $, $ X \in \U(\n_Q), Y \in \U(\g) $ and $ \Sfrak := \Sfrak_{std, Q, t, \omega, \omega'} \in \V_{std, \Gamma_Q} $ with $ \omega' $ containing $ \{ n \in N^{\Gamma, Q} \mid |\log(n)| < s' \} $.
    We call $ \Sfrak $ as above \textit{admissible} (for $f$).\index{admissible generalised standard Siegel set}
\end{dfn}

\begin{rem}
    Every element of $ \CS(\Gamma_Q \bs G, \phi) $ is clearly a Schwartz function near $eQ$.
\end{rem}

\begin{lem}\label{Lem estimate of the little constant term}
    Let $ f $ be a Schwartz function near $eQ$, $ X \in \U(\n_Q) $ and $ Y \in \U(\g) $.
    \\ If $ \Sfrak := \Sfrak_{std, Q, t, \omega, \omega'} $ is an admissible generalised standard Siegel set in $G$ with respect to $ Q $ and if $ r \geq 0 $, then there exists a constant $ C > 0 $ (depending on $ \omega $ and $ r $) such that
    \begin{multline}\label{eq Lem estimate of the little constant term}
        |L_X R_Y f^{Q, lc}(g)| \\
        \leq C \LO{p}{\Sfrak'}_{r, X, Y}(f)
        h_Q(g)^{\rho_{\Gamma_Q}} (1 + \log a_{p_{N^{\Gamma, Q}}(\nu_Q(g)) h_Q(g)})^{-r} a_{p_{N^{\Gamma, Q}}(\nu_Q(g)) h_Q(g)}^{-\rho^{\Gamma_Q}}
    \end{multline}
    for all $ g \in \Sfrak $. Thus, $ f^{Q, lc} $ is a Schwartz function near $eQ$. 
\end{lem}

\begin{proof}
    Let $f$ and $ \Sfrak $ be as above, $ X \in \U(\n_Q) $, $ Y \in \U(\g) $ and $ g \in \Sfrak $.
    Then,
    \[
        \int_{\Gamma_Q \bs N_{\Gamma, Q} M_{\Gamma, Q}}{ | L_X R_Y (f (nm \, \cdot))(g)| \, dn \, dm}
    \]
    is equal to
    \[
        \int_{\Gamma_Q \bs N_{\Gamma, Q} M_{\Gamma, Q}}{ | L_{\Ad(p_{N_{\Gamma, Q}}(\nu_Q(g)))^{-1}X} R_Y (f (nm \ \cdot))(p_{N^{\Gamma, Q}}(\nu_Q(g)) h_Q(g) k_Q(g))| \, dn \, dm}
    \]
    as $ \Gamma_Q \bs N_{\Gamma, Q} M_{\Gamma, Q} $ is unimodular. Then, this is less or equal than
    \[
        \int_{M_{\Gamma, Q}} \int_{\omega}{ | L_{\Ad(nm) \Ad(p_{N_{\Gamma, Q}}(\nu_Q(g)))^{-1}X} R_Y f (nm p_{N^{\Gamma, Q}}(\nu_Q(g)) h_Q(g) k_Q(g))| \, dn \, dm} \ .
    \]
    Since we can write $ \Ad(n m) X $ ($ n \in \omega, m \in M_{\Gamma, Q} $) as finite linear combination of elements in $ \U(\n_Q) $ with bounded coefficients, it suffices to show that
    \[
        \int_{M_{\Gamma, Q}} \int_{\omega}{ | L_X R_Y f (nm p_{N^{\Gamma, Q}}(\nu_Q(g)) h_Q(g) k_Q(g))| \, dn \, dm}
    \]
    has the required estimate.
    By Lemma \ref{Lemma omega p_(N^Gamma)(n(g)) h(g) K in Sfrak'} and as $ f $ is a Schwartz function near $eQ$, there exists $ c > 0 $ such that this divided by $ \vol(\Gamma_Q \bs N_{\Gamma, Q} M_{\Gamma, Q}) $ is less or equal than
    \begin{multline*}
          c \LO{p}{\Sfrak'}_{r, X, Y}(f)
          \int_{M_{\Gamma, Q}} \int_{\omega} h_Q(g)^{\rho_{\Gamma_Q}} (1 + \log a_{n m p_{N^{\Gamma, Q}}(\nu_Q(g)) h_Q(g)})^{-r} \\
          a_{n m p_{N^{\Gamma, Q}}(\nu_Q(g)) h_Q(g)}^{-\rho^{\Gamma_Q}} \, dn \, dm \ .
    \end{multline*}
    But as $ \vol(M_{\Gamma, Q}) = 1 $, this is less or equal than
    \[
        c' \vol(\omega) \LO{p}{\Sfrak'}_{r, X, Y}(f) h_Q(g)^{\rho_{\Gamma_Q}} (1 + \log a_{p_{N^{\Gamma, Q}}(\nu_Q(g)) h_Q(g)})^{-r} a_{p_{N^{\Gamma, Q}}(\nu_Q(g)) h_Q(g)}^{-\rho^{\Gamma_Q}}
    \]
    for some constant $ c' > 0 $ as $ \sup_{n \in \omega} a_n $ is finite.
    \eqref{eq Lem estimate of the little constant term} follows.
    Thus, $ f^{Q, lc} $ is a Schwartz function near $eQ$.
\end{proof}

\begin{prop}
    Let $ f \in \CS(\Gamma \bs G, \phi) $. Then, $ f^Q = 0 $ if $ f^{Q, lc} $ vanishes.
\end{prop}

\begin{proof}
    Let $ f \in \CS(\Gamma \bs G, \phi) $ be such that $ f^{Q, lc} = 0 $ and $ g \in G $. By Fubini's theorem applied to $ R_g f $ which is integrable over $ \Gamma_Q \bs N_Q M_{\Gamma, Q} $ by Proposition \ref{Prop f^P well-defined}, we have
    \begin{multline*}
        f^Q(g) = \frac1{\vol(\Gamma_Q \bs N_{\Gamma, Q} M_{\Gamma, Q})} \int_{\Gamma_Q \bs N_Q M_{\Gamma, Q}}{ f(n m g) \, dn \, dm } \\
        = \frac1{\vol(\Gamma_Q \bs N_{\Gamma, Q} M_{\Gamma, Q})} \int_{N_{\Gamma, Q} \bs N_Q}{ \underbrace{\int_{\Gamma_Q \bs N_{\Gamma, Q} M_{\Gamma, Q}}{ f(n m n' g) \, dn \, dm } }_{= 0} \, dn' } = 0 \ .
    \end{multline*}
\end{proof}

\newpage

\subsubsection{\texorpdfstring{$ \SU(1, 2) $}{SU(1, 2)}-reduction}

In the following, we generalise some formulas obtained in Appendix \ref{sec:Sp(1,n)} for $ G = \Sp(1, n) $ by using $ \SU(1, 2) $-reduction.

For $ t > 0 $, set $ a_t = \exp(\log(t) H_\alpha) $.
This is a generalisation of the $ a_t $ defined in Appendix \ref{sec:Sp(1,n)}.

For $ X \in \g_\alpha $, $ Y \in \g_{2\alpha} $ and $ r \in \RR $, set $ X \cdot r = r X $ and $ X \cdot Y = -\frac12[\theta X, Y] \in \g_\alpha $.

Note that $ X \cdot Y $ is simply the natural multiplication induced from the multiplication on $ \HH $ if $ G = \Sp(1, n) $ ($ n \geq 2 $).

\begin{thm}\label{Thm explicit formulas}
    We have:
    \begin{enumerate}
    \item Write $ \bar{n} \in \bar{N} $ as $ \bar{n} = \exp(X + Y) \ ( X \in \g_{-\alpha} $, $ Y \in \g_{-2\alpha} ) $. Then,
        \begin{enumerate}
        \item $ \nu(\bar{n}) = \exp\bigl(- X \cdot \frac{ 1 + \frac12|X|^2 - 2 Y }{(1 + \frac12|X|^2)^2 + 2|Y|^2} - \frac{ Y }{(1 + \frac12|X|^2)^2 + 2|Y|^2} \bigr) $, and
        \item $ h(\bar{n}) = \frac1{\sqrt{(1 + \frac12|X|^2)^2 + 2 |Y|^2}} \ $.
        \end{enumerate}
    \item Let $ n \in N $ and $ a \in A $. Write $ n \in N $ as $ n = \exp(X + Y) $ ($ X \in \g_{\alpha} $, $ Y \in \g_{2\alpha} $) and $ a = a_t $ ($ t > 0 $). Then,
        \[
            a_{n a} = e^{ \arccosh\Bigl(\sqrt{ \big(\tfrac12(t + \tfrac1t) + \tfrac1{4t} |X|^2\big)^2 + \tfrac{|Y|^2}{2t^2}}\Bigr) } \ .
        \]
    \end{enumerate}
\end{thm}

\begin{rem} \nlenum
    \begin{enumerate}
    \item The second formula was already established by S.~Helgason and the third formula is a generalisation of a formula obtained by S.~Helgason (cf. Proposition \ref{Theorem 3.8 of Helgason, p.414}).
    \item One can show that the formula for $ a_n $ ($ n \in N $) is indeed equal to the one provided by Proposition \ref{Theorem 3.8 of Helgason, p.414}.
    \end{enumerate}
\end{rem}

\begin{proof}
    The subgroups $ N $, $ A $, $M$ and $K$ of $ \Sp(1, n) $, defined in Appendix \ref{sec:Sp(1,n)}, induce Lie groups $ N_0 $, $ A_0 $, $ M_0 $ and $ K_0 $ of $ \SU(1, 2) $.
    Let $ \n_0 $, $ \a_0 $, $ \m_0 $, $ \k_0 $ be the Lie algebras of $ N_0 $, $ A_0 $, $ M_0 $ respectively $ K_0 $.
    Let $ \theta_0 $ the Cartan involution associated to $ \k_0 $. The Lie algebra $ \n_0 $ is given by $ \{ X_{z, it} \mid z \in \CC , t \in \RR \} $ ($ X_{v, r} $ is defined in Appendix \ref{sec:Sp(1,n)}).
    So, $ \n_0 = \CC X_{1, 0} \oplus \RR X_{0, i} $.
    We denote $ \nu_{\theta_0, P_0} $, $ h_{\theta_0, P_0} $, $ k_{\theta_0, P_0} $ simply by $ \nu_0 $, $ h_0 $, $ k_0 \ $.

    Let $ Z^* = Z^*_1 + Z^*_2 \in \n $ ($ Z^*_i \in \g_{i \alpha}, Z^*_1 \neq 0 $). If $ Z^*_i \neq 0 $, set $ X_{i\alpha} = \frac1{|Z^*_i|} \cdot Z^*_i $. Otherwise, set $ X_{i\alpha} = 0 $.
    By construction, $ Z^* \in \RR X_\alpha \oplus \RR X_{2\alpha} $.

    Let $ a \in A $ and let $ H^* \in \a $ be such that $ a = \exp(H^*) $.

    \smallskip

    We distinguish the following two cases:
    \begin{enumerate}
    \item \underline{$ X_{2\alpha} \neq 0 $:} Let $ \l = \su(1, 2) $ and let $ L = \SU(1, 2) $. Then, it follows from Theorem 3.1, Chapter IX, of \cite[p.409]{Helgason} ($ \SU(2, 1) $-reduction) that the Lie subalgebra $ \g^* \subset \g $ generated by $ X_\alpha $, $ X_{2\alpha} $, $ \theta(X_\alpha) $, $ \theta(X_{2 \alpha}) $ is isomorphic to $ \l $ via a map $ \Phi \colon \l \to \g^* $ (The isomorphism from $ \g^* $ to $ \su(2, 1) $ is explicitly described in \cite[p.413]{Helgason}. Note that it is norm preserving.) satisfying $ \Phi(\log(a_e)) = H_\alpha $ ($ \alpha(H_\alpha) = 1 $), $ \Phi(X_{1, 0}) = \sqrt{2} X_\alpha $, $ \Phi(X_{0, i}) = \sqrt{2} X_{2\alpha} $ and $ \Phi \circ \theta_0 = \theta \circ \Phi $.
    \item \underline{$ X_{2\alpha} = 0 $:}
        Let $ \l = \so(1, 2) \subset \su(1, 2) $ and let $ L = \SU(1, 2) $. Then, $\l$ is the Lie subalgebra of $ \su(1, 2) $ which is generated by $ X_{1, 0} $ and $ \theta(X_{1, 0}) $ ($\l$ is isomorphic to $ \su(1, 1) $).
        The Lie subalgebra $ \g^* \subset \g $ generated by $ X_\alpha, \theta(X_\alpha) $ is isomorphic to $ \l $ via a map $ \Phi \colon \l \to \g^* $ satisfying $ \Phi(X_{1, 0}) = \sqrt{2} X_\alpha $, $ \Phi(\log(a_e)) = H_\alpha $ and $ \Phi \circ \theta_0 = \theta \circ \Phi $.
    \end{enumerate}
    Let us show in the following that the isomorphism $ \Phi $ is norm-preserving when we put on $ \g^* $ the norm $ |\cdot| $ of $\g$ restricted to $ \g^* $.

    Let $ G^* $ be the analytic subgroup of $G$ with Lie algebra $ \g^* $ and let $ B^* $ be the Killing form of $ \g^* $. Then, $ B^* $ is up to a positive constant equal to the restriction of the Killing form on $ \g $ to $ \g^* \times \g^* $.
    Since the map
    \[
        \g^* \times \g^* \to \RR , \quad (X, Y) \mapsto B(\Phi^{-1}(X), \Phi^{-1}(Y)) := \tr(\ad_\g\bigl(\Phi^{-1}(X)\bigr) \ad_\g\bigl(\Phi^{-1}(Y)\bigr))
    \]
    is an $ \Ad(G^*) $-invariant bilinear form on $ \g^* $ and as
    \begin{multline*}
        B(\Phi^{-1}(H_\alpha), \Phi^{-1}(H_\alpha)) = B(\log(a_e), \log(a_e)) = \tfrac1{c_{d_X}} \tilde{B}(\log(a_e), \log(a_e)) = \tfrac1{c_{d_X}} \\
        = \tfrac1{c_{d_X}} \alpha(H_\alpha) = \tfrac1{c_{d_X}} \tilde{B}(H_\alpha, H_\alpha) = B(H_\alpha, H_\alpha) \ ,
    \end{multline*}
    this map is equal to the Killing form on $ \g $ restricted to $ \g^* \times \g^* $. Thus, $ |\Phi(X)| = |X| $ for all $ X \in \l $.

    Let $ K^* = G^* \cap K $, $ A^* = G^* \cap A = A $, $ M^* = G^* \cap M $, $ N^* = G^* \cap N $ and $ \bar{N}^* = G^* \cap \bar{N} $. Then, it follows from Lemma 3.7, Chapter IX, of \cite[p.413]{Helgason} and the construction of his Lie groups $ K^* $, $ A^* $, $ M^* $ and $ N^* $ that $ P^* := M^* A^* N^* \subset P $ is a parabolic subgroup of $ G^* $ and that $ G^* = N^* A^* K^* $ is an Iwasawa decomposition of $ G^* $ and that $ G^* = K^* A_+^* K^* $ is a Cartan decomposition of $ G^* $.

    Moreover, $ g \in G^* $ is equal to $ \nu(g) h(g) k(g) $ (resp. $ k_g a_g h_g $) with $ \nu(g) \in N^* $, $ h(g) \in A^* $ and $ k(g) \in K^* $ (resp. $ k_g, h_g \in K^* $, $ a_g \in A^* $).

    The Lie algebra homomorphism $ \Phi \colon \l \to \g^* \hookrightarrow \g $ induces a Lie group isomorphism $ \tilde{\Phi} $ from the universal cover $ \tilde{L} $ of $ L $ to $ G $ such that $ \tilde{\Phi}_* = \Phi $. Since $ N_0 $, $ \bar{N}_0 := \theta_0(N_0) $ and $ A_0 $ are simply connected,
    \[
        \tilde{\Phi}(\exp(Z)) = \exp(\tilde{\Phi}_{*}(Z)) = \exp(\Phi(Z)) \in N                     \qquad (Z \in \n_0) \ ,
    \]
    \[
        \tilde{\Phi}(\exp(Z)) = \exp(\tilde{\Phi}_{*}(Z)) = \exp(\Phi(Z)) \in \bar{N}               \qquad (Z \in \bar{\n}_0 := \theta_0(\n_0)) \ ,
    \]
    and
    \[
        \tilde{\Phi}(\exp(H)) = \exp\big(\tilde{\Phi}_{*}(H)\big) = \exp(\Phi(H)) \in A             \qquad (H \in \a_0) \ .
    \]
    Let $ \tilde{K}_0 $ be the universal cover of $ K_0 $. Then, $ \tilde{\Phi}(\tilde{K}_0) = K $. Thus, $ N_0 A_0 \tilde{K}_0 $ is an Iwasawa decomposition of $\tilde{L}$ and $ \tilde{K}_0 (A_0)_+ \tilde{K}_0 $ is a Cartan decomposition of $\tilde{L}$. The Lie group homomorphism $ \tilde{\Phi} $ respects these Iwasawa decompositions and these $ K A K $-decompositions.

    Thus, $ \nu(\bar{n}) = \tilde{\Phi}\big(\nu_0(\tilde{\Phi}^{-1}(\bar{n}))\big) $, $ h(\bar{n}) = \tilde{\Phi}\big(h_0(\tilde{\Phi}^{-1}(\bar{n}))\big) $ for all $ \bar{n} \in \bar{N}^* $ and
    \[
        a_{na} = \tilde{\Phi}\big(a_{\tilde{\Phi}^{-1}(n) \tilde{\Phi}^{-1}(a)}\big)
    \]
    for all $ n \in N^* $ and $ a \in A^* $.

    There are $ z \in \CC $, $ u \in \RR $ and $ t > 0 $ such that
    $
        Z := \Phi^{-1}(Z^*) = z X_{1, 0} + u X_{0, i} = X_{z, i u}
    $
    ($ Z^* \in \n^* $) and $ H^* = \log(t) H_\alpha $.

    Let $ X = z X_{1, 0} $, $ Y = u X_{0, i} $ and $ H = \log(t) \log(a_e) $. Then, $ |H| = |\log(t)| \cdot |\log(a_e)| = |\log(t)| $.
    Since
    \[
        |X_{z, i u}| = \tfrac{\sqrt{2}}2 \|X_{z, i u}\| = \sqrt{2} \sqrt{|z|^2 + u^2} \ ,
    \]
    \[
        |X| = |z| \cdot |X_{1, 0}| = \sqrt{2} |z| \quad \text{and} \quad
        |Y| = |u| \cdot |X_{0, i}| = \sqrt{2} |u| \ .
    \]
    As $ h_0(\bar{n}_{z, iu}) = a_{ \frac1{\sqrt{(1 + |z|^2)^2 + 4u^2} } } $, $ \nu_0(\bar{n}_{z, i u}) = n_{(-z \cdot \frac{ 1 + |z|^2 - 2 i u }{(1 + |z|^2)^2 + 4u^2 }, -\frac{ i u }{(1 + |z|^2)^2 + 4u^2 }) } $
    and as
    \[
        a_{n_{z, iu} a_t} = e^{ \arccosh\bigl(\sqrt{ \big(\tfrac12(t + \tfrac1t) + \tfrac1{2t} |z|^2\big)^2 + \tfrac{u^2}{t^2}}\bigr) } \ ,
    \]
    the theorem follows.
\end{proof}

\newpage

\subsubsection{Two geometric estimates}

Let us prove now the following statement:
\\ \textit{There exist $ c, c' > 0 $ such that
\[
    a_{n_1 n_2 a} \geq c \, a_{n_1 a} , \qquad a_{n_1 n_2 a} \geq c' a_{n_2 a}
\]
for all $ n_1 \in N_{\Gamma, Q} $, $ n_2 \in N^{\Gamma, Q} $ and $ a \in A_Q $.}

\smallskip

We need the second estimate in order to be able to show that the right translation induces also in the geometrically finite case a representation on the Schwartz space.

We start by doing some preparations.

\begin{lem}\label{Lem about asymp}
    Let $X$ be a set. Let $ f $, $ g $ be two real-valued functions on $X$.

    \begin{enumerate}
    \item If $ f $ is nonnegative (resp. positive) and if $ f \prec g $, then $ g $ is nonnegative (resp. positive).
    \item If $ f \asymp g $ ($ \asymp $ was defined in Definition \ref{Dfn asymp}), then $ f $ is nonnegative (resp. positive) if and only if $ g $ is nonnegative (resp. positive).
    \item If $f$ and $g$ are nonnegative and if $ f \asymp \tilde{f} $ and $ g \asymp \tilde{g} $ for some real-valued functions $ \tilde{f} $, $ \tilde{g} $ on $X$, then
        $
            f + g \asymp \tilde{f} + \tilde{g} .
        $
    \item If $f$ and $g$ are positive, then $ f \asymp g $ implies that
        $
            \frac1f \asymp \frac1g .
        $
    \end{enumerate}
\end{lem}

\begin{proof}
    Obvious.
\end{proof}

\begin{lem}\label{Lem sup product over sum less or equal than a power of eps}
    Let $ N \geq 2 $, $ \eps > 0 $, $ p_i \geq 0 $ and $ q_i > 0 $ $ (i \in \{1, \dotsc, N\}) $ be such that $ \sum_{k=1}^N \frac{p_k}{q_k} \leq 1 $. Let $ D_\eps = \{ x \in \RR^N \mid \sum_{i=1}^N x_i^{q_i} \geq \eps^2 \} $ and let $ D'_\eps = D_\eps \cap \{ x \in \RR^N \mid x_i \geq 0 \quad \forall i \} $.
    Then,
    \[
        \sup_{x \in D'_\eps} \frac{ x_1^{p_1} \cdots x_N^{p_N} }{ \sum_{i=1}^N x_i^{q_i} }
    \]
    is less or equal than $ \eps^{2 (\sum_{k=1}^N \frac{p_k}{q_k} - 1)} $.
\end{lem}

\begin{proof}
    Let $ N \geq 2 $ and $ \eps > 0 $. Let $ p_i $, $ q_i $ and $ D'_\eps $ be as above.
    \\ Let $ \tilde{N} = \#\{ i \in \{1, \dotsc, N \} \mid p_i \neq 0 \} $

    We may assume without loss of generality that all the $ q_i $'s are equal to 2. Indeed, if this is not the case set $ q'_i = 2 $ and $ p'_i = \frac{2p_i}{q_i} $.
    
    By reordering the variables if necessary, we may assume without loss of generality that $ p_i = 0 $ for all $ i > \tilde{N} $ and $ p_i \neq 0 $ for all $ i \leq \tilde{N} $.

    By the change of variables $ (x_1, \dotsc, x_N) = r \Phi_N(\varphi, \psi_1, \dotsc, \psi_{N-2}) $ (spherical coordinates),
    \[
        \sup_{x \in D'_\eps} \frac{ x_1^{p_1} \cdots x_{\tilde{N}}^{p_{\tilde{N}}} }{ \sum_{i=1}^N x_i^2 }
    \]
    is equal to
    \[
        \sup_{r \in [\eps, \infty), (\varphi, \psi_j) \in E } \frac{ ( r \Phi_N(\varphi, \psi_1, \dotsc, \psi_{N-2})_1 )^{p_1} \cdots (r \Phi_N(\varphi, \psi_1, \dotsc, \psi_{N-2})_{\tilde{N}})^{p_{\tilde{N}}} }{ r^2 } ,
    \]
    where $ E = [0, \pi] \times [0, \tfrac{\pi}2]^{N-2} $.
    As $ \sum_{k=1}^N p_k \leq 2 $ and as $ \Phi_N(\varphi, \psi_1, \dotsc, \psi_{N-2})_i \leq 1 $ for all $i$, this is less or equal than
    \[
        \eps^{\sum_{k=1}^N p_k - 2}
    \]
    The lemma follows.
\end{proof}

Next we give a criterion for polynomials in $N$ variables to have the same growth behaviour.

\begin{cor}\label{Cor Estimate for Schwartz space = representation space}
    Let $ c > 0 $, $ N \in \NN $ and $ N' \in \NN $.
    \begin{enumerate}
    \item Let $ q_i > 0 $ and let $ p_{j, i} \geq 0 $ $ (i \in \{1, \dotsc, N\}, j \in \{1, \dotsc, N'\}) $ be such that
        $
            \sum_{k=1}^N \frac{p_{j, k}}{q_k} \leq 1
        $
        for all $j$. Let $ D' = \{ x \in \RR^N \mid x_i \geq 0, \ \sum_{i=1}^N x_i^{q_i} \geq 1 \} $.
        \\ If $ c_j \in (0, +\infty) $, then
        \[
            \big( c \sum_{i=1}^N x_i^{q_i} + \sum_{j=1}^{N'} c_j \ x_1^{p_{j, 1}} \cdots x_N^{p_{j, N}} \big)
            \asymp \sum_{i=1}^N x_i^{q_i} \qquad ( x \in D' ) \ .
        \]
    \item Let $ q_i \in 2 \NN $ and let $ p_{j, i} \in \NN_0 $ $ (i \in \{1, \dotsc, N\}, j \in \{1, \dotsc, N'\}) $ be such that
        $
            \sum_{k=1}^N \frac{p_{j, k}}{q_k} \leq 1
        $
        for all $j$. Let $ D_\eps = \{ x \in \RR^N \mid \sum_{i=1}^N x_i^{q_i} \geq \eps^2 \} $.
        Let
        \[
            S = \{ j \in \{1, \dotsc, N'\} \mid (\exists \, k : p_{j, k} \text{ is odd or } c_j < 0) \text{ and } \sum_{k=1}^N \frac{p_{j, k}}{q_k} = 1 \} \ .
        \]
        If $ c_j \in \RR $ satisfies $ |c_j| < \frac{ c }{ \# S } $ for all $ j \in S $, then there exists $ \eps \geq 1 $ such that
        \[
            \big( c \sum_{i=1}^N x_i^{q_i} + \sum_{j=1}^{N'} c_j \ x_1^{p_{j, 1}} \cdots x_N^{p_{j, N}} \big)
            \asymp \sum_{i=1}^N x_i^{q_i} \qquad (x \in D_\eps) \ .
        \]
        Moreover, we can choose $ \eps = 1 $ if for all $j$ either $ \sum_{k=1}^N \frac{p_{j, k}}{q_k} = 1 $ or ($ p_{j, k} $ is even for all $k$ and $ c_j \geq 0 $).
    \end{enumerate}
\end{cor}

\begin{proof}
    Let $ N \in \NN $, $ N' \in \NN $, $ c > 0 $ and $ \eps > 0 $. Let $ p_{j, i} $, $ q_i $, $m_j$, $c_j$, $ D_\eps $, $ D' $ and $S$ be as above and let $ x \in D_\eps $.
    We have
    \[
        \sum_{i=1}^N x_i^{q_i} + \sum_{j=1}^{N'} \frac{c_j}{c} \ x_1^{p_{j, 1}} \cdots x_N^{p_{j, N}}
        = \big(\sum_{i=1}^N x_i^{q_i} \big) \big(1 + \sum_{j=1}^{N'} \frac{c_j}{c} \cdot \frac{ \ x_1^{p_{j, 1}} \cdots x_N^{p_{j, N}} }{ \sum_{i=1}^N x_i^{q_i} } \big) \ .
    \]
    If $ x $ vary in $ D' $, then it follows from the previous lemma that
    \[
        \frac{c_j}{c} \cdot \frac{ \ x_1^{p_{j, 1}} \cdots x_N^{p_{j, N}} }{ \sum_{i=1}^N x_i^{q_i} }
    \]
    belongs to a compact subset of $ [0, +\infty) $ for all $ j $. The corollary follows in this case.

    Let us consider now the other case. We have the following subcases.

    \begin{enumerate}
    \item \underline{Case:} $ j \in S $ \\
    It follows from the previous lemma that
    \[
        \frac{c_j}{c} \cdot \frac{ \ x_1^{p_{j, 1}} \cdots x_N^{p_{j, N}} }{ \sum_{i=1}^N x_i^{q_i} } \qquad (x \in D_1)
    \]
    belongs to a compact subset of $ (-\tfrac1{ \# S}, \tfrac1{\# S}) $.
    \item \underline{Case:} $ j \in \{1, \dotsc, N'\} $ such that $ p_{j, k} $ are all even and $ c_j \geq 0 $ \\
    It follows from the previous lemma that
    \[
        \frac{c_j}{c} \cdot \frac{ \ x_1^{p_{j, 1}} \cdots x_N^{p_{j, N}} }{ \sum_{i=1}^N x_i^{q_i} } \qquad (x \in D_1)
    \]
    belongs to a compact subset of $ [0, +\infty) $.
    \item \underline{Case:} $ j \in \{1, \dotsc, N'\} $ such that $ \sum_{k=1}^N \frac{p_{j, k}}{q_k} < 1 $ \\
    Let $ \delta_j > 0 $. It follows from the previous lemma that
    \[
        \frac{c_j}{c} \cdot \frac{ \ x_1^{p_{j, 1}} \cdots x_N^{p_{j, N}} }{ \sum_{i=1}^N x_i^{q_i} } \qquad (x \in D_\eps)
    \]
    belongs to $ [-\delta_j, \delta_j] $ for $ \eps \geq 1 $ sufficiently large.
    \end{enumerate}
    The corollary follows by combining the different cases and by choosing the $ \delta_j $'s sufficiently small.
\end{proof}

Assume that $P$ is a $ \Gamma $-cuspidal parabolic subgroup of $G$.

Let $ \n_{2, 2}^{\Gamma, P} $ be the orthogonal complement of $ \n_{2, 1}^{\Gamma, P} := p_{\g_{2\alpha_P}}(\n_{\Gamma, P}) \cap \n_2^{\Gamma, P} $ in $ \n_2^{\Gamma, P} $.

It follows from the above that $ p_{\g_{2\alpha_P}}(\n_{\Gamma, P}) $ is equal to $ \n_{2, 1}^{\Gamma, P} \oplus (\g_{2\alpha_P} \cap \n_{\Gamma, P}) $.

We denote by $p$ the canonical projection from $ p_{\g_{2\alpha_P}}(\n_{\Gamma, P}) $ to $ \g_{2\alpha_P} \cap \n_{\Gamma, P} $.

Let $ \{ Z_1, Z_2, \dotsc, Z_{m_{\alpha_P} + m_{2\alpha_P}}\} $ be an orthonormal basis of
\[
    \n = \g_{\alpha_P} \oplus \big((\g_{2\alpha_P} \cap \n_{\Gamma, P}) \oplus \n_{2, 1}^{\Gamma, P} \oplus \n_{2, 2}^{\Gamma, P}\big) \ ,
\]
which respects this decomposition.
For $ v \in \RR^{m_{\alpha_P}} $ and $ r \in \RR^{m_{2\alpha_P}} $, set
\[
    X_{v, r} = \sum_{j=1}^{m_{\alpha_P}} v_j Z_j + \sum_{j=1}^{m_{2\alpha_P}} r_j Z_{m_{\alpha_P} + j}
\]
and $ n_{v, r} = \exp(X_{v, r}) $.
Let $ \iota \colon \RR^{m_{\alpha_P}} \times \RR^{m_{2\alpha_P}} \to \n , \; (v, r) \mapsto X_{v, r} \ $.
Then, $ \iota $ is an isomorphism of vector spaces. For $ v, v' \in \RR^{m_{\alpha_P}} $ and $ r, r' \in \RR^{m_{2\alpha_P}} $, set
\[
    [(v, r), (v', r')] = \iota^{-1}([X_{v, r}, X_{v', r'}]) \ .
\]
To simplify notations, we don't write the isomorphism $ \iota $ in the following.

\begin{lem}\label{Lem a_(n_(v_1, phi(v_1) + r_1) n_(v_2, r_2) a) asymp a_(n_(v_1 + v_2, r_1 + r_2) a)}
    Under the above assumptions, there exists a compact set $ C $ in $ \n $ such that
    \[
        a_{n_{v, r} n_{v', r'} a} \asymp a_{n_{v + v', p(r) + r'} a}
    \]
    for $ a $ varying in $ A $ and $ n_{v, r} $ (resp. $ n_{v', r'} $) varying in $ N_{\Gamma, P} $ (resp. $ N^{\Gamma, P} $) such that $ \log(n_{v, r}) + \log(n_{v', r'}) \not \in C $.
\end{lem}

\begin{proof}
    If $eP$ has full rank, then the assertion of the lemma is trivial. So, we may assume without loss of generality that $eP$ is of smaller rank.

    Let $ n_1 = n_{v, r} \in N_{\Gamma, P} $, $ n_2 = n_{v', r'} \in N^{\Gamma, P} $ and $ a = a_t \in A $. Then
    \[
        n_1 n_2 = n_{v, r} n_{v', r'} = n_{v + v', r + r' + \tfrac12 [(v, 0), (v', 0)]}
    \]
    by Lemma \ref{Lem using the Baker-Campbell-Hausdorff formula}. For $ x \geq 1 $, we have
    \[
        e^{\arccosh(x)} = x + \sqrt{x^2 - 1} = x \cdot \Big(1 + \sqrt{1 - \frac1{x^2}}\Big) \asymp x
    \]
    as $ \big\{ 1 + \sqrt{1 - \frac1{x^2}} \mid x \geq 1 \big\} $ is bounded by positive constants.

    Let $ \|v\| = \frac{\sqrt{2}}2 |X_{v, 0}| $ and $ |r| = \frac{\sqrt{2}}2 |X_{0, r}| $. Then, these norms are up to a constant (the constant is 1 for $ G = \Sp(1, n) $, $ n \geq 2 $) equal to the 2-norm and
    \[
        a_{n_{v, r} a_t} = e^{ \arccosh\bigl(\sqrt{\big(\tfrac12(t + \tfrac1t) + \tfrac1{2t} \|v\|^2\big)^2 + \tfrac{|r|^2}{t^2}}\bigr) }
    \]
    (Theorem \ref{Thm explicit formulas} or \eqref{eq a_(n_(v, r) a_t) for Sp(1,n)} when $ G = \Sp(1, n) $, $ n \geq 2 $). 
    So, we must show that
    \begin{multline}\label{Proof of Lem a_(n_(v_1, phi(v_1) + r_1) n_(v_2, r_2) a) asymp a_(n_(v_1 + v_2, r_1 + r_2) a) asymp to show}
        \big(\tfrac12(t + \tfrac1t) + \tfrac1{2t} \|v+v'\|^2\big)^2 + \tfrac{|r + r' + \tfrac12 [(v, 0), (v', 0)] |^2}{t^2} \\
        \asymp \big(\tfrac12(t + \tfrac1t) + \tfrac1{2t} \|v+v'\|^2\big)^2 + \tfrac{|p(r)+r'|^2}{t^2}
    \end{multline}
    for $ t $ varying in $ (0, \infty) $ and $ (v + v', r + r') $ varying in a subset of $ \n $ having compact complement (we still assume that $ n_{v, r} \in N_{\Gamma, P} $ and $ n_{v', r'} \in N^{\Gamma, P} $). Note that a constant times the left hand side or the right hand side of \eqref{Proof of Lem a_(n_(v_1, phi(v_1) + r_1) n_(v_2, r_2) a) asymp a_(n_(v_1 + v_2, r_1 + r_2) a) asymp to show} is then greater or equal than 1 for all $ t \in (0, +\infty) $ and for $ (v + v', r + r') $ varying in an appropriate subset of $ \n $ having compact complement. Since $ p_{\g_{\alpha_P}}(\n_{\Gamma, P}) $ is orthogonal to $ \n_1^{\Gamma, P} $, $ \|v+v'\|^2 = \|v\|^2 + \|v'\|^2 $.
    Thus, the above is equivalent to show that
    \begin{multline*}
        \tfrac14 \big(t^2 + 1 + \|v\|^2 + \|v'\|^2\big)^2 + |r + r' + \tfrac12 [(v, 0), (v', 0)]|^2 \\
        \asymp \big(\tfrac12(t^2 + 1) + \tfrac12 (\|v\|^2 + \|v'\|^2)\big)^2 + |p(r)|^2 + |r'|^2
    \end{multline*}
    for $ t $ varying in $ (0, \infty) $ and $ (v + v', r + r') $ varying in a subset of $ \n $ having compact complement. If $ \g_{2\alpha_P} $ is trivial (real hyperbolic case), then the lemma follows already.
    
    Let $ c > 0 $. If $ \|v\|^4 + \|v'\|^4 \geq 1 $, then this is again equivalent to show that
    \begin{multline*}
        t^4 + c(\|v\|^4 + \|v'\|^4) + 1 + |r + r' + \tfrac12 [(v, 0), (v', 0)]|^2 \\
        \asymp 1 + t^4 + \|v\|^4 + \|v'\|^4 + |p(r)|^2 + |r'|^2
    \end{multline*}
    for $ t $ varying in $ (0, \infty) $ and $ (v + v', r + r') $ varying in a subset of $ \n $ having compact complement, by Lemma \ref{Lem about asymp} combined with Corollary \ref{Cor Estimate for Schwartz space = representation space}. 

    Write $ v = (v_1, v_2, \dotsc, v_{m_\alpha}) $ and $ v' = (v'_1, v'_2, \dotsc, v'_{m_\alpha}) $.
    Write $ r = (r_1, \dotsc, r_{m_{2\alpha}}) $ and $ r' = (r'_1, \dotsc, r'_{m_{2\alpha}}) $.
    For $ \eps > 0 $, set
    \begin{multline*}
        D_\eps = \{ (v, r, v', r') \in \RR^{m_\alpha} \times \RR^{m_{2\alpha}} \times \RR^{m_\alpha} \times \RR^{m_{2\alpha}} \\
        \mid n_{v, r} \in N_{\Gamma, P}, \ n_{v', r'} \in N^{\Gamma, P}, \ |v_i|^4 + |v'_i|^4 + |r_j|^2 + |r'_j|^2 > \eps \quad \forall i, j \} \ .
    \end{multline*}
    Since the Lie bracket is bilinear, $ [(v, 0), (v', 0)] $ is a polynomial in
    \[
        v_1, v_2, \dotsc, v_{m_\alpha}, v'_1, v'_2, \dotsc, v'_{m_\alpha} \ .
    \]
    Moreover, for every $ k \in \{1, \dotsc, m_{2\alpha}\} $,
    \[
        [(v, 0), (v', 0)]_k = \sum_{i=1}^{m_\alpha} \sum_{j=1}^{m_\alpha} a^{(k)}_{i, j} v_i v'_j
    \]
    for some constants $ a^{(k)}_{i, j} \in \RR $. 
    For $ c > 0 $ sufficiently large and $ (v, r, v', r') $ varying in $ D_c $, we have
    \begin{align*}
        & c (\|v\|^4 + \|v'\|^4) + |r + r' + \tfrac12 [(v, 0), (v', 0)]|^2 \\
        \asymp& c \big(\sum_{i=1}^{m_\alpha} v_i^4 + \sum_{i=1}^{m_\alpha} ({v'_i})^4 \big) + |r + r' + \tfrac12 [(v, 0), (v', 0)]|^2 \\
        \asymp& c \big(\sum_{i=1}^{m_\alpha} v_i^4 + \sum_{i=1}^{m_\alpha} ({v'_i})^4 \big) + \sum_{j=1}^{m_{2\alpha}} \big(r_j + r'_j + \tfrac12 \big([(v, 0), (v', 0)]\big)_j \big)^2 \\ 
        =& c \Big(\sum_{i=1}^{m_\alpha} v_i^4 + \sum_{i=1}^{m_\alpha} ({v'_i})^4 + \sum_{j=1}^{m_{2\alpha}} \big(\tfrac{r_j}{\sqrt{c}} + \tfrac{r'_j}{\sqrt{c}})^2 \\
        & \qquad + \sum_{j=1}^{m_{2\alpha}} (\tfrac{1}{\sqrt{c}} \cdot \tfrac{r_j + r'_j}{\sqrt{c}} \cdot \big([(v, 0), (v', 0)]\big)_j + \tfrac14 \big([(v, 0), (v', 0)]\big)^2_j \big) \Big) \\
        \asymp & \sum_{i=1}^{m_\alpha} \big( v_i^4 + ({v'_i})^4 \big) + \sum_{j=1}^{m_{2\alpha}} (r_j + r'_j)^2
        \asymp \|v\|^4 + \|v'\|^4 + |r + r'|^2
    \end{align*}
    by Lemma \ref{Lem about asymp} combined with Corollary \ref{Cor Estimate for Schwartz space = representation space}.
    Thus, for $ \eps > 0 $ sufficiently large and $ (v, r, v', r') \in D_\eps $, we have
    \begin{multline*}
        t^4 + c(\|v\|^4 + \|v'\|^4) + 1 + |r + r' + \tfrac12 [(v, 0), (v', 0)]|^2 \\
        \asymp t^4 + \|v\|^4 + \|v'\|^4 + 1 + |r + r'|^2 \ .
    \end{multline*}
    For each $ v \in p_{\g_{\alpha_P}}(\n_{\Gamma, P}) $, there exists $ r \in p_{\g_{2\alpha_P}}(\n_{\Gamma, P}) $ such that $ n_{v, r} \in N_{\Gamma, P} $. Write $ r = u + z $ with $ u \in \n_{2, 1}^{\Gamma, P} $ and $ z \in \g_{2\alpha_P} \cap \n_{\Gamma, P} $. Set $ \phi(v) = u $ and set $ \phi = 0 $ on $ \n_1^{\Gamma, P} $. This defines an $ \RR $-linear map $ \phi $ from $ \RR^{m_\alpha} $ to $ \RR^{m_{2\alpha}} $. One checks easily that this map is well-defined.
    So, there exist $ c_i \in \RR^{m_{2\alpha}} $ such that
    \[
        \phi\big((x_1, x_2, \dotsc, x_{m_\alpha})\big) = \sum_{i=1}^{m_\alpha} c_i x_i  \qquad (x_i \in \RR) \ .
    \]
    It follows from Lemma \ref{Lem about asymp} combined with Corollary \ref{Cor Estimate for Schwartz space = representation space} that
    \[
        \|v\|^4 + |r + r'|^2 \asymp \|v\|^4 + |p(r)|^2 + |r'|^2
    \]
    for $ (v, r, r') $ such that $ \|v\|^4 + |p(r)|^2 + |r'|^2 > 1 $.
    Indeed, the terms that are on the left but not on the right hand side are either of the form $ \lambda v_a^2 $ ($ \lambda \in \RR $) or of the form $ \lambda v_a r'_b $ ($ \lambda \in \RR $, $ r'_b \in \n_{2, 2}^{\Gamma, P} $).

    Since moreover, for $ \eps \geq 1 $,
    \[
        \{ (v + v', r + r') \mid n_{v, r} \in N_{\Gamma, P}, \ n_{v', r'} \in N^{\Gamma, P}, \ (v, r, v', r') \not \in D_\eps \}
    \]
    is contained
    \[
        \{ (v + v', r + r') \mid \max\{ |v_i|, |v'_i|, |r_j|, |r'_j|\} \leq \eps \quad \forall i, j \}
    \]
    (compact), the lemma follows.
\end{proof}

\begin{prop}\label{Prop a_(n_1 n_2 a) >= c a_(n_2 a)}
    There exist $ c, c' > 0 $ such that
    \[
        a_{n_1 n_2 a} \geq c \, a_{n_1 a} , \qquad a_{n_1 n_2 a} \geq c' a_{n_2 a}
    \]
    for all $ n_1 \in N_{\Gamma, Q} $, $ n_2 \in N^{\Gamma, Q} $ and $ a \in A_Q $.
\end{prop}

\begin{proof}
    If $eQ$ has full rank, then the proposition follows from Lemma \ref{Lem a_(na) estimates}.

    Without loss of generality, we may assume that $ Q = P $ and that $P$ is a $ \Gamma $-cuspidal parabolic subgroup with associated cusp of smaller rank.

    It follows from the proof of the previous lemma that there exists a compact set $ C $ in $ \n $ such that
    \[
        a_{n_{v, r} n_{v', r'} a_t} \asymp \frac1t \sqrt{1 + t^4 + \|v\|^4 + \|v'\|^4 + |p(r)|^2 + |r'|^2}
    \]
    for $t$ varying in $ (0, +\infty) $ and $ n_{v, r} $ (resp. $ n_{v', r'} $) varying in $ N_{\Gamma, P} $ (resp. $ N^{\Gamma, P} $) such that $ \log(n_{v, r}) + \log(n_{v', r'}) \not \in C $. In particular,
    \[
        a_{n_{v, r} a_t} \asymp \frac1t \sqrt{1 + t^4 + \|v\|^4 + |p(r)|^2}
    \]
    and
    \[
        a_{n_{v', r'} a_t} \asymp \frac1t \sqrt{1 + t^4 + \|v'\|^4 + |r'|^2}
    \]
    for $t$ varying in $ (0, +\infty) $ and $ n_{v', r'} $ varying in $ N^{\Gamma, P} $ such that $ \log(n_{v', r'}) \not \in C $. The proposition follows.
\end{proof}

\newpage

\subsubsection{The right regular representation of \texorpdfstring{$G$}{G} on the Schwartz space}

In this section, we show that the right translation by elements of $G$ defines a representation on the Schwartz space (see Theorem \ref{Thm representations geometrically finite case}). Some more preparations are needed to be able to show this result. After this, we are finally also able to show that the Schwartz space does not depend on choices (see Proposition \ref{Prop Schwartz space is independent of choices}).

\medskip

For a parabolic subgroup $ P' $ of $G$, set $ N_{\alpha_{P'}} = \exp(\g_{\alpha_{P'}}) $ and $ N_{2\alpha_{P'}} = \exp(\g_{2\alpha_{P'}}) $.\index[n]{N_alpha N_2alpha@ $ N_\alpha $, $ N_{2\alpha} $}
Let
\[
    p_i \colon N_Q \to N_{i \alpha_Q} \, , \quad n \mapsto \exp\bigl(p_{\g_{i \alpha_Q}}(\log(n))\bigr)
    \qquad (i \in \{1, 2\}) \ .
\]
Then, $ \log \circ \, p_i = p_{\g_{i\alpha_Q}} \circ \log $ and for every $ n \in N_Q $
\[
    n = p_1(n) p_2(n)
\]
by Lemma \ref{Lem using the Baker-Campbell-Hausdorff formula}. Hence, $ p_1 $ and $ p_2 $ are projections.

Let $ p_{N_i^{\Gamma, Q}} = p_i \circ \, p_{N^{\Gamma, Q}} $. Then, this defines a projection from $ N_Q $ to $ N_i^{\Gamma, Q} $.
\\ Clearly, $ p_i \circ \, p_{N_i^{\Gamma, Q}} = p_{N_i^{\Gamma, Q}} $.
Let $ \n_{2, 1}^{\Gamma, Q} = p_{\g_{2\alpha_Q}}(\n_{\Gamma, Q}) \cap \n_2^{\Gamma, Q} $. We denote the projection on $ \n_{2, 1}^{\Gamma, Q} $ by $ p_{ \n_{2, 1}^{\Gamma, Q} } $.
Let $ X \in p_{\g_{\alpha_Q}}(\n_{\Gamma, Q}) $. Then, there exists $ Y \in \n_{\Gamma, Q} $ such that $ p_{\g_{\alpha_Q}}(Y) = X $.
Let $ \phi(X) = p_{ \n_{2, 1}^{\Gamma, Q} }( Y ) $.
Then, $ \phi $ is a well-defined $ \RR $-linear map from $ p_{\g_{\alpha_Q}}(\n_{\Gamma, Q}) $ to $ \n_{2, 1}^{\Gamma, Q} $. Indeed, let $ Y' \in \n_{\Gamma, Q} $ be such that $ p_{\g_{\alpha_Q}}(Y') = X $. Then, $ p_{\g_{\alpha_Q}}(Y - Y') = 0 $. Thus, $ Y - Y' \in \n_{\Gamma, Q} \cap \g_{2\alpha_Q} $. Hence, $ p_{ \n_{2, 1}^{\Gamma, Q} }( Y - Y' ) = 0 $. So, $ p_{ \n_{2, 1}^{\Gamma, Q} }( Y ) = p_{ \n_{2, 1}^{\Gamma, Q} }( Y' ) $.

By construction, $ p_{ \n_{2, 1}^{\Gamma, Q} }(X) = \phi\big(p_{\g_{\alpha_Q}}(X)\big) $ for all $ X \in \n_{\Gamma, Q} $.

\begin{lem}\label{Lem N = p_1(N_Gamma) N_1^Gamma N_2 decomposition} \nlenum
    \begin{enumerate}
    \item $ N_Q = p_1(N_{\Gamma, Q}) N_1^{\Gamma, Q} N_{2 \alpha_Q} $ is a decomposition of $ N_Q $. Moreover, for every $ n \in N_Q $, there exists a unique $ n' \in N_{2\alpha_Q} $ such that
        \[
            n = (p_1 \circ p_{N_{\Gamma, Q}})(n) p_{N_1^{\Gamma, Q}}(n) n' \ .
        \]
        It follows that $ p_1 \circ p_{N_{\Gamma, Q}} \colon N_Q \to p_1(N_{\Gamma, Q}) $ is a projection.
        \\ Moreover, the above decomposition induces an $A_Q$-invariant Lie group decomposition on $ N_Q/N_{2\alpha_Q} $.
    \item For every $ n \in N_Q $, $ \log(n) $ is equal to
        \[
            \log\bigl( (p_1 \circ \, p_{N_{\Gamma, Q}})(n) \bigr) + \log\bigl( p_{N_1^{\Gamma, Q}}(n) \bigr) + \log\bigl( p_2(n) \bigr) \ .
        \]
        Thus,
        \[
            \log \circ \, p_1 \circ \, p_{N_{\Gamma, Q}} = p_{\g_{\alpha_Q}} \circ \, p_{\n_{\Gamma, Q}} \circ \log \, , \quad \log \circ \, p_{N_1^{\Gamma, Q}} = p_{\n_1^{\Gamma, Q}} \circ \log \ .
        \]
    \item $ p_1 \circ \, p_{N_{\Gamma, Q}} $ is a projection from $ N_Q $ to $ p_1(N_{\Gamma, Q}) = \exp\bigl(p_{\g_{\alpha_Q}}(\n_{\Gamma, Q})\bigr) $.
    \end{enumerate}
\end{lem}

\begin{rem}
    In general, we don't have $ \log \circ \, p_{N_{\Gamma, Q}} = p_{\n_{\Gamma, Q}} \circ \log $, $ \log \circ \, p_{N^{\Gamma, Q}} = p_{\n^{\Gamma, Q}} \circ \log $, $ \log \circ \, p_2 \circ \, p_{N_{\Gamma, Q}} = p_{\g_{2 \alpha_Q}} \circ \, p_{\n_{\Gamma, Q}} \circ \log $ or $ \log \circ \, p_{N_2^{\Gamma, Q}} = p_{\n_2^{\Gamma, Q}} \circ \log $.
\end{rem}

\begin{proof}
    Let $ n \in N_Q $. By Proposition \ref{Prop Decomposition of N}, by Lemma \ref{Lem using the Baker-Campbell-Hausdorff formula}, we have
    \[
        n = p_1(n) p_2(n) = p_1(p_{N_{\Gamma, Q}}(n) p_{N^{\Gamma, Q}}(n)) p_2(n)
        = (p_1 \circ p_{N_{\Gamma, Q}})(n) p_{N_1^{\Gamma, Q}}(n) n' \ ,
    \]
    where $ n' := p_2(n) \exp(-\frac12[\log(p_1(p_{N_{\Gamma, Q}}(n))), \log(p_{N_1^{\Gamma, Q}}(n))]) \in N_{2\alpha_Q} $.
    \\ It follows that $ N_Q = p_1(N_{\Gamma, Q}) N_1^{\Gamma, Q} N_{2\alpha_Q} $.
    \\ Let $ n_1, n'_1 \in p_1(N_{\Gamma, Q}) $, $ n_2, n_2' \in N_1^{\Gamma, Q} $, $ n_3, n_3' \in N_{2\alpha_Q} $ be such that $ n_1 n_2 n_3 = n'_1 n'_2 n'_3 $.

    Let $ X_1, X'_1 \in p_{\g_{\alpha_Q}}(\n_{\Gamma, Q}) $ be such that $ n_1 = \exp(X_1) $, $ n'_1 = \exp(X'_1) $.
    \\ Let $ X_2, X'_2 \in \n_1^{\Gamma, Q} $ be such that $ n_2 = \exp(X_2) $, $ n'_2 = \exp(X'_2) $.

    It follows from the Baker-Campbell-Hausdorff formula and the fact that every element of $ \n_Q = \g_{\alpha_Q} \oplus \g_{2\alpha_Q} $ admits a unique decomposition that
    \[
        X_1 + X_2 = X'_1 + X'_2 \ .
    \]
    Thus, $ X_1 = X'_1 $ and $ X_2 = X'_2 $ as $ \g_{\alpha_Q} = p_{\g_{\alpha_Q}}(\n_{\Gamma, Q}) \oplus \n_1^{\Gamma, Q} $. Hence, $ n_1 = n'_1 $ and $ n_2 = n'_2 $. So, $ n_3 $ is also equal to $ n'_3 $.

    The first assertion follows now easily from this.

    The first part of the second assertion follows now from the Baker-Campbell-Hausdorff formula and the second part follows from the fact that $ (p_{\g_{\alpha_Q}} \circ \, p_{\n_{\Gamma, Q}})(X) + p_{\n_1^{\Gamma, Q}}(X) + p_{\g_{2\alpha_Q}}(X) $ gives a unique decomposition of $ X \in \n_Q $.

    We have
    \[
        p_1(N_{\Gamma, Q}) = \exp( p_{\g_{\alpha_Q}}\big( \log(N_{\Gamma, Q}) \big) ) = \exp( p_{\g_{\alpha_Q}} \n_{\Gamma, Q}) \ .
    \]
    Let $ n \in p_1(N_{\Gamma, Q}) $. Then, by the previous assertion and by uniqueness of the decomposition of $ \n_Q = p_{\g_{\alpha_Q}}(\n_{\Gamma, Q}) \oplus \n_1^{\Gamma, Q} \oplus \g_{2\alpha_Q} $, we have
    \[
        (p_1 \circ \, p_{N_{\Gamma, Q}})(n) = \exp\bigl( \log( (p_1 \circ \, p_{N_{\Gamma, Q}})(n) ) \bigr)
        = \exp\bigl( (p_{\g_{\alpha_Q}} \circ \, p_{\n_{\Gamma, Q}})(\log(n)) \bigr) = n \ .
    \]
    The last assertion follows.
\end{proof}

\begin{lem}\label{Lem |[X, Y]| <= c |X| |Y|}
    There exists a constant $ c > 0 $ such that $ |[X, Y]| \leq c |X| \cdot |Y| $ for all $ X, Y \in \g $.
\end{lem}

\begin{proof}
    Let $ X, Y \in \g $. Without loss of generality, we may assume that $ X $ and $Y$ are nonzero. Then,
    \[
        [X, Y] = |X| \cdot |Y| \cdot [\tfrac{X}{|X|}, \tfrac{Y}{|Y|}] \ .
    \]
    Let $ S = \{ Z \in \g \mid |Z| = 1 \} $ (compact) and let $ c = \sup_{Z_1, Z_2 \in S} |[Z_1, Z_2]| $.
    Since $ [\cdot, \cdot] \colon \g \times \g \to \g $ is continuous and as $ S \times S $ is compact, $ c $ is finite. The lemma follows.
\end{proof}

\begin{prop}\label{Prop gamma g h in Siegel set with estimates}
    Let $ \Sfrak_{std, Q, t, \omega, \omega'} $ be a generalised standard Siegel set with respect to $ Q $ and let $ S $ be a relatively compact subset of $G$. Then, there exists a generalised standard Siegel set $ \Sfrak_{std, Q, \tilde{t}, \omega, \tilde{\omega}'} $ such that
    \[
        \gamma_{g, h} g h \in \Sfrak_{std, Q, \tilde{t}, \omega, \tilde{\omega}'}
    \]
    for some $ \gamma_{g, h} \in \Gamma_Q $, for all $ g \in \Sfrak_{std, Q, t, \omega, \omega'} $ and all $ h \in S $.
    Moreover, we have the following:
    \begin{enumerate}
    \item There exist constants $ c > 0 $, $ c' > 0 $ such that
        \[
            \sup_{h \in S} |\log(p_1(\nu_Q(\gamma_{g, h})))| \leq c + c' \, h_Q(g)^{\alpha_Q}
        \]
        for all $ g \in \Sfrak_{std, Q, t, \omega, \omega'} $.
    \item There exist constants $ c > 0 $, $ c' > 0 $, $ c'' > 0 $ such that
        \[
            \sup_{h \in S} |\log\bigl(\nu_Q(\gamma_{g, h})\bigr)| \leq c
            + c' \, h_Q(g)^{\alpha_Q} (1 + |\log\bigl( p_{N_1^{\Gamma, Q}}(\nu_Q(g)) \bigr)|)
            + c'' \, h_Q(g)^{2\alpha_Q}
        \]
        for all $ g \in \Sfrak_{std, Q, t, \omega, \omega'} $.
    \end{enumerate}
\end{prop}

\begin{rem}\label{Rem of Prop gamma g h in Siegel set} \nlenum
    \begin{enumerate}
    \item It follows from the above lemma that
        \[
            \{ |\log(p_1(\nu_Q(\gamma_{g, h})))| \mid h \in S, g \in \Sfrak_{std, Q, t, \omega, \omega'} : h_Q(g)^{\alpha_Q} \leq T \}
        \]
        is bounded for all $ T > 0 $.
    \item \begin{sloppypar}
        Let $ \{X_j\} $ be an orthonormal basis of $ \n_Q $. Then, the proof of the proposition provides estimates of the absolute value of the coefficients of $ \log(p_1(\nu_Q(\gamma_{g, h}))) $ and of $ \log\bigl(\nu_Q(\gamma_{g, h})\bigr) $ ($ h \in S $) with respect to the basis $ \{X_j\} $.
        \end{sloppypar}
    \end{enumerate}
\end{rem}

\begin{proof}
    Let $ \Sfrak := \Sfrak_{std, Q, t, \omega, \omega'} $ and $ S $ be as above and let $ h \in S $. Let $ g \in G $.
    \\ We have
    \begin{align}\label{eq gh}
        g h &= \nu_Q(g h) h_Q(gh) k_Q(gh) \\
        &= \big( \nu_Q(g) h_Q(g) \nu_Q(k_Q(g)h) h_Q(g)^{-1} \big) h_Q(g) h_Q(k_Q(g)h) k_Q(gh) \ . \nonumber
    \end{align}
    By Lemma \ref{Lem gamma g h in Siegel set}, there exists a generalised Siegel standard set $ \Sfrak_{std, Q, \tilde{t}, \omega, \tilde{\omega}'} $ such that
    \[
        \gamma_{g, h} g h \in \Sfrak_{std, Q, \tilde{t}, \omega, \tilde{\omega}'}
    \]
    for some $ \gamma_{g, h} \in \Gamma_Q $, for all $ g \in \Sfrak_{std, Q, t, \omega, \omega'} $ and all $ h \in S $.

    Note that \eqref{eq gh} is again equal to
    \begin{multline*}
        \big( p_1(\nu_Q(g)) h_Q(g)^{\alpha_Q} p_1(\nu_Q(k_Q(g)h)) p_2(\nu_Q(g)) h_Q(g)^{2 \alpha_Q} p_2(\nu_Q(k_Q(g)h)) \big) \\
        h_Q(g) h_Q(k_Q(g)h) k_Q(gh) \ .
    \end{multline*}
    Thus, by Lemma \ref{Lem using the Baker-Campbell-Hausdorff formula}, $ \nu_Q(gh) $ is equal to
    \begin{align*} 
        & p_1(p_{N_{\Gamma, Q}}(\nu_Q(g))) p_{N_1^{\Gamma, Q}}(\nu_Q(g))
        \exp\bigl(-\tfrac12 [p_1(p_{N_{\Gamma, Q}}(\nu_Q(g))), p_{N_1^{\Gamma, Q}}(\nu_Q(g))] \bigr) \\
        & \qquad h_Q(g)^{\alpha_Q} p_1\big( \nu_Q(k_Q(g)h) \big)
        p_2(\nu_Q(g)) h_Q(g)^{2 \alpha_Q} p_2(\nu_Q(k_Q(g)h)) \ . 
    \end{align*}
    Hence, by the Baker-Campbell-Hausdorff formula, $ \log\bigl(\nu_Q(g h)\bigr) $ is equal to
    \begin{align*}
        & \log\bigl(p_1(p_{N_{\Gamma, Q}}(\nu_Q(g)))\bigr) + \log\bigl( p_{N_1^{\Gamma, Q}}(\nu_Q(g)) \bigr)
        + h_Q(g)^{\alpha_Q} \log\bigl( p_1( \nu_Q(k_Q(g)h) ) \bigr) \\
        & \qquad + \tfrac12 h_Q(g)^{\alpha_Q}[\log\bigl(p_1(p_{N_{\Gamma, Q}}(\nu_Q(g)))\bigr), \log\bigl( p_1( \nu_Q(k_Q(g)h) ) \bigr)] \\
        & \qquad + \tfrac12 h_Q(g)^{\alpha_Q}[\log\bigl( p_{N_1^{\Gamma, Q}}(\nu_Q(g)) \bigr), \log\bigl( p_1( \nu_Q(k_Q(g)h) ) \bigr)] \\
        & \qquad + \log\bigl( p_2(\nu_Q(g)) \bigr)
        + h_Q(g)^{2 \alpha_Q} \log\bigl( p_2(\nu_Q(k_Q(g)h)) \bigr)
        \big) \ .
    \end{align*}
    By Lemma \ref{Lem N = p_1(N_Gamma) N_1^Gamma N_2 decomposition}, this yields
    \begin{align}\label{eq log(nu(g h))}
        & \log\bigl(\nu_Q(g)\bigr)
        + h_Q(g)^{\alpha_Q} \log\bigl( p_1( \nu_Q(k_Q(g)h) ) \bigr) \\ \nonumber
        & \qquad + \tfrac12 h_Q(g)^{\alpha_Q}[\log\bigl( p_1(p_{N_{\Gamma, Q}}(\nu_Q(g))) \bigr), \log\bigl( p_1( \nu_Q(k_Q(g)h) ) \bigr)] \\ \nonumber
        & \qquad + \tfrac12 h_Q(g)^{\alpha_Q}[\log\bigl( p_{N_1^{\Gamma, Q}}(\nu_Q(g)) \bigr), \log\bigl( p_1( \nu_Q(k_Q(g)h) ) \bigr)] \\ \nonumber
        & \qquad + h_Q(g)^{2 \alpha_Q} \log\bigl( p_2(\nu_Q(k_Q(g)h)) \bigr) \ .
    \end{align}
    So, $ \log\bigl(p_{N_1^{\Gamma, Q}}(\nu_Q(g h))\bigr) = p_{\g_{\alpha_Q}}\Bigl(p_{\n^{\Gamma, Q}} \log\bigl(\nu_Q(g h)\bigr)\Bigr) $ is equal to
    \[
        \log\bigl(p_{N_1^{\Gamma, Q}}(\nu_Q(g))\bigr) + h_Q(g)^{\alpha_Q} \log\bigl( p_{N_1^{\Gamma, Q}}(\nu_Q(k_Q(g)h)) \bigr) \ .
    \]
    One can also use these computations in order to prove the first assertion.

    \medskip

    \needspace{2\baselineskip}
    Let us estimate in the following $ |\log(p_1(\nu_Q(\gamma_{g, h})))| = |p_{\g_{\alpha_Q}}\big(\log(\nu_Q(\gamma_{g, h}))\big)|$.

    By \eqref{eq log(nu(g h))} and by Lemma \ref{Lem using the Baker-Campbell-Hausdorff formula}, $ \log\bigl(p_1(p_{N_{\Gamma, Q}}(\nu_Q(\gamma_{g, h} g h)))\bigr) $ is equal to
    \begin{multline*}
        \log\bigl(p_1(p_{N_{\Gamma, Q}}(\nu_Q(\gamma_{g, h})))\bigr) + \log\bigl(p_1(p_{N_{\Gamma, Q}}(\nu_Q(l_{M_Q}(\gamma_{g, h}) g)))\bigr) \\
        + h_Q(g)^{\alpha_Q} \log\bigl( p_1( p_{N_{\Gamma, Q}}( \nu_Q(k_Q(g)h) ) ) \bigr) \ .
    \end{multline*}
    Since $ \omega $ is relatively compact, there exists a constant $ c > 0 $ such that
    \[
        \sup_{ h \in S, \, m \in M_{\Gamma, Q} } |\log\bigl(p_1(p_{N_{\Gamma, Q}}(\nu_Q(\gamma_{g, h} g h)))\bigr)
        - \log\bigl(p_1(p_{N_{\Gamma, Q}}(\nu_Q(m g)))\bigr)| \leq c
    \]
    for all $ g \in \Sfrak $.
    Thus,
    \[
        \sup_{ h \in S } | \log\bigl(p_1(p_{N_{\Gamma, Q}}(\nu_Q(\gamma_{g, h})))\bigr) | \leq c + c' h_Q(g)^{\alpha_Q} \ ,
    \]
    where $ c' := \sup_{ k \in K_Q, \, h \in S } |\log\bigl( p_1( p_{N_{\Gamma, Q}}( \nu_Q(kh) ) ) \bigr)| $.

    \medskip

    Let us estimate now $ | \log\bigl(p_{N_{2\alpha_Q} \cap N_{\Gamma, Q}}(\nu_Q(\gamma_{g, h}))\bigr) | $.

    By \eqref{eq log(nu(g h))} and by Lemma \ref{Lem using the Baker-Campbell-Hausdorff formula},
    $
        p_{\g_{2\alpha_Q} \cap \n_{\Gamma, Q}} \log\bigl(\nu_Q(\gamma_{g, h} g h)\bigr)
    $
    is equal to
    \begin{align*}
        & p_{\g_{2\alpha_Q} \cap \n_{\Gamma, Q}}\big(\log(\nu_Q(\gamma_{g, h}))\big) + \log\bigl(p_{N_{2\alpha_Q} \cap N_{\Gamma, Q}}(\nu_Q(l_{M_Q}(\gamma_{g, h}) g))\bigr) \\
        & \qquad + h_Q(g)^{2 \alpha_Q} p_{\g_{2\alpha_Q} \cap \n_{\Gamma, Q}}\big( \log( \nu_Q(k_Q(g)h) ) \big) \\
        & \qquad + \tfrac12 h_Q(g)^{\alpha_Q} p_{\g_{2\alpha_Q} \cap \n_{\Gamma, Q}}([\log\bigl( p_1(\nu_Q(\gamma_{g, h})) \bigr), \log\bigl( p_1( \nu_Q(k_Q(g)h) ) \bigr)]) \\
        & \qquad + \tfrac12 h_Q(g)^{\alpha_Q} p_{\g_{2\alpha_Q} \cap \n_{\Gamma, Q}}([\log\bigl( p_1( \Ad(l_{M_Q}(\gamma_{g, h})) p_{N_{\Gamma, Q}}(\nu_Q(g)) ) \bigr), \log\bigl( p_1( \nu_Q(k_Q(g)h) ) \bigr)]) \\
        & \qquad + \tfrac12 h_Q(g)^{\alpha_Q} p_{\g_{2\alpha_Q} \cap \n_{\Gamma, Q}}([\log\bigl( p_{N_1^{\Gamma, Q}}(\nu_Q(g)) \bigr), \log\bigl( p_1( \nu_Q(k_Q(g)h) ) \bigr)]) \ .
    \end{align*}
    Since $ \omega $ is relatively compact, there exists a constant $ c > 0 $ such that
    \[
        \sup_{h \in S, \, m \in M_{\Gamma, Q}} | p_{\g_{2\alpha_Q} \cap \n_{\Gamma, Q}}\big(\log(\nu_Q(\gamma_{g, h} g h))\big)
        - p_{\g_{2\alpha_Q} \cap \n_{\Gamma, Q}}\big(\log(\nu_Q(m g))\big)| \leq c
    \]
    for all $ g \in \Sfrak $. Similar as above and by using that result and Lemma \ref{Lem |[X, Y]| <= c |X| |Y|}, one shows that there exist constants $ c' > 0 $, $ c'' > 0 $ such that
    \begin{multline}\label{eq Estimate of sup_h |p_(n_2 cap n_Gamma) log(gamma)|}
        \sup_{h \in S} | p_{\g_{2\alpha_Q} \cap \n_{\Gamma, Q}} \bigl(\log(\nu_Q(\gamma_{g, h}))\bigr) | \\
        \leq c + c' h_Q(g)^{\alpha_Q} (1 + |\log\bigl( p_{N_1^{\Gamma, Q}}(\nu_Q(g)) \bigr)|) + c'' h_Q(g)^{2\alpha_Q} \ .
    \end{multline}
    Since $ | p_{ \n_{2, 1}^{\Gamma, Q} }\bigl( \log\bigl( \nu_Q(\gamma_{g, h}) \bigr) \bigr) | = | \phi\bigl( p_{\g_{\alpha_Q}}\bigl( \log( \nu_Q(\gamma_{g, h}) ) \bigr) \bigr) | $ and since any linear map between finite-dimensional vector spaces is bounded, $ \sup_{h \in S} | p_{ \n_{2, 1}^{\Gamma, Q} }\bigl( \log\bigl( \nu_Q(\gamma_{g, h}) \bigr) \bigr) | $ has an estimate as in \eqref{eq Estimate of sup_h |p_(n_2 cap n_Gamma) log(gamma)|}. The last assertion follows now from the previous one's as
    \[
        X = p_{\g_{\alpha_Q}}(X) + p_{\g_{2\alpha_Q} \cap \n_{\Gamma, Q}}(X) + p_{\n_{2, 1}^{\Gamma, Q}}(X) \ .
    \]
    for all $ X \in \n_{\Gamma, Q} $.
\end{proof}

\begin{lem}\label{Lem Ad(gamma)X - X}
    Let $ \gamma \in \Gamma_Q $.
    \begin{enumerate}
    \item If $ X \in \n_Q $, then
        \[
            \Ad(\gamma)X - X = [\log(p_1(\nu_Q(\gamma))), \Ad(l_{M_Q}(\gamma))X] \in \g_{2\alpha_Q} \subset \n_Q \ .
        \]
    \item If $ X \in \a_Q $, then $ \Ad(\gamma)X - X $ is equal to
        \[
            [\log( \nu_Q(\gamma) ), X]
            + \tfrac12 [\log(p_1(\nu_Q(\gamma))), [\log(p_1(\nu_Q(\gamma))), X]]
            \in \n_{\Gamma, Q} \oplus \g_{2\alpha_Q} \subset \n_Q \ .
        \]
    \item If $ X \in \m_Q $, then $ \Ad(\gamma)X - X $ is equal to
        \begin{multline*}
            [\log( \nu_Q(\gamma) ), \Ad(l_{M_Q}(\gamma))X] \\
            + \tfrac12 [\log(p_1(\nu_Q(\gamma))), [\log(p_1(\nu_Q(\gamma))), \Ad(l_{M_Q}(\gamma))X]] \in \n_Q \ .
        \end{multline*}
    \end{enumerate}
\end{lem}

\begin{proof}
    The lemma follows by direct computation.
\end{proof}

For $ f \in \CS(\Gamma \bs G, \phi) $, define
\[
    \pi(g) f(x) = R_g f(x) = f(xg)  \qquad (g, x \in G) \ .
\]
\begin{prop}\label{Prop convergence of Sfrak p_(r, X, Y)(R_h f - f)}
    Let $ \Sfrak $ be a generalised standard Siegel set in $G$ with respect to $ Q $, let $ X \in \U(\n_Q) $, let $ Y \in \U(\g) $, let $ f \in \CS(\Gamma \bs G, \phi) $ and let $ r \geq 0 $. Then, for each relatively compact subset $S$ of $G$, there exists a seminorm $ p \colon \CS(\Gamma \bs G, \phi) \to \RR $ such that $ \sup_{h \in S} \LO{p}{\Sfrak'}_{r, X, Y}(\pi(h) f) $ is less or equal than $ p(f) < \infty $ and $ \LO{p}{\Sfrak'}_{r, X, Y}(R_h f - f) $ converges to zero when $ h $ tends to $e$.
\end{prop}

\begin{proof}
    To prove this proposition, we do a similar argument as in Proposition \ref{Prop convergence of Up_(r, X, Y)(R_h f - f)}, which itself was proven in a similar way as Lemma 14 of \cite[p.21]{HC66}.

    Let $ h \in G $, $ x \in G $, $ X \in \U(\n_Q) $, $ Y \in \U(\g) $, $ \Sfrak' $ as above and let $ f \in \CS(\Gamma \bs G, \phi) $.
    By Proposition \ref{Prop a_(n_1 n_2 a) >= c a_(n_2 a)}, there exists $ C > 0 $ such that
    \[
        a_{\gamma g} \geq C \, a_{p_{N^{\Gamma, Q}}(\nu_Q(g)) h_Q(g)} \qquad (\gamma \in \Gamma_Q, g \in \Sfrak) \ .
    \]
    This is again greater or equal than $ C' a_g $ for all $ g \in \Sfrak $. We have
    \begin{align*}
        \LO{p}{\Sfrak'}_{r, X, Y}(\pi(h)f) &= \sup_{gK_Q \in \Sfrak'}{ a_Q(g^{-1})^{\rho_{\Gamma_Q}} (1 + \log a_g)^r a_g^{\rho^{\Gamma_Q}} |L_X R_Y (R_h f)(g)| } \\
        &\leq \sup_{gK_Q \in \Sfrak'}{ a_Q(g^{-1})^{\rho_{\Gamma_Q}} (1 + \log a_g)^r a_g^{\rho^{\Gamma_Q}} |L_X R_{\Ad(h^{-1})Y} f(gh)| } \ .
    \end{align*}
    Let $S$ be a relatively compact subset of $G$. By Proposition \ref{Prop gamma g h in Siegel set with estimates}, there exists a generalised standard Siegel set $ \tilde{\Sfrak} $ with respect to $Q$ such that
    \[
        \gamma_{g, h} g h \in \tilde{\Sfrak}
    \]
    for some $ \gamma_{g, h} \in \Gamma_Q $, for all $ g \in \Sfrak $ and all $ h \in S $.

    By Lemma \ref{Lem 1 + log a_h <= (1 + log a_gh)(1 + log a_g)}, Lemma \ref{Lem about a(g)} and the estimates $ C' a_g \leq a_{\gamma_{g, h} g} $ ($ g \in \Sfrak $), $ a_x \leq a_{x h} a_h $ ($ x \in G $), the above is again less or equal than
    \[
        C'' \sup_{g K_Q \in \Sfrak'} a_Q((\gamma_{g, h} g h)^{-1})^{\rho_{\Gamma_Q}} (1 + \log a_{\gamma_{g, h} g h})^r a_{\gamma_{g, h} g h}^{\rho^{\Gamma_Q}} 
        |L_{\Ad(\gamma_{g, h})X} R_{\Ad(h^{-1})Y} f(\gamma_{g, h} g h)|
    \]
    for some constant $ C'' > 0 $.
    Assume now that $ X \in \n_Q $. The case $ X \in \U(\n_Q)_l $ will then follow by induction. Let $ t > 0 $ and let $ g \in \Sfrak $ be such that $ h_Q(g) \geq t $.
    \\ Set $ Z_{g, h} = [\log(p_1(\nu_Q(\gamma_{g, h}))), \Ad(l_{M_Q}(\gamma_{g, h}))X] \in \g_{2\alpha_Q} $.
    Then, for $ x \in G $, we have
    \begin{align*}
        L_{\Ad(\gamma_{g, h})X - X} R_{\Ad(h^{-1})Y} f(x) &= L_{Z_{g, h}} R_{\Ad(h^{-1})Y}  f(x) \\
        &= -R_{\Ad(k_Q(g)^{-1}) h_Q(g)^{-2\alpha_Q} Z_{g, h}} R_{\Ad(h^{-1})Y} f(x) \ .
    \end{align*}
    Since $ K_Q $ is compact and since the coefficients of $ h_Q(g)^{-2\alpha_Q} Z_{g, h} $ with respect to $ \{X_i\} $ are bounded by Proposition \ref{Prop gamma g h in Siegel set with estimates} and Lemma \ref{Lem |[X, Y]| <= c |X| |Y|} and since
    \[
        \{ |\log(p_1(\nu_Q(\gamma_{g, h})))| \mid h \in S, g \in \Sfrak_{std, Q, t, \omega, \omega'} : h_Q(g) \leq t \}
    \]
    is finite by Remark \ref{Rem of Prop gamma g h in Siegel set} of Proposition \ref{Prop gamma g h in Siegel set with estimates},
    \begin{multline*}
        \sup_{h \in S} \sup_{g K_Q \in \Sfrak'} a_Q((\gamma_{g, h} g h)^{-1})^{\rho_{\Gamma_Q}} (1 + \log a_{\gamma_{g, h} g h})^r a_{\gamma_{g, h} g h}^{\rho^{\Gamma_Q}} \\
        |L_{\Ad(\gamma_{g, h})X - X} R_{\Ad(h^{-1})Y} f(\gamma_{g, h} g h)|
    \end{multline*}
    is less or equal than $ l_1(f) $ for some seminorm $ l_1 \colon \CS(\Gamma \bs G, \phi) \to \RR $.

    Thus, $ \sup_{h \in S} \LO{p}{\Sfrak'}_{r, X, Y}(\pi(h)f) $ is less or equal than
    \[
        l_2(f) := l_1(f) + \sup_{h \in S} \sup_{xK_Q \in \tilde{\Sfrak}'}{ a_Q(x^{-1})^{\rho_{\Gamma_Q}} (1 + \log a_x)^r a_x^{\rho^{\Gamma_Q}} |L_X R_{\Ad(h^{-1})Y} f(x)| } \ .
    \]
    In particular, $ \pi(h)f $ belongs to $ \CS(\Gamma \bs G, \phi) $.
    The remaining argument is very similar to the one in Proposition \ref{Prop convergence of Up_(r, X, Y)(R_h f - f)}.
\end{proof}

\begin{thm}\label{Thm representations geometrically finite case}
    $ \big(\pi, \CS(\Gamma \bs G, \phi)\big) $ and $ \big(\pi, \c{\CS}(\Gamma \bs G, \phi)\big) $ are representations of $G$.
\end{thm}

\begin{proof}
    The theorem follows from Proposition \ref{Prop convergence of Sfrak p_(r, X, Y)(R_h f - f)} for the $ \LO{p}{\Sfrak'}_{r, X, Y} $ seminorms and the proof of Theorem \ref{Thm representations convex-cocompact case} for the $ \LOI{p}{U}_{r, X, Y} $ seminorms. Indeed, the cocompact action of $ \Gamma $ on $ X \cup \Omega_\Gamma $ is only needed to show that the Schwartz space lies in $ L^2(\Gamma \bs G, \phi) $.
    The detailed proof is very similar to the one of Theorem \ref{Thm representations convex-cocompact case}.
\end{proof}

\begin{prop}\label{Prop Schwartz space is independent of choices}
    The Schwartz space does not depend on made choices.
\end{prop}

\begin{proof}
    Since the Schwartz space is a representation space by Proposition \ref{Prop convergence of Sfrak p_(r, X, Y)(R_h f - f)}, it suffices to prove that alternative seminorms are finite for Schwartz functions.

    Let $ f \in \CS(\Gamma \bs G, \phi) $ and $ r \geq 0 $.
    \begin{enumerate}
    \item We can choose an other Cartan involution $ \theta'_Q $ of $ \g $ in order to get $ A_Q $ and $ M_Q $.
        Indeed, let $ K'_Q $ be the maximal compact subgroup of $G$ corresponding to $ \theta'_Q $. 
        Let $ \tilde{h} \in N_Q A_Q $ be such that $ K'_Q = \tilde{h} K_Q \tilde{h}^{-1} $.

        Let $ \Sfrak'_{\theta'_Q, Q} $ be a generalised Siegel set in $G/K_{\theta'_Q} $ with respect to $ Q $. Set $ \Sfrak'_Q = \{ g \tilde{h} K_Q \mid gK_{\theta'_Q} \in \Sfrak'_{\theta'_Q, Q} \} $.
        Then, this is a generalised Siegel set in $G/K_Q$ with respect to $ Q $.

        Let $ X \in \U(\n_Q) $ and $ Y \in \U(\g) $.
        By \eqref{eq Compare a(g)}, \eqref{eq Compare a_g}, Lemma \ref{Lem about a(g)} and as
        \[
            gK'_Q \in \Sfrak'_{\theta'_Q, Q} \iff g \tilde{h} K_Q \in \Sfrak'_Q \ ,
        \]
        \[
            \sup_{gK'_Q \in \Sfrak'_{\theta'_Q, Q}}{ a_{\theta'_Q, Q}(g^{-1})^{\rho_{\Gamma_Q}} (1 + \log a_g)^r a_g^{\rho^{\Gamma_Q}} |L_X R_Y f(g)| }
        \]
        is finite if and only if
        \[
            \sup_{gK_Q \in \Sfrak'_Q}{ a_{\theta_Q, Q}(g^{-1})^{\rho_{\Gamma_Q}} (1 + \log a_g)^r a_g^{\rho^{\Gamma_Q}} |L_X R_{\Ad(\tilde{h}^{-1})Y} (R_{\tilde{h}^{-1}} f)(g)| }
        \]
        is finite.
        It follows from Proposition \ref{Prop convergence of Sfrak p_(r, X, Y)(R_h f - f)} that this is finite.
    \item We can choose an other representative of the $ \Gamma $-cuspidal parabolic subgroup of $G$.
        Indeed, let $ \gamma \in \Gamma $ and let $ Q' = \gamma Q \gamma^{-1} $. Then, $ N_{Q'} = \gamma N_Q \gamma^{-1} $, $ A_{Q'} = \gamma A_Q \gamma^{-1} $, $ M_{Q'} = \gamma M_Q \gamma^{-1} $ and $ K_{Q'} = \gamma K_Q \gamma^{-1} $.
        
        Let $ \Sfrak_{Q'} $ be a generalised Siegel set in $G$ with respect to $ Q' $ ($ \theta_{Q'} = \LOI{\theta_Q}{\gamma} $). Set $ \Sfrak'_Q = \gamma^{-1} \Sfrak'_{Q'} \gamma $.
        Let $ K_{Q'} = \gamma K_Q \gamma^{-1} $. Let $ X \in \U(\Ad(\gamma)\n_Q \oplus \Ad(\gamma)\a_Q \oplus \Ad(\gamma)\m_Q) $ and $ Y \in \U(\g) $.
        By \eqref{eq Compare a(g)}, \eqref{eq Compare a_g} and as $ gK_{Q'} \in \Sfrak'_{Q'} \iff g \gamma K_Q \gamma^{-1} \in \gamma \Sfrak'_Q \gamma^{-1} \iff \gamma^{-1} g \gamma K_Q \in \Sfrak'_Q $, $ \LO{p}{\Sfrak_{Q'}}_{r, X, Y}(f) $ is equal to
        \begin{align*}
            & \sup_{gK_{Q'} \in \Sfrak'_{Q'}}{ a_{Q'}(g^{-1})^{\rho_{\Gamma_{Q'}}} (1 + \log a_g)^r a_g^{\rho^{\Gamma_{Q'}}} |L_X R_Y f(g)| } \\
            =& \sup_{g K_Q \in \Sfrak'_Q}{ a_Q(g^{-1})^{\rho_{\Gamma_Q}} (1 + \log a_g)^r a_g^{\rho^{\Gamma_Q}} |L_{\Ad(\gamma^{-1})X} R_{\Ad(\gamma^{-1})Y} (R_{\gamma^{-1}} f)(g)| } \ .
        \end{align*}
        Thus, it follows from Proposition \ref{Prop convergence of Sfrak p_(r, X, Y)(R_h f - f)} that this is finite as $ \Ad(\gamma^{-1})X \in \U(\n_Q) $.
    \end{enumerate}
    The proposition follows.
\end{proof}

\newpage

\subsubsection{Partition of unity indexed over \texorpdfstring{$ \Gamma $}{Gamma}}

In this section, we construct a partition of unity indexed over $ \Gamma $. See Theorem \ref{Thm Cut-off function in the geometrically finite case} on p.\pageref{Thm Cut-off function in the geometrically finite case}.
In the convex-cocompact case, Lemma \ref{Lem Cut-off function in the convex-cocompact case} (lemma providing the $ \chi $-function) provides such a partition of unity.

We need this function in order to prove that the compactly supported smooth functions are dense in the Schwartz space (cf. Proposition \ref{Prop Ccinf(Gamma|G) is dense in the geometrically finite case} on p.\pageref{Prop Ccinf(Gamma|G) is dense in the geometrically finite case}).

In the following, we do some preparations for the proof of this result.

\medskip

\begin{lem}\label{Lem Estimate terms between R and R + 1}
    Let $ g $ be a positive, bounded function on a topological space $Y$ and $ n \in \NN, n' \in \NN_0 $, $ k \in \NN $, $ k' \in \NN $. By convention, $ \RR^0 = \ZZ^0 = \{0\} $.
    Let $ f \colon \RR^{n} \times \RR^{n'} \times Y \to \RR $ be a positive, measurable function such that
    \[
        f((x, x'), y) \asymp 1 + (\|x\|_n^{k} + \|x'\|_{n'}^{k'}) \cdot g(y)
    \]
    when $ x $ varies in $ \RR^n $, $ x' \in \RR^{n'} $ and $ y $ varies in $Y$.
    Here, $ \|\cdot\|_n $ (resp. $ \|\cdot\|_{n'} $) is a norm on $ \RR^{n} $ (resp. $ \RR^{n'} $).
    Let $ R > 0 $ be sufficiently large so that
    \[
        \#\{ (\underline{m}, \underline{m}') \in \ZZ^{n} \times \ZZ^{n'} \mid f((\underline{m}, \underline{m}'), y) \leq R\} > 0 \ .
    \]
    Then,
    \[
        \dfrac{ \#\{ (\underline{m}, \underline{m}') \in \ZZ^{n} \times \ZZ^{n'} \mid f((\underline{m}, \underline{m}'), y) \in [R, R+1] \} }{ \#\{ (\underline{m}, \underline{m}') \in \ZZ^{n} \times \ZZ^{n'} \mid f((\underline{m}, \underline{m}'), y) \leq R\} }
    \]
    is uniformly bounded in $ y \in Y $ and converges uniformly in $ y \in Y $ to zero when $ R $ tends to $ \infty $.
\end{lem}

\begin{proof}
    Let $Y$ be a topological space and let $ f, g $ be as above. Let $ n, n', k, k' $ be as above and let $ y \in Y $.
    Without loss of generality, we may assume that $ \|\cdot\|_n = \|\cdot\|_\infty $ and that $ \|\cdot\|_{n'} = \|\cdot\|_\infty $.
    Let $ d = \frac{n}{k} + \frac{n'}{k'} > 0 $ and let $ \eps > 0 $. Then,
    \begin{multline*}
        \#\{ (\underline{m}, \underline{m}') \in \ZZ^{n} \times \ZZ^{n'}  \mid \|\underline{m}\|^{k} + \|\underline{m}'\|^{k'} \leq r \} \\
        \asymp \#\{ (\underline{m}, \underline{m}') \in \ZZ^{n} \times \ZZ^{n'} \mid \max\{\|\underline{m}\|^{k}, \|\underline{m}'\|^{k'}\} \leq r \}
        \asymp r^d
    \end{multline*}
    when $ r $ varies in $ [\eps, \infty) $.
    Since moreover $ g $ is bounded,
    \[
        \#\{ (\underline{m}, \underline{m}') \in \ZZ^{n} \times \ZZ^{n'} \mid f((\underline{m}, \underline{m}'), y) \leq R\} \asymp \frac{R^{d}}{g(y)^{d}}
    \]
    when $ R > 0 $ is sufficiently large so that the left hand side is positive.
    Hence,
    \[
        \#\{ (\underline{m}, \underline{m}') \in \ZZ^{n} \times \ZZ^{n'} \mid f((\underline{m}, \underline{m}'), y) \in [R, R+1]\} \asymp \frac{\big((R+1)^{d} - R^{d}\big)}{g(y)^{d}} \ .
    \]
    So,
    \[
        \dfrac{ \#\{ (\underline{m}, \underline{m}') \in \ZZ^{n} \times \ZZ^{n'} \mid f((\underline{m}, \underline{m}'), y) \in [R, R+1] \} }{ \#\{ (\underline{m}, \underline{m}') \in \ZZ^{n} \times \ZZ^{n'} \mid f((\underline{m}, \underline{m}'), y) \leq R\} } \asymp \frac{\big((R+1)^{d} - R^{d}\big)}{R^{d}}
    \]
    is uniformly bounded in $ y \in Y $ and converges uniformly in $ y \in Y $ to zero when $ R $ tends to $ \infty $.
\end{proof}

For every multi-index $ \ubar{k} = (k_1, \dotsc, k_l) \in \NN_0^l $, set $ |\ubar{k}| = k_1 + \dotsb + k_l $ and
\[
    D_{\ubar{k}} = \frac{ \del^{|\ubar{k}|} }{ \del x_{k_1} \cdots \del x_{k_l} } \ ,
\]
where $ \del_0 = \Id $ by convention.

\begin{lem}\label{Lem about e^(arccosh(sqrt(x)))}
    For any $ n \in \NN $,
    \[
        \sup_{x \in (1, \infty)} (x-1)^{\frac{2n-1}2} \frac{\del^n}{\del x^n} \big(e^{\arccosh(\sqrt{x})}\big)
    \]
    is finite.
\end{lem}

\begin{proof}
    One can easily prove by induction that
    \[
        \lim_{x \to 1^{+}} (x-1)^{\frac{2n-1}2} \frac{\del^n}{\del x^n} \big(e^{\arccosh(\sqrt{x})}\big)
    \]
    and
    \[
        \lim_{x \to \infty} (x-1)^{\frac{2n-1}2} \frac{\del^n}{\del x^n} \big(e^{\arccosh(\sqrt{x})}\big)
    \]
    are finite for any $ n \in \NN $.
    The lemma follows as
    \[
        x \in (1, \infty) \mapsto (x-1)^{\frac{2n-1}2} \frac{\del^n}{\del x^n} \big(e^{\arccosh(\sqrt{x})}\big)
    \]
    is continuous for any $ n \in \NN $.
\end{proof}

\begin{prop}\label{Prop for derivatives of a_(n_2a)}
    Let $ m = \dim(\g_\alpha) $ and $ n = \dim(\g_{2\alpha}) $. Identify $ \g_\alpha $ with $ \RR^m $ and $ \g_{2\alpha} $ with $ \RR^n $ such that $ |X| = \|X\|_2 $ and $ |Y| = \|Y\|_2 $ for all $ X \in \g_\alpha $ and $ Y \in \g_{2\alpha} $.

    Let $ f(x_1, \dotsc, x_m, y_1, \dotsc, y_n, s) $ be defined by
    \[
        \big(\frac12(e^s + e^{-s}) + \frac1{4e^s}(x_1^2 + \dotsb + x_m^2)\big)^2 + \frac1{2e^{2s}} (y_1^2 + \dotsb + y_n^2) \ .
    \]
    Let $ A(X, Y, s) = a_{\exp(X + Y) a_{e^s}} = e^{\arccosh\bigl(\sqrt{f(X, Y, s)}\bigr)} $.
    Let $ \eps > 0 $. Then, for all $ l, l', l'' \in \NN_0 $,
    \begin{equation}\label{Cor for derivatives of a_(n_2a) eq1}
         \frac{e^{\frac{(l + 2l')s}2} A(X, Y, s)^{\frac{l + 2l'}2}}{ A(X, Y, s) } \cdot |\frac{ \del^{l} }{ \del X_{k_1} \cdots \del X_{k_{l}} }
         \frac{ \del^{l'} }{ \del Y_{k'_1} \cdots \del Y_{k'_{l'}} } \frac{ \del^{l''} }{ \del s^{l''} } A(X, Y, s)|
    \end{equation}
    $ (\ubar{k} \, : \, |\ubar{k}| = l, \ubar{k}' \, : \, |\ubar{k}'| = l' $, $ X, Y, s : A(X, Y, s) > 1 + \eps) $ is uniformly bounded.
\end{prop}

\begin{rem}\label{Rem of Prop for derivatives of a_(n_2a)}
    As $ a_{na} \geq a_a $ for all $ n \in N $ and $ a \in A $, we have the estimate
    \[
        A(X, Y, s) \geq \max\{e^s, e^{-s}\} \geq 1 \qquad (X \in \g_\alpha , Y \in \g_{2\alpha}, s \in \RR) \ .
    \]
    In particular, $ e^{\pm s} A(X, Y, s) \geq 1 $ for all $ X \in \g_\alpha , Y \in \g_{2\alpha}, s \in \RR $.
\end{rem}

\begin{proof}
    Let $ f(x, y, s) = \big(\frac12(e^s + e^{-s}) + \frac{\|x\|^2}{4e^s}\big)^2 + \frac{\|y\|^2}{2e^{2s}} \ $.
    Then,
    \[
        \frac{ \del }{ \del x_j } f(x, y, s) = \big(\frac12(e^s + e^{-s}) + \frac{\|x\|^2}{4e^s}\big) \cdot \frac{x_j}{e^s} \ ,
    \]
    \[
        \frac{ \del }{ \del y_j } f(x, y, s) = e^{-2s} y_j
    \]
    and
    \[
        \frac{ \del }{ \del s } f(x, y, s) = -\big(\frac12(e^s + e^{-s}) + \frac{\|x\|^2}{4e^s}\big) \cdot \frac{\|x\|^2}{2e^s} - \frac{\|y\|^2}{e^{2s}} \ ,
    \]
    Thus,
    \[
        |\frac{ \del }{ \del x_j } f(x, y, s)| \leq \sqrt{f(x, y, s)} \cdot e^{-\frac{s}2} \cdot \frac{|x_j|}{e^{\frac{s}2}}
        \prec e^{-\frac{s}2} f(x, y, s)^{\frac34} \ ,
    \]
    \[
        |\frac{ \del }{ \del y_j } f(x, y, s)| \leq e^{-2s} (1 + |y|) \prec e^{-s} \sqrt{f(x, y, s)}
    \]
    and
    \[
        |\frac{ \del }{ \del s } f(x, y, s)| \prec f(x, y, s) \ .
    \]
    It follows easily from the above, Lemma \ref{Lem about e^(arccosh(sqrt(x)))} and the fact that $ e^{\arccosh(x)} \asymp x $ when $ x $ varies in $ [1, \infty) $ that \eqref{Cor for derivatives of a_(n_2a) eq1} is uniformly bounded for $ l = l' = l'' = 1 $. The general case is proven similarly. For example,
    \begin{align*}
        \frac{ \del^2 }{ \del x_i \del x_j } f(x, y, s)
        &= \frac{ \del }{ \del x_i } \big(\frac12(e^s + e^{-s}) \cdot \frac{x_j}{e^s} + \frac{x_i^2}{4e^s} \cdot \frac{x_j}{e^s}\big) \\
        &= e^{-s} \Big(\delta_{i, j} \big(\frac12(e^s + e^{-s}) + \frac{x_i^2}{4e^s} \big)
        + \frac{x_i}{2e^{\frac{s}2}} \cdot \frac{x_j}{e^{\frac{s}2}}\Big) \ .
    \end{align*}
    Thus,
    \[
        |\frac{ \del^2 }{ \del x_i \del x_j } f(x, y, s)| \prec e^{-s} \sqrt{f(x, y, s)} \ .
    \]
    It follows again easily from this that $ |\frac{ \del^2 }{ \del X_i \del X_j } A(X, Y, s)| $ satisfies the needed estimate.
\end{proof}

\begin{lem}
    Let $ r_1, r_2 \in \RR \cup \{\pm\infty\} $ be such that $ r_1 < r_2 $ (we use the usual convention).
    Then, there is a smooth function $ h_{r_1, r_2} \colon \RR \to [0, 1] $ such that $ h_{r_1, r_2}(x) = 1 $ for $ x \leq r_1 $, $ 0 < h_{r_1, r_2}(x) < 1 $ for $ x \in (r_1, r_2) $ and $ h_{r_1, r_2}(x) = 0 $ for $ x \geq r_2 $.
\end{lem}

\begin{rem} \nlenum
    \begin{enumerate}
    \item The lemma is a slight generalisation of Lemma 2.21 of \cite[p.42]{Lee}.
    \item Let $ 0 < r_1 < r_2 $ ($ r_2 = +\infty $ is also allowed here). Let $ H_{r_1, r_2} \colon \RR^k \to \RR $ ($ k \in \NN $) be defined by $ H_{r_1, r_2}(x) = h_{r_1, r_2}(\|x\|_2) $. Then, this is also a smooth function by the proof of Lemma 2.22 of \cite[p.42]{Lee}.
    \item $ \sup_{x \in \RR} |h_{r_1, r_2}^{(n)}(x)| $ is finite for all $ n \in \NN_0 $.
        \\ Indeed, for $ n = 0 $, this is trivial. Let $ n \in \NN $. As $ h_{r_1, r_2}^{(n)} $ is a continuous function with compact support, $ \sup_{x \in \RR} |h_{r_1, r_2}^{(n)}(x)| $ is finite, too.
    \end{enumerate}
\end{rem}

\begin{proof}
    Define $ f \colon \RR \cup \{\pm \infty\} \to \RR $ by
    \[
        f(x) = \begin{cases}
            1                   & \text{if } x = \infty \\
            e^{-\tfrac1x}       & \text{if } x > 0 \\
            0                   & \text{if } x \leq 0
        \end{cases} \ .
    \]
    This function is smooth on $ \RR $ by Lemma 2.20 of \cite[p.41]{Lee}.

    Let $ r_1, r_2 \in \RR \cup \{\pm \infty\} $ be such that $ r_1 < r_2 $. For $ x \in \RR $, define
    \[
        h_{r_1, r_2}(x) = \frac{ f(r_2 - x) }{ f(r_2 - x) + f(x - r_1) } \ .
    \]
    Then, $ h_{r_1, r_2} \colon \RR \to \RR $ is a well-defined, smooth function such that $ h_{r_1, r_2}(x) = 1 $ for $ x \leq r_1 $, $ 0 < h_{r_1, r_2}(x) < 1 $ for $ x \in (r_1, r_2) $ and $ h_{r_1, r_2}(x) = 0 $ for $ x \geq r_2 $.

\end{proof}

\begin{lem}\label{Lem a a_(na) extends analytically}
    For $ n \in N $ and $ a \in A $, set $ F(naK) = a a_{na} $ and $ F(n w P) = a(n w) $. Then, $F$ is analytic on $ (X \smallsetminus \{eK\}) \cup \{n w P \mid n \in N\} $.
\end{lem}

\begin{proof}
    Let $ n \in N $ and $ a \in A $. Write $ n \in N $ as $ n = \exp(X + Y) $ ($ X \in \g_{\alpha} $, $ Y \in \g_{2\alpha} $) and $ a = a_t $ ($ t > 0 $). Then,
    \[
        a_{n a} = e^{ \arccosh\Bigl(\sqrt{ \big(\tfrac12(t + \tfrac1t) + \tfrac1{4t} |X|^2\big)^2 + \tfrac{|Y|^2}{2t^2}}\Bigr) }
    \]
    by Theorem \ref{Thm explicit formulas}.
    Since $ \arccosh(x) = \log(x + \sqrt{x^2 - 1}) $ for all $ x \geq 1 $, $ a a_{na} $ is again equal to
    \[
        \sqrt{ \big(\tfrac12(t^2 + 1) + \tfrac14 |X|^2\big)^2 + \tfrac{|Y|^2}{2}} + \sqrt{\big(\tfrac12(t^2 + 1)
        + \tfrac14 |X|^2\big)^2 + \tfrac{|Y|^2}{2} - t^2} \ .
    \]
    This function extends analytically to $ (\RR \times N) \smallsetminus \{(1, e)\} $. At $ t = 0 $, it is equal to
    \[
        2\sqrt{ \big(\tfrac12 + \tfrac14 |X|^2\big)^2 + \tfrac{|Y|^2}{2}}
        = \sqrt{ \big(1 + \tfrac12 |X|^2\big)^2 + 2|Y|^2 } \ .
    \]
    By Theorem \ref{Theorem 3.8 of Helgason, p.414}, this is again equal to $ a(\theta n) = a(n w) $. The lemma follows.
\end{proof}

\needspace{11\baselineskip}
Let $ \Delta_X $ be the Laplace-Beltrami operator of $X$.

\begin{thm}\label{Thm Cut-off function in the geometrically finite case}
    There exists a cut-off function $ \chi \in \Cinf(X \cup \Omega, [0, 1]) $ with the following properties:
    \begin{enumerate}
    \item $ \sum_{\gamma \in \Gamma} L_\gamma \chi $ is locally finite;\label{Lem Cut-off function in the geometrically finite case eq1}
    \item $ \sum_{\gamma \in \Gamma} L_\gamma \chi = 1 $;\label{Lem Cut-off function in the geometrically finite case eq2}
    \item $ \sup_{g \in G} |L_X R_Y \chi(g)| < \infty $ for all $ X \in \U(\n_{P_i}) \, (i \in \{1, \dotsc, m\}) $, $ Y \in \U(\g) $;\label{Lem Cut-off function in the geometrically finite case eq3}
    \item $ \sup_{gK \in U} |L_X R_Y \chi(g)| < \infty $ for all $ X, Y \in \U(\g) $ and all $ U \in \U_\Gamma $;\label{Lem Cut-off function in the geometrically finite case eq3b}
    \item $ \sup_{gK \in U} a_g |d\chi(gK)| < \infty $ for every $ U \in \U_\Gamma $;\label{Lem Cut-off function in the geometrically finite case eq3c}
    \item $ \sup_{gK \in U} a_g |\Delta_X \chi(gK)| < \infty $ for every $ U \in \U_\Gamma $;\label{Lem Cut-off function in the geometrically finite case eq3d}
    \item $ \{\gamma \in \Gamma \mid \gamma^{-1} U \cap \supp \chi \neq \emptyset \} $ is finite for every $ U \in \U_\Gamma $;\label{Lem Cut-off function in the geometrically finite case eq4}
    \item $ \{ \gamma\Gamma_{P_i} \in \Gamma/\Gamma_{P_i} \mid \big(\bigcup_{\gamma' \in \Gamma_{P_i}} {\gamma'}^{-1} \gamma^{-1} \Sfrak'_{X, std, i}\big) \cap \supp \chi \neq \emptyset \} $ is finite for all $i$;\label{Lem Cut-off function in the geometrically finite case eq5}
    \item there is a constant $ C > 0 $ such that
        \[
            a_{{\gamma'}^{-1} g} \leq C a_g
        \]
        for all $ gK \in \bigcup_{\gamma \in \Gamma} \gamma \Sfrak'_{X, std, i} \, (i \in \{1, \dotsc, m\}) $, $ \gamma' \in \Gamma_{P_i} $ such that $ {\gamma'}^{-1} g K \in \supp \chi $.\label{Lem Cut-off function in the geometrically finite case eq6}
    \end{enumerate}
    We denote the restriction of $ \chi $ to $ \Omega $ by $ \chi_\infty $.
\end{thm}

\begin{rem}
    Compare this $ \chi $-function with the $ \chi $-function from Lemma \ref{Lem Cut-off function in the convex-cocompact case} (convex-cocompact case).
\end{rem}

\begin{proof}
    Since the theorem follows from Lemma \ref{Lem Cut-off function in the convex-cocompact case} if $ \Gamma \bs X $ has no cusp, we may assume without loss of generality that $ \Gamma \bs X $ has a cusp. Then, \eqref{Lem Cut-off function in the geometrically finite case eq3b} follows from \eqref{Lem Cut-off function in the geometrically finite case eq3} by Lemma \ref{Lem X in U(n) enough}.

    Let us show now that the properties \eqref{Lem Cut-off function in the geometrically finite case eq3c} and \eqref{Lem Cut-off function in the geometrically finite case eq3d} hold if we have constructed a cut-off function $ \chi \in \Cinf(X \cup \Omega, [0, 1]) $ satisfying Property \eqref{Lem Cut-off function in the geometrically finite case eq4}.

    Let $ U \in \U_\Gamma $. Then, $ E := \{\gamma \in \Gamma \mid \gamma^{-1} U \cap \supp \chi \neq \emptyset \} $ is finite by Property \eqref{Lem Cut-off function in the geometrically finite case eq4}. Then, by using a cut-off function which is identically one on $ \bigcup_{\gamma \in E} \gamma^{-1} U $ (relatively compact subset of $ X \cup \Omega $) that $ \chi $ is equal to a cut-off function in $ \Cinf(\bar{X}, [0, 1]) $ on $ \bigcup_{\gamma \in E} \gamma^{-1} U $. The properties \eqref{Lem Cut-off function in the geometrically finite case eq3c} and \eqref{Lem Cut-off function in the geometrically finite case eq3d} follow now from the proof of Lemma \ref{Lem Cut-off function in the convex-cocompact case}.

    We use the following convention: If $ f $ is a function on $ \bar{X} $, then $ f(g) $ ($ g \in G $) is by definition equal to $ f(gK) $ (and in general not equal to $ f(gP) $).

    Let $ \omega_{i, R} = \{ n \in N_{\Gamma, P_i} \mid a_n \leq R \} $ and
    $
        \omega'_{i, R'} = \{ n \in N^{\Gamma, P_i} \mid a_n < R' \}
    $.
    \\ For $ t, s, s' > 0 $, set $ \Sfrak_{i, t, s, s'} = \omega_{i, s} N^{\Gamma, P_i} A_i K_i \smallsetminus \omega_{i, s} \omega'_{i, s'} A_{i, < t} K_i $.

    Recall that $ w_i \in N_{K_i}(\a_{P_i}) $ denotes a representative of the nontrivial Weyl group element.

    Let $ \lambda_i \geq 1 $ be such that $ \frac{ a_{P_i}(n w_i) }{ a_{P_i}(p_{N^{\Gamma, P_i}}(n) w_i) } \leq \lambda_i \cdot \frac{ a_n }{ a_{p_{N^{\Gamma, P_i}}(n)} } $ for all $ n \in N_i $.
    \\ Let $ \lambda = \max_i \lambda_i $.
    Choose $ t > 2, R > 1, R' > 2 $ sufficiently large so that
    \begin{enumerate}
    \item $ \omega_i \subsetneq \omega_{i, \frac{R}{\lambda}} $, $ \omega'_{std, i} \subsetneq \omega'_{i, R'-1} $ for all $ i \in \{1, \dotsc, m\} $ \\
        (recall that $ \Sfrak_{std, i} := \omega_{std, i} N^{\Gamma, P_i} A_i K_i \smallsetminus \omega_{std, i} \omega'_{std, i} (A_i)_{< 1} K_i $);
    \item $ \{\gamma \in \Gamma \mid \gamma \clo(\Sfrak'_{X, i, t - 1, R, R' - 1}) \cap \clo(\Sfrak'_{X, i, t - 1, R, R' - 1}) \neq \emptyset \} $ is finite and contained in $ \Gamma_{P_i} $ for all $i$ and $ \Gamma \clo(\Sfrak'_{X, i, t - 1, R, R' - 1}) \cap \clo(\Sfrak'_{X, j, t - 1, R, R' - 1}) = \emptyset $ for all $ i \neq j $. It follows from the proof of Proposition \ref{Prop Decomp. of Gamma|(X cup Omega) in the geometrically finite case} that this is possible.
    \end{enumerate}
    Then, there is $ U \in \U_\Gamma $ such that
    \begin{equation}\label{eq Sfrak'_i subset Sfrak'_(i, t - 1, R, R' - 1) cup U}
        \Sfrak'_{X, std, i} \subset \Sfrak'_{X, i, t - 1, R, R' - 1} \cup U
    \end{equation}
    for all $i$. By Corollary \ref{Cor Decomp. of Gamma|X in the geometrically finite case}, there is $ V \in \U_\Gamma $ such that
    \begin{equation}\label{Proof of Lem Cut-off function in the geometrically finite case eq1}
        \Gamma \bs X = \Gamma V \cup \bigcup_{i=1}^m \Gamma \Sfrak'_{X, i, t, R, R'} \ .
    \end{equation}
    Let $ W_0 $ be the closure of $ V $ in $ X \cup \Omega $.
    Let $ \tilde{W}_0 $ be an open, relatively compact subset of $X \cup \Omega$ which contains $ W_0 $.

    By Lemma 2.26 of \cite[p.55]{Lee}, there exists a smooth cut-off function $ \psi_0 \in \Cinf(\bar{X}, [0, 1]) $ such that $ \psi_0 $ is identically one on $ W_0 $ and identically zero on $ X \smallsetminus \tilde{W}_0 $. So, $ \supp \psi_0 $ is contained in $ X \cup \Omega $. Thus, $ \sup_{gK \in X} |L_X \psi_0(gK)| $ is finite for all $ X \in \U(\g) $.
    It follows from Lemma \ref{Lem Right derivatives in k-direction suffice when we have all the left derivatives} that $ \sup_{g \in G} |L_X R_Y \psi_0(g)| $ is also finite for all $ X, Y \in \U(\g) $.

    Let $ i \in \{1, \dotsc, m\} $. By the proof of Lemma 2.26 of \cite[p.45]{Lee}, there is a cut-off function $ \varphi_i \in \Cinf(G/K_i \cup G/P_i, [0, 1]) $ such that
    \begin{enumerate}
    \item $ \varphi_i = 0 $ on $ \omega_{std, i} \omega'_{i, R'-1} A_{i, \leq t - 1} K_i $,
    \item $ \varphi_i = 1 $ on $ \omega_{i, R} \omega'_{i, R'} A_{i, \geq t} K_i $, and
    \item $ \varphi_i $ takes values in $ (0, 1) $ otherwise, on $ G/K_i $.
    \end{enumerate}

    \begin{clm}
        $ \prod_{\gamma \in \Gamma_{P_i}} L_\gamma \varphi_i $ is a well-defined function in
        \[
            \Cinf(G/K_i \cup \{ n w_i P_i \mid n \in N_i \}, [0, 1]) \ .
        \]
        Moreover, on $ G/K_i $, it has support on
        \[
            G/K_i \smallsetminus N_{\Gamma, P_i} \omega'_{i, R'-1} A_{i, < t - 1} K_i = \bigcup_{\gamma \in \Gamma_{P_i}} \gamma \Sfrak'_{i, t-1, R, R'-1} \ .
        \]
    \end{clm}

    \begin{proof}
        Note that the last assertion holds.

        Let $ U $ be an open subset of $G/K_i$ which is relatively compact in $ G/K_i \cup \{ n w_i P_i \mid n \in N_i \} $. Then, $ \prod_{\gamma \in \Gamma_{P_i}} \varphi_i(\gamma^{-1} g K_i) = \prod_{\gamma \in \Gamma_{P_i}} \varphi_i(\gamma g K_i)  $ has only finitely many factors which are not equal to 1 when $ gK_i $ varies in $U$. Indeed, if $ \varphi_i(\gamma g K_i) \neq 1 $ for some $ \gamma \in \Gamma_{P_i} $ and some $ gK_i \in U $, then
        \[
            n_\gamma m_\gamma p_{N_{\Gamma, P_i}}(\nu_{P_i}(g)) m_\gamma^{-1} = p_{N_{\Gamma, P_i}}(\nu_{P_i}(\gamma g)) \in \omega_{i, R} \ .
        \]
        But as $ U $ is relatively compact in $ G/K_i \cup \{ n w_i P_i \mid n \in N_i \} $,
        \[
            \{ p_{N_{\Gamma, P_i}}(\nu_{P_i}(h)) \mid hK_i \in U\}
        \]
        is a relatively compact subset of $ N_{\Gamma, P_i} $. Hence, $ \gamma $ is contained in a compact subset of $ N_{\Gamma, P_i} M_{\Gamma, P_i} $ (independently of the choice of $ \gamma $ and $g$). Since moreover $ \Gamma_{P_i} $ is discrete,
        \[
            \{ \gamma \in \Gamma_{P_i} \mid \varphi_i(\gamma g K_i) \neq 1 \text{ for some } gK_i \in U\}
        \]
        is finite.
        Similarly, one shows that $ \prod_{\gamma \in \Gamma_{P_i}} \varphi_i(\gamma^{-1} g P_i) $ has only finitely many factors which are not equal to 1 when $ gP_i $ varies in a relatively compact subset of $ \{ n w_i P_i \mid n \in N_i \} $.
        The claim follows.
    \end{proof}

    For $ n \in N_i $ and $ a \in A_i $, set $ F_i(naK_i) = a a_{na} $ and $ F_i(n w_i P_i) = a_{P_i}(n w_i) $. Then, $F_i$ is analytic on $ (G/K_i \smallsetminus eK_i) \cup \{n w_i P_i \mid n \in N_i\} $ by Lemma \ref{Lem a a_(na) extends analytically}.

    For $ g \in G $, set $ a^{(i)}_g = \restricted{a_h}{ h = p_{N^{\Gamma, P_i}}(\nu_{P_i}(g)) h_{P_i}(g) } $.

    For $ n \in N_i $ and $ a \in A_i $, set $ G_i(naK_i) = a a^{(i)}_{na} $ and $ G_i(n w_i P_i) = a_{P_i}(p_{N^{\Gamma, P_i}}(n) w_i) $. Then, $G_i$ is analytic on $ G/K_i \cup \{n w_i P_i \mid n \in N_i\} $ by the proof of Lemma \ref{Lem a a_(na) extends analytically} and as $ p_{N^{\Gamma, P_i}} $ is analytic.

    Let $ R_i = h_{R, R + 1}\big( \tfrac{ F_i }{ G_i } \big) $. Note that there is a neighbourhood of $ eK_i $ in $ G/K_i $ on which $ R_i $ vanishes. So, $ R_i $ belongs to $ \Cinf(G/K_i \cup \{n w_i P_i \mid n \in N_i\}, [0, 1]) $.

    Define $ \psi_i \in \Cinf(G/K_i \cup \{n w_i P_i \mid n \in N_i\}, [0, 1]) $ by
    \[
        \psi_i = \big(\prod_{\gamma \in \Gamma_{P_i}} L_\gamma \varphi_i \big) h_{R, R + 1}\big( \tfrac{ F_i }{ G_i } \big) \ .
    \]
    Let $ g \in G $. Then,
    \[
        \psi_i(g K_i) = \prod_{\gamma \in \Gamma_{P_i}} \varphi_i(\gamma^{-1} g K_i) h_{R, R + 1}\big( \tfrac{ a_{g} }{ a^{(i)}_g } \big) \in [0, 1] \ .
    \]
    On $ G/K_i $, $ \psi_i $ has support on $ G/K_i \smallsetminus N_{\Gamma, P_i} \omega'_{i, R'-1} A_{i, < t - 1} K_i = \bigcup_{\gamma \in \Gamma_{P_i}} \gamma \Sfrak'_{i, t-1, R, R'-1} $.

    Let $ h_i \in N_i A_i $ be such that $ K_i = h_i K h_i^{-1} $. Then, there is $ k_i \in K $ such that $ P_i = h_i k_i P k_i^{-1} h_i^{-1} $ by Lemma \ref{Lem A_Q = k A k^{-1}}.
    Set $ h_0 = e $, $ k_0 = e $, $ K_0 = K $ and $ P_0 = P $. For $ g \in G $, set
    \[
        \tilde{\psi}_j(gK) = \psi_j(g h_j^{-1} K_j) \quad \text{and} \quad \tilde{\psi}_j(gP) = \psi_j(g k_j^{-1} h_j^{-1} P_j) \ .
    \]
    Let $ \Gamma'_i $ ($ i \in \{1, \dotsc, m\} $) be the subgroup of finite index of $ \Gamma_{P_i} $ which is provided by Auslander's theorem.

    We have:
    \begin{enumerate}
    \item Let $ g \in \bigcup_{\gamma \in \Gamma} \gamma \Sfrak_{i, t, R, R'} $.
        Then, $ \psi_i(\gamma g K_i) = 1 $ for some $ \gamma \in \Gamma $.
        \\ Indeed,
        \begin{multline*}
            \#\{\gamma \in \Gamma \mid \psi_i(\gamma^{-1} g K_i) = 1 \} \geq \#\{\gamma \in \Gamma'_i \mid a_{p_{N_{\Gamma, P_i}}(\nu_{P_i}(\gamma^{-1} g))} \leq R \} \\
            = \#\{\gamma \in \Gamma'_i \mid p_{N_{\Gamma, P_i}}(\nu_{P_i}(\gamma^{-1} g)) \in \omega_{i, R} \}
        \end{multline*}
        as $ \tfrac{ a_{g} }{ a^{(i)}_g } \leq a_{p_{N_{\Gamma, P_i}}(\nu_{P_i}(g))} $. This is greater or equal than 1 since
        \[
            a_{p_{N_{\Gamma, P_i}}(\nu_{P_i}(\gamma^{-1} g))} = a_{n_\gamma^{-1} p_{N_{\Gamma, P_i}}(\nu_{P_i}(g))}
        \]
        for all $ \gamma \in \Gamma_{P_i} $ and since $ \bigcup_{\gamma \in \Gamma'_i} n_\gamma \omega_{std, i} = N_{\Gamma, P_i} $.
    \item Let $ n w_i P_i \in \bigcup_{\gamma \in \Gamma} \gamma \bigl(\clo(\Sfrak'_{i, t, R, R'}) \cap \{gP_i \mid g h_i k_i P \in \Omega_\Gamma\}\bigr) $.
        \\ Then, $ \psi_i(\gamma n w_i P_i) = 1 $ for some $ \gamma \in \Gamma $.
        \\ Indeed, $ \psi_i(n w_i P_i) = \frac{ a_{P_i}(n w_i) }{ a_{P_i}(p_{N^{\Gamma, P_i}}(n) w_i) } \leq \lambda \cdot \frac{ a_n }{ a_{p_{N^{\Gamma, P_i}}(n)} } \leq \lambda a_{p_{N_{\Gamma, P_i}}(n)} $.
        \\ Now we can show the rest similarly as above (use that $ \omega_{std, i} $ is contained in $ \omega_{i, \frac{R}{\lambda}} $).
    \item $ \supp(\tilde{\psi}_j) \cap \supp(\tilde{\psi}_k) = \emptyset $ for all $ j, k $ such that $ 1 \leq j < k \leq m $.
    \item $ \{ \gamma \in \Gamma \mid \supp(\tilde{\psi}_i) \cap \gamma \clo(\Sfrak'_{X, i, t-1, R, R'-1}) \neq \emptyset \} \subset \Gamma_{P_i} $: \\
        As $ \supp(\tilde{\psi}_i) \subset \Gamma_{P_i} \clo(\Sfrak'_{i, t-1, R, R'-1}) $,
        \[
            \{ \gamma \in \Gamma \mid \supp(\tilde{\psi}_i) \cap \gamma \clo(\Sfrak'_{X, i, t-1, R, R'-1}) \neq \emptyset \} \subset \Gamma_{P_i} \ .
        \]
    \item $ \{ \gamma \in \Gamma \mid \supp(\tilde{\psi}_i) \cap \gamma \supp(\tilde{\psi}_i) \neq \emptyset \} \subset \Gamma_{P_i} $.
    \end{enumerate}
    Let $ g \in G $ be such that $ \psi_0(g) \neq 0 $. By the above, there is at most one $ j \in \{1, \dotsc, m\} $ such that $ \tilde{\psi}_j(g) \neq 0 $.

    Set $ \psi = \sum_{j=0}^m \tilde{\psi}_j $. Then, $ \psi \in \Cinf(X \cup (\dX \smallsetminus \bigcup_{j=1}^m \{eP_j\}), [0, 1]) $.

    Define $ \chi = \frac{\psi}{\sum_{\gamma \in \Gamma} L_\gamma \psi} $ on $ X \cup \Omega $. It follows from \eqref{Proof of Lem Cut-off function in the geometrically finite case eq1}, the construction of our maps, Proposition \ref{Prop [gamma in Gamma | gamma (U cap cup_i Sfrak'_i) cap (U cap cup_i Sfrak'_i) neq emptyset]}, Proposition \ref{Prop Decomp. of Gamma|(X cup Omega) in the geometrically finite case} and \eqref{eq Sfrak'_i subset Sfrak'_(i, t - 1, R, R' - 1) cup U} that $ \chi $ is well-defined, that the first two and the last three properties hold. Moreover, $ \chi \in \Cinf(X, [0, 1]) $.
    
    \smallskip

    Let us prove now that \eqref{Lem Cut-off function in the geometrically finite case eq3} holds.

    \begin{enumerate}
    \item Let $ \eps > 0 $. It follows from Proposition \ref{Prop L_X R_Y a_g} that for all $ X, Y \in \U(g) $, there is $ c_{X, Y} > 0 $ such that
        \begin{equation}\label{eq in proof |L_X R_Y a_g| <= c_(X, Y) a_g (g in G : a_g > 1 + eps)}
             |L_X R_Y a_{g}| \leq c_{X, Y} a_{g} \qquad (g \in G : a_g > 1 + \eps) \ .
        \end{equation}
    \item If $ g \in G $ belongs to
        \begin{equation}\label{eq Domain for g}
            D_i := \{g \in G \mid \tfrac{ a_{g} }{ a^{(i)}_g} \leq R + 1 \}
        \end{equation}
        (e.g. if $ g \in \Sfrak_{i, t, R, R'} $) then it follows from the proof of Proposition \ref{Prop a_(n_1 n_2 a) >= c a_(n_2 a)} that there is a constant $ c_R > 0 $ such that
        \begin{equation}\label{eq Estimate of |log p_(N_1)(p_(N_(Gamma, P_i))(n_(P_i)(g)))|}
            1 + |\log p_1\big(p_{N_{\Gamma, P_i}}(\nu_{P_i}(g))\big)|
            \leq c_R \sqrt{h_{P_i}(g)} a^{\frac12}_{p_{N^{\Gamma, P_i}}(\nu_{P_i}(g)) h_{P_i}(g)} \ .
        \end{equation}
    \item For $ k \in \NN_0 $, set $ \U^k(\n_{P_i}) = \{ X \in \U(\n_{P_i}) \mid \ad(H)X = k X \quad \forall H \in \a_{P_i} \} $. Then,
        \[
            \U(\n_{P_i}) = \bigoplus_{k=0}^\infty \U^k(\n_{P_i}) \ .
        \]
        \begin{clm}
            Let $ l \in \NN_0 $. For any $ X \in \U^l(\n_{P_i}) $, $ Y \in \U(\g) $ and any $ \eps > 0 $,
            \[
                \underbrace{ h_{P_i}(g)^{\frac{l}2} (a^{(i)}_g)^{\frac{l}2} }_{ \geq 1 }
                \cdot \frac{ |L_X R_Y a^{(i)}_g| }{ a^{(i)}_g }
            \]
            is uniformly bounded when $ g $ varies in
            \begin{equation}\label{eq domain for g}
                \{ g \in G \mid a_g > 1 + \eps , \ \tfrac{ a_{g} }{ a^{(i)}_g} \leq R + 1 \} \ .
            \end{equation}
        \end{clm}

        \begin{proof}
            Without loss of generality, we may assume that $ k_i(g) = e $.

            For $ X \in \g_{\alpha_i} $, $ Y \in \g_{2\alpha_i} $ and $ s \in \RR $, set $ A(X, Y, s) = a_{\exp(X + Y) a_{e^s}} $.  

            Let now $ X \in \U^l(\n_{P_i}) $, $ Y \in \U(\g) $.

            Since the function is right $K_i$-invariant, we may also assume without loss of generality that $ Y \in \U(\n_{P_i}) \U(\a_{P_i}) $.
            Let $ Z \in \U^{l'}(\n_{P_i}) $ and $ H \in \U(\a_{P_i}) $.

            Note that there are finitely many $ X_{j, g} \in \U^l(\n_{P_i}) $, $ Z_{k, g} \in \U^{l'}(\n_{P_i}) $ such that
            \begin{multline*}
                L_X R_Z R_H a^{(i)}_g \\
                = \sum_{j, k} D_{X_{j, g}} D_{Z_{k, g}} D_H A(\log p_{N_1^{\Gamma, P_i}}(\nu_{P_i}(g)),
                    \log p_{N_2^{\Gamma, P_i}}(\nu_{P_i}(g)), \log h_{P_i}(g))
            \end{multline*}
            (the notations are self-explaining) with
            \[
                |X_{j, g}| \leq c_1 (1 + |\log p_1(p_{N_{\Gamma, P_i}}(g))|)^{l}
                + c_2 (1 + |\log p_{N_1^{\Gamma, P_i}}(g)|)^{l}
            \]
            and
            \[
                |Z_{k, g}| \leq c_1 h_{P_i}(g)^{l'} (1 + |\log p_1(p_{N_{\Gamma, P_i}}(g))|)^{l'} + c_2 h_{P_i}(g)^{l'} (1 + |\log p_{N_1^{\Gamma, P_i}}(g)|)^{l'}
            \]
            for some constants $ c_1, c_2 > 0 $.

            Since $ a(\bar{n}) \leq a_{na} a $ for all $ n \in N_i $ and $ a \in A_i $,
            \[
                1 + |\log p_{N_1^{\Gamma, P_i}}(g)| \prec a^{\frac12}_{p_{N^{\Gamma, P_i}}(g)} \leq (a^{(i)}_g)^{\frac12} h_{P_i}(g)^{\frac12} \ .
            \]
            Thus,
            \[
                |X_{j, g}| \prec h_{P_i}(g)^{\frac{l}2} (a^{(i)}_g)^{\frac{l}2}
            \]
            and
            \[
                |Z_{k, g}| \prec h_{P_i}(g)^{\frac{3 l'}2} (a^{(i)}_g)^{\frac{l'}2}
                \prec h_{P_i}(g)^{l'} (a^{(i)}_g)^{l'} \ .
            \]
            The claim follows now from Proposition \ref{Prop for derivatives of a_(n_2a)} and Remark \ref{Rem of Prop for derivatives of a_(n_2a)}.
        \end{proof}
    \item We can now easily conclude from the above that for all $ X \in \U(\n_{P_i}) $, $ Y \in \U(\g) $, there is $ c'_{X, Y} > 0 $ such that 
        \[
             |L_X R_Y \tfrac{ a_{g} }{ a^{(i)}_g }| \leq c'_{X, Y} \cdot \tfrac{ a_{g} }{ a^{(i)}_g }
        \]
        when $ g $ varies in \eqref{eq domain for g}.
    \end{enumerate}

    Since moreover $ \sup_{x \in \RR} |h_{r_1, r_2}^{(n)}(x)| $ is finite for all $ n \in \NN_0 $, it follows from the Leibniz rule and the chain rule that $ \sup_{g \in G} |L_X R_Y \psi_i(g)| $ is finite for all $ X \in \U(\n_{P_i}) $ and $ Y \in \U(\g) $.
    Let $ l \in \NN_0 $. Then,
    \begin{equation}\label{eq Estimate for L_X R_Y psi_i(g)}
        \sup_{g \in G} h_{P_i}(g)^{\frac{l}2} (a^{(i)}_g)^{\frac{l}2} |L_X R_Y \psi_i(g)|
    \end{equation}
    is even finite for any $ X \in \U^l(\n_{P_i}) $, $ Y \in \U(\g) $.

    \medskip

    By the proof of Proposition \ref{Prop a_(n_1 n_2 a) >= c a_(n_2 a)}, Lemma \ref{Lem Estimate terms between R and R + 1} and as $ \Gamma'_i $ has finite index in $ \Gamma_{P_i} $,
    \begin{equation}\label{eq estimate of cardinalities}
        \#\{\gamma \in \Gamma_{P_i} \mid \psi_i(\gamma^{-1} g) \in (0, 1) \} \leq \#\{\gamma \in \Gamma_{P_i} \mid \psi_i(\gamma^{-1} g) = 1 \} \quad (g \in G) \ .
    \end{equation}
    Let $ g \in G $. If $ X \in \g \U(\g) $ and if $ \psi_i(g) \in \{0, 1\} $, then $ L_X \psi_i(g) = 0 $ by construction of $ \psi_i $.
    Thus, $ L_X \chi(g) = 0 $ for all $ X \in \U(\g) $ if $ \psi(g) = 0 $.
    
    Let $ X \in \g_\CC $. Then,
    \begin{align*}
        L_X \chi(g)
        &= \restricted{ \frac{d}{dt} }{ t = 0 } \frac{\psi(\exp(-tY)g)}{\sum_{\gamma \in \Gamma} \psi(\gamma^{-1}\exp(-tY)g)} \\
        &= \frac{ L_X \psi(g) \big(\sum_{\gamma \in \Gamma} \psi(\gamma^{-1}g)\big)
        - \sum_{\gamma \in \Gamma} L_X (L_\gamma \psi)(g) \cdot \psi(g) }{ \big(\sum_{\gamma \in \Gamma} \psi(\gamma^{-1}g)\big)^2 } \\
        &= \frac{ L_X \psi(g) }{ \sum_{\gamma \in \Gamma} \psi(\gamma^{-1}g) }
        - \frac{ \sum_{\gamma \in \Gamma} L_X (L_\gamma \psi)(g) }{ \sum_{\gamma \in \Gamma} \psi(\gamma^{-1}g) } \cdot \chi(g) \ .
    \end{align*}
    Let $ X_1, X_2 \in \g_\CC $. Then, $ L_{X_2} L_{X_1} \chi(g) $ is equal to
    \begin{align*}
        & \restricted{ \frac{d}{dt} }{ t = 0 } \Big(\frac{ L_{X_1} \psi(\exp(-t X_2) g) }{ \sum_{\gamma \in \Gamma} \psi(\gamma^{-1} \exp(-t X_2)g) }
        - \frac{ \sum_{\gamma \in \Gamma} L_{X_1} (L_\gamma \psi)(\exp(-t X_2)g) }{ \sum_{\gamma \in \Gamma} \psi(\gamma^{-1} \exp(-t X_2)g) } \cdot \chi(\exp(-t X_2)g)\Big) \\
        =& \frac{ L_{X_2} L_{X_1} \psi(g) }{ \sum_{\gamma \in \Gamma} \psi(\gamma^{-1} g) }
        - \frac{ L_{X_1} \psi(g) \big(\sum_{\gamma \in \Gamma} L_{X_2} (L_\gamma \psi)(g)\big) }{ \big(\sum_{\gamma \in \Gamma} \psi(\gamma^{-1} g)\big)^2 } \\
        & \; \Big(- \frac{ \sum_{\gamma \in \Gamma} L_{X_2} L_{X_1} (L_\gamma \psi)(g) }{\sum_{\gamma \in \Gamma} \psi(\gamma^{-1} g) } 
        + \frac{ \big(\sum_{\gamma \in \Gamma} L_{X_1} (L_\gamma \psi)(g)\big) \big(\sum_{\gamma \in \Gamma} L_{X_2} (L_\gamma \psi)(g)\big)}{ \big(\sum_{\gamma \in \Gamma} \psi(\gamma^{-1} g)\big)^2 }
        \Big) \cdot \chi(g) \\
        & \; - \frac{ \sum_{\gamma \in \Gamma} L_{X_1} (L_\gamma \psi)(g) }{ \sum_{\gamma \in \Gamma} \psi(\gamma^{-1}g) } \cdot L_{X_2} \chi(g) \ .
    \end{align*}
    By induction, one can show that the appearing factors in the terms of $ L_X R_Y \chi(g) $ ($ X, Y \in \g \U(\g) $) have one of the following forms:
    \[
        L_{X'} R_{Y'} \chi(g) \quad \text{or} \quad \frac{ L_{X'} R_{Y'} \psi(g) }{ \sum_{\gamma \in \Gamma} \psi(\gamma^{-1}g) }
        \quad \text{or} \quad \frac{ \sum_{\gamma \in \Gamma} L_{X'} R_{Y'} (L_\gamma \psi)(g) }{ \sum_{\gamma \in \Gamma} \psi(\gamma^{-1}g) }
    \]
    ($ X', Y' \in \g \U(\g) $).

    Let $ I_g $ be the subset of $ \{0, \dotsc, m\} $ such that $ j \in I_g $ if and only if $ \psi_j(\gamma^{-1}g) \neq 0 $ for some $ \gamma \in \Gamma $.
    Recall that, by construction, $ \#I_g \in \{1, 2\} $.

    Let $ N_{i, gK} = \#\{\gamma \in \Gamma \mid \psi_i(\gamma^{-1}g h_i^{-1}) = 1 \} $ and let $ N'_{i, gK} = \#\{\gamma \in \Gamma \mid \psi_i(\gamma^{-1}g h_i^{-1}) \in (0, 1) \} $.
    If $ j \in I_g $, then $ N_{j, gK} \geq 1 $.

    We have the following estimates:
    \[
        \frac{ |L_{X'} R_{Y'} \psi_i(g)| }{ \sum_{\gamma \in \Gamma} \psi_i(\gamma^{-1}g) } \leq |L_{X'} R_{Y'} \psi_i(g)|
    \]
    and
    \begin{multline*}
        \frac{ \sum_{\gamma \in \Gamma} |L_{X'} R_{Y'} (L_\gamma \psi)(g)| }{ \sum_{\gamma \in \Gamma} \psi(\gamma^{-1}g) }
        \leq \sum_{j \in I_g} \frac{ \sum_{\gamma \in \Gamma} |L_{X'} R_{Y'} (L_\gamma \psi_j)(g h_j^{-1})| }{ \sum_{\gamma \in \Gamma} \psi_j(\gamma^{-1}g h_j^{-1}) } \\
        \leq 2 \sup_{j=0, \dotsc, m : j \in I_g} \frac{ \sum_{\gamma \in \Gamma} |L_{X'} R_{Y'} (L_\gamma \psi_j)(g h_j^{-1})| }{ \sum_{\gamma \in \Gamma} \psi_j(\gamma^{-1}g h_j^{-1}) }\ .
    \end{multline*}
    Thus, in order to prove \eqref{Lem Cut-off function in the geometrically finite case eq3}, it remains to prove the following assertion:

    \begin{clm}\label{Clm varphi_i has finite seminorms}
        For $ X \in \U(\n_{P_j}) $ if $ j \in \{1, \dotsc, m\} $ (resp. $ X \in \U(\g) $ if $ j = 0 $), $ Y \in \U(\g) $, we have
        \[
            \sup_{gK \in X \, : \, N_{j, gK} \geq 1} \frac{ \sum_{\gamma \in \Gamma} |L_{X} R_{Y} (L_\gamma \psi_j)(g h_j^{-1})| }{ \sum_{\gamma \in \Gamma} \psi_j(\gamma^{-1}g h_j^{-1}) } < \infty \ .
        \]
    \end{clm}

    \begin{proof}
        Let $ gK \in X $ be such that $ N_{j, gK} \geq 1 $. Then,
        $
            \frac{ \sum_{\gamma \in \Gamma} |L_{X} R_{Y} (L_\gamma \psi_j)(g h_j^{-1})| }{ \sum_{\gamma \in \Gamma} \psi_j(\gamma^{-1}g h_j^{-1}) }
        $
        is less or equal than
        \[
            \frac{ N'_{j, gK} }{ N_{j, gK} } \cdot \sup_{\gamma \in \Gamma \, : \, \psi_j(\gamma^{-1}g h_j^{-1}) \in (0, 1) }
            |L_{\Ad(\gamma^{-1})X} R_{Y} \psi_j(\gamma^{-1} g h_j^{-1})| \ .
        \]
        Since $ N'_{j, gK} \leq N_{j, gK} $ for all $ j \geq 1 $ by \eqref{eq estimate of cardinalities}, and as $ \sup_{g \in G} N'_{0, gK} $ is finite, there is a constant $ c > 0 $ such that this is again less or equal than
        \[
            c \sup_{\gamma \in \Gamma \, : \, \psi_j(\gamma^{-1} g h_j^{-1}) \in (0, 1) } |L_{\Ad(\gamma^{-1})X} R_{Y} \psi_j(\gamma^{-1} g h_j^{-1})| \ .
        \]
        Recall that $ \{ \gamma \in \Gamma \mid \supp(\tilde{\psi}_i) \cap \gamma \Sfrak'_{X, i, t-1, R, R'-1} \neq \emptyset \} $ is contained in $ \Gamma_{P_i} $ for all $ i \in \{1, \dotsc, m\} $. Thus, we may assume without loss of generality that
        \[
            g K \in \bigcup_{\gamma \in \Gamma} \gamma V \cup \bigcup_{\gamma \in \Gamma_j} \gamma \Sfrak'_{X, j, t, R, R'} \ .
        \]
        Let
        $
            S = \{ \gamma \in \Gamma \mid \gamma V \cap (V \cup \Sfrak'_{X, i, t, R, R'}) \neq \emptyset \quad \text{for some } i \in \{1, \dotsc, m\}\}
        $
        (finite set) and let $ V' = \bigcup_{\gamma \in S} \gamma V \in \U_\Gamma $. Then, $ \sup_{gK \in V'} (N_{j, gK} + N'_{j, gK}) $ is finite.

        Since moreover $ \sup_{g \in G} |L_X R_Y \psi_j(g)| $ ($ j \in \{0, \dotsc, m\} $) is finite for all $ X \in \U(\n_{P_j}) $ and $ Y \in \U(\g) $, it remains to show that
        \[
            \sup_{\gamma' \in \Gamma_{P_i}} \sup_{gK \in \Sfrak'_{X, i, t, R, R'}, \gamma \in \Gamma_{P_i} \, : \, \psi_i(\gamma^{-1} \gamma' g h_i^{-1}) \in (0, 1) } |L_{\Ad(\gamma^{-1})X} R_{Y} \psi_i(\gamma^{-1} \gamma' g h_i^{-1})|
        \]
        is finite for all $ X \in \U(\n_{P_i}) $ ($ i \in \{1, \dotsc, m\} $) and $ Y \in \U(\g) $.
        It follows from \eqref{eq Estimate of |log p_(N_1)(p_(N_(Gamma, P_i))(n_(P_i)(g)))|} that
        \begin{equation}\label{eq Estimate of |log p_(N_1)(gamma)|}
            1 + |\log p_1(\nu_{P_i}(\gamma))| \prec \sqrt{h_{P_i}(g h_i^{-1})} a^{\frac12}_{p_{N^{\Gamma, P_i}}(\nu_{P_i}(g h_i^{-1})) h_{P_i}(g h_i^{-1})}
        \end{equation}
        for all $ \gamma \in \Gamma_{P_i} $ and all $ g \in D_i $ (in particular for all $ gK \in \Sfrak'_{X, i, t, R, R'} $). The claim follows now from Lemma \ref{Lem Ad(gamma)X - X}, \eqref{eq Estimate of |log p_(N_1)(gamma)|}, \eqref{eq Estimate for L_X R_Y psi_i(g)} and the fact that the right-hand side of \eqref{eq Estimate of |log p_(N_1)(gamma)|} is left $ \Gamma_{P_i} $-invariant. 
    \end{proof}
    This completes the proof of the theorem.
\end{proof}

\needspace{5\baselineskip}
We also have a partition of unity of $ \Gamma \bs X $ having nice properties:
\begin{prop}\label{Prop f = sum_(i=0)^m varphi_i f for all f in CS(Gamma|G)}
    There exist $ \varphi_i \in \Cinf(\Gamma \bs (X \cup \Omega_\Gamma), [0, 1]) $ $ (i \in \{0, \dotsc, m\}) $ and $ U \in \U_\Gamma $ such that
    \begin{enumerate}
    \item $ \sum_{i=0}^m \varphi_i = 1 $,
    \item $ \Gamma gK = \Gamma hK $ implies $ \Gamma_{P_i} gK = \Gamma_{P_i} hK $ on $ \supp(\varphi_i) $ $ (i \in \{1, \dotsc, m\}) $,
    \item $ \supp(\varphi_0) \subset \Gamma U $, $ \supp(\varphi_i) \subset \Gamma \Sfrak'_{X, std, i} $ $ (i \in \{1, \dotsc, m\}) $,
    \item $ \sup_{gK \in X} |L_X R_Y \varphi_0(g)| < \infty $ for all $ X, Y \in \U(\g) $,
    \item $ \sup_{gK \in X} |L_X R_Y \varphi_i(g)| < \infty $ for all $ X \in \U(\n_{P_i}) $ $ (i \in \{1, \dotsc, m\}) $ and $ Y \in \U(\g) $.
    \end{enumerate}
    Thus, $ \varphi_0 f \in \CS(\Gamma \bs G, \phi) $ and $ \varphi_i f \in \CS(\Gamma_{P_i} \bs G, \phi) $ $ (i \in \{1, \dotsc, m\}) $ for all $ f \in \CS(\Gamma \bs G, \phi) $.
\end{prop}

\begin{rem}\label{Rem of Prop f = sum_(i=0)^m varphi_i f for all f in CS(Gamma|G)} \nlenum
    \begin{enumerate}
    \item Let $ f \in \CS(\Gamma \bs G, \phi) $. Then, $ f = \sum_{i=0}^m \varphi_i f $ with $ \varphi_0 f \in \CS(\Gamma \bs G, \phi) $ and $ \varphi_i f \in \CS(\Gamma_{P_i} \bs G, \phi) $ $ (i \in \{1, \dotsc, m\}) $.
    \item The proposition can be used to construct Schwartz functions:
        \begin{enumerate}
        \item Let $ f \in \CS(G, V_\phi) $ be such that $ \supp(f) $ is contained in some $ U \in \U_\Gamma $ (e.g. $ f = \frac{\psi_0}{\sum_\gamma L_\gamma \psi} \cdot g $ for some $ g \in \CS(G, V_\phi) $). Then,
            \[
                \sum_{\gamma \in \Gamma} \phi(\gamma) L_\gamma f \in \CS(\Gamma \bs G, \phi) \ .
            \]
        \item Let $ f_i \in \CS(\Gamma_{P_i} \bs G, \phi) $ $ (i \in \{1, \dotsc, m\}) $. Then, $ \sum_{i=1}^m \varphi_i f \in \CS(\Gamma \bs G, \phi) $.
        \end{enumerate}
    \end{enumerate}
\end{rem}

\begin{proof}
    We use the notation of the proof of the previous theorem. For $ g \in G $, set $ \varphi_i = \frac{\sum_{\gamma \in \Gamma} L_\gamma \tilde{\psi}_i}{\sum_{\gamma \in \Gamma} L_\gamma \psi} $.
    Then, the proposition follows from the proof of the mentioned theorem. The two last properties follow in particular from Claim \ref{Clm varphi_i has finite seminorms}.
\end{proof}

\subsubsection{\texorpdfstring{Density of $ \Ccinf(\Gamma \bs G, \phi) $ in $ \CS(\Gamma \bs G, \phi) $}{Density of the compactly supported smooth functions}}

\begin{prop}\label{Prop Ccinf(Gamma|G) is dense in the geometrically finite case}
    The inclusion of $ \Ccinf(\Gamma \bs G, \phi) $ into $ \CS(\Gamma \bs G, \phi) $ is continuous with dense image.
\end{prop}

\begin{proof}
    The proof of this proposition is similar to the one of Proposition \ref{Prop Ccinf(Gamma|G) is dense} (use Theorem \ref{Thm Cut-off function in the geometrically finite case} instead of Lemma \ref{Lem Cut-off function in the convex-cocompact case}).
\end{proof}

\newpage

\subsubsection{Examples of noncompactly supported Schwartz functions}

Let $ G = \SO(1, n)^0 $, $ n \geq 2 $, and let $ M $, $ A $, $ N $, $K$ be the subgroups of $G$ induced by the corresponding one of $ \Sp(1, n) $ computed in Appendix \ref{sec:Sp(1,n)} (then $ N \subset \exp(\g_\alpha) $). Let $ P = MAN $. Assume that $ \Gamma = \Gamma_P $.

We have: $ \frac{k}2 = \rho_{\Gamma} $ and $ \frac{n - k - 1}2 = \rho^{\Gamma} $.

Let $ \phi $ be a $ \Gamma $-invariant smooth function on $ N_\Gamma $. Set $ V_\phi = \CC $.

We denote by $ \S(\SO(1, n-k)^0) $ the space of rapidly decreasing functions on $ \SO(1, n-k)^0 $ (cf. \cite[p.230]{Wallach}).

Fix $ \gamma \in \hat{K} $. For $ f \in \S(\SO(1, n-k)^0) $, $ n_1 \in N_\Gamma $, $ n_2 \in N^\Gamma $, $ a \in A $ and $ h \in K $, define
\[
    F(n_1 n_2 a h) = \phi(n_1) \gamma(h)^{-1} f(n_2 a) \ .
\]

\begin{clm}
    The function $ F $ belongs to $ \CS(Y, V_Y(\gamma)) $, where $ Y := \Gamma \bs X $.
\end{clm}

\begin{proof}
    Let $ g \in G $, let $ \eps > 0 $, $ d \geq \frac{n - 1}2 + \eps $ and $ r \geq 0 $. As $ a_{na}^\alpha \geq a_a^\alpha = \max\{a^\alpha, a^{-\alpha}\} $ for all $ n \in N $ and $ a \in A $ by Lemma \ref{Lem a_(na) estimates},
    \[
        |F(g)| = |f(p_{N^\Gamma}(\nu(g)) h(g))|
    \]
    is less or equal than
    \[
        c_d \, h(g)^{\frac{k}2} a_{p_{N^\Gamma}(\nu(g)) h(g)}^{-\frac{n - k - 1}2 - \eps}
    \]
    for some constant $ c_d > 0 $. Let $ \Sfrak $ be a generalised standard Siegel set with respect to $P$. Then, $ \LO{p}{\Sfrak'}_{r, 1, 1}(f) $ is finite.

    Let $ X \in \U(\n) = \U(\n_\Gamma)\U(\n^\Gamma) $ and $ Y \in \U(\g) $. Without loss of generality, we may assume that $ k(g) $ is trivial and that $ Y \in \U(\n_\Gamma)\U(\n^\Gamma)\U(\a) $.

    We have
    \begin{enumerate}
    \item $ R_H \big(h(g)^m\big) = m \alpha(H) h(g)^m $ for all $ H \in \a $ and $ m \in \ZZ $;\label{eq Property1 for F in C(Y, V_Y(gamma))}
    \item $ \sup_{n \in N_\Gamma} |R_Z \phi(n)| < \infty $ for all $ Z \in \U(\n_\Gamma) $;\label{eq Property2 for F in C(Y, V_Y(gamma))}
    \item $ L_{Y_1} R_{Y_2} f(p_{N^\Gamma}(\nu(\cdot)) h(\cdot))(g) = 0 $ for all $ Y_1, Y_2 \in \U(\n_\Gamma) $.
    \end{enumerate}
    Let $ \l $ be a Lie algebra over $ \RR $. For $ m \in \NN_0 $, set $ \U^m(\l) = \{ X \in \U(\l) \mid \ad(H)X = m X \} $. Then,
    \[
        \U(\l) = \bigoplus_{m=0}^\infty \U^m(\l) \ .
    \]
    Let $ X_1 \in \U(\n_\Gamma) $, $ Z_1 \in \U^l(\n_\Gamma) $, $ X_2, Z_2 \in \U(\n^\Gamma) $ and $ H \in \U(\a) $.
    Since
    \[
        n_1 n_2 a n'_1 n'_2 = n_1 n_2 (a n'_1 a^{-1}) (a n'_2 a^{-1}) a
        = \underbrace{n_1 (a n'_1 a^{-1})}_{\in N_\Gamma} \underbrace{n_2 (a n'_2 a^{-1})}_{\in N^\Gamma} a
    \]
    ($ n_1, n'_1 \in N_\Gamma, n_2, n'_2 \in N^\Gamma $),
    $ L_{X_1} L_{X_2} R_{Z_1} R_{Z_2} F(g) $ is equal to
    \begin{equation}\label{eq L_(X_1) L_(X_2) R_(Z_1) R_(Z_2) F(g)}
         h(g)^{-l} \big(L_{X_1} R_{Z_1}
         \phi(p_{N_\Gamma}(\nu(g)))\big) L_{X_2} R_{Z_2} f(p_{N^\Gamma}(\nu(g)) h(g)) \ .
    \end{equation}
    Let now $ d \geq \frac{n - 1}2 + l + \eps $. As moreover $ f \in \S(\SO(1, n-k)^0) $, it follows from \eqref{eq Property1 for F in C(Y, V_Y(gamma))}, \eqref{eq Property2 for F in C(Y, V_Y(gamma))} and \eqref{eq L_(X_1) L_(X_2) R_(Z_1) R_(Z_2) F(g)} that $ |L_{X_1} L_{X_2} R_{Z_1} R_{Z_2} R_H F(g)| $ is less or equal than
    \[
        c_d \, h(g)^{-l} a_{ p_{N^\Gamma}(\nu(g)) h(g) }^{-d}
    \]
    for some constant $ c_d > 0 $. This is again less or equal than
    \[
        c_d \, h(g)^{\frac{k}2} a_{p_{N^\Gamma}(\nu(g)) h(g)}^{-\frac{n - k - 1}2 - \eps}
    \]
    as $ a_{na} \geq a_a $ for all $ n \in N $ and $ a \in A $. The claim follows.
\end{proof}

\newpage

\subsubsection{Examples of cusp forms}\label{ssec:Examples of cusp forms}

In this section, we determine the cusp forms on $ \Gamma \bs G $ which are induced from cusp forms on $G$.

\begin{lem}\label{Lemma about sum_(gamma in Gamma) a_(gamma^(-1) g)^(-(lambda + rho))}
    Let $ U \in \U_\Gamma $. If $ \lambda > \delta_\Gamma $, then there is $ c > 0 $ such that
    \[
        \sum_{\gamma \in \Gamma} a_{\gamma^{-1} g}^{-(\lambda + \rho)} \leq c \, a_g^{-(\lambda + \rho)}
    \]
    for all $ gK \in U $.
\end{lem}

\begin{proof}
    Let $ U \in \U_\Gamma $. Assume that $ \lambda > \delta_\Gamma $.
    Since the assertion of the lemma is trivial when $ U $ is relatively compact in $X$, we may assume without loss of generality that there is an open relatively compact subset $ V $ of $ \Omega $ such that
    \[
        U = \{ k a K \mid kM \in V , \, a \in A_+ \} \ .
    \]
    Let $W$ be a compact subset of $ (\dX \smallsetminus \clo(V)M)M $ which contains a neighbourhood of $ \Lambda $.
    Since $ \Gamma $ acts properly discontinuously on $ X \cup \Omega $, $ S := \{ \gamma \in \Gamma \mid k_\gamma M \not \in W \} $ is finite. Thus, we may also assume without loss of generality that we sum over $ \Gamma \smallsetminus S $.
    By Corollary 2.4 of \cite[p.85]{BO00}, there is $ c > 0 $ such that
    \[
        a_{\gamma^{-1} g} = a_{\gamma^{-1} k_g a_g} \geq c \, a_\gamma a_g
    \]
    for all $ \gamma \in \Gamma \smallsetminus S $ ($ k_\gamma M \in W $) and $ g \in U $ ($ kM \in V $).
    The lemma follows now from the definition of the critical exponent.
\end{proof}

Define $ \CS_{weak}(\Gamma \bs G, \phi) $ by
\begin{multline*}
    \{ f \in \Cinf(\Gamma \bs G, \phi) \mid \LOI{p}{U}_{r, 1, Y}(f) < \infty \; \forall r \geq 0, Y \in \U(\g), \\
    \LO{p}{\Sfrak'_Q}_{r, 1, Y}(f) < \infty \; \forall U \in \U_\Gamma, \Sfrak'_Q \in \V_{std, \Gamma}, Y \in \U(\g), r \geq 0 \} \ .
\end{multline*}
Then, $ \CS_{weak}(\Gamma \bs G, \phi) $ is a Fréchet space when we equip it with the topology induced by the seminorms. Moreover,
\[
    \CS(\Gamma \bs G, \phi) \subset \CS_{weak}(\Gamma \bs G, \phi) \subset L^2(\Gamma \bs G, \phi) \ .
\]
If $ f \in \CS_{weak}(\Gamma \bs G, \phi) $, then $ f^\Omega $ and $ f^Q $ are well-defined for all $ [Q] \in \P_\Gamma $.
Set
\[
    \c{\CS_{weak}}(\Gamma \bs G, \phi) = \{ f \in \CS_{weak}(\Gamma \bs G, \phi) \mid f^\Omega = 0 , \,
    f^Q = 0 \quad \forall [Q]_\Gamma \in \P_\Gamma \} \ . 
\]
We believe that
\begin{equation}
    \overline{\c{\CS_{weak}}(\Gamma \bs G, \phi)} = \overline{\c{\CS}(\Gamma \bs G, \phi)} \ .
\end{equation}

\begin{prop}\label{Prop Examples of cusp forms}
    Assume that $ \delta_\Gamma < 0 $. Let $ f \in \CS(G, V_\phi) $. For $ g \in G $, set
    \[
        F_f(g) = \sum_{\gamma \in \Gamma} \phi(\gamma) f(\gamma^{-1}g) \ .
    \]
    Then, $F_f$ is well-defined and belongs to $ \CS_{weak}(\Gamma \bs G, \phi) $.

    If $ f \in \c{\CS}(G, V_\phi) $, then $ F_f \in \c{\CS_{weak}}(\Gamma \bs G, \phi) $.

    \smallskip

    If $ G = \SO(1, n)^0 $ $ (n \geq 2) $ and if $ \Gamma \bs X $ has at least one cusp, then $ F_f \in \CS(\Gamma \bs G, \phi) $ for all $ f \in \CS(G, V_\phi) $ and $ F_f \in \c{\CS}(\Gamma \bs G, \phi) $ for all $ f \in \c{\CS}(G, V_\phi) $.
\end{prop}

\begin{rem}
    The proof shows that we can also take smooth matrix coefficients
    \[
        c_{v, \tilde{v}}     \qquad (v \in V_{\pi, \infty} \otimes V_\phi, \tilde{v} \in V_{\pi', K})
    \]
    of an integrable discrete series representation $ \pi $ for $f$ when $ \delta_\Gamma \geq 0 $.
\end{rem}

\begin{proof}
    Assume that $ \delta_\Gamma < 0 $. Let $ f \in \CS(G, V_\phi) $ and let $ F_f $ be defined as above.
    Note that $ F_f $ is left $ \Gamma $-equivariant.

    Let $ U \in \U_\Gamma $, $ Y \in \U(\g) $ and $ r \geq 0 $. Since $ f \in \CS(G, V_\phi) $, $ \sum_{\gamma \in \Gamma} |R_Y f(\gamma^{-1}g)| $ is less or equal than
    \[
        c \sum_{\gamma \in \Gamma} a_{\gamma^{-1}g}^{-\rho} (1 + \log a_{\gamma^{-1}g})^{-r}
    \]
    for some constant $ c > 0 $. By Lemma \ref{Lemma about sum_(gamma in Gamma) a_(gamma^(-1) g)^(-(lambda + rho))} and as $ \delta_\Gamma < 0 $, there is $ c' > 0 $ such that this is less or equal than
    \[
        c' a_g^{-\rho} (1 + \log a_g)^{-r}
    \]
    for all $ g \in G $ such that $ gK \in U $. Thus, $ \LOI{p}{U}_{r, 1, Y}(F_f) $ is finite and $F_f$ is well-defined on $ \{ g \in G \mid gK \in U \} $. Moreover, $ F_f $ is smooth on this set.

    Let $Q$ be a $ \Gamma $-cuspidal parabolic subgroup of $G$. Let $ \Sfrak $ be a generalised standard Siegel set in $G$ with respect to $Q$.

    Let $ X, Y \in \U(\g) $ and $ r \geq 0 $. Without loss of generality, we may assume that $ r \geq 1 $.
    Since $ f \in \CS(G, V_\phi) $, $ \sum_{\gamma \in \Gamma_Q} |L_X R_Y f(\gamma^{-1}g)| $ is less or equal than
    \begin{equation}\label{Proof of Prop Examples of cusp forms}
        c'' \sum_{\gamma \in \Gamma_Q} a_{\gamma^{-1}g}^{-\rho} (1 + \log a_{\gamma^{-1}g})^{-2r}
    \end{equation}
    for some constant $ c'' > 0 $.
    By Proposition \ref{Prop a_(n_1 n_2 a) >= c a_(n_2 a)},
    \begin{equation}\label{eq a_(n_1 n_2 a) >= c a_(n_2 a), a_(n_1 n_2 a) >= c a_(n_1 a)}
        a_{n_1 n_2 a} \succ a_{n_2 a} \qquad \text{and} \qquad a_{n_1 n_2 a} \succ a_{n_1 a}
    \end{equation}
    for all $ n_1 \in N_{\Gamma, Q} $, $ n_2 \in N^{\Gamma, Q} $ and $ a \in A_Q $.
    As moreover $ a_{na} \geq a(\bar{n}) a^{-1} \succ a_n a^{-1} $ for all $ n \in N $ and $ a \in A $,
    \begin{equation}\label{eq a_(gamma^(-1) g) <= c a_gamma^(-1) h_Q(g)}
        a_{\gamma^{-1}g}^{-1} \prec a_{\gamma^{-1} \nu_Q(g)}^{-1} h_Q(g)
        \prec a_\gamma^{-1} a_Q(g^{-1})^{-1} \qquad (g \in \Sfrak) \ .
    \end{equation}
    It follows also from \eqref{eq a_(n_1 n_2 a) >= c a_(n_2 a), a_(n_1 n_2 a) >= c a_(n_1 a)} that
    \begin{equation}\label{eq a_(gamma^(-1) g) >= a_g}
        a_{\gamma^{-1}g} \succ a_g
    \end{equation}
    for all $ \gamma \in \Gamma_Q $ and $ g \in \Sfrak $ and it follows from \eqref{eq a_(n_1 n_2 a) >= c a_(n_2 a), a_(n_1 n_2 a) >= c a_(n_1 a)} and the estimate $ a_{na}^2 \geq a_n $ (see Lemma \ref{Lem a_(na) estimates}) that
    \[
        a_{n_1 n_2 a} \succ \sqrt{a_{n_1}}
    \]
    for all $ n_1 \in N_{\Gamma, Q} $, $ n_2 \in N^{\Gamma, Q} $ and $ a \in A_Q $.
    Thus,
    \begin{equation}\label{eq (1 + log a_(gamma^(-1)g))^(-r) <= c (1 + log a_gamma)^(-r)}
        (1 + \log a_{\gamma^{-1}g})^{-r} \prec (1 + \log a_\gamma)^{-r}
    \end{equation}
    for all $ \gamma \in \Gamma_Q $ and $ g \in \Sfrak $. Since $ \rho_{\Gamma_Q} = \delta_{\Gamma_Q} $ by Lemma \ref{Lem rho_(Gamma_Q) = delta_(Gamma_Q)},
    \begin{equation}\label{eq sum_(gamma in Gamma_Q) a_gamma^(-rho_(Gamma_Q)) (1 + log a_gamma)^(-r) < infty}
        \sum_{\gamma \in \Gamma_Q} a_\gamma^{-\rho_{\Gamma_Q}} (1 + \log a_\gamma)^{-r}
    \end{equation}
    is finite for all $ r \geq 1 $.

    Let $ g \in \Sfrak $. It follows from \eqref{eq a_(gamma^(-1) g) <= c a_gamma^(-1) h_Q(g)}, \eqref{eq a_(gamma^(-1) g) >= a_g}, \eqref{eq (1 + log a_(gamma^(-1)g))^(-r) <= c (1 + log a_gamma)^(-r)} and \eqref{eq sum_(gamma in Gamma_Q) a_gamma^(-rho_(Gamma_Q)) (1 + log a_gamma)^(-r) < infty} that there is $ c''' > 0 $ such that \eqref{Proof of Prop Examples of cusp forms} is less or equal than
    \[
        c''' a_g^{-\rho^{\Gamma_Q}} (1 + \log a_g)^{-r}
        \sum_{\gamma \in \Gamma_Q} a_{\gamma^{-1}g}^{-\rho_{\Gamma_Q}} (1 + \log a_\gamma)^{-r}
        \prec a_Q(g^{-1})^{-\rho_{\Gamma_Q}} a_g^{-\rho^{\Gamma_Q}} (1 + \log a_g)^{-r}
    \]
    as moreover $ \rho^{\Gamma_Q} + \rho_{\Gamma_Q} = \rho $.
    Thus, $ \sum_{\gamma \in \Gamma_Q} |L_X R_Y f(\gamma^{-1}g)| $ is finite for all $ g \in G $.

    Assume from now on that $ X \in \U(\n_Q) $ if $ G = \SO(1, n)^0 $ and that $ X = 1 $ otherwise. Then, $ \Ad(\gamma)^{-1}X = X $ for all $ \gamma \in \Gamma_Q $.
    Hence, $ g \in G \mapsto \sum_{\gamma \in \Gamma_Q} L_X R_Y f(\gamma^{-1}g) $ is smooth and is equal to $ L_X R_Y F_f $.
    So, $ \LO{p}{\Sfrak'}_{r, X, Y}(F_f) $ is finite.

    The part of the proposition, where we assume only that $f$ is a Schwartz function, follows now from Proposition \ref{Prop f = sum_(i=0)^m varphi_i f for all f in CS(Gamma|G)}. Consequently, $ F_f^{\Omega} $ and $ F_f^Q $ are well-defined. 

    Assume from now on that $ f \in \c{\CS}(G, V_\phi) $. Then,
    \[
        F_f^{\Omega}(g, h) = \sum_{\gamma \in \Gamma} \phi(\gamma) \int_N f(\gamma^{-1}g n h) \, dn = 0
    \]
    for all $ g \in G(\Omega) $, $ h \in G $ and $ F_f^Q(g) $ is, up to a positive multiplicative constant, equal to
    \[
        \sum_{\Gamma/\Gamma_Q} p_Q\big(\phi(\gamma) \int_{N_Q M_{\Gamma, Q}}{ f(\gamma^{-1} nmg) \, dn \, dm}\big)
        = 0
    \]
    for all $ g \in G $. This completes the proof of the proposition.
\end{proof}

\newpage

\subsubsection{Series expansion of \texorpdfstring{$f$}{f}}\label{ssec:Conv of the constant terms: Fourier expansion of f}

Let $ G = \SO(1, n)^0 $, $ n \geq 2 $, and let $ M $, $ A $, $ N $, $K$ be the subgroups of $G$ induced by the corresponding one of $ \Sp(1, n) $ computed in Appendix \ref{sec:Sp(1,n)} (then $ N \subset \exp(\g_\alpha) $).

Let $ Q $ be a $ \Gamma $-cuspidal parabolic subgroup of $G$. Without loss of generality, we may assume that $ Q = P := MAN $.

Let $ N_\Gamma = N_{\Gamma, P} $, $ M_\Gamma = M_{\Gamma, P} $, $ N^\Gamma = N^{\Gamma, P} $ and $ k = \dim(N_\Gamma) \geq 1 $.

Let $ \Gamma' $ be the normal subgroup of finite index of $ \Gamma_P $ provided by Auslander's theorem (see Section \ref{ssec:Geometrically finite groups} --- also for what follows).

Let $ \varphi \colon \{ n_\gamma \mid \gamma \in \Gamma' \} \to M_{\Gamma'} , \, n_\gamma \mapsto m_\gamma $. This defines a group homomorphism from $ \{ n_\gamma \mid \gamma \in \Gamma' \} $ to $ M_{\Gamma'} $.

By Lemma 3.2 of \cite[p.26]{BO07}, there exists a subgroup $V$ of finite index in $ \{ n_\gamma \mid \gamma \in \Gamma' \} $ such that the restriction of $ \varphi $ to $V$ extends to a group homomorphism from $ N_\Gamma $ to $M_\Gamma$, which we denote by abuse of notation also by $ \varphi $.

By 3.1.5 of \cite[p.26]{BO07}, we can choose $ \Gamma' $ with the following properties:
\begin{enumerate}
\item $ \{ n_\gamma \mid \gamma \in \Gamma' \} $ is contained in $V$;
\item $ M_{\Gamma'} $ coincides with $ \overline{\varphi(N_\Gamma)} $;
\item $ M_{\Gamma'} $ has finite index in $ M_\Gamma $.
\end{enumerate}

Let $ Z_1, \dotsc, Z_k $ be a basis of $ \n_\Gamma $ such that
\[
    \{ n_\gamma \mid \gamma \in \Gamma' \}
    = \{ \exp(n_1 Z_1 + \dotsb + n_k Z_k) \mid n_1, \dotsc, n_k \in \ZZ \} \ .
\]
For $ z \in \RR^k $, set $ n^{N_\Gamma, \{Z_j\}}_z = \exp(z_1 Z_1 + \dotsb + z_k Z_k) $.

Let $ \underline{n} \in \ZZ^k $. Then, $ n^{N_\Gamma, \{Z_j\}}_{\underline{n}} \varphi(n^{N_\Gamma, \{Z_j\}}_{\underline{n}}) \in \Gamma' $.
Set $ \tau(\underline{n}) = \phi( n^{N_\Gamma, \{Z_j\}}_{\underline{n}} \varphi(n^{N_\Gamma, \{Z_j\}}_{\underline{n}}) ) $.
Then, $ \tau $ is a unitary representation of $ \ZZ^k $ on $ V_\phi $ since
\[
    \underline{n} \in \ZZ^k \mapsto n^{N_\Gamma, \{Z_j\}}_{\underline{n}} \varphi(n^{N_\Gamma, \{Z_j\}}_{\underline{n}}) \in \Gamma'
\]
is a group isomorphism. If we identify $ \Gamma' $ as a group with $ \ZZ^k $, then $ \tau $ is equal to $ \phi $.

Since $ \Gamma' $ is isomorphic to $ \ZZ^k $ (abelian) as a group and since every unitary finite-dimensional representation on an abelian group can be written as a direct sum of 1-dimensional representations, it suffices to have a series expansion for $ \dim_\CC V_\phi = 1 $. So, let us assume that we can identify $ V_\phi $ with $ \CC $.

Then, there exist $ \lambda_j \in [0, 1) $ such that $ \tau(v) = e^{2\pi i \lambda_1 v_1} \cdots e^{2\pi i \lambda_k v_k} $ for all $ v \in \ZZ^k $. Thus, we can extend $ \tau $ to a representation of $ \RR^k $, which we denote by abuse of notation also by $ \tau $.

Set $ \lambda = (\lambda_1, \dotsc, \lambda_k)^t $. Let $ g \in G $ and let $ f \in \Cinf(\Gamma \bs G, \phi) $.

Since $ z \in \RR^k \mapsto \tau(-z) f(n^{N_\Gamma, \{Z_j\}}_z \varphi(n^{N_\Gamma, \{Z_j\}}_z) g) $ belongs to $ \Cinf(\ZZ^k \bs \RR^k) $, we have
\[
    e^{ -2 \pi i \langle z, \lambda \rangle } f(n^{N_\Gamma, \{Z_j\}}_z \varphi(n^{N_\Gamma, \{Z_j\}}_z) g) = \sum_{x \in \ZZ^k} c_{x, g}(f) e^{ 2 \pi i \langle z, x \rangle }
    \qquad (z \in \RR^k) \ ,
\]
where $ c_{x, g}(f) := \int_{\ZZ^k \bs \RR^k} f(n^{N_\Gamma, \{Z_j\}}_u \varphi(n^{N_\Gamma, \{Z_j\}}_u) g) e^{ -2 \pi i \langle u, \lambda + x \rangle } \, du $. Thus,
\[
    f(n^{N_\Gamma, \{Z_j\}}_z \varphi(n^{N_\Gamma, \{Z_j\}}_z) g) = \sum_{x \in \ZZ^k} c_{x, g}(f) e^{ 2 \pi i \langle z, \lambda + x \rangle }  \qquad (z \in \RR^k)
\]
and
\[
    f(g) = \sum_{x \in \ZZ^k} c_{x, g}(f) \ .
\]
Since $ \vol(\Gamma' \bs N_\Gamma M_{\Gamma'}) = \#(\Gamma' \bs \Gamma_P) \vol(\Gamma_P \bs N_\Gamma M_\Gamma) \#(M_\Gamma/M_{\Gamma'}) $ and since
\[
    \int_{\Gamma' \bs N_\Gamma M_{\Gamma'}}{ f^0(nmg) \, dn \, dm} = \#(\Gamma' \bs \Gamma_P) \#(M_\Gamma/M_{\Gamma'}) \int_{\Gamma_P \bs N_\Gamma M_\Gamma}{ f^0(nmg) \, dn \, dm}
\]
for all $ g \in G $,
\[
    \frac1{\vol(\Gamma' \bs N_\Gamma M_{\Gamma'})} \int_{\Gamma' \bs N_\Gamma M_{\Gamma'}}{ f^0(nmg) \, dn \, dm} = f^{P, lc}(g) \ .
\]
It follows that
\[
    (f - f^{P, lc})(g) = \sum_{x \in \ZZ^k \, : \, \lambda + x \neq 0 } c_{x, g}(f) \ .
\]
Indeed:
\begin{enumerate}
\item If $ \lambda $ is zero, then $ f = f^0 $ and $ f^{P, lc}(g) = c_{0, g}(f) $.
\item If $ \lambda $ is nonzero, then $ \lambda + x \neq 0 $ for all $ x \in \ZZ^k $, $ f^0 = 0 $ and $ (f - f^{P, lc})(g) $ is equal to $ f(g) $.
\end{enumerate}
Moreover, $ \vol(\ZZ^k \bs \RR^k) = 1 $.

\newpage

\subsubsection{Rapid decay on Siegel sets}

\begin{lem}\label{Lem sum_(0 neq x in Z^k) 1/( |x|_2^(N+1) ) }
    Let $ k \in \NN $, let $ \lambda \in \RR^k $ and let $ \|\cdot\| $ be a norm on $ \RR^k $. Then,
    \[
        \sum_{x \in \ZZ^k \, : \, x \neq -\lambda} \frac1{ \|\lambda + x\|^{N+1} }
    \]
    converges if $ N \geq k $.
\end{lem}

\begin{proof}
    Without loss of generality, $ \|\cdot\| = \|\cdot\|_2 $ and $ \lambda_j \in [0, 1) $. Then, $ x_j + \lambda_j = 0 $ if and only if $ x_j = 0 $ and $ \lambda_j = 0 $.

    Then, $ x \neq -\lambda $ for all $ x \in \ZZ^k \smallsetminus \{0\} $. Moreover, there are constants $ c_1, c_2 > 0 $ such that
    \[
        c_1 \|x\|_2 \leq \|\lambda + x\|_2 \leq c_2 \|x\|_2
    \]
    for all $ x \in \ZZ^k \smallsetminus \{0\}$. So, the above series converges if and only if $ \sum_{0 \neq x \in \ZZ^k} \frac1{ \|x\|_2^{N+1} } $ converges.

    Let $ N \geq k $. Since $ 0 \neq x \in \RR^k \mapsto \frac1{ \|x\|_2^{N+1} } $ continuous, nonnegative function that is decreasing in each variable,
    \[
        \sum_{x \in \ZZ^k \, : \, x_j \neq 0} \frac1{ \|x\|_2^{N+1} }
        = \sum_{x \in \ZZ^k \, : \, |x_j| \geq 1} \frac1{ \|x\|_2^{N+1} }
    \]
    converges by the integral test for convergence applied repeatedly as
    \begin{multline*}
        \int_{ \{x \in \RR^n \, : \, |x_j| \geq 1\} } \frac1{\|y\|_2^{m+1}} \, dy
        \leq \int_{ \{x \in \RR^n \, : \, \|x\|_2 \geq 1\} } \frac1{\|y\|_2^{n+1}} \, dy \\
        = \vol(S^{n-1}) \int_1^\infty \frac1{ r^{n+1} } \cdot r^{n-1} \, dr
        = \vol(S^{n-1}) \int_1^\infty \frac1{ r^2 } \, dr < \infty
    \end{multline*}
    for all $ m \geq n \geq 1 $. This shows in particular that the lemma holds for $ k = 1 $.
    \\ Let $ E_k = \ZZ^k \smallsetminus \{0\} $ and let
    \[
        E_{k, i} = \{ (x_1, \dotsc, x_{i-1}, x_{i+1}, \dotsc, x_k) \in \RR^{k-1} \mid
            (x_1, \dotsc, x_{i-1}, x_{i+1}, \dotsc, x_k) \neq 0 \} \ .
    \]
    Since
    \[
        \sum_{x \in E_k} \frac1{ \|x\|_2^{N+1} }
        = \sum_{x \in \ZZ^k \, : \, x_j \neq 0} \frac1{ \|x\|_2^{N+1} }
        + \sum_{i = 1}^k \sum_{x \in E_{k, i}} \frac1{ \|(x_1, \dotsc, x_{i-1}, x_{i+1}, \dotsc, x_k) \|_2^{N+1} } \ ,
    \]
    the lemma follows by induction.
\end{proof}

\begin{dfn}
    Let $ f $ be a $ \Gamma $-equivariant smooth function from $ G $ to $ V_\phi $ and let $ \Sfrak $ be a (nongeneralised) Siegel set in $G$ with respect to $ Q $.
    We say that $f$ is \textit{rapidly decreasing} on $ \Sfrak $ if for all $ r \geq 0 $, $ X \in \U(\n_Q) $, $ Y \in \U(\g) $ and all $ N \in \NN $ there exists a constant $ C > 0 $ such that
    \[
        |L_X R_Y f(g)| \leq C h_Q(g)^{-N \alpha_Q} a_g^{-\rho^{\Gamma_Q}} (1 + \log a_g)^{-r}
    \]
    for all $ g \in \Sfrak $.\index{rapidly decreasing}
\end{dfn}

\begin{rem}
    We restrict ourselves here to Siegel sets as the estimate is only interesting when $ h_Q(g) $ is big.
\end{rem}

\begin{thm}\label{Thm f - f^(Q, lc) is rapidly decreasing}
    Let $ Q $ be a $ \Gamma $-cuspidal parabolic subgroup of $G$ and let $ f $ be a Schwartz function near $eQ$. Let $ \Sfrak $ be an admissible (nongeneralised) standard Siegel set in $G$ with respect to $Q$. 
    If $ \n_{\Gamma, Q} $ is an ideal in $ \n_Q $, then $ f - f^{Q, lc} $ is rapidly decreasing on $ \Sfrak $.
\end{thm}

\begin{rem} \nlenum
    \begin{enumerate}
    \item $ \n_{\Gamma, Q} $ is an ideal in $ \n_Q $ if and only if $ \g_{2\alpha_Q} $ is contained in $ \n_{\Gamma, Q} $ or $ \n_{\Gamma, Q} $ is contained in $ \g_{2\alpha_Q} $.
    \item If $eQ$ has full rank or if we are in the real hyperbolic case, then $ \n_{\Gamma, Q} $ is an ideal in $ \n_Q $.
    \item Compare with, e.g., Theorem 7.5 of \cite[p.72]{Borel} and Theorem 6.9 of \cite[p.27]{Borel2}.
    \end{enumerate}
\end{rem}

\begin{proof}
    Let $ Y \in \U(\n_Q) $ and $ Y' \in \U(\g) $. Without loss of generality, we may assume that $ Y' = 1 $.
    Without loss of generality, we may assume that $Q$ is equal to $P$. Let us consider first the real hyperbolic case: $ G = \SO(1, n)^0 $ ($ n \geq 2 $).
    Let $ k = \dim(N_{\Gamma, P}) $.
    
    Since $ \Gamma' $ is isomorphic to $ \ZZ^k $ (abelian) as a group and since every unitary finite-dimensional representation on an abelian group can be written as a direct sum of one-dimensional representations, we may assume without loss of generality that $ \dim_\CC V_\phi = 1 $. Let us identify $ V_\phi $ with $ \CC $.

    In the following, we use the notations seen in Section \ref{ssec:Conv of the constant terms: Fourier expansion of f}. Then,
    \[
        f(g) = \sum_{x \in \ZZ^k} c_{x, g}(f) \ .
    \]
    Let $ x = (x_1, \dotsc, x_k) \in \ZZ^k $. Let $ j \in \{1, \dotsc, k\} $ be such that $ |x_j| = \max_i |x_i| $.
    By doing integration by parts repeatedly, we have
    \begin{multline}\label{eq partial integration of Fourier coefficients}
        \int_{[0, 1]^k} (L_Y f)(n^{N_\Gamma, \{Z_i\}}_u \varphi(n^{N_\Gamma, \{Z_i\}}_u) \, \cdot)(g) e^{-2 \pi i \langle u, \lambda + x \rangle} \, du \\
        = \big(\frac{-1}{2\pi i (\lambda_j + x_j)}\big)^p \int_{[0, 1]^k} \frac{ \del^p}{\del u_j^p} \big((L_Y f)(n^{N_\Gamma, \{Z_i\}}_u \varphi(n^{N_\Gamma, \{Z_i\}}_u) \, \cdot)(g)\big) e^{-2 \pi i \langle u, \lambda + x \rangle} \, du \ .
    \end{multline}
    To simplify notations, we continue the proof with the assumption $ Y = 1 $.
    Without this further assumption, one must use at some moment that we can write $ \Ad(n m) X $ ($ n \in \omega, m \in M_\Gamma $) as finite linear combination of elements in $ \U(\n) $ with bounded coefficients.

    Let $ \{X_l\}_{l=1, \dotsc, \dim(\g)} $ be a basis of $ \g $. Then, there exist smooth functions $ c_{j, l} $ on $K$ such that
    \[
        \Ad(k(g)^{-1})(Z_j) = \sum_{l=1}^{\dim(\g)} c_{j, l}(k(g)^{-1}) X_l \ .
    \]
    Assume that either $ \lambda $ or $ x $ is nonzero. Since $ \frac{ \del^p f}{\del u_j^p}(n^{N_\Gamma, \{Z_i\}}_u \varphi(n^{N_\Gamma, \{Z_i\}}_u) g) $ is equal to
    \[
        R_{Z_j}^p (R_g f)(n^{N_\Gamma, \{Z_i\}}_u \varphi(n^{N_\Gamma, \{Z_i\}}_u))
    \]
    by the chain rule.
    This is again equal to
    \begin{align*}
        & \restricted{ \frac{ \del^p }{ \del s_1 \cdots \del s_p } }{ s_1 = \dotsm = s_p = 0 } f(n^{N_\Gamma, \{Z_i\}}_u \varphi(n^{N_\Gamma, \{Z_i\}}_u) g \\
        & \qquad \exp(h(g)^{-\alpha_P} (s_1 \Ad(k(g)^{-1})(Z_j) + \dotsb + s_p \Ad(k(g)^{-1})(Z_j)))  \\
        =& h(g)^{-p\alpha_P} \sum_{m=1}^p \sum_{l_m=1}^{\dim(\g)} \prod_{i=1}^p c_{j, l_i}(k(g)^{-1}) (R_{X_{l_1} \dotsm X_{l_p}} f)(n^{N_\Gamma, \{Z_i\}}_u \varphi(n^{N_\Gamma, \{Z_i\}}_u) g) \ ,
    \end{align*}
    \eqref{eq partial integration of Fourier coefficients} is equal to
    \begin{multline*}
        \big(\frac{-1}{2\pi i (\lambda_j + x_j)}\big)^p h(g)^{-p\alpha_P} \sum_{m=1}^p \sum_{l_m=1}^{\dim(\g)} c_{j, l_m}(k(g)^{-1}) \\
        \int_{[0, 1]^k} R_{X_{l_1} \dotsm X_{l_p}} f(n^{N_\Gamma, \{Z_i\}}_u \varphi(n^{N_\Gamma, \{Z_i\}}_u) g) e^{-2 \pi i \langle u, \lambda + x \rangle } \, du \ .
    \end{multline*}
    As $ f $ is a Schwartz function near $eP$ and as $ \sup_{k \in K, j, l} |c_{j, l}(k)| $ is finite, the absolute value of this is less or equal than
    \[
        C \sum_{m=1}^p \sum_{l_m = 1}^{\dim(\g)} \LO{p}{\Sfrak'}_{r, 1, X_{l_1} \dotsm X_{l_p}}(f) \cdot \frac{1}{(\max_i |\lambda_i + x_i|)^p} \cdot h(g)^{\rho_{\Gamma_P} - p\alpha_P} a_g^{-\rho^{\Gamma_P}} (1 + \log a_g)^{-r}
    \]
    for some constant $ C > 0 $. The theorem follows in the real hyperbolic case as
    \[
        (f - f^{P, lc})(g) = \sum_{x \in \ZZ^k \, : \,  x \neq -\lambda } c_{x, g}(f)
    \]
    is equal to
    \begin{align*}
        & \sum_{x \in \ZZ^k \, : \,  x \neq -\lambda} \int_{[0, 1]^k} f(n^{N_\Gamma, \{Z_i\}}_u \varphi(n^{N_\Gamma, \{Z_i\}}_u) g) e^{-2 \pi i \langle u, \lambda + x \rangle } \, du \\
        =& \sum_{x \in \ZZ^k \, : \,  x \neq -\lambda} \big(\frac{-1}{2\pi i (\lambda_j + x_j)}\big)^p h(g)^{-p\alpha_P} \sum_{m=1}^p \sum_{l_m=1}^{\dim(\g)} c_{j, l_m}(k(g)^{-1}) \\
        & \qquad \int_{[0, 1]^k} (R_{X_{l_1} \dotsm X_{l_p}} f)(n^{N_\Gamma, \{Z_i\}}_u \varphi(n^{N_\Gamma, \{Z_i\}}_u) g) e^{-2 \pi i \langle u, \lambda + x \rangle } \, du
    \end{align*}
    and as $ \sum_{x \in \ZZ^k \, : \,  x \neq -\lambda} \frac1{\|\lambda + x\|_\infty^p} $ converges for $p$ sufficiently large by Lemma \ref{Lem sum_(0 neq x in Z^k) 1/( |x|_2^(N+1) ) }.

    \medskip

    \needspace{2\baselineskip}
    Let us consider now the general case. To simplify notations, we still assume that $ Y = 1 $.

    Since $ \Gamma' \cap N_{2\alpha} M_{\Gamma'} $ is abelian, $ V_\phi $ can be decomposed into an orthogonal direct sum of 1-dimensional $ (\Gamma' \cap N_{2\alpha} M_{\Gamma'} ) $-invariant spaces.

    Without loss of generality, we may assume that $ f $ takes values in such a 1-dimensional space $ W_\phi \equiv \CC $.

    As $ L := \{ n_\gamma \mid \gamma \in \Gamma' \cap N_{2\alpha} M_{\Gamma'} \} $ is a cocompact lattice of $ N_\Gamma \cap N_{2\alpha} $ (abelian), $ L $ is isomorphic to some $ \ZZ^l $.

    Let now $ Z_1, \dotsc, Z_l $ be a basis of $ \n_\Gamma \cap \g_{2\alpha} $ such that
    \[
        L = \{ \exp(n_1 Z_1 + \dotsb + n_l Z_l) \mid n_1, \dotsc, n_l \in \ZZ \} \ .
    \]
    For $ z \in \RR^l $, set $ n^{N_\Gamma \cap N_{2\alpha}, \{Z_j\}}_z = \exp(z_1 Z_1 + \dotsb + z_l Z_l) $.

    Let $ p $ the orthogonal projection from $ W_\phi $ to $ W_\phi^{\Gamma' \cap N_{2\alpha} M_{\Gamma'}} $.
    
    Similarly as before, one can show that $f$ has a series expansion:
    \[
        f(g) = \sum_{x \in \ZZ^l} c_{x, g}(f) \ ,
    \]
    where $ c_{x, g}(f) := \int_{[0, 1]^l} f(n^{N_\Gamma \cap N_{2\alpha}, \{Z_i\}}_u g) e^{ -2 \pi i \langle u, \lambda + x \rangle } \, du $.

    Let $ M_L = \overline{\{ m_\gamma \mid \gamma \in \Gamma' \cap N_{2\alpha} M_{\Gamma'} \}} $.

    Let $ c_1 = \frac1{ \vol((\Gamma' \cap N_{2\alpha} M_{\Gamma'}) \bs (N_\Gamma \cap N_{2\alpha})M_L) } $ and $ c_2 = \frac1{\vol((\Gamma' \cap N_{2\alpha} M_{\Gamma'} \bs \Gamma') \bs ((N_\Gamma \cap N_{2\alpha}) M_L \bs N_\Gamma M_{\Gamma'}))} $. Then,
    \[
        c_1 c_2 = \frac1{\vol(\Gamma' \bs N_\Gamma M_{\Gamma'})} \ .
    \]
    Let
    \[
        f^{P, lc, 1}(g) = c_1 \int_{(\Gamma' \cap N_{2\alpha} M_{\Gamma'}) \bs (N_\Gamma \cap N_{2\alpha}) M_L} (p(f))(n m g) \, dn \, dm
    \]
    (called ``\textit{partial little constant term}'').

    If $ W_\phi = W_\phi^{\Gamma' \cap N_{2\alpha} M_{\Gamma'}} $, then $ \lambda = 0 $. Otherwise, $ W_\phi^{\Gamma' \cap N_{2\alpha} M_{\Gamma'}} = 0 $ and then $ \lambda \neq 0 $. Hence, $ \lambda + x \neq 0 $ for all $ x \in \ZZ^l $. Moreover, $ f^{P, lc, 1} = 0 $ as $ p(f) = 0 $ and therefore $ f^{P, lc} $ vanishes, too.

    As previously, one can show that
    \[
        f(g) - f^{P, lc, 1}(g) = \sum_{x \in \ZZ^l \, : \, x \neq -\lambda} c_{x, g}(f)
    \]
    has rapid decay. Thus, $ f = f - f^{P, lc, 1} = f - f^{P, lc} $ decays rapidly if $ W_\phi^{\Gamma' \cap N_{2\alpha} M_{\Gamma'}} = 0 $.
    It remains to consider the case $ W_\phi = W_\phi^{\Gamma' \cap N_{2\alpha} M_{\Gamma'}} $.
    Then, $ \phi $ induces a representation of $ \Gamma' \cap N_{2\alpha} M_{\Gamma'} \bs \Gamma' $ on $ W_\phi^{\Gamma' \cap N_{2\alpha} M_{\Gamma'}} $.

    Since $ f(g) - f^{P, lc}(g) = f(g) - c_1 c_2 \int_{\Gamma' \bs N_\Gamma M_{\Gamma'}} f^0(n m g) \, dn \, dm $ is equal to
    \begin{align*}
        & \big(f(g) - f^{P, lc, 1}(g)\big)
        + \big(f^{P, lc, 1}(g) - c_1 c_2 \int_{\Gamma' \bs N_\Gamma M_{\Gamma'}} f^0(n m g) \, dn \, dm \big) \\
        =& \big( f(g) - f^{P, lc, 1}(g) \big) \\
        & \qquad + \Big( f^{P, lc, 1}(g) - c_2 \int_{(\Gamma' \cap N_{2\alpha} \bs \Gamma') \bs (N_\Gamma \cap N_{2\alpha} \bs N_\Gamma M_{\Gamma'})} (f^{P, lc, 1})^0(n m g) \, dn \, dm \Big) \ ,
    \end{align*}
    it remains to show that the last two terms decay rapidly.

    The group $ N_\Gamma \cap N_{2\alpha} \bs N_\Gamma $ is abelian ($ [N_\Gamma, N_\Gamma] \subset N_\Gamma \cap N_{2\alpha} $) and can be identified with $ \RR^{l'} $, where
    \[
        l' := \dim_\RR(N_\Gamma \cap N_{2\alpha} \bs N_\Gamma)
        = \dim_\RR\big((\n_\Gamma \cap \g_{2\alpha}) \bs \n_\Gamma\big) \ .
    \]
    Let $ C $ be the orthogonal complement of $ \n_\Gamma \cap \g_{2\alpha} $ in $ \n_\Gamma $.

    Since $ \{ n_\gamma \mid \gamma \in \Gamma'\} $ is a cocompact lattice in $ N_\Gamma $ and since $ L $ is a cocompact lattice in $ N_\Gamma \cap N_{2\alpha} $, $ L \bs \{ n_\gamma \mid \gamma \in \Gamma'\} $ is a cocompact lattice in $ (N_\Gamma \cap N_{2\alpha}) \bs N_\Gamma $.

    Let $ Z'_1, \dotsc, Z'_{l'} $ be a basis of $ C $ such that
    \[
        \{ L n_\gamma \mid \gamma \in \Gamma' \}
        = \{ L \exp(n_1 Z'_1 + \dotsb + n_{l'} Z'_{l'}) \mid n_1, \dotsc, n_{l'} \in \ZZ \} \ .
    \]
    For $ z \in \RR^{l'} $, set $ n^{(N_\Gamma \cap N_{2\alpha}) \bs N_\Gamma, \{Z'_j\}}_z = \exp(z_1 Z'_1 + \dotsb + z_{l'} Z'_{l'}) $.

    Hence, we have a series expansion for $ c_{0, g}(f) = f^{P, lc, 1}(g) $ ($ \lambda = 0 $):
    \[
        f^{P, lc, 1}(g) = \sum_{x \in \ZZ^{l'}} c'_{x, g}(f) \ ,
    \]
    where $ c'_{x, g}(f) $ is given by
    \[
         \int_{[0, 1]^{l'}} f^{P, lc, 1}(n^{(N_\Gamma \cap N_{2\alpha}) \bs N_\Gamma, \{Z'_i\}}_u \varphi(n^{(N_\Gamma \cap N_{2\alpha}) \bs N_\Gamma, \{Z'_i\}}_u) g)
        e^{ -2 \pi i \langle u, \lambda' + x \rangle } \, du \ .
    \]
    Again, one can show similarly as in the real hyperbolic case that
    \begin{align*}
        & f^{P, lc, 1}(g) - f^{P, lc}(g) \\
        =& f^{P, lc, 1}(g) -  c_2 \int_{((\Gamma' \cap N_{2\alpha} M_{\Gamma'}) \bs \Gamma') \bs ((N_\Gamma \cap N_{2\alpha}) M_L \bs N_\Gamma M_{\Gamma'})} (f^{P, lc, 1})^0(n m g) \, dn \, dm \\
        =& \sum_{x \in \ZZ^{l'} \, : \, x \neq -\lambda'} c'_{x, g}(f)
    \end{align*}
    has rapid decay. Indeed, as $ [\n^\Gamma, \n_\Gamma] \subset \n_\Gamma $ (recall that $ \n_\Gamma $ is an ideal in $ \n $), there is $ Z' \in U(\n_\Gamma \cap \g_{2\alpha}) $ such that
    \[
        \Ad(p_{N^\Gamma}(\nu(g)^{-1}))Z = Z + Z'
    \]
    for every $ g \in \Sfrak $ and every $ Z \in \U(\n_\Gamma) $.
    The theorem follows.
\end{proof}

\newpage

\subsection{Relation between the constant term along \texorpdfstring{$ \Omega $}{Omega} and the one along a smaller rank cusp}\label{ssec:Convergence of the constant terms}

Since we can approach a cusp of smaller rank by ordinary points, the question whether one can relate also the constant terms at these points by a limit formula rises naturally. We show in this section that this is indeed possible in the real hyperbolic case. The formula is given in Theorem \ref{Thm Conv of the constant terms} on p.\pageref{Thm Conv of the constant terms}. With the help of this formula we can then show that the constant terms of a Schwartz function $f$ along a cusp of smaller rank vanishes if the constant term of $f$ along $ \Omega $ is zero (Theorem \ref{Thm f^Omega = 0 implies f^P = 0}, p.\pageref{Thm f^Omega = 0 implies f^P = 0}).

\medskip

Let $ G = \SO(1, n)^0 $, $ n \geq 2 $, and let $ M $, $ A $, $ N $, $K$ be the subgroups of $G$ induced by the corresponding one of $ \Sp(1, n) $ computed in Appendix \ref{sec:Sp(1,n)} (then $ N \subset \exp(\g_\alpha) $).

Let $ \Gamma $ be a geometrically finite subgroup of $G$.
Recall that
\[
    w = \begin{pmatrix} 1&0&0&0\\0&-1&0&0\\0&0&-1&0\\0&0&0&I_{n-2} \end{pmatrix} \ ,
\]
\[
    h(w n_v) = a_{ \frac1{1 + \|v\|^2} }    \qquad (v \in \RR^{n-1})
\]
and that
\[
    \nu(w n_v) = n_{\frac1{1+\|v\|^2}(-v_1, v_2, \dotsc, v_{n-1})} \qquad (v \in \RR^{n-1}) \ .
\]
Let $ Q $ be a $ \Gamma $-cuspidal parabolic subgroup of $G$ with associated cusp of smaller rank. 
Without loss of generality, we may assume that $ Q = P := MAN $.

Let $ N_\Gamma = N_{\Gamma, P} $, $ M_\Gamma = M_{\Gamma, P} $, $ N^\Gamma = N^{\Gamma, P} $ and $ k = \dim(N_\Gamma) \geq 1 $.

Since $eP$ is not of full rank, $ n-k-1 \geq 1 $. This is only possible if $ n \geq 3 $.

Let $ M_k = \{ m \in M \mid m n m^{-1} = n \quad \forall n \in N^\Gamma , \, m n' m^{-1} \in N_\Gamma \quad \forall n' \in N_\Gamma \} $.

As $ m_A n_v m_A^{-1} = n_{Av} $, we may assume without loss of generality that
\[
    N_\Gamma = \{ n^{N_\Gamma}_{(v_{n-k}, \dotsc, v_{n-1})} := n_{(0_{1, n-k-1}, v_{n-k}, \dotsc, v_{n-1})} \mid v_i \in \RR \}
\]
(by replacing $ \Gamma_P $ by $ m \Gamma_P m^{-1} $ if necessary).
Then,
\[
    N^\Gamma = \{ n^{N^\Gamma}_{(v_1, \dotsc, v_{n-k-1})} := n_{(v_1, \dotsc, v_{n-k-1}, 0_{1, k})} \mid v_i \in \RR \}
\]
and
\[
    M_k = \Big\{ m_{\begin{pmatrix} I_{n-k+1} & 0_{n-k+1, k} \\ 0_{k, n-k+1} & A \end{pmatrix}} \mid A \in \SO(k) \Big\} \ .
\]

Define
\[
    \Phi_2(\psi_0) = (\cos \psi_0, \sin \psi_0)
\]
and
\[
    \Phi_d(\psi_0, \psi_1, \dotsc, \psi_{d-2}) = (\cos \psi_{d-2} \cdot \Phi_{d-1}(\psi_0, \psi_1, \dotsc, \psi_{d-3}), \sin \psi_{d-2}) \ .
\]

\subsubsection{Several small results}

Let $ d \in \NN $, $ d \geq 2 $. Let us show now how the spherical coordinates and $ \so(d) $ are related.

\begin{lem}\label{Lem Relation between spherical coordinates and so(d)}
    There exist elements $ Y_0, Y_1, \dotsc, Y_{d-2} \in \so(d) $ such that
    \[
        (\Phi_d(\psi_0, \psi_1, \dotsc, \psi_{d-2}))^t = \prod_{j=0}^{d-2} \exp(-\psi_j Y_j) (1, 0_{1, d-1})^t
    \]
    for all $ \psi_j \in \RR $.
\end{lem}

\begin{proof}
    Let $ d \in \NN $, $ d \geq 2 $. For $ j \in \{0, \dotsc, d-2\} $, define
    \[
        Y_j = \begin{pmatrix}
            0               & 0_{1, j}      & 1             & 0_{1, d-j-2} \\
            0_{j, 1}        & 0_{j, j}      & 0_{j, 1}      & 0_{j, d-j-2} \\
            -1              & 0_{1, j}      & 0             & 0_{1, d-j-2} \\
            0_{d-j-2, 1}    & 0_{d-j-2, j}  & 0_{d-j-2, 1}  & 0_{d-j-2, d-j-2}
        \end{pmatrix} \in \so(d) \ .
    \]
    The lemma follows easily by induction on $d$.
\end{proof}

Let us denote the partial derivative with respect to the $j$-th Cartesian coordinate in $ \RR^d $ by $ \del_j $ and the partial derivatives relative to the spherical coordinates in $ \RR^d $ by $ \del_r $, $ \del_{\psi_0} $, $ \ldots $, $ \del_{\psi_{d-2}} $.

\medskip

The following lemma shows how $ \del_j $ ($ j \in \{1, \dotsc, d\} $) is expressed in spherical coordinates.

\begin{lem}\label{Lem cartesian partial derivatives in terms spherical partial derivatives}
    There exist continuous functions $ p_{i, j}^{(d)} \colon (0, 2\pi) \times (-\tfrac{\pi}2, \tfrac{\pi}2)^{d-2} \to \RR $ such that
    \[
        \cos \psi_1 (\cos \psi_2)^2 \cdots (\cos \psi_{d-2})^{d-2} p_{i, j}^{(d)}(\psi_0, \dotsc, \psi_{d-2})
    \]
    are polynomials in
    \[
        \cos(\psi_0), \dotsc, \cos(\psi_{d-2}), \sin(\psi_0), \dotsc, \sin(\psi_{d-2})
    \]
    for all $ i \in \{0, \dotsc, d-2\} $, $ j \in \{1, \dotsc, d\} $ and
    \[
        \del_j = \Phi_d(\psi_0, \dotsc, \psi_{d-2})_j \del_r + \frac1r \sum_{i=0}^{d-2} p_{i, j}^{(d)}(\psi_0, \dotsc, \psi_{d-2}) \del_{\psi_i} \ .
    \]
\end{lem}

\begin{proof}
    By induction, one can show that
    \begin{equation}\label{eq det d Phi_d}
        \det\bigl(\d (r\Phi_d(\psi_0, \psi_1, \dotsc, \psi_{d-2}))\bigr) = r^{d-1} \cos \psi_1 (\cos \psi_2)^2 \cdots (\cos \psi_{d-2})^{d-2} \ .
    \end{equation}
    Let $ A_j $ be the matrix formed by replacing the $ j $-th column of
    \[
        \d (r\Phi_d(\psi_0, \psi_1, \dotsc, \psi_{d-2}))^t
    \]
    by the column vector $ (\del_r, \del_{\psi_0}, \del_{\psi_1}, \dotsc, \del_{\psi_{d-2}}) $. By the chain rule, we have
    \[
        \begin{pmatrix} \del_r \\ \del_{\psi_0} \\ \del_{\psi_1} \\ \vdots \\ \del_{\psi_{d-2}} \end{pmatrix}
        = \big(\d (r\Phi_d(\psi_0, \psi_1, \dotsc, \psi_{d-2}))\big)^t \begin{pmatrix} \del_1 \\ \del_2 \\ \vdots \\ \del_d \end{pmatrix}
    \]
    By Cramer's rule and \eqref{eq det d Phi_d}, we get
    \[
        \del_j = \frac{ \det(A_j) }{ r^{d-1} \cos \psi_1 (\cos \psi_2)^2 \cdots (\cos \psi_{d-2})^{d-2} } \qquad (j \in \{1, \dotsc, d\}) \ .
    \]
    Let $ q_{i, j}^{(d)}(r, \psi_0, \psi_1, \dotsc, \psi_{d-2}) $ be the $ \del_r $-coordinate, (resp. $ \del_{\psi_{i-2}} $-coordinate) if $ i = -1 $ (resp. $ i \in \{0, \dotsc, d-2\} $) of $ \frac{ \det(A_j) }{ r^{d-1} \cos \psi_1 (\cos \psi_2)^2 \cdots (\cos \psi_{d-2})^{d-2} } $.

    Note that $ q_{-1, j}^{(d)}(r, \psi_0, \psi_1, \dotsc, \psi_{d-2}) $ and $ r q_{i, j}^{(d)}(r, \psi_0, \psi_1, \dotsc, \psi_{d-2}) $ ($ i \in \{0, \dotsc, d-2\} $) are independent of $r$ for all $ j $.
    For $ j \in \{1, \dotsc, d\} $, set
    \[
        p_{-1, j}^{(d)}(\psi_0, \psi_1, \dotsc, \psi_{d-2}) = q_{-1, j}^{(d)}(r, \psi_0, \psi_1, \dotsc, \psi_{d-2})
    \]
    and
    \[
        p_{i, j}^{(d)}(\psi_0, \psi_1, \dotsc, \psi_{d-2}) = r q_{i, j}^{(d)}(r, \psi_0, \psi_1, \dotsc, \psi_{d-2}) \qquad (i \in \{0, \dotsc, d-2\}) \ .
    \]
    Then, the functions $ p_{i, j}^{(d)} \colon (0, 2\pi) \times (-\tfrac{\pi}2, \tfrac{\pi}2)^{d-2} \to \RR $ are as required.

    Let $ f \colon \RR^d \to \RR, \; v \mapsto \|v\|_2 $. Then, $ \del_r f = 1 $ and $ \del_{\psi_0} f = \del_{\psi_1} f = \dotsm = \del_{\psi_{d-2}} f = 0 $.
    If $ v = r \Phi_d(\psi_0, \psi_1, \dotsc, \psi_{d-2}) $, then
    \[
        p_{-1, j}^{(d)}(\psi_0, \psi_1, \dotsc, \psi_{d-2}) = \del_j(\|v\|_2) = \frac{v_j}{\|v\|_2} = \Phi_d(\psi_0, \psi_1, \dotsc, \psi_{d-2})_j
    \]
    for all $ j \in \{1, \dotsc, d\} $. The lemma follows.
\end{proof}

\begin{lem}\label{Lem relation between Haar measures}
    Let $ f $ be a continuous function on $ S^{d-1} $ ($ d \geq 2 $). Then,
    \[
        \int_{\SO(d)} f(A (1, 0_{1, d-1})^t) \, dA = \frac1{\vol(S^{d-1})} \int_{S^{d-1}} f(x) \, dx \ .
    \]
\end{lem}

\begin{rem}
    If $ f $ is a function on $ S^0 = \{-1, 1\} $, then
    \[
        \int_{S^0} f(x) \, dx = \frac12(f(1) + f(-1)) \ .
    \]
\end{rem}

\begin{proof}
    Let $ d \geq 2 $. Let $ g \colon \SO(d)/\SO(d-1) \to S^{d-1} , \ A\SO(d-1) \mapsto A (1, 0_{1, d-1})^t $.
    \\ Then, $g$ is a diffeomorphism which identifies these two spaces.

    Since $ \SO(d)/\SO(d-1) $ has up to a multiplicative constant a unique invariant measure, the push-forward of the spherical measure on $ S^{d-1} $ by $ g $ divided by $ \vol(S^{d-1}) $ is equal to the normalised invariant measure on $ \SO(d)/\SO(d-1) $.

    Let $ f $ be a continuous function on $ S^{d-1} $. Then, $ \int_{\SO(d)} f(A (1, 0_{1, d-1})^t) \, dA $ is equal to
    \[
        \int_{\SO(d)/\SO(d-1)} f(A (1, 0_{1, d-1})^t) \, dA = \frac1{\vol(S^{d-1})}  \int_{S^{d-1}} f(x) \, dx
    \]
    as $ \vol(\SO(d-1)) = 1 $.
\end{proof}

Let $ n \geq 3 $. For $ A = (a_{i, j})_{i, j} \in \SO(n-1) $, we have
\[
    w m_A w
    = \begin{pmatrix}
        1   &   0   &   0           &   0           & \cdots    &   0 \\
        0   &   1   &   0           &   0           & \cdots    &   0 \\
        0   &   0   &   a_{1, 1}    &   -a_{1, 2}   & \cdots    & -a_{1, n-1} \\
        0   &   0   &   -a_{2, 1}   &   a_{2, 2}    & \cdots    & a_{2, n-1} \\
        0   &   0   &   \vdots      &   \vdots      & \ddots    & \vdots \\
        0   &   0   &   -a_{n-1, 1} &   a_{n-1, 2}  & \cdots    & a_{n-1, n-1}
    \end{pmatrix} \ .
\]
In particular, $ w m_A w = m_A $ if $ A = \begin{pmatrix}
    1   &   0_{1, n-2} \\
    0_{n-2, 1}   &   B
\end{pmatrix} $ for some $ B \in \SO(n-2) $.

Thus, $ w m w = m $ for all $ m \in M_k $.

\begin{lem}
    For all $ n \in N $ and $ m \in M $, we have
    \begin{enumerate}
    \item $ h(w m n m^{-1}) = h(w n) $;
    \item $ \nu(w m n m^{-1}) = (w m w^{-1}) \nu(w n) (w m w^{-1})^{-1} $;
    \item $ k(w m n m^{-1}) = w m w^{-1} k(w n) m^{-1} $.
    \end{enumerate}
\end{lem}

\begin{proof}
    Let $ n \in N $ and $ m \in M $. As $ w M w^{-1} \subset M $,
    \begin{align*}
        w m n m^{-1} &= (w m w^{-1}) w n m^{-1} \\
        &= (w m w^{-1}) \nu(w n) h(w n) k(w n) m^{-1} \\
        &= ((w m w^{-1}) \nu(w n) (w m w^{-1})) h(w n) (w m w^{-1} k(w n) m^{-1}) \ .
    \end{align*}
    The lemma follows.
\end{proof}

\needspace{3\baselineskip}
Let $ \bar{n}_v = \theta(n_v) \in \bar{N} $.

\begin{lem}\label{Lem about Iwasawa k}
    For $ n \in N $, we have
    \[
        n(\theta n) = \nu(\theta n) \ ,
    \]
    \[
        a(\theta n) = h(\theta n)^{-1}
    \]
    and
    \[
        \kappa(\theta n) = k(\theta n) = w k(w n) \ .
    \]
    Moreover, for any $ v \in \RR^{n-1} $ ($ n \geq 2 $), we have
    \[
        k(\bar{n}_v) = n_v a_{\tfrac1{1 + \|v\|^2}} \bar{n}_v = \bar{n}_v a_{1 + \|v\|^2} n_v \ .
    \]
\end{lem}

\begin{proof}
    For $ v \in \RR^{n-1} $, we have
    \[
        n(\bar{n}_v) = \nu(\bar{n}_{-v})^{-1} = n_{-\tfrac1{1+\|v\|^2}v} = \nu(\bar{n}_v) \ ,
    \]
    \[
        a(\bar{n}_v) = h(\bar{n}_{-v})^{-1} = a_{1 + \|v\|^2} = h(\bar{n}_v)^{-1}
    \]
    and
    \[
        k(\bar{n}_v) = h(\bar{n}_v)^{-1} \nu(\bar{n}_v)^{-1} \bar{n}_v
        = a_{1 + \|v\|^2} n_{\tfrac1{1+\|v\|^2}v} \bar{n}_v
        = n_v a_{\tfrac1{1 + \|v\|^2}} \bar{n}_v \ .
    \]
    Thus, $ k(\bar{n}_v) $ is also equal to
    \[
        \theta(k(\bar{n}_v)) = \bar{n}_v a_{1 + \|v\|^2} n_v
    \]
    and $ \kappa(\bar{n}_v) $ is equal to
    \[
        k(\bar{n}_{-v})^{-1} = \theta(k(\bar{n}_{-v})^{-1})
        = \theta(\bar{n}_v a_{1 + \|v\|^2} n_v )
        = n_v a_{\tfrac1{1 + \|v\|^2}} \bar{n}_v = k(\bar{n}_v) \ .
    \]
    Hence,
    \begin{align*}
        w k(w n_v) &= w k(\bar{n}_{(v_1, -v_2, \dotsc, -v_{n-1})}) w \\
        &= w n_{(v_1, -v_2, \dotsc, -v_{n-1})} a_{\tfrac1{1 + \|v\|^2}} \bar{n}_{(v_1, -v_2, \dotsc, -v_{n-1})} w \\
        &= \bar{n}_v a_{1 + \|v\|^2} n_v = k(\bar{n}_v) \ .
    \end{align*}
    This completes the proof of the lemma.
\end{proof}

\begin{cor}\label{Cor Iwasawa k}
    We have
    \begin{enumerate}
    \item $ k(w n_v)^{-1} = k(w n_{(-v_1, v_2, \dotsc, v_{n-1})}) $;\label{enum formulas with k eq1}
    \item $ w k(w n_v) w = k(w n_{(v_1, -v_2, \dotsc, -v_{n-1})}) $.\label{enum formulas with k eq2}
    \end{enumerate}
\end{cor}

\begin{proof}
    This corollary follows immediately from the previous lemma.

\end{proof}

\begin{lem}\label{Lem formula bar(n)_v n_v bar(n)_v = k(bar(n)_v)^2 = n_v bar(n)_v n_v}
    Let $ v \in \RR^{n-1} $ ($ n \geq 2 $). Then,
    \[
        \bar{n}_v n_v \bar{n}_v = k(\bar{n}_v) n_{(\|v\|^2 - 1)v} k(\bar{n}_v) \ .
    \]
    If $ \|v\| = 1 $, then
    \[
        \bar{n}_v n_v \bar{n}_v = k(\bar{n}_v)^2 = n_v \bar{n}_v n_v \ .
    \]
\end{lem}

\begin{proof}
    For $ v \in \RR^{n-1} $, we have
    \[
        \bar{n}_v n_v \bar{n}_v = \kappa(\bar{n}_v) a(\bar{n}_v) n(\bar{n}_v) n_v \nu(\bar{n}_v) h(\bar{n}) k(\bar{n}_v) \ .
    \]
    By Lemma \ref{Lem about Iwasawa k}, this is equal to
    \begin{align*}
        k(\bar{n}_v) h(\bar{n})^{-1} \nu(\bar{n}_v) n_v \nu(\bar{n}_v) h(\bar{n}) k(\bar{n}_v)
        &= k(\bar{n}_v) a_{1 + \|v\|^2} n_{v -\tfrac2{1+\|v\|^2}v } a_{1 + \|v\|^2}^{-1} k(\bar{n}_v) \\
        &= k(\bar{n}_v) n_{(\|v\|^2 - 1)v} k(\bar{n}_v) \ .
    \end{align*}
    The lemma follows.
\end{proof}

\begin{lem}
    Let $ n \geq 2 $ and let $ v \in \RR^{n-1} $ be of norm 1. Let $ A \in \SO(n-1) $ be such that $ A (1, 0_{1, n-2})^t = v $. Then,
    \[
        k(w n_v) w k(w n_v)
        = w n_v w n_{(v_1, -v_2, \dotsc, -v_{n-1})} w n_v
        = w m_A w m_A^{-1} \in M \ .
    \]
\end{lem}

\begin{proof}
    Let $ v $ and $A$ be as above. Then,
    \[
        w n_v w n_{(v_1, -v_2, \dotsc, -v_{n-1})} w n_v
    \]
    is equal to
    \begin{align*}
        w n_v \bar{n}_v n_v
        &= w (m_A n_{e_1} m_A^{-1}) (m_A \bar{n}_{e_1} m_A^{-1}) (m_A n_{e_1} m_A^{-1}) \\
        &= w m_A n_{e_1} \bar{n}_{e_1} n_{e_1} m_A^{-1} \ .
    \end{align*}
    A direct computation shows that $ n_{e_1} \bar{n}_{e_1} n_{e_1} $ is equal to $w$. Thus, the above is equal to $ w m_A w m_A^{-1} $.
    The lemma follows now from Lemma \ref{Lem formula bar(n)_v n_v bar(n)_v = k(bar(n)_v)^2 = n_v bar(n)_v n_v}.
\end{proof}

\begin{lem}\label{Lem lim_(t to infty) k(w n_(v(t)))}
    Let $ n \geq 3 $ and $ v(t) \in \RR^{n-1} $ be such that $ \|v(t)\| $ converges to $ \infty $ when $t$ tends to $ \infty $. Assume that $ v' := \lim_{t \to \infty} \frac{v(t)}{\|v(t)\|} $ exists and has norm 1. Let $ A \in \SO(n-1) $ be such that $ A (1, 0_{1, n-2})^t = v' $. Then, $ \lim_{t \to \infty} k(w n_{v(t)}) $ exists and is equal to
    \begin{align*}
        & k(w n_{v'}) w k(w n_{v'}) \\
        =& w n_{(v'_1, \dotsc, v'_{n-1})} w n_{(v'_1, -v'_2, \dotsc, -v'_{n-1})} w n_{(v'_1, \dotsc, v'_{n-1})} \\
        =& w m_A w m_A^{-1} \ .
    \end{align*}
    In particular, $ \lim_{t \to \infty} k(w n_{t (1, 0_{1, n-1})}) = e $.
\end{lem}

\begin{rem}
    Let $ v(t) $ be as above. Assume that there exists $ c(t) \in (0, +\infty) $ such that $ v' := \lim_{t \to \infty} \frac{v(t)}{c(t)} $ exists and has norm 1.
    Then,
    \[
        \lim_{t \to \infty} \frac{v(t)}{\|v(t)\|} = \lim_{t \to \infty} \frac{v(t)}{c(t)} \frac{c(t)}{\|v(t)\|}
    \]
    exists also, is equal to $ v' $ and has norm 1.
\end{rem}

\begin{proof}
    Let $ v(t) $, $ v' $ and $A$ be as above. Let $ \tilde{v}(t) = \tfrac{v(t)}{\|v(t)\|} $.

    Since $ w n_{\tilde{v}(t)} w n_{(\tilde{v}_1(t), -\tilde{v}_2(t), \dotsc, -\tilde{v}_{n-1}(t))} w n_{\tilde{v}(t)} $ belongs to $M$ by the previous lemma,
    \[
        k(w n_{v(t)}) = k(w n_{\tilde{v}(t)} a_{\|v(t)\|}^{-1}) 
    \]
    is equal to
    \begin{align*}
        & k(n_{(-\tilde{v}_1(t), \tilde{v}_2(t), \dotsc, \tilde{v}_{n-1}(t))} w n_{-\tilde{v}(t)} w a_{\|v(t)\|}^{-1}) w n_{\tilde{v}(t)} w n_{(\tilde{v}_1(t), -\tilde{v}_2(t), \dotsc, -\tilde{v}_{n-1}(t))} w n_{\tilde{v}(t)} \\
        =& k(a_{\|v(t)\|} (w n_{-\tilde{v}(t)} w) a_{\|v(t)\|}^{-1}) w n_{\tilde{v}(t)} w n_{(\tilde{v}_1(t), -\tilde{v}_2(t), \dotsc, -\tilde{v}_{n-1}(t))} w n_{\tilde{v}(t)} \ .
    \end{align*}
    Since $ w n_{-\tilde{v}(t)} w \in \bar{N} $, since $ \lim_{t \to \infty} \|\tilde{v}(t)\| = 1 < \infty $ and since $ \lim_{t \to \infty} \|v(t)\| = +\infty $, $ a_{\|v(t)\|} (w n_{-\tilde{v}(t)} w) a_{\|v(t)\|}^{-1} $ converges to $e$ when $t$ tends to $ \infty $, the above converges to
    \[
        w n_{v'} w n_{(v'_1, -v'_2, \dotsc, -v'_{n-1})} w n_{v'}
    \]
    when $t$ tends to $ \infty $. The lemma follows by the previous one.
\end{proof}

\newpage

\subsubsection{The diffeomorphism \texorpdfstring{$ \Upsilon_{1, t, s} $}{Upsilon\_(1, t, s)}}\label{sssec:diffeo Upsilon_(1, t, s)}

The map $ \Upsilon_{1, t, s}^G \colon N \to N^\Gamma A $, which we define in this section, will satisfy
\[
    \kappa(n_t w) n_{-t} n = n' \Upsilon_{1, t, s}^G(n) k' \qquad (n \in N)
\]
for some $ n' \in N_\Gamma $ and $ k' \in K $ depending on $t$ and $n$.
Moreover, we show that it induces a diffeomorphism.

\medskip

As $ w = w^{-1} $, we have
\begin{equation}\label{eq kappa(n_t w)}
    \kappa(n_t w) = n_t w a(n_t w)^{-1} (a(n_t w) n(n_t w)^{-1} a(n_t w)^{-1})
    = n_t w h(w n_{-t}) n_t \ .
\end{equation}
For $ n \in N $, we have
\begin{align*}
    n_t w h(w n_{-t}) n a_s &= n_t h(w n_{-t})^{-1} a_s^{-1} w n^{a_{\frac1s}} \\
    &= n_t h(w n_{-t})^{-1} a_s^{-1} \nu(w n^{a_{\frac1s}}) h(w n^{a_{\frac1s}}) k(w n^{a_{\frac1s}}) \\
    &= n_t \nu(w n^{a_{\frac1s}})^{h(w n_{-t})^{-1} a_s^{-1}} h(w n_{-t})^{-1} a_s^{-1} h(w n^{a_{\frac1s}}) k(w n^{a_{\frac1s}}) \ .
\end{align*}
The canonical projection of the horosphere
\[
    (\kappa(n_t w) N \kappa(n_t w)^{-1}) \kappa(n_t w) a_s K
\]
on $ N^{\Gamma} A $ gives the preimage of a horoball under the natural map from $ N^{\Gamma} A $ to $ N^{\Gamma} A K $ which we denote by $ B^G_{t, s} $.

Let us now look at an example:
When $X$ is the upper half-space, then the projection of the 3-dimensional sphere on the upper half-plane gives us a disc (the hatched area):
\begin{center}
    \includegraphics[width = 12cm]{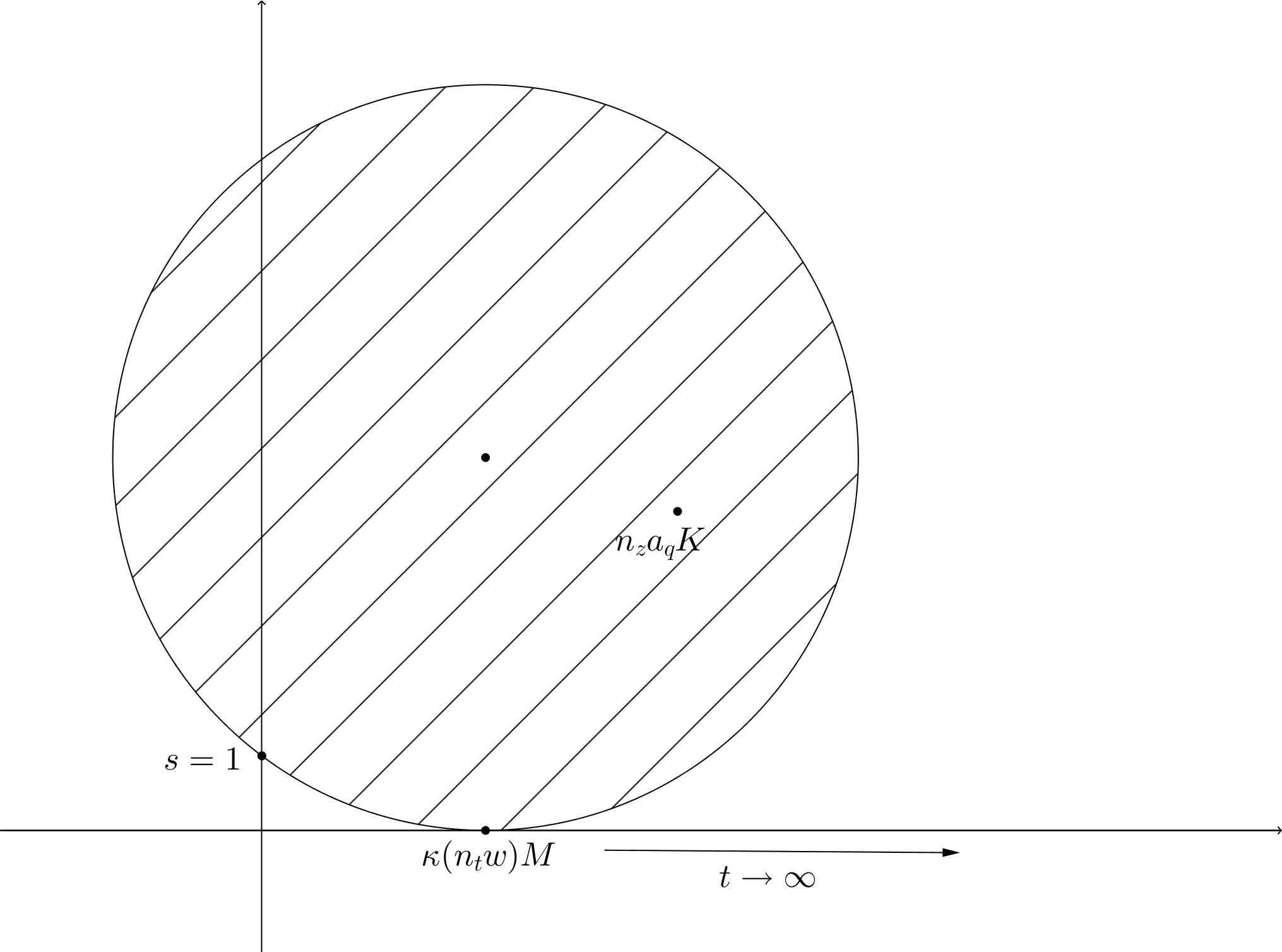}
\end{center}
The disc in the upper half-plane converges to the area $ \{ x + i y \mid x \in \RR, \, y > 1 \} $.

For $ g \in G $, set $ \nu_{N_\Gamma}(g) = p_{N_\Gamma}(\nu(g)) $ and $ \nu_{N^\Gamma}(g) = p_{N^\Gamma}(\nu(g)) $.

Define
\begin{align*}
    \Upsilon_{1, t, s}^G \colon & N \to N^\Gamma A, \\
    & 
    n \mapsto \nu_{N^\Gamma}(n_t w h(w n_{-t}) n a_s) h(n_t w h(w n_{-t}) n a_s) \\
    & \qquad = \nu_{N^\Gamma}(n_t \nu(w n^{a_{\frac1s}})^{h(w n_{-t})^{-1} a_s^{-1}}) h(w n_{-t})^{-1} a_s^{-1} h(w n^{a_{\frac1s}}) \ .
\end{align*}

By construction, $ \im(\Upsilon_{1, t, s}^G) = B^G_{t, s} $.

Set $ \Upsilon_{1, t, s} = (n_{\cdot} a_{\cdot})^{-1} \circ \Upsilon_{1, t, s}^G \circ n_{\cdot} \colon \RR^{n-1} \to \RR^{n-k-1} \times (0, +\infty) $.

Set $ B_{t, s} = (n_{\cdot} a_{\cdot})^{-1}(B^G_{t, s}) $. So, $ B_{t, s} = \im(\Upsilon_{1, t, s}) $.

Since $ n_t w h(w n_{-t}) a_s = n_t a_{\frac1s(1+t^2)} w $, $ (t, 0_{1, n-k-2}, \frac1s(1+t^2)) \in B_{t, s} $.

$ B_{t, s} $ is the ball in $ \RR^{n-k-1} \times (0, +\infty) $ whose boundary contains $ (t, 0_{1, n-k-1}) $ and is tangent to $ \RR^{n-k-1} \times \{ 0 \} $.
Thus, the center $ C_{t, s} $ of $ B_{t, s} $ is given by $ C_{t, s} = (t, 0_{1, n-k-2}, c_{t, s}) $, where $ c_{t, s} := \frac1{2s}(1+t^2) $.
Since 
\begin{align*}
    & \nu_{N^\Gamma}(n_t \nu(w n_v^{a_{\frac1s}})^{h(w n_{-t})^{-1} a_s^{-1}}) \\
    = & \nu_{N^\Gamma}(n_t \nu(w n_{\frac1s v})^{(a_{\frac1{t^2+1}})^{-1} a_s^{-1}}) \\
    = & \nu_{N^\Gamma}(n_t (n_{\frac1{1+\|\frac1s v\|^2} \cdot \frac1s (-v_1, v_{2..n-1})})^{a_{\frac{t^2+1}{s}}}) \\
    = & \nu_{N^\Gamma}(n_t n_{\frac{t^2+1}{s} \cdot \frac1{1+\frac1{s^2} \|v\|^2} \cdot \frac1s (-v_1, v_{2..n-1})}) \\
    = & n^{N^\Gamma}_{(t - \frac{t^2+1}{s^2+\|v\|^2} v_1, \frac{t^2+1}{s^2+\|v\|^2} v_{2..n-k-1})}
\end{align*}
and
\[
    h(w n_{-t})^{-1} a_s^{-1} h(w n_v^{a_{\frac1s}}) = a_{t^2+1} a_{\frac1s} h(w n_{ \frac1s v})
    = a_{\frac{t^2+1}s} a_{\frac1{1 + \frac1{s^2}\|v\|^2 }}
    = a_{\frac{t^2+1}{s + \frac1s \|v\|^2 }} \ ,
\]
we have
\[
    \Upsilon_{1, t, s}(v)
    = \Big( \big(t - \frac{t^2+1}{s^2+\|v\|^2} v_1, \frac{t^2+1}{s^2+\|v\|^2} v_{2..n-k-1} \big), \frac{t^2+1}{s + \frac1s \|v\|^2} \Big) \ .
\]
Define
\[
    S^{k-1} \times \RR^{n-1} \to \RR^{n-1}, \quad u.(v_1, \dotsc, v_{n-1}) \mapsto \left( v_1, \dotsc, v_{n-k-1}, \sqrt{ \sum_{j=n-k}^{n-1} v_j^2 } \ u \right) \ .
\]
As $ \Upsilon_{1, t, s} $ is uniquely defined if we know $ v_1, \dotsc, v_{n-k-1} $ and $ \sqrt{ \sum_{j=n-k}^{n-1} v_j^2 } $,
\[
    \Upsilon_{1, t, s}(u.v) = \Upsilon_{1, t, s}(v)
\]
for all $ u \in S^{k-1} $ and $ v \in \RR^{n-1} $.
Let $ x \in \RR^{n-1} \smallsetminus \{0\} $. Then, there is a unique $ u \in S^{k-1} $ such that $ u.x \in \RR^{n-k-1} \times \{ (r, 0_{1, k-1}) \mid r \in (0, \infty) \} $.
By abuse of notation, we will denote the map
\[
    \RR^{n-k-1} \times (0, \infty) \to \RR^{n-k-1} \times (0, \infty), \; (v, r) \mapsto \Upsilon_{1, t, s}((v, r, 0_{1, k-1}))
\]
also by $ \Upsilon_{1, t, s} $. It is given by
\[
    \Upsilon_{1, t, s}(v, r)
    = \Big( \big(t - \frac{t^2+1}{s^2+\|v\|^2+r^2} v_1, \frac{t^2+1}{s^2+\|v\|^2+r^2} v_{2..n-k-1} \big), \frac{t^2+1}{s + \frac1s (\|v\|^2 + r^2)} \Big) \ .
\]

Let us now try to compute the inverse of $ \Upsilon_{1, t, s} $. We have
\begin{align*}
    q := \frac{t^2+1}{s + \frac1s (\|v\|^2 + r^2)} \in (0, +\infty) \iff & q (s + \frac1s \big(\|v\|^2 + r^2)\big) = t^2+1 \\
    \iff & \|v\|^2 + r^2 = \frac{ s(t^2+1 - q s) }{ q } \ .
\end{align*}
Thus,
\[
    z_1 := t - \frac{t^2+1}{s^2+\|v\|^2+r^2} v_1 \iff t - z_1  = \frac{q}{s} \cdot v_1
    \iff v_1 = \frac{s(t - z_1)}{q}
\]
and, for $ j = 2, \dotsc, n-k-1 $,
\[
    z_j := \frac{t^2+1}{s^2+\|v\|^2+r^2} v_j \iff z_j = \frac{q}{s} \cdot v_j
    \iff v_j = \frac{s z_j}{q} \ .
\]
Hence,
\[
    r = \sqrt{ \frac{ s(t^2+1 - q s) }{ q } - \|v\|^2 }
    = \sqrt{ -s^2 + \frac{ s(t^2+1) }{ q } - \frac{s^2}{q^2} (t - z_1)^2 - \frac{s^2}{q^2} \sum_{j=2}^{n-k-1} z_j^2 } \ .
\]
Note that this is equal to
\[
    \frac{s}{q} \cdot \sqrt{ c_{t, s}^2 - (q - c_{t, s})^2 - (t - z_1)^2 - \sum_{j=2}^{n-k-1} z_j^2 } \ .
\]
For $ (z, q) \in \im(\Upsilon_{1, t, s}) = B_{t, s} $, define
\begin{multline*}
    \Psi_{1, t, s}(z, q) = \Bigg( \frac{s(t - z_1)}{q}, \frac{s z_2}{q}, \dotsc, \frac{s z_{n-k-1}}{q}, \\
    \sqrt{ -s^2 + \frac{ s(t^2+1) }{ q } - \frac{s^2}{q^2} (t - z_1)^2 - \frac{s^2}{q^2} \sum_{j=2}^{n-k-1} z_j^2 } \Bigg) \in \RR^{n-k-1} \times (0, +\infty) \ .
\end{multline*}
Thus, $ \Upsilon_{1, t, s} \colon \RR^{n-k-1} \times (0, \infty) \to B_{t, s} $ is bijective and $ \Psi_{1, t, s} $ is its inverse. As $ \Upsilon_{1, t, s} $ and $ \Psi_{1, t, s} $ are smooth, $ \Upsilon_{1, t, s} $ is a diffeomorphism.

\newpage

\subsubsection{The diffeomorphism \texorpdfstring{$ \Upsilon_{2, t, s} $}{Upsilon\_(2, t, s)}}\label{ssec:diffeo Upsilon_(2, t, s)}

The diffeomorphism $ \Upsilon_{2, t}^G \colon N^\Gamma A \to N^\Gamma A $ (inducing a diffeomorphism $ \Upsilon_{2, t} \colon \RR^{n-k-1} \times (0, +\infty) \to \RR^{n-k-1} \times (0, +\infty) $), which we define in this section, will satisfy
\[
    n a = \kappa(n_t w) n_{-t} (\Upsilon_{2, t}^G)^{-1}(na) k' \qquad (n \in N^\Gamma, a \in A)
\]
for some $ k' \in K $ depending on $t$ and $n$.

As $ w = w^{-1} $, we have
\[
    \kappa(n_t w) = n_t w a(n_t w)^{-1} (a(n_t w) n(n_t w)^{-1} a(n_t w)^{-1})
    = n_t w h(w n_{-t}) n_t \ .
\]
For $ n \in N^\Gamma $, we have 
\begin{align*}
    n_t w h(w n_{-t}) n a_q &= n_t h(w n_{-t})^{-1} a_q^{-1} w n^{a_{\frac1q}} \\
    &= n_t h(w n_{-t})^{-1} a_q^{-1} \nu(w n^{a_{\frac1q}}) h(w n^{a_{\frac1q}}) k(w n^{a_{\frac1q}}) \\
    &= n_t \nu(w n^{a_{\frac1q}})^{h(w n_{-t})^{-1} a_q^{-1}} h(w n_{-t})^{-1} a_q^{-1} h(w n^{a_{\frac1q}}) k(w n^{a_{\frac1q}}) \ .
\end{align*}
Let
\[
    \Upsilon_{2, t}^G \colon N^\Gamma A \to N^\Gamma A ,
    \quad n a \mapsto \nu(n_t w h(w n_{-t}) n a) h(n_t w h(w n_{-t}) n a) \ .
\]
Set $ \Upsilon_{2, t} = (n^{N^\Gamma}_{\cdot} a_{\cdot})^{-1} \circ \Upsilon_{2, t}^G \circ (n^{N^\Gamma}_{\cdot} a_{\cdot}) \colon \RR^{n-k-1} \times (0, +\infty) \to \RR^{n-k-1} \times (0, +\infty) $.

Since
\[
    \nu(n_t \nu(w n_z^{a_{\frac1q}})^{h(w n_{-t})^{-1} a_q^{-1}})
    = n^{N^\Gamma}_{(t - \frac{t^2+1}{q^2+\|z\|^2} z_1, \frac{t^2+1}{q^2+\|z\|^2} z_{2..n-k-1})}
\]
and
\[
    h(w n_{-t})^{-1} a_q^{-1} h(w n_z^{a_{\frac1q}})
    = a_{\frac{t^2+1}{q + \frac1q \|z\|^2 }} \ ,
\]
we have
\[
    \Upsilon_{2, t}(z, q)
    = \Big( \big(t - \frac{t^2+1}{q^2+\|z\|^2} z_1, \frac{t^2+1}{q^2+\|z\|^2} z_{2..n-k-1} \big), \frac{t^2+1}{q + \frac1q \|z\|^2} \Big) \ .
\]
Let us now try to compute the inverse of $ \Upsilon_{2, t} $. We have
\[
    q' := \frac{t^2+1}{q + \frac1q \|z\|^2} \iff q' (q + \frac1q \|z\|^2) = t^2+1 \iff \|z\|^2 = \frac{ q(t^2+1 - q' q) }{ q' } \ .
\]
Thus,
\[
    z'_1 := t - \frac{t^2+1}{q^2+\|z\|^2} z_1 \iff t - z'_1  = \frac{q'}{q} \cdot z_1
    \iff z_1 = \frac{q(t - z'_1)}{q'}
\]
and, for $ j = 2, \dotsc, n-k-1 $,
\[
    z'_j := \frac{t^2+1}{q^2+\|z\|^2} z_j \iff z'_j = \frac{q'}{q} \cdot z_j
    \iff z_j = \frac{q z'_j}{q'} \ .
\]
Hence,
\begin{align*}
    & \frac{ q(t^2+1 - q' q) }{ q' } = \frac{q^2(t - z'_1)^2}{{q'}^2}
    + \sum_{j=2}^{n-k-1} \frac{q^2 {z'_j}^2}{{q'}^2} \\
    \iff & q'(\frac{t^2+1}q - q') = t^2 - 2 t z'_1 + \|z'\|^2 \\
    \iff & q = \frac{q'(t^2+1)}{ t^2 - 2 t z'_1 + \|z'\|^2 + {q'}^2 } \\
    \iff & q = \frac{q'(t^2+1)}{ (t - z'_1)^2 + \sum_{j=2}^{n-k-1} {z'_j}^2 + {q'}^2 } \ .
\end{align*}
So,
\[
    z_1 = \frac{(t - z'_1)(t^2+1)}{ (t - z'_1)^2 + \sum_{j=2}^{n-k-1} {z'_j}^2 + {q'}^2 }
\]
and, for $ j = 2, \dotsc, n-k-1 $,
\[
    z_j  = \frac{z'_j (t^2+1)}{ (t - z'_1)^2 + \sum_{j=2}^{n-k-1} {z'_j}^2 + {q'}^2 } \ .
\]
For $ (z', q') \in \im(\Upsilon_{2, t}) $, define $ \Psi_{2, t}(z', q') $ by
$
    ( z, q ) \in \RR^{n-k-1} \times (0, +\infty)
$.

Thus, $ \Upsilon_{2, t} \colon \RR^{n-k-1} \times (0, +\infty) \to \RR^{n-k-1} \times (0, +\infty) $ is a bijective map and $ \Psi_{2, t} $ is its inverse. As $ \Upsilon_{2, t} $ and $ \Psi_{2, t} $ are smooth, $ \Upsilon_{2, t} $ is a diffeomorphism.

Since
\begin{align*}
    & (z', q') = \Upsilon_{2, t}(z, q) \in B_{t, s} \\
    \iff & (z'_1 - t)^2 + \sum_{j=1}^{n-k-1} {z'_j}^2 + (q' - c_{t, s})^2 < c_{t, s}^2 \\
    \iff & \big(\frac{t^2+1}{q^2+\|z\|^2}\big)^2 \|z\|^2
    + (\frac{t^2+1}{q + \frac1q \|z\|^2} - c_{t, s})^2 < c_{t, s}^2 \\
    \iff & \frac{ \|z\|^2 }{(q^2+\|z\|^2)^2}
    + (\frac{q}{q^2 + \|z\|^2} - \frac1{2s})^2 < \frac1{4s^2} \\
    \iff & \frac{ \|z\|^2 + q^2 }{(q^2+\|z\|^2)^2} < \frac{q}{s(q^2 + \|z\|^2)} \\
    \iff & q > s \ ,
\end{align*}
$ \Upsilon_{2, t, s} := \restricted{ \Upsilon_{2, t} }{ \RR^{n-k-1} \times (s, +\infty) } \colon \RR^{n-k-1} \times (s, +\infty) \to B_{t, s} $ is also a diffeomorphism.

Thus, $ \Upsilon_{2, t, s}^G := \restricted{ \Upsilon_{2, t}^G }{ N^\Gamma A_{>s} } $ is a diffeomorphism from $ N^\Gamma A_{>s} $ to $ B_{t, s}^G $.

\newpage

\subsubsection{The diffeomorphism \texorpdfstring{$ \Upsilon_{t, s} $}{Upsilon\_(t, s)}}\label{sssec:diffeo Upsilon_(t, s)}

For $ \tilde{z} \in \RR^{n-k-1} $, define $ T_{\tilde{z}} \colon \RR^{n-k-1} \times (0, + \infty) \to \RR^{n-k-1} \times (0, +\infty) , \; (z, q) \mapsto (z + \tilde{z}, q) $ (diffeomorphism).
Since $ \Upsilon_{1, t, s} $ is a diffeomorphism from $ \RR^{n-k-1} \times (0, +\infty) $ to $ B_{t, s} $ and $ \restricted{\Upsilon_{2, t, s}}{ \RR^{n-k-1} \times (s, +\infty)} $ is a diffeomorphism from $ \RR^{n-k-1} \times (s, +\infty) $ to $ B_{t, s} $,
\[
    \Upsilon_{t, s} := \restricted{ T_{(-t, 0_{1, n-k-2})} \circ \Upsilon_{1, t, s}^{-1} \circ \Upsilon_{2, t, s} \circ T_{(t, 0_{1, n-k-2})} }{ \RR^{n-k-1} \times (s, +\infty)}
\]
is a diffeomorphism from $ \RR^{n-k-1} \times (s, +\infty) $ to $ \RR^{n-k-1} \times (0, +\infty) $.

Let $ (v_1, \dotsc, v_{n-k-1}, r) = \Upsilon_{t, s}(z, q) $ ($ z \in \RR^{n-k-1} $, $ q > s $). Direct computation shows that $ (v_1, \dotsc, v_{n-k-1}, r) $ is equal to
\[
    ( -t + \frac{s (z_1 + t)}q, \dotsc, \frac{s z_{n-k-1}}q, \sqrt{s} \cdot \sqrt{1 - \frac{s}{q}} \cdot \sqrt{ q + \frac1q \|(z_1 + t, z_2, \dotsc, z_{n-k-1})\|^2 } \ ) \ .
\]
Put $ r(t, s, z, q) = \sqrt{s} \cdot \sqrt{1 - \frac{s}{q}} \cdot \sqrt{ q + \frac1q \|(z_1 + t, z_2, \dotsc, z_{n-k-1})\|^2 } $.
\\ Let $ (z', q') = \Upsilon_{1, t, s}(T_{(t, 0_{1, n-k-2})}((v, r))) $.
Let us treat now first the case when $ k \geq 2 $.
Then,
\[
    \kappa(n_t w) n_{(v_1, \dotsc, v_{n-k-1}, r \Phi_k(\psi_0, \psi_1, \dotsc, \psi_{k-2}))}
\]
is equal to
\begin{align*}
    & \nu_{N_\Gamma}(w n_{(v_1 - t, \dotsc, v_{n-k-1}, r \Phi_k(\psi_0, \psi_1, \dotsc, \psi_{k-2}))}^{a_{\frac1s}})^{h(w n_{-t})^{-1} a_s^{-1}} n_{z'} a_{q'} \\
    & \qquad k(w n_{(v_1 - t, \dotsc, v_{n-k-1}, r \Phi_k(\psi_0, \psi_1, \dotsc, \psi_{k-2}))}^{a_{\frac1s}}) \\
    =& \nu_{N_\Gamma}(w n_{(v_1 - t, \dotsc, v_{n-k-1}, r \Phi_k(\psi_0, \psi_1, \dotsc, \psi_{k-2}))}^{a_{\frac1s}})^{h(w n_{-t})^{-1} a_s^{-1}} \kappa(n_t w) n_z a_q \\
    & \qquad k(w n_{(z_1 + t, z_2, \dotsc, z_{n-k-1})}^{a_{\frac1{q}}})^{-1} k(w n_{(v_1, \dotsc, v_{n-k-1}, r \Phi_k(\psi_0, \psi_1, \dotsc, \psi_{k-2}))}^{a_{\frac1s}}) \\
    =& n^{N_\Gamma}_{\big(\frac{t^2+1}{s (q + \frac1q \|(z_1 + t, z_2, \dotsc, z_{n-k-1})\|^2)} \cdot r(t, s, z, q) \Phi_k(\psi_0, \psi_1, \dotsc, \psi_{k-2}) \big)} \kappa(n_t w) n_z a_q k(w n_{\frac{ (z_1 + t, z_2, \dotsc, z_{n-k-1})}{ q }})^{-1} \\
    & \qquad k(w n_{(\frac{ z_1 + t }{ q}, \frac{ z_2 }{ q}, \dotsc, \frac{ z_{n-k-1} }{ q}, \frac1{\sqrt{s}} \cdot \sqrt{1 - \frac{s}{q}} \cdot \sqrt{ q + \frac1q \|(z_1 + t, z_2, \dotsc, z_{n-k-1})\|^2 } \Phi_k(\psi_0, \psi_1, \dotsc, \psi_{k-2}))})
\end{align*}
as
\begin{align*}
    &\nu_{N_\Gamma}(w n_{(v_1 - t, v_2, \dotsc, v_{n-k-1}, r \Phi_k(\psi_0, \psi_1, \dotsc, \psi_{k-2}))}^{a_{\frac1s}})^{h(w n_{-t})^{-1} a_s^{-1}} \\
    =& n^{N_\Gamma}_{\big(\frac{t^2+1}{s (q + \frac1q \|(z_1 + t, z_2, \dotsc, z_{n-k-1})\|^2)} \cdot r(t, s, z, q) \Phi_k(\psi_0, \psi_1, \dotsc, \psi_{k-2}) \big)}
\end{align*}
and as
\begin{align*}
    &k(w n_{(v_1 - t, v_2, \dotsc, v_{n-k-1}, r \Phi_k(\psi_0, \psi_1, \dotsc, \psi_{k-2}))}^{a_{\frac1s}}) \\
    =& k(w n_{(\frac{ z_1 + t}{ q}, \frac{ z_2}{q}, \dotsc, \frac{ z_{n-k-1} }{q}, \frac1s r(t, s, z, q) \Phi_k(\psi_0, \psi_1, \dotsc, \psi_{k-2}))}) \\
    =& k(w n_{(\frac{ z_1 + t}{ q}, \frac{ z_2}{q}, \dotsc, \frac{ z_{n-k-1} }{q}, \frac1{\sqrt{s}} \cdot \sqrt{1 - \frac{s}{q}} \cdot \sqrt{ q + \frac1q \|(z_1 + t, z_2, \dotsc, z_{n-k-1})\|^2 } \Phi_k(\psi_0, \psi_1, \dotsc, \psi_{k-2}))}) \ .
\end{align*}
For $ u \in [0, 1] $, set
\begin{align*}
    m_u &:= \lim_{t \to \infty} k(w n_{t (\sqrt{u}, 0_{1, n-k-2}, \sqrt{ 1 - u }, 0_{1, k-1})}) \\
    &= w n_{(\sqrt{u}, 0_{1, n-k-2}, \sqrt{ 1 - u }, 0_{k-1})} w n_{(\sqrt{u}, 0_{1, n-k-2}, -\sqrt{ 1 - u }, 0_{1, k-1})} w n_{(\sqrt{u}, 0_{1, n-k-2}, \sqrt{ 1 - u }, 0_{1, k-1})} \\
    &= \begin{pmatrix}
        1 & 0 & 0 & 0 & 0 & 0_{1, k-1} \\
        0 & 1 & 0 & 0 & 0 & 0_{1, k-1} \\
        0 & 0 & 2u - 1 & 0_{1, n-k-2} & 2\sqrt{u(1-u)} & 0_{1, k-1} \\
        0_{n-k-2, 1} & 0_{n-k-2, 1} & 0_{n-k-2, 1} & I_{n-k-2} & 0_{n-k-2, 1} & 0_{1, k-1} \\
        0 & 0 & -2\sqrt{u(1-u)} & 0_{1, n-k-2} & 2u - 1 & 0_{1, k-1} \\
        0_{k-1, 1} & 0_{k-1, 1} & 0_{k-1, 1} & 0_{k-1, 1} & 0_{k-1, 1} & I_{k-1, k-1}
    \end{pmatrix} \ .
\end{align*}
Note that $ \lim_{t \to \infty} k(w n_{t (\sqrt{u}, 0_{1, n-k-2}, -\sqrt{ 1 - u }, 0_{1, k-1})}) $ is equal to $ m_u^{-1} $.

Let $ Y_0, Y_1, \dotsc, Y_{k-2} \in \so(k) $ be as in the proof of Lemma \ref{Lem Relation between spherical coordinates and so(d)} (with $ d = k $).
\\ Set $ A(\psi_0, \dotsc, \psi_{k-2}) = \prod_{j=0}^{k-2} \exp(-\psi_j Y_j) \in \SO(k) $. For $ A \in \SO(k) $, set
\[
    m_A^{M_k} = m_{\diag(I_{n-k-1}, A)} \ .
\]
By Lemma \ref{Lem lim_(t to infty) k(w n_(v(t)))} and as $ w m w = m $ for all $ m \in M_k $,
\begin{align*} 
    & k(w n_{\frac{(z_1 + t, z_2, \dotsc, z_{n-k-1})}q})^{-1} \\
    & \qquad k(w n_{(\frac{ z_1 + t }{q}, \frac{ z_2 }{q}, \dotsc, \frac{ z_{n-k-1} }{q}, \frac1{\sqrt{s}} \cdot \sqrt{1 - \frac{s}{q}} \cdot \sqrt{ q + \frac1q \|(z_1 + t, z_2, \dotsc, z_{n-k-1})\|^2 } \Phi_k(\psi_0, \psi_1, \dotsc, \psi_{k-2}))}) \\
    =& k(w n_{\frac{(z_1 + t, z_2, \dotsc, z_{n-k-1})}q})^{-1} m^{M_k}_{A(\psi_0, \dotsc, \psi_{k-2})} \\
    & \qquad k(w n_{(\frac{ z_1 + t }{q}, \frac{ z_2 }{q}, \dotsc, \frac{ z_{n-k-1} }{q}, \frac1{\sqrt{s}} \cdot \sqrt{1 - \frac{s}{q}} \cdot \sqrt{ q + \frac1q \|(z_1 + t, z_2, \dotsc, z_{n-k-1})\|^2 })}) (m^{M_k}_{A(\psi_0, \dotsc, \psi_{k-2})})^{-1}
\end{align*}
converges to $ m^{M_k}_{A(\psi_0, \dotsc, \psi_{k-2})} m_{\frac{s}{q}} (m^{M_k}_{A(\psi_0, \dotsc, \psi_{k-2})})^{-1} $.

\medskip

If $ k = 1 $, then
$
    \kappa(n_t w) n_{(v_1, \dotsc, v_{n-k-1}, \pm r)}
$
is equal to
\begin{multline*}
    n^{N_\Gamma}_{\pm \big(\frac{t^2+1}{s (q + \frac1q \|(z_1 + t, z_2, \dotsc, z_{n-k-1})\|^2)} \cdot r(t, s, z, q) \big)} \kappa(n_t w) n_z a_q k(w n_{\frac{ (z_1 + t, z_2, \dotsc, z_{n-k-1})}{ q }})^{-1} \\
    k(w n_{(\frac{ z_1 + t }{ q}, \frac{ z_2 }{ q}, \dotsc, \frac{ z_{n-k-1} }{ q}, \pm \frac1{\sqrt{s}} \cdot \sqrt{1 - \frac{s}{q}} \cdot \sqrt{ q + \frac1q \|(z_1 + t, z_2, \dotsc, z_{n-k-1})\|^2 } )})
\end{multline*}
and
\[
    k(w n_{(\frac{ z_1 + t }{q}, \frac{ z_2 }{q}, \dotsc, \frac{ z_{n-k-1} }{q}, \pm \frac1{\sqrt{s}} \cdot \sqrt{1 - \frac{s}{q}} \cdot \sqrt{ q + \frac1q \|(z_1 + t, z_2, \dotsc, z_{n-k-1})\|^2 })})
\]
converges to $ m_{\frac{s}{q}} $ respectively to $ m_{\frac{s}{q}}^{-1} $.

\medskip

Let us compute now $ |\det(\d(\Upsilon_{1, t, s}^{-1} \circ \Upsilon_{2, t, s} )(z, q))| = |\det(\d(\Psi_{1, t, s} \circ \Upsilon_{2, t, s})(z, q))| $ for $ z \in \RR^{n-k-1} $ and $ q \in (s, +\infty) $:
\[
    |\det\begin{pmatrix}
        \diag(\frac{s}q, \dotsc, \frac{s}q) & (-\frac{s z_1}{q^2}, \dotsc, -\frac{s z_{n-k-1}}{q^2})^t \\
        (\sqrt{s} \cdot \sqrt{1 - \frac{s}{q}} \cdot \frac{ \frac{2 z_1}q }{ 2\sqrt{ q + \frac1q \|z\|^2 } }, \dotsc, \sqrt{s} \cdot \sqrt{1 - \frac{s}{q}} \cdot \frac{ \frac{2 z_{n-k-1}}q }{ 2\sqrt{ q + \frac1q \|z\|^2 } })
        & \frac{ s - \frac{s}{q^2} \|z\|^2 + \frac{2 s^2}{q^3} \cdot \|z\|^2 }{2 \sqrt{s} \cdot \sqrt{1 - \frac{s}{q}} \cdot \sqrt{ q + \frac1q \|z\|^2 } }
    \end{pmatrix}| \ .
\]
This is again equal to
\begin{multline*}
    \big(\frac{s}{q}\big)^{n-k-1} \cdot \frac{\sqrt{s}}{2 \sqrt{1 - \frac{s}{q}} \cdot \sqrt{ q + \frac1q \|z\|^2 }} \cdot \\
    |\det\begin{pmatrix}
        I_{n-k-1} & (-\frac{z_1}{q}, \dotsc, -\frac{z_{n-k-1}}{q})^t \\
        (\big(1 - \frac{s}{q}\big) \cdot \frac{2 z_1}q, \dotsc, \big(1 - \frac{s}{q}\big) \cdot \frac{2 z_{n-k-1}}q)
        & 1 - \frac1{q^2} \|z\|^2 + \frac{2 s}{q^3} \cdot \|z\|^2
    \end{pmatrix}| \ .
\end{multline*}
Set
\[
    d_2 = \det\begin{pmatrix}
        1 & -\frac{z_1}{q} \\
        \big(1 - \frac{s}{q}\big) \cdot \frac{2 z_1}q &
        1 - \frac1{q^2} \|z\|^2 + \frac{2 s}{q^3} \cdot \|z\|^2
    \end{pmatrix}
\]
and
\[
    d_m = \det\begin{pmatrix}
        I_{m-1} & (-\frac{z_1}{q}, \dotsc, -\frac{z_{m-1}}{q})^t \\
        (\big(1 - \frac{s}{q}\big) \cdot \frac{2 z_1}q, \dotsc, \big(1 - \frac{s}{q}\big) \cdot \frac{2 z_{m-1}}q)
        & 1 - \frac1{q^2} \|z\|^2 + \frac{2 s}{q^3} \cdot \|z\|^2
    \end{pmatrix}
\]
for $ m = 3, \dotsc, n-k $.
Then,
\begin{align*}
    d_2
    &= 1 - \frac1{q^2} \|z\|^2 + \frac{2 s}{q^3} \cdot \|z\|^2 + \frac{2 z_1^2}{q^2} \cdot \big(1 - \frac{s}{q}\big) \\
    &= 1 - \frac{\|z\|^2}{q^2} + \frac{2 s}{q^3} \cdot (\|z\|^2 - z_1^2) + \frac{2 z_1^2}{q^2}
\end{align*}
and
\begin{align*} 
    d_m &= d_{m-1} - \big(1 - \frac{s}{q}\big) \cdot \frac{2 z_{m-1}}q \cdot
    \det\begin{pmatrix}
        I_{m-2} & (-\frac{z_1}{q}, \dotsc, -\frac{z_{m-2}}{q})^t \\
        0_{1, m-2} & -\frac{z_{m-1}}{q}
    \end{pmatrix} \\
    &= d_{m-1} + \big(1 - \frac{s}{q}\big) \cdot \frac{2 z_{m-1}^2}{q^2} \\
    &= \dotsm = d_2 + \frac{2 \cdot \big(1 - \frac{s}{q}\big)}{q^2} \cdot \sum_{j=2}^{m-1} z_j^2 \\
    &= 1 - \frac{\|z\|^2}{q^2} + \frac{2 s}{q^3} \cdot (\|z\|^2 - z_1^2) + \frac{2 z_1^2}{q^2} + \frac{2 \big(1 - \frac{s}{q}\big)}{q^2} \cdot \sum_{j=2}^{m-1} z_j^2 \\
    &= 1 + \frac{\|z\|^2}{q^2} + \frac{2 s}{q^3}( \|z\|^2 - \sum_{j=1}^{m-1} z_j^2 ) \ .
\end{align*}
Thus, $ |\det(\d\Upsilon_{t, s}(z, q))| $ is equal to
\begin{multline*}
    \big(\frac{s}{q}\big)^{n-k-1} \cdot \frac{\sqrt{s}}{2 \sqrt{1 - \frac{s}{q}} \cdot \sqrt{ q + \frac1q \|(z_1 + t, z_2, \dotsc, z_{n-k-1})\|^2 }} \cdot |d_{n-k}| \\
    = s^{n-k-\frac12} q^{-n+k} \cdot \frac{ \sqrt{ q + \frac1q \|(z_1 + t, z_2, \dotsc, z_{n-k-1})\|^2 } }{2 \sqrt{1 - \frac{s}{q}}} \ .
\end{multline*}
Put 
\[
    \mu(t, s, z, q) = r(t, s, z, q)^{k-1} \cdot |\det(\d\Upsilon_{t, s}(z, q))| \ .
\]
This is equal to
\[
    \frac{s^{n - \frac{k}2 - 1}(q + \frac1q \|(z_1 + t, z_2, \dotsc, z_{n-k-1})\|^2)^{\frac{k}2}}{2 q^{n-k} } \cdot \big(1 - \frac{s}{q}\big)^{\frac{k}2-1} \ .
\]
Thus,
\begin{equation}\label{eq limit mu/t^k}
    \lim_{t \to \infty} \frac{\mu(t, s, z, q)}{t^k} = \tfrac1{2s} (\tfrac{s}{q})^{n - \frac{k}2}(1 - \tfrac{s}q)^{\frac{k}2 - 1}
\end{equation}
for all $ z \in \RR^{n-k-1} $ and $ q > s $ fixed.

\newpage

\subsubsection{The formula for left \texorpdfstring{$ N_\Gamma M_\Gamma $}{N\_Gamma M\_Gamma}-invariant Schwartz functions}

Let $ g \in G $ and let $ f $ be a left $ N_\Gamma M_\Gamma $-invariant Schwartz function near $eP$. Let $T$ (depending on $g$) be sufficiently large so that for all $ t \geq T $ there exists an admissible generalised standard Siegel set $ \Sfrak $ with respect to $P$ (depending on $g$ and $t$) in $G$ such that $ U_{\kappa(n_t w), g} $ (see Definition \ref{Def of U_(g, h)}) is contained in $ \Sfrak' $. Thus, $ f^\Omega(\kappa(n_t w), g) $ for $t$ sufficiently large is well-defined by Proposition \ref{Prop seminorms are equivalent} and by the proof of Proposition \ref{Prop f^Omega well-defined}.

\begin{lem}\label{Lem Conv Of Constant Terms dominating function}
    Let $ f $ be a Schwartz function near $eP$ and let $ d > 1 $. Then, there exists a constant $ C > 0 $ such that
    $ |f(n' \kappa(n_t w) n_z a_q h)| \cdot \frac{\mu(t, s, z, q)}{t^k} $ is less or equal than
    \begin{equation}\label{Lem Conv Of Constant Terms dominating function eq1}
        C s^{n - \frac{k}2 - 1} a(\theta n_z)^{-\frac{n-k-1}2} (1 + \log a_{\theta n_z})^{-d}
        q^{- \frac{n-k+1}2} \big(1 - \frac{s}{q}\big)^{\frac{k}2-1} \big(1 + |\log(q)|\big)^{-d}
    \end{equation}
    for all $ n' \in N_\Gamma $ and $ h \in K $.
\end{lem}

\begin{rem}
    \eqref{Lem Conv Of Constant Terms dominating function eq1} is integrable over $ \RR^{n-k-1} \times \{ q \in \RR \mid q > s \} $ by Proposition \ref{Prop N bar integral is finite}.
\end{rem}

\begin{proof}
    Let $ \xi = \xi_{P, N^\Gamma} $. Let $ d > 1 $. By the proof of Lemma \ref{Lem Estimate for int_N |f(g n h)| dn} applied to $ E = N^\Gamma $ and as
    \[
        h(\kappa(n_t w) n_z a_q)^{(1-\xi) \rho_P} = \frac{(t^2+1)^{\frac{k}2}}{ (q + \frac1q \cdot \|(z_1 + t, z_2, \dotsc, z_{n-k-1})\|^2)^{\frac{k}2} } \ ,
    \]
    there exist $ C > 0 $ and $ T > 0 $ such that
    \begin{multline*}
        |f(n' \kappa(n_t w) n_z a_q h)| \leq C a(\theta n_z)^{-\xi \rho_P}
        \cdot \frac{(t^2+1)^{\frac{k}2}}{ (q + \frac1q \cdot \|(z_1 + t, z_2, \dotsc, z_{n-k-1})\|^2)^{\frac{k}2} } \\
        \cdot q^{\frac{n-k-1}2} (1 + \log a_{\theta n_z})^{-d} \big(1 + |\log(q)|\big)^{-d}
    \end{multline*}
    for all $ n' \in N_\Gamma $, $ z \in \RR^{n-k-1} $ and $ q > s $ with $ t \geq T $.
    Since moreover $ \sup_{t \geq T} \frac{ (t^2+1)^{\frac{k}2} }{2t^k} $ is finite, there exists a constant $ C' > 0 $ such that
    \begin{multline*}
        |f(n' \kappa(n_t w) n_z a_q h)| \cdot \frac{\mu(t, s, z, q)}{t^k} \\
        = |f(n' \kappa(n_t w) n_z a_q h)| \cdot \frac{s^{n - \frac{k}2 - 1}(q + \frac1q \cdot \|(z_1 + t, z_2, \dotsc, z_{n-k-1})\|^2)^{\frac{k}2}}{2 t^k q^{n-k} } \cdot \big(1 - \frac{s}{q}\big)^{\frac{k}2-1}
    \end{multline*}
    is less or equal than
    \[
        C' s^{n - \frac{k}2 - 1} a(\theta n_z)^{-\xi \rho_P} (1 + \log a_{\theta n_z})^{-d} q^{- \frac{n-k+1}2}
        \big(1 - \frac{s}{q}\big)^{\frac{k}2-1} \big(1 + |\log(q)|\big)^{-d} \ .
    \]
\end{proof}

\newpage

\begin{prop}\label{Prop formula Conv of constant terms}
    Let $ f $ be a left $ N_\Gamma M_\Gamma $-invariant Schwartz function near $eP$ and let $ h \in K $.
    If $ k \geq 2 $, then $ \lim_{t \to \infty} \frac1{t^k} \int_N f(\kappa(n_t w) n a_s h) \, dn $ is equal to
    \[
        \frac{\vol(S^{k-1})}2 \cdot
        \int_0^1 u^{n - \frac{k}2 - 2}(1 - u)^{\frac{k}2 - 1}
        \int_{M_k} \int_{N_\Gamma \bs N} f(n a_{\tfrac{s}{u}} m m_u m^{-1} h) \, dn \, dm \, du \ ,
    \]
    where $ m_u $ is defined as in Section \ref{sssec:diffeo Upsilon_(t, s)}.

    If $ k = 1 $, then $ \lim_{t \to \infty} \frac1{t^k} \int_N f(\kappa(n_t w) n a_s h) \, dn $ is equal to
    \begin{multline*}
        \frac12 \Big(\int_0^1 u^{n - \frac{k}2 - 2}(1 - u)^{\frac{k}2 - 1} \int_{N_\Gamma \bs N} f(n a_{\tfrac{s}{u}} m_u h) \, dn \, du \\
        + \int_0^1 u^{n - \frac{k}2 - 2}(1 - u)^{\frac{k}2 - 1} \int_{N_\Gamma \bs N} f(n a_{\tfrac{s}{u}} m_u^{-1} h) \, dn \, du \Big) \ .
    \end{multline*}
\end{prop}

\begin{rem}\label{Rem of Prop formula Conv of constant terms}
    Let $ \Sfrak' $ be an admissible generalised standard Siegel set with respect to $ P $ such that
    \[
        n' \kappa(n_t w) n a K \in \Gamma_P \Sfrak'
    \]
    for all $ n' \in N_\Gamma $, $ n \in N^\Gamma $, $ a \in A_{> s} $ and all $ t $ greater or equal than some $ T > 0 $.
    Let $ f \in \Cinf(\Gamma \bs G, \phi) $ be a Schwartz function near $eP$.
    Then, the proof of the proposition, combined with Proposition \ref{Prop a_(n_1 n_2 a) >= c a_(n_2 a)}, shows that there exists a nonnegative integrable function $ g_s $ on $ \RR^{n-k-1} \times (s, +\infty) $ such that
    \[
        \sup_{n' \in N_\Gamma, h \in K} |f(n' \kappa(n_t w) n_z a_q h)| \cdot \frac{\mu(t, s, z, q)}{t^k}
    \]
    is less or equal than
    \[
        \LO{p}{\Sfrak'}_{2d, 1, 1}(f) g_s(z, q)
    \]
    for all $ t \geq T $, $ z \in \RR^{n-k-1} $ and $ q \in (s, +\infty) $.
\end{rem}

\begin{proof}
    Let $f$ be as above and let $ h \in K $.
    Without loss of generality we may assume that $ h = e $. To simplify notations we assume further that $f$ is right $K$-invariant. The general case works analogously. The limit of the $K$-part has already been computed in Section \ref{sssec:diffeo Upsilon_(t, s)}. To get the integral over $ M_k $, one has to use Lemma \ref{Lem relation between Haar measures}.

    By \eqref{eq kappa(n_t w)},
    \[
        \frac1{t^k} \int_N f(\kappa(n_t w) n a_s) \, dn
        = \frac1{t^k} \int_{\RR^{n-1}} f(n_t w h(w n_{-t}) n_t n_v a_s) \, dv
    \]
    for all $ t > 0 $. As $f$ is left $ N_\Gamma $-invariant and right $K$-invariant, this is equal to
    \[
        \frac1{t^k} \int_{\RR^{n-k-1} \times \RR^k} f(\kappa(n_t w) n_{-t} (n_\cdot a_\cdot)(\Upsilon_{2, t, s}^{-1}(\Upsilon_{1, t, s}((v_1 + t, v_2, \dotsc, v_{n-1})))) \, dv \ .
    \]
    By passing to spherical coordinates, we get
    \begin{align}\label{eq integral with mu}
        & \frac{\vol(S^{k-1})}{t^k} \int_{\RR^{n-k-1}} \int_0^\infty f(\kappa(n_t w) (n_\cdot a_\cdot)(\Upsilon_{t, s}^{-1}(v, r))) r^{k-1} \, dr \, dv \\
        =& \vol(S^{k-1}) \int_s^{+\infty} \int_{\RR^{n-k-1}} f(\kappa(n_t w) n_z a_q) r(t, s, z, q)^{k-1} \cdot \frac{|\det(\d\Upsilon_{t, s})(z, q)|}{t^k} \, dz \, dq \nonumber \\
        =& \vol(S^{k-1}) \int_s^{+\infty} \int_{\RR^{n-k-1}} f(\kappa(n_t w) n_z a_q) \cdot \frac{\mu(t, s, z, q)}{t^k} \, dz \, dq \ . \nonumber
    \end{align}
    It follows from Lemma \ref{Lem Conv Of Constant Terms dominating function} that
    \[
        |f(\kappa(n_t w) n_z a_q)| \cdot \frac{\mu(t, s, z, q)}{t^k}
    \]
    is dominated by an integrable function on $ \RR^{n-k-1} \times \{ q \in \RR \mid q > s \} $.
    Since moreover
    \[
        \lim_{t \to \infty} \frac{\mu(t, s, z, q)}{t^k} = \frac1{2s} (\frac{s}{q})^{n - \frac{k}2}(1 - \frac{s}q)^{\frac{k}2 - 1}
    \]
    for all $ z \in \RR^{n-k-1} $ and $ q > s $ fixed, by \eqref{eq limit mu/t^k}, \eqref{eq integral with mu} converges to
    \[
        \frac{\vol(S^{k-1})}{2s} \int_s^\infty (\tfrac{s}{q})^{n - \frac{k}2}(1 - \tfrac{s}q)^{\frac{k}2 - 1} \int_{N_\Gamma \bs N} f(n a_q) \, dn \, dq
    \]
    when $t$ tends to $ +\infty $, by Lebesgue's theorem of dominated convergence and as $ \RR^{n-k-1} \equiv N^\Gamma \equiv N_\Gamma \bs N $ (as groups).

    Put $ u = \frac{s}{q} $. Then, $ q = \frac{s}{u} $ and $ dq = - \frac{s}{u^2} du $. Thus, the above is equal to
    \[
        \frac{\vol(S^{k-1})}2 \cdot \int_0^1 u^{n - \frac{k}2 - 2}(1 - u)^{\frac{k}2 - 1} \int_{N_\Gamma \bs N} f(n a_{\tfrac{s}{u}}) \, dn \, du \ .
    \]
    The proposition follows.
\end{proof}

\newpage

\subsubsection{The formula for derivatives of Schwartz functions in the \texorpdfstring{$ \n_\Gamma $}{n\_Gamma}-direction}

The following, including the next lemma, is a preparation for Proposition \ref{Prop formula Conv of constant terms for derivative}, which says that the constant term along $ \Omega $ divided by $ t^k $ converges to zero for derivatives of Schwartz functions in the $ \n_\Gamma $-direction when we approach a cusp of smaller rank.

\medskip

Let $ t \geq 1 $ and let
\[
    R(t, s, z, q) = \frac{t^2+1}{s (q + \frac1q \|(z_1 + t, z_2, \dotsc, z_{n-k-1})\|^2)} \cdot r(t, s, z, q) \ .
\]
By Section \ref{sssec:diffeo Upsilon_(t, s)}, $R$ is the Euclidean norm of
\[
    \log(\nu_{N_\Gamma}( \kappa(n_t w) n_{(v_1, \dotsc, v_{n-k-1}, r \Phi_k(\psi_0, \psi_1, \dotsc, \psi_{k-2}))} )) \in \n_\Gamma \equiv \RR^k
\]
for $ (v_1, \dotsc, v_{n-k-1}, r) = \Upsilon_{t, s}(z, q) $.
Let $ z(t) = (z_1+t, z_2, \dotsc, z_{n-k-1}) $.
\\ Since $ r(t, s, z, q) = \sqrt{s} \cdot \sqrt{1 - \frac{s}{q}} \cdot \sqrt{ q + \frac1q \|z(t)\|^2 } $, $ R(t, s, z, q) $ is equal to
\[
     \frac{(t^2+1) \cdot \sqrt{1 - \frac{s}{q}}}{\sqrt{s} \sqrt{q + \frac1q \|z(t)\|^2}} \ .
\]
Direct computation shows that $ \frac{ \del R }{ \del q }(t, s, z, q) $ is equal to
\[
    \frac{ t^2+1 }{ 2\sqrt{s} } \cdot
    \frac{ q + \frac1q \|z(t)\|^2 - 2q(1 - \frac{s}{q}) }{ q \sqrt{1 - \frac{s}{q}} (q + \frac1q \|z(t)\|^2)^{\frac32} } \ .
\]
If $ q \in (s, \frac12 t] $ and $ \|z(t)\| \in [t, 2t] $, then this is greater or equal than
\[
    \frac{ t^2+1 }{ 2\sqrt{s} } \cdot \frac{ \|z(t)\|^2 - q^2  }{ \sqrt{q} \sqrt{1 - \frac{s}{q}} (q^2 + \|z(t)\|^2)^{\frac32} }
    \geq \frac{ t }{ 2\sqrt{s} } \cdot \frac{ \frac34 }{ \sqrt{q} \sqrt{1 - \frac{s}{q}} (\frac14 + 4)^{\frac32} }  \ .
\]
Let $ c = \big( \frac{ 3 }{ 8 (\frac14 + 4)^{\frac32} } \big)^{-1} > 0 $. Then, this is again greater or equal than
\[
    \frac{ t }{ c \sqrt{s} \sqrt{q} \sqrt{1 - \frac{s}{q}} } > 0 \ .
\]
Since
\[
    \mu(t, s, z, q) = \frac{s^{n - \frac{k}2 - 1}(q + \frac1q \|z(t)\|^2)^{\frac{k}2}}{2 q^{n-k} } \cdot \big(1 - \frac{s}{q}\big)^{\frac{k}2-1} \ ,
\]
$ \mu_1(t, s, z, q) := \Big( \frac{ \del R }{ \del q }(t, s, z, q) \Big)^{-1} \mu(t, s, z, q) $ is equal to
\[
    \frac{ s^{n - \frac{k+1}2} }{ t^2+1 } \cdot
    \frac{ 1 }{ q + \frac1q \|z(t)\|^2 - 2q(1 - \frac{s}{q}) }
    \cdot q^{-n+k+1} (q + \frac1q \|z(t)\|^2)^{\frac{k+3}2} \cdot \big(1 - \frac{s}{q}\big)^{\frac{k-1}2} \ .
\]
Hence, $ \frac{ \del }{ \del q } \mu_1(t, s, z, q) $ is equal to
\begin{align*}
    & \Big( \frac{ \del R }{ \del q }(t, s, z, q) \Big)^{-1} \mu(t, s, z, q) \Big(
    \frac{ 1 + \frac1{q^2} \|z(t)\|^2 }{ q + \frac1q \|z(t)\|^2 - 2q(1 - \frac{s}{q}) }
    + \frac{-n+k+1}{q} \\
    & \qquad + \frac{ k + 3 }2 \cdot \frac{ 1 - \frac1{q^2} \|z(t)\|^2 }{ q + \frac1q \|z(t)\|^2 }
    + \frac{ k - 1 }2 \cdot \frac{ s }{ q^2 (1 - \frac{s}{q}) }
    \Big) \ .
\end{align*}
Set $ \mu_2(t, s, z, q) = \frac1{ R(t, s, z, q) } \cdot \mu(t, s, z, q) $.

\medskip

\needspace{5\baselineskip}
The following lemma allows us to compare the measure $ \mu(t, s, z, q) \, dz \, dq $ with other measures that will arise in the proof of Proposition \ref{Prop formula Conv of constant terms for derivative}.

\begin{lem}\label{Lem estimates of mu'/mu}
    Let $ k \geq 1 $ and let $ t \geq 1 $. Let $ c > 0 $ be as above, $ q \in (s, \frac12 t] $ and $ z \in \RR^{n-k-1} $ be such that $ \|z(t)\| \in [t, 2t] $. Then: \nlenum
    \begin{enumerate}
    \item\label{Lem estimates of mu'/mu eq1}
        $ \Big( \frac{ \del R }{ \del q }(t, s, z, q) \Big)^{-1} $ is less or equal than $ \frac{c}t \sqrt{s} \sqrt{q} \sqrt{1 - \frac{s}{q}} $.
    \item\label{Lem estimates of mu'/mu eq2}
        $ |\frac{ \frac{ \del }{ \del q } \mu_1(t, s, z, q) }{ \mu(t, s, z, q) }| $ is less or equal than
        \[
            \frac{c}t \sqrt{s} \sqrt{q} \sqrt{1 - \frac{s}{q}} \Big( \frac{n - \frac{k}2 + \frac{13}2 }{q} + \frac{ k - 1 }2 \cdot \frac{ s }{ q^2 (1 - \frac{s}{q}) } \Big) \ .
        \]
    \item\label{Lem estimates of mu'/mu eq3}
         $ \frac1{ R(t, s, z, q) } $ is less or equal than $ \frac{ 3\sqrt{s} }{ t \sqrt{q} \sqrt{1 - \frac{s}{q}} } $.
    \end{enumerate}
\end{lem}

\begin{proof} Let $ t \geq 1 $, $ q $, $ z $ be as above. \nlenum
    \begin{enumerate}
    \item By the above, $ \Big( \frac{ \del R }{ \del q }(t, s, z, q) \Big)^{-1} \leq \frac{c}t \sqrt{s} \sqrt{q} \sqrt{1 - \frac{s}{q}} \ $.
        The first assertion of the lemma follows.
    \item As $ \big( \frac{ \del R }{ \del q }(t, s, z, q) \big)^{-1} \leq \frac{c}t \sqrt{s} \sqrt{q} \sqrt{1 - \frac{s}{q}} $,
        $ |\frac{ 1 + \frac1{q^2} \|z(t)\|^2 }{ q + \frac1q \|z(t)\|^2 - 2q(1 - \frac{s}{q}) }| \leq \frac{ \frac14 + 4 }{ 1 - \frac14 }  \cdot \frac1q
        \leq \frac6q $ and as
        $ |\frac{ 1 - \frac1{q^2} \|z(t)\|^2 }{ q + \frac1q \|z(t)\|^2 }| \leq \frac1q $,
        \[
            |\frac{ \frac{ \del }{ \del q } \mu_1(t, s, z, q) }{ \mu(t, s, z, q) }|
        \]
        is less or equal than
        \[
            \frac{c}t \sqrt{s} \sqrt{q} \sqrt{1 - \frac{s}{q}} \Big( \frac{n - \frac{k}2 + \frac{13}2 }{q} + \frac{ k - 1 }2 \cdot \frac{ s }{ q^2 (1 - \frac{s}{q}) } \Big) \ .
        \]
        This shows the second assertion of the lemma.
    \item Note that $ \frac1{ R(t, s, z, q)} $ is less or equal than $ \frac{\sqrt{ 4 + \frac14 } }{ t } \cdot \frac{\sqrt{s}}{ \sqrt{q} \sqrt{1 - \frac{s}{q}} } $.
        So, the third assertion of the lemma holds, too.
    \end{enumerate}
\end{proof}

Let $ D = \{ (\psi_0, \psi_1, \dotsc, \psi_{k-2}) \in \RR^{k-1} \mid \psi_0 \in (0, 2\pi) \text{ and } \psi_l \in (-\tfrac{\pi}2, \tfrac{\pi}2) \text{ for all } l \geq 1 \} $ and let
\[
    E_t = \{ z \in \RR^{n-k-1} \mid \|z(t)\| \in [t, 2t] \} \ .
\]

\begin{prop}\label{Prop formula Conv of constant terms for derivative}
    Let $ f $ be a Schwartz function near $eP$, $ h \in K $ and let $ X \in \n_\Gamma $. Then,
    \[
        \frac1{t^k} \int_N L_X f(\kappa(n_t w) n a_s h) \, dn
    \]
    converges to zero when $t$ tends to $ \infty $.
\end{prop}

\begin{proof}
    Let $ f $ be as above and let $ X \in \n_\Gamma $. Without loss of generality, we may assume that $ h = e $ and that $ X = \log(n_{e_{n-k-1+j}}) $ for some $ j \in \{1, \dotsc, k\} $. Assume that $ k \geq 2 $.

    Let $ Y_0, Y_1, \dotsc, Y_{k-2} \in \so(k) $ be as in the proof of Lemma \ref{Lem Relation between spherical coordinates and so(d)} (with $ d = k $).
    \\ Set $ A(\psi_0, \dotsc, \psi_{k-2}) = \prod_{j=0}^{k-2} \exp(-\psi_j Y_j) \in \SO(k) $. For $ A \in \SO(k) $, set $ m_A^{M_k} = m_{\diag(I_{n-k-1}, A)} $.

    It follows from Remark \ref{Rem of Prop formula Conv of constant terms}, Section \ref{sssec:diffeo Upsilon_(t, s)}, Proposition \ref{Prop convergence of Sfrak p_(r, X, Y)(R_h f - f)} and the fact that the indicator functions $ 1_{ q \not \in (s, \frac12 t] } $, $ 1_{ \|z(t)\| \not \in [t, 2t] } $ converge pointwise to zero when $t$ tends to $ \infty $ that it suffices to show that
    \begin{align*}
        & \int_D \int_s^{\frac12 t} \int_{E_t}  \\
        & \quad L_X f(n^{N_\Gamma}_{(R(t, s, z, q) \Phi_k(\psi_0, \psi_1, \dotsc, \psi_{k-2}))} \kappa(n_t w) n_z a_q m^{M_k}_{A(\psi_0, \dotsc, \psi_{k-2})} m_{\frac{s}{q}} (m^{M_k}_{A(\psi_0, \dotsc, \psi_{k-2})})^{-1}) \\
        & \quad \cdot \frac{\mu(t, s, z, q)}{t^k} \, dz \, dq \ \cos \psi_1 (\cos \psi_2)^2 \cdots (\cos \psi_{k-2})^{k-2} \, d(\psi_0, \psi_1, \dotsc, \psi_{k-2})
    \end{align*}
    converges to zero when $t$ tends to $ \infty $.
    Let
    \[
        Z_t = \begin{pmatrix}
            0 & 0 & 0 & 0 & 0 & 0_{1, k-1} \\
            0 & 0 & 0 & 0 & 0 & 0_{1, k-1} \\
            0 & 0 & 0 & 0_{1, n-k-2} & t & 0_{1, k-1} \\
            0_{n-k-2, 1} & 0_{n-k-2, 1} & 0_{n-k-2, 1} & I_{n-k-2} & 0_{n-k-2, 1} & 0_{1, k-1} \\
            0 & 0 & -t & 0_{1, n-k-2} & 0 & 0_{1, k-1} \\
            0_{k-1, 1} & 0_{k-1, 1} & 0_{k-1, 1} & 0_{k-1, 1} & 0_{k-1, 1} & 0_{k-1, k-1}
        \end{pmatrix} \in \m \ .
    \]
    Then, $ \exp(Z_{\arccos(2u-1)}) = m_u $.
    Let
    \begin{multline*}
        \tilde{f}_t(R, \psi_0, \psi_1, \dotsc, \psi_{k-2}, z, q, \tilde{q}, x, y) \\
        = f\big(n^{N_\Gamma}_{R \Phi_k(\psi_0, \psi_1, \dotsc, \psi_{k-2})} \kappa(n_t w) n_z a_q m^{M_k}_{A(x_0, \dotsc, x_{k-2})} \exp(Z_{\tilde{q}}) (m^{M_k}_{A(y_0, \dotsc, y_{k-2})})^{-1} \big) \ .
    \end{multline*}
    Define $ \varphi(t, s, \psi_0, \dotsc, \psi_{k-2}, z, q) $ by
    \[
        \big(R(t, s, z, q), z, q, m^{M_k}_{A(\psi_0, \dotsc, \psi_{k-2})} \exp(Z_{\arccos(2 \cdot \frac{s}{q}-1)}) (m^{M_k}_{A(\psi_0, \dotsc, \psi_{k-2})})^{-1}\big) \ .
    \]
    Then,
    \[
        L_X f(n^{N_\Gamma}_{(R(t, s, z, q) \Phi_k(\psi_0, \psi_1, \dotsc, \psi_{k-2}))} \kappa(n_t w) n_z a_q m^{M_k}_{A(\psi_0, \dotsc, \psi_{k-2})} m_{\frac{s}{q}} (m^{M_k}_{A(\psi_0, \dotsc, \psi_{k-2})})^{-1})
    \]
    is equal to
    \[
        \del_{n-k-1+j} \tilde{f}_t(\varphi(t, s, \psi_0, \dotsc, \psi_{k-2}, z, q)) \ .
    \]
    Let $ p_{i, j}^{(k)} \colon (0, 2\pi) \times (-\tfrac{\pi}2, \tfrac{\pi}2)^{k-2} \to \RR $ be as in Lemma \ref{Lem cartesian partial derivatives in terms spherical partial derivatives} when we apply it to $ \{0_{1, n-k-1}\} \times \RR^k \equiv \RR^k $. Then,
    \begin{equation}\label{eq del_(n-k-1+j) w.r.t. spherical partial derivatives}
        \del_{n-k-1+j} = \Phi_k(\psi_0, \dotsc, \psi_{k-2})_j \del_{R} + \frac1{R} \sum_{i=0}^{k-2} p_{i, j}^{(k)}(\psi_0, \dotsc, \psi_{k-2}) \del_{\psi_i} \ .
    \end{equation}
    Let $ H_t = \log(a_t) $. Since $ m a = a m $ for all $ a \in A $ and $ m \in M $,
    \[
        \big(\frac{ \del \tilde{f}_t }{ \del q } \circ \varphi\big)(t, s, \psi_0, \dotsc, \psi_{k-2}, z, q)
    \]
    is equal to
    \begin{multline}\label{eq del_q f}
        R_{H_1} f(n^{N_\Gamma}_{(R(t, s, z, q) \Phi_k(\psi_0, \psi_1, \dotsc, \psi_{k-2})) } \\
        \kappa(n_t w) n_z a_q m^{M_k}_{A(\psi_0, \dotsc, \psi_{k-2})} m_{\frac{s}{q}} (m^{M_k}_{A(\psi_0, \dotsc, \psi_{k-2})})^{-1})
    \end{multline}
    and $ \big(\frac{ \del \tilde{f}_t }{ \del \tilde{q} } \circ \varphi\big)(t, s, \psi_0, \dotsc, \psi_{k-2}, z, q) $ is equal to
    \begin{multline}\label{eq del_tilde q f}
        R_{m^{M_k}_{A(\psi_0, \dotsc, \psi_{k-2})} Z_1 (m^{M_k}_{A(\psi_0, \dotsc, \psi_{k-2})})^{-1}} f(n^{N_\Gamma}_{(R(t, s, z, q) \Phi_k(\psi_0, \psi_1, \dotsc, \psi_{k-2})) } \\
        \kappa(n_t w) n_z a_q m^{M_k}_{A(\psi_0, \dotsc, \psi_{k-2})} m_{\frac{s}{q}} (m^{M_k}_{A(\psi_0, \dotsc, \psi_{k-2})})^{-1}) \ .
    \end{multline}
    When one estimates the absolute value of these functions, one has to use that $ \Ad(x^{-1})Z_1 $ depends smoothly on $x$ and that $ M $ is compact.
    
    By the chain rule and as $ \tfrac{ \del }{ \del q } \arccos(\tfrac{2s}{q} - 1) = \frac{\sqrt{s}}{ q^{\frac32} \sqrt{1- \tfrac{s}{q}} } $, we have
    \[
        \frac{\del (\tilde{f}_t \circ \varphi)}{\del q} = \Big(\frac{ \del \tilde{f}_t }{ \del R } \circ \varphi\Big) \frac{ \del R }{ \del q }
        + \Big(\frac{ \del \tilde{f}_t }{ \del \tilde{q} } \circ \varphi\Big) \cdot \frac{\sqrt{s}}{ q^{\frac32} \sqrt{1- \tfrac{s}{q}} } + \frac{ \del \tilde{f}_t }{ \del q } \circ \varphi \ .
    \]
    Thus,
    \begin{equation}\label{eq chain rule for del_q (f circ varphi)}
        \frac{ \del \tilde{f}_t }{ \del R } \circ \varphi = \Big(\frac{\del (\tilde{f}_t \circ \varphi)}{\del q} - \Big(\frac{ \del \tilde{f}_t }{ \del \tilde{q} } \circ \varphi\Big) \cdot \frac{\sqrt{s}}{ q^{\frac32} \sqrt{1- \tfrac{s}{q}} } - \frac{ \del \tilde{f}_t }{ \del q } \circ \varphi\Big) \Big(\frac{ \del R }{ \del q }\Big)^{-1} \ .
    \end{equation}
    Hence, it follows from Lemma \ref{Lem Conv Of Constant Terms dominating function} and Lemma \ref{Lem estimates of mu'/mu}, partial integration and Fubini that
    \begin{align*}
        & \int_D \int_s^{\frac12 t} \int_{E_t} \frac{\del (\tilde{f}_t \circ \varphi)}{\del q}(t, s, \psi_0, \dotsc, \psi_{k-2}, z, q)
        \cdot \frac{\mu_1(t, s, z, q)}{t^k} \, dz \, dq \\
        & \quad \Phi_k(\psi_0, \dotsc, \psi_{k-2})_j \cos \psi_1 (\cos \psi_2)^2 \cdots (\cos \psi_{k-2})^{k-2} \, d(\psi_0, \psi_1, \dotsc, \psi_{k-2})
    \end{align*}
    converges to zero when $t$ tends to $ \infty $. Because of Lemma \ref{Lem Conv Of Constant Terms dominating function} and Lemma \ref{Lem estimates of mu'/mu}, the boundary term appearing from partial integration converges to zero.

    It follows from this, \eqref{eq del_q f}, \eqref{eq del_tilde q f}, \eqref{eq chain rule for del_q (f circ varphi)}, Lemma \ref{Lem Conv Of Constant Terms dominating function} and Lemma \ref{Lem estimates of mu'/mu} that
    \begin{align*}
        & \int_D \int_s^{\frac12 t} \int_{E_t} \big( \frac{ \del \tilde{f}_t }{ \del R } \circ \varphi \big)(t, s, \psi_0, \dotsc, \psi_{k-2}, z, q)
        \cdot \frac{\mu(t, s, z, q)}{t^k} \, dz \, dq \\
        & \quad \Phi_k(\psi_0, \dotsc, \psi_{k-2})_j \cos \psi_1 (\cos \psi_2)^2 \cdots (\cos \psi_{k-2})^{k-2} \, d(\psi_0, \psi_1, \dotsc, \psi_{k-2})
    \end{align*}
    converges also to zero when $t$ tends to $ \infty $.

    By the chain rule, we have
    \[
        \frac{ \del (\tilde{f}_t \circ \varphi) }{ \del \psi_i } = \frac{ \del \tilde{f}_t }{ \del \psi_i } \circ \varphi
        + \sum_{j=0}^{k-2} \Big( \frac{ \del \tilde{f}_t }{ \del x_j } \circ \varphi \Big)
        + \sum_{j=0}^{k-2} \Big( \frac{ \del \tilde{f}_t }{ \del y_j } \circ \varphi \Big) \ .
    \]
    Thus,
    \[
        \frac{ \del \tilde{f}_t }{ \del \psi_i } \circ \varphi
        = \frac{ \del (\tilde{f}_t \circ \varphi) }{ \del \psi_i } - \sum_{j=0}^{k-2} \Big( \frac{ \del \tilde{f}_t }{ \del x_j } \circ \varphi \Big)
        - \sum_{j=0}^{k-2} \Big( \frac{ \del \tilde{f}_t }{ \del y_j } \circ \varphi \Big) \ .
    \]
    Note that $ (\frac{ \del \tilde{f}_t }{ \del y_i } \circ \varphi)(\psi_0, \psi_1, \dotsc, \psi_{k-2}, z, q) $ is equal to
    \begin{multline*}
        R_{-\Ad\big( (\prod_{j=i+1}^{k-2} m^{M_k}_{\exp(-\psi_j Y_j)})^{-1} \big)Y_i} f(n^{N_\Gamma}_{(R(t, s, z, q) \Phi_k(\psi_0, \psi_1, \dotsc, \psi_{k-2}))} \\
        \kappa(n_t w) n_z a_q m^{M_k}_{A(\psi_0, \dotsc, \psi_{k-2})} m_{\frac{s}{q}} (m^{M_k}_{A(\psi_0, \dotsc, \psi_{k-2})})^{-1}) \ .
    \end{multline*}
    and that $ (\frac{ \del \tilde{f}_t }{ \del x_i } \circ \varphi)(\psi_0, \psi_1, \dotsc, \psi_{k-2}, z, q) $ is equal to
    \begin{multline*}
        R_{-\Ad\big( (\prod_{j=i+1}^{k-2} m^{M_k}_{\exp(-\psi_j Y_j)} m_{\frac{s}{q}} \prod_{j=1}^{k-2} m^{M_k}_{\exp(-\psi_j Y_j)})^{-1} \big)Y_i} f(n^{N_\Gamma}_{(R(t, s, z, q) \Phi_k(\psi_0, \psi_1, \dotsc, \psi_{k-2}))} \\
        \kappa(n_t w) n_z a_q m^{M_k}_{A(\psi_0, \dotsc, \psi_{k-2})} m_{\frac{s}{q}} (m^{M_k}_{A(\psi_0, \dotsc, \psi_{k-2})})^{-1}) \ .
    \end{multline*}
    When one estimates the absolute value of these functions, one has to use as above that $ \Ad(x^{-1})Y_i $ depends smoothly on $x$ and that $ M $ is compact.
    Similar as above, one shows that
    \begin{align*}
        & \int_D \int_s^{\frac12 t} \int_{E_t} \big( \frac{ \del \tilde{f}_t }{ \del \psi_i } \circ \varphi \big)(t, s, \psi_0, \dotsc, \psi_{k-2}, z, q)
        \cdot \frac{\mu(t, s, z, q)}{t^k} \, dz \, dq \\
        & \quad p_{i, j}^{(k)}(\psi_0, \psi_1, \dotsc, \psi_{k-2}) \cos \psi_1 (\cos \psi_2)^2 \cdots (\cos \psi_{k-2})^{k-2} \, d(\psi_0, \psi_1, \dotsc, \psi_{k-2})
    \end{align*}
    converges to zero when $t$ tends to $ \infty $. To show this, one needs of course also Lemma \ref{Lem cartesian partial derivatives in terms spherical partial derivatives}. Let us check now that indeed no boundary term appears from partial integration.

    Note first that $ \d( R\Phi_k(\psi_0, \psi_1, \dotsc, \psi_{k-2}) )^t $ is given by
    \[
        \blockmat{\cos \psi_{k-2} \ \d( R\Phi_{k-1}(\psi_0, \psi_1, \dotsc, \psi_{k-3}) )^t}{\begin{matrix} \sin \psi_{k-2} \\ 0_{k-2, 1}\end{matrix}}{-\sin \psi_{k-2} \ R \Phi_{k-1}(\psi_0, \psi_1, \dotsc, \psi_{k-3})}{R \cos \psi_{k-2}} \ .
    \]
    Thus, if $ \psi_0 \in \{0, 2\pi\} $ and $ \psi_l \in \{-\frac{\pi}2, \frac{\pi}2\} $ for all $ l \geq 1 $, then $ \d( R\Phi_k(\psi_0, \psi_1, \dotsc, \psi_{k-2}) )^t $ is equal to
    \[
        \begin{pmatrix} 1 & 0 \\ 0 & R \end{pmatrix} \quad \text{if $ k = 2 $}
        \quad \text{and} \quad
        \blockmat{0_{k-1,k-1}}{ (\pm 1, 0_{1, k-2})^t }{ (0_{1, k-2}, -R)}{0_{1, k-1}} \quad \text{if $ k > 2 $}
        \ .
    \]
    It follows now from the proof of Lemma \ref{Lem cartesian partial derivatives in terms spherical partial derivatives} that
    \[
        \cos \psi_1 (\cos \psi_2)^2 \cdots (\cos \psi_{k-2})^{k-2} p_{i, j}^{(k)}(\psi_0, \psi_1, \dotsc, \psi_{k-2})
    \]
    is equal to zero for all $ i \in \{0, \dotsc, k-2\} $, $ j \in \{1, \dotsc, k\}) $ if $ \psi_0 \in \{0, 2\pi\} $ and $ \psi_l \in \{-\frac{\pi}2, \frac{\pi}2\} $ for all $ l \geq 1 $ such that $ k \neq 2 $ or $ j \neq 2 $.
    Moreover, $ p_{0, 2}^{(2)}(\psi_0) $ is equal to 1 if $ \psi_0 \in \{0, 2\pi\} $. Thus, $ p_{0, 2}^{(2)}(2\pi) - p_{0, 2}^{(2)}(0) = 0 $.
    Hence, no boundary term appears.

    We leave the case $ k = 1 $, which can be proven similarly, to the reader. 
    The proposition follows from \eqref{eq del_(n-k-1+j) w.r.t. spherical partial derivatives}.
\end{proof}

\newpage

\subsubsection{The general formula}

We use here the notations seen in Section \ref{ssec:Conv of the constant terms: Fourier expansion of f}.

\begin{lem}\label{Lem f = Delta F}
    If $ f $ is a Schwartz function near $eP$ with $ f^{P, lc} = 0 $, then there exists a Schwartz function $ F $ near $eP$ such that
    \[
        f = \Delta F \ ,
    \]
    where $ \Delta := \sum_{j=1}^k L_{Z_j}^2 $.
\end{lem}

\begin{proof}
    Let $ f $ be as above. Let $ g \in N^\Gamma A K $.

    As in the proof of Theorem \ref{Thm f - f^(Q, lc) is rapidly decreasing}, we may assume without loss of generality that $ \dim_\CC V_\phi = 1 $.
    Using the notations seen in Section \ref{ssec:Conv of the constant terms: Fourier expansion of f}, we have
    \[
        f(n^{N_\Gamma, \{Z_j\}}_z g) = \sum_{-\lambda \neq x \in \ZZ^k} c_{x, g}(f) e^{ 2 \pi i \langle z, \lambda + x \rangle }  \qquad (z \in \RR^k) \ .
    \]
    Set
    \[
        F(n^{N_\Gamma, \{Z_j\}}_z g) = - \sum_{-\lambda \neq x \in \ZZ^k} \frac{ c_{x, g}(f) }{ 4 \pi^2 \|\lambda + x\|_2^2 } e^{2 \pi i \langle z, \lambda + x \rangle} \qquad (z \in \RR^k, g \in N^\Gamma A K) \ .
    \]
    Then,
    \[
        (\Delta F)(n^{N_\Gamma, \{Z_j\}}_z g) = \sum_{-\lambda \neq x \in \ZZ^k} c_{x, g}(f) e^{2 \pi i \langle z, \lambda + x \rangle} = f(n^{N_\Gamma, \{Z_j\}}_z g)
    \]
    since $ \frac{ \del^2 }{ \del z_l^2 } e^{2 \pi i \langle z, \lambda + x \rangle} = - 4 \pi^2 (\lambda_l + x_l)^2 e^{2 \pi i \langle z, \lambda + x \rangle} $.

    Let $ \Sfrak $ be an admissible generalised Siegel set relative to $P$. Let $ r \geq 0 $, $ X \in \U(\n) $ and $ Y \in \U(\g) $. Let $ g \in \Sfrak $. Then,
    \[
        F(g) = -\sum_{-\lambda \neq x \in \ZZ^k} \frac{ c_{x, g}(f) }{ 4 \pi^2 \|\lambda + x\|_2^2 } \ .
    \]
    Since $f$ is a Schwartz function near $eP$, there exists $ c' > 0 $ such that
    \begin{multline*} 
        |\int_{[0, 1]^k} ((L_X R_Y f)(n^{N_\Gamma, \{Z_j\}}_u \varphi(n^{N_\Gamma, \{Z_j\}}_u) \ \cdot))(g) e^{ -2 \pi i \langle u, \lambda + x \rangle } \, du| \\
        \leq c' h(g)^{(1-\xi)\rho_P} a_g^{-\xi \rho_P} (1 + \log a_g)^{-r}    \qquad (x \in \ZZ^k, g \in \Sfrak) \ ,
    \end{multline*}
    where $ \xi := \xi_{P, N^\Gamma} $. Since $ \sum_{-\lambda \neq x \in \ZZ^k} \frac{ 1 }{ 4 \pi^2 \|\lambda + x\|_2^2 }$ converges by Lemma \ref{Lem sum_(0 neq x in Z^k) 1/( |x|_2^(N+1) ) }, there is a constant $ c'' > 0 $ such that $ |L_X R_Y F(g)| $ is less or equal than
    \[
        c'' h(g)^{(1-\xi)\rho_P} a_g^{-\xi \rho_P} (1 + \log a_g)^{-r}    \qquad (g \in \Sfrak) \ .
    \]
    Hence, $ \LO{p}{\Sfrak'}_{r, X, Y}(F) $ is finite. The lemma follows.
\end{proof}

\begin{prop}\label{Prop formula Conv of constant terms if f^(P, lc) is zero}
    Let $ f \in \Cinf(\Gamma \bs G, \phi) $ be a Schwartz function near $eP$ with $ f^{P, lc} = 0 $ and let $ h \in K $. Then, $ \lim_{t \to \infty} \frac1{t^k} \int_N f(\kappa(n_t w) n a_s h) \, dn $ is zero.
\end{prop}

\begin{proof}
    Let $ h \in K $, $ f $ be as above and let $F$ be as in Lemma \ref{Lem f = Delta F}.
    Then,
    \[
        \frac1{t^k} \int_N f(\kappa(n_t w) n a_s) \, dn = \sum_{j=1}^k \frac1{t^k} \int_N \big(\sum_{j=1}^k L_{Z_j} L_{Z_j} F\big)(\kappa(n_t w) n a_s) \, dn \ .
    \]
    This converges to zero by Proposition \ref{Prop formula Conv of constant terms for derivative} applied to $ L_{Z_j} F $ (Schwartz function near $eP$).
\end{proof}

\needspace{4\baselineskip}
\begin{thm}\label{Thm Conv of the constant terms}
    Let $ f \in \Cinf(\Gamma \bs G, \phi) $ be a Schwartz function near $eP$ and let $ h \in K $. If $ k \geq 2 $, then
    \[
        \lim_{t \to \infty} \frac1{t^k} f^\Omega(\kappa(n_t w), a_s h)   \qquad (s \in (0, +\infty))
    \]
    is up to a constant equal to
    \[
        \int_0^1 u^{n - \frac{k}2 - 2}(1 - u)^{\frac{k}2 - 1}
        \int_{M_k} f^P(a_{\tfrac{s}{u}} m' m_u {m'}^{-1} h) \, dm' \, du \ ,
    \]
    where
    \begin{align*}
        m_u &:= \lim_{t \to \infty} k(w n_{t (\sqrt{u}, 0_{1, n-k-2}, \sqrt{ 1 - u }, 0_{1, k-1})}) \\
        &= w n_{(\sqrt{u}, 0_{1, n-k-2}, \sqrt{ 1 - u }, 0_{k-1})} w n_{(\sqrt{u}, 0_{1, n-k-2}, -\sqrt{ 1 - u }, 0_{1, k-1})} w n_{(\sqrt{u}, 0_{1, n-k-2}, \sqrt{ 1 - u }, 0_{1, k-1})} \\
        &= \begin{pmatrix}
            1 & 0 & 0 & 0 & 0 & 0_{1, k-1} \\
            0 & 1 & 0 & 0 & 0 & 0_{1, k-1} \\
            0 & 0 & 2u - 1 & 0_{1, n-k-2} & 2\sqrt{u(1-u)} & 0_{1, k-1} \\
            0_{n-k-2, 1} & 0_{n-k-2, 1} & 0_{n-k-2, 1} & I_{n-k-2} & 0_{n-k-2, 1} & 0_{1, k-1} \\
            0 & 0 & -2\sqrt{u(1-u)} & 0_{1, n-k-2} & 2u - 1 & 0_{1, k-1} \\
            0_{k-1, 1} & 0_{k-1, 1} & 0_{k-1, 1} & 0_{k-1, 1} & 0_{k-1, 1} & I_{k-1, k-1}
        \end{pmatrix} \ .
    \end{align*}
    If $ k = 1 $, then
    \[
        \lim_{t \to \infty} \frac1{t^k} f^\Omega(\kappa(n_t w), a_s h)     \qquad (s \in (0, +\infty))
    \]
    is up to a constant equal to
    \[
        \int_0^1 u^{n - \frac{k}2 - 2}(1 - u)^{\frac{k}2 - 1}
        \big(f^P(a_{\tfrac{s}{u}} m_u h) + f^P(a_{\tfrac{s}{u}} m_u^{-1} h)\big) \, du \ .
    \]
\end{thm}

\begin{proof}
    Let $ f $ be as above. Then, $ f^{P, lc} $ is a $ N_\Gamma M_\Gamma $-invariant Schwartz function near $eP$ and $ (f - f^{P, lc})^{P, lc} = 0 $ by Proposition \ref{Prop Little constant term properties} and Lemma \ref{Lem estimate of the little constant term}.

    The theorem follows now directly from Proposition \ref{Prop formula Conv of constant terms} applied to $ f^{P, lc} $ and Proposition \ref{Prop formula Conv of constant terms if f^(P, lc) is zero} applied to $ f - f^{P, lc} $.
\end{proof}

\newpage

\subsubsection{The main theorem}

Let $ (\gamma, V_\gamma) \in \hat{K} $. Let $ I $ be the canonical isomorphism from $ \SO(2) $ to $ S := \{ m_u \mid u \in [0, 1] \} $.

Since all irreducible, unitary representations of $ \SO(2) $ are one-dimensional and since $ V_\gamma $ is finite-dimensional, there exist $ v_j \in V_\gamma $ such that
\[
    V_\gamma = \bigoplus_{j=1}^l \CC v_j
\]
as representation spaces of $S$. We identify $ V_\gamma $ with $ \CC^l $ via this basis. 

For $ \phi \in \RR $, set $ k_\phi = \begin{pmatrix} \cos(\phi) & \sin(\phi) \\ -\sin(\phi) & \cos(\phi) \end{pmatrix} \in \SO(2) $. Let $ \pi $ be an irreducible representation of $ \SO(2) $. Then, there exists $ m \in \ZZ $ such that $ \pi(k_\phi) = e^{i m \phi} $ for all $ \phi \in \RR $.
Thus, there are $ m_j \in \ZZ $ such that
\[
    \gamma(I(k_\phi)) = D_{m_1, \dotsc, m_l}(\phi) \ ,
\]
where $ D_{m_1, \dotsc, m_l}(\phi) := \begin{pmatrix} e^{i m_1 \phi} & 0 & \cdots & 0 \\ 0 & e^{i m_2 \phi} & \ddots & \vdots \\ \vdots & \ddots & \ddots & 0 \\ 0 & \cdots & 0 & e^{i m_l \phi} \end{pmatrix} $.

Let $ Y \in \so(k) $ be such that $ m_x = m_{\diag(I_{n-k-1}, \exp(\arccos(2x-1) Y))} $.

Then, $ \gamma(m_x)^{\pm 1} = \gamma(I(k_{\pm \arccos(2x-1)})) $. 

\bigskip

\begin{thm}\label{Thm f^Omega = 0 implies f^P = 0}
    If $ f \in \CS(\Gamma \bs G, \phi) $ with $ f^\Omega = 0 $, then $ f^Q = 0 $ for all $ \Gamma $-cuspidal parabolic subgroups $Q$ of $G$ with associated cusp of smaller rank.
\end{thm}

\begin{proof}
    Assume that $ k \geq 2 $. We leave the proof for $ k = 1 $ to the reader.

    For a measurable function $ \varphi $ on $ (0, +\infty) $ and $ z \in \CC $ such that $ \varphi(x) x^z \in L^1((0, +\infty)) $, set $ M(\varphi)(z) = \int_0^\infty \varphi(x) x^{z-1} \, dx $.
    The function $ M(\varphi) $ is called the \textit{Mellin transform} of $ \varphi $.

    The Mellin convolution operator of two measurable functions on $ (0, +\infty) $, $ \varphi_1 $ and $ \varphi_2 $, is defined by
    \[
        (\varphi_1 * \varphi_2)(x) := \int_0^\infty \varphi_1(\tfrac{x}{y}) \varphi_2(y) \, \tfrac{dy}{y}
    \]
    whenever the integral exists.

    Let $Q$ be a $ \Gamma $-cuspidal parabolic subgroups of $G$ with associated of smaller rank. Without loss of generality, $ Q = P $.

    Since the space $ \CS(\Gamma \bs G, \phi)_K $ of $K$-finite Schwartz function is dense in $ \CS(\Gamma \bs G, \phi) $, it suffices to show that $ f^Q = 0 $ for all $ \Gamma $-cuspidal parabolic subgroups of $G$ with associated cusp of smaller rank if $ f \in \CS(Y, V_Y(\gamma, \phi)) $ ($ \gamma \in \hat{K} $) with $ f^\Omega = 0 $.

    Let $ \gamma \in \hat{K} $, let $ f \in \CS(Y, V_Y(\gamma, \phi)) $ with $ f^\Omega = 0 $.

    Let $ \eps > 0 $. Choose a function $ \chi_\eps \in \Ccinf((0, +\infty), [0, 1]) $ such that
    \[
        (\eps, +\infty) \subset \supp(\chi_\eps) \ . 
    \]
    Define
    \[
        g_{1, \eps} \colon x \in (0, +\infty)
        \mapsto \chi_\eps(x) \int_{\Gamma_P \bs N M_\Gamma} f^0(n m a_x) \, dn \, dm
    \]
    and
    \[
        g_2 \colon x \in (0, +\infty) \mapsto 1_{(0, 1)}(x) x^{n - \frac{k}2 - 1}(1 - x)^{\frac{k}2 - 1}
        \int_{M_k} \gamma(m') \gamma(m_x)^{-1} \gamma(m')^{-1} \, dm' \ .
    \]
    
    It follows from Proposition \ref{Prop f^P well-defined} that $ x \in (0, +\infty) \mapsto |g_{1, \eps}(x)| x^{-\tfrac{n-1}2} $ belongs to $ L^1((0, +\infty)) $. So, the Mellin transform of $ g_{1, \eps} $ at $ -\tfrac{n-1}2 + i r $ ($ r \in \RR $) is well-defined.
    \\ Since $ n - \frac{k}2 - 1 - \tfrac{n-1}2 - 1 = \tfrac{n-k-1}2 - 1 $, the integral $ \int_0^\infty |g_2(x)| x^{-\tfrac{n-1}2} \, \frac{dx}x $ is less or equal than
    \[
        \int_0^1 x^{\tfrac{n-k-1}2 - 1}(1 - x)^{\frac{k}2 - 1} \, dx = B(\tfrac{n-k-1}2, \tfrac{k}2) \ ,
    \]
    where $ B(x, y) = \int_0^1 u^{x-1} (1 - u)^{y-1} \, du $ ($ \R(x) > 0 $, $ \R(y) > 0 $) denotes the $ \beta $-function. Since $ eP $ is not of full rank, $ k < n-1 $. Thus, $ \tfrac{n-k-1}2 > 0 $. Hence, the above integral is finite. So, $ x \in (0, +\infty) \mapsto |g_2(x)| x^{-\tfrac{n-1}2} $ belongs to $ L^1((0, +\infty)) $. Consequently, the Mellin transform of $ g_2 $ at $ -\tfrac{n-1}2 + i r $ ($ r \in \RR $) is well-defined.

    Thus, the Mellin transform
    \[
        M(g_{1, \eps} * g_2)(-\tfrac{n-1}2 + ir) := \int_0^\infty g_2(y) g_{1, \eps}(\tfrac{x}{y}) \, \tfrac{dy}{y}
    \]
    of $ g_{1, \eps} * g_2 $ at $ -\tfrac{n-1}2 + ir $ ($ r \in \RR $) is well-defined and equal to
    \[
        M(g_2)(-\tfrac{n-1}2 + ir) M(g_{1, \eps})(-\tfrac{n-1}2 + ir) \ .
    \]
    Let $ D = \{ z \in \CC \mid \R(z) > -\frac{n-k-1}2 \} $ (connected, open subset of $ \CC $). Note that the function
    \begin{multline*}
        \Phi \colon z \in D \mapsto M(g_2)(-\tfrac{n-1}2 + z) \\
        = \int_{M_k} \int_0^1 \gamma(m') x^{\tfrac{n-k-1}2 - 1 + z} (1 - x)^{\frac{k}2 - 1} \gamma(m_x)^{-1} \gamma(m')^{-1} \, dx \, dm'
    \end{multline*}
    is complex analytic. Indeed, since $\Phi$ is measurable, since $ B(\cdot, \tfrac{k}2) $ satisfies the assumptions of the theorem about complex differentiation of parameter dependent integrals on $ \{ z \in \CC \mid 0 < r < \R(z) < R \} $ for every $ R > r > 0 $ and since
    \begin{multline*}
        \int_{M_k} \int_0^1 |\gamma(m') x^{\tfrac{n-k-1}2 - 1 + z} (1 - x)^{\frac{k}2 - 1} \gamma(m_x)^{-1} \, \gamma(m')^{-1}| \, dx \, dm' \\
        = B(\tfrac{n-k-1}2 + \R(z), \tfrac{k}2) \ ,
    \end{multline*}
    $ \Phi $ satisfies also the assumptions of that theorem on
    \[
        \{ z \in \CC \mid R > \R(z) > -\tfrac{n-k-1}2 + r \}
    \]
    for every $ R > r > 0 $. Thus, $ \Phi $ is holomorphic and hence complex analytic.
    
    Let $ \{v_j\} $ be a basis of $ V_\phi^{\Gamma_P} $ and let $ p_j $ be the orthogonal projection on $ \CC v_j \equiv \CC $.
    By abuse of notation, we denote $ \Id \otimes \, p_j \colon V_\gamma \otimes V_\phi \to V_\gamma \otimes V_\phi $ also by $ p_j $.

    As by assumption $ f^\Omega = 0 $,
    \[
        M(p_j(g_{1, \eps}) * g_2) = M(g_2) M(p_j(g_{1, \eps}))
    \]
    vanishes for all $ j $, by Proposition \ref{Prop formula Conv of constant terms}.

    For all $ x \in [0, 1] $, $ |\gamma(m_x)^{-1} - \Id| \leq 2 $. Let $ \delta \in (0, 1] $ be such that $ |\gamma(m_x)^{-1} - \Id| \leq \frac12 $ for all $ x \in [0, \delta] $.
    For $ z \in D $, set $ \psi(z) = \int_0^1 x^{\tfrac{n-k-1}2 - 1 + z} (1 - x)^{\frac{k}2 - 1} \, dx $. Then, $ \psi $ is real and positive on $ D \cap \RR $.
    For $ z \in D \cap \RR $, consider
    \begin{align*}
        & \frac1{ \psi(z) } \Phi(z) - \Id \\
        = & \frac1{ \psi(z) } \int_{M_k} \int_0^1 \gamma(m') x^{\tfrac{n-k-1}2 - 1 + z} (1 - x)^{\frac{k}2 - 1} (\gamma(m_x)^{-1} - \Id) \gamma(m')^{-1} \, dx \, dm' \ .
    \end{align*}
    In absolute values, this is less or equal than
    \[
        \frac12 + \frac2{ \psi(z) } \int_\delta^1 x^{\tfrac{n-k-1}2 - 1 + z} (1 - x)^{\frac{k}2 - 1} \, dx \ .
    \]
    Choose $ z \in D \cap \RR $ sufficiently close to $ -\frac{n-k-1}2 $ so that this is less than $ 1 $. This is indeed possible as $ \psi(z) $ tends to $ \infty $ when $ z \in D \cap \RR $ tends to $ -\frac{n-k-1}2 $ and as $ \sup_{ z \in [-\frac{n-k-1}2, -\frac{n-k-1}2 + 1] } \int_\delta^1 x^{\tfrac{n-k-1}2 - 1 + z} (1 - x)^{\frac{k}2 - 1} \, dx $ is finite.

    It follows that $ \frac1{ \psi(z) } \Phi(z) $ and hence $ \Phi(z) $ are invertible for every $ z \in D \cap \RR $ sufficiently close to $ -\frac{n-k-1}2 $.
    Consequently, $ \det(\Phi(z)) $ is a nonzero complex analytic function. It follows from the Identity theorem that $ \{ z \in D \mid \det(\Phi(z)) = 0 \} $ has no accumulation point in $D$.

    Hence, $ M(g_2)(-\tfrac{n-1}2 + ir) $ is invertible for almost all $ r \in \RR $.
    So, $ M(p_j(g_{1, \eps}))(-\tfrac{n-1}2 + ir) $ vanishes for all $j$ and almost all $ r \in \RR $.
    Since $ g_{1, \eps} $ is smooth and since the Mellin transform is injective on $ L^1((0, +\infty), V_\gamma \otimes V_\phi^{\Gamma_P}) $, it follows that $ g_{1, \eps} $ is identically zero.
    As this holds for every $ \eps > 0 $, we conclude that $ f^P(a) = 0 $ for all $ a \in A $. Since $ f^P $ is left $ N $-invariant and right $K$-equivariant, $ f^P $ vanishes, too.
    The theorem follows.
\end{proof}

\cleardoublepage

\begin{appendices}

\standardsection{Some results about direct sum decompositions}

\begin{lemSec}\label{Lem X = Y oplus F topological direct sum if Y is closed and if F is finite-dimensional}
    Let $ X $ be a locally convex space and let $ Y \subset X $ be a closed finite-codimensional subspace. Let $F$ be a finite-dimensional subspace complementing $Y$. Then,
    \[
        X = Y \oplus F
    \]
    as topological vector spaces.
\end{lemSec}

\begin{proof}
    For a proof, see Lemma B.1 of \cite[p.179]{Wur04}.
\end{proof}

\begin{lemSec}\label{Lem V = (S cap V) oplus (S^perp cap V)}
    Let $ H $ be a Hilbert space and let $ V \subset H $ be a Fréchet space such that the injection $ V \hookrightarrow H $ is continuous. Let $ S $ be a closed subspace of $H$.
    Then, $ S \cap V $ and $ S^\perp \cap V $ are closed in $V$ and hence also Fréchet spaces.
    We denote the canonical projection on $ S $ by $ p_S $.
    If $ p_S(V) $ is contained in $V$, then
    \[
        V = (S \cap V) \oplus (S^\perp \cap V)
    \]
    as topological direct sum. 
\end{lemSec}

\begin{proof}
    Let $ H $, $ V $ be as above. Let $ S $ be a closed subspace of $H$.
    As $ \iota \colon V \hookrightarrow H $ is continuous and as $ S $ and $ S^\perp $ are closed in $H$, $ S \cap V = \iota^{-1}(S) $ and $ S^\perp \cap V = \iota^{-1}(S^\perp)$ are closed in $V$.

    Assume moreover that $ p_S(V) $ is contained in $V$. Then, $ p_{S^\perp}(V) $ is also contained in $V$.
    Furthermore, $ p_S(V) = S \cap V $ and $ p_{S^\perp}(V) = S^\perp \cap V $.
    The lemma follows now by the open mapping theorem.
\end{proof}

\begin{lemSec}\label{Lem V = (S cap V) oplus (S^perp cap V) generalisation}
    Let $ H $ be a Hilbert space and let $ V \subset H $ be a Fréchet space that can be continuously injected into $H$. Let $ S_1, \dotsc, S_d $ be closed subspaces of $H$ such that $ H = \bigoplus_{j=1}^d S_d $ (Hilbert direct sum).
    Then, $ S_j \cap V $ is closed in $V$ for all $j$ and hence also a Fréchet space.
    We denote the canonical projection on $ S_j $ by $ p_{S_j} $.
    If $ p_{S_j}(V) $ is contained in $V$ for all $j$, then
    \[
        V = \bigoplus_{j=1}^d (S_j \cap V)
    \]
    as topological direct sum. 
\end{lemSec}

\begin{proof}
    The proof works completely similarly to the one of Lemma \ref{Lem V = (S cap V) oplus (S^perp cap V)}.
    
\end{proof}

\begin{lemSec}\label{Lem dual map is continuous}
    Let $ V $ and $W$ be topological vector spaces and let $ u \colon V \to W $ be a continuous linear map. Let $ V' $ (resp. $ W' $) be (conjugate-)linear strong dual of $V$ (resp. $ W $). Then, the dual map $ u' \colon W' \to V' $ is continuous.
\end{lemSec}

\begin{proof}
    This can easily be checked. In the case where $ V $ and $W$ are locally convex, one can find the proof (for the linear strong duals) of this statement for example in Yosida's book ``Functional Analysis'' (see Theorem 2 of Chapter VII in \cite[p.194]{Yosida}).
\end{proof}

\begin{lemSec}\label{Lem V' simeq V'_1 oplus V'_2}
    Let $ V $ be a topological vector space. Let $ V_1, \dotsc, V_d $ ($ d \geq 2 $) be subspaces of $V$ such that $ V = \bigoplus_{j=1}^d V_j $ is a topological direct sum. Let $ V' $ (resp. $ V'_j $) be the (conjugate-)linear strong dual of $ V $ (resp. $ V_j $).
    For $ j = 1, \dotsc, d $, let $ \iota_j \colon V'_j \to V' $ denote the natural topological embedding. Then,
    \[
        V' = \bigoplus_{j=1}^d \iota_j(V'_j)
    \]
    (topological direct sum).
\end{lemSec}

\begin{proof}
    This lemma can easily be proven with the help of Lemma \ref{Lem dual map is continuous}.

\end{proof}

\section{\texorpdfstring{The group $ \Sp(1, n) $}{The group Sp(1, n)}}\label{sec:Sp(1,n)}

In the following, we do some explicit computations in the group $ \Sp(1, n) $ ($ n \geq 1 $).

\medskip

Let $ n \geq 2 $. We denote the standard norm on $ \HH $ by $ |\cdot| $ and the 2-norm on $ \HH^{n-1} $ by $ \|\cdot \| $.
We identify $ \HH $ with $ \RR^4 $ as vector spaces.
The quaternionic hyperbolic space is defined by
\[
    X = \HH H^n := \{ [1 : x_1 : \ldots : x_n] \mid \sum_{i=1}^n |x_i|^2 < 1 \} \ .
\]
Put $ z_0 = [1  : 0 : \ldots : 0] $. Put $ G = \Sp(1, n) $ and $ K = G_{z_0} = \Sp(n) $ (compact symplectic group). Then, $ X = G/K $. Define
\[
    N = \left\{ n_{v, r} :=
    \begin{pmatrix}
        1+r+\frac12 \|v\|^2 & -r-\frac12 \|v\|^2 & \bar{v}^t \\[0.2em]
        r+\frac12 \|v\|^2 & 1-r-\frac12 \|v\|^2 & \bar{v}^t \\[0.2em]
        v & -v & I_{n-1}
    \end{pmatrix} \mid v \in \HH^{n-1}, \, r \in \I(\HH) \right\} \ ,
\]
\[
    n_v = n_{v, 0} \qquad (v \in \HH^{n-1}) \ .
\]
Put $ n_{v^t, r^t} = n_{v, r} $ ($ v \in \HH^{n-1}, \, r \in \I(\HH) $). We use similar conventions for all the functions that are defined with the help of $ n_{v, r} $ as for example $ n_v $.

Set
\[
    A = \left\{ a_t :=
    \begin{pmatrix}
        \frac12(t + \tfrac1t) & \frac12(t - \tfrac1t) & 0 \\[0.2em]
        \frac12(t - \tfrac1t) & \frac12(t + \tfrac1t) & 0 \\[0.2em]
        0 & 0 & I_{n-1}
    \end{pmatrix}
    \mid t > 0 \right\} \ .
\]
We have
\[
    a_t n_{v, 0} a_t^{-1} = n_{t v, 0} \quad \text{and} \quad a_t n_{0, r} a_t^{-1} = n_{0, t^2r} \ .
\]
Thus, $ a_t^\alpha = t $.
Then, $ n_{v, r} a_t z_0 $ is equal to
\begin{align*}
     &
     \begin{pmatrix}
        1+r+\frac12 \|v\|^2 & -r-\frac12 \|v\|^2 & \bar{v}^t \\[0.2em]
        r+\frac12 \|v\|^2 & 1-r-\frac12 \|v\|^2 & \bar{v}^t \\[0.2em]
        v & -v & I_{n-1}
    \end{pmatrix}
    \begin{pmatrix}
        \frac12(t + \tfrac1t) & \frac12(t - \tfrac1t) & 0 \\[0.2em]
        \frac12(t - \tfrac1t) & \frac12(t + \tfrac1t) & 0 \\[0.2em]
        0 & 0 & I_{n-1}
    \end{pmatrix}
    [1  : 0 : \ldots : 0] \\
    =&
    \begin{pmatrix}
        1+r+\frac12 \|v\|^2 & -r-\frac12 \|v\|^2 & \bar{v}^t \\[0.2em]
        r+\frac12 \|v\|^2 & 1-r-\frac12 \|v\|^2 & \bar{v}^t \\[0.2em]
        v & -v & I_{n-1}
    \end{pmatrix}
    [\tfrac12(t + \tfrac1t) : \tfrac12(t - \tfrac1t) : 0 : \ldots : 0] \\
    =& [\tfrac12(t + \tfrac1t)(1+r+\tfrac12 \|v\|^2) - \tfrac12(t - \tfrac1t)(r + \tfrac12 \|v\|^2) : \\
    & \qquad \tfrac12(t + \tfrac1t)(r+\tfrac12 \|v\|^2) + \tfrac12(t - \tfrac1t)(1-r-\tfrac12 \|v\|^2) : \tfrac12(t + \tfrac1t)v^t - \tfrac12(t - \tfrac1t)v^t] \\
    =& [\tfrac12(t + \tfrac1t) + \tfrac1t(r+\tfrac12 \|v\|^2) :
    \tfrac12(t - \tfrac1t) + \tfrac1t (r+\tfrac12 \|v\|^2) : \tfrac1t v^t] \ .
\end{align*}
For $ g \in \Sp(1, n) $, write
\[
    g z_0 = [a_0 : a_1 : \ldots : a_n] = [1 : a_1 a_0^{-1} : \ldots : a_n a_0^{-1}] \ .
\]
Then,
\begin{align*}
    & a_1 a_0^{-1} = \tfrac12(t - \tfrac1t) + \tfrac1t (r+\tfrac12 \|v\|^2) \big(\tfrac12(t + \tfrac1t) + \tfrac1t(r+\tfrac12 \|v\|^2)\big)^{-1} \\
    \iff & a_1 a_0^{-1} \cdot \tfrac12(t + \tfrac1t) + a_1 a_0^{-1} \cdot \tfrac1t(r+\tfrac12 \|v\|^2) = \tfrac12(t - \tfrac1t) + \tfrac1t (r+\tfrac12 \|v\|^2) \\
    \iff & (a_1 a_0^{-1} - 1) \cdot \tfrac1{2t} \|v\|^2 = \tfrac12(t - \tfrac1t) + (1 - a_1 a_0^{-1})\tfrac{r}t - a_1 a_0^{-1} \cdot \tfrac12(t + \tfrac1t) \ .
\end{align*}
Thus,
\begin{equation}\label{eq norm v^2 quaternionic case}
    \|v\|^2 = 2t (a_1 a_0^{-1} - 1)^{-1} \big( \tfrac12(t - \tfrac1t) + (1 - a_1 a_0^{-1})\tfrac{r}t - a_1 a_0^{-1} \cdot \tfrac12(t + \tfrac1t) \big) \ .
\end{equation}
For $ j = 2, \dotsc, n $, we have
\[
    a_j a_0^{-1} = \frac1t v_{j-1} \cdot \big(\tfrac12(t + \tfrac1t) + \tfrac1t(r+\tfrac12 \|v\|^2)\big)^{-1} \ .
\]
Thus,
\begin{align*}
    v_{j-1} &= t a_j a_0^{-1} \cdot \big(\tfrac12(t + \tfrac1t) + \tfrac1t(r+\tfrac12 \|v\|^2)\big) \\
    &= a_j a_0^{-1} \Big(\tfrac12(t^2 + 1) + r + (a_1 a_0^{-1} - 1)^{-1} \\
    & \qquad \cdot \big( \tfrac12(t^2 - 1) + (1 - a_1 a_0^{-1})r - a_1 a_0^{-1} \cdot \tfrac12(t^2 + 1) \big) \Big) \\
    &= a_j a_0^{-1} (a_1 a_0^{-1} - 1)^{-1} \big( (a_1 a_0^{-1} - 1)(\tfrac12(t^2 + 1) + r) + \tfrac12(t^2 - 1) \\
    & \qquad + (1 - a_1 a_0^{-1})r - a_1 a_0^{-1} \cdot \tfrac12(t^2 + 1) \big) \\
    &= -a_j a_0^{-1} (a_1 a_0^{-1} - 1)^{-1} = a_j (a_0 - a_1)^{-1} \ .
\end{align*}
By \eqref{eq norm v^2 quaternionic case}, we have
\begin{align*}
    (a_1 a_0^{-1} - 1) \|v\|^2
    &= t^2 - 1 + (1 - a_1 a_0^{-1}) \cdot 2r  - a_1 a_0^{-1} \cdot (t^2 + 1) \\
    &= (1 - a_1 a_0^{-1}) t^2 - 1 - a_1 a_0^{-1} + (1 - a_1 a_0^{-1}) \cdot 2r \ .
\end{align*}
Hence,
\begin{align*}
    t^2 &= (1 - a_1 a_0^{-1})^{-1} \big(1 + a_1 a_0^{-1} - (1 - a_1 a_0^{-1}) \cdot 2r + (a_1 a_0^{-1} - 1) \|v\|^2 \big) \\
    &= (1 - a_1 a_0^{-1})^{-1} \Big(1 + a_1 a_0^{-1} - (1 - a_1 a_0^{-1}) \cdot 2r + (a_1 a_0^{-1} - 1) \big(\sum_{j=2}^n \frac{|a_j|^2}{|a_0 - a_1|^2} \big) \Big) \\
    &= a_0 (a_0 - a_1)^{-1}(a_0 + a_1)a_0^{-1} - 2r - \sum_{j=2}^n \frac{|a_j|^2}{|a_0 - a_1|^2} \ .
\end{align*}
Since this belongs to $ \RR $ and since $ \I(x y x^{-1}) = x \I(y) x^{-1} $ for all $ x, y \in \HH $,
\[
    r = \tfrac12 a_0 \I\big((a_0 - a_1)^{-1}(a_0 + a_1)\big) a_0^{-1} \ .
\]
So,
\[
    h(g) = a_{\sqrt{\R\big((a_0 - a_1)^{-1}(a_0 + a_1)\big) - {\displaystyle \sum_{j=2}^n \frac{|a_j|^2}{|a_0 - a_1|^2}}}}
\]
as $ \R(x y x^{-1}) = x \R(y) x^{-1} = \R(y) $ for all $ x, y \in \HH $,
and
\[
    \nu(g) = n_{(a_2 (a_0 - a_1)^{-1}, \dotsc, a_n (a_0 - a_1)^{-1}, \tfrac12 a_0 \I\big((a_0 - a_1)^{-1}(a_0 + a_1)\big) a_0^{-1})} \ .
\]
Let $ m \in M = \left\{ m_{A, q} := \begin{pmatrix}
    q   &   0   &   0 \\
    0   &   q   &   0 \\
    0   &   0   &   A
\end{pmatrix}
\mid A \in \Sp(n-1), q \in \Sp(1)
\right\}
$.
\\ Set $ m_A = m_{A, 1} $. Then,
\begin{equation}\label{eq m_(A, q) n_(v, r) m_(a, q)^(-1) = n_(Av bar(q), q r bar(q))}
    m_{A, q} n_{v, r} m_{A, q}^{-1} = n_{Av \bar{q}, q r \bar{q}}
\end{equation}
and
\begin{equation}\label{eq m_A n_v m_A^(-1) = n_(Av)}
    m_A n_v m_A^{-1} = n_{Av} \ .
\end{equation}

For $ g \in \Sp(1, n) $ (resp. $ X \in \smallSp(1, n) $), define
\[
    \theta(g) = (\bar{g}^t)^{-1} \qquad (\text{resp. } \theta(X) = - \bar{X}^t ) \ .
\]
Let $ \bar{n}_{v, r} = \theta(n_{v, r}) \in \bar{N} $ ($ v \in \HH^{n-1}, \, r \in \I(\HH) $). Since $ \bar{r} = -r $,
\[
    g := \bar{n}_{v, r} = \begin{pmatrix}
        1+r+\frac12 \|v\|^2 & r+\frac12 \|v\|^2 & -\bar{v}^t \\[0.2em]
        -r-\frac12 \|v\|^2 & 1 - r-\frac12 \|v\|^2 & \bar{v}^t \\[0.2em]
        -v & -v & I_{n-1}
    \end{pmatrix} \ .
\]
Thus,
\[
    (a_0 - a_1)^{-1} = (1 + 2r + \|v\|^2)^{-1} = \frac{ 1 + \|v\|^2 - 2r }{(1 + \|v\|^2)^2 + 4|r|^2 } \ .
\]
As $ a_0 + a_1 = 1 $, we have
\begin{align*}
    h(\bar{n}_{v, r}) &= a_{\sqrt{\R\big((a_0 - a_1)^{-1}\big) - {\displaystyle \sum_{j=1}^{n-1} \tfrac{|v_j|^2}{|a_0 - a_1|^2}}}} \\
    &= a_{\sqrt{\frac{ 1 + \|v\|^2 }{(1 + \|v\|^2)^2 + 4|r|^2 } - \tfrac{\|v\|^2}{(1 + \|v\|^2)^2 + 4|r|^2}}}
    = a_{ \frac1{\sqrt{(1 + \|v\|^2)^2 + 4|r|^2} } }
\end{align*}
and, as in addition $ a_0 $ and $ \I((a_0 - a_1)^{-1}) = \I((2a_0 - 1)^{-1}) $ commute,
\begin{align*}
    \nu(\bar{n}_{v, r}) &= n_{(-v_1 (a_0 - a_1)^{-1}, \dotsc, -v_{n-1} (a_0 - a_1)^{-1}, \tfrac12\I((a_0 - a_1)^{-1}))} \\
    &= n_{(-(v_1, v_2, \dotsc, v_{n-1})\frac{ 1 + \|v\|^2 - 2r }{(1 + \|v\|^2)^2 + 4|r|^2 }, -\frac{ r }{(1 + \|v\|^2)^2 + 4|r|^2 }) } \ .
\end{align*}
The formula for $ h(\bar{n}_{v, r}) $ follows also from Proposition \ref{Theorem 3.8 of Helgason, p.414}.

Put
$
    w = \begin{pmatrix} 1&0&0&0\\0&-1&0&0\\0&0&-1&0\\0&0&0&I_{n-2} \end{pmatrix}.
$
Then,
\[
    w n_{v, r} = (w n_{v, r} w) w = \bar{n}_{(v_1, -v_2, \dotsc, -v_{n-1}), -r} w \ .
\]
Thus,
\[
    h(w n_{v, r}) = a_{ \frac1{\sqrt{(1 + \|v\|^2)^2 + 4|r|^2} } }
\]
and
\[
    \nu(w n_{v, r}) = n_{((-v_1, v_2, \dotsc, v_{n-1})\frac{ 1 + \|v\|^2 + 2r }{(1 + \|v\|^2)^2 + 4|r|^2 }, \frac{ r }{(1 + \|v\|^2)^2 + 4|r|^2 }) } \ .
\]
Let $ b \colon \HH^{n+1} \times \HH^{n+1} \to \HH $ be the following $ \HH $-Hermitian form
\[
    b(x, y) = \overline{x_0} y_0 - \sum_{i=1}^n \overline{x_i} y_i \qquad (x = (x_0, \dotsc, x_n), y = (y_0, \dotsc, y_n) \in \HH^{n+1}) \ .
\]
Note that $ b(y, x) = \overline{b(x, y)} $. Define a metric on $ \HH H^n $ by the formula
\[
    \cosh^2\big(d(x, y)\big) = \frac{ b(x, y) b(y, x) }{b(x, x) b(y, y)} \qquad (x, y \in \HH H^n) \ .
\]
For $ g \in \Sp(1, n) $, write
\[
    g z_0 = [a_0 : a_1 : \ldots : a_n] \ .
\]
Since $ b(z_0, z_0) = 1 $, $ b(z_0, g z_0) = a_0 $ and $ b(g z_0, g z_0) = |a_0|^2 - \sum_{i=1}^n |a_i|^2 $,
\[
    a_g = e^{d(z_0, g z_0)} = e^{ \arccosh\bigl(\sqrt{\frac{|a_0|^2}{|a_0|^2 - \sum_{i=1}^n |a_i|^2 }}\bigr) } \ .
\]
Since
\begin{multline*}
    \cosh^2\big(d(z_0, a_t z_0)\big) = \frac{(t^2+1)^2}{(t^2+1)^2 - (t^2-1)^2} = \frac{(t^2+1)^2}{4t^2} \\
    = \big( \frac{t^2+1}{2t} \big)^2 = \big( \frac{t+\frac1t}{2} \big)^2 = \cosh^2(\log t) \ ,
\end{multline*}
$ a_{a_t} = t $ for all $ t \geq 1 $. Moreover,
\[
    a_g = e^{\arccosh\bigl(\sqrt{\frac{|a_0|^2}{|a_0|^2 - \sum_{i=1}^n |a_i|^2 }}\bigr) } \ .
\]
Recall that $ n_{v, r} a_t z_0 $ is equal to $ [\tfrac12(t + \tfrac1t) + \tfrac1t(r+\tfrac12 \|v\|^2) : \tfrac12(t - \tfrac1t) + \tfrac1t (r+\tfrac12 \|v\|^2) : \tfrac1t v^t] $.
Since $ |\tfrac12(t \pm \tfrac1t) + \tfrac1t(r+\tfrac12 \|v\|^2)|^2 = \big(\tfrac12(t \pm \tfrac1t) + \tfrac1{2t} \|v\|^2\big)^2 + \tfrac{|r|^2}{t^2} $, $ \cosh(\log(a_{n_{v, r} a_t}))^2 $ is equal to
\begin{multline*}
    \frac{\big(\tfrac12(t + \tfrac1t) + \tfrac1{2t} \|v\|^2\big)^2 + \tfrac{|r|^2}{t^2} }{\big(\tfrac12(t + \tfrac1t) + \tfrac1{2t} \|v\|^2\big)^2 - \big(\tfrac12(t - \tfrac1t) + \tfrac1{2t} \|v\|^2\big)^2 - \frac1{t^2} \|v\|^2 } \\
    = \frac{\big(\tfrac12(t + \tfrac1t) + \tfrac1{2t} \|v\|^2\big)^2 + \tfrac{|r|^2}{t^2} }{\tfrac1t \big(t + \tfrac1{t} \|v\|^2 \big) - \frac1{t^2} \|v\|^2 }
    = \big(\tfrac12(t + \tfrac1t) + \tfrac1{2t} \|v\|^2\big)^2 + \tfrac{|r|^2}{t^2} \ .
\end{multline*}
Thus,
\begin{equation}\label{eq a_(n_(v, r) a_t) for Sp(1,n)}
    a_{n_{v, r} a_t} = e^{ \arccosh\bigl(\sqrt{ \big(\tfrac12(t + \tfrac1t) + \tfrac1{2t} \|v\|^2\big)^2 + \tfrac{|r|^2}{t^2}}\bigr) } \ .
\end{equation}
In particular, $ a_{n_v a_t} $ is equal to
\[
     e^{ \arccosh\bigl(\tfrac12(t + \tfrac1t) + \tfrac1{2t} \|v\|^2\bigr) } \ .
\]

For $ v \in \HH^{n-1} $ and $ r \in \I(\HH) $, set
$
    X_{v, r} = \begin{pmatrix}
        r & -r & \bar{v}^t \\[0.2em]
        r & -r & \bar{v}^t \\[0.2em]
        v & -v & 0_{n-1}
    \end{pmatrix} .
$
\\ Then, $ \exp(X_{v, r}) = n_{v, r} $.

\smallskip

For $ X, Y \in \g $, $ B_0(X, Y) = \tr(X Y) $. This defines an $ \Ad(G) $-invariant inner product on $ \g $. As these are unique up to a constant for simple groups and as both are positive defined on $ \s := \{ X \in \g \mid \theta(X) = -X\} $, there exists a positive constant $ c $ such that $ B_0(X, Y) = c B(X, Y) $.
For $ X \in \g $, set $ \|X\| = B_0(X, \theta(X)) $.

Let us compute now the value of $c$.

Let $ H \in \a $. Then, $ \|H\|^2 = B_0(H, H) = c B(H, H) = c \tr\bigl(\ad(H) \ad(H)\bigr) $.

Since $ [H, Y] = 0 $ for all $ Y \in \m \oplus \a $ and $ [H, X] = \mu(H) X $ for all $ X \in \g_\mu $,
\[
    \tr\bigl(\ad(H) \ad(H)\bigr) = 2 m_\alpha \alpha(H)^2 + 2 m_{2\alpha} \big(2 \alpha(H)\big)^2 = 2 (m_\alpha + 4 m_{2\alpha}) \alpha(H)^2 \ .
\]
Let $ H_1 \in \a $ be such that $ \alpha(H_1) = 1 $. Then, $ H_1 = \log(a_e) $ as $ a_t^\alpha = t $. 
Since
\[
    \|\log(a_t)\|^2 = \tr(\log(a_t) \log(a_t)) = 2 (\log t)^2
    \ ,
\]
$ \|H_1\|^2 = 2 $.
Thus,
\[
     2 = \|H_1\|^2 = c \tr\bigl(\ad(H_1) \ad(H_1)\bigr) = 2 c \big( m_\alpha + 4 m_{2\alpha} \big)
\]
as $ \alpha(H_1) = 1 $. Hence, $ c = \frac1{m_\alpha + 4 m_{2\alpha}} $.
So,
\[
    |X_{v, r}| = \frac{\sqrt{c_{d_X}}}{\sqrt{c}} \|X_{v, r}\| = \frac{\sqrt{2}}2 \|X_{v, r}\| \ .
\]
Let us check now that
\[
    \|X_{v, r}\| = \|\log n_{v, r}\| = 2\sqrt{\|v\|^2 + |r|^2}
    \qquad \big(v \in \HH^{n-1}, \, r \in \I(\HH)\big) \ .
\]
Since
\begin{align*}
    X_{v, r} \theta(X_{v, r}) &= \begin{pmatrix}
        r & -r & \bar{v}^t \\[0.2em]
        r & -r & \bar{v}^t \\[0.2em]
        v & -v & 0_{n-1}
    \end{pmatrix}
    \begin{pmatrix}
        \bar{r} & \bar{r} & \bar{v}^t \\[0.2em]
        -\bar{r} & -\bar{r} & -\bar{v}^t \\[0.2em]
        v & v & 0_{n-1}
    \end{pmatrix} \\
    &= \begin{pmatrix}
        \|v\|^2 + 2|r|^2 & \|v\|^2 + 2|r|^2 & 2r\bar{v}^t \\[0.2em]
        \|v\|^2 + 2|r|^2 & \|v\|^2 + 2|r|^2 & 2r\bar{v}^t \\[0.2em]
        -2vr & 2vr & \diag(2 |v_1|^2, \dotsc, 2 |v_{n-1}|^2)
    \end{pmatrix} \ ,
\end{align*}
\[
    \|\log n_{v, r}\|^2 = \|X_{v, r}\|^2 = \tr(X_{v, r} \theta(X_{v, r})) = 4\|v\|^2 + 4|r|^2 \ .
\]

\section{Results of independent interest about discrete series representations}\label{sec: Results of independent interest about discrete series representations}

Let $ \Gamma $ be a convex-cocompact, non-cocompact, discrete subgroup of $G$.
\\ Let $ (\phi, V_\phi) $ be a finite-dimensional unitary representation of $ \Gamma $ and let $ \sigma \in \hat{M} $.

Assume that $ X \neq \OO H^2 $ or that $ \delta_\Gamma < 0 $ and $ \lambda \in (0, -\delta_\Gamma) $.

Let us introduce the following short notations: For $ \mu \in \a^*_\CC $, set
\begin{multline*}
    H^{\sigma, \mu}_{-\infty} = \C^{-\infty}(\dX, V(\sigma_\mu)) , \ H^{\sigma, \mu}_{\Omega} = \C^{-\infty}_\Omega(\dX, V(\sigma_\mu)) , \\
    H^{\sigma, \mu}_{-\infty}(\phi) = \C^{-\infty}(\dX, V(\sigma_\mu, \phi))
    \quad \text{and} \quad
    H^{\sigma, \mu}_{\Omega}(\phi) = \C^{-\infty}_\Omega(\dX, V(\sigma_\mu, \phi)) \ .
\end{multline*}
Let $ (\pi, V_\pi) $ be a discrete series representation. Recall that $ \Lambda_{V_{\pi, K}} \in \a^* $ denotes the real part of the highest leading exponent appearing in an asymptotic expansion for $ a \to \infty $ of $ c_{v, \tilde{v}}(ka) $ ($ v \in V_{\pi, \infty} $, $ \tilde{v} \in V_{\pi', K} $).

Let $ \lambda \in \a^* $ be such that $ -(\lambda + \rho) = \Lambda_{V_{\pi, K}} $. As $ \pi $ is a discrete series representation, $ 0 < \lambda \in \a^* $.

Recall that we defined in Section \ref{ssec: Discrete series representations} an embedding map $ \beta $ from $ V_{\pi, -\infty} $ to $ H^{\sigma, -\lambda}_{-\infty} $.
With the help of this map, we constructed then a map
$
    A \colon H^{\sigma, \lambda}_{-\infty} \to H^{\sigma, -\lambda}_{-\infty}
$.
By (19) of \cite[p.100]{BO00}, we have
\[
    \hat{J}_{\sigma, -\mu} \hat{J}_{\sigma, \mu}
    = \hat{J}_{\sigma, \mu} \hat{J}_{\sigma, -\mu} = \frac1{p_\sigma(\mu)} \Id \ . 
\]
Here $ p_\sigma $ is the Plancherel density. Since $ p_\sigma $ is meromorphic, it has a pole of finite order at $ \mu = \lambda $.

Let $ \mu \in \a^*_\CC $ be such that $ \R(\mu) > 0 $. Let
\[
    J^{-}_{-\mu} := \hat{J}_{\sigma, -\mu} := \sigma(w) \hat{J}^w_{\sigma, -\mu}
    \colon H^{\sigma, -\mu}_{-\infty} \to H^{\sigma, \mu}_{-\infty} \ ,
\]
let $ J_{-} = J^{-}_{-\lambda} $ and let $ I^{\sigma, \lambda} $ denote the image of $ J_{-} $.

Let $ J^{+}_\mu \colon H^{\sigma, \mu}_{-\infty} \to H^{\sigma, -\mu}_{-\infty} $ be the renormalised operator $ \sigma(w) \hat{J}^w_{\sigma, \mu} $ which is regular and nonzero in $ \mu = \lambda $ and which satisfies $ J^{+}_\mu J^{-}_{-\mu} = (\mu - \lambda)^k \Id $ for some $ k \in \NN $. Let $ J_{+} = J^{+}_\lambda $.

Then, the operators $J_{-}$ and $J_{+}$ induce nonzero intertwining operators
\[
    J_{-} \colon H^{\sigma, -\lambda}_{-\infty}(\phi) \to H^{\sigma, \lambda}_{-\infty}(\phi)
    \quad \text{and} \quad
    J_{+} \colon H^{\sigma, \lambda}_{-\infty}(\phi) \to H^{\sigma, -\lambda}_{-\infty}(\phi) \ .
\]
Let $ V_{\pm \infty} = \beta(V_{\pi, \pm \infty}) $, $ V = (V_{\pm \infty})_K = \beta(V_{\pi, K}) $, $ W_{-\infty} = \ker J_{-} $, $ W_\infty = \ker(J_{-}) \cap H^{\sigma, -\lambda}_\infty $, $ W = (\ker J_{-})_K $, $ V_\Omega = V_{-\infty} \cap H^{\sigma, -\lambda}_\Omega $, $ V_\Omega^\Gamma = V_{-\infty}^\Gamma \cap H^{\sigma, -\lambda}_\Omega $, $ W_\Omega = W_{-\infty} \cap  H^{\sigma, -\lambda}_{\Omega} $ and $ W_\Omega^\Gamma = W_{-\infty}^\Gamma \cap H^{\sigma, -\lambda}_{\Omega} $.

Let $ J_{-}^\Gamma $ be the restriction of $ J_{-} $ to $ (H^{\sigma, -\lambda}_{-\infty})^\Gamma $. Let $ W_{-\infty}^\Gamma = \ker J_{-}^\Gamma $.

In an obvious manner, we define also $ V_{\pm \infty}(\phi) $, $ V(\phi) $, $ V_\Omega^\Gamma(\phi) $ and $ W_\Omega^\Gamma(\phi) $.

\medskip

\needspace{3\baselineskip}
We want to prove the following:
\begin{thmSec}\label{Thm V_Omega subset V_(-infty)^Gamma is dense}
    The space $ V_\Omega^\Gamma(\phi) $ is dense in $ V_{-\infty}^\Gamma(\phi) $. Moreover, $ A((H^{\sigma, \lambda}_\Omega)^\Gamma(\phi)) $ $(\text{resp.}$ $ A((H^{\sigma, \lambda}_{-\infty})^\Gamma(\phi)) )$ has finite codimension in $ V_\Omega^\Gamma(\phi) $ $(\text{resp.}$ $ V_{-\infty}^\Gamma(\phi) ) $.
\end{thmSec}

The proof contains several statements of independent interest.

\begin{thmSec}\label{Thm W = V oplus Z}
    There exists a second discrete series module $Z$, $ Z \not \simeq V $, such that $ W = V \oplus Z $.
\end{thmSec}

\begin{remSec}\label{Rem of Thm W = V oplus Z}\nlenum
    \begin{enumerate}
    \item It follows that we have the following nonsplit exact sequence:
        \[
            0 \to V \oplus Z \to H^{\sigma, -\lambda}_K \to I^{\sigma, \lambda} \to 0 \ ,
        \]
        where $ I^{\sigma, \lambda} $ is the Langlands quotient.
    \item As the proof shows, this result is also always true when $ X = \OO H^2 $.
    \end{enumerate}
\end{remSec}

\begin{propSec}\label{Prop W_(pm infty) = J_+(H^(sigma, lambda)_(pm infty)}
    The space $ W_{\pm \infty} $ is equal to $ J_{+}(H^{\sigma, \lambda}_{\pm \infty}) $.
\end{propSec}

\begin{remSec}\label{Rem J^+_mu J^-_(-mu) = (mu - lambda) Id}
    It follows from the proof that $ J^{+}_{\mu} J^{-}_{-\mu} = J^{-}_{\mu} J^{+}_{-\mu} = (\mu - \lambda) \Id $.
\end{remSec}

\begin{corSec}\label{Cor A = c p_V circ J_+}
    The map $ A $ is up to a constant equal to $ p_V \circ J_{+} $, where $ p_V \colon W_{-\infty} \to V_{-\infty} $ is the natural projection.
\end{corSec}

\begin{propSec}\label{Prop W_Omega is dense in W_(-infty)^Gamma}
    The space $ W_\Omega^\Gamma(\phi) $ is dense in $ W_{-\infty}^\Gamma(\phi) $. Moreover,
    \[
        J_+((H^{\sigma, \lambda}_\Omega)^\Gamma(\phi)) \qquad \big(\text{resp. } J_+((H^{\sigma, \lambda}_{-\infty})^\Gamma(\phi)) \big)
    \]
    has finite codimension in $ W_\Omega^\Gamma(\phi) $ $(\text{resp.}$ $ W_{-\infty}^\Gamma(\phi) ) $.
\end{propSec}

\begin{propSec}\label{Prop p_V(W_Omega) = V_Omega}
    We have: $ p_V(W_\Omega) = V_\Omega $.
\end{propSec}

\begin{remSec}\label{Rem of Prop p_V(W_Omega) = V_Omega}
    As $ W_{-\infty} = V_{-\infty} \oplus Z_{-\infty} $ by Theorem \ref{Thm W = V oplus Z}, it follows from this proposition that $ p_{V(\phi)}(W_\Omega^\Gamma(\phi)) = V_\Omega^\Gamma(\phi) $.
\end{remSec}

\begin{proof}[Proof of Theorem \ref{Thm V_Omega subset V_(-infty)^Gamma is dense} using the other statements] \nlenum \\
    By Proposition \ref{Prop W_Omega is dense in W_(-infty)^Gamma}, $ W_\Omega^\Gamma(\phi) $ is dense in $ W_{-\infty}^\Gamma(\phi) $. Moreover, by its proof, there is a finite-dimensional vector space $F$ in $ W_\Omega^\Gamma(\phi) $ such that $ J_+((H^{\sigma, \lambda}_\Omega)^\Gamma(\phi)) \oplus F = W_\Omega^\Gamma(\phi) $ and $ J_+((H^{\sigma, \lambda}_{-\infty})^\Gamma(\phi)) \oplus F = W_{-\infty}^\Gamma(\phi) $.
    By Remark \ref{Rem of Prop p_V(W_Omega) = V_Omega}, $ V_\Omega^\Gamma(\phi) = p_{V(\phi)}(W_\Omega^\Gamma(\phi)) $.
    Thus, $ V_\Omega^\Gamma(\phi) $ is dense in $ p_{V(\phi)}(W_{-\infty}^\Gamma(\phi)) = V_{-\infty}^\Gamma(\phi) $. Moreover, it follows from Corollary \ref{Cor A = c p_V circ J_+} that $ V_\Omega^\Gamma(\phi) = p_{V(\phi)}(W_\Omega^\Gamma(\phi)) $ (resp. $ V_{-\infty}(\phi) = p_{V(\phi)}(W_{-\infty}(\phi)) $) is equal to $ A(H^{\sigma, \lambda}_\Omega(\phi)) + p_{V(\phi)}(F) $ (resp. $ A(H^{\sigma, \lambda}_{-\infty}(\phi)) + p_{V(\phi)}(F) $). The theorem follows.
\end{proof}

\begin{proof}[Theorem \ref{Thm W = V oplus Z} implies Proposition \ref{Prop p_V(W_Omega) = V_Omega}]
    Since $ V_\Omega $ is contained $ W_\Omega $, $ V_\Omega $ is also contained in $ p_V(W_\Omega) $.
    We have to show that $ p_V(W_\Omega) $ is also contained in $ V_\Omega $.

    Let $ \gamma \in \hat{K} $ and $ V_\gamma \subset V_\pi(\gamma) $ as in Section \ref{ssec: Discrete series representations} and let $ \tilde{t} \in \Hom_K(V_{\tilde{\gamma}}, V_{\pi'}) $ (unique up to a constant). We choose the same as in the definition of the embedding $ \beta $. Set $ S(\beta v) = \t{\tilde{t}}(v) \in V_\gamma $. Then, $ S $ is a nonzero element of $ \Hom_K(V, V_\gamma) $.

    By Theorem \ref{Thm W = V oplus Z}, we have
    $
        W = V \oplus Z
    $
    as $K$-representations. 
    This decomposition provides us with a projection $ p_V $ from $ W $ to $V$.
    So, we can view $ S \circ p_V $ also as an element of $ \Hom_K(W, V_\gamma) $ which vanishes on $Z$ but not on $V$.

    Let $ \tilde{p} $ be a $K$-equivariant projection from $ H^{\sigma, -\lambda}_K $ to $ W $.

    Define a $K$-homomorphism $ \tilde{S} \in \Hom_K(H^{\sigma, -\lambda}_K, V_\gamma) $ by $ S \circ p_V \circ \tilde{p} $.

    By Frobenius reciprocity, we get $ T' \in \Hom_M(V_\sigma, V_\gamma) $. By construction,
    \[
        P^{T'}_{-\lambda} f = 0
    \]
    for all $ f \in Z $.
    Since $ \Hom_K(V_{\tilde{\gamma}}, V_{\pi'}) $ is 1-dimensional, the nonzero map $ v \in V_{\tilde{\gamma}} \mapsto \t{\beta} v_{T', \lambda} \in V_{\pi'} $ is up to a constant equal to $ \tilde{t} $. We choose $ T' \in \Hom_K(V_{\tilde{\gamma}}, V_{\pi'}) $ so that we have equality. Thus, $ \t{A} v_{T', \lambda} = \t{q} \tilde{t} v = v_{T, -\lambda} $, where $ T \in \Hom_M(V_\sigma, V_\gamma) $ is defined by \eqref{eq Def of T} (see p.\pageref{eq Def of T}). Hence, $ P^{T'} A f = P^{T} f $ for all $ f \in H^{\sigma, \lambda}_{-\infty} $.

    Let $ t \in \Hom_M(V_\gamma, V_\sigma) $ be as in the definition of $ \beta $. For $ f \in W_\infty $ and $ k \in K $, define $ p(f)(k) = t \lim_{a \to \infty} a^{\lambda + \rho} P^{T'}_{-\lambda}(f)(ka) $. Since $ f \in W_\infty $ is smooth, $ P^T_{-\lambda} f(k a) $, $ k \in K $, has an asymptotic expansion for $ a \to \infty $ whose coefficients depend smoothly on $ k \in K $ (look at the proof of Lemma \ref{Lem Asymptotics of c (Af, v)(kah) for -lambda}). Thus, $ p $ is well-defined. Extend $p$ to a map from $ W_{-\infty} $ to $ V_{-\infty} $.

    If $ f \in W_\Omega $, then $ p(f) $ is smooth on $ \Omega $ and $ p(f)(k) $ ($ k \in K(\Omega) $) is equal to
    \[
        t \lim_{a \to \infty} a^{\lambda + \rho} P^{T'}_{-\lambda} f(ka)
    \]
    ($ t P^{T'}_{-\lambda} f(ka) $ has an asymptotic expansion when $ k \in K(\Omega) $ and the leading term depends smoothly on $ k \in K(\Omega) $, compare with Corollary \ref{Cor of Lem Asymptotics of c (Af, v)(kah) for -lambda}).
    Let $ \beta v \in V_\Omega $. Then,
    \[
        p(\beta v)(k) = t \lim_{a \to \infty} a^{\lambda + \rho} P^{T'}_{-\lambda} \beta v(ka) = t(\tilde{\beta} v(k)) = (\beta v)(k) \ .
    \]
    It follows that $ p $ is the identity on $ V_{-\infty} $. Thus, $ p = p_V $ by uniqueness of the projection.
    Hence, $ p_V(f) = p(f) $ is smooth on $ \Omega $. So, $ p_V(W_\Omega) $ is contained in $ V_\Omega $.
\end{proof}

Let $ \star $ be $ \infty $ (or $ -\infty $, resp.) and let $ \sharp $ be $ \emptyset $ (or $ \Omega $, resp.).

\begin{proof}[Theorem \ref{Thm W = V oplus Z} and Proposition \ref{Prop W_(pm infty) = J_+(H^(sigma, lambda)_(pm infty)} imply Proposition \ref{Prop W_Omega is dense in W_(-infty)^Gamma}]
    Let $ E^1(\sigma_\mu, \phi) $ be as in Chapter 4 of \cite[p.37]{Olb02} ($ k = 1 $) and let
    \[
        E^1_\sharp(\sigma_\mu, \phi) = E^1(\sigma_\mu, \phi) \cap \C^{-\infty}_\sharp(\dX, V(\sigma_\mu, \phi))^\Gamma \ .
    \]
    Then, $ E^1_\sharp(\sigma_\mu, \phi) $ is equal to the space of distributions $ f \in \C^{-\infty}_\sharp(\dX, V(\sigma_\mu, \phi))^\Gamma $ which can be extended to a holomorphic family
    \[
        \nu \mapsto f_\nu \in \C^{-\infty}_\sharp(\dX, V(\sigma_\nu, \phi))^\Gamma
    \]
    in a neighbourhood of $ \mu $.

    Since $ \ext $ has at most finite-dimensional singularities by Theorem 5.10 of \cite[p.103]{BO00} and since $ \Cinf(B, V_B(\sigma_{-\lambda}, \phi)) $ is dense in $ \C^{-\infty}(B, V_B(\sigma_{-\lambda}, \phi)) $, $ E^1(\sigma_{-\lambda}, \phi) $ has finite codimension and there exists a finite-dimensional space $F$ in $ (H_\Omega^{\sigma, -\lambda})^\Gamma(\phi) $ such that
    \[
        \C^{-\infty}_\sharp(\dX, V(\sigma_{-\lambda}, \phi))^\Gamma = E^1_\sharp(\sigma_{-\lambda}, \phi) \oplus F \ .
    \]
    \begin{clmSec}
        The map
        \[
            \restricted{\rest \circ J_{-}}{ E^1_\sharp(\sigma_{-\lambda}, \phi) } \colon E^1_\sharp(\sigma_{-\lambda}, \phi) \to \C^{-\star}(B, V_B(\sigma_\lambda, \phi))
        \]
        is surjective.
    \end{clmSec}

    \begin{proof}
        Let $ \varphi \in \C^{-\star}(B, V_B(\sigma_\lambda, \phi)) $. Define $ \tilde{e} \in E^1_\sharp(\sigma_{-\lambda}, \phi) $ by the evaluation of
        \[
            \tilde{e}_{-\mu} := J^{+}_\mu \ext(\varphi_\mu)
        \]
        at $ \mu = \lambda $. Then,
        \[
            \rest(J_{-}(\tilde{e})) = \rest(\ext(\varphi)) = \varphi
        \]
        by \eqref{eq res o ext = Id}, as $ \ext $ is regular at $ \lambda $ by Lemma \ref{Lem ext regular} and as $ J^{-}_{\mu} J^{+}_{-\mu} = (\mu - \lambda) \Id $ by Remark \ref{Rem J^+_mu J^-_(-mu) = (mu - lambda) Id}. The claim follows.
    \end{proof}
    Let now $ f \in \C^{-\infty}_\sharp(\dX, V(\sigma_{-\lambda}, \phi))^\Gamma $.
    By the previous claim, there exists $ \tilde{e} \in E^1_\sharp(\sigma_{-\lambda}, \phi) $ such that
    $ \rest(J_{-}(f)) = \rest(J_{-}(\tilde{e})) $. Let $ e = J_{-}(\tilde{e}) $. Thus, there is $ g \in \C^{-\infty}(\Lambda, V(\sigma_\lambda, \phi))^\Gamma $ such that $ J_{-}(f) = e + g $. Let $ h = f - \tilde{e} $. Then, $ f = \tilde{e} + h $ and
    \[
        J_{-}(h) = J_{-}(f) - J_{-}(\tilde{e}) = J_{-}(f) - e = g \in \C^{-\infty}(\Lambda, V(\sigma_\lambda, \phi))^\Gamma \ ,
    \]
    i.e. $ h \in \ker(\rest \circ J_{-}) $. Moreover, $ h \in \C^{-\infty}_\sharp(\dX, V(\sigma_{-\lambda}, \phi))^\Gamma $.

    It follows that there exists a finite-dimensional subspace of
    $
        \ker(\rest \circ J_{-}) \cap (H^{\sigma, -\lambda}_\Omega)^\Gamma ,
    $
    which we denote by abuse of notation also by $F$, such that
    \[
        \C^{-\infty}_\sharp(\dX, V(\sigma_{-\lambda}, \phi))^\Gamma = E^1_\sharp(\sigma_{-\lambda}, \phi) \oplus F \ .
    \]
    Let $ F_0 = F \cap \ker(J_{-}^\Gamma) \subset W_{-\infty}^\Gamma(\phi) $ and let $ F_1 $ be an algebraic complement of $ F_0 $ in $F$.

    Let now $ f \in \ker(J_{-}^\Gamma) \cap \C^{-\infty}_\sharp(\dX, V(\sigma_{-\lambda}, \phi))^\Gamma $. Write $ f = \tilde{e} + h_0 + h_1 $ with $ \tilde{e} \in E^1_\sharp(\sigma_{-\lambda}, \phi) $, $ h_0 \in F_0 $ and $ h_1 \in F_1 $. Thus,
    \[
        J_{-}(\tilde{e}) = -J_{-}(h_1) \in E^1_\sharp(\sigma_\lambda, \phi) \cap \C^{-\infty}(\Lambda, V(\sigma_\lambda, \phi))^\Gamma
        = \{0\}
    \]
    as $ \rest \circ \ext = \Id $ by \eqref{eq res o ext = Id}.
    Hence, $ J_{-} h_1 = 0 $, i.e. $ h_1 \in F_0 $. So, $ h_1 = 0 $ and $ \tilde{e} \in \ker J_{-} $. Consequently, $ f \in (E^1_\sharp(\sigma_{-\lambda}, \phi) \cap \ker J_{-}) \oplus F_0 $. We have just proven that $ \ker(J_{-}^\Gamma) \cap \C^{-\infty}_\sharp(\dX, V(\sigma_{-\lambda}, \phi))^\Gamma $ is equal to $ (E^1_\sharp(\sigma_{-\lambda}, \phi) \cap \ker J_{-}^\Gamma) \oplus F_0 $.

    \begin{clmSec}\label{Clm E^1 cap W^Gamma_(-infty) = J_+(E^1)}
        We have
        \[
            E^1_\sharp(\sigma_{-\lambda}, \phi) \cap W_{-\infty}^\Gamma(\phi) = J_{+}(E^1_\sharp(\sigma_{\lambda}, \phi)) \ .
        \]
    \end{clmSec}

    \begin{proof}
        It follows from Proposition \ref{Prop W_(pm infty) = J_+(H^(sigma, lambda)_(pm infty)} that $ J_{+}(E^1_\sharp(\sigma_{\lambda}, \phi)) $ is contained in
        \[
            E^1_\sharp(\sigma_{-\lambda}, \phi) \cap W_{-\infty}^\Gamma(\phi) \ .
        \]
        Let $ f \in E^1_\sharp(\sigma_{-\lambda}, \phi) \cap W_{-\infty}^\Gamma(\phi) $. Extend $f$ to a holomorphic family
        \[
            \mu \mapsto f_{-\mu} \in \C^{-\infty}_\sharp(\dX, V(\sigma_{-\mu}, \phi))^\Gamma \ ,
        \]
        defined in a neighbourhood of $ \lambda $.
        Since $ f \in \ker(J_{-}) $, $ J^{-}_{-\mu} f_{-\mu} $ has a zero in $ \lambda $. Thus,
        \[
            g_\mu := \frac1{ \mu - \lambda } J^{-}_{-\mu} f_{-\mu}
        \]
        is well-defined and belongs to $ \C^{-\infty}_\sharp(\dX, V(\sigma_\mu, \phi))^\Gamma $. Moreover, $ J^{+}_{\mu} g_\mu $ evaluated at $ \lambda $ is equal to $ f $ as $ J^{+}_{\mu} J^{-}_{-\mu} = (\mu - \lambda) \Id $ by Remark \ref{Rem J^+_mu J^-_(-mu) = (mu - lambda) Id}. The claim follows.
    \end{proof}

    Since $ J_{+}(\C^{-\infty}(\Lambda, V(\sigma_\lambda, \phi))^\Gamma) = \{0\} $ and since $ \ext $ is regular at $ \lambda $, $ J_{+}(E^1_\sharp(\sigma_{\lambda}, \phi)) $ is equal to $ J_{+}(\C^{-\infty}_\sharp(\dX, V(\sigma_\lambda, \phi))^\Gamma) $.

    It follows from this and Claim \ref{Clm E^1 cap W^Gamma_(-infty) = J_+(E^1)} that
    \[
        W_{-\infty}^\Gamma(\phi) = \ker(J_{-}^\Gamma) = J_{+}(\C^{-\infty}(\dX, V(\sigma_\lambda, \phi))^\Gamma) \oplus F_0
    \]
    and that
    \[
        W_\Omega^\Gamma(\phi) = \ker(J_{-}^\Gamma) \cap \C^{-\infty}_\Omega(\dX, V(\sigma_{-\lambda}, \phi))^\Gamma
        = J_{+}(\C^{-\infty}_\Omega(\dX, V(\sigma_\lambda, \phi))^\Gamma) \oplus F_0 \ .
    \]
    As $ J_{+}( (H_\Omega^{\sigma, \lambda})^\Gamma(\phi) ) $ is dense in $ J_{+}( (H^{\sigma, \lambda})^\Gamma(\phi) ) $ and as $ F_0 $ is finite-dimensional, the proposition follows.
\end{proof}

\begin{proof}[Theorem \ref{Thm W = V oplus Z} implies Proposition \ref{Prop W_(pm infty) = J_+(H^(sigma, lambda)_(pm infty)}]
    Let us now introduce certain extensions of the bundles $ V(\sigma_\mu, \phi) $ with itselves (see Chapter 4 of \cite{Olb02} for further details). Let $\Pi$ be the space of polynomials on $\a$. Then, the group $A$ acts on $ \Pi $ by translations and this action extends to a representation $ 1^+ \colon MAN \rightarrow \GL(\Pi) $:
    \\ For $ f \in \Pi $, $ m \in M $, $ a \in A $, $ n \in N $ and $ H \in \a $, set
    \[
        1^+(man)f(H) = f(H - \log a) \ .
    \]
    For any $ k \in \NN $, we denote the finite dimensional subspace of polynomials of degree at most $ k-1 $ by $ \Pi^k $. It is invariant under the above action. Let us denote the restriction of $ 1^+ $ to $ \Pi^k $ by $ 1^k $. Set $ V^k(\sigma_{\pm \lambda}) = V(\sigma_{\pm \lambda} \otimes 1^k) $ and set $ V^0(\sigma_{\pm \lambda}) = \{0\} $.
    Let $ H^{(k)}_{\sigma, \pm \lambda} = \C^{-\infty}(\dX, V^k(\sigma_{\pm \lambda})) $. We can identify $ H^{(k)}_{\sigma, \pm \lambda} $ with germs at $ \pm\lambda $ of holomorphic families $ \mu \mapsto f_\mu \in H^{\sigma, \mu}_{-\infty} $ modulo families vanishing of order $k$. In general,
    \begin{equation}\label{eq short exact sequence for H^(2)_(sigma, mu)}
        0 \to H^{\sigma, \pm \lambda}_{-\infty} \hookrightarrow H^{(2)}_{\sigma, \pm \lambda} \twoheadrightarrow H^{\sigma, \pm \lambda}_{-\infty} \to 0
    \end{equation}
    is a nonsplitting short exact sequence of $G$-representations.
    
    In short notation, we write $ H^{(2)}_{\sigma, \pm \lambda} = \Quotient{H^{\sigma, \pm \lambda}_{-\infty}}{H^{\sigma, \pm \lambda}_{-\infty}} $ (notation as in \cite{Collingwood}).
    The image of the embedding of $ H^{\sigma, \mu}_{-\infty} $ is equal to $ H^{(1)}_{\sigma, \mu} $.

    By Chapter 4 of \cite[p.30]{Olb02}, there is a natural operator $ \rho $ from $ H^{(2)}_{\sigma, -\lambda} $ to $ H^{(1)}_{\sigma, -\lambda} $, which is called \textit{shift operator}. With the above identification, $ \rho $ can be identified with the multiplication by $ \mu - \lambda $.

    The intertwining operator $ J_{-} $ (resp. $ J_{+} $) induces an intertwining operator
    \[
        J^{(2)}_{-} \colon H^{(2)}_{\sigma, -\lambda} \to H^{(2)}_{\sigma, \lambda} \quad
        (\text{resp. } J^{(2)}_{+} \colon H^{(2)}_{\sigma, \lambda} \to H^{(2)}_{\sigma, -\lambda} ) \ .
    \]
    We have the following commutative diagram:
    \[
        \xymatrix@R4pc@C4pc@M2.5mm{
            0 \ar[r] & H^{\sigma, -\lambda}_{-\infty} \ar@{->}[d]_{J_{-}} \ar@{^{(}->}[r]_{\iota_{-}} & H^{(2)}_{\sigma, -\lambda} \ar[d]_{J^{(2)}_{-}} \ar@{->>}[r]_{p_{-}} & H^{\sigma, -\lambda}_{-\infty} \ar[d]_{J_{-}} \ar[r] & 0 \\
            0 \ar[r] & H^{\sigma, \lambda}_{-\infty} \ar@{^{(}->}[r]_{\iota_{+}} \ar@<-1ex>[u]_{J_{+}} & H^{(2)}_{\sigma, \lambda} \ar@<-1ex>[u]_{J^{(2)}_{+}} \ar@{->>}[r]_{p_{+}} & H^{\sigma, \lambda}_{-\infty} \ar[r] \ar@<-1ex>[u]_{J_{+}} & 0
        }
    \]
    Recall that $ J^{+}_\mu J^{-}_{-\mu} = (\mu - \lambda)^l \Id $ for some $ l \in \NN $. Hence, $ J^{(2)}_{+} J^{(2)}_{-} $ is equal to $ (- \rho)^l $. This map would be zero if $ l $ were greater or equal than $ 2 $.

    It follows from Theorem \ref{Thm W = V oplus Z} that $ (H_{\sigma, -\lambda}^{(2)})_K = \DoubleQuotient{ I }{ V \oplus Z }{ I }{ V \oplus Z } $ (here every division of the outer box into two parts corresponds to an exact sequence), where $ I := I^{\sigma, \lambda} $ is the Langlands quotient and $Z$ is a second discrete series module which is not isomorphic to $V$.

    Since moreover $ J^{(2)}_{-} $ does not vanish on $ I \subset (H_{\sigma, -\lambda}^{(1)})_K $, $ (\ker J^{(2)}_{-})_K $ has only composition factors $ V, Z $ appearing once or twice.

    Let $M$ be a Harish-Chandra module having composition factors only discrete series modules (e.g. $ M = (\ker J^{(2)}_{-})_K $), then $M$ is a direct sum of discrete series modules.

    Indeed, the exponents of $ H^0(\n, M) $ ($ \a $-weights) are only those that appear in $ H^0(\n, V) $ or $ H^0(\n, Z) $.

    Thus, a suitable matrix coefficient map gives an embedding of $M$ into $ L^2(G) \oplus \dotsm \oplus L^2(G) $ which is a unitary representation. Hence, the image can be decomposed.

    Consequently, $ (\ker J^{(2)}_{-})_K = V \oplus Z \oplus R $, where $ R \cap H_{\sigma, -\lambda}^{(1)} = \{0\} $.

    But by the proof of Proposition 4.21 of \cite[p.44]{Olb02} and as $ \lambda \neq 0 $, every nonzero submodule of $ (H_{\sigma, -\lambda}^{(2)})_K $ meets $ (H_{\sigma, -\lambda}^{(1)})_K $ nontrivially. So, $ R = \{0\} $.

    It follows that $ (\ker J^{(2)}_{-})_K = W \subset (H_{\sigma, -\lambda}^{(1)})_K $.
    Let $ S = \ThreeQuotient{W}{I}{W} $ (submodule of $ H^{(2)}_{\sigma, -\lambda} $). Then,
    \[
        \rho J^{(2)}_{-}(S) = J^{(2)}_{-}(\rho S) = J^{(2)}_{-}(W) = \{0\}
    \]
    by definition of the $ \rho $-shift and as $ (\ker J^{(2)}_{-})_K = W $. Thus, $ J^{(2)}_{-}(S) $ is contained in $ \ker \rho = H_{\sigma, \lambda}^{(1)} $. It follows that $ J^{(2)}_{-}(S) $ is equal to $ H_{\sigma, \lambda}^{(1)} $ as $ (\ker J^{(2)}_{-})_K = W $.
    Hence, $ (H_{\sigma, \lambda}^{(1)})_K \subset \im(J^{(2)}_{-}) $.
    Moreover, $ J^{(2)}_{+} (H_{\sigma, \lambda}^{(1)})_K \neq 0 $ as $ J_+ H^{\sigma, \lambda}_K \neq 0 $.
    So, $ J^{(2)}_{+} J^{(2)}_{-} \neq 0 $. So, $ l = 1 $.

    Consequently, $ J^{(2)}_{-}(\rho^{-1}(W)) \subset (H_{\sigma, \lambda}^{(1)})_K $ by definition of the $ \rho $-shift. This implies that $ W = J^{(2)}_{+} J^{(2)}_{-}(\rho^{-1}(W)) \subset J^{(2)}_{+}((H_{\sigma, \lambda}^{(1)})_K) $.
    Thus, $ W $ is contained in $ J_{+}((H_{\sigma, \lambda})_K) $.

    Since $ J^{-}_{-\mu} J^{+}_\mu = J^{+}_\mu J^{-}_{-\mu} = (\mu - \lambda) \Id $, $ J_{+}(H^{\sigma, \lambda}_K) $ is also contained in $W$.
    The proposition follows now from the Casselman-Wallach globalisation theory (cf. Chapter 11 in \cite{Wallach2}).
\end{proof}

\begin{proof}[Proof of Theorem \ref{Thm W = V oplus Z}]
    The space of distribution vectors  $ V_{\pi, -\infty} $ can be embedded into some principal series $ H^{\sigma, -\lambda}_{-\infty} $ by using the leading terms of asymptotic expansions of matrix coefficients if and only if $ V_{\pi, K} $ appears as a submodule of $ H^{\sigma, -\lambda}_K $ and $ \lambda \in \a^* $ is minimal among all $ \mu \in \a^* $ with $ V_{\pi, K} \hookrightarrow H^{\sigma, -\mu} $.

    This is Casselman's subrepresentation theorem combined with Frobenius reciprocity and a theorem of Mili\v{c}i\'{c} (stated in \cite[Theorem 1.5.5, p.32]{Collingwood}).

    The assumption that $G$ is linear and has a discrete series implies that the set of infinitesimal characters of discrete series representations coincides with the set of infinitesimal characters of all finite-dimensional representations of $G$ (which are of the form $ \mu_\tau + \rho_\g $, where $ \mu_\tau $ is the highest weight of an irreducible finite-dimensional representation and $ \rho_\g $ is the half-sum of positive roots (with respect to the positive system $ \Phi^+ $, defined below)). Both are the regular integral infinitesimal characters.

    The translation functor sending $ V $ to $ (V \otimes F)^{\chi_F} $ gives an equivalence of categories between the category of Harish-Chandra modules with trivial (generalised) infinitesimal characters and the ones with (generalised) infinitesimal characters $ \chi_F $ (\cite[p.66]{Collingwood}).

    This functor sends discrete series representations to discrete series representations and principal series representations to principal series representations (see, e.g., Proposition 6.1 of \cite[p.57]{Olb02}).

    Since $G$ is covered by one of the groups $ \Spin(1, 2m), \SU(1, n), \Sp(1, n), F_4 $ (if it has a discrete series), we may assume that $G$ is one of these groups.

    The above discussion reduces the question of embeddings of discrete series representations into principal series representation to the case of trivial infinitesimal character $ \chi_{\rho_\g} $.

    These embeddings have been determined by Collingwood: See Proposition 6.2.7--6.2.10 in \cite[pp.140-141]{Collingwood}. 
    There the discrete series representations are denoted by $ \pi_i $ (only one index $i$).

    We want to show that there exists a second discrete series module $Z$, $ Z \not \simeq V $, such that $ W = V \oplus Z $ and that $ H^{\sigma, -\lambda}_K $ is isomorphic to $ \Quotient{I^{\sigma, \lambda}}{V \oplus Z} $. We say then that $ H^{\sigma, -\lambda}_K $ ``decomposes according to Schmid's character identity'' (see 5.1.7 in \cite[p.107]{Collingwood}). We shortly say that the embedding $ V_\pi \hookrightarrow H^{\sigma, -\lambda} $ is of ``Schmid type''.

    \underline{Case $ G = \Spin(1, 2m) $:} \\
    By Proposition 6.2.7 of \cite[p.140]{Collingwood}, each discrete series representation admits only one embedding and the corresponding principal series ``decomposes according to Schmid's character identity'' by Theorem 5.2.4 of \cite[p.114]{Collingwood}. (This also holds for $ m = 1 $, i.e. $ G \simeq \SL(2, \RR) $, compare \cite[Chapter 5]{Wallach} or \cite[Remark 5.2.5]{Collingwood}).
    It remains to deal with the other cases.

    By Proposition 6.2.8 -- 6.2.10 of \cite[pp.140-141]{Collingwood}, some discrete series representations have more than one embedding. We have to show that for every embedding $ V_\pi \hookrightarrow H^{\sigma_1, -\mu} $ which is not of ``Schmid type'' there is one embedding $ V_\pi \hookrightarrow H^{\sigma_2, -\lambda} $ of ``Schmid type'' with $ 0 < \lambda < \mu $. The last condition has to be checked for every infinitesimal character, not only for the trivial one.

    We want to parametrise principal series representations for fixed (regular integral) infinitesimal character.

    Let $ \tfrak \subset \m $ be the Lie algebra of a maximal torus. Then, $ \a \oplus \tfrak $ is a Cartan subalgebra of $ \g $.

    We consider $ \h := \a \oplus i \tfrak \subset \g_\CC $. Then, $ \h_\CC $ is a Cartan subalgebra of $ \g_\CC $. 

    Let $ \Phi = \Phi(\g_\CC, \h_\CC) \subset \h^* $ be a root system.

    Let $ |\cdot| $ denote the norm on $ \h^*_\CC $ induced by the Killing form.

    Fix $ H_0 \in \a \smallsetminus \{0\} $ such that $ \alpha(H_0) > 0 $ for the short root $ \alpha $ of $ \a $ in $ \n $. Then,
    \[
        \Phi = \{ \mu \in \Phi \mid \mu(H_0) \neq 0 \} \cup \Phi_\m \ ,
    \]
    where $ \Phi_\m := \{ \mu \in \Phi \mid \mu(H_0) = 0 \} $. 
    Choose a positive root system $ \Phi_\m^{+} \subset \Phi_\m $ and set
    \[
        \Phi^+ := \{ \mu \in \Phi \mid \restricted{\mu}{\a} \in \a^*_+ \} \cup \Phi_\m^{+} \subset \Phi
    \]
    (compatible positive root system).

    It follows from Lemma 2.3.5 of \cite[p.58]{Wallach} that the existence of a compact Cartan subgroup of $G$ implies that there is a real root $ \mu_\RR \in \Phi^{+} $, i.e. $ \restricted{ \mu_\RR }{ i\tfrak } = 0 $ and $ \restricted{\mu_\RR}{\a} \in \a^*_+ $. Moreover, it is unique as $G$ is of real rank one.
    Observe that $ \mu_\RR $ is a dominant root. 

    Let $ \Lambda \in \h^* $ be dominant, regular and integral: $ \frac{ 2 \langle \Lambda, \mu \rangle }{ \langle \mu, \mu \rangle } \in \NN $ for all $ \mu \in \Phi^{+} $.

    Let $ W_\m $ be the Weyl group of $ \Phi_\m $ and $ W := W_\g $ be the Weyl group of $ \Phi $. The set of principal series representations with infinitesimal character $ \chi_\Lambda = \chi_F $ is in one-to-one correspondence with the quotient of Weyl groups $ W^\m := W_\m \bs W $ via the map $ [w] \mapsto (\sigma, -\lambda) $ (see \cite[p.48]{Collingwood}), 
    where $ \lambda = \restricted{w\Lambda}{\a} $ and $ \sigma $ is a finite-dimensional irreducible representation of $M$ with infinitesimal character $ \restricted{w\Lambda}{i \tfrak} $. Since $M$ is connected ($ G \neq \SL(2, \RR) $), $ \sigma $ is uniquely determined.

    We are only interested in those $ [w] \in W^\m $ with $ \restricted{ w \Lambda }{ \a } \in \a^*_+ $ ($ \lambda \in \a^*_+ $) --- this condition is in fact independent of $ \Lambda $ as it is equivalent to
    \begin{multline*}
        \langle \restricted{\Lambda}{\a}, \restricted{w^{-1} \mu_\RR}{\a} \rangle > 0
        \iff \langle \alpha, \restricted{w^{-1} \mu_\RR}{\a} \rangle > 0 \\
        \iff \restricted{ w^{-1} \mu_\RR }{\a} \in \a^*_+ \iff w^{-1} \mu_\RR \in \Phi^+ \ .
    \end{multline*} 
    This defines a subset $ W^\m_+ \subset W^\m $.

    By Theorem 2.4 of \cite[p.321]{Procesi}, the stabiliser of a dominant element $ \nu $ is generated by the simple reflections $ s_\mu $ for which $ \langle \mu, \nu \rangle = 0 $.

    Applying this result for $ \nu = \mu_\RR $ (dominant root) yields that $ \Stab_W(\mu_\RR) = W_\m $ since $ \langle \mu, \mu_\RR \rangle = 0 $ ($ \mu \in \Phi $) if and only if $ \mu \in \Phi_\m $.

    Thus, we can identify $ W^\m $ with $ \{ w^{-1} \mu_\RR \mid w \in W \} \subset \Phi $ and $ W^\m_{+} $ with $ \{ w^{-1} \mu_\RR \mid w \in W \} \cap \Phi^{+} $.

    Moreover, $W$ acts transitively (for $ \g $ simple) on roots of equal lengths. We obtain an identification of $ W^\m_{+} $ with $ \{ \mu \in \Phi^+ \mid |\mu| = |\mu_\RR| \} $.

    Let $ [w_1], [w_2] \in W^\m_+ $. We say that $ [w_1] \leq [w_2] $ if $ w_1^{-1} \mu_\RR - w_2^{-1} \mu_\RR $ is a nonnegative linear combination of positive roots.

    Let $ \mu, \tau \in \a^* $. Then, we write $ \mu > \tau $ if $ \mu - \tau \in \a^*_+ \iff \mu(H_0) > \tau(H_0) $.

    \begin{lemSec}\label{Lem [w_1] < [w_2] => ...}
        $ [w_1] < [w_2] \Rightarrow w_1 \restricted{ \Lambda }{ \a } > w_2 \restricted{ \Lambda }{ \a } $ for all regular, dominant $ \Lambda $.
    \end{lemSec}

    \begin{proof}
        Let $ [w_1], [w_2] \in W^\m_+ $. Write $ w_1^{-1} \mu_\RR - w_2^{-1} \mu_\RR = \sum_{\mu \in \Phi^+} n_\mu \mu $.
        Then,
        \[
            \langle w_1 \Lambda - w_2 \Lambda , \mu_\RR \rangle = \langle \Lambda, w_1^{-1} \mu_\RR - w_2^{-1} \mu_\RR \rangle
            = \sum_{\mu \in \Phi^+} n_\mu \underbrace{ \langle \Lambda, \mu \rangle }_{ > 0 }
        \]
        If $ [w_1] < [w_2] $, then $ n_\mu \geq 0 $ and they are not all zero. Thus,
        \[
            \langle w_1 \Lambda - w_2 \Lambda , \mu_\RR \rangle > 0 \ .
        \]
    \end{proof}

    \needspace{3\baselineskip}
    Let us now discuss the cases
    \begin{enumerate}
    \item \underline{$ G = \SU(1, n) $:} \\
        Write $ \h^* = \{ \sum_{i=0}^n x_i e_i \mid \sum_{i=0}^n x_i = 0 \} $, $ \Phi^+ = \{ e_i - e_j \mid 0 \leq i < j \leq n \} $ and $ \mu_\RR = e_0 - e_n $ (compare with \cite[p.14]{Collingwood}).

        We parametrise $ W^\m $ by $ \{ w_{ij} \mid 0 \leq i, j \leq n - 1 , i + j \leq n - 1 \} $, where the $ w_{ij} $'s are determined by $ w_{ij}^{-1} \mu_\RR = e_i - e_{n-j} $.

        Here $ w_{ij} $ corresponds to $ \gamma_{i, j+1} $ in \cite[p.83]{Collingwood} ($ \gamma_{i, j+1} = w_{ij} \rho_\g $, $ \rho_\g $ is the half-sum over $ \Phi^{+} $).

        As $ e_a - e_b \in \Phi^{+} \cup \{0\} $ for all $ a, b \in \NN_0 $ such that $ 0 \leq a \leq b \leq n $, $ i \leq k $ implies that $ e_i - e_k \in \Phi^{+} \cup \{0\} $ and $ j \leq l  $ implies that $ e_{n-l} - e_{n-j} \in \Phi^{+} \cup \{0\} $. It follows that
        \[
            i \leq k, j \leq l \, \Rightarrow \, w_{ij} \leq w_{kl} \qquad (*) \ .
        \]
        By Theorem 5.3.1 of \cite[p.115]{Collingwood}, the Schmid type principal series are parametrised by $ w_{ij} $ with $ i + j = n - 1 $.

        By Proposition 6.2.8 of \cite[p.140]{Collingwood}, the discrete series representations with infinitesimal character $ \chi_\Lambda $, indexed by $ \pi_0, \pi_1, \dotsc, \pi_n $ embed as follows
        \begin{enumerate}
        \item $ \pi_0 \to w_{0, n-1} $ (of Schmid type)
        \item $ \pi_n \to w_{n-1, 0} $ (of Schmid type)
        \item $ \pi_l $ ($ 0 < l < n $) $ \to $ $ w_{l, n-1-l} $ (of Schmid type), $ w_{l-1, n-l} $ (of Schmid type), $ w_{l-1, n-1-l} $ (not of Schmid type)
        \end{enumerate}
        By $ (*) $ and the above lemma, we are done.
    \item \underline{$ G = \Sp(1, n) $:} \\
        Write $ \h^* = \{ \sum_{i=0}^n x_i e_i \mid x_i \in \RR \} $,
        \[
            \Phi^+ = \underset{\text{short roots}}{\{ e_i \pm e_j \mid 0 \leq i < j \leq n \}} \cup \underset{\text{long roots}}{\{ 2 e_i \mid i = 0, \dotsc, n \}}
        \]
        and $ \mu_\RR = e_0 + e_1 $ (short root). Compare with \cite[p.14]{Collingwood}.

        We parametrise $ W^\m $ by $ \{ w_{ij} \mid 0 \leq i \leq j , i + j \leq 2 n - 1\} $, where the $ w_{ij} $'s are determined by $ w_{ij}^{-1} \mu_\RR = \begin{cases}
            e_i + e_{j+1}   & \text{if } j \leq n - 1 \\
            e_i - e_{2n -j} & \text{if } j \geq n
        \end{cases} $.
        \\ Here $ w_{ij} $ corresponds to $ \gamma_{i, j+1} $ in \cite[p.85]{Collingwood}.

        For all $ a, b \in \NN_0 $ such that $ 0 \leq a \leq b \leq n $, $ e_a \pm e_b \in \Phi^{+} \cup \{0\} $. It follows that
        \[
            i \leq k, j \leq l \, \Rightarrow \, w_{ij} \leq w_{kl} \qquad (**) \ .
        \]
        By Theorem 5.4.1 of \cite[p.120]{Collingwood}, the Schmid type principal series arise if and only if $ i + j = 2n - 1 $.

        The discrete series are parametrised by $ \pi_0, \dotsc, \pi_n $.
        By Proposition 6.2.9 of \cite[p.141]{Collingwood}, we have
        \begin{enumerate}
        \item $ \pi_0 \to w_{0, 2n-1} $ (of Schmid type)
        \item $ \pi_n \to w_{n-1, n} $ (of Schmid type), $ w_{n-2, n-2} $ (not of Schmid type)
        \item $ \pi_l $ ($ 0 < l < n $) $ \to $ $ w_{l, 2n-1-l} $ (of Schmid type), $ w_{l-1, 2n-l} $ (of Schmid type), $ w_{l-1, 2n-1-l} $ (not of Schmid type).
        \end{enumerate}

        The assertion follows by $ (**) $ and Lemma \ref{Lem [w_1] < [w_2] => ...}.
    \item \underline{$ G = F_4 $:} \\
        Write $ \h^* = \{ \sum_{i=0}^3 x_i e_i \mid x_i \in \RR \} $,
        \begin{multline*}
            \Phi^+ = \underset{\text{long roots}}{\{ e_i \pm e_j \mid 0 \leq i < j \leq 3 \}} \cup \underset{\text{short roots}}{\{ e_i \mid i = 0, \dotsc, 3 \}} \\
                \cup \underset{\text{short roots}}{\{ \tfrac12(e_0 \pm e_1 \pm e_2 \pm e_3) \mid 0 \leq i < j \leq 3 \}}
        \end{multline*}
        and $ \mu_\RR = e_0 $ (short root). Compare with \cite[p.14]{Collingwood}.
        The $ \pm $ signs are taken here independent.

        Thus, $ W^\m_+ $ is in bijection with the short positive roots:
        \[
            W^\m_+ \simeq \{ w^{-1} \mu_\RR \mid w \in W^\m_+ \} \simeq \{ \text{short positive roots} \} \ .
        \]
        We want to compare this parametrisation with the one in \cite[p.87]{Collingwood}:
        \[
            \{ w \rho_\g \mid [w] \in W^\m_+ \} = \{ \gamma_{ij} \} \ .
        \]
        Define a map from $ W^\m_+ $ to $ \h^* $ by
        \[
            [w] \mapsto w \rho_\g \ ,
        \]
        where $ w $ is chosen such that $ \restricted{w\rho_\g}{i \tfrak} $ is dominant with respect to $ \Phi_\m^+ $.

        We have $ \rho_\g = \frac12 (11, 5, 3, 1) \in \h^* \equiv \RR \oplus \RR^3 $. For every short root $ \mu = w^{-1} \mu_\RR $, we compute
        \[
            \langle \mu, 2 \rho_\g \rangle = \langle \mu_\RR, 2 w \rho_\g \rangle = \langle e_0, 2 w \rho_\g \rangle \ .
        \]
        \needspace{3\baselineskip}
        Using this and the values of $ \langle e_0, 2 \gamma_{ij} \rangle $ provided by \cite[p.87]{Collingwood}, we are able to construct most of the following tables:
        \begin{center}
            \begin{tabular}{c|c|c}
                $ \mu $                     & $ \langle \mu, 2 \rho_\g \rangle $        & $ \gamma_{ij} $   \\ \hline
                $ e_0 $                     & $ 11 $                                    & $ \gamma_{01} $ \tikzmark{e_0}  \\
                \tikzmark{e_1L}$ e_1 $      & $ 5 $                                     & $ \gamma_{24} $ \tikzmark{e_1}  \\
                \tikzmark{e_2L}$ e_2 $      & $ 3 $                                     & $ \gamma_{25} $ \tikzmark{e_2}  \\
                \tikzmark{e_3L}$ e_3 $      & $ 1 $                                     & $ \gamma_{26} $ \tikzmark{e_3}
            \end{tabular}
            \hspace{1.5cm}
            \begin{tabular}{c|c|c}
                $ \mu $                                                                          & $ \langle \mu, 2 \rho_\g \rangle $        & $ \gamma_{ij} $   \\ \hline
                \tikzmark{gamma_02L}$ \tfrac12(e_0 + e_1 + e_2 + e_3) $ \tikzmark{gamma_02R}        & $ 10 $                                       & $ \gamma_{02} $   \\
                $ \tfrac12(e_0 + e_1 + e_2 - e_3) $ \tikzmark{gamma_03R}                            & $ 9 $                                        & $ \gamma_{03} $   \\
                $ \tfrac12(e_0 + e_1 - e_2 + e_3) $ \tikzmark{gamma_04R}                            & $ 7 $                                        & $ \gamma_{04} $\tikzmark{gamma_04} \\
                \tikzmark{gamma_14L}$ \tfrac12(e_0 + e_1 - e_2 - e_3) $ \tikzmark{gamma_14R}        & $ 6 $                                        & $ \gamma_{14} $\tikzmark{gamma_14} \\
                $ \tfrac12(e_0 - e_1 + e_2 + e_3) $ \tikzmark{gamma_05R}                            & $ 5 $                                        & $ \gamma_{05} $\tikzmark{gamma_05} \\
                \tikzmark{gamma_15L}$ \tfrac12(e_0 - e_1 + e_2 - e_3) $ \tikzmark{gamma_15R}        & $ 4 $                                        & $ \gamma_{15} $\tikzmark{gamma_15} \\
                \tikzmark{gamma_16L}$ \tfrac12(e_0 - e_1 - e_2 + e_3) $ \tikzmark{gamma_16R}        & $ 2 $                                        & $ \gamma_{16} $   \\
                $ \tfrac12(e_0 - e_1 - e_2 - e_3)$ \tikzmark{gamma_17R}                             & $ 1 $                                        & $ \gamma_{17} $
            \end{tabular}
            \begin{tikzpicture}[overlay, remember picture, yshift=.25\baselineskip, shorten >=.5pt, shorten <=.5pt]
                \draw [->] ([xshift=2pt, yshift=3pt]{pic cs:gamma_04}) [bend left] to ([xshift=2pt, yshift=3pt]{pic cs:gamma_05});
                \draw [->] ([xshift=2pt, yshift=3pt]{pic cs:gamma_14}) [bend left] to ([xshift=2pt, yshift=3pt]{pic cs:gamma_15});
                \draw [->] ([xshift=-2pt, yshift=3pt]{pic cs:gamma_14L}) -- ([yshift=3pt]{pic cs:e_1});
                \draw [->] ([xshift=-2pt, yshift=3pt]{pic cs:gamma_15L}) -- ([yshift=3pt]{pic cs:e_2});
                \draw [->] ([xshift=-2pt, yshift=3pt]{pic cs:gamma_16L}) -- ([yshift=3pt]{pic cs:e_3});
                \draw [->] ([yshift=3pt]{pic cs:e_0}) -- ([xshift=-2pt, yshift=3pt]{pic cs:gamma_02L});
                \draw [->] ([yshift=3pt]{pic cs:gamma_02R}) [bend left] to ([yshift=3pt]{pic cs:gamma_03R});
                \draw [->] ([yshift=3pt]{pic cs:gamma_03R}) [bend left] to ([yshift=3pt]{pic cs:gamma_04R});
                \draw [->] ([yshift=3pt]{pic cs:gamma_04R}) [bend left] to ([yshift=3pt]{pic cs:gamma_14R});
                \draw [->] ([yshift=3pt]{pic cs:gamma_05R}) [bend left] to ([yshift=3pt]{pic cs:gamma_15R});
                \draw [->] ([yshift=3pt]{pic cs:gamma_15R}) [bend left] to ([yshift=3pt]{pic cs:gamma_16R});
                \draw [->] ([yshift=3pt]{pic cs:gamma_16R}) [bend left] to ([yshift=3pt]{pic cs:gamma_17R});
                \draw [->] ([xshift=-2pt, yshift=3pt]{pic cs:e_1L}) [bend right] to ([xshift=-2pt, yshift=3pt]{pic cs:e_2L});
                \draw [->] ([xshift=-2pt, yshift=3pt]{pic cs:e_2L}) [bend right] to ([xshift=-2pt, yshift=3pt]{pic cs:e_3L});
            \end{tikzpicture}
        \end{center}
        Only two values ($ 1 $ and $ 5 $) of $ \langle \mu, 2 \rho_\g \rangle $ appear twice. For these we need in an additional argument.
        Since we have in the left table precisely those roots for which $w$ can be represented by a permutation, we distinguish roots for the values $ 1 $ and $ 5 $ by looking at $ w $.
        
        Again one checks (check the arrows above) that
        \[
            i \leq k, j \leq l \Longrightarrow \gamma_{ij} \leq \gamma_{kl} \qquad (***) \ .
        \]
        In fact not all of these implications are actually needed.

        By Theorem 5.5.12 of \cite[pp.130-131]{Collingwood}, $ \gamma_{26} $ and $ \gamma_{17} $ are the Schmid type principal series.

        Denote the discrete series representations by $ \pi_1, \pi_2, \pi_3 $. By Proposition 6.2.9 of \cite[p.141]{Collingwood}, we have
        \begin{enumerate}
        \item $ \pi_1 \to \gamma_{26} $ (of Schmid type), $ \gamma_{05} $ (not of Schmid type)
        \item $ \pi_2 \to \gamma_{26} $ (of Schmid type), $ \gamma_{17} $ (of Schmid type), $ \gamma_{16} $ (not of Schmid type)
        \item $ \pi_3 \to \gamma_{17} $ (of Schmid type), $ \gamma_{01} $ (not of Schmid type).
        \end{enumerate}
        The assertion follows by $ (***) $ and Lemma \ref{Lem [w_1] < [w_2] => ...}.
    \end{enumerate}
\end{proof}

\end{appendices}

\cleardoublepage
\phantomsection

\printindex[n] 

\addcontentsline{toc}{section}{Index of notation}

\cleardoublepage
\phantomsection

\printindex[idx] 

\addcontentsline{toc}{section}{Index of terminology}

\cleardoublepage
\phantomsection

\eject\addcontentsline{toc}{section}{Bibliography}

\bibliographystyle{plain}

\end{document}